\documentclass[11pt]{article}
\usepackage[latin1]{inputenc} 
\usepackage{latexsym,amsmath,amsfonts,euscript} 
\usepackage{amssymb}
\usepackage{amsthm}

\newtheorem{Th}{Theorem} 

\newtheorem{Theorem}{Theorem}[section] 
\newtheorem{Prop}[Theorem]{Proposition} 
\newtheorem{Def}[Theorem]{Definition} 
\newtheorem{Lemma}[Theorem]{Lemma} 
\newtheorem{Coro}[Theorem]{Corollary} 
\newtheorem{Remark}[Theorem]{Remark} 
 
\begin{document}

\newcommand{\finishproof}{\hfill $\Box$ \vspace{5mm}} 
\newcommand{\cI}{{\cal I}} 
\newcommand{\al}{\alpha} 
\newcommand{\be}{\beta} 
\newcommand{\ga}{\gamma} 
\newcommand{\Ga}{\Gamma} 
\newcommand{\ep}{\epsilon} 
\newcommand{\de}{\delta} 
\newcommand{\De}{\Delta} 
\newcommand{\ka}{\kappa} 
\newcommand{\la}{\lambda} 
\newcommand{\te}{\theta} 
\newcommand{\om}{\omega} 
\newcommand{\si}{\sigma} 
\newcommand{\bg}{{\bar g}} 
\newcommand{\bs}{{\bar s}} 
\newcommand{\bL}{{\bar L}} 
\newcommand{\tL}{{\tilde L}}
\newcommand{\tK}{{\tilde K}}
\newcommand{\tg}{\tilde g} 
\newcommand{\eqdef}{\stackrel{\rm def}{=}} 
\newcommand{\const}{\mathop{\rm const}} 
\newcommand{\sg}{\mathop{\rm sgrad}} 
\newcommand{\id}{{\bf 1}} 
\newcommand{\dds}{\frac{d}{ds}|_{s=0}} 
\newcommand{\n}{{\|\!|}} 
\newcommand{\p}{{\partial}}
\newcommand{\C}{\mathbb C}
\newcommand{\LL}{\mathbb L}

\newcommand{\A}{{\mathbb A}}
\newcommand{\B}{\mathbb B}
\newcommand{\D}{\mathbb D}
\newcommand{\Z}{\mathbb Z} 
\newcommand{\N}{\mathbb N} 
\newcommand{\R}{\mathbb R}
\newcommand{\DD}{\mathbb D} 
\newcommand{\Q}{\mathbb Q} 
\newcommand{\T}{\mathbb T}
\def\1{{1\hspace{-.23cm}1}}

\newcommand{\CC}{{\mathbb C}}
\newcommand{\tI}{\tilde I}

\renewcommand{\mod}{\mathop{\rm mod}} 
\renewcommand{\theequation}{\thesection .\arabic{equation}} 
\renewcommand{\arraystretch}{1.3}

\renewcommand{\theequation}{\thesection .\arabic{equation}}  
\renewcommand{\arraystretch}{1.3}

\title{From KAM Tori to Isospectral Invariants and Spectral Rigidity of Billiard Tables. } 
\author{G.Popov, P.Topalov} 
\maketitle 

\begin{abstract} 
\noindent
This article is a part of a project  investigating the relationship between the dynamics of completely integrable or ``close'' to completely integrable billiard tables, the integral geometry on them, and the spectrum of the corresponding Laplace-Beltrami operators. 
It is concerned with new isospectral invariants and with the spectral rigidity problem for the Laplace-Beltrami operators $\Delta_t$, $t\in [0,1]$,  with Dirichlet, Neumann or Robin boundary conditions, associated with  $C^1$ families  of billiard tables $(X,g_t)$. 
We introduce a notion of weak isospectrality for such deformations. 

The main dynamical assumption on  $(X,g_0)$  is that the corresponding billiard ball map $B_0$ or an iterate $P_0= B_0^m$  of it  posses a Kronecker invariant torus  with a Diophantine frequency  $\omega_0$ and that the corresponding Birkhoff Normal Form  is nondegenerate in  Kolmogorov sense. 
Then we prove that there exists $\delta_0>0$ and  a set $\Xi$ of Diophantine frequencies containing $\omega_0$ and of full Lebesgue measure around  $\omega_0$   such that for each $\omega\in \Xi$ and $0<\delta<\delta_0$ there exists a $C^1$ family of Kronecker tori $\Lambda_t(\omega)$ of $P_t$ for  $t\in [0,\delta]$. If the family $\Delta_t$, $t\in [0,1]$,  satisfies the weak isospectral condition we prove that the average action $\beta_t(\omega)$ on  $\Lambda_t(\omega)$ and  the Birkhoff Normal Form of $P_t$ at  $\Lambda_t(\omega)$ are  independent of $t\in [0,\delta]$  for each $\omega\in \Xi$. 

As an application we obtain  infinitesimal spectral rigidity for Liouville billiard tables in dimensions 2 and 3. In particular  infinitesimal spectral rigidity for the ellipse and the ellipsoid is obtained under the weak isospectral condition. Applications are obtained also for strictly convex billiard tables in $\R^2$ as well as in the case when $(X,g_0)$ admits an elliptic periodic billiard trajectory with no resonances of order $\le 4$. 

In particular we obtain spectral rigidity (under the weak isospectral condition) of  elliptical billiard tables in the class of analytic and $\Z_2\times \Z_2$ symmetric billiard tables in $\R^2$. We prove also that billiard tables with boundaries close to ellipses   are spectrally rigid in this class. 

The results are based on a construction of  $C^1$ families of quasi-modes associated with the Kronecker tori $\Lambda_t(\omega)$ and on suitable KAM  theorems for $C^1$ families of Hamiltonians. We propose a new iteration schema (a modified iterative lemma) in the proof of the KAM theorem with parameters, which avoids the Whitney extension theorem for $C^\infty$ jets and allows one to obtain global  estimates  of  the corresponding canonical transformations and Hamiltonians in the scale of all H\"older norms. The classical and quantum Birkhoff Normal Forms for $C^1$ or analytic families of symplectic mappings (Hamiltonians) obtained here can be used as well in order to investigate  problems related to the quantum non-ergodicity of  $C^\infty$-smooth KAM systems.

\end{abstract} 

\pagebreak

\tableofcontents

 !TEX root = iso-inv-main.tex

\section{Introduction} \label{Sec:introduction}

This article is a part of a project (cf. \cite{PT1}-\cite{PT4}) investigating the relationship between the dynamics of completely integrable or ``close'' to completely integrable billiard tables, the integral geometry on them, and the spectrum of the corresponding Laplace-Beltrami operators.
It is concerned with new isospectral invariants and the spectral rigidity of the Laplace-Beltrami operator associated with $C^1$ deformations  $(X,g_t)$, $0\le t\le 1$, of a billiard table $(X,g)$, where  $X$ is a $C^\infty$ smooth compact manifold with a connected boundary $\Gamma:=\partial X$ of dimension ${\rm dim}\,  X=n\ge 2$ and 
$t\to g_t$ is a $C^1$    family of smooth  Riemannian metric  on $X$.

Substantial progress in the inverse spectral geometry has been made 
by means of the  wave-trace formula 
 \cite{G}, \cite{G-M1}, \cite{IS, ISZ}, \cite{MM}, \cite{Z1}-\cite{Z5}, and by its semi-classical analogue - the Gutzwiller trace formula \cite{Mein}, \cite{IS, ISZ}, \cite{GPU, GU}. The wave-trace formula, known in  physics as the Balian-Bloch formula  and treated rigorously by Y. Colin de Verdi\`ere \cite{CV1}, J. Duistermaat and V. Guillemin \cite{D-G}, V. Guillemin and R. Melrose \cite{G-M2} and S. Zelditch \cite{Z4}  (see also \cite{CV3}, \cite{MM}, \cite{P1}, \cite{S-V}, \cite{SZ}), as well as the Gutzwiller trace formula  relate the spectrum of the operator with certain invariants of the corresponding closed geodesics  such as their  lengths and the spectrum of the linear Poincaré map. 
  
 It has been proved  in \cite{G}, \cite{Z1, Z2}, \cite{IS, ISZ}, that for certain nondegenerate closed geodesics one can extract the Birkhoff Normal Form (BNF) from the singularity expansions of the wave-trace. S. Zelditch \cite{Z3} - \cite{Z6} and H. Hezari and S. Zelditch \cite{H-Z} have   reconstructed the boundary for a large class of  analytic  domains on $\R^n$ having certain symmetries. 
Hezari and Zelditch \cite{H-Z2}  have proven infinitesimal rigidity of isospectral deformations of the ellipse. 

Spectral rigidity of   closed Riemannian manifolds of negative sectional curvature has been obtained by V. Guillemin and D. Kazhdan \cite{G-K} (in dimension two), C. Croke and V. Sharafutdinov \cite{C-S} (in any dimension) and by G.  Paternain, M. Salo, and G.  Uhlmann \cite{P-S-U} for closed oriented Anosov surfaces.   In order to  link the spectrum of the Laplace-Beltrami operator with the length spectrum of the manifold   the wave-trace formula is used. The wave-trace  formula is especially useful for $C^1$-deformations $(X,g_t)$ of a closed  Riemannian manifold  $(X,g_0)$ with an Anosov  geodesic flow  since every closed geodesic of $(X,g_0)$ is hyperbolic, hence, nondegenerate and it gives rise to a $C^1$
family of closed hyperbolic geodesics of $(X,g_t)$ for $|t|$ small enough. This reduces the problem of the infinitesimal spectral rigidity of Anosov manifolds to the injectivity of a  geodesic ray transform which has been proved for negatively curved closed manifolds of any dimension  \cite{C-S} and for closed oriented Anosov surfaces \cite{P-S-U}. 
Moreover,   infinitesimal rigidity  implies spectral rigidity because of the structural stability. 
Non of these properties is valid for deformations of a billiard table ``close'' to  an integrable billiard table which makes the spectral rigidity problem much more difficult in that case. The wave-trace method requires  certain technical assumptions such as simplicity of the length spectrum (a non-coincidence condition) and non-degeneracy of the corresponding closed geodesic and its iterates which are not fulfilled in general. 

The main dynamical assumption on  $(X,g_0)$ in the present work is that the corresponding billiard ball map or an iterate of it  posses a Kronecker invariant torus (see Definition \ref{Def:Kronecker}) with a Diophantine frequency vector and that the corresponding Birkhoff Normal Form (BNF) is nondegenerate in a Kolmogorov sense. Such Hamiltonian systems are said to be of Kolmogorov-Arnold-Moser (KAM) type. 
The dynamics of such systems  is quite complex. In particular, the non-coincidence  and the non-degeneracy conditions  may not hold for the corresponding closed geodesics.  This makes the wave-trace method useless for such systems in general.  
 On the other hand, the Kronecker tori with Diophantine frequencies  survive under small perturbations which makes them the right objects to look for. For this reason we propose another method which is based on the construction of {\em $C^1$-families of quasi-modes} associated with these tori.

Du to the Kolmogorov-Arnold-Moser (KAM) theory, if the initial Hamiltonian system  ($t=0$) is completely integrable and if it satisfies the Kolmogorov  nondegeneracy condition, then  a large part of the invariant tori of the initial system having Diophantine frequencies $\omega$ survive under the perturbation for $t$ in a small interval  $[0,\delta_0)$  and give rise to cylinders of invariant tori $t\to \Lambda_t(\omega)$, $0\le t < \delta_0$ along the perturbation.
The positive number $\delta_0$ depends on the small constant $\kappa$ and on the exponent $\tau$ in the Diophantine condition \eqref{eq:sdc}.
 The aim of this paper is to prove  that  the invariant tori form $C^1$-families with respect to $t$ and that the value at $\omega$ of the corresponding  Mather's  $\beta$-function does not depend on $t$, or equivalently that 
the Birkhoff Normal Form (BNF) of the system at each torus $\Lambda_t(\omega)$ does not depend on $t$ for any $C^1$-smooth isospectral deformation. 
Applications will be obtain in the following three cases: for deformations of Liouville billiard tables, in the case of deformations of strictly convex domains and in the case  when $g_0$ admits an elliptic (broken) geodesic which has no resonances of order $\le 4$  and has a nondegenerate BNF.

Let us formulate the main problems that we are going to investigate. 
Denote by $\Delta_t$ the ``geometric'' Laplace-Beltrami operator corresponding to the Riemannian manifold  $(X,g_t)$ with Dirichlet, Neumann or Robin boundary conditions.  This is a self-adjoint operator in $L^2(X)$ with discrete spectrum accumulating at $+\infty$. The corresponding  eigenvalues $\lambda$   solve the spectral problem  
\begin{equation}  
\left\{  
\begin{array}{rcll}  
\Delta_t\,  u\ &=& \ \la\,  u \, \quad \mbox{in}\ X\, , \\  
\displaystyle {\mathcal B}_t \, u \ &=&\ 0  \, , 
\end{array}  
\right.  
                                         \label{thespectrum}  
\end{equation}
where ${\mathcal B}_tu = u|_{\Gamma}$ in the case of Dirichlet boundary conditions, ${\mathcal B}_t u = \frac{\partial u}{\partial \nu_t}|_{\Gamma}$ in the case of Neumann boundary conditions, and ${\mathcal B}_tu = 
\frac{\partial u}{\partial \nu_t}|_{\Gamma} - f\, u|_{\Gamma}$ in the case of Robin boundary conditions,  where $\nu_t(x)$, $x\in \Gamma_t$, is the outward unit normal to  
$\Gamma$ with respect to the metric $g_t$ and $f$ is a smooth real valued function on $\Gamma$. 

The  method we use is based on the construction of  $C^1$ smooth with respect to $t$ quasi-modes. This method has been applied in \cite{PT4} in order to investigate the spectral rigidity of the problem \eqref{thespectrum} with Robin boundary conditions in the case when the metric $g$ is fixed and $t\to f_t$ is a continuous deformation of the function appearing  in the Robin boundary condition. Let us formulate  the isospectral condition. 

Consider a union ${\mathcal I}$ 
of infinitely many disjoint intervals $[a_k,b_k]$ going to infinity, of  length $o\left(\sqrt{a_k}\right)$, and which are polynomially separated. 

More precisely, fix  two positive constants $d\ge 0$ and $c>0$, and suppose that
\begin{enumerate}
\item[$(\mbox{H}_1)$] {\em  $\mathcal{I}\subset (0,\infty) $ 
 is a union of infinitely many disjoint  intervals $[a_k,b_k]$, $k\in\N$,  such that 
\begin{enumerate}
\item[$\bullet$]
$\displaystyle{\lim\,  a_k\, =\, \lim\,  b_k\,  =\,  +\infty}$;
\item[$\bullet$] 
$\displaystyle\lim\,   \frac{b_k -a_k}{\sqrt{a_k}}\,  =\,  0$; 
\item[$\bullet$]
$a_{k+1} - b_{k}\, \ge \, c b_k^{-d}$ \quad for  any 
$k\in \N$. 
\end{enumerate}}
\end{enumerate}
Given a set ${\mathcal I}$ satisfying $(\mbox{H}_1)$, 
we impose the following {\em ``weak isospectral assumption''} 
\begin{enumerate}
\item[$(\mbox{H}_2)$]
{\em There is  $a\ge 1$ such that  \quad 
${\rm Spec}\left(\Delta_t
\right)\, \cap [a,+\infty)\  
\subset \   {\mathcal I}\quad  \forall\,  t \in [0,1]   \, .$}
\end{enumerate}
Note that the length of the intervals $[a_k,b_k]$ can increase and even go  to infinity as $k\to \infty$ but not faster than  $o\left(\sqrt{a_k}\right)$. Physically this means that we allow noise in the system. 
Using the asymptotic behavior  of the eigenvalues $\lambda_j$ as $j\to \infty$ one can  show that  conditions (H$_1$)-(H$_2$)
are ``natural'' for any $d> n/2$ and $c>0$. By ``natural'' we mean that for any $d> n/2$ and $c>0$ the usual isospectral condition 
\[
 {\rm Spec}\left(\Delta_t \right)\, = \, 
{\rm Spec}\left(\Delta_0 \right) \quad \forall\,  t \in [0,1] 
\]
implies that  
there exists $a\ge 1$ and  a family of infinitely many disjoint intervals $[a_k,b_k]$ such that (H$_1$)-(H$_2$) are satisfied -- see \cite{PT4}, Lemma 2.2,  
for  details. The exponent $d$ depends on the level spacing of the spectrum of $\Delta_0$. 

The elastic reflection of geodesics of  $(X,g_t)$ at $\Gamma$ determines
continuous curves on $X$ called {\em billiard trajectories}  as well as a discontinuous dynamical system on the corresponding coshere bundle $S_t^\ast X$ --
the ``{\em billiard flow}'' consisting of broken bicharacteristics of the Hamiltonian $h_t$ associated to $g_t$  via the Legendre transform. The latter 
 induces  a discrete dynamical system $B_t$ defined  on an open subset $\widetilde {\bf B}_t^\ast \Gamma$ (depending on $t$) of the open  coball bundle ${\bf B}_t^\ast \Gamma$ of 
$\Gamma$ called billiard ball map (see Section \ref{subsec:billiard-ball} for a definition). 
The  map $B_t:\widetilde {\bf B}_t^\ast \Gamma\to {\bf B}_t^\ast \Gamma$  is exact symplectic.  
Fix an integer $m\ge 1$ and consider the exact symplectic map  
\[
P_t=B_t^m: U_t\to  {\bf B}_t^\ast \Gamma
\]
where $U_t$ is an open subset of ${\bf B}_t^\ast \Gamma_t$ such that $B_t^j(U_t) \subset \widetilde {\bf B}_t^\ast \Gamma$ for any $0\le j<m$. Given an interval $J\subset [0,1]$ we say that $P_t=B_t^m$, $t\in J$, is a $C^1$ family of exact symplectic maps if for every $t_0\in J$ and $\rho^0\in U_{t_0}$ there exist neighborhoods $J_0\subset J$ of $t_0$ and $V\subset U_{t_0}$ of $\rho^0$ such that $V\subset U_{t}$ for every $t\in J_0$ and the map $J_0 \ni t\mapsto P_t\big|_V \in C^\infty(V, T^\ast \Gamma)$ is $C^1$.

We are interested in Kronecker  invariant 
tori of $P_t$  of  Diophantine frequencies, which are defined as follows. 

We denote  by $\T^{d}$ the torus $\T^{d}:=\R^{d}/2\pi \Z^{d}$ of dimension $d\ge 1$ and by ${\rm pr\,}:\, \R^d \to \T^d$ the canonical projection and we consider  $\T^{d}$ as a $\Z$-module. A ``distance" from a given   $\alpha\in\T:=\T^1$ to $0$ can be defined by
\[
|\alpha|_\T\, :=\, \inf \{|a|:\, a\in {\rm pr\,}^{-1}(\alpha)\}.
\]
Fix $\kappa \in (0,1)$ and $\tau > n-1$ and denote by  $D(\kappa,\tau)$ the set all  $\omega\in\T^{n-1}=\R^{n-1}/2\pi \Z^{n-1}$ satisfying the ``strong''  $(\kappa,\tau)$-Diophantine  condition 
\begin{equation}                    
|\langle k,\omega\rangle|_\T = | k_1\omega_1 + \cdots +k_{n-1}\omega_{n-1}|_\T  \ \ge \ 
\frac{\kappa}{\big(\sum_{j=1}^{n-1}|k_j|\big)^{\tau}}\quad \forall\, k=(k_1,\ldots,k_{n-1})\in \Z^{n-1}\setminus \{0\}.
                       \label{eq:sdc}                                        
\end{equation}
The condition \eqref{eq:sdc} is equivalent to the following one.  There exists $\widetilde\omega^\prime \in {\rm pr\,}^{-1}(\omega)$ such that 
\[                   
|\langle\widetilde\omega^\prime,k\rangle + 2\pi k_{n}|\ \ge \ \frac{\kappa}{\big(\sum_{j=1}^{n-1}|k_j|\big)^{\tau}}\quad \forall\, (k,k_{n})\in \Z^{n-1}\times\Z,\ 
k\neq 0. 
\]
Obviously, if this condition is satisfied for one $\widetilde\omega^\prime \in {\rm pr\,}^{-1}(\omega)$ then it holds for each $\widetilde\omega^\prime \in {\rm pr\,}^{-1}(\omega)$. 
We denote by  $\widetilde D(\kappa,\tau)$ the set all  $\widetilde\omega\in\R^{n}$ satisfying the following  ``weak''   $(\kappa,\tau)$-Diophantine  condition: 
\begin{equation}   
|\langle\widetilde\omega,k\rangle|\ \ge \ \frac{\kappa}{\big(\sum_{j=1}^{n}|k_j|\big)^{\tau}}\quad \forall\, k\in \Z^{n},\ 
k\neq 0. 
                       \label{eq:sdc-1-intro}                                        
\end{equation}
Thus the relation $\omega\in D(\kappa,\tau)$ implies that $\widetilde \omega:= (\widetilde \omega^\prime,2\pi) \in \widetilde D(\kappa,\tau)$ for at least one (and then for all) $\widetilde\omega^\prime \in {\rm pr\,}^{-1}(\omega)$. 

The set $D(\kappa,\tau)$ ($\widetilde D(\kappa,\tau)$) is closed and nowhere dense in $\R^{n-1}$ ($\R^n$), respectively. Moreover, 
the union $\cup_{0<\kappa\le 1} D(\kappa,\tau)$ of  $(\kappa,\tau)$-Diophantine vectors is of full Lebesgue measure in $\R^{n-1}$ for $\tau>n-1$ fixed, and $D(\kappa',\tau) \subset D(\kappa,\tau)$ for $0<\kappa < \kappa'$. Denote by $D^0(\kappa,\tau)$ the set of \emph{points of positive Lebesgue density} in $D(\kappa,\tau)$, i.e.  $\omega\in D^0(\kappa,\tau)$  if the Lebesgue measure ${\rm meas\, }(D(\kappa,\tau)\cap V) > 0$ for any neighborhood $V$ of $\omega_0$ in $\R^{n-1}$. By definition, the complement of $D^0(\kappa,\tau)$ in $D(\kappa,\tau)$ is of zero Lebesgue measure. In the same way we define the subset $\widetilde D^0(\kappa,\tau)$ of points of positive Lebesgue density in $\widetilde D(\kappa,\tau)$. 
\begin{Def}\label{Def:Kronecker} 
A Kronecker torus of $P_t$ of a frequency $\omega$ is an embedded submanifold $\Lambda_t(\omega)$ of $ {\bf B}_t^\ast \Gamma$  diffeomorphic
to $\T^{n-1}$  such that 
\begin{enumerate}
	\item[(i)]  $B_t^j(\Lambda_t(\omega))$ is a subset of  $\widetilde {\bf B}_t^\ast \Gamma$ for each  $0\le j\le m-1$; 
	\item[(ii)] $\Lambda_t(\omega)$ is   invariant with respect to $P_t=B_t^m$;
	\item[(iii)] The restriction of $P_t$ to $\Lambda_t(\omega)$ is
$C^\infty$ conjugated to the  translation   $R_{ \omega}: \T^{n-1}\to \T^{n-1}$ given by
$R_{ \omega}(\varphi) =\varphi + \omega $. 
\end{enumerate}
\end{Def} 
This means that there is a smooth  embedding  $f_{t,\omega}: \T^{n-1}\to {\bf B}_t^\ast \Gamma$ such that  $\Lambda_t(\omega)= f_{t,\omega}(\T^{n-1})$ 
and   the    diagram
\begin{equation}
\label{eq:diagram} 
\displaystyle{\begin{array}{cccl} 
\T^{n-1}&\stackrel{R_{ \omega}}{\longrightarrow}&\T^{n-1}\cr
\downarrow\lefteqn{f_{t,\omega}}& &\downarrow\lefteqn{f_{t,\omega}} \cr
\Lambda_t(\omega)&\stackrel{P_t}{\longrightarrow}&\Lambda_t(\omega)&  
\end{array} }
\end{equation}
is commutative.
\begin{Def}\label{Def:Kronecker-C-1}
By a $C^1$-smooth family of Kronecker tori $\Lambda_t(\omega)$ of $P_t$, $t\in [0,\delta]$, with a  frequency $\omega$   we mean  a $C^1$ family of smooth embeddings  $[0,\delta] \ni t \mapsto f_{t,\omega}\in C^\infty (\T^{n-1}, T^\ast\Gamma)$ satisfying (i)-(iii) of Definition \ref{Def:Kronecker}. 
\end{Def}
For each Diophantine frequency $\omega\in D(\kappa,\tau)$ the embedding  $f_{t,\omega}: \T^{n-1}\to {\bf B}_t^\ast \Gamma$ is a {\em Lagrange embedding} (see \cite{Her}, Sect. I.3.2). We simply say that each Kronecker torus 
$\Lambda_t(\omega)\subset {\bf B}_t^\ast  \Gamma$ is Lagrangian for such frequencies. Note that the map $P_t: \Lambda_t(\omega)\to \Lambda_t(\omega)$ is uniquely  ergodic for $\omega\in D(\kappa,\tau)$, i.e. there  is a unique probability measure $\mu_t$ on $\Lambda_t(\omega)$ which is $P_t$ invariant. Evidently, its pull-back $f_{t,\omega}^\ast (d\mu_t)$  by  the diffeomorphism $f_{t,\omega}$ coincides with 
the Lebesgue-Haar measure $d\theta$ of $\T^{n-1}$. The automorphism  $x\to \frac{x}{2\pi}$ of $\R^{n-1}$ induces an isomorphism of  groups $\jmath:\R^{n-1}/2\pi \Z^{n-1} \to \R^{n-1}/\Z^{n-1}$ assigning to any frequency vector $\omega$ the corresponding rotation vector which will be denoted by $\omega/2\pi$. Hereafter we will deal mainly  with frequency vectors which is motivated by the extensive use of the Fourier analysis.

To any Kronecker torus $\Lambda_t(\omega)$ with a Diophantine  frequency $\omega\in D(\kappa,\tau)$ one can associate three dynamical invariants as follows. 

The first one is 
the {\em average action} on the torus, which corresponds to the Mather's $\beta$-function in the case of twist maps. 
Given $\varrho\in \widetilde {\bf B}_t^\ast \Gamma$ we denote by 
$$
A_t(\varrho):= \int_{\tilde\gamma_t(\varrho)} \xi dx 
$$
the action on the broken bicharacteristic $\tilde\gamma_t(\varrho)$  ``issuing from'' $\varrho_0:=\varrho$ and ``having  endpoint'' at $\varrho_m:= P_t(\varrho)$, where $\xi dx$ is the fundamental one-form on $T^\ast X$. Denote by $X_{h_t}$ the Hamiltonian vector field  where  $h_t$ is  the Legendre transform of the metric tensor $g_t$. 
  The broken bicharacteristic  
$\tilde\gamma_t(\varrho)$ is a disjoint union of   integral curves $\gamma_t(\varrho_j)$  of the Hamiltonian vector field $X_{h_t}$  ``issuing'' from  $\varrho_j:=B_t^j(\varrho)$ and ``ending'' at $\varrho_{j+1}=B_t^{j+1}(\varrho)$ and lying on the coshere bundle $\Sigma_t:=S_t^\ast X=\{h_t=1\}$ (for a more precise definition see Section \ref{subsec:billiard-ball}). The vertices of  $\tilde\gamma_t(\varrho)$ can be identified with  $\varrho_j$, $0\le j\le m$, and we have 
$$
A_t(\varrho) = \sum_{j=0}^{m-1}\int_{\gamma_t(\varrho_j)} \xi dx. 
$$
Notice that $2A_t(\varrho)$ is just the length of the broken geodesic in $(X,g_t)$ obtained by projecting   the broken bicharacteristic $\tilde\gamma_t(\varrho)$  to $X$. In particular, $A_t(\varrho)>0$.
By  Birkhoff's ergodic theorem 
\begin{equation}
\displaystyle \beta_t(\omega):=- 2\lim_{N\to+\infty} \frac{1}{2N} \sum_{k=-N}^{N-1} A_t(P_t^k\varrho) = -2\int_{\Lambda_t(\omega)} A_t d\mu_t <0
\label{eq:beta function}
\end{equation}
does not depend on the choice of $\varrho\in \Lambda_t(\omega)$.  The function $\beta_t$ can be extended as a convex function in the case when $n=2$ and  $P_t$ is a monotone twist map.  It   can be related to the Mather's $\beta$-function $\beta_t^M$ \cite{Sor} (cf. also \cite{Sib}) by the isomorphism $\jmath: \R^{n-1}/2\pi \Z^{n-1} \to \R^{n-1}/\Z^{n-1}$, i.e. $\beta_t=\beta_t^M\circ\jmath$. 

Another invariant of a Kronecker torus $\Lambda_t(\omega)$ with a Diophantine frequency is the {\em Liouville class} on it which is defined as the cohomology class $ [f_{t,\omega}^\ast(\xi dx)]\in H^1(\T^{n-1},\R)$, where $\xi dx$ stands for  the fundamental one-form of $T^\ast \Gamma$  (recall that $f_{t,\omega}: \T^{n-1}\to {\bf B}_t^\ast \Gamma$ is a Lagrange embedding). Let $e_1, \ldots, e_{n-1}$ be the canonical basis of $\R^{n-1}$ and $s\to c_j(s)={\rm pr\, }(se_j)$, $j=1,\ldots,n-1$,  be the corresponding loops on $\T^{n-1}$. Then 
$\gamma_{t,\omega}^j:=f_{t,\omega}\circ c_j$, $j=1,\ldots,n-1$,  provide  a basis of loops  of $H_1(\Lambda_t(\omega),\Z)$. In  the dual basis of $H^1(\T^{n-1},\R)$ we write $[f_{t,\omega}^\ast(\xi dx)]$ as
\begin{equation}
I_t(\omega) = \left( \int_{\gamma_{t,\omega}^1} \xi dx,\cdots,  \int_{\gamma_{t,\omega}^{n-1}} \xi dx\right).
\label{eq:momentum-I}
\end{equation}
The {\em  Birkhoff Normal Form} (BNF) of $P_t$ is another invariant related to a Kronecker torus. 
To each Kronecker torus $\Lambda_t(\omega)$ with a Diophantine frequency $\omega$ one can associate a BNF of $P_t$ as follows. There exist an exact symplectic map $\chi_t$ from a neighborhood of $\T^{n-1}\times \{I_t(\omega)\}$ in $T^\ast \T^{n-1}$  to a neighborhood of $\Lambda_t(\omega) $ in $T^\ast\Gamma$   a smooth function $L_t$ and a map $R_t$ such that $\Lambda_t(\omega)= \chi_t\big(\T^{n-1}\times \{I_t(\omega)\}\big)$  and 
\begin{equation}\label{eq:birkhoff-L}
\big(\chi_t^{-1}\circ P_t \circ \chi_t \big)(\varphi,I) = (\varphi+\nabla L_t(I),I) +R_t(\varphi,I),\quad \partial_I^\alpha R_t(\varphi,I_t(\omega))=0 \ \forall\, \alpha\in \N^{n-1}, 
\end{equation}
(see Sect. \ref{Sec:BBM-BNF}). The BNF of $P_t$ at the torus $\Lambda_t(\omega)$ is said to be {\em nondegenerate} if the Hessian matrix of $L_t$ at $I=I_t(\omega)$ is nondegenerate, i.e. 
\begin{equation}\label{eq:BNF-intro}
\det \partial_I^2 L_t (I_t(\omega))\,  \neq\,  0.
\end{equation}
One can choose $L_t$ so that
\begin{equation}
\beta_t(\omega) + L_t(I_t(\omega)) = \langle \omega,I_t(\omega)\rangle \quad {\rm and}\quad \nabla L_t(I_t(\omega))=\omega
\label{eq:beta-alpha}
\end{equation}
(see Lemma \ref{lemma:beta}). 

Given an interval $J\subset [0,1]$ and a $C^1$ family of Kronecker tori $J\ni t \mapsto \Lambda_t(\omega)$, we say that \eqref{eq:birkhoff-L} provides a $C^1$ family of BNFs in $J$ if $t\to \chi_t$, $t\to L_t$ and $t\to R_t$ are $C^1$ families with values in the corresponding $C^\infty$ spaces (see Definition \ref{Def:BNF}). 

Let  $[0,\bar\delta)\ni t\to P_t$ be a $C^1$ family of exact symplectic maps and $0<\bar\delta\le 1$.
We are interested in the following problems. \\

\noindent
\underline{Problem I.}\quad {\em Let $\Lambda_0(\omega_0)$ be a Kronecker torus  of $P_0$ with a $(\kappa_0,\tau)$-Diophantine frequency $\omega_0\in D(\kappa_0,\tau)$, where $0<\kappa_0\le 1$ and $\tau>n-1$.
Suppose that the BNF of $P_0$ at $\Lambda_0(\omega_0)$ is nondegenerate. Do there exist  $\Xi\subset \T^{n-1}$  and $0<\delta\le \bar\delta$ such that 
\begin{enumerate}
	\item $\omega_0\in \Xi$ and $\Xi$ is a set of Diophantine frequencies of full Lebesgue measure at $\omega_0$ which means that ${\rm meas\,} (\Xi\cap W)={\rm meas\,} (W) + o({\rm meas\,} (W))$ as ${\rm meas\,} (W) \to 0$ for any open neighborhood $W$ of $\omega_0$;
	\item For each $\omega\in \Xi$ there exists a $C^1$ family of Kronecker tori $[0,\delta] \ni t \mapsto \Lambda_t(\omega)$ of $P_t$.  
\end{enumerate}}

\vspace{0.3cm}
A positive answer of this question is given by    Theorem \ref{Th:main1}, item {\em 2}, and Theorem \ref{Theo:soft-BNF}. The set $\Xi$ is of the form (see Section \ref{Sec:BBM-BNF}) 
\begin{equation}\label{eq:the-set-Xi}
\Xi\, =\ \bigcup_{0<\kappa\le \kappa_1}\, \Omega_\kappa^0
\end{equation}
where $\kappa_1\le \kappa_0$,  $\omega_0\in \Omega_\kappa^0$ and the set $\Omega_\kappa^0$ consists only of points of positive Lebesgue density for any  $\kappa$ fixed. 
Moreover, Theorem \ref{Theo:soft-BNF} gives a $C^1$ family of simultaneous BNFs associated with the $C^1$ families of invariant tori
\[
[0,\delta] \ni t \to \Lambda_t(\omega) \quad \forall\, \omega \in \Omega_\kappa^0.
\]
which means that the family of symplectic maps $\chi_t$, $t\in [0,\delta]$, is $C^1$ and for any fixed $t$ the map $\chi_t$ provides   a BNF  \eqref{eq:birkhoff-L} of $P_t$ at $\Lambda_t(\omega) $ for all $\omega \in \Omega_\kappa^0$ at once. 
These families are analytic in $t$ if the map $t\to P_t$ is analytic.  We apply that to  the following three situations \\
\begin{enumerate}
	\item $(X, g_0)$ is a nondegenerate Classical Liouville Billiard table as defined in Section \ref{sec:LBT}. Then the  billiard ball map $P_0=B_0$ is  completely integrable and the  Kolmogorov  non-degeneracy condition is fulfilled. Hence Theorem \ref{Th:main1}, {\em 1}-{\em 2}, and Theorem \ref{Theo:soft-BNF} hold for every invariant torus of Diophantine frequency in this case. 
	\item $(X, g_0)$ has an elliptic closed broken geodesic with $m\ge 2$ vertices of no resonances of order $\le 4$ and with a nondegenerate BNF (see Section \ref{Sec:elliptic}). The corresponding return map $P_0=B_0^m$ has a large family of Kronecker tori with Diophantine frequencies, the BNF of each of them is nondegenerate and one can apply Theorem \ref{Th:main1}, {\em 1}-{\em 2}, and Theorem \ref{Theo:soft-BNF}. 
	\item $(X, g_0)$ is a locally strictly geodesically convex billiard table of dimention two (see Section  \ref{Sec:convex}). Then there exists a large family of Kronecker tori with Diophantine frequencies of the billiard ball map $P_0=B_0$. These invariant circles  accumulate  at the boundary $S_0^\ast \Gamma$ of the coball bundle $\textbf{B}_0^\ast \Gamma$ and give rise of the so  called Lazutkin caustics in the interior of $X$. Close to $S_0^\ast \Gamma$ the  map $P_0=B_0$ is twisted. This implies that the BNF of $P_0$ is nondegenerate at each Kronecker torus sufficiently close to $S_0^\ast \Gamma$  and we can apply Theorem \ref{Th:main1}, {\em 1}-{\em 2}, as well as Theorem \ref{Theo:soft-BNF}. 
\end{enumerate}
From now on we denote by $\Xi$ a set of Diophantine frequencies of the form \eqref{eq:the-set-Xi} such that  the items 1. and 2. in Problem I are satisfied in $\Xi$. 

\vspace{0.3cm}
\noindent
\underline{ Problem II.}\quad  {\em Let $(X,g_t)$, $t\in [0,1]$, be a $C^1$ family of billiard tables satisfying the weak isospectral condition  $(\mbox{H}_1)$-$(\mbox{H}_2)$. Consider  a $C^1$ family of Kronecker tori $[0,\delta] \ni t \mapsto \Lambda_t(\omega)$ of $P_t$ for each $\omega\in \Xi$.  
Are the functions $t\to \beta_t(\omega)$, $t\to I_t(\omega)$ and $t\to L_t(I_t(\omega))$ independent of $t\in [0,\delta]$ for each $\omega \in \Xi$?}\\

Affirmative answer of this question is given in  Theorem \ref{Th:main1} and Theorem \ref{Th:main2}.  This result can be applied in the cases $1. -  3.$ listed above. The proof is based on the construction of $C^1$ families of quasi-modes of the spectral problem \eqref{thespectrum} in  Theorem \ref{Theorem:quasimodes}. We  present below the main idea of the proof.

\vspace{0.3cm}
\noindent
\underline{ Problem III.}\quad {\em Does the weak isospectral condition 
$(\mbox{H}_1)$-$(\mbox{H}_2)$ imply the existence of a $C^1$ family of Kronecker tori $[0,1] \ni t \mapsto \Lambda_t(\omega)$ of $P_t$ all along the perturbation  for each $\omega \in \Xi$?}\\

This problem is closely related with a mysterious phenomena  in the Hamiltonian dynamics of close to integrable systems - the destruction of Kronecker tori with Diophantine frequencies along a perturbation. The $C^1$ family of Kronecker tori $ t \mapsto \Lambda_t(\omega)$ exists in a certain interval $[0,\delta_0)$ but it may cease to exist at $t=\delta_0$. Does the ``weak isospectral condition'' prevent the tori from destroying? 
We give a positive answer of Problem III in the following two cases - in the case 2. mentioned  above if the elliptic periodic broken geodesic has no resonances of order $\le 12$ (see Theorem \ref{theo:main-elliptic} and Proposition \ref{prop:interval}) and for the Lazutkin caustics in the case of a $C^1$ deformation of a strictly convex billiard table in $\R^2$ (see Theorem \ref{Th:convex-main}). 
The proof of these two results is rather involved. It requires a KAM theorem and BNF theorem where the  constant $\epsilon$ appearing in the smallness condition essentially  depends only on the dimension $n$ and the exponent $\tau>n-1$ but not on the particular completely integrable Hamiltonian (see Theorem \ref{Theo:KAM} and Theorem \ref{Theo:BNF}). 
We need as well suitable uniform with respect to $t$ \emph{global} estimates of the  H\"{o}lder $C^\ell$-norms $(\ell\ge 1)$ of the functions $L_t$  in the BNF \eqref{eq:birkhoff-L}. These estimates are obtained in Theorem \ref{Theo:BNF}.

\vspace{0.3cm}
\noindent
\underline{ Problem IV.}\quad {\em Spectral rigidity under the weak isospectral condition.}\\

We show in Proposition \ref{Prop:variation-alpha} that the variation  $\dot\beta_t(\omega)$, $\omega\in\Xi$, can be written by means of a suitable Radon transform at the family of Kronecker tori $\Lambda_t(\omega)$, $\omega\in\Xi$,  applied to the ``vertical component'' of the variation of the boundary $\Gamma_t$. In particular, the equality $\beta_t(\omega)=\beta_0(\omega)$, $t\in [0,\delta]$, $\omega\in \Xi$, obtained in Theorem \ref{Th:main2} implies that the image of  the Radon transform is zero for any weakly isospectral family (see Theorem \ref{Th:main3}). Hence, to prove infinitesimal rigidity one has to obtain injectivity of that Radon transform. In this way we obtain  infinitesimal spectral rigidity under the weak isospectral conditions  for classical Liouville Billiard Tables of dimension $2$ and $3$ in Theorem \ref{Th:Liouville-dim2} and Theorem \ref{Th:Liouville-dim3}. We obtain in particular that the billiard tables inside the ellipse in $\R^2$ and inside the ellipsoid in $\R^3$ are  infinitesimally spectrally rigid under the weak isospectral conditions $(\mbox{H}_1)-(\mbox{H}_2)$.  
Infinitesimal spectral rigidity of the billiard table inside the ellipse has been obtained by Hezari and Zelditch \cite{H-Z} under the usual isospectral condition using the wave-trace method.  Unfortunately infinitesimal spectral rigidity does not always  apply  spectral rigidity as in the case of negatively curved manifolds because of the phenomena of destruction of Kronecker tori with Diophantine frequencies.  

As an application of   Theorem \ref{Th:main3} we prove in Theorem \ref{theo:elliptic} spectral rigidity of analytic $\Z_2\oplus\Z_2$ symmetric billiard tables $(X_t,g)$ of dimension two if one of the corresponding bouncing ball trajectories is elliptic, it has no resonances of order $\le 4$ and the Poincar\'e map is nondegenerate.

\vspace{0.3cm}
\noindent
\underline{ Problem V.}\quad {\em Are classical Liouville Billiard Tables spectrally rigid?}\\

It turns out (see Corollary \ref{coro:ellipse})   that the map $P=B^2$  is always Kolmogorov nondegenerate (twisted)  at the elliptic point  for  elliptical billiard tables (bounded by an ellipse). 
 Moreover, except of five confocal families of ellipses given explicitly by \eqref{eq:rotation_ellipse'}, the  geodesic $\gamma_1$ is $4$-elementary. 
 These two conditions  are open in the $C^5$ topology, and applying Theorem  \ref{theo:elliptic} we obtain spectral rigidity not only of such elliptical billiard tables but also of analytic $\Z_2\oplus\Z_2$ symmetric billiard tables close to them. 

		\vspace{0.3cm}
		\noindent
		\underline{ Problem VI.}\quad {\em Estimates of the canonical transformation and the transformed Hamiltonian  in the KAM theorem in the scale of $C^\alpha$ norms.}\\

		In order to prove the main theorems in the first part of the article we need certain global estimates (in the whole domain of frequencies) of the canonical transformations and the transformed Hamiltonian in the KAM theorem in the whole scale of H\"older norms. Such estimates are obtained in Theorem \ref{Theo:A}, (iii), and in Theorem \ref{Theo:Holder}, using a new iteration schema, which allows us to avoid the Whitney extension theorem for $C^\infty$ jets.

\vspace{0.3cm}
Before giving the structure of the paper we would like to compare  different features of the spectral rigidity problem in the cases of negatively curved closed manifolds and  of close to integrable Hamiltonian systems.\\

\begin{tabular}{r l c}
\underline{\Large negative curvature}    &         
\multicolumn{2}{c}{\underline{\Large close to integrable} }\\[0.3cm]
Anosov  system                  &         KAM system  \\
$C^1$ families of hyperbolic closed geodesics & $C^1$ families of Kronecker tori\\
Labeled length spectrum                 &   Average action on the Kronecker tori,\\
                                &   Mather's $\beta$ function\\
Wave trace formula              &  $C^1$ families of quasi-modes \\
Geodesic ray transform          & Radon transform on Kronecker tori \\
Structural stability  of Anosov dynamics   & Phenomena of destruction of invariant tori\\
                                            & with Diophantine frequencies\\
Infinitesimal rigidity easily  & Passing from infinitesimal rigidity to \\
implies spectral rigidity &   spectral rigidity is a hard problem
\end{tabular}

\vspace{0.5cm}
We are going to describe now the structure of the paper. 

In Section \ref{Sec:results} we recall first the definition of the billiard ball map and then we formulate some of the main results. We give as well a proof of Theorem \ref{Th:main3} which reduces the spectral rigidity problem to the injectivity of a suitable Radon transform.

In Section \ref{Sec:BBM-BNF} we obtain by Theorem \ref{Theo:soft-BNF}   a $C^1$ family of BNFs of $P_t$ associated with  $C^1$ families of Kronecker tori $\Lambda_t(\omega)$, where $\omega\in \Omega_{ \kappa}^0$, $t\in J$, and $J\subset [0, \delta_0]$ is an interval. 
This family is analytic in $t$ if the map $t\to P_t$ is analytic. 
The theorem is based on the BNF obtained in Theorem \ref{Theo:BNF}.

Section \ref{sec:LBT} is devoted to Liouville billiard tables of dimension $n=2$ or $n=3$.  Liouville billiard tables of dimension two were defined in \cite[Sec. 2]{PT1} by using a branched double covering map. We give here an invariant definition of Liouville billiard tables in dimension two and we prove the equivalence of the two definitions  in Appendix \ref{sec:Appendix B}. Then we recall the definition of Liouville billiard tables of classical type in dimension two. Infinitesimal spectral rigidity of such billiard tables under the ``weak isospectral condition'' is obtained in Theorem \ref{Th:Liouville-dim2}. Infinitesimal spectral rigidity of nondegenerate Liouville billiard tables of classical type in dimension three is obtained in 
Theorem \ref{Th:Liouville-dim3}. Here we essentially use the injectivity of the corresponding Radon transform which has been proven in \cite{PT3}. In particular we obtain infinitesimal spectral rigidity of the ellipse in $\R^2$ and the ellipsoid in $\R^3$ under the ``weak isospectral conditions''. 

In Section \ref{Sec:elliptic} we consider $C^1$ isospectral deformations $[0,1] \ni t \to (X,g_t)$
of a given billiard table $(X,g_0)$ admitting an elliptic closed broken geodesic $\gamma$ 
with $m\ge 2$ vertices and we denote by $P_0:=B_0^m$ the corresponding Poincar\'e map. We suppose that $\gamma$ admits no resonances of order $\le 4$ and that the BNF of $P_0$ is nondegenerate. By the implicit function theorem  there exixts $\bar\delta>0$ and a $C^1$ family  of elliptic closed broken geodesic $\gamma_t$, $t\in [0,\bar\delta)$, 
with $m\ge 2$ vertices in $(X,g_t)$ having no resonances of order $\le 4$ and such that the BNF    of the corresponding Poincar\'e maps $P_t$ are nondegenerate. Theorem \ref{theo:elliptic} gives spectral rigidity for a class of analytic $\Z_2\oplus\Z_2$ symmetric billiard tables of dimension two. Then we  address the following questions. Suppose that the $C^1$ family of billiard tables $(X,g_t)$, $0\le t\le 1$, is weakly isospectral. Assume that  $(X,g_0)$ admits a periodic  elliptic broken geodesic $\gamma_0$ and that the corresponding local Poincaré map is twisted. Does there exist a $C^1$ family of periodic  elliptic broken geodesics $[0,1]\ni t\to \gamma_t$ in $(X,g_t)$ \emph{along the whole perturbation}? Do the corresponding local Poincaré map remain twisted?  Do the invariant tori $\Lambda_0(\omega)$ associated to $\gamma_0$ give rise to $C^1$ families of invariant tori $[0,1]\ni t\to \Lambda_t(\omega)$ \emph{along the whole perturbation}? We give an answer of these questions in  Theorem \ref{theo:main-elliptic} and Proposition \ref{prop:interval}. 

Section \ref{Sec:convex} is devoted to isospectral deformations of  locally strictly 
geodesically convex billiard tables of dimension two. Firstly we obtain in Proposition \ref{Prop:appr-inter-hamiltonian} a $C^1$ family of BNFs for the billiard ball maps $B_t$ in a neighborhood of $S^\ast_t\Gamma$ in terms of  the interpolating Hamiltonian $\zeta_t$ introduced  by Marvizi and Melrose \cite{MM}. Then Theorem \ref{Th:convex} gives an affirmative answer of  Problems I-III in the case of  $C^1$ families  $(X,g_t)$, $t\in [0,\delta]$,  of  locally strictly geodesically convex  billiard tables of dimension two satisfying the weak isospectral condition $(\mbox{H}_1)-(\mbox{H}_2)$.  Moreover, if $(X_t,g)$, $t\in [0,1]$, is a $C^1$ family of  billiard tables in $\R^2$ equipped with the Euclidean metric and satisfying  the weak isospectral condition and if $X_0$ is strictly convex then we prove in that $X_t$ remains strictly convex for each $t\in [0,1]$ and we get an affirmative answer of Problem III for  $t\in [0,1]$ (see Theorem \ref{Th:convex-main}). 

In Section \ref{Sec:QBNF} we reduce the problem \eqref{thespectrum} microlocally at the boundary. The main idea is explained in the beginning of Section \ref{Sec:QBNF}. Let $f_t=f_t(\cdot,\lambda)$ be a $\frac{1}{2}$-density   on $\Gamma$ depending on a large parameter $\lambda$ and with a frequency support  contained in a small neighborhood of the union of the invariant tori $\Lambda_t(\omega)$.  We consider the corresponding outgoing solution $u_t$ of the reduced wave equation (the Helmholtz equation) in $X_t$ with initial data $f_t$ and we ``reflect it at the boundary'' $m-1$ times if $m\ge 2$. To this end we use the outgoing parametrix of the reduced wave equation which is a Fourier Integral Operator with a large parameter $\lambda$ ($\lambda$-FIO). Taking the pull-back  to $\Gamma$ of the last branch of the solution $u_t(\cdot,\lambda)$ we get a $\frac{1}{2}$-density  $M_t^0(\lambda)f_t$, where  $M_t^0(\lambda)$ is a $\lambda$-FIO of order zero at $\Gamma$ the canonical relation of which is just the graph of $P_t$. We call $M_t^0(\lambda)$ a monodromy operator. In this way obtain that 
\[
(-\Delta_t +\lambda^2)u_t= O_N(|\lambda|^{-N})f_t, \quad {\mathcal B}_t u_t=O_N(|\lambda|^{-N})f_t
\]
if and only if
\[
M_t^0(\lambda)f_t=f_t + O_N(|\lambda|^{-N})f. 
\]
We are looking for couples $(\lambda, f_t)$ solving the last equation. 
To this end using the BNF of $P_t$  we obtain   a suitable microlocal (quantum) Birkhoff normal form  of $M_t^0(\lambda)$ for $t\in J$ (see Proposition \ref{operator-T} and  Proposition \ref{prop:commutator}  in Sect. \ref{subsec:QBNF}).
This enables us to  ''separate the variables'' microlocally near the whole family of invariant tori $\Lambda_t(\omega)$, $\omega\in\Omega_\kappa^0$, and to obtain a microlocal spectral decomposition of $M_t^0(\lambda)$   in Proposition \ref{prop:spectral-decomposition}. Then the  problem of finding $\lambda$ is reduced to an algebraic equation $\mu_t(\lambda)=1 + O_N(|\lambda|^{-N})$ where $\mu_t(\lambda)$ are suitable eigenvalues of $M_t^0(\lambda)$ with eigenfunctions $f_t$.  In this way we obtain that 
$\lambda$  should satisfy \eqref{eq:quantization-condition} and \eqref{eq:new-equation} and we solve that system of equations recursively. This is done in 
 Section \ref{Sec:Quasi}, where  we obtain  $C^1$ families of quasi-modes. Using these quasi-modes   we prove item 3 of Theorem \ref{Th:main1}
which claims that the function $t\to \beta_t(\omega)$ is independent of $t$ for  Diophantine frequencies $\omega$ provided that the weak isospectral condition $(H_1)$-$(H_1)$ is fulfilled. We are going to give the idea of the proof.\\

\noindent
\underline{\em Sketch of the Proof.}\\ 

\noindent
We fix $t\in [0,\delta)$, $\kappa>0$ and $\omega\in \Omega_\kappa^0$ and we impose the following \\

\noindent
{\em Strong Quantization Condition on the torus $\Lambda_t(\omega)$}. \quad 
There exists  an infinite sequence $\widetilde {\mathcal M}(\omega)$  of $(q,\lambda)\in \Z^n \times [1,\infty)$ such that
$q=(k,k_n)\in \Z^{n-1}\times \Z$ and  $\lambda=\mu_q^0\ge 1$ satisfy
\begin{equation}\label{eq:strong-quantization}
\left\{
\begin{array}{clrr}
\displaystyle c_0^{-1}|q| \, \le \, \mu_q^0 \, \le\,  c_0 |q|\\[0.3cm]
\displaystyle \lim_{|q|\to \infty}\, \Big|\mu_q^0 \, \Big(I_t(\omega), L_t(I_t(\omega))\Big)\  -\ \Big(k+\frac{\vartheta_0}{4}, 2\pi \Big(k_n+\frac{\vartheta}{4}\Big)\Big)\Big|\,  =\,  0.
\end{array}
\right.
\end{equation}
Here $c_0>1$ is a constant, $(\vartheta,\vartheta_0)\in \Z^n$,  $\vartheta$ is related with the Maslov class of $\Lambda_t(\omega)$ and $\vartheta_0$ is a Maslov index. It turns out that  condition \eqref{eq:strong-quantization}  is satisfied for each $\omega$ in a subset $\Xi^t_\kappa\subset \Omega_\kappa^0$ of full Lebesgue measure in $\Omega_\kappa^0$ (see Lemma \ref{Lemma:quantization}). 

Now we fix $\omega\in \Xi_\kappa^t$ and 
denote by  $\mathcal M\subset \Z^n$ the image of the projection of $\widetilde {\mathcal M}(\omega)\subset \Z^n \times  [1,\infty)$ on the first factor. The set $\mathcal M$ will be the index set of the $C^1$ family of quasi-modes that we are going to construct.  To obtain a quantization condition for the tori $\Lambda_s(\omega)$ for $s$  close to $t$ we introduce 
for any $q\in {\mathcal M}$ the interval 
\[
J_q:=\left[t,t+2|q|^{-1}  \right] .
\] 
Since the  maps  $s\mapsto L_s$ and 
$s\mapsto I_s$ are $C^1$ in a neighborhood of $t$ with values in the corresponding $C^\infty$ spaces (see Theorem \ref{Theo:soft-BNF} and Definition  \ref{Def:BNF}), the following quantization condition of the tori $\Lambda_s(\omega)$ is satisfied\\

\noindent
 {\em Quantization Condition}. \quad There exists a constant $C=C(\omega)>0$ independent of  $q\in {\mathcal M}$ and  $s\in J_q$ such that 
\begin{equation}
\Big| \mu_q^0 \Big(I_s(\omega), L_s(I_s(\omega))\Big)\  -\ \Big(k+\frac{\vartheta_0}{4}, 2\pi \Big(k_n+\frac{\vartheta}{4}\Big)\Big)\Big|\, \le\,  C \quad \forall\, q\in {\mathcal M},\  s\in J_q. 
\label{eq:usual-quantization}
\end{equation}
Using this condition we construct $C^1$ families quasi-modes $(\mu_q(s)^2,u_{s,q})$, $q\in {\mathcal M}$, $s\in J_q$,  of order $M$ for the problem \eqref{thespectrum} such that (see Theorem \ref{Theorem:quasimodes})
\begin{enumerate}
\item[(1)] $ u_{s,q}\in  D(\Delta_s)$, $\|u_{s,q}\|_{ L^2(X)}=1$, and there exists a constant $C_M>0$ such that
\[
\left\{  
\begin{array}{lcr}  
\left\|\Delta\,  u_{s,q}\ - \ \mu_q^2(s)\,  u_{s,q}\right\|\ \le \   
C_M\,  \mu_q^{-M}(s)  \, \quad \mbox{in}\ L^2(X)\, , \\ [0.3cm] 
\displaystyle {\mathcal B}\,  u_{s,q}|_\Gamma \   =\ 0 \,   
\end{array}  
\right.  
\]
for every $q\in {\cal M}$ and  $s\in J_q$;   
\item[(2)] We have 
\[
\mu_q(s) = \mu_q^0 +c_{q,0}(s) + c_{q,1}(s) \frac{1}{\mu_q^0}  + \cdots + 
c_{q,M}(s) \frac{1}{(\mu_q^0)^M}
\]   
where the functions $s\mapsto c_{q,j}(s)$ are real valued and $C^1$ on the interval $J_q$ and  
there exists  a constant $C_M'>0$ such that
$|c_{q,j}(s)| \le C_M'$ for every  $q\in {\mathcal M}$, $0\le j\le M$,  and any $s\in J_q$;
\item[(3)] There exists $C>0$ such that 
\[
\Big|\mu_q(s) L_s\Big(\frac{k + \vartheta_0/4}{\mu_q(s)}\Big) -  2\pi \Big(k_n +
\frac{\vartheta}{4}\Big)\Big| \, \le\,  \frac{C}{\mu_q(s)}
\]
for every $q\in  {\cal M}$ and $s\in J_q$;
\item[(4)]  We have
\[
\frac{k + \vartheta_0/4}{\mu_q(t)}\, =\, I_t(\omega) + o\Big(\frac{1}{|q|}\Big)\quad \mbox{as}\ |q|\to \infty. 
\]
\end{enumerate}
We point out that the strong quantization condition \eqref{eq:strong-quantization} is needed only in the proof of item (4). Items (1)-(3) follow from the weaker quantization condition \eqref{eq:usual-quantization}. The estimate in (3) is related to the nullity of the subprincipal symbol of the Laplace-Beltrami operator. Note also  that the quasi-eigenvalues $\mu_q(s)^2$ are defined only in the intervals $J_q$ which shrink to $t$ as $q\to \infty$. 

Consider now the self-adjoint operator $\Delta_s$ with Dirichlet, Neumann or Robin boundary conditions and the corresponding spectral problem \eqref{thespectrum}. The relation between the spectrum of $\Delta_s$ and the quasi-eigenvalues $\mu_q(s)^2$ is given by 
\[
\mbox{dist} \Big(\mbox{Spec} \left(\Delta_s\right), \mu_q(s)^2\Big) \, \le \, C_M\,  \mu_q^{-M}(s)
\]
where $C_M$ is the constant in (1) and $M$ is the order of the quasi-mode. We fix $M>2d+2$ where $d\ge 0$ is the exponent in $(\mbox{H}_1)$. It  follows then from   (H$_2$)  that the quasi-eigenvalues  $\mu_q(s)^2$, $|q|\ge q_0 \gg 1$,  $s\in J_q$, belong
to the  union of intervals 
\[
A_k:=\Big[a_k-\frac{c}{2}a_k^{-d-1},b_k+\frac{c}{2}a_k^{-d-1}\Big] \quad k\ge k_0\gg 1
\]
 where $c$ is the positive constant in the third hypothesis of (H$_1$). These intervals do not 
 intersect each other  for $k\ge k_0\gg 1$ in view of the third hypothesis of  (H$_1$). The function $J_q\ni s\mapsto \mu_q(s)^2$  being  continuous in $J_q$ 
can not jump from one interval to another, hence, it is trapped in a certain interval $A_k$.  Then using the first and second hypothesis of (H$_1$) we obtain
\[
|\mu_q(s) - \mu_q(t)| \le \frac{1}{\mu_q(t)}|\mu_q(s)^2 - \mu_q(t)^2|\le 2a_k^{-\frac{1}{2}}(b_k-a_k) = o(a_k^\frac{1}{2})=o(\mu_q(t))= o\left(\frac{1}{|q|}\right)
\]
for $|q|\ge q_0$ where $q_0\gg 1$ does not depend on the choice of $s\in J_q$. 
This implies
\[
\mu_q(s) = \mu_q(t)\left(1+ o\left(\frac{1}{|q|}\right)\right) \quad \mbox{as}\ q \to \infty
\]
uniformly with respect to $s\in J_q$. A detailed proof of this statement is given in   Lemma \ref{Lemma:isospectral}. Now using (4) we get
\[
\zeta_q(s): = \frac{k + \vartheta_0/4}{\mu_q(s)} = \frac{k + \vartheta_0/4}{\mu_q(t)(1 + o(1/|q|))}= \frac{k + \vartheta_0/4}{\mu_q(t)} + o\left(\frac{1}{|q|}\right)= I_t(\omega) + o\left(\frac{1}{|q|}\right),\quad q\to \infty
\]
uniformly with respect to $s\in J_q$.  In the same way we get from (3) 
\[
L_s\left(\zeta_q(s)\right)=  2\pi\frac{k_n-\vartheta/4}{\mu_q(s)} 
 +   O(|q|^{-2})= 2\pi\frac{k_n-\vartheta/4}{\mu_q(t)} + o\left(\frac{1}{|q|}\right)= L_t\left(\zeta_q(t)\right)  + o\left(\frac{1}{|q|}\right),\quad q\to \infty
\]
uniformly with respect to $s\in J_q$.  Setting $\eta:=1/|q|$ we obtain from the above   equalities that 
\[
\begin{array}{lcrr}
L_{t+\eta}\big(I_t(\omega)\big)=L_{t+\eta}\big(\zeta_q(t+\eta)+ o(\eta)\big)=L_{t+\eta}\big(\zeta_q(t+\eta)\big)+ o(\eta)\\[0.3cm]
= L_{t}\big(\zeta_{q}(t)\big)+ o(\eta)=  
L_t\big(I_t(\omega)\big)  + o(\eta), \quad \eta=1/|q| \to 0.
\end{array}
\]
Recall that   the map $[0,\delta]\to L_s(\cdot)$ is $C^1$ with values in the corresponding $C^\infty$ space.
Hence, 
\[
 \dot{L_t}(I_t(\omega))=\frac{d}{ds}L_s(I_t(\omega))\big|_{s=t} = 0 \quad \forall\,  \omega\in \Xi_\kappa^t\, . 
\]
On the other hand, $\Xi_\kappa^t$ is dense in $\Omega_\kappa^0$ since any point of $\Omega_\kappa^0$ is of positive Lebesgue density and   $\Omega_\kappa^0\setminus\Xi_\kappa^t$ has measure zero, and by continuity (the function $I\to \dot L_t(I)$ is smooth) we get  this equality for each   $\omega\in\Omega_\kappa^0$.  The point $t$ has been  fixed arbitrary in $[0, \delta)$, hence,  $\dot{L_t}(I_t(\omega))=0$  for every $t\in [0, \delta)$ and $\omega\in\Omega_\kappa^0$. 
Now differentiating the first equation of \eqref{eq:beta-alpha} with respect to $t$   we obtain
\[
 \dot{\beta}_t(\omega) = \langle \omega, \dot{I}_t(\omega)\rangle - \dot{L}_t(I_t(\omega)) - \langle \nabla L_t(I_t(\omega)), \dot{I}_t(\omega)\rangle = 0 \quad \forall\,  \omega\in \Omega_\kappa^0 
\]
since $\nabla L_t(I_t(\omega))=\omega$. Hence, $\beta_t(\omega)=\beta_0(\omega)$ for every $t\in [0,\delta)$ and $\omega\in \Omega_\kappa^0$.
By continuity we get  the last equality for every $t\in [0,\delta]$ as well. \finishproof

We point out that the classical and quantum BNFs are analytic in $t$ if the perturbation is analytic in $t$ which leads to analytic in $t$ quasi-modes. This can be used as in   \cite{Go, Go-H} to extend the results of S. Gomes and A. Hassell about  the quantum non-ergodicity of $C^\infty$-smooth KAM systems. Moreover, using Theorem \ref{Theo:Holder} and the pseudodifferential calculus of operators with symbols of finite smoothness \cite{PT4}, one may extend them to KAM systems of finite smoothness. \\

The second part of the manuscript  is devoted to the KAM theorem and the BNF around families of invariant tori in both the continuous and discrete cases. The main novelty in it can be briefly summarized as follows 
\begin{itemize}
	\item[--] the constant  $\epsilon$ in the smallness condition essentially depends only on the dimension of the configuration space and on the exponent in the Diophantine condition;
	\item[--] $C^k$ smooth (analytic) families of invariant tori $t\to \Lambda_t(\omega)$ with Diophantine frequencies are obtained;
	\item[--] $C^k$ smooth  (analytic)  with respect to the parameter $t$ BNF is obtained around the union of $\Lambda_t(\omega)$;
	\item[--] global estimates in the whole scale of H\"older norms with universal constants are obtained.  To this end a new  iterative schema is proposed. The Modified Iterative  Lemma proven in Sect. \ref{Prop:IterativeLemma} provides in a limit smooth functions in the whole domain $\Omega$  (not only smooth Whitney jets on the Cantor set $\Omega_{ \kappa}$) with a good control of the H\"older  norms. 
\end{itemize}
We need all these properties un the first part of the manuscript. 
The KAM theorems here are based on Theorem \ref{Theo:A} which is a KAM theorem for  $C^k$ (k=0; 1) or analytic families of $C^\infty$ smooth Hamiltonians $H_t$ in $\T^n\times D$ with parameters $\omega\in \Omega$ where $H_t$ are small perturbations of the  normal form  ${\mathcal N}(I;\omega):= \langle \omega, I\rangle$. The proof of the theorem, especially of the so called KAM step follows that of J. P\"{o}schel \cite{Poe1} in the case of analytic Hamiltonians but it requires additional work in order to adapt it to the case of $C^k$ families of Hamiltonians $[0,1]\ni t \to C^\infty(\T^n\times D; \Omega)$. For this reason we give a complete proof of the KAM step. Next we adapt the  Iterative Lemma in \cite{Poe1} to the case of smooth Hamiltonians. 
In this way one obtains  an iteration schema which gives in a limit $C^\infty$ Whitney jets on $\Omega_{ \kappa}$. The Whitney extension theorem in the $C^\infty$ case does not provide in general global estimates of the H\"older norms in $\Omega$ without loss of derivatives \cite{F-J-W}. For this reason we provide another iteration scheme based on a Modified Iteration Lemma given by Proposition \ref{Prop:ModifiedIterativeLemma}, which involves Gevrey almost analytic extensions of certain cut-off Gevrey functions. The $\bar \partial$ derivatives of such functions are exponentially small near the reals and one can use Cauchy (Green's) formula. This allows one to obtain a convergent iteration schema on the the whole space of frequencies and to obtain the desired global (in $\Omega$) H\"older estimates of any   order.

Using Theorem  \ref{Theo:A} we obtain a KAM theorem for $C^1$ families of Hamiltonians $H_t$ which are perturbations of a given $C^1$ family of completely integrable nondegenerate in Kolmogorov sense  Hamiltonians $H_t^0$. The  family $H_t^0$ is given as follows. We consider a $C^1$ family of nondegenerate real valued functions $K_t\in C^\infty(\Omega)$, $t\in [0,\delta]$, in a domain $\Omega\subset \R^n$, where by nondegenerate we mean that the gradient map 
\[
\Omega\ni \omega \mapsto \nabla K_t(\omega) \in D_t:= \nabla K_t(D_t)
\]
is a diffeomorphism for each $t\in [0,\delta]$. We denote by $H_t^0\in C^\infty(D_t)$ the Legendre transform $K_t^\ast$ of $K_t$ given by
\[
H_t^0(I)\, =\, K_t^\ast(I) \, :=\,  \mbox{Crit.val.\,}_{\omega\in \Omega}\{\langle \omega, I\rangle - K_t(\omega)\}.
\]
Then $H_t^0\in C^\infty(D_t)$ is nondegenerate and $(H_t^0)^\ast= K_t$. Theorem \ref{Theo:KAM2} provides a result of KAM type for $C^k$ families of Hamiltonians $[0,\delta]\ni t \mapsto H_t\in C^\infty (\T^n\times D_t)$ which is a small perturbation of the family $H_t^0$. The  constant $\epsilon$ in the smallness condition essentially depends only on the dimension $n$ and on the exponent $\tau>n-1$ in the Diophantine condition \eqref{eq:sdc-1-intro}. To  this end, given $\omega\in \Omega$  we linearize $H_t^0$ at $I=\nabla H_t^{0\ast}(\omega)$ applying Taylor's formula and sending the nonlinear part of it to the functions with perturbation. The smallness condition and the  estimates  in  Theorem \ref{Theo:KAM2}  are given in terms of suitable weighted H\"{o}lder norms. In order to obtain  estimates in H\"{o}lder norms with universal constants  of a composition of functions with  $ H_t^{0\ast}$,     we suppose that $\Omega$ is a strictly convex bounded domain  in $\R^n$ and that $t \to  H_t^{0\ast}$  is a $C^k$ family with values in $C^\infty(\overline{\Omega}, \R)$, $\overline{\Omega}$ being the closure of $\Omega$. 
The idea is to use the interpolation inequalities for H\"{o}lder norms in $\R^n$ or in $\T^n\times \R^n$ which simplify a lot the estimates of higher order  H\"{o}lder norms for the inverse function and for the composition of functions. 
The problem about the composition of functions in   H\"{o}lder spaces is quite delicate. It has been investigated recently by  R. de la Llave and  R. Obaya
\cite{L-O}. We can not use directly their results here since we need estimates with universal constants. These estimates are obtained in Appendix \ref{subsec:composition}. 

Theorem \ref{Theo:A-for-maps}  is a counterpart of Theorem \ref{Theo:A} in the discrete case for $C^1$ families of exact symplectic maps $P_t$. Theorem \ref{Theo:A-for-maps} is obtained  from Theorem \ref{Theo:A} using an idea of R. Douady \cite{Dou}. The BNF of $C^k$ families of exact symplectic maps at $C^k$ families of Kronecker tori is obtained in Section \ref{Subsec:BNF-for-maps}. We point out that the constant $\epsilon$ in the corresponding smallness conditions \eqref{eq:smallness-condition-maps-2} and \eqref{eq:smallness-condition-1} essentially depend only on the dimension and on the exponent $\tau$. Moreover, the constants $C_m$ in the corresponding H\"{o}lder estimates \eqref{eq:estimates-F-1}, \eqref{eq:estimates-chi}, \eqref{eq:estimates} and in Theorem \ref{Theo:Holder} are universal. This makes these results especially useful in the case when the symplectic maps $P_t$ have singularities. They can be applied for  example for  the billiard ball map $B_t$ near the singular set $S_t^\ast \Gamma$ in the case when $(X,g_t)$ is locally strictly geodesically convex and ${\rm dim}\, X=2$ (see Theorem \ref{Th:convex}). 

\part{Isospectral invariants and rigidity}

\section{Main Results}\label{Sec:results}

Before formulating the main results we recall from Birkhoff \cite{Birkhoff} (see also \cite{Tabach}) the definition 
of the billiard ball map $B$ associated to a billiard table
$(X,g)$ with a smooth boundary $\Gamma$.

\subsection{Billiard ball map}\label{subsec:billiard-ball}
Let $(X,g)$ be a smooth billiard table wioth boundary $\Gamma$. 
The ``broken geodesic flow'' given by the elastic reflection of geodesics hitting transversely the boundary induces a discrete dynamical system at the boundary which can be described as follows. 

Denote by $h$  the Hamiltonian on $T^\ast X$ 
corresponding to the Riemannian metric $g$ on $X$ via the Legendre
transformation and by $h^0$ the Hamiltonian on $T^\ast\Gamma$ corresponding to the induced
Riemannian metric on $\Gamma$. 
The billiard ball map $B$ lives in an open  subset of  the open 
coball bundle 
$\textbf{B}^\ast \Gamma = \{(x,\xi)\in T^\ast \Gamma:\, h^0(x,\xi) < 1\}$. 
It is
defined as follows.  Denote by   $S^\ast X:=\{(x,\xi)\in T^\ast X:\, h(x,\xi) =  1\}$ the cosphere bundle,  and set
\[
\Sigma = S^\ast X|_{\Gamma}:= \{(x,\xi)\in S^\ast X:\, x\in
\Gamma\}\ \mbox{and}\ \Sigma^{\pm}:=  \{(x,\xi)\in \Sigma :\,
\pm \langle \xi,\nu(x)\rangle >0\}\, 
\]
where $\nu(x)\in T_x X$, $x\in\Gamma$,  is the outward unit normal to $\Gamma$ with respect to  the metric $g$. 
Let $\overline{\textbf{B}^\ast \Gamma}$ be  closed coball bundle, i.e. the closure  of $\textbf{B}^\ast\Gamma$ in $T^\ast \Gamma$.
Consider the natural projection $\pi_\Sigma:\, 
\Sigma \, \rightarrow \, \overline{\textbf{B}^\ast \Gamma}$ assigning to each $(x,\eta)\in
\Sigma$ the covector $(x,\eta|_{T_x\Gamma})$. Its restriction to $\Sigma^+\cup \Sigma^-$ admits two smooth inverses 
\begin{equation}\label{eq:outgoing-vector}
\pi_\Sigma^{\pm}:\,  \textbf{B} ^\ast \Gamma\,  \rightarrow\, 
\Sigma^{\pm}\, ,\ \pi_\Sigma^{\pm}(x,\xi) = (x,\xi^\pm) \, .
\end{equation}
The maps $\pi_\Sigma^{\pm}$  can be extended continuously on the  closed coball bundle  $ \overline{{\bf B} ^\ast \Gamma}$.  
Given $\varrho^\pm\in \Sigma^\pm$ we consider  the integral  
curve  $\exp(sX_{h})(\varrho)$  of the Hamiltonian vector field   
$X_{h}$ starting at $\varrho$. 
If it intersect $\Sigma$ transversely    
 at a time $T=T(\varrho)$ and lies entirely in the   
interior  of $S^\ast  X$   
for $t$ between $0$ and $T$ we set   
\[
J(\varrho^\pm) =\exp(T X_{h})(\varrho^\pm)\in \Sigma^{\mp}\, .
\] 
Notice that $J$ is a smooth involution defined  in an open dense  subset of $\Sigma$. 
 In this way we obtain a smooth exact symplectic map   
$B:\widetilde {\bf B}^\ast \Gamma \rightarrow {\bf B}^\ast \Gamma $, given by 
$B=\pi_\Sigma\circ J \circ \pi_\Sigma^+$, where $ \widetilde B^\ast \Gamma$ is an open dense subset of ${\bf B}^\ast \Gamma$. 
The map $J$ can be extended to a smooth involution of $\Sigma$ in the case when $X$ is a strictly convex billiard table in $\R^n$. In this case the billiard ball map is well defined and smooth  in $  {\bf B}^\ast \Gamma$ and can be extended by continuity as the identity map on its boundary $  S^\ast \Gamma$. This case will be considered in more details in Sect. \ref{Sec:convex}.

Suppose now that $t\to g_t$ is a $C^1$ family of Riemannian metrics on $X$. For any $t$ we denote   the corresponding cosphere bundle by $S_t^\ast X:=\{h_t =  1\}$ and the corresponding open coball bundle of $\Gamma$ by ${\bf B}_t^\ast \Gamma:=\{h_t^0 <  1\}$. Let 
$\pi_{t}:\, \Sigma_t \, \rightarrow \, \overline{B_t^\ast \Gamma}$ be the natural projection and $\pi_{t}^{\pm}:\,  {\bf B}_t ^\ast \Gamma\,  \rightarrow\, \Sigma_t^{\pm}\, ,\ \pi_{\Sigma_t}^{\pm}(x,\xi) = (x,\xi_t^\pm)$ 
its inverses.  Denote by $J_t$ the corresponding involution in $\Sigma_t$ and consider the billiard ball map $B_t: \widetilde {\bf B}_t^\ast \Gamma \to {\bf B}_t^\ast \Gamma$. If $(x,\xi)\in  \widetilde {\bf B}_{t_0}^\ast \Gamma$, then $(x,\xi)\in  \widetilde {\bf B}_{t}^\ast \Gamma$ for any $t$ in a neighborhood of $t_0$ because of the transversality and one can show that the map $t\to  B_t\in C^\infty(\widetilde {\bf B}_{t}^\ast \Gamma, {\bf B}_{t}^\ast \Gamma)$ is $C^1$. 
In this way we obtain a $C^1$ family of symplectic mappings 
$t\to (B_t: \widetilde {\bf B}_t^\ast \Gamma \to {\bf B}_t^\ast \Gamma)$.

\subsection{Main Results}\label{subsec:main-results}
Recall that $P_0$ admits a $C^\infty$ Birkhoff normal form at any Kronecker invariant torus $\Lambda_0(\omega)$ with a Diophantine frequency $\omega$ (cf \cite{La}, Proposition 9.13).  
The Birkhoff normal form of $P_0$ at $\Lambda_0(\omega)$ is said to be nondegenerate if the quadratic part of it is a nondegenerate quadratic form. The non-degeneracy of the Birkhoff normal form enables one to apply the KAM theorem. Recall that  $D^0(\kappa,\tau)$ is the set of points of positive Lebesgue measure in $D(\kappa,\tau)$, defined in the Introduction. 
\begin{Th}\label{Th:main1}
Let $(X,g_t)$, $t\in [0,1]$,  be a $C^1$ family of billiard tables. Let $\Lambda_0(\omega_0)\subset {\bf B}^\ast_0\Gamma$ be a Kronecker invariant torus of $P_0:= B_0^m$ of a	Diophantine frequency $\omega_0\in D^0(\kappa_0,\tau)$. Suppose that the Birkhoff normal form of $P_0$ at $\Lambda_0(\omega_0)$ is nondegenerate. Then there exists $\delta_0 = \delta_0(\kappa_0)>0$ such that the following holds.
\begin{enumerate}
	\item There exists a $C^1$-family of Kronecker invariant tori  
	\[
	[0,\delta_0)\ni t\to\Lambda_t(\omega_0)\subset {\bf B}^\ast_t\Gamma
	\]
	of $P_t:= B_t^m$ of a frequency $\omega_0$.  
	\item For any $0<\delta<\delta_0$ there exists a set $\Xi\subset \T^{n-1}$ of Diophantine frequencies such that $\omega_0\in\Xi$, 
	\[
	\frac{{\rm meas\,}(B(\omega_0,\varepsilon)\cap \Xi)}{{\rm meas\,}(B(\omega_0,\varepsilon))} = 1 - O_\delta(\varepsilon) \quad \mbox{as}\ \varepsilon\to 0,
	\]
	and for any $\omega\in\Xi$ there exists a $C^1$-family of Kronecker invariant tori 
	\[
	[0,\delta]\ni t\to\Lambda_t(\omega)\subset {\bf B}^\ast_t \Gamma
	\]
	of $P_t:= B_t^m$ of a frequency $\omega$.
	\item If the the billiard tables satisfy the weak isospectral condition $(\mbox{H}_1)-(\mbox{H}_2)$ then $\beta_t(\omega)$, $I_t(\omega)$ and  $L_t(I_t(\omega))$  
are independent of $t\in [0,\delta]$ for any $\omega\in\Xi$. 
\end{enumerate}
\end{Th}
We are going to apply this result for $C^1$ deformations of the boundary keeping the Riemannian metric fixed. Let $X$ be a smooth compact manifold of dimension $n\ge 2$ with non-empty boundary $\Gamma:=\partial X$ which is smoothly embedded into a Riemannian manifold $(\widetilde X,g)$ of the same dimension and without boundary.
We say that $(X_t,g)$, $t\in [0,\varepsilon]$, is a $C^1$ variation of $(X,g)$  if  $(X_t,g)$ is a billiard table in $(\widetilde X,g)$  with boundary $\Gamma_t=\partial X_t$, $X_0=X$, and if there exists a $C^1$ family of embeddings  
\begin{equation}
[0,\varepsilon]\ni t\to\psi_t\in C^\infty(\Gamma,\widetilde X)
\label{eq:variation-boundary}
\end{equation}
such that $\psi_t(\Gamma)=\Gamma_t$. In this case we say that  $(X_t,g)$ is a $C^1$ family of billiard tables. 
Then there exists a $C^1$ family of  diffeomorphisms onto their images
$[0,\varepsilon]\ni t\to\Psi_t\in C^\infty(X,\widetilde X)$   such that $\Psi_t(X)= X_t$, $\Psi_t|_{\Gamma}=\psi_t$,
and $\Psi_t$ is identity outside an open neighborhood of $\Gamma$ in $X$. 
The family $\Psi_t$ can be constructed  parameterizing a neighborhood of $\Gamma$ in $\widetilde X$ by the exponential
map $(x,s)\mapsto \exp_x(s\nu(x))$ corresponding to the Riemannian metric $g=g_0$ in $\widetilde X$, where $\nu(x)$ is
the outward unit normal to $\Gamma$.
In particular, we get a family of billiard tables $(X_t,g)$, $t\in [0,\varepsilon]$, which are isometric to $(X,g_t)$,
$g_t:=\Psi_t^\ast g$. We say that a family of Kronecker invariant tori 
$[0,\delta)\ni t\to\Lambda_t(\omega)\subset {\bf B}^\ast \Gamma_t$ is $C^1$-smooth if the family 
$[0,\delta)\ni t\to\psi_t^\ast(\Lambda_t(\omega))\subset {\bf B}^\ast \Gamma$ is $C^1$. 
Consider the corresponding Laplace-Beltrami operator $\Delta_t$ in $(X_t,g)$ with Dirichlet, Neumann or Robin boundary
conditions on $\Gamma_t$. As a corollary of the main theorem we obtain 
\begin{Th}\label{Th:main2}
Let $(X_t,g)$, $t\in [0,1]$,  be a $C^1$ family of billiard tables. Let $\Lambda_0(\omega_0)\subset {\bf B}^\ast\Gamma$
be a Kronecker invariant torus of $P_0:= B_0^m$ of a Diophantine frequency $\omega_0\in D^0(\kappa_0,\tau)$. Suppose that
the Birtkhoff normal form of $P_0$ at $\Lambda_0(\omega_0)$ is non-degenerate. Then there exists 
$\delta_0 = \delta_0(\kappa_0)>0$ such that the following holds.
\begin{enumerate}
	\item There exists a $C^1$-family of Kronecker invariant tori  
	\[
	[0,\delta_0)\ni t\to\Lambda_t(\omega_0)\subset {\bf B}^\ast\Gamma_t
	\]
	of $P_t:= B_t^m$ of a frequency $\omega_0$.  
	\item For any $0<\delta<\delta_0$ there exists a set $\Xi\subset \T^{n-1}$ of Diophantine frequencies such that $\omega_0\in\Xi$, 
	\[
	\frac{{\rm meas\,}(B(\omega_0,\varepsilon)\cap \Xi)}{{\rm meas\,}(B(\omega_0,\varepsilon))} = 1 - O_\delta(\varepsilon) \quad \mbox{as}\ \varepsilon\to 0,
	\]
	and for any $\omega\in\Xi$ there exists a $C^1$-family of Kronecker invariant tori 
	\[
	[0,\delta]\ni t\to\Lambda_t(\omega)\subset {\bf B}^\ast \Gamma_t
	\]
	of $P_t:= B_t^m$ of a frequency $\omega$.
	\item If the the billiard tables satisfy the weak isospectral condition $(\mbox{H}_1)-(\mbox{H}_2)$ then $\beta_t(\omega)$, $I_t(\omega)$ and  $L_t(I_t(\omega))$  
are independent of $t\in [0,\delta]$ for any $\omega\in\Xi$. 
\end{enumerate}
\end{Th}

We shall denote by $\delta_\nu \Gamma_t:\Gamma_t\to \R$ the vertical component of the variation $\Gamma_s$ of 
$\Gamma_t$ which is defined by 
\begin{equation}\label{eq:vertical-component}
\forall \, x\in \Gamma_t,\quad \delta_\nu \Gamma_t(x):= \left\langle \dot\psi_t(\psi_t^{-1}(x)),
\nu_t(x)\right\rangle=\left\langle \frac{d \psi_s}{ds}|_{s=t} (\psi_t^{-1}(x)),\nu_t(x)\right\rangle,
\end{equation}
where $\nu_t(x)$ is the outward unit normal to $\Gamma_t$ at $x$ with respect to the metric $g$.

Let $\pi_t: T^\ast \Gamma_t \to \Gamma_t$ be  the natural projection. Given $\zeta=(x,\xi) \in B ^\ast  \Gamma_t$, we
denote by $\xi_t^+(\zeta)\in T_x^\ast X$ the corresponding outgoing unit covector which means that the restriction of 
the covector $\xi_t^+$ to $T_x\Gamma_t$ equals $\xi$, $h_t( \xi_t^+)=1$, and 
$\langle \xi_t^+(\zeta),\nu_t( \pi_t(\zeta))\rangle_x\ge 0$, where $\langle \cdot,\cdot\rangle_x$ stands for the 
paring between  covectors in $T_x^\ast \widetilde X$ with vectors in $T_x\widetilde X$  
(see Sect. \ref{subsec:billiard-ball}). 
In this case the second part of Corollary \ref{Th:main2} can be stated as follows. 
\begin{Th}\label{Th:main3}
Let $(X_t,g)$, $t\in [0,\varepsilon]$,  be a $C^1$-family of billiard tables satisfying the isospectral condition
$(\mbox{H}_1)-(\mbox{H}_2)$.  Let $[0,\delta)\ni t\to \Lambda_t(\omega)$,  $0<\delta\le \varepsilon$, be a $C^1$
family of invariant tori of $P_t= B_t^m$ with a Diophantine vector of rotation $\omega$. Then  
\begin{equation}
\sum_{j=0}^{m-1}\int_{\Lambda_t(\omega)}\,\left(\Big\langle \xi_t^+(\zeta),\nu_t(\pi_t(\zeta))\Big\rangle\,  
\left(\delta_\nu \Gamma_t\right)(\pi_t(\zeta))\right)\Big|_{\zeta=B^j_t(\rho)}\,  d\mu_t(\rho)\ =\ 0
\label{eq:integral-invariants}
\end{equation}
for any  $t\in [0,\delta)$. 
\end{Th}

There are three particular cases we will focus our attention at, namely, $C^1$ deformations of nondegenerate Liouville billiard tables, deformations of a manifold  having non-degenerate elliptic periodic geodesics and deformations of strictly convex planar domains.

We shall  prove  a spectral rigidity result for analytic billiard tables of dimension two having the symmetries of the ellipse  provided that  one of the bouncing ball rays of the initial billiard table is elliptic,  $4$-elementary, and has a nondegenerate BNF. 
Define a class of billiard tables as follows.
Let $(\widetilde X,g)$, $\mbox{dim}\, \widetilde X= 2$ be a Riemannian manifold of dimension two. Suppose that it admits   two commuting involutions ${\mathcal J}_k$, $k=1,2$, acting as isometries. Consider the family ${\mathcal B}$ of  billiard tables    $(X,g)$ in $(\widetilde X,g)$ such that the boundary $\Gamma=\partial X$ of $X$ is connected and  invariant with respect to ${\mathcal J}_k$, $k=1,2$. Then the set of fixed points of ${\mathcal J}_k$, $k=1,2$, in $X$ defines a bouncing ball geodesic $\gamma_k$ of any $(X,g)\in{\mathcal B}$.   Denote by ${\mathcal B}_{\rm an}$ the set of analytic billiard tables which belong to  ${\mathcal B}$. 
\begin{Th}\label{theo:elliptic} Let $(X,g)\in{\mathcal B}$. Assume that the  broken geodesic $\gamma_1$ given by the set of fixed points of ${\mathcal J}_1$ in $X$ is  elliptic $4$-elementary and that  the corresponding Poincare map admits a non-degenerate BNF. 
	Suppose that $(X_t,g)\in{\mathcal B}$, $t\in [0,1]$,  is a $C^1$ deformation of $(X,g)$ satisfying the weak isospectral condition $(\mbox{H}_1)-(\mbox{H}_2)$. Then   $\gamma_1$ is a bouncing ball geodesic of $(X_t,g)$ for each $t\in [0,1]$ and $\Gamma_t$ has a contact of infinite order to $\Gamma_0$ at the vertexes  of $\gamma_1$. In particular, $X_1=X_0$ if the boundaries $\Gamma_1=\partial X_1$ and $\Gamma_2=\partial X_2$ are both analytic.
\end{Th}

It turns out (see Corollary \ref{coro:ellipse})   that the Poincar\'e map associated with the elliptic bouncing ball geodesic $\gamma_1$ is always non-degenerate (twisted)  for Liouville Billiard Tables in surfaces of constant curvature. Let us fix the foci $F_1\neq F_2$ and consider the corresponding confocal ellipses. Then,    except of five confocal families of ellipses given explicitly by \eqref{eq:rotation_ellipse'}, the  geodesic $\gamma_1$ is $4$-elementary.    Denote by $\mathcal{E}$ the set of ellipses which do not belong to these families. A billiard table in $\R^2$ is said to be elliptical if its boundary is an ellipse. 
\begin{Th}\label{theo:ellipse1} We have the following.
	\begin{enumerate}
		\item 
	Each elliptical billiard table with boundary in  $\mathcal{E}$  is spectrally rigid in the class ${\mathcal B}_{\rm an}$  under the weak isospectral condition $(\mbox{H}_1)-(\mbox{H}_2)$;
	\item Each billiard table $X$ in ${\mathcal B}_{\rm an}$ with  boundary  sufficiently close to an ellipse $\Gamma \in \mathcal{E}$ in the $C^5$ topology is  spectrally rigid in the class ${\mathcal B}_{\rm an}$  under the weak isospectral condition $(\mbox{H}_1)-(\mbox{H}_2)$. 
		\end{enumerate}
\end{Th}
Theorems   \ref{theo:elliptic}  and \ref{theo:ellipse1} are proved in Section \ref{Sec:elliptic}.  

\begin{Remark}\label{rem:low-order-terms}
	To obtain the preceding results we use only the leading term of the quasi-mode.  	Replacing  the second condition of $(\mbox{H}_1)$ with a stronger one 
	\[
	\mathrm{lim}\,  a^{s/2}(b_k-a_k) =0 \quad \mbox{as} \quad k\to +\infty,
	\]
	where $s\in\N$ is fixed, and using an analogue 
	of Lemma 2.5 \cite{PT4}, one can  obtain further isospectral invariants, which could  be used to remove at least one of the symmetries.  
\end{Remark} 

Consider a strictly convex billiard table  $X$ in $\R^2$. Lazutkin has proved that 
any fixed $\tau>1$ and any $0<\kappa<\kappa_0\ll 1$ there is a subset $\Omega_\kappa\subset D(\kappa,\tau)$ of positive Lebesgue measure consisting $(\kappa,\tau)$-Diophantine frquences such that for any $\omega\in \Omega_\kappa$ there is a Kronecker invariant circle $\Lambda(\omega)$ of the billiard ball map $B$ of a frequency $\omega$ (see \cite{La} and the references there).  Morreover, the corresponding caustic ${\mathcal C}(\omega)$ - the envelope of the rays issuing from $\Lambda(\omega)$ - is a closed  smooth convex curve lying in the interior of $X$. As $\kappa$ tends to $0$ the invariant curves accumulate at the boundary $S^\ast(\Gamma)$ of ${\bf B}^\ast (\Gamma)$. 
\begin{Th}\label{Th:convex-main} Let $X_t\subset \R^2$, $t\in [0,1]$, be a $C^1$ family of  compact   billiard tables  satisfying the weak isospectral condition $(\mbox{H}_1)-(\mbox{H}_2)$. Suppose that $X_0$ is strictly convex. Then 
\begin{enumerate}
\item $X_t$ is strictly convex for each $t\in [0,1]$
\item There is a Cantor set $\Xi\subset (0,1]$ consisting of Diophantine numbers such that 
\[
{\rm meas}\, \left(\Xi\cap (0,\varepsilon)\right) = \varepsilon(1 - O(\varepsilon))\quad  \mbox{as}\quad \varepsilon\to 0^+ 
\]
and such that $\forall\, \omega\in \Xi$ there is a $C^1$ family of Kronecker invariant circles 
$[0,1]\ni t \to \Lambda_t(\omega)$ of $B_t$ of frequency $\omega$, 
\item $\forall \omega\in \Xi$ and $t\in [0,1]$, $\beta_t(\omega)= \beta_0(\omega)$, $L_t(I_0(\omega))=L_0(I_0(\omega))$, and $I_t(\omega)=I_0(\omega)$. 
 
\end{enumerate}
\end{Th}

\subsection{Proof of Theorem \ref{Th:main3}}\label{subsec:proof-Th3}
Theorem \ref{Th:main3} follows from Theorem \ref{Th:main2} and the following statement. 
\begin{Prop}\label{Prop:variation-alpha}
Let $[0,\delta]\ni t\to \Lambda_t(\omega)\subset{\bf B}^\ast\Gamma_t$ be a $C^1$ family of invariant tori of 
$P_t= B_t^m$ with a frequency $\omega\in \Omega_\kappa^0$. Then
\begin{equation}
\frac{d}{dt}\beta_t(\omega) = -\frac{2}{\mu_t(\Lambda_t(\omega))}\sum_{j=0}^{m-1}\int_{\Lambda_t(\omega)} 
\left( \delta_\nu \Gamma_t \right)(\pi_t(\zeta))\,\Big\langle \xi_t^+(\zeta),\nu_t(\pi_t(\zeta))\Big\rangle\Big|_{\zeta=B^j_t(\rho)}
\,d\mu_t(\rho)
\label{eq:variation-alpha}
\end{equation}
for each  $t\in [0,\delta]$ and   $\omega\in \Omega_\kappa^0$. 
\end{Prop}
We are going to use the following 
\begin{Lemma}\hspace{-2mm}{\bf .}\label{Lem:length-differentiation} 
Consider a $C^1$ family of curves $c : (-\ep,\ep)\times[0,l]\to X$, 
$c_s(\cdot)\equiv c(s,\cdot) : [0,l]\to X$, 
such that $c(\theta):=c_0(\theta)$ is a geodesic of the metric $g$ and denote by $l(s):=l_g(c_s)$  the length of the curve $c_s(\cdot) : [0,l]\to X$ with respect to the metric $g$. 
 Then 
\[ 
\frac{dl}{ds}(0) =  -g\left(\frac{{\dot c}_0(0)}{||{\dot c}_0(0)||_g},\frac{\partial c}{\partial s}(0,0)\right)+ 
g\left(\frac{{\dot c}_0(l)}{||{\dot c}_0(l)||_g},\frac{\partial c}{\partial s}(0,l)\right), 
\] 
where ${\dot c}_0(\theta):= \frac{dc}{d\theta}(0,\theta)$. 
\end{Lemma} 
{\em Proof.} 
The Lemma follows from a straightforward differentiation of the length function $l(s)$ and the Euler-Lagrange
equation. 
\finishproof

\noindent{\em Proof of Proposition \ref{Lem:length-differentiation}. } Using the notations introduced just before
Theorem \ref{Th:main2} we set 
$p(\theta,t):= f_t(\theta)$, $\theta\in \T^{n-1}$, where $f_t\equiv f_{t,\omega}: \T^{n-1}\to{\bf B}^\ast\Gamma_t$ 
is the embedding of the Kronecker torus $\Lambda_t(\omega)\subset{\bf B}^\ast\Gamma_t$.
Consider the function $x: \T^{n-1}\times [0,\delta] \to \Gamma_t$ given by
\begin{equation}
x(\theta,t):= \pi_{t}(p(\theta,t))=\psi\big(\pi_0({\tilde f}_t(\theta)),t\big). 
\label{eq:x}
\end{equation} 
where ${\tilde f}_t(\theta):=\psi_t^\ast(f_t(\theta))$.
Clearly, $x\in C^1(\T^{n-1}\times [0,\delta],{\widetilde X})$ and 
$x(\theta,0) = \pi_0(f_0(\theta))$ for any $\theta\in \T^{n-1}$. Suppose first that  $m=1$, that means that $P_t=B_t$.
By \eqref{eq:beta function} we have
\[
\beta_s(\omega) = -\frac{1}{(2\pi)^{n-1}}\int_{\T^{n-1}} l\Big(x(\theta,s),  
x\big(R_{2\pi\omega}(\theta),s\big)\Big)\,d\theta
\]
where $l(x,x')$ is the corresponding length function  which is well defined and smooth in a neighborhood of
the projection of $\Lambda_t(\omega)\times \Lambda_t(\omega)$ in ${\widetilde X}\times{\widetilde X}$. 

First we prove \eqref{eq:variation-alpha} for $t=0$.
Setting $\widetilde \theta := R_{2\pi\omega}(\theta)$  and using  Lemma \ref{Lem:length-differentiation} we get 
\[
\begin{array}{lcrr}
-\displaystyle \frac{d}{ds}\Big|_{s=0} \beta_s(\omega) = \frac{1}{(2\pi)^{n-1}} \int_{\T^{n-1}} 
\left[\frac{d}{ds} l \left(x(\theta,s),  x(R_{2\pi\omega}(\theta),s)\right)\right]\Big|_{s=0} d\theta \\[0.3cm]
\displaystyle= \frac{1}{(2\pi)^{n-1}} \int_{\T^{n-1}} \left[\left\langle \pi^+_{\Sigma_s}(p(\widetilde\theta,s)),
\frac{\partial x}{\partial s}(\widetilde\theta,s)\right\rangle - \left\langle \pi^-_{\Sigma_s}(p(\theta,s)), 
\frac{\partial x}{\partial s}(\theta,s)\right\rangle \right]\Big|_{s=0}  d\theta \\[0.3cm]
\displaystyle= \frac{1}{(2\pi)^{n-1}} \int_{\T^{n-1}} \left\langle \left(\pi^+_{\Sigma}-\pi^-_{\Sigma}
\right)(p(\theta,0)), \frac{\partial x}{\partial t}(\theta,0)\right\rangle d\theta \\[0.3cm]
\displaystyle= \frac{2}{(2\pi)^{n-1}} \int_{\T^{n-1}}  \left\langle \pi^+_{\Sigma}(f_0(\theta)), 
\nu(x(\theta,0))\right\rangle 
g\left(\frac{\partial\psi}{\partial t}\big(\pi_0(f_0(\theta)),0\big),\nu\big(\pi_0(f_0(\theta)\big)\right)d\theta .
\end{array}
\]
where $\pi^\pm_{\Sigma_s} : {\bf B}^\ast\Gamma_s\to\Sigma_s^\pm$ is the map defined in \eqref{eq:outgoing-vector}.
In the variables $(x,\xi)=f_0(\theta)$,  we get
\[
\displaystyle \frac{d}{ds}\Big|_{s=0} \beta_s(\omega) = -\frac{2}{\mu_0(\Lambda_0(\omega))} 
\int_{\Lambda_0(\omega)} \left\langle\xi^+(\xi),\nu(\pi_0(\xi))\right\rangle(\delta_\nu\Gamma)(\pi_0(\xi))\,d\mu_0, 
\]
where $\delta_\nu\Gamma(x) = \langle\frac{\partial\psi}{\partial t}(x,0),\nu(x)\rangle $. 
The same argument holds for $m\ge 2$. Finally, to prove \eqref{eq:variation-alpha} for any $t\in [0,\delta]$ 
we replace $\Gamma_0$ by $\Gamma_t$ and apply the same arguments.
\finishproof

\section{Birkhoff Normal Forms of $C^k$ deformations.} \label{Sec:BBM-BNF}

The aim of this section is to prove  items {\em 1} and {\em 2} of Theorem \ref{Th:main1}. We shall obtain 
a $C^k$, $k\in\{0;1\}$, family of Birkhoff Normal Forms for the $C^k$ family of exact symplectic maps $t\to P_t$ around the corresponding invariant tori. This BNF will be used  to construct a $C^k$ family of Quantum Birkhoff Normal Forms. In order to obtain the BNF we will apply Theorem \ref{Theo:BNF}  to a suitable $C^k$ family of exact symplectic mappings $\widetilde P_t$ which will be constructed below.  

Let $[0,1]\ni t\mapsto P_t\in C^\infty(U, U)$ be a $C^k$ family of exact symplectic maps where $U\subset T^\ast\Gamma$ is an open set. Suppose that $P_0$ has a Kronecker torus $\Lambda_0(\omega_0)$ with frequency $\omega_0\in D(\kappa_0,\tau)$ where $0<\kappa_0< 1$ and $\tau>n-1$. 
Then $P_0$ admits a Birkhoff Normal Form (BNF) at $\Lambda_0(\omega_0)$ (cf \cite{La}, Proposition 9.13 and \cite{PT4}, Proposition 3.3) which means the following. There exists an exact symplectic transformation $\widetilde \chi: \A \to T^\ast\Gamma$, where $\A:= \T^{n-1}\times D\subset T^\ast\T^{n-1} $ and $D$ is a neighborhood of $I_0(\omega_0)$ given by \eqref{eq:momentum-I}, $\widetilde \chi(\A)$ is a neighborhood of $\Lambda_0(\omega_0)$, $\widetilde \chi(\T^{n-1}\times \{I_0(\omega_0)\}=\Lambda_0(\omega_0)$ and the exact symplectic map $\widetilde P_0:= \widetilde \chi^{-1}\circ P_0 \circ \widetilde \chi$ has the form 
\begin{equation}\label{eq:BNF-torus}
\left\{
\begin{array}{lcrr}
\widetilde P_0 (\theta, r) \, =\, (\theta + \nabla K(r),r) \, +\, R(\theta,r),\\[0.3cm]
\partial_r^\alpha R(\theta,I_0(\omega_0))=0,\ \forall\, \theta\in \T^{n-1},\ \forall\, \alpha\in\N^{n-1}.
\end{array}
\right.
\end{equation}
We are going to use the following  definition of a generating function of a symplectic map.  Denote by $\rm{pr}:\R^{n-1}\to \T^{n-1}$ the canonical projection. 
\begin{Def}\label{def:generating-function}
	Let $D\subset \R^{n-1}$ be an open set and $F\in C^\infty(\T^{n-1}\times D)$. 
	The function $S\in C^\infty(\R^{n-1}\times D)$ given by $S(x,r)= \langle x,r\rangle - F(\rm{pr}(x),r)$ is said to be a generating function of a  symplectic map $P$ in $\T^{n-1}\times D$ if 
	\begin{itemize}
		\item the map 
		$x \to \nabla_r S(x,r)=x -\nabla_r F(\rm{pr}(x),r)$ projects to a diffeomorphism of $\T^{n-1}$ homothope to the identity for any fixed $r\in D$ 
		\item for any $(\theta,r)\in \T^{n-1}\times D$  
		\[
		P\big(\theta - \nabla_r F(\theta,r),r\big)\, =\, \big(\theta,r-\nabla_\theta F(\theta,r)\big) 
		\]
	\end{itemize}
\end{Def}
Hereafter,  we  slightly abuse the notations identifying $y=\nabla_r F(\theta,r)\in \R^{n-1}$ with its image $\rm{pr}(y)\in \T^{n-1}$. 

Shrinking $U$ if necessary we  set $U=\widetilde \chi(\A)$.
We suppose that the BNF is nondegenerate, which means that the Hessian matrix $\partial^2 K(r_0)$ is nondegenerate at $r_0=I_0(\omega_0)$. Shrinking $D$ if necessary we suppose that  the map $\nabla K: D \to \nabla K(D)\subset \R^{n-1}$ is a diffeomorphism.
Then there exists $\bar\delta>0$ such that the exact symplectic maps $\widetilde P_t= \widetilde \chi^{-1}\circ P_t \circ \widetilde \chi$ admit for $0\le t\le \bar\delta$ a $C^1$ family of generating functions $[0,\bar\delta]\ni t \rightarrow \widetilde G_t\in C^\infty (\R^{n-1}\times D)$ such that
\begin{equation}\label{eq:BNF-generating-functions1}
\widetilde G_t(x,r) = \langle x, r \rangle  - K(r) - G_t({\rm pr}(x),r)
\end{equation}
where the map $t\to G_t\in C^\infty(\T^{n-1}\times D)$ is $C^1$, ${\rm pr}: \R^{n-1} \to \T^{n-1}$ is the canonical projection  and 
\begin{equation}\label{eq:BNF-generating-functions2}
 \partial_r^\alpha G_0(\theta,I_0(\omega_0))=0\quad \forall\,  \theta\in \T^{n-1},\   \alpha\in\N^{n-1}.
\end{equation}  
On the other hand,
\[
\big|\partial_\theta^\alpha\partial_r^\beta \big( G(\theta,r)-G_0(\theta,r)\big)\big| \le C_{\alpha,\beta}\,  t\quad \forall\, (\theta,r)\in \T^{n-1}\times D(\kappa),\ t\in [0,\bar\delta]. 
\] 
These inequalities allow us to apply Theorem \ref{Theo:KAM}. 
 Consider the Legendre transform $K^\ast$ of $K$ in   a  ball $B(\omega_0,\varepsilon)$ centered at $\omega_0$ and with sufficiently small radius $0<\varepsilon\ll 1$. Then $\nabla K^\ast : B(\omega_0,\varepsilon)\to \nabla K^\ast (B(\omega_0,\varepsilon))$ becomes  a diffeomorphism. Using Theorem \ref{Theo:KAM} with $\kappa=\varrho =\kappa_1$, where $\kappa_1\le \kappa_0$ is sufficiently small, we obtain a $C^k$ family of 
 Kronecker tori $[0,\delta_1] \ni t \mapsto \Lambda_t(\omega_0)$ of $P_t$ with a frequency $\omega_0$, where
 $\delta_1=\delta_1(\kappa_1)>0$. Following Lazutkin  (cf \cite{La}, Proposition 9.13 and \cite{PT4}, Proposition 3.3) we obtain a $C^k$ family of Birkhoff Normal Forms of $\widetilde P_t$ around the tori $\Lambda_t(\omega_0)$ up to any order, which means the following. Fix $N\ge 4$. There exist $C^k$ families of exact symplectic mappings $t\to \chi_t^0\in C^\infty(\T^{n-1}\times D,\T^{n-1}\times D)$ and vector valued  functions $ t\to I_t(\omega_0)\in D$,  $t\in [0,\delta_1] $  (analytic in $t$ if the map $t\to P_t$ is analytic)  with the properties
 \begin{itemize}
 	\item 
  $\chi_t^0(\T^{n-1}\times \{I_t(\omega_0)\}) = \Lambda_t(\omega_0)$ for each $t\in [0,\delta_1] $;
  \item the exact symplectic map $\widetilde P_t^0= (\chi_t^0)^{-1}\circ \widetilde P_t \circ  \chi_t^0$ admits for each $0\le t\le \delta_1$ a  generating function
 \begin{equation}\label{eq:BNF-generating-functions3}
 \widetilde G_t^0(x,r) = \langle x, r \rangle  - K_t(r) - G_t^0({\rm pr}(x),r), 
 \end{equation}
such that  the maps $t\to G_t^0\in C^\infty(\T^{n-1}\times D)$ and $t\to K_t\in C^\infty(D)$ are  $C^k$ (analytic in $t$ if the map $t\to P_t$ is analytic)  and 
 \begin{equation}\label{eq:BNF-generating-functions4}
  \partial_\theta^\alpha \partial_r^\beta G_t^0(\theta,I_t(\omega_0))=0\quad \forall\,  \theta\in \T^{n-1},\   \forall\,   \alpha, \beta\in\N^{n-1}, \ \mbox{with}\ |\beta|\le 2N ,
 \end{equation}  
and for each $0\le t\le \delta_1$. Moreover, $K_t$ is a polynomial of degree $2N$ for each $t$ fixed. 
\end{itemize}

Let $\omega_0$ be  a point of positive Lebesgue density in $D(\kappa_0,\tau)$ ( $\omega_0\in D^0(\kappa_0,\tau)$ ), which means that the Lebesgue measure ${\rm meas\, }(D(\kappa_0,\tau)\cap V) > 0$ for any neighborhood $V$ of $\omega_0$ in $\R^{n-1}$. Then $\omega_0\in D^0(\kappa,\tau)$ for each $0<\kappa\le \kappa_0$. 
We  suppose that $\kappa_1\le \varepsilon^2$ and for every $0<\kappa\le\kappa_1<1$ we  set
\begin{equation}\label{eq:Omega-kappa}
\Omega(\kappa)\, :=\, B(\omega_0,\sqrt{\kappa}),\quad \Omega_\kappa\, :=\, D(\kappa,\tau)\cap B(\omega_0,\sqrt{\kappa}-\kappa),\quad \Omega_\kappa^0\, :=\, D^0(\kappa,\tau)\cap B(\omega_0,\sqrt{\kappa}-\kappa)
\end{equation}
It follows from \cite{La}, Proposition 9.9,  that 
\begin{equation}\label{eq:measure-omega-BNF}
\frac{{\rm meas}\, \big(\Omega(\kappa)\setminus\Omega_{\kappa}\big)}{ {\rm meas}\, (\Omega(\kappa)) }\,  \le \,  C\,  \kappa.
\end{equation}
Moreover, ${\rm meas}\, (\Omega_\kappa^0) ={\rm meas}\, (\Omega_\kappa) $. Set $D_t(\kappa):= \nabla K_t^\ast(\Omega(\kappa))$ and $\A_t=\T^{n-1}\times D_t(\kappa)$. 
Notice that there exists $c>0$ such that $D_t(\kappa) \subset B(I_t(\omega_0), c \sqrt{\kappa})$ and \eqref{eq:BNF-generating-functions4} implies that for every $\alpha,\, \beta\in \N^{n-1}$  there exists $C_{\alpha,\beta, N}$ such that
\[
\big|\partial_\theta^\alpha(\kappa\partial_r)^\beta G_t(\theta,r)\big| \le C_{\alpha,\beta,N} \, \kappa^N\quad \forall\, (\theta,r)\in \T^{n-1}\times D_t(\kappa),\quad t\in [0,\delta_1(\kappa_1)].
\]
This inequality allows one to apply  Theorem \ref{Theo:BNF} taking $\varrho=\varepsilon\kappa$ and $0<\kappa\le \kappa_1$, where $\varepsilon$ and $\kappa_1$ are sufficiently small in order  
In this way we obtain the following   
\begin{Theorem}\label{Theo:soft-BNF} (Birkhoff Normal Form)\quad Let $[0, \delta]\ni t\to P_t\in C^\infty (U, U)$, $U \subset T^\ast\Gamma$,  be a $C^1$ family of symplectic maps. Let $P_0$ have a Kronecker torus $\Lambda_0(\omega_0)$ with a frequency $\omega_0\subset D^0(\kappa_0,\tau)$, where  $\tau>n-1$. Suppose that the BNF of $P_0$ at $\Lambda_0(\omega_0)$ is nondegenerate. Then there exists  $0<\kappa_1\le \kappa_0$ and  $\delta_1=\delta_1(\kappa_1)>0$ such that  for each   $M>0$  fixed and $0<\kappa\le \kappa_1$ the following holds
\begin{enumerate}
\item[(i)]
 For each  $\omega\in \Omega_{\kappa}^0$ 
there exists a $C^1$ family of Kronecker invariant tori $[0,\delta_1] \ni t \to \Lambda_t(\omega)$ of $P_t$ with a frequency $\omega$;  
\item[(ii)] There exists a  
$C^1$-smooth with respect to $t\in [0,\delta_1]$ family of exact symplectic maps 
$\chi_t:\A_t \to U$ and of real valued functions  $L_t\in C^\infty(D_t(\kappa))$ and $R_t^0\in C^\infty(\A_t)$   (analytic in $t$ if the map $t\to P_t$ is analytic)  such that 
\begin{enumerate}
\item[1.] $\Lambda_t(\omega)=\chi_t(\T^{n-1}\times \{I_t(\omega)\})$ for each $t\in [0,\delta_1]$ and $\omega\in  \Omega_{\kappa}^0$, where $I_t(\omega)$ is given by \eqref{eq:momentum-I}; 
\item[2.] The function $\R^{n-1}\times D\ni (x,I)\mapsto \phi_t(x,I):=\langle x,  I\rangle -L_t(I) - R_t^0(x,I)$ is a generating function in the sense of Definition \ref{def:generating-function} of the exact symplectic map 
\[
P_t^0:= \chi_t^{-1}\circ P_t\circ \chi_t: \A\to\A ; 
\]
\item[3.] $\nabla L_t: D_t \to \Omega$ is a diffeomorphism,   $L_t=K_t$  outside $D_t^1:= \nabla K_t^\ast\big(B\big(\omega_0,\sqrt{\kappa}-\frac{1}{2}\kappa\big)\big)$  and  $\nabla L_t^\ast(\omega)=I_t(\omega)$ is given by \eqref{eq:momentum-I} for each $\omega\in \Omega_{ \kappa}^0$;
\item[4.] $R_t^0$ is flat at $\T^{n-1}\times E_t^{\kappa}$, where $E_t^{\kappa}:=\nabla L_t^\ast( \Omega_{\kappa}^0)$. 
 \item[5.]  $\|\nabla L_t - \nabla K_t\|_{m,D_t;\kappa}   + \| \nabla R_t^0\|_{m,D_t;\kappa} \le C_{m,M}\kappa^M$ for each $\alpha,\beta\in \N^{n-1}$ and $m\in \N$. 
\end{enumerate}
Moreover, if the map $t\to P_t$ is analytic in a disc $B(0,\delta)$ in $\C$, then the maps $t\to \chi_t$, $t\to L_t$, $t\to R^0_t$, are analytic in a disc $B(0,\delta_1)$ and the estimate in 5. holds for $t\in B(0,\delta_1)$. 
\end{enumerate}
\end{Theorem}
To prove the Theorem we apply Theorem \ref{Theo:BNF} taking $N\gg M$, $\varrho=\varepsilon\kappa$ and $0<\kappa\le \kappa_1$, where $\varepsilon$ and $\kappa_1$ are sufficiently small in order to satisfy \eqref{eq:smallness-condition-1}.

Observe that  each $I\in E_t^{\kappa}$ is an element of positive Lebesgue density of $E_t^{\kappa}$ since the map $\nabla L_t^\ast : \Omega(\kappa) \to D_t(\kappa)$ is a diffeomorphism. For any $0<\kappa\le \kappa_1$ fixed, we  extend $I_t$ as a $C^1$ family of smooth functions setting
\[
I_t(\omega):= \nabla L_t^\ast(\omega)\quad \forall\, \omega\in \Omega. 
\]
Then we have 
\[
P_t^0(\varphi,I)= (\varphi + \nabla L_t(I),I) + O_N(|I- I_t(\omega)|^N)
\] 
for each $\omega\in \Omega_\kappa^0$ and  $N\in\N$. This formula can be differentiated  with respect to $(\varphi,I)$ as many times as we want. To summarize we give the following
\begin{Def}\label{Def:BNF} We say  that the $C^1$ family of exact symplectic maps $P_t$, $t\in [0,\delta]$,  admits a $C^1$-smooth family of nondegenerate  Birkhoff Normal Forms associated with  a $C^1$ family of invariant tori $\Lambda_t(\omega)$ with  frequencies  $\omega\in \Omega_{\kappa}^0$ if item (ii) of 
Theorem \ref{Theo:BNF} holds true.  
\end{Def}
Recall as well that the complement of  $\Omega_\kappa^0$ in  $\Omega_\kappa$ is of Lebesgue measure zero.\\

\noindent Setting 
\begin{equation}\label{eq:Xi}
\Xi\, :=\, \displaystyle  \bigcup_{0<\kappa\le \kappa_1} \Omega_\kappa^0
\end{equation}
we prove  items {\em 1} and  {\em 2} of Theorem \ref{Th:main1}. 

The advantage of working with $\Omega_{\kappa}^0$ instead of $\Omega_{\kappa}$ is given by the following
\begin{Lemma}\label{Lemma:flat}
Let $\Omega$ be an open subset of $\R^d$, $d\ge 1$. Let $E\subset \Omega$ be a measurable set of positive Lebesgue measure and let $E^0\subset E$ be the set of points of positive Lebesgue density in $E$. Then  any 
 smooth function $f$ on $\Omega$ which  is zero on $E^0$ is flat at $E^0$, i.e. the equality $f|_{E^0}=0$ implies  $\partial^k f|_{E^0}=0$ for any $k\in\N^d$. Moreover, the Lebesgue measure of $E\setminus E^0$ is zero by construction.
\end{Lemma}
{\em Proof}. The result is evident when $d=1$. Suppose that $d\ge 2$. Let $\omega=(\omega_1,\omega')\in E^0$. By Fubini's theorem, for any neighborhood $U_1\subset \R$ of $\omega_1$ and $U'\subset \R^{d-1}$ of $\omega'$ there is $z'\in U'$ and a set of positive Lebesgue measure $V_1\subset U_1$ such that $(z_1,z')\in E^0$ for any $z_1\in V_1$. Then $f(z_1,z')=0$ for any $z_1\in V_1$ and there is $y_1\in U_1$ such that $\partial_1 f(y_1,z')=0$. By continuity we obtain $\partial_1 f(\omega)=0$. In the same way we prove   by induction that the restriction of $\partial^\alpha f$ to $ \Omega_E^0$ is zero for any $\alpha\in \N^{d}$. \finishproof

We are going to give a relation between the function $\beta_t$ defined by \eqref{eq:beta function} and the restrictions on $\Omega_\kappa^0$  of the functions $I_t$ and $L_t$ given by Theorem \ref{Theo:soft-BNF}. As $\chi_t:\A\to U\subset T^\ast\Gamma$ is exact symplectic for each $t\in J:=[0,\delta(\kappa)]$ there exists a $C^1$ family of functions $\Psi_t\in C^\infty(\A)$ such that
\[
\chi_t^\ast(\xi dx) \, =\, I d\varphi + d\Psi_t . 
\]
 Notice that the generating functions $\phi_t$ of $P_t^0$ are uniquely defined up to  additive constants $C_t$ such that the function $J\ni t\to C_t$ is $C^1$. 
\begin{Lemma}\label{lemma:beta} Choosing appropriately the $C^1$ function $J\ni t\to C_t$ we obtain the following 
\begin{enumerate}
	\item [(i)] We have 
	\[
	\langle I, \nabla L_t(I)\rangle  - L_t(I) = A_t(\chi_t(\varphi,I)) + \Psi_t(\varphi,I) - \Psi_t(P_t^0(\varphi,I)) + R_t^1(\varphi,I)
	\]
	where the function $R_t^1$ is flat at $\T^{n-1}\times E_t^\kappa$ for each $t\in J$;
	\item [(ii)] $\beta_t(\omega) + L_t(I_t(\omega)) = \langle \omega, I_t(\omega)\rangle \quad \forall\, \omega\in \Omega_\kappa^0,\ t\in J$. 
\end{enumerate}
\end{Lemma}
{\em Proof.} \quad The Poincar\'e identity implies  
\begin{equation}
P_t^\ast  (\xi dx) = \xi dx + dA_t,
\label{poincare}
\end{equation}
where $\xi dx$ is the fundamental one-form on $T^\ast\Gamma$ and
$A_t(\rho)$, $\rho =\chi(\varphi, I)\in U$ 
is the action 
\[
A_t(\rho) = \int_{\widetilde\gamma_{t}(\rho)}\, \xi dx
\]
along the broken bicharacteristic
$\widetilde \gamma_{t}(\rho)$ strating at $\pi_t^+(\rho)$ and ending at $\pi_t^-(P_t(\rho))$. Then we obtain 
$(P^0)^\ast  (I d\varphi) - I d\varphi = d((A\circ \chi) +  \Psi -
\Psi\circ P^0)$. On the other hande, $(P^0)^\ast  (I d\varphi) = d( L_t(I) - \langle I, \nabla L_t(I)\rangle + R_t^1(\varphi,I))$ where $R_t^1$ is a flat function at $\T^{n-1}\times E_t^\kappa$ and we obtain (i). 
To prove (ii) we use \eqref{eq:beta function}. 
\finishproof

\section{Infinitesimal spectral rigidity of Liouville billiard tables}\label{sec:LBT}
A Liouville billiard table  of dimension $n\ge 2$,  is a completely integrable billiard table $(X,g)$ admitting $n$
functionally independent and Poisson commuting integrals of the billiard flow on $T^*X$
which are quadratic forms in the momentum. It can be viewed as a $2^{n-1}$-folded branched covering of a disk-like
domain in $\R^n$ by the cylinder $ \T^{n-1}\times [-N,N]$, $N>0$. 

Liouville billiard tables of dimension two were defined in \cite[Sec. 2]{PT1} by using a branched double covering map.
Here we give an invariant definition of Liouville billiard tables in dimension two.
The equivalence of the two definitions is proven in Appendix \ref{sec:Appendix B}.

\begin{Def}\label{def:LBT_invariant_definition}
A Liouville billiard table is a smooth oriented compact and connected Riemannian manifold of dimension two 
$(X,g)$ with connected boundary $\Ga\equiv\partial X$ such that the following two conditions are satisfied:
\begin{itemize}
\item[$(a)$] There exists a smooth quadratic in velocities integral of the geodesic flow $I : TX\to \R$ that is 
invariant with respect to the reflection at the boundary $TM|_{\Ga}\to TM|_{\Ga}$, $\xi\mapsto \xi - 2g(\nu,\xi)$, 
where $\nu$ is the outward unit normal to $\Ga$. In addition, we assume that the metric $g$ does {\em not} allow global 
Killing  symmetries;
\item[$(b)$] There is no point $x_0\in\Ga$ and a constant $c\in\R$ such that $g_{x_0}(\xi,\xi)= c I_{x_0}(\xi,\xi)$
for any  $\xi\in T_{x_0}X$.
\end{itemize}
\end{Def}
\noindent In view of Theorem \ref{th:invariant_definition} in Appendix \ref{sec:Appendix B} there exists a double
covering map with two branched points,
\begin{equation}
\tau : C\to X,
\end{equation}
where $C$ denotes the cylinder ${(\R/\Z)}\times [-N, N]$, $N>0$, coordinatized by the variables $x$ and $y$
respectively, so that the metric $\tau^*(g)$ and the integral $\tau^*(I)$ have the following form on $C$,
\begin{eqnarray}\label{eq:the_metric_1}
dg^2&=&\big(f(x)-q(y)\big)(dx^2+dy^2)\\
dI^2&=&\alpha\,dF^2+\beta\,dg^2\nonumber
\end{eqnarray}
where $\alpha\ne 0$ and $\beta$ are real constants and
\begin{equation}\label{eq:the_integral_1}
dF^2:=\big(f(x)-q(y)\big)\big(q(y)\,dx^2+f(x)\,dy^2\big)\,.
\end{equation}
In other words, the integral $dI^2$ belongs to the pencil of $dg^2$ and $dF^2$.
Here $f\in C^\infty(\R)$ is 1-periodic, $q\in C^\infty([-N,N])$, and
\begin{itemize}  
\item[(i)] $f$ is  even, $f>0$ if $x\notin\frac{1}{2}{\Z}$, and        
$f(0)=f(1/2)=0$;  
\item[(ii)] $q$ is even, $q<0$ if $y\ne 0$,  $q(0)=0$ and $q^{''}(0)<0$;  
\item[(iii)] $f^{(2k)}(l/2)=(-1)^kq^{(2k)}(0)$,  $l=0,1$,    
for every natural $k\in{\N}$.  
\end{itemize}  
In particular, if  $f\sim \sum_{k=1}^{\infty}\ f_kx^{2k}$ is the Taylor expansion of $f$  
at $0$,  then,   by  (iii),  the Taylor expansion of $q$ at $0$  
is $q\sim \sum_{k=1}^{\infty}\ (-1)^k f_kx^{2k}$.  
A Liouville billiard table is said to be  {\em of classical type} if it satisfies the following additional
conditions,
\begin{itemize}  
\item[(iv)] the boundary $\Gamma$ of $X$  is   
locally geodesically convex which amounts  to $q'(N)<0$;   
\item[(v)] $f(x)=f(1/2-x)$ for any $x$ and $f$ is strictly increasing   
on the interval $[0,1/4]$;   
\end{itemize}    
\begin{Remark}
Note that in contrast to \cite[Sec. 2]{PT1} we do {\em not} assume that the functions $f$ and $g$ are analytic
Morse functions.
\end{Remark}
The points $F_1:=\tau(0,0)$ and $F_2:=\tau(1/2,0)$ on $X$ are the two branched points of the covering map
$\tau : C\to X$. All other points on $X$ are regular values of $\tau$. The preimage of any regular value consists of two
points. Note also that $\tau : C\to X$ commutes with the involution on the cylinder $C$ induced by the map
$\si : (x,y)\mapsto (-x,-y)$. The fixed points of this involution are precisely the singular points $(0,0)$ and $(1/2,0)$
of the covering map $\tau$. One can see that any Liouville billiard table  possesses the {\em string property} which 
means that any broken geodesic starting from the singular point $F_1\, (F_2)$ passes through 
$F_2\, (F_1)$ after the first reflection at the boundary and the sum of distances from any
point of $\Gamma$ to $F_1$ and $F_2$ is constant \cite{PT1}. In particular, the only Liouville billiard table 
in $\R^2$ equipped with the Euclidean metric is the interior of the  ellipse. Thus Liouville billiard tables can be 
regarded as a natural generalization of elliptic billiards to curved space. In view of condition $(v)$ in the
definition of the Liouville billiard tables of  classical type,  there is a group  $I(X)\cong{\Z}_2\oplus{\Z}_2$ 
acting on $(X,g)$ by isometries. This group is generated by the involutions $\sigma_1(x,y)= (x,-y)$ and 
$\sigma_2(x,y) = (\pi - x,y)$ of the cylinder.

The construction of Liouville billiards of dimension two involving the covering map $\tau$ was generalized to 
any dimension in \cite[\S 5.3]{PT2} (cf. also \cite[\S 3]{PT3}). As in the two dimensional case, one defines the
subclass of Liouville billiards of classical type in a similar way.
An important example of a  Liouville billiard table of classical type is the interior of the $n$-axial ellipsoid
in $\R^n$ equipped with the Euclidean metric. More generally, there is a non-trivial two-parameter family
of  analytic Liouville billiard tables of classical type  of constant scalar curvature $K$ having the same broken
geodesics (considered as non-parameterized curves) as the ellipsoid \cite[Theorem 3]{PT2}. 
This family includes the ellipsoid ($K=0$), and  Liouville billiard tables on the sphere ($K=1$) and 
in the hyperbolic space ($K=-1$). 

In what follows we will apply Theorem \ref{Th:main3} to Liouville billiard tables of classical type in dimensions two 
and three for obtaining several new isospectral results. 
The main idea is to interpret the integrals in \eqref{eq:integral-invariants} as values of a suitable Radon transform which is 
one-to-one. 

Let $(X,g)$ be a Liouville billiard table of classical type. The  group of isometries of $(X,g)$ has a subgroup $I(X,g)$ 
isomorphic to $(\Z/2\Z)^n$.  One can extend $(X,g)$ to an open Riemannian manifold $(M,g)$ so that any isometry in 
$I(X,g)$ can be extended to an isometry of $(M,g)$. In this way, the group of isometries of $(M,g)$ contains a subgroup 
$I(M,g)$ isomorphic to $I(X,g)$. 
Denote by ${\rm Symm\, }(M,g)$ the class of $C^\infty$-smooth billiard tables $(Y,g)$, $Y$ isometrically embedded in
$M$, so that any isometry in $I(M,g)$ is an isometry of $(Y,g)$. 
Recall from Sec. \ref{Sec:results} that $\delta_\nu \Gamma:\Gamma\to \R$ is the vertical component of the variation 
$\dot \psi_0:\Gamma\to T M|_{\Gamma}$, where $\psi_0=id$ and $\Gamma_t=\psi_t(\Gamma_0)$ is a 
$C^1$-deformation of $\Gamma=\Gamma_0$ and $(X_t,g)$ is the billiard table with boundary $\Gamma_t=\partial X_t$. 

\begin{Theorem}\label{Th:Liouville-dim2}
Let $(X,g)$ be a Liouville billiard table of classical type dimension 2 and let $(X_t,g)$, $t\in (-\varepsilon,\varepsilon)$,  
be a $C^1$-family of billiard tables in ${\rm Symm\, }(M,g)$ satisfying the weak isospectral condition 
$(\mbox{H}_1)-(\mbox{H}_2)$ and such that $X_0=X$.   Then $\delta_\nu \Gamma\equiv 0$.
\end{Theorem}
This means that any Liouville billiard table of classical type $(X,g)$ is {\em infinitesimally spectrally rigid} in ${\rm Symm\, }(M,g)$ under the weak-isospectral condition $(\mbox{H}_1)-(\mbox{H}_2)$. \\

\noindent{\em Proof of Theorem \ref{Th:Liouville-dim2}. } The theorem follows from Theorem \ref{Th:main2} as in 
the proof of Corollary 1.4 in \cite{PT4}.
A first integral of $B$ in $B^\ast \Gamma$ is the function ${\cal I}(x,\xi) = f(x)-\xi^2$ the regular values $h$ of which 
belong to  $(q(N),0)\cup (0,f(1/4))$  (see \cite{PT1}, Lemma 4.1 and Proposition 4.2). Moreover, for each regular value 
$h\in (q(N),0)$ the corresponding level set $L_h: \{\mathcal I=h\}$ consists of two connected circles which are invariant
with respect to $B$, having rotation numbers  $\omega=\pm \rho(h)$.  By  Proposition 4.4 \cite{PT1}, the rotation function 
$\rho$ is smooth and strictly increasing in an interval  $(q(N), q(N) +\varepsilon) $, and we obtain a diffeomorphism 
$\rho: (q(N), q(N) +\varepsilon) \to (0,\omega_0)$. Then the Kolmogorov non-degeneracy condition is satisfied in that 
interval. Hence, one can apply Corollary \ref{Th:main2} to any Kronecker invariant circle $\Lambda_0(\omega)$,  
with a Diophantine vector of rotation  $\omega=\rho(h)\in (0,\omega_0)$. 

We are going to interpret \eqref{eq:integral-invariants}  as a value of a suitable Radon transform. 
The Leray form on the invariant circle $\Lambda_0(\omega)\subset L_h$ is 
\[
 \lambda_h \ =\ 
 \left\{   
 \begin{array}{ccc} 
 \frac{dx}{\sqrt{f(x)-h}},\ \xi > 0\, ,\\  
-\frac{dx}{\sqrt{f(x)-h}}, \ \xi < 0\, .
\end{array}
\right.
\]
Since  $\omega\in (0,\omega_0)$ is Diophantine and 
 the  Leray form  is invariant with respect to $B$, there exists a constant $c(h)\neq 0$ such that $\lambda_h=c(h)d\mu_0$, where $d\mu_0$ is the unique probability measure on $\Lambda_0(\omega)$ which is invariant with respect to $B$. 
Setting $K:=\pi_{\Gamma_t}^\ast \left( \delta_\nu \Gamma_t \right)$ we
consider the corresponding Radon transform which  assigns to each circle $\Lambda_0(\omega):=\{(x,y(h)): x\in\T\}$, $h=\rho^{-1}(\omega)$,  the integral
\[
R_K(\Lambda_0(\omega)) = \int_{\Lambda_0(\omega)}\, \langle \xi^\pm,\nu\rangle \, K\circ \pi_\Gamma\,  \lambda_h\, .
\]
We have
\[
\langle \xi^\pm(x,y(h)),\nu(x)\rangle \ =\ \sqrt{\frac{h-q(N)}{f(x)-q(N)}}\ ,
\]
hence,
\[ 
R_{K}(\Lambda_0(\omega))  = \pm c(h)\sqrt{h-q(N)}  
\int\limits_0^{1}\, \frac{K_1(x)}{\sqrt{f(x)-h}}\;dx  \ ,
\ h\in (q(N),0)\, ,
\] 
where $K_1(x)= K(x)/\sqrt{f(x)-q(N)}$.
Since $K$ is invariant with respect to the group of isometries $I(X)$ then so is $K_1$ and applying \eqref{eq:integral-invariants} for $t=0$ we get
\begin{equation}
0=\int\limits_0^{1/4}\, \frac{K_1(x)}{\sqrt{f(x)-h}}\;dx =  \int\limits_0^{f(1/4)}\, \frac{K_2(s)}{\sqrt{s-h}}\;dx\, 
\label{eq:radon}
\end{equation}
for any $h\in (q(N),q(N)+\varepsilon)$ such that $\rho(h)\in D(\kappa,\tau)$, where $K_2 = (K_1/f')\circ f^{-1}\in L^1(0, f(1/4))$. 
On the other hand,  the right hand side of \eqref{eq:radon} is analytic in $h\in (q(N),0)$ and the set of $h=\rho^{-1}(\omega)$, $\omega \in D(\kappa,\tau)\cap(0,\omega_0)$ is of positive measure, and we obtain  
(\ref{eq:radon}) for any $h\in (q(N), 0)$. Differentiating (\ref{eq:radon}) with respect to $h$ at $h= q(N)$ we get 
\[
\int\limits_0^{f(1/4)}\, \frac{K_2(s)}{\sqrt{s-q(N)}}(s-q(N))^{-k}\;ds=0
\]
for any $k\in\N$, which implies $K_2=0$ since the set $\{(s-q(N))^{-k}:\, k\in\N\}$ is dense in $L^1(0, f(1/4))$ and $K_2$ is continuous in $(0, f(1/4))$. This completes the proof of the theorem. 
\finishproof

In order to apply Theorem \ref{Th:main2} to Liouville billiard tables of dimension three we need to ensure that
the following Kolmogorov nondegeneracy condition is satisfied. 
Consider a Liouville-Arnold chart which consists of an open set $U$ of the phase space of 
the billiard ball map $B$ and a symplectic map $(\varphi,I):U\to \T^n\times D$, $D$ being an open subset of $\R^n$,  
giving ``action-angle'' coordinates on $U$, i.e. $B$ is given by the map $(\varphi,I)\mapsto (\varphi + \nabla K(I),I)$ in 
these coordinates. The Kolmogorov  condition means that the map $\nabla K: D \to D^\ast:= \nabla K(D)$ is a 
diffeomorphism. We are interested in maximal charts with this property. It turns out that Liouville billiard tables of 
classical type of dimension two are always non-degenerate in a Kolmogorov sense close to the boundary \cite{PT1}. 
The non-degeneracy property of  Liouville billiard tables of classical type of dimension three has been investigated in 
\cite{PT3}. 

It is proved in \cite{PT3} that any Liouville Billiard Table of classical type of dimension 3 admits
four not necessarily connected maximal Liouville-Arnold charts $U_j$, $1\le j\le 4$, of action-angle variables in $B^\ast\Gamma$.
Two of them, say $U_1$ and $U_2$, have the property that any unparameterized geodesic in $S^\ast\Gamma$ can be
obtained as a limit of orbits of $B$ lying either in $U_1$ or in $U_2$
(then the corresponding broken geodesics approximate geodesics of the boundary). Moreover, in any connected component of $U_1$ and $U_2$ there is such a sequence of orbits of $B$, while any orbit of $B$ in $U_3$ and $U_4$ is far away from $S^\ast\Gamma$. In other words,   the charts  $U_1$ and $U_2$  can be  characterized by the property that there is a family of {\em ``whispering gallery rays''}  issuing  from any of their connected components.  For this reason the  two cases $j=1,2$ are referred  as to boundary cases. 
Denote by ${\mathcal F}_b$ the set of all regular invariant tori $\Lambda\in {\mathcal F}$ lying either in $U_1$ or in $U_2$. 
 We say that a Liouville Billiard Table is Kolmogorov nondegenerate if $B$ satisfies the Kolmogorov nondegeneracy condition in $U_1$ and $U_2$. It is shown in \cite{PT3}, Theorem 5.1,   that \emph{any  analytic $3$-dimensional Liouville billiard table of classical type having  at least one non-periodic geodesic on the boundary is Kolmogorov nondegenerate}. An example of such billiard tables is the ellipsoid. 

It is proved in \cite{PT3}, Theorem \ref{Th:Problem_B}, that the Radon transform is one-to-one for  Liouville billiard tables  of classical type  of dimension $3$. More precisely, we have 
\begin{Theorem}\label{Th:Problem_B}\cite{PT3}
	Let $(X,g)$, $\mbox{dim}\, X=3$, be a Liouville billiard table of classical  type. 
	If $K\in C(\Gamma)$ is invariant under the group of symmetries $G$ of $\Gamma$ and the Radon transform ${\mathcal R}_{K}(\Lambda) = 0$ for any  $\Lambda\in {\mathcal F}_b$, then  $K\equiv 0$. 
\end{Theorem}
We point out that Liouville billiard tables  of classical type are smooth by construction but \emph{they are not supposed to be analytic}. 

Using Corollary \ref{Th:main2}  and Theorem \ref{Th:Problem_B}, we obtain as above the following 
\begin{Theorem}\label{Th:Liouville-dim3}
Any nondegenerate Liouville billiard table of dimension 3 of classical type $(X,g)$ is infinitesimally spectrally rigid in ${\rm Symm\, }(M,g)$ under the weak-isospectral condition $(\mbox{H}_1)-(\mbox{H}_2)$.
\end{Theorem}
This means that if $(X_t,g)$, $t\in (-\varepsilon,\varepsilon)$ is a $C^1$-family of billiard tables in ${\rm Symm\, }(M,g)$ satisfying the weak isospectral condition $(\mbox{H}_1)-(\mbox{H}_2)$ and such that $X_0=X$, then  $\delta_\nu \Gamma\equiv 0$.

A smooth deformation $(X_t,g)$, $t\in (-\varepsilon,\varepsilon)$, is said to be  flat at $t=0$ if \eqref{eq:variation-boundary} is $C^\infty$ smooth with respect to $t$ in an interval $(-\varepsilon,\varepsilon)$ and  the vertical component of $k$-th variation $\delta_\nu^k \Gamma$ is zero for any integer $k\ge 1$, where
\[
\delta_\nu^k \Gamma(x):= \left\langle\left(\frac{d}{dt}\right)^k\,  \psi_t\big|_{t=0}(x),\nu(x)\right\rangle\, , \quad x\in \Gamma. 
\]
\begin{Coro}\label{Th:Liouville-flat-dim2}
Let $(X,g)$ be a classical Liouville billiard table of dimension 2 or a classical non-degenerate Liouville billiard table of dimension 3 and let $(X_t,g)$, $t\in (-\varepsilon,\varepsilon)$,  be a $C^\infty$-family of billiard tables in ${\rm Symm\, }(M,g)$ satisfying the weak isospectral condition $(\mbox{H}_1)-(\mbox{H}_2)$ and such that $X_0=X$. 
Then the deformation is flat at $t=0$. In particular, $\Gamma_t=\Gamma_0$ for $t\in(-\varepsilon,\varepsilon)$ if the family is analytic with respect to $t$.
\end{Coro}
In the case of the ellipse similar results have been obtained by Hezari and Zelditch \cite{H-Z} under the usual isospectral condition using the wave-trace method. 

\section{Isospectral deformations in the presence of elliptic geodesics.}\label{Sec:elliptic}
Let  $(X,g)$ be a smooth billiard table with  boundary $\Gamma:=\partial X$.   Consider  a $C^1$-smooth family 
\begin{equation}
[0,1] \ni t \to (X,g_t)
\label{eq:perturbation}
\end{equation}
of Riemannian metrics on $X$ with $g_0=g$.  Suppose that $(X,g)$ admits 
 an elliptic closed broken geodesic $\gamma$ 
with $m\ge 2$ vertices. 
Denote by $\{B^{j}(\rho):\ 0\le j\le m-1\}$ the corresponding periodic trajectory of the billiard ball map $B$. Then $ \rho= (x,\xi)\in {\bf B}^\ast\Gamma$ is a fixed point of the local Poincare map $P=B^m$ which is symplectic. Recall that $\gamma$ is  {\em elliptic} if $ \rho$ is an {\em elliptic} fixed point of $P$ which means that the eigenvalues  of the linear Poincaré map $dP( \rho):T_{ \rho}\Gamma \to T_{ \rho}\Gamma$ are all distinct, different from one and of modulus one, hence, 
\[
{\rm Spec}\, (dP( \rho))=\{e^{\pm  i\phi_j}:\,  1\le j\le n-1\}, 
\]
where $0 <\phi_1< \cdots <\phi_{n-1}\le \pi$. 
Set $\phi=(\phi_1,\ldots,\phi_{n-1})$ and fix a positive  integer $ N$. The linear Poincaré map  $dP(\rho)$ is said to be {\em $N$-elementary}   if the scalar product 
\[
 \langle \phi,k\rangle \, \notin \, 2\pi \Z
\]
for any integer vector $k=(k_1,\ldots,k_{n-1})\in \Z^{n-1}$  such that $0<|k_1|+\cdots+|k_{n-1}|\le N$. We say as well that $\gamma$ admits no resonances of order less or equal to $N$. 

From now on we fix $N\ge 4$ and suppose that $\gamma_0=\gamma$ is elliptic in $(X,g_0)$ and that it admits no resonances of order less or equal to $N$. 
By the implicit function theorem  there exists $\bar\delta>0$ such that the following holds.  There is a unique $C^1$ curve $[0,\bar\delta) \ni t \mapsto \rho_t\in T^\ast\Gamma$ starting from  $\rho_0=\rho$ such that $\rho_t\in {\bf B}^\ast_t\Gamma$ is an elliptic  fixed point of $P_t= B_t^m$  for any $t\in [0,\bar\delta)$. Moreover, the linear Poincare map $dP_t(\rho_t)$ is $N$-elementary. The eigenvalues of $P_t$ have the form $\exp(i \phi_{j}(t))$, where $0 <\phi_{1}(t)< \cdots <\phi_{n-1}(t) <\pi$ and the map $t\to \phi(t):= (\phi_{1}(t), \ldots,\phi_{n-1}(t))$ is $C^1$ in $[0,\bar\delta)$. 
Moreover,  $P_t$ admits a Birkhoff normal form of order $[N/2]\ge 2$ ($[a]$ denotes the integer part of $a\in\R$)  in suitable polar symplectic coordinates which will be described below.

In order to avoid eventual singularities at $r_j=0$, $j=1,\ldots,n-1$,  we fix  $0<c_0\ll 1$ and $r_0>0$, and set 
\begin{equation}\label{eq:D-b}
\D=\D(c_0):=\{r=(r_1,\ldots,r_{n-1})\in\R^{n-1}:\ 0 < c_0|r|< |r_j| < r_0,\  1\le j\le n-1\} \, .
\end{equation}
and  $\A:=\T^{n-1}\times \DD$. Denote by $\rm{pr}:\R^{n-1}\to \T^{n-1}$ the canonical projection.

\begin{Prop}\label{Prop:BNF-elliptic} (Birkhoff Normal Form). 
For any $0<\delta<\bar\delta$ there exists 
\begin{itemize}
	\item[-] a $C^1$-family of exact symplectic transformation $[0,\delta]\ni t\to \big(\widetilde \chi_t: \A \to U_t:=\chi_t(\A)\big)$ , where $U_t\subset {\bf B}^\ast_t\Gamma$ is an open set
   \item[-]  a $C^1$-family of polynomials $K_t\in \R_{[N/2]}[\xi_1,\ldots,\xi_{n-1}]$ with real coefficients of $n-1$ variables $\xi_1,\ldots,\xi_{n-1}$ and of degree $\left[\frac{N}{2}\right]$
   \item[-] a $C^1$-family of  real valued functions $G_t\in C^\infty(\A)$
\end{itemize}
such that the following holds
\begin{enumerate}
	\item the function $\widetilde G_t\in C^\infty(\R^{n-1}\times \D)$ defined  by $\widetilde G_t(x,r):= \langle x,r\rangle - K_t(r)-G_t(\rm{pr}(x),r)$ is a generating function of the symplectic map $\widetilde P_t:=\widetilde \chi_t^{-1}\circ P_t\circ\widetilde \chi_t$ in $\T^{n-1}\times \DD$ , 
	\item for any  $\alpha,\beta\in \N^{n-1}$  there exists $C_{\alpha,\beta}>0$ such that 
\[
|\partial_\theta^\alpha \partial_r^\beta  G_{t}(\theta,r)|\, \le \, C_{\alpha,\beta}\,  |r|^{\frac{N+1}{2}-|\beta|} 
\]
for any $t\in [0,\delta]$, $(\theta,r)\in \A$, and
	\item $\lim_{r\to 0}\widetilde\chi_t(\theta,r)=\rho_t$. 
\end{enumerate} 
\end{Prop}
The first condition is satisfied for $N\ge 2$ and $|r|\le r_0\ll 1$ in view of the estimate in {\em 2.} and the inverse function theorem.  Multiplying  $G_t$ by a  smooth  cut-off function of the form $f_0(r)=f(r/r_0)$ where $f$ is compactly supported  in the unit ball and $f(r)=1$ for $|r|\le 1/2$  we obtain a smooth function with support contained in the ball of radius $r_0$. Notice that the estimates of {\em 2} still hold for the function $G_tf_0$ with constants $C_{\alpha,\beta}$ depending on $f$ but not on $r_0$. 

The polynomial  $K_t$  is the so called Birkhoff polynomial, 
\[
\nabla K_t(0) = \phi(t)
\]
and 
$(r,\theta)$ are local polar symplectic coordinates. The construction of the Birkhoff normal form  follows from that of Moser \cite{M3} (see also \cite{Klin}, Lemma 3.3.2). 

The  Birkhoff normal form of $P_t$ is said to be {\em nondegenerate} if   the Hessian of $K_t$ at $0$ does not vanish, i.e. $\det \partial^2 K_t(0)\neq 0$. 
We say as well that $P_t$ is a {\em twisted map} at $\rho_t$ in this case. Suppose now that $P_0$ is twisted. Choosing $\delta>0$ sufficiently small we obtain by continuity that $P_t$ is twisted for any $t\in [0,\delta]$, i.e. 
\begin{equation}\label{eq:twisted}
  \det \partial^2 K_t(0)\neq 0 \quad \forall\, t\in [0,\delta].
\end{equation}
Then $\nabla K_t:\D \to \D^\ast_t:= \nabla K_t(\D)$ is a diffeomorphism for any $t\in [0,\delta]$ provided that $r_0\ll 1$. Denote by $Q_t:\A\to\A$ the corresponding non-degenerate completely integrable map defined by $Q_t(\theta,r):=(\theta+\nabla K_t(r),r)$.
The  set of frequency vectors of $Q_t$ is $\Omega_t:=\D^\ast_t$. This is an open cone-like set in $\R^{n-1}$ with vertex at $\phi(t)$.

The remaining of this Section is devoted to the spectral rigidity of the Kronecker  tori in a vicinity of an elliptic geodesic. We address the following questions. Suppose that the $C^1$ family of billiard tables $(X,g_t)$, $0\le t\le 1$, is weakly isospectral. Assume that  $(X,g_0)$ admits a periodic  elliptic broken geodesic $\gamma_0$ and that the corresponding local Poincaré map is twisted. Does there exist a $C^1$ family of periodic  elliptic broken geodesics $[0,1]\ni t\to \gamma_t$ in $(X,g_t)$ along the whole perturbation? Do the corresponding local Poincaré map remain twisted?  Do the invariant tori $\Lambda_0(\omega)$ associated to $\gamma_0$ give rise to $C^1$ families of invariant tori $[0,1]\ni t\to \Lambda_t(\omega)$ along the whole perturbation? We give an answer of these questions in the following Theorem. 

 Denote by 
$\B(\alpha,\epsilon)=\B^{n-1}(\alpha, \epsilon)$ the ball of radius $\epsilon$ and center $\alpha$ in $\R^{n-1}$. Recall that the functions $\beta_t(\omega)$, $I_t(\omega)$ and $L_t(I)$ are defined by \eqref{eq:beta function}, \eqref{eq:momentum-I} and \eqref{eq:birkhoff-L} respectively.
\begin{Theorem}\label{theo:main-elliptic}
Let $(X,g_t)$, $t\in [0,1]$,  be a $C^1$ family of billiard tables of dimension $n\ge 2$  satisfying the weak isospectral condition $(\mbox{H}_1)-(\mbox{H}_2)$. Suppose that $(X,g_0)$ admits a closed  elliptic   broken billiard trajectory $\gamma_0$ with $m\ge 2$ vertices. Suppose as well that the corresponding linear Poincare map $dP(\rho_0)$ is $N\ge 8$ elementary and that  $P=B^m$  is twisted at $\rho_0$. Then there exists $\delta_0>0$ such that for any interval $I=[0,\delta]$, $0<\delta<\delta_0$   the following holds.
\begin{enumerate}
	\item[(i)] There exist a $C^1$-family of   elliptic  fixed points $I\ni t\to \rho_t\in {\bf B}^\ast_t\Gamma$  of $P_t=B_t^m$,  the corresponding linear Poincaré map $dP_t(\rho_t)$ is $N$-elementary and $P_t$ is twisted at $\rho_t$.  Moreover for any $t\in I$ and  $|\alpha| \le \frac{N}{4} - 1$ 
\begin{equation}
\partial^\alpha\nabla K_t(0)=\partial^\alpha\nabla K_0(0).
\label{eq:Birkhoff-polynomials}
\end{equation}
\item[(ii)] There is a set $\Xi$ of positive Lebesgue measure consisting of Diophantine frequencies such that 
\[
\lim_{\epsilon \searrow 0}\,  \frac{{\rm meas}(\Xi\cap \B(\phi(0),\epsilon))}{{\rm meas}(\B(\phi(0),\epsilon))}\, =\,  0 
\]
and 
for any $\omega\in\Xi$ there is a $C^1$ family of Kronecker invariant tori $I\ni t\to \Lambda_t(\omega)\subset {\bf B}^\ast_t\Gamma$ of $P_t$ of a frequency  $\omega$.
\item[(iii)]  
 $\beta_t(\omega)=\beta_0(\omega)$, $L_t(I_0(\omega))=L_0(I_0(\omega))$, and $I_t(\omega)=I_0(\omega)$ for any  $t\in I$ and $\omega\in\Xi$.
\end{enumerate}
\end{Theorem}
It follows from \eqref{eq:Birkhoff-polynomials} that the function $[0,\delta_0)\ni t\to  \phi(t)=\nabla K_t(0)$ is constant, i.e. 
\begin{equation}
\forall\, t\in [0,\delta_0),\quad  {\rm Spec}\, (dP_t( \rho_t))={\rm Spec}\, (dP_0( \rho_0)).
\label{eq:spectrum-of-dP}
\end{equation}
A natural question  is to describe the largest interval $I$,  if it exists,  for which  Theorem \ref{theo:main-elliptic} holds. The answer is given by 
\begin{Prop}\label{prop:interval} Let $N\ge 12$. Then there are only two possibilities that may occur. 
\begin{enumerate}
	\item[(i)]
the conclusion $(i)-(iii)$ of Theorem \ref{theo:main-elliptic} holds with $I=[0,1]$
\item[(ii)] there is $0<\delta_0\le 1$ such that Theorem \ref{theo:main-elliptic} holds in any interval  $I=[0,\delta]$ with $0<\delta<\delta_0$ and  the limit set 
\[
\Sigma_{\delta_0}:= \displaystyle\lim_{t\to\delta_0}\{\rho_t, B(\rho_t), \cdots, B^{m-1}(\rho_t)\}
\] 
intersects the boundary of $B_{\delta_0}^\ast\Gamma$.  
\end{enumerate}
\end{Prop}
The case $(ii)$  means that there is a generalized (glancing  to the boundary at certain point)  geodesic  issuing from the limit set $\Sigma_{\delta_0}$. If the billiard tables $(X,g_t)$, $t\in[0,1]$, are locally strictly geodesically convex then each generalized geodesics of $(X,g_t)$  lies entirely on $\Gamma$ and the second case can not occur. \\

We are going to prove Theorem \ref{theo:elliptic}, Theorem \ref{theo:main-elliptic} and Proposition \ref{prop:interval}. 
Firstly using  Theorem \ref{Theo:BNF} we will obtain 
 a KAM theorem for the $C^1$ family of symplectic maps $\widetilde P_t$ given by Proposition \ref{Prop:BNF-elliptic}. 
To this end  we will determine the convex set $\Omega$, fix the parameters $\kappa$ and $\varrho$, and then estimate the corresponding quantities ${\mathcal B}_\ell$ which appear in Theorem \ref{Theo:BNF}. 

Consider the $C^1$ family of exact symplectic mappings $\widetilde P_t$ in $\A=\T^{n-1}\times \D$ with generating functions 
\begin{equation}\label{eq:generating-function-elliptic}
\widetilde G_t(x,r):= \langle x,r\rangle - K_t(r)-G_t(\rm{pr}(x),r),\quad (x,r)\in \R^{n-1}\times \D,
\end{equation} 
given by Proposition \ref{Prop:BNF-elliptic}. Recall that $t\to K_t$ is a $C^1$-family of polynomials with real coefficients of $n-1$ variables and of degree $\left[\frac{N}{2}\right]$, i.e. $K_t\in \R_{[N/2]}(\xi_1,\ldots,\xi_{n-1})$, while $t\to G_t\in C^\infty(\A)$ is a $C^1$-family of  real valued functions with support in $\B(0,r_0)$ such that
\begin{equation}\label{eq:estimates-G-elliptic}
|\partial_\theta^\alpha \partial_r^\beta  G_{t}(\theta,r)|\, \le \, C_{\alpha,\beta}\,  |r|^{\frac{N+1}{2}-|\beta|} 
\end{equation}
for any $t\in [0,\delta]$, $(\theta,r)\in \A$,
and  $\alpha,\beta\in \N^{n-1}$. 

There exists  a constant $A\ge 1$ such that 
\begin{equation}
\forall\, t,s\in [0,\delta],\ \|K_t - K_s\|_{C^{[N/2]}} \le A|t-s|,
\label{eq:K-Lipschitz}
\end{equation}
where the norm is taken in $C^{[N/2]}(\B(0,1))$.  
The map $P_0$ is twisted, then by continuity  $P_t$ remains twisted for any $t\in [0,\delta]$ provided that $\delta>0$ is sufficiently small. Choosing $\delta>0$ small enough,  there exists $\varepsilon>0$ such that the  Legendre transform $K_t^\ast$ of $K_t$ given by  \eqref{eq:L-transform} is well defined in  $\B(\phi(0),\varepsilon)$ and 
\[
\nabla K_t^\ast: \B(\phi(0),\varepsilon) \to  V_t:= \nabla K_t^\ast\big(\B(\phi(0),\varepsilon)\big)
\] 
is a $C^1$ family of diffeomorphisms with respect to $t\in [0,\delta]$, where $V_t$ is a neighborhood of $0$. Moreover,  the corresponding inverse maps are 
$\nabla K_t: V_t \to  \B(\phi(0),\varepsilon)$, hence, $\nabla K_t\circ \nabla K_t^\ast = {\rm id}$ on  $\B(\phi(0),\varepsilon)$  for any $t\in [0,\delta]$. In particular the inverse map of $d\nabla K_t(0):\R^{n-1}\to \R^{n-1}$ is $d\nabla K_t^\ast(\phi(t))$. 

We are ready to define  suitable convex sets of frequencies $\Omega$. Choose $e=(e_1,\ldots,e_{n-1})\in \R^{n-1}$ such that
\begin{equation}\label{eq:vector-e-norm}
2c_0<|e_j|<\frac{1}{2n},\quad j=1,\ldots,n-1,
\end{equation}
where $0<c_0<1/4n$ is fixed in \eqref{eq:D-b}
and set $e^\ast:= d\nabla K_t(0)e$. 
Given  $0<a_0<\varepsilon$ and $0<\eta_0<1$  we consider for any $0<a\le a_0$ the cube of center 
$\phi(t) + ae^\ast$  with sides of length $2\eta_0 a$ defined by 
\[
\Omega= \Omega(t,a):=\Big\{\omega\in\R^{n-1}:\, |\omega_j-\phi_j(t)-ae_j^\ast| < \eta_0 a,\, 1\le j\le n-1 \Big\},
\]
Obviously, $\overline \Omega(t,a) \subset \B(\phi(0),\varepsilon)$ for $a_0\ll \varepsilon$, hence, $\nabla K_t^\ast$ is well defined and smooth on the convex set  $\overline \Omega(t,a)$. Denote by $I_{t,a}$ the set of all $s\in [0,\delta]$ such that  $|t-s|\le \eta_0 a$. Denote by $\D_a$ the connected component of the set 
\[
\Big\{r=(r_1,\ldots,r_{n-1})\in\R^{n-1}:\, c_0 < | r_j| <   a/n,\,  1\le j\le n-1\Big\}
\]
containing $ae$. Then $\D_a$ is a {\em convex} open set and $ae\in \D_a$. 
We claim that there exist $0<a_0< 1$ and  $0<\eta_0<1$ such that for any $0<a\le a_0$, $t\in [0,\delta]$ and  $s\in I_{t,a}$ the following relation holds
\begin{equation}\label{eq:inclusion}
\D^s(t,a)\, :=\, \nabla K_s^\ast(\Omega(t,a)) \, \subset \, \D_a \, \subset \, \D\, ,
\end{equation}
where $\D$ is defined by \eqref{eq:D-b}. 
Indeed, for any $\omega\in \Omega(t,a)$, using Taylor's formula up to order three for the function $\omega\to \nabla K_t^\ast(\omega)$ at $\omega=\phi(t)$ and the identity  $\nabla K_t^\ast(\phi(t))= \nabla K_t^\ast(\nabla K_t(0))=0$ we obtain
\[
\begin{array}{rcll}
|\nabla K_s^\ast(\omega)-ae|&\le& A \eta_0 a + |\nabla K_t^\ast(\omega)-a(d \nabla K_t^\ast)(\phi(t))e^\ast|\\[0.3cm] 
&\le& a C_n\big((A+B) \eta_0  + B a\big) 
\end{array}
\]
where $C_n$ depends only on $n$, $A>0$ is the constant in \eqref{eq:K-Lipschitz} and  
\[
B:= 1+\sup_{0\le t\le \delta}\| \nabla K_t^\ast\|_{C^2(\B(\phi(t), \varepsilon)}.
\]
Then the inclusion $\D^s(t,a) \subset \, \D_a$ follows from \eqref{eq:vector-e-norm} choosing $\eta_0$ and $a_0$ so that 
\[
C_n\big((A+B) \eta_0  + B a_0\big) < c_0<1/2n. 
\]
On the other hand, the inequalities $c_0a < | r_j| <   a/n$, $j=1,\ldots,n-1$,  imply $c_0|r|<c_0 a< | r_j|$ and we obtain the second inclusion in \eqref{eq:inclusion}, which proves the claim.

 Set $\A^s(t,a):=\T^{n-1}\times \D^s(t,a)$ for $s\in I_{t,a}$. The relation \eqref{eq:inclusion} allows one to apply Proposition \ref{Prop:BNF-elliptic} in $\D^s(t,a)$ for any $t\in [0,\delta]$ fixed, where the parameter of the deformation $s$ varies in $I_{t,a}$. We point out that both $\Omega(t,a)$ and $\D_a$ are {\em convex} open sets which allows us to apply Theorem \ref{Theo:BNF}  and to obtain the corresponding H\"{o}lder estimates. 

Fix $\tau>n-1$  and  choose $\kappa = \eta a$ 
in the Diophantine condition \eqref{eq:sdc}, where $0<\eta<\eta_0$. Denote by $\Omega_{t,\kappa}$ the set of all   $\omega\in \Omega(t,a)\cap D(\tau,\kappa)$ such that ${\rm dist}(\omega,\R^{n-1}\setminus \Omega)\ge \kappa$. 
There exists  $0<\eta_1 =c(n,\tau) \eta_0$, where $0<c(n,\tau)<1$ depends only on $n$ and $\tau$ such that the Lebesgue measure of $\Omega_{t,\kappa}$ is positive for any $t\in [0,\delta]$, $0<\eta <\eta_1$ and $a\in(0,a_0]$. Indeed, it follows from \cite{La}, Proposition 9.9 that 
\begin{equation}\label{eq:measure-omega}
{\rm meas}\, (\Omega(t,a)\setminus\Omega_{t,\eta a})\,  \le \,  C\, \frac{\eta}{\eta_0} \, {\rm meas}\, (\Omega(t,a)) \, ,
\end{equation}
where the positive constant $C$ depends only on $n$ and $\tau$, and we take $c=1/C$.  Let us  fix $0<\eta<\eta_1$ and denote by $\Omega_{t,\kappa}^0$ the set of points of positive Lebesgue density in $\Omega_{t,\kappa}$ (see Sect. \ref{Subsec:BNF-for-maps}).
\begin{Theorem}\label{theo:KAM-elliptic} 
Let $[0,\delta]\ni t\to \widetilde P_t$ be a $C^1$ family of symplectic mappings with generating functions $\widetilde G_t$ given by \eqref{eq:generating-function-elliptic}, where $K_t\in \R_{[N/2]}[\xi_1,\ldots,\xi_{n-1}]$ while $G_t\in C^\infty(\A)$ satisfies \eqref{eq:estimates-G-elliptic} with $N\ge 4$. 
Then for any $t\in [0,\delta]$ there exists  a  $C^1$-family  of exact symplectic maps 
\[
I_{t,a}\ni s \to \big(\chi_s: \A^s(t,a) \to \A^s(t,a)\big)
\]
 and of real valued functions  $L_s\in C^\infty(\D^s(t,a))$ and $R_t\in C^\infty(\A^s(t,a))$  such that for any $s\in I_{t,a}$ the following holds
\begin{enumerate}
\item $\widetilde G^0_s(x, I)= \langle x,I\rangle -L_s(I) - R_s({\rm pr}(x),I)$ is a generating function of $P_s^0:= \chi_s^{-1}\circ \widetilde P_s\circ \chi_s$ 
\item $\nabla L_s:  \D^s(t,a) \to \Omega(t,a)$ is a diffeomorphism
\item $R_s$ is flat at $\T^{n-1}\times\nabla L_s^\ast( \Omega_{t,\kappa}^0)$ 
\item  for any  $m\in\N$ there exists a constant  $C_m>0$ independent of  $a\in(0,a_0]$ and $t\in [0,\delta]$ such that the following estimates hold
\[
|\partial_\varphi^\alpha(a\partial_\omega)^\beta \sigma_a^{-1}(\chi_s - {\rm id})| + |\partial_\varphi^\alpha(a\partial_\omega)^\beta \sigma_a^{-1}(\chi_s^{-1} - {\rm id})| \le C_m  a^{\frac{ N-3}{4}}
\]
on $\A^s(t,a)$, and 
\begin{equation}\label{eq:estimates-birkhoff}
|(a\partial_I)^\beta  (\nabla L_s(I) - \nabla K_s(I))| \le C_m a^{{\frac{ N+1}{4}}}
\end{equation}
on $\D^s(t,a)$ for any $s\in I_{t,a}$ and $|\alpha|+|\beta|\le m$. 
\end{enumerate}
Moreover, if $\nabla K_t(0) = \nabla K_0(0)$ for any  $t\in[0,\delta]$, then for any $0<a\le a_0$ there exists  a  $C^1$-family  of exact symplectic maps 
$\chi_s:\A^0(0,a) \to \A^0(0,a)$ in $I=[0,\delta]$ and of real valued functions  $L_s\in C^\infty(\D^0(0,a))$ and $R_t\in C^\infty(\A^0(0,a))$  such that  1.) - 4.) hold for any $s\in I$ and $t=0$. 
\end{Theorem}

\noindent
{\em Proof.} To prove the first part of the Theorem  we apply Theorem  \ref{Theo:BNF} to the $C^1$ family of symplectic mappings $I_{t,a}\ni s\to \widetilde P_s\in  C^\infty(\A^s(t,a), \A^s(t,a))$.  

Let us  estimate the corresponding quantities ${\mathcal B}_l$ for $\ell\ge 1$ defined by \eqref{eq:B-0-sequence}-\eqref{eq:B-sequence}. First of all the constant $\lambda$ in \eqref{eq:non-degeneracy-BNF}   can be fixed by
\[
\lambda =  \sup_{t\in [0,\delta]} \|\partial^2 K_t\|_{C^{[N/2]}(\B(0,1))} .
\]
 Given $\ell=m+\mu$ with $m\in\N_\ast$ and $0\le \mu<1$ we get by \eqref{eq:inclusion} that for any $s\in I_{t,a}$  
\[
\|G_s\|_{\ell, \A^s(t,a);\kappa} \le \|G_s\|_{\ell, \T^{n-1}\times \D_a;\kappa} \le \|G_s\|_{m+1, \T^{n-1}\times \D_a;\kappa}.
\]
The second inequality follows from the fact that $\D_a$ is convex. On the other hand,
\[
\|G_s\|_{m+1, \T^{n-1}\times \D_a;\kappa} = \sup_{|\alpha|+ |\beta|\le m+1}\, \|\partial_\theta^\alpha (\kappa\partial_r)^\beta G_t\|_{C^0(\T^{n-1}\times \D_a)}
\]
and using Proposition \ref{Prop:BNF-elliptic}, {\em 2}, we obtain
$\|G_s\|_{\ell, \A^s(t,a);\kappa}\,  \le\, C_m \, a^{\frac{N+1}{2}}$ 
for any $s\in I_{t,a}$. 
We choose  $\varrho=  \kappa^{\frac{N+1}{4} }<\kappa$, where $N\ge 4$. Then 
\[
\|G_s\|_{\ell, \A^s(t,a);\kappa}\,  \le\,  \kappa\varrho \, C_m' \,a^{\frac{N-3}{4}}
\]
for any $s\in I_{t,a}$, where $C_m'= C_m \eta^{-\frac{N+5}{4}}$. Moreover, 
\[
|\!|\!| \partial^2 K_s |\!|\!|_{\ell, \A^s(t,a);\kappa} \, \le\,  |\!|\!| \partial^2 K_s |\!|\!|_{\ell, \T^{n-1}\times \B(0,1);\kappa} \le  C\| \partial^2 K_t \|_{C^{[N/2]}(\B(0,1))}  
\]
and 
\[
S_{\ell}(\nabla K^{\ast})\, \le \, \,\sup_{0\le t\le \delta}\, \big(1+ \|\nabla K^{\ast}_t\|_{C^1(\B(\phi(0),\varepsilon))}\big)^{\ell-1} \big(1+ \|\nabla K^{\ast}_t\|_{C^\ell(\B(\phi(0),\varepsilon))}\big)  .
\]
Thus for any $\ell\ge 1$ we obtain
\[
{\mathcal B}_\ell\, \le \,  \kappa\varrho \, C_\ell \,a^{\frac{N-3}{4}}
\]
where $C_\ell>0$ depends neither on $a\in (0,a_0]$ nor on  $t\in [0,\delta]$. Choosing $a_0\ll 1$ we get ${\mathcal B}_2\le \epsilon \kappa\varrho \lambda^{-4}$, which gives \eqref{eq:smallness-condition-1}. Applying Theorem  \ref{Theo:BNF} we obtain {\em 1}-{\em 4}.

The equality $\nabla K_t(0) = \nabla K_0(0)$ means that $\phi(t)=\phi(0)$ which implies  $\Omega(t,a)=\Omega(0,a)$. Then one can take $I_{t,a}=[0,\delta]$ in \eqref{eq:inclusion} which implies {\em 1}-{\em 4} in $[0,\delta]$.
\finishproof
\\

\noindent 
{\em Proof of  Theorem \ref{theo:elliptic}.}\quad  The set of fixed points of ${\mathcal J}_1$ in $X_t\in {\mathcal B}$ defines a bouncing ball geodesic $\gamma_{t}$ in $(X_t,g)$ which is preserved by ${\mathcal J}_2$.
We are going to apply Theorem \ref{theo:KAM-elliptic} to the Birkhoff Normal Forms of the local Poincaré maps $P_t$ associated to $\gamma_{t}$. 

 Let $t\to \rho_t$ be a $C^1$ family of fixed points of $P_t=B_t^2$. Denote by $ \rho_{t,j}=(x_{t,j},0) = B_t^{j}(\rho_t)$, $j= 1,\, 2$,  
the corresponding periodic orbit   of $B_t$.
Fix $t\in [0,\delta]$.    Denote by $U\subset \widetilde X$ a neighborhood  of the vertices $x_{t,0}$ and $x_{t,1}$ of $\gamma_{t,1}$ such that ${\mathcal J}_k(U)=U$, $k=1,2$. Denote the restrictions of the two involutions to $\Gamma_t\cap U$  by $J_1$ and $J_2$ and
 by $\widetilde J_j:T^\ast(\Gamma_t\cap U)\to T^\ast(\Gamma_t\cap U)$ the corresponding lifts.  The set $\Gamma_t\cap U$ has two connected components $\Gamma_{t}^j$, $j=1,\, 2$,  and   $J_1(\Gamma_{t}^1)=\Gamma_{t}^1$ while  $J_2(\Gamma_{t}^1)=\Gamma_{t}^2$. Since  $J_1$ and $J_2$ act as isometries and commute with each other, using the definition of $B_t$ in Sect. \ref{subsec:billiard-ball},  we obtain that the involutions $\tilde J_j$, $j=1,2$,  commute with each other and also with $B_t$. 
 
Denote by $\Xi$ the union of $\Omega_{t,\eta a}$,   $0<a\le a_0$.  For any $\omega\in \Xi$ we set  $\Lambda_t^1(\omega) =\Lambda_t(\omega)$  and $\Lambda_t^2(\omega) =B_t(\Lambda_t(\omega))$. 
Then $\widetilde J_1(\Lambda_t^j(\omega))$, $j=1,2$, are also invariant circles of $P_t= B_t^2$ of frequency  $\omega\in \Xi$  and  
 $\Lambda_t^j(\omega) =\widetilde J_1(\Lambda_t^j(\omega))$ for $j=1,2$, while $\Lambda_\omega^2 =\widetilde J_2(\Lambda_t^1(\omega))$.
To prove it we use the following argument. Since $\mbox{dim}\, T^\ast\Gamma_{t,j}=2$ the KAM circle    $\Lambda_t^j(\omega)$ divides $T^\ast\Gamma_{t,j}$ into two connected components,  and it contains the elliptic fixed point $ \rho_{t,j}=(x_{t,j},0)$ of $P_t$ in its interior $D_j$. Moreover, $\widetilde J_1( \rho_j)= \rho_j$,  hence, $\widetilde J_1(\Lambda_t^j(\omega))$ contains 
$ \rho_{t,j}$ in its interior $\tilde J_1(D_j)$  as well. 
On the other hand, $\tilde J_1$ preserves the volume form of $T^\ast\Gamma_{t,1}$, hence, $\Lambda_t^j(\omega)$ intersects  $\widetilde J_1(\Lambda_t^j(\omega))$. This implies  $\Lambda_t^1(\omega) =\widetilde J_1(\Lambda_t^1(\omega))$, since $P_t$  acts  transitively on both of them. In the same way we prove that  $\Lambda_t^2(\omega) =\widetilde J_2(\Lambda_t^1(\omega))$.

Recall that the family $\Gamma_t$, $t\in [0,1]$, is given by a $C^1$ family of embeddings $\psi_t\in C^\infty(\Gamma,\widetilde X)$, where $\psi_t(\Gamma)=\Gamma_t$. 
Without loss of generality we suppose that $\gamma_0= {\rm id}_\Gamma$ is the identity at $\Gamma$.  
Notice that that the vectors $\frac{\partial \psi_t}{\partial x}(x)$ and $\nu_t(\psi_t(x))$ provide a base of $T_{\psi_t(x)}\widetilde X$ for any $x\in\Gamma$, hence, 
\begin{equation}
\forall\, x\in \Gamma_t, \quad \dot{\psi}_t(x)=\lambda(t,x)\frac{\partial \psi_t}{\partial x}(x) + \delta_\nu\Gamma_t(\psi_t(x))\nu_t(\psi_t(x))
\label{eq:horizontal-vertical-component}
\end{equation}
where $t\to \lambda(t,\cdot)\in C^\infty(\Gamma)$ is continuous on $[0,\delta]$ and    the function $\delta_\nu \Gamma_t$ is defined by \eqref{eq:vertical-component} and it belongs to $C^\infty(\Gamma_t)$. We are going to show that the function 
$\delta_\nu \Gamma_t$ is flat at $x_{t,1}$.

Using Corollary \ref{Th:main3} and the symmetry with respect to $\widetilde{\mathcal J}_1$ given above we obtain
\begin{equation}
\int_{\Lambda_t(\omega)}\,    \langle \xi_t^+(\rho),\nu_t(\pi_{t}(\rho))\rangle\, \delta_\nu \Gamma_t (\pi_{t}(\rho))\,  d\mu_t(\rho)\ =\ 0
\label{eq:theequality}
\end{equation}
for any $\omega\in\Xi$, with $0<a\le a_0$. Moreover, the functions
$$
\rho\to  f_t(\rho):=\langle \xi_t^+(\rho),\nu_t(\pi_{t}(\rho))\rangle\, ,\quad \rho\to h_t(\rho):= \delta_\nu \Gamma_t (\pi_{t}(\rho))
$$ 
are invariant with respect to the involution $\widetilde{\mathcal J}_1$. 
Let us parametrize $\Gamma_{t}^1$ by its arclength $y\in [-c,c]$ so that $y(x_{t,1})=0$ and denote by $(y,\eta)$ the corresponding local coordinates in $T^\ast\Gamma_{t,1}$. Then  $J_1(y)=-y$ for any $y$.   For any invariant circle  $\Lambda_t^1(\omega)$, $\omega\in \Xi$, there is $y(\omega) >0$ such that    $\pi_{t}(\Lambda_t^1(\omega)) = [-y(\omega),y(\omega)]$.  Notice that $f_t(y,\eta)=\langle \xi_t^+(y,\eta),\nu_t(y)\rangle >0$ for $(y,\eta)\in \Lambda_t^1(\omega)$ since 
$\Lambda_t^1(\omega)$ is contained in ${\bf B}^\ast_t(\Gamma_t)$. On the other hand, $h_t(y,\eta)=h_t(y)$ depends only on $y$. 
We are going to show that there exists an infinite sequence $(y_j)_{j\in \N} \subset (0,c)$ such that $\lim y_j=0$ and  $h_t(y_j)=0$.
Indeed, suppose that   $h_t(y)\neq 0$ in $(0,b)$ for some $b>0$. Take $\omega\in \Xi$  such that $0<y(\omega)<b$. 
The function $h_t(y)$ is {\em even} because it is invariant with respect to ${\mathcal J}_1$, hence  it will  not change its sign in the interval 
$[-y(\omega),y(\omega)]$.  
Then $h_t(y,\eta)f_t(y)$ will not change its sign on $\Lambda_t^1(\omega)$ and it is not identically null, which contradicts (\ref{eq:theequality}). This proves the existence of an infinite sequence $\{y_j\}_{j\in\N}$ such that $h_t(y_j)=0$, $y_j\neq 0$ for any $j\in\N$ and $\lim y_j=0$.  Now there exists an infinite sequence $(y_j')_{j\in \N} \subset (0,b)$ such that $y_j\le y_j'\le y_{j+1}$ and $\frac{dh_t}{dy}(y_j')=0$, and so on. This implies that the Taylor series of  $h_t(y)$ vanishes at $y=0$. In particular we obtain that $x_{t,j}=x_{0,j}$ since $\lambda(t,x_{t,j})=0$. Hence, the function  $\Gamma_t\ni x\to\delta_\nu \Gamma_t(x)$ is flat at $x=x_{0,1}$.

Take local coordinates $x:\Gamma_0^1\to \R$ in the neighborhood $\Gamma_0^1$ of $x_{0,1}$ in $\Gamma=\Gamma_0$ such that $x(x_{0,1})=0$  
and consider the equation
\begin{equation}
\frac{d}{dt}u_t(x) = -\lambda(t,u_t(x))
\label{eq:equation-w}
\end{equation}
with initial data $u_0(x) = x$. This problem has a unique solution $u_t(x)$ for $t$ in a neighborhood of $0$ and $x$ in an open interval 
$V\subset \R$ containing  $x=0$. Moreover, $u_t: V\to \R$ is a $C^1$ family of local diffeomorphisms.    
Consider the $C^1$ family of  embeddings  $v_t=\psi_t\circ u_t:V\to \widetilde X$.  The set  $v_t(V)$ is an open neighborhood of $x_{0,1}$ in $\Gamma_t$, $v_t$ gives a local parametrization of $\Gamma_t$ in $v_t(V)$ and $v_0(x)=x$. 
Using \eqref{eq:horizontal-vertical-component} and \eqref{eq:equation-w} one obtains  that the map
\[
V\ni x \to \dot{v}_t(x)= \delta_\nu\Gamma_t(v_t(x))\nu_t(v_t(x))
\]
is flat at $x=0$ for any $s\in [0,\delta]$. Then  for any $\varphi\in C^\infty(\widetilde X)$, the function
\[
V\ni x\to\varphi(v_t(x))-\varphi(x)= \int_0^t d\varphi(v_s(x))\dot v_s(x)ds
\]
is flat at $x=0$ which means 
 that $\Gamma_t$ is tangent to infinite order to $\Gamma_0$ at $x_{0,1}$ for $t>0$ sufficiently small. Replacing $\Gamma$ by $\Gamma_t$, $t\in [0,1]$, we complete the proof of the Theorem.  
\finishproof

\noindent 
\emph{Proof of Theorem \ref{theo:ellipse1}.}  Corollary \ref{coro:ellipse} implies    that the Poincar\'e map associated with the elliptic bouncing ball geodesic $\gamma_1$ is always non-degenerate (twisted)  for elliptical billiard tables. Fix  the foci $F_1\neq F_2$. Except of five confocal families of ellipses given explicitly by \eqref{eq:rotation_ellipse'}, the  geodesic $\gamma_1$ is $4$-elementary.    
The two conditions are open in the $C^5$ topology and the Theorem follows from Theorem \ref{theo:elliptic}. 
\finishproof

\noindent 
{\em Proof of Theorem \ref{theo:main-elliptic}.} \quad It follows from Theorem \ref{theo:KAM-elliptic}  that for any $0<a\le a_0\ll 1$ and 
 $\omega\in \Omega_{t,\kappa}$ with $\kappa=\eta a$  there is a $C^1$ family of Kronecker invariant tori $I_{t,a}\ni s\to \Lambda_s(\omega)$ of $P_s$. 
 Moreover, Corollary \ref{Th:main2} implies that $I_s(\omega)=I_t(\omega)$ and $\nabla L_s(I_t(\omega))= \nabla L_t(I_t(\omega))$, hence,  $\partial^\alpha_I\nabla L_s(I_t(\omega))= \partial^\alpha_I \nabla L_t(I_t(\omega))$ for any $\alpha\in \N^{n-1}$ in view of Lemma \ref{Lemma:flat}. 
Then for any $\omega\in \Omega_{t,\kappa}$, $s\in I_{t,a}$ and  $\alpha\in\N^{n-1}$ of length $|\alpha|\le N/4-1$  using \eqref{eq:estimates-birkhoff} we obtain
\[
\begin{array}{rcll}
|\partial^\alpha_I \nabla K_s(I_t(\omega)) - \partial^\alpha_I \nabla K_t(I_t(\omega))| 
&\le& |\partial^\alpha_I \nabla K_s(I_t(\omega)) - \partial^\alpha_I \nabla L_s(I_t(\omega))|\\[0.3cm]
&+& |\partial^\alpha_I \nabla K_t(I_t(\omega)) - \partial^\alpha_I \nabla L_t(I_t(\omega))|\\[0.3cm]
&\le&  C a^{\frac{N+1}{4}- |\alpha|}\le C a^{\frac{5}{4}} .
\end{array}
\]
Taking $s=t+a$ there is  $t(a)\in [t,t+a]$ such that
\[
\left|\frac{d}{ds}\Big|_{s=t(a)}\partial^\alpha_I \nabla K_s(I_t(\omega))\right| \le C  a^{\frac{1}{4}},
\]
where the positive constant $C$ is independent of $a\in(0,a_0)$ and $\omega\in \Omega_{t,\kappa}^0$. Let  $a\to 0$.  
Then  $\Omega_{t,\kappa}^0\subset \Omega(t,a)$ shrinks to $\phi(t)$, hence, 
\[
\lim_{a\to 0}\, \sup \{|I_t(\omega)|:\, \omega\in \Omega_{t,\kappa}^0\} =0
\] 
and we get 
\[
\frac{d}{dt}\partial^\alpha \nabla  K_t(0) = 0
\]
for any $t\in [0,\delta]$. This implies \eqref{eq:Birkhoff-polynomials} in that interval.
By assumption $N\ge 8$, hence \eqref{eq:spectrum-of-dP} holds  for $t\in [0,\delta]$. 

We are going to define the set of frequencies $\Xi$ as a union of $\Omega_{0,\kappa}^0$. Recall that $\kappa=\eta a$, where $0<\eta\le  \eta_1=c(n,\tau)\eta_0$. Moreover, $\Omega(0,a)$ depends on the choice of $e^\ast:= d\nabla K_t(0)e$, where $d\nabla K_t(0)$ is an isomorphisme of $\R^{n-1}$ and $e$ satisfies \eqref{eq:vector-e-norm}, hence, $\Omega(0,a)$ depends  as well as on the parameter $0<c_0<1/4n$ defined in \eqref{eq:D-b}. Now varying the parameters $0<c_0<1/2n$, $e$ satisfying \eqref{eq:vector-e-norm},  $0<a\le a_0$ and $0<\eta\le \eta_1=c(n,\tau)\eta_0$ we denote by   $\Xi$ the union of the corresponding sets $\Omega_{0,\kappa}^0$. The set of frequencies  $\Xi$ satisfies $(ii)$ by construction in view of \eqref{eq:measure-omega}. 
Using the second part of Theorem \ref{theo:KAM-elliptic} and Corollary~\ref{Th:main2}  in $I=[0,\delta]$ we complete the proof of the Theorem.  \finishproof\\

\noindent
{\em Proof of Proposition \ref{prop:interval}.}\quad 
Denote by $0<\delta_0\le 1$ the supremum of all $\delta>0$ such that $(i)$, Theorem \ref{theo:main-elliptic} holds in $[0,\delta]$. 
  Suppose that $\delta_0<1$ and that $\Sigma_{\delta_0}\subset {\bf B}^\ast_{\delta_0}\Gamma$. Fix $\rho$ in the limit set $\displaystyle\lim_{t\to \delta_0}\rho_t$. Then $P_{\delta_0}= B^m_{\delta_0}$ is well-defined and smooth in a neighborhood of $\rho$. By continuity, $\rho$ is a fixed  point of $P_{\delta_0}$ and  \eqref{eq:spectrum-of-dP}  holds true for $t\in [0,\delta_0]$.  Hence, $\rho$ is an elliptic fixed point  of $P_{\delta_0}$ and there are no resonances of order less or equal to $N$.   Now  Proposition \ref{Prop:BNF-elliptic} provides a $C^1$ family of Birkhoff normal forms of  $P_{t}$ in an interval $t\in [0,\delta_0+\varepsilon]$, where $\varepsilon>0$. On the other hand, 
\eqref{eq:Birkhoff-polynomials} implies that $\partial^2 K_t(0)=\partial^2 K_0(0)$ for $t\in [0,\delta_0[$ since $N\ge 12$. By continuity, this equality is true for $t\in [0,\delta_0]$.  Then $P_t$ is twisted for any $t\in [0,\delta_0+\varepsilon]$, provided that $\varepsilon>0$ is sufficiently small. 
Applying Theorem \ref{theo:KAM-elliptic}   we show as above that $(i)$ and $(ii)$ in  Theorem \ref{theo:main-elliptic} hold in 
$[0,\delta_0 +\epsilon]$.   This contradicts the choice of $\delta_0$. If $\delta_0=1$ and $\Sigma_{1}\subset {\bf B}^\ast_{\delta_0}\Gamma$, then (i) holds in $[0,1]$ and  Theorem \ref{theo:KAM-elliptic}  holds in $I=[0,1]$. 
 \finishproof

\section{Isospectral deformation of  locally strictly 
geodesically convex billiard tables of dimension two.}\label{Sec:convex}

The aim of this Section is to prove Theorem \ref{Th:convex-main}. More generally we consider isospectral deformations of a billiard table $(X,g)$   in an ambient Riemannian manifold $(\widetilde X, g)$ with a locally strictly geodesically convex (with respect to the outward normal) boundary  $\Gamma=\partial X$. This means that if a geodesic $s\to \gamma(s)$ of $(\widetilde X, g)$ is tangent to  $\Gamma$ at $s=0$ then the order of the tangency is exactly two and $\gamma(s)\notin X$  for $0<|s|\ll 1$. The behavior of the billiard ball map near $S^\ast \Gamma$ is investigated by Melrose \cite{Mel} and Marvizi and Melrose \cite{MM} in the more general context of pairs of glancing surfaces. 

Consider the hypersurfaces $\Sigma_1:=S^\ast \widetilde X$ and $\Sigma_2:= T^\ast \widetilde X\big|_{\Gamma}$ in $T^\ast \widetilde X$. Set $f_1= h-1$, where the Hamiltonian $h$ is just the Legendre transform $h(x,\xi)= g^{ij}(x)\xi_i\xi_j$ of the Riemannian metric defined locally by $g(x,v)= g_{ij}(x)v^i v^j$ and denote by $f_2\in C^\infty(T^\ast \widetilde X)$ a smooth function which is constant on the fibers ($f_2(x,\xi)=f_2(x)$) and such that $f_2(x)>0$ for $x$ in the interior of $X$, $f_2(x)<0$  in the exterior  of $X$ and $f_2(x)=0$, $df_2(x)\neq 0$ for $x\in\Gamma$. Then the hypersurfaces $\Sigma_j$, $j=1,2$, are just the zero level sets of the non-degenerate Hamiltonians $f_j$ ($df_j\neq 0$ on $\Sigma_j=\{f_j=0\}$). One can show that 
 $\Gamma$ is locally strictly geodesically convex with respect to the outward normal to $\Gamma$ if and only if the following relation holds
\begin{equation}\label{eq:glancing}
f_1(\varrho)=f_2(\varrho)=\{f_1,f_2\}(\varrho)=0 \quad \Longrightarrow \quad \{f_1,\{f_1,f_2\}\}(\varrho)<0 \ \mbox{and} \ \{f_2,\{f_2,f_1\}\}(\varrho)>0 ,
\end{equation}
where $\{,\}$ is the Poisson bracket related to the canonical  symplectic two-form $\widetilde \omega$ of $T^\ast\widetilde X$. In particular 
\[
{\mathcal K}:= \{\varrho\in T^\ast \widetilde X:\, f_1(\varrho)=f_2(\varrho)=\{f_1,f_2\}(\varrho)=0\}
\]
is a smooth submanifold of $T^\ast \widetilde X$ of co-dimension two. The characteristic foliations of the two-form $\widetilde \omega\big|_{\Sigma_{j}}$ given by the non-parametrized integral curves of the hamiltonian vector fields of $f_{j}$ define two involutions ${\mathcal J}_{j}$ in a neighborhood $U$ of the glancing manifold ${\mathcal K}$ in $\Sigma:=\Sigma_{1}\cap\Sigma_{2}$.  
For any $\varrho\in U\setminus {\mathcal K}$, the point ${\mathcal J}_{j}(\varrho)\in U$ is just the second point of intersection of the characteristic of $\widetilde \omega\big|_{\Sigma_{j}}$ passing through $\varrho$ with  $U$. The set of fixed points of ${\mathcal J}_{j}$ is just ${\mathcal K}$. Moreover, \eqref{eq:glancing} implies that the differentials of ${\mathcal J}_{j}$ are linearly independent at any point of  ${\mathcal K}$. The billiard ball map is given by the composition ${\mathcal J}:={\mathcal J}_{2}\circ{\mathcal J}_{1}: U\to U$.  Moreover, ${\mathcal J}_{j}^\ast \omega_\Sigma=\omega_\Sigma$, where $\omega_\Sigma:=\widetilde\omega\big|_{\Sigma}$.  Then the billiard ball map ${\mathcal J}$ preserves $\omega\big|_\Sigma$ as well. Notice that the map ${\mathcal J}$ is smooth in $U$ but the two-form $\omega\big|_\Sigma$ is degenerate at 
${\mathcal K}$. To make the later symplectic one considers the quotient space $U/{\mathcal J}_{2}$ of $U$ by the action of ${\mathcal J}_{2}$. In our case it is given by the closed co-ball bundle $\overline {B}^\ast\Gamma \subset T^\ast\Gamma$ equipped with the canonical symplectic two-form. Let $\pi: U \to U/{\mathcal J}_{2} \cong\overline {B}^\ast\Gamma$ be the canonical projection. Then the  billiard ball map is represented by the boundary map $B= \pi\circ {\mathcal J}_{1} \circ \pi^+$, where $\pi^+$ is defined by \eqref{eq:outgoing-vector}. We call $B$ a billiard ball map as well. 

A local normal form of the pair of involutions ${\mathcal J}_{j}$, $j=1,2$, and of the two form $\omega$ has been obtained by Melrose  \cite{Mel} in a neighborhood of any point of  the glancing manifold 
${\mathcal K}$.   This normal form  leads to a local symplectic normal form of the billiard ball map $B$ at any point of the projection $\pi({\mathcal K})= S^\ast\Gamma$ (see also \cite{Hor}, Theorem 21.4.8). 

Consider a $C^1$ family of Riemannian metrics $[0,\delta_0] \ni t\to g_t$ in $\widetilde X$  and suppose that $\Gamma$ is locally strictly geodesically convex in $(\widetilde X, g_0)$ with respect to the outward normal field at $\Gamma$. Choosing $0<\delta\le \delta_0$ sufficiently small we obtain by \eqref{eq:glancing} that $\Gamma$ remains locally strictly geodesically convex in $(\widetilde X, g_t)$ for any $t\in [0,\delta]$. We denote by $\Sigma_{j,t}=\{f_{j,t}=0\}$, $j=1,2$, the  corresponding pairs of glancing hypersurfaces.   Here  $f_{1,t}+1$ is the Hamiltonian corresponding to the Riemannian metric $g_t$ via the Legendre transform and $f_{2,t}=f_2$, hence, both families of Hamiltonians are $C^1$ smooth with respect to $t$. Moreover, $f_{j,t}$ satisfy \eqref{eq:glancing} for $t\in [0,\delta]$ and we denote  by ${\mathcal K}_t$ the corresponding glancing manifolds.
Consider the  corresponding $C^1$ family of billiard ball maps $B_t:U_t \to {\bf B}_t^\ast\Gamma$, $t\in [0,\delta]$, where $U_t$ are suitable open subsets of ${\bf B}_t^\ast\Gamma$.  The   map $B_t$ is exact symplectic and smooth in  $U_t$, and it is extended by continuity as  the identity map on $S^\ast\Gamma$. 
Using the construction of the local symplectic  normal form of $B_t$ at $S^\ast\Gamma$ in \cite{Mel} and the interpolating Hamiltonian introduced  by Marvizi and Melrose \cite{MM} we obtain below a $C^1$ family of Birkhoff Normal Forms of $B_t$, $t\in [0,\delta]$. 

From now on we suppose that ${\rm dim}\, \widetilde X=2$  and we denote by $2\pi l_t$ the length of $\Gamma$ with respect to the Riemannian metric $g_t$.    Set $\A_t:= \T\times (l_t-\varepsilon,l_t+\varepsilon)$, where $\varepsilon>0$ will be chosen bellow small enough. Denote by $\rm{pr}:\R \to \T$ the canonical projection. 
\begin{Prop}\label{Prop:appr-inter-hamiltonian}
Let $(X,g_t)$, $t\in [0,\delta]$, be a $C^1$ family of connected locally strictly geodesically convex  billiard tables in $\widetilde X$. Then there exists 
\begin{itemize}
	\item[(1)] a $C^1$-family of exact symplectic transformation $[0,\delta]\ni t\to \widetilde\chi_t\in C^\infty(\A_t, V_t)$, where $V_t:=\widetilde\chi_t(\A_t)\subset T^\ast \Gamma$ is a neighborhood of $S^\ast_t\Gamma$,  $\widetilde \chi_t(\T\times \{l_t\})=S^\ast_t\Gamma$ and  
	$\widetilde \chi_t(\T\times (l_t-\varepsilon,l_t))\subset {\bf B}^\ast_t\Gamma$, 
   \item[(2)]  a $C^1$-family of real valued functions  $\zeta_t\in C^\infty (\R)$ and  $G_t\in C_0^\infty(\A_t)$, $t\in [0,\delta]$,  with $\zeta_t(l_t)=0$  and  $\zeta_t^\prime(l_t)<0$ 
\end{itemize}
such that the following holds
\begin{enumerate}
	\item[(i)] the function $\widetilde G_t\in C^\infty(\R\times \R)$ defined  by 
\[
\widetilde G_t(x,r):= x r - \frac{2}{3} \zeta_t(r)^\frac{3}{2}-G_t(\rm{pr}(x),r)
\]
is a generating function of the symplectic map $\widetilde P_t:=\widetilde \chi_t^{-1}\circ B_t\circ\widetilde \chi_t$ in  $\T\times (l_t-\varepsilon, l_t)$, 
	\item[(ii)] $G_{t}$ is flat at $r=l_t$, which means that $\partial_r^\alpha   G_{t}(\theta, l_t)\, =\, 0$
	for any $\theta\in\T$ and $\alpha\in \N$.  
\end{enumerate} 
\end{Prop} 
The function $ \widetilde\zeta_t:= \zeta_t\circ \chi_t^{-1}$ is an interpolating Hamiltonian of $B_t$ in the sense of Marvizi and Melrose \cite{MM}, which means that for any $\varphi\in C^\infty(T^\ast\Gamma)$, the function
\[
\varphi\circ B_t\, -\, \varphi\circ \exp\left(\widetilde \zeta^{\frac{1}{2}}X_{\widetilde \zeta}\right)
\]
is {\em flat} at $S^\ast_t\Gamma$, where $t\to \exp(tX_{\widetilde \zeta})$ is the flow of the Hamiltonian vector field $X_{\widetilde \zeta}$ of $\widetilde \zeta$ on $T^\ast\Gamma$. \\

\noindent
{\em Proof}. \quad The proof of the Proposition is based on a local  normal form of pairs of glancing hypersurfaces near the glancing manifold obtained by Melrose in \cite{Mel}. Consider the two involutions ${\mathcal J}_{j,t}$ associated to the characteristic foliations of $\Sigma_{j,t}=\{f_{j,t}=0\}$ in a neighborhood $U_t$ of the glancing manifold ${\mathcal K_t}$ in $\Sigma_t:=\Sigma_{1,t}\cap\Sigma_{2,t}$ for $t\in [0,\delta]$. In this way one obtains a $C^1$ family of
 billiard ball maps  given by the compositions ${\mathcal J}_{t}={\mathcal J}_{2,t}\circ {\mathcal J}_{1,t}: U_t\to U_t$ and 
 $B_t= \pi_t\circ {\mathcal J}_{1,t}\circ \pi_t^+$. Arguing as in the proof of \cite{Hor}, Theorem C.4.8, we first obtain a  $C^1$ family of normal forms of the two involutions ${\mathcal J}_{j,t}$. More precisely, following the first part of the proof of that theorem we get a $C^1$ family of diffeomorfisms $\widetilde\Psi_t: V \to U_t$, where $V$ is a neighborhood of a point 
$z^0=(z_1^0,0,z^{0\prime})\in \R^{2n-2}$ such that 
\[
\begin{array}{lcrr}
\widetilde\Psi_t^{-1}\circ {\mathcal J}_{1,t}\circ\widetilde\Psi_t(z_1,z_2,z') = (z_1+z_2,-z_2,z') + O_N(z_2^N)\, ,\\[0.3cm] 
\widetilde\Psi_t^{-1}\circ {\mathcal J}_{2,t}\circ \widetilde\Psi_t(z_1,z_2,z') = (z_1,-z_2,z') 
\end{array}
\]
for any $N\in\N$. In order to do this we consider an  asymptotic expansion of $\widetilde\Psi_t$ in formal power series  $\widetilde\Psi_t(z_1,z_2,z') \approx \sum \widetilde\Psi_{t,k}(z_1,z')z_2^k$. The functions $\widetilde\Psi_{t,k}(z_1,z')$ are obtained by solving linear systems of ordinary differential equations (see  the proof of \cite{Hor}, Theorem C.4.8). In this way we obtain  that the maps $t\to \widetilde\Psi_{t,k}\in C^\infty(V)$ are  $C^1$ smooth, and then using Borel's extension theorem we get a $C^1$ smooth family of maps $\psi_t:V\to U_t$ such that
\[
|\widetilde\Psi_t(z_1,z_2,z') - \sum_{k=0}^n \widetilde\Psi_{t,k}(z_1,z')z_2^k|\, \le \, C_N |z_2|^N.
\]
Then following the proof of Theorem 21.4.4 in \cite{Hor} (see also \cite{G-P}) one finds a $C^1$ family of diffeomorphisms $\Psi_t$ defined by  $\Psi_t^{-1}: W \to U_t$, where $W$ is a neighborhood of a point 
$\varrho^0=(x_1^0,x^{0\prime},0,\xi^{0\prime})\in T^\ast\R^{n-1}$, $U_t:=\Psi_t^{-1}(W)$ is an open neighborhood of the glancing manifold ${\mathcal K_t}$ in $\Sigma_{1,t}\cap\Sigma_{2,t}$  and such that 
\begin{equation}\label{eq:glancing-normal-form}
\left\{
\begin{array}{lcrr}
(\Psi_t^{-1})^\ast(\omega\big|_{\Sigma_t})=  dx_1\wedge d(\xi_1^2) + \sum_{k=2}^{n-1} dx_k\wedge d\xi_k  \\[0.3cm]
\Psi_t\circ {\mathcal J}_{1,t}\circ\Psi_t^{-1}(x_1,x',\xi_1,\xi') = (x_1+\xi_1,x',-\xi_1,\xi') + \widetilde R_t(x,\xi)\, ,\\[0.3cm] 
\Psi_t\circ {\mathcal J}_{2,t}\circ\Psi_t^{-1}(x_1,x',\xi_1,\xi') = (x_1,x',-\xi_1,\xi') \\[0.3cm] 
\forall\, j\in\N,\ \partial_{\xi_1}^j \widetilde R_t(x,0,\xi')=0.
\end{array}
\right.
\end{equation}
We mention just for information that  the formal power series are not convergent in general even when the hypersurfaces are analytic. 
It has been proved in \cite{G-P} that for any $t$ fixed  the corresponding functions in \eqref{eq:glancing-normal-form} belong to the Gevrey class $G^2$ of index two if the glancing hypersurfaces are analytic.

 Suppose now that $\Sigma_{1,t} = S_t^\ast \widetilde X$ and $\Sigma_{2,t}= T^\ast \widetilde X\big|_\Gamma$. Then $U_t$ is an open subset of $ S_t^\ast \widetilde X\big|_\Gamma$. Choosing normal to $\Gamma$ coordinates with respect to the metric $g_t$ one can assume that locally $f_{1,t}(y,y_n,\eta,\eta_n)=y_n$ and $f_{2,t}(y,y_n,\eta,\eta_n)= \eta_n^2 + q_t(y,y_n,\eta)$, where $t\to q_t(y,y_n,\eta)$ is a $C^1$ family of quadratic forms with respect to $\eta=(\eta_1,\ldots,\eta_{n-1})$ and $q_t(y,0,\eta)$ is the Hamiltonian corresponding to the induced metric on $\Gamma$ via the Legendre transform. In these coordinates $U_t$
 can be identified with  the set of $(y,y_n,\eta,\eta_n)$, where $y_n=0$, $(y,\eta)$ are  local coordinates in $T^\ast \Gamma$ near a  point $\varrho_t^0 \in S_t^\ast\Gamma$ and $\eta_n^2 + q_t(y,0,\eta,\eta_n)=1$,  while ${\mathcal K}_t\subset U_t$ is given by $y_n=\eta_n=0$. Moreover, ${\mathcal J}_{2,t}(y,0,\eta,-\eta_n)= - {\mathcal J}_{2,t}(y,0,\eta,\eta_n)$ and $\pi_t^{\pm}(y,\eta) = \pm \sqrt{1- q_t(y,0,\eta)}$ for $(y,\eta)\in {\bf B}_t^\ast\Gamma$. Setting 
\[
\Psi_{t}(y,\eta,\eta_n)=(x_t(y,\eta,\eta_n),\xi_{t1}(y,\eta,\eta_n),\xi'_t(y,\eta,\eta_n))
\]
where $\eta_n^2 + q_t(y,0,\eta,\eta_n)=1$
one obtains from the third relation of \eqref{eq:glancing-normal-form}  that $x_t$ and $\xi'_t$ are even functions of $\eta_n$ while $\xi_{t1}$ is odd. Then there exists a $C^1$ family of functions $t\to (\widetilde x_t,\widetilde\xi_t)\in C^\infty(T^\ast\R^{n-1})$ such that 
\[
\Psi_{t}(y,0,\eta,\eta_n)=(\widetilde x_t(y,\eta,\eta_n^2), \eta_n\widetilde\xi_{t1}(y,\eta,\eta_n^2),\widetilde\xi_t'(y,\eta,\eta_n^2)), 
\]
where $\eta_n^2+q_t(y,0,\eta)=1$.

We define  a $C^1$ family of diffeomorphisms $\widetilde \chi_t$ by
\[
\widetilde \chi_t^{-1}(y,\eta)=(\widetilde x_t(y,\eta, \widetilde \eta_n),\widetilde \eta_n\widetilde\xi_{t1}(y,\eta,\widetilde \eta_n)^2,\widetilde\xi_t'(y,\eta,\widetilde \eta_n)), 
\]
where $\widetilde \eta_n:=1-q_t(y,0,\eta)$. 
Then $\widetilde \chi_t: V\to V_t:=\widetilde \chi_t(V)\subset T^\ast\Gamma$ is a $C^1$ family of symplectic mappings, i.e. 
\[
\widetilde \chi_t^\ast(\sum_{j=1}^{n-1}dy_j\wedge d\eta_j)=\sum_{j=1}^{n-1}dx_j\wedge d\xi_j,  
\]
where $V\subset T^\ast\R^{n-1}$ is an open neighborhood of a given point $(x^0,0,\xi^{0\prime})$ and 
   we get the following symplectic  normal form of the billiard ball maps
\[
\widetilde\chi_t^{-1}\circ B_{t}\circ \widetilde\chi_t(x_1,x',\xi_1,\xi') = (x_1+\sqrt{\xi_1} ,x',\xi_1,\xi') + R_t(x,\xi),
\]  
where $t\to R_t\in C^\infty(T^\ast\R^{n-1}, T^\ast\R^{n-1})$ is a $C^1$ family of maps such that $\partial_{\xi_1}^j R_t(x,0,\xi')=0$ for any $j\in\N$. 
The interpolating Hamiltonian $\widetilde \zeta_t$  is defined by the $\xi_1$ component of  $\widetilde\chi_t^{-1}$, i.e. 
\[
\widetilde \zeta_t(y,\eta)= (1-q_t(y,0,\eta))\widetilde\xi_{t1}(y,\eta,1-q_t(y,0,\eta))^2.
\] 
As in \cite{MM} and \cite{P3} one obtains that $\widetilde \zeta_t$ is uniquely defined modulo a flat function on $S^\ast_t\Gamma$.

We suppose now that ${\rm dim\,} \Gamma =2$. 
To obtain the Hamiltonian $\zeta_t$   we find action-angle coordinates of $\widetilde \zeta_t$ as in \cite{P3}. To simplify the notations we drop the index $t$. Denote
by $M_{u}$ the closed curve $\{\varrho\in T^\ast\Gamma:\, \widetilde\zeta(\varrho) = u\}$ in $T^{*}\Gamma$
where $u$ varies in a small neighborhood of the origin. For any $\varrho \in
M_u$ consider
the map $\R \ni \ t\longrightarrow \ \exp(tX_{\widetilde\zeta})(\varrho ) \in\
M_{u}$ and
denote by $2\pi \Pi (u)$ its period. Let $S$ be
a section transversal to $M_0$ in  $T^\ast\Gamma$. It is equipped with local
coordinates $S\ni \varrho \rightarrow u=\widetilde\zeta (\varrho )$.
Denote by ${\mathcal O}$ the discrete group in $\R\times S$
generated by
\[
\R\times S\ \ni\ (t,u) \longrightarrow
(t + 2\pi \Pi (u)),u),\ u=\widetilde\zeta (\varrho ).
\]
Let $(\R\times S)/{\mathcal O}$ be the corresponding factor space. It is a symplectic manifold  equipped with the symplectic two-form 
$ dt \wedge du$ and the mapping
\[
\R\times S\ \ni\ (t,\varrho) \longrightarrow
\exp(tX_{\widetilde\zeta })(\varrho ) \in\ T^{*}\Gamma
\]
lifts to a symplectic diffeomorphism  from
$(\R\times S)/{\mathcal O}$ to a neighborhood of $M_{0}$.
Making suitable symplectic change of the variables
\[
\theta \ = t/\Pi (u),\ r = g(u),
\]
in $\R\times S$ we can suppose that ${\mathcal O}$ is generated by
$(\theta,r) \longrightarrow\ (\theta + 2\pi ,r)$ while the symplectic
two-form becomes $d\theta \wedge dr$. Then $ g'(u) = - \Pi (u)$ which yields
\begin{equation}
r(u) = l- \int_{0}^{u} \Pi (t)\,dt
 ,                   \label{eq:u-zeta}
\end{equation}
where $l={\rm length(\Gamma)}/2\pi$. 
Denote by $\zeta (r)$ the function inverse to
$r(u )$. 

We have obtained 
symplectic coordinates $(\theta _{t}(x,\xi), r_{t}(x,\xi))$, $t\in [0,\delta]$, in a
neighborhood of the boundary $S_t^\ast\Gamma$ in the co-ball bundle of $\Gamma$ with values in $\T\times \R $ such that
$S_t^\ast\Gamma = \{r_{t} = l_t\}$ and
${\bf B}^\ast_t  \subset \{r_{t} < l_t\}$. The map $t\to (\theta_t,r_t)\in C^\infty(T^\ast\Gamma))$ is $C^1$ by construction. 
The exact symplectic map $B_{t}$ is
generated in this coordinates by the function $\widetilde G_t$.
\finishproof

Recall that the functions $\beta_t(\omega)$, $I_t(\omega)$ and $L_t(I)$ are defined by \eqref{eq:beta function}, \eqref{eq:momentum-I} and \eqref{eq:birkhoff-L} respectively. 

\begin{Theorem}\label{Th:convex} Let $(X,g_t)$, $t\in [0,\delta]$,   be a $C^1$ family of  compact locally strictly geodesically convex  billiard tables of dimension two satisfying the weak isospectral condition $(\mbox{H}_1)-(\mbox{H}_2)$.  Then 
\begin{enumerate}
\item[(i)] There is a Cantor set $\Xi\subset (0,1]$ consisting of Diophantine numbers such that 
\[
\frac{{\rm meas}\, \left(\Xi\cap (0,\varepsilon)\right)}{\varepsilon}\,  =\,  1 - O(\varepsilon^2)\quad  \mbox{as}\quad \varepsilon\to 0^+ 
\]
and for any $\omega\in \Xi$ there exists a $C^1$ family of Kronecker invariant circles 
$[0,\delta]\ni t \to \Lambda_t(\omega)$ of $B_t$ of frequency $\omega$, 
\item[(ii)] $\forall \omega\in \Xi$ and $t\in [0,\delta]$, $\beta_t(\omega)= \beta_0(\omega)$, $I_t(\omega)=I_0(\omega)$ and $L_t(I_0(\omega))=L_0(I_0(\omega))$, 
\item[(iii)] $l_t=r_0$ and the function $\zeta_t-\zeta_0$ is flat at $r_0$ for any $t\in [0,\delta]$. 
\end{enumerate}
\end{Theorem}

We are going to prove Theorem \ref{Th:convex}.
Firstly, using  Theorem \ref{Theo:BNF} we will obtain 
 a suitable KAM theorem and a BNF at the corresponding family of invariant circles for the $C^1$ family of symplectic maps $\widetilde P_t$ given by Proposition \ref{Prop:appr-inter-hamiltonian}. 
To this end  we will determine the convex set $\Omega$, fix the parameters $\kappa$ and $\varrho$, and then estimate the corresponding quantities ${\mathcal B}_\ell$ and $\lambda$ which appear in Theorem \ref{Theo:BNF}. 

Consider the function $K_t\, := \, -\frac{2}{3} \zeta_t^\frac{3}{2}$ in $[l_t-\varepsilon,l_t]$. Fix $\varepsilon >0$ so that 
$\zeta_t^\prime(r) < 0$ for $(r,t)\in [l_t-\varepsilon, l_t]\times [0,\delta]$,
and denote by $K_t^\ast$ the Legendre transform  of $K_t$ in an interval $[0,a_0]$, $0<a_0\ll 1$. One can easily show that the family $t\to K_t^\ast$ can be extended as a $C^1$ family of smooth odd functions  $[0,\delta]\ni t\mapsto K_t^\ast\in C^\infty([-a_0,a_0])$. Indeed, the function $K_t$ admits an asymptotic expansion of the form
\[
K_t(r)\, \approx \, -\sum_{k=1}^\infty \,  \frac{1}{{\frac{1}{2}+k}}\, (l_t-r)^{\frac{1}{2}+k}\,  u_k(t) \quad \mbox{as} \quad r \nearrow l_t\, ,
\]
where $u_k\in C^1([0,\delta])$ and  $\displaystyle u_1(t) = (-\zeta_t^\prime(l_t))^\frac{3}{2}>0$. Moreover, this asymptotic expansion can differentiated infinitely many times with respect to $r$.
Recall that for any $t\in [0,\delta]$ fixed the derivative $\displaystyle K_t^{\ast \, \prime}$ of $K_t^{\ast}$  satisfies the identity $K_t^\prime(K_t^{\ast \, \prime}(\omega))=\omega$ for any $\omega\in (0,a_0]$. Moreover, $K_t(K_t^{\ast \, \prime}(\omega)) + K_t^{\ast}(\omega) = \omega K_t^{\ast \, \prime}(\omega)$, $K_t^{\ast \, \prime}(0) = l_t$,  and we easily obtain the asymptotic expansion 
\[
K_t^\ast(r)\, \approx \,  \sum_{k=0}^\infty \,  \frac{1}{2k+1} \, \omega^{2k+1}\,  v_k(t) \quad \mbox{as} \quad \omega\searrow 0\, ,
\]
where $v_k\in C^1([0,\delta])$,  
\[
v_0(t)=l_t\, ,\quad  v_1(t)= -u_1(t)^{-2}=-(-\zeta_t^\prime(l_t))^{-3}<0, 
\]
and so on. 

Fix $\tau>1$, set  $\Omega(a):= [a/2, 2a]$, choose $\kappa_a = a^2$ and denote by $\Omega_\kappa(a)$ the set of Diophantine frequencies $[a/2+ a^2, 2a-a^2]\cap D(\kappa,\tau)$. It follows from \cite{La}, Proposition 9.9,  that 
\begin{equation}\label{eq:measure-omega-convex}
\frac{{\rm meas}\, (\Omega(a)\setminus\Omega_{\kappa}(a))}{ {\rm meas}\, (\Omega(a)) }\,  \le \,  C\,  \kappa\, =\, C\, a^2 .
\end{equation}
Choose $a_0>0$ so that the Lebesgue measure of $\Omega_{\kappa}(a_0)$ is positive and denote by $\Omega_{\kappa}^0(a)$ the set of points of $\Omega_{\kappa}(a)$ of positive Lebesgue density. We have 
\begin{equation}\label{eq:D-convex}
\D(t,a):= K_t^{\ast \, \prime}(\Omega(a)) = [l_t + v_1(t)a^2/4 + O(a^4), l_t + 4v_1(t) a^2 + O(a^4)] \subset [l_t-\varepsilon, l_t]
\end{equation}
for $0<a \le a_0\ll 1$. 
Set $\A(t,a)=\T\times \D(t,a)$. 
We are ready to announce   the corresponding KAM theorem for the $C^1$ family of symplectic mappings $[0,\delta]\ni t\mapsto \widetilde P_t$ with generating functions $\widetilde G_t$ satisfying $(i)$ and $(ii)$ in Proposition \ref{Prop:appr-inter-hamiltonian}. 

\begin{Theorem}\label{theo:KAM-convex} 
For any $a\in (0,a_0]$ there exists  a  $C^1$-family  of exact symplectic maps 
\[
[0,\delta]\ni t \mapsto \big(\chi_t: \A(t,a) \to \A(t,a)\big)
\]
and of real valued functions  $[0,\delta]\ni t \mapsto L_t\in C^\infty(\D(t,a))$ and $[0,\delta]\ni t \mapsto R_t\in C^\infty(\A(t,a))$  such that for any $t\in [0,\delta]$ the following holds
\begin{enumerate}
\item $\widetilde G^0_t(x, I)= xI -L_t(I) - R_t({\rm pr}(x),I)$ is a generating function of $P_t^0:= \chi_t^{-1}\circ \widetilde P_t\circ \chi_t$ 
\item $L_t^\prime:  \D(t,a) \to \Omega(a)$ is a diffeomorphism with inverse $L_t^{\ast\, \prime}:  \Omega(a)\to \D(t,a) $, where $L_t^{\ast}$ is the Legendre transform of $L_t$
\item $R_t$ is flat at $\T\times L_t^{\ast \, \prime}( \Omega_{\kappa}^0(a))$ 
\item  for any integer $N\ge 1$ and $m\in\N$  there exists a constant  $C= C_{m,N}>0$  independent of  $a\in(0,a_0]$ and $t\in [0,\delta]$ such that the following estimates hold
\[
|\partial_\varphi^\alpha(\kappa\partial_\omega)^\beta \sigma_\kappa^{-1}(\chi_s - {\rm id})| + 
|\partial_\varphi^\alpha(\kappa\partial_\omega)^\beta \sigma_\kappa^{-1}(\chi_s^{-1} - {\rm id})|\,  \le\,  C   \kappa^{2N-m-\frac{3}{4}}
\]
on $\A(t,a)$, and 
\begin{equation}\label{eq:estimates-convex}
\left|\left(\kappa \frac{d}{d I}\right)^m  (L_t'(I) - K_t'(I))\right|\,  \le \, C  \kappa^{2N-m+\frac{1}{4}}
\end{equation}
on $\D(t,a)$ for any $t\in [0,\delta]$ and $m\in \N$. 
\end{enumerate}
\end{Theorem}

\noindent
{\em Proof.} We are going to apply Theorem  \ref{Theo:BNF} to the $C^1$ family of symplectic mappings $[0,\delta]\ni t\mapsto \widetilde P_t\in  C^\infty(\A(t,a), \A(t,a))$ given by Proposition \ref{Prop:appr-inter-hamiltonian}.  

Let us  estimate the corresponding quantities  ${\mathcal B}_m$ for $m\in\N$ and $\lambda$ defined by \eqref{eq:B-0-sequence} - \eqref{eq:B-sequence} and \eqref{eq:non-degeneracy-BNF}. The constant $\lambda$    can be fixed by
\begin{equation}\label{eq:lambda}
\lambda = \lambda_a= \sup_{t\in [0,\delta]} \| K''_t\|_{\D(t,a);\kappa} = C_0 a^{-1}= C_0 \kappa^{-1/2},
\end{equation}
where $C_0$ is a positive constant independent of $a$.  Given $\ell=m+\mu$ with $m\in\N_\ast$  we get by \eqref{eq:inclusion} that for any $t\in [0,\delta)$ the following inequality holds 
\[
\|G_t\|_{\ell, \A(t,a);\kappa}  \le \|G_t\|_{m+1, \T\times \D(t,a);\kappa}.
\]
since  $\D(t,a)$ is an interval. Moreover, 
\[
\|G_t\|_{m+1, \T\times \D_a;\kappa} = 
\sup_{|\alpha|+ |\beta|\le m+1}\, \|\partial_\theta^\alpha (\kappa\partial_r)^\beta G_t\|_{C^0(\T\times \D(t,a))}.
\]
Fix $N\ge 2$.  It follows from   Proposition \ref{Prop:appr-inter-hamiltonian}, $(ii)$, and the definition of $\D(t,a)$ that 
\[
\|G_t\|_{\ell, \A(t,a);\kappa}\,  \le\, C_{m,N} \, a^{8N+4} \,  = \, C_{m,N} \, \kappa^{4N+2}
\]
for any $t\in [0,\delta]$, where $C_{m,N}$ is a positive constant. 
Choosing  $\varrho=  \kappa^{2N +5/4}<\kappa$ we obtain
\[
\|G_t\|_{\ell, \A(t,a);\kappa}\,  \le\,  C_m'\, \kappa\varrho  \, \kappa^{2N-1/4}
\]
for any $t\in [0,\delta]$. Moreover, 
\begin{equation}\label{eq:norm-K-second}
|\!|\!| K_t'' |\!|\!|_{\ell, \D(t,a);\kappa} \, \le\,  \| K_t'' \|_{m+1, \D(t,a);\kappa}\,  \le \, C_m a^{-1} \, = C_m\,  \kappa^{-1/2}
\end{equation}
and 
\[
S_{\ell}(\nabla K^{\ast})\, \le \, \,\sup_{0\le t\le \delta}\, \big(1+ \|\nabla K^{\ast}_t\|_{C^1([-a_0,a_0])}\big)^{\ell-1} \big(1+ \|\nabla K^{\ast}_t\|_{C^\ell([-a_0,a_0])}\big)  \, \le \, C_m
\]
where $C_m$ is a positive constant. 
Thus for any $m\in \N$ we obtain from \eqref{eq:B-0-sequence} - \eqref{eq:B-sequence} that 
\begin{equation}\label{eq:estimate-B}
{\mathcal B}_m\, \le \, C_m' {\mathcal B}_m^0\, \le \,  C_\ell \, \varrho \,  \kappa^{2N +3/4} \, = \,  C_\ell \, \kappa\varrho \,  \kappa^{2N -1/4}
\end{equation}
where $C_m,\, C_m'>0$ depends neither on $a\in (0,a_0]$ nor on  $t\in [0,\delta]$. 
Choosing $a_0\ll 1$ we get ${\mathcal B}_2\le \epsilon \kappa\varrho \lambda^{-4}$ for any $a\in (0,a_0]$ since $N\ge 2$ and $\lambda= C_0 \kappa^{-1/2}$,  which gives \eqref{eq:smallness-condition-1}. Applying Theorem  \ref{Theo:BNF} we obtain {\em 1}-{\em 4}. In particular, taking into account \eqref{eq:lambda} - \eqref{eq:estimate-B} we obtain from \eqref{eq:estimates} the estimate
\[
\sup_{t\in [0,\delta]}\, \sup_{I\in \D(t,a)}\, \left|\left(\kappa \frac{d}{d I}\right)^m  (L_t'(I) - K_t'(I))\right|\,  \le \, C_m'  \kappa^{2N+\frac{3}{4}}\lambda^{2m}(\lambda + \kappa^{-\frac{1}{2}})
\,  \le \, C_m  \kappa^{2N-m+\frac{1}{4}}
\]
for any $m\in \N$, where $C_m'$ and $C_m$ are positive constants. 
\finishproof
\\

\noindent 
{\em Proof of Theorem \ref{Th:convex}.} \quad It follows from Theorem \ref{theo:KAM-convex}  that for any $0<a\le a_0\ll 1$ and 
 $\omega\in \Omega_{\kappa}^0(a)$ with $\kappa= a^2$  there exists a $C^1$ family of Kronecker invariant tori $[0,\delta]\ni t\to \Lambda_t(\omega)$ of $B_t$. 
Corollary \ref{Th:main2} implies that $I_t(\omega)=I_0(\omega)$ and $L_t'(I_0(\omega))= L_0'(I_0(\omega))$. Notice that
\[
\displaystyle\limsup_{a\to 0}\{I_t(\omega)-l_t:\, \omega\in [a/2, 2a]\}=0.
\]
Then 
\[
|l_t-l_0|\le \displaystyle\limsup_{a\to 0}\{I_t(\omega)-l_t:\, \omega\in [a/2, 2a]\} + \displaystyle \limsup_{a\to 0}\{I_0(\omega)-r_0:\, \omega\in [a/2, 2a]\} = 0
\]
hence, $l_t=l_0$ for any $t\in [0,\delta]$.
Moreover, the function $\omega\mapsto L_t'(I_0(\omega)) - L_0'(I_0(\omega))$ is flat at the set $\Omega_{\kappa}^0(a)$ in view of Lemma \ref{Lemma:flat}. 
Then for any $\omega\in \Omega_{\kappa}^0(a)$, $t\in [0,\delta]$, and  any  $m\in N$  
using  the equality $I_t(\omega)=I_0(\omega)$ and the estimate \eqref{eq:estimates-convex}  with $N=m$ we obtain 
\[
\begin{array}{rcll}
\left|(d/d I)^m \left(K_t' -  K_0'\right)(I_0(\omega))\right| 
&\le&  \left|(d/d I)^m \left(K_t' -  L_t'\right)(I_t(\omega))\right|\\[0.3cm]
&+& \left|(d/d I)^m \left(K_0' - L_0'\right)(I_0(\omega))\right|\\[0.3cm]
&\le&  C \kappa^{1/4} = C a^{1/2}.
\end{array}
\]
Taking the limit as $a\searrow 0$ we obtain that the function $K_t -K_0$ is smooth in $[r_0,r_0+\varepsilon]$ and flat at $r_0$. Then 
$\zeta_t -\zeta_0= (3 K_t/2)^{2/3}-(3 K_0/2)^{2/3}$ is also flat at $r_0$. The set of frequencies $\Xi$ is defined as the union of $\Omega_{\kappa}^0(a)$.   \finishproof\\

\noindent
{\em Proof of Theorem \ref{Th:convex-main}.} \quad It remains to show that $\Gamma_t$ is strictly convex for any $t\in
[0,1]$. To do this we  are going to use an argument from \cite{P3}. To simplify the notations we will omit the index $t\in [0,1]$. 

Consider the  interpolating Hamiltonian $\widetilde \zeta(x,\xi ) =
\zeta(\widetilde\chi ^{-1}(x,\xi ))$ of $B$, where the function $\zeta$ and the symplectic transformation $\widetilde\chi$ are given by Proposition \ref{Prop:appr-inter-hamiltonian}. For any $r$ with $|r|$ 
small enough  the level set
$$
M(r)\ =\ \{(x,\xi)\in T^\ast\Gamma :\, \widetilde\zeta(x,\xi ) = r \}
$$
is an ``circle" and we set
\begin{equation}
\nu(r)\ =\ \int_{M(r)} d\, \big(z \big|_{M(r)}\big)
                                                     \label{eq:integral-inv}
\end{equation}
where $r\to d \big(z\big|_{M(r)} \big)$ is a smooth family of $1$-forms on $M(r)$ such that $d \big(z\big|_{M(r)} \big)(X_{\widetilde\zeta})=1$. One can consider $z\big|_{M(r)}$ as a multivalued function on the circle $M(r)$ which is well defined on the corresponding covering space $\R\to M(r)$ so that
\begin{equation}
\{\widetilde\zeta,z\}\ =\ dz(X_{\widetilde\zeta}) =\ 1.
                                                      \label{eq:poisson}
\end{equation}
It is easy to show that the set of Taylor coefficients of $\nu(r)$
at $r=0$ is algebraically equivalent to the set of Taylor's coefficients
of $\zeta(I)$ at $I=l$. Indeed, performing the symplectic change of the
variables $(x,\xi ) = \widetilde\chi (\varphi ,I),\ (\varphi ,I)\in \A
$, and using (\ref{eq:poisson}) we easily get
\[
\zeta'(I) \nu(\zeta(I))\ =\ 2\pi .
\]
Denote by ${\cal R}_r$ the function inverse to $I \rightarrow \zeta(I)$. Then  
(\ref{eq:integral-inv}) implies
\begin{equation}
\nu(r)\ =\ 2\pi {\cal R}'(r)
\label{eq:integral-inv1}
\end{equation}
and we obtain  that the Taylor coefficients of $\nu(r)$ at $r=0$ determine those of $\zeta$ at $I=\ell$ and vice versa. 

 The Taylor coefficients of $\nu(r)$ at $r=0$, also called
integral invariants, have been investigated by Sh.~Marvizi and
R.~Melrose
\cite{MM}. They are given by integrals on $\Gamma$ of certain
polynomials of the curvature $\kappa (x)$ of $\Gamma$ and its
derivatives. In particular,  (4.6) in \cite{MM} and (\ref{eq:integral-inv1})
yield together
\begin{eqnarray}
\label{eq:inv1}
{\cal R}'(0)\ & = & -\frac{1}{\pi} \int_{0}^{\ell}
\kappa (x)^{2/3} dx\\   
\label{eq:inv2}                 
{\cal R}''(0)\ & = & \frac{1}{2160\pi} \int_{0}^{\ell}
(9\kappa (x)^{4/3} \ +\ 8\kappa (x)^{-8/3}
\kappa '(x)^{2})dx
\end{eqnarray}
(see also \cite{Sor}).

Suppose now  that $X_t$ is strictly convex for $0\le t<\delta$ but only
convex for $t=\delta$. Consider the function ${\cal R}_t(r)$ inverse
to $r = \zeta_t(I)$. Then
Theorem \ref{Th:convex}, $(iii)$, yields 
$$
{\cal R}_t(r) = {\cal R}_0(r) + O_N(r^N) \quad \mbox{as} \ r\to 0
$$
for any $N\in\N$ and  we obtain
\begin{equation}
{\cal R}'_t(0)\ =\ {\cal R}'_0(0),\
{\cal R}''_t(0)\ =\ {\cal R}''_0(0),\ s \in [0,b_0).
                                        \label{eq:inv3}
\end{equation}
Denote by $\kappa _t(x)>0,\ x \in \Gamma _t$ the curvature of $
\Gamma _t$ and define $f_\delta$ by $f_t(x) = \kappa _t (x)^{-1/3}$ for
$t<\delta$ and 
$f_\delta(x) = \kappa _\delta (x)^{-1/3}$ if $\kappa_\delta(x) \neq 0$ and $f_\delta(x) = 0$ if $\kappa _\delta(x) = 0$. The second
equality of (\ref{eq:inv3}) and (\ref{eq:inv2}) yield together
\begin{equation}
\int_{\Gamma_t} \mid f'_t(x) \mid ^2 dx \ \leq \ C,\
s \in [0,\delta).
                                          \label{eq:inv4}
\end{equation}
where $C$ is a positive constant.
On the other hand, the first equality of (\ref{eq:inv3}) and
(\ref{eq:inv1}) imply that for any $t
\in [0,\delta)$ there exists $x_t \in \Gamma_t$ such that
$$
\kappa _t(x_t)\,  \ge\,  C_1:=\left( -\frac{\pi}{l_0}{\cal R}'_0(0)\right)^{3/2} >0.
$$
Then $f_t(x_t) \leq C_1^{-1/3}$ for $t\in [0,\delta)$, and using Taylor's formula and (\ref{eq:inv4}) we obtain
the estimate
$$
\int_{\Gamma_t} (\mid f_t(x) \mid ^2  \ +\
\mid f'_t(x) \mid ^2) dx \ \leq\ C_2,\
s \in [0,\delta),
$$
where $C_2$ is a positive constant. Let $[0,\delta]\ni t\to \psi_t: \Gamma\to \R^2$ be a $C^1$ family of embeddings such that $\Gamma_0=\Gamma$ and $\psi_t(\Gamma)=\Gamma_t$. 
Then $\{f_t\circ \psi_t:\ t\in [0,\delta)\}$ is a compact subset of
$L^2(\Gamma)$ and we obtain that $f_{\delta} \circ \psi_{\delta}\in L^2(\Gamma)$ as well. On the other hand, $\Gamma_{\delta}$ is convex but not strictly convex, hence the
curvature its curvature $k_{\delta}$ is a non-negative function and it has a zero of at least second
order at a point $x_0 \in \Gamma$. Then
$$
\mid f_{\delta}(x) \mid \ \geq\ C\mid x-x_0\mid ^{-2/3}
$$
in any local coordinates in a neighborhood  of $x_0$ in $\Gamma_{\delta}$. Hence
 $f_{\delta} \notin L^2(\Gamma)$ which leads to a contradiction. This implies that 
$\Gamma_t$ is strictly convex for any $t\in [0,1]$. \finishproof

\section{Microlocal Birkhoff Normal Form of the monodromy operator} \label{Sec:QBNF}

Starting from the BNF in Theorem \ref{Theo:soft-BNF} we are going to find a microlocal (quantum) Birkhoff normal form (shortly QBNF) at the union of the invariant tori $\Lambda_t(\omega)$, $\omega\in\Omega_\kappa^0$, of the  corresponding microlocal monodromy operator for the  family of Laplace-Beltrami operators $\Delta_t$ in $X$ with Dirichlet boundary conditions. A similar  QBNF has been obtained in \cite{PT4} for perturbations of the function in the Robin boundary conditions around a single Kronecker torus. In contrast to \cite{PT4}  the BNF of the tori here is nondegenerate which simplifies the construction. 
 
Let us present  the main steps in the construction. 
At first we reduce the problem to the boundary and   introduce the corresponding microlocal monodromy operator $M_t^0(\lambda)(\lambda)$, $t\in J$. The reduction to the boundary is obtained by a variant of  the reflection method for the wave equation which consists in the following. Given a suitable function $f(\cdot;\lambda)$ on $\Gamma$ depending on a large parameter $\lambda$ the frequency support of which is contained in a small neighborhood of the union of the invariant tori $\Lambda_t(\omega)$, we consider the corresponding outgoing solution of the reduced wave equation (the Helmholtz equation) in $X$ and we reflect it at the boundary $m-1$ times if $m\ge 2$. After each reflection at the boundary we consider the corresponding branch of the solution $u_t$ of the Helmholtz equation given by the outgoing parametrix. We denote by  $M_t(\lambda)f$ the restriction at $\Gamma$ of the last  branch of the  solution $u_t$. We call $M_t(\lambda)$  a monodromy operator.  By construction, the function $f(\cdot,\lambda)$ on $\Gamma$ gives rise to an asymptotic solution $u_t(\cdot,\lambda)$  of the Dirichlet problem of the Helmholtz equation 
\[
(-\Delta_t +\lambda^2)u_t= O_N(|\lambda|^{-N})f_t, \quad u_t|_\Gamma = O_N(|\lambda|^{-N})f_t,
\]
 or a quasi-mode $(\lambda, u_t)$ of the Laplace-Beltrami operator with Dirichlet boundary conditions when $\|u_t\|_{L^2} = 1$
 if and only if 
\[
M_t^0(\lambda)(\lambda)f=f + O_N(|\lambda|^{-N})f.
\] 
The family of operators   $M_t^0(\lambda)(\lambda)$, $t\in J$, is a $C^1$ family of Fourier Integral Operator with a large parameter $\lambda$ ($\lambda$-FIO) the canonical relation of each of them being  the graph of $P_t$. For this reason we recall in Sect. \ref{subsec:FIO} some properties of the $\lambda$-FIOs associated with a $C^1$ family of Lagrange immersions. The reduction to the boundary and the construction of the microlocal monodromy operator is done in Sect. \ref{subsec:reduction}. 

 Our next goal is to ''separate the variables''
microlocally near the whole family of invariant tori $\Lambda_t(\omega)$, $\omega\in\Omega_\kappa^0$. This is done in Sect. \ref{subsec:QBNF}. To this end we use the  Birkhoff normal form of $P_t$ given by Theorem \ref{Theo:soft-BNF}.  First we conjugate   $M_t^0(\lambda)$ with a microlocally unitary  $\lambda$-FIO $T_t(\lambda)$ the canonical relation of which is the graph of the symplectic transformation
$\chi_t$ given by Theorem \ref{Theo:soft-BNF}. In this way we obtain 
 a $\lambda$-FIO  $W_t(\lambda)$ the canonical relation of which is just the graph of $P_t^0$ (see Proposition \ref{operator-T}).  Then we obtain a microlocal  Birkhoff normal form  $W_t^0(\lambda)$ of $W_t(\lambda)$ by conjugating it  with a suitable  $\lambda$-PDO  and  solving at any step the corresponding homological equation.  
In this way we separate microlocally the variables near the whole family of invariant tori. This means that the amplitude of $W_t^0(\lambda)$ does not depend on the angular variables but only on the action variables at the family of invariant tori, which allows us to obtain  a microlocal ''spectral decomposition''   of $W_t(\lambda)$ near the family $\Lambda_t(\omega)$, $\omega\in\Omega_\kappa^0$. At any step  the corresponding phase functions and amplitudes are $C^1$ with respect to the parameter $t$. 

\subsection{$C^1$ families of PDOs and FIOs with a large parameter $\lambda$.}\label{subsec:FIO}
\subsubsection{$C^1$ families of symbols and $\lambda$-PDOs .}\label{subsec:PDOs}
Let   $M^d$  be a smooth paracompact  manifold of  dimension $d$. We are going to define  a class of $C^1$ families of pseudo-differential operators depending on a large parameter $\lambda$ (shortly $\lambda$-PDOs) acting on the half-density bundle 
$\Omega^{\frac{1}{2}}(M^d)$ of $M^d$. The large  parameter $\lambda$ will belong to  the  set 
\begin{equation}\label{eq:D-lambda}
{\mathcal D}\, := \{\lambda\in \C:\, |{\rm Re}\, \lambda| \ge C_0,\, |{\rm Im}\, \lambda| \le C_1\} ,\quad \sup_{\lambda\in {\mathcal D}}\, |\lambda|\, =\, +\infty\, ,
\end{equation}
where  $C_0,\ C_1>0$.  One can switch to   the semi-classical setting by introducing $\hbar := 1/\lambda$. 

Let us first define the symbols we are going to deal with. Given an interval $J\subset \R$ we define a  $C^1$ family of symbols $J\ni\to a_t$  of order $0$ in $T^\ast \R^d$ as a map
\[
J\times {\mathcal D}\, \longrightarrow \, C_0^\infty(T^\ast \R^d)\, ,\quad (t,\lambda)\, \longmapsto \, a_t(\cdot,\lambda),
\]
such that 
\begin{itemize}
	\item[--] The map $J\ni t \to a_t(\cdot,\lambda)\in C^\infty(T^\ast \R^d)$ is $C^1$ for any $\lambda\in {\mathcal D}$ fixed;
	\item[--] The support ${\rm supp\, } a_t(\cdot,\lambda)$ is contained in a fixed compact subset of $T^\ast \R^d$ independent of $(t,\lambda)\in I\times {\mathcal D}$;
	\item[--] For any $\alpha,\, \beta\in \N^d$ there exists a positive constant $C_{\alpha,\beta}$ such that
\[
|\partial_t^k\partial_x^\alpha\partial_\xi^\beta a_t(x,\xi,\lambda)|\, \le\, C_{\alpha,\beta}
\]
for every $(t,\lambda)\in I\times {\mathcal D}$, $(x,\xi)\in T^\ast\R^d$ and $k\in \{0, 1\}$.
\end{itemize}
 In this case we say that $a_t$ is a $C^1$ family of symbols in $S^0(T^\ast\R^d\times {\mathcal D})$ with respect to the parameter  $t\in J$. We set $S^p(T^\ast\R^d\times {\mathcal D})= \lambda^p S^0(T^\ast\R^d\times {\mathcal D})$ for $p\in\R$ and denote by $S^{-\infty}(T^\ast\R^d\times {\mathcal D})$ the residual set $\cap_{p\ge 0}S^{-p}(T^\ast\R^d\times {\mathcal D})$. 
We say that 
\begin{equation}\label{eq:formal-symbols}
J\ni \, t\, \longmapsto\,  \sum_{j\in\N}\,  a_{t ,j} \lambda^{-j}
\end{equation}
is a $C^1$ family of formal symbols of order $0$ if for any $j\in \N$ the map $J\ni t\mapsto  a_{t, j}\in C^\infty(T^\ast\R^d)$ is $C^1$ smooth and the support   ${\rm supp\, } a_{t, j}$ is contained in a fixed compact subset of $T^\ast \R^d$ independent of $(t,j)\in I\times \N$. A $C^1$ family of symbols $J\ni t\mapsto a_t\in S^0(T^\ast\R^d\times {\mathcal D})$ is said to be a  realizations of the $C^1$ family of formal symbols \eqref{eq:formal-symbols} if for any $N\in \N$ and  $\alpha,\, \beta\in \N^d$ there exists a positive constant $C_{N,\alpha,\beta}$ such that
\begin{equation}\label{eq:amplitude}
\sup_{(t,x,\xi,\lambda)\in J\times T^\ast\R^d\times {\mathcal D} }\, 
\Big|\partial_t^k\partial_x^\alpha \partial_\xi^\beta\Big(a_t(x, \xi,\lambda) - \sum_{j=0}^{N-1} a_{t,j}(x, \xi)\lambda^{-j}\Big)\Big| \, \le \, 
C_{N,\alpha,\beta} |\lambda|^{-N}
\end{equation}
for  $k\in \{0, 1\}$. Symbols admitting an asymptotic expansion of the form \eqref{eq:amplitude} for any $N$ are said to be classical. We denote by $S^0_{\rm cl}(T^\ast\R^d\times {\mathcal D})$ the class of the classical symbols.  Any $C^1$ family of formal symbols of order zero  admits a $C^1$ family of realizations by Borel's theorem. 
\begin{Prop}\label{prop:realization}
Any $C^1$ family of formal symbols \eqref{eq:formal-symbols} of order $0$ admits a  realization as a $C^1$ family of symbols $J\ni t\mapsto a_t\in S^0(T^\ast\R^d\times {\mathcal D})$. Moreover, if  $a_t$ and $a_t^\prime$  are two  $C^1$ family of realizations of \eqref{eq:formal-symbols} then $J\ni t\mapsto a_t-a_t^\prime \in S^{-p}(T^\ast\R^d\times {\mathcal D})$ is  a $C^1$ family for every $p\ge 0$. 
\end{Prop}
We give a prove of Borel's theorem in Appendix \ref{Sec:Whitney}.

We say that the family of operators $J\ni t\mapsto {\rm Op\, }(a_t)$ with Schwartz kernels
\begin{equation}\label{eq:kernel-of-pdo1}
K_{{\rm Op\, }(a_t)}(x,y,\lambda):=  \left(\frac{\lambda}{2\pi}\right)^d \, \left(\int_{\R^{d}}\, e^{i\lambda\langle x-y,\xi\rangle }a_t(x,\xi,\lambda)\, d\xi\right)\, |dx|^\frac{1}{2}|dy|^\frac{1}{2}
\end{equation}
is a $C^1$ family of  $\lambda$-PDOs of order zero acting on $\frac{1}{2}$-densities if $J\ni t\mapsto a_t\in S^0(T^\ast\R^d\times {\mathcal D})$ is a $C^1$ family of symbols. We say that a   family of operators $J\ni t\mapsto A_t$ acting on the smooth sections of the  half-density bundle $\Omega^\frac{1}{2}(M^d)$ of the manifold $M^d$ is a $C^1$ family of  $\lambda$-PDOs if it is given  by a  $C^1$ family of  $\lambda$-PDOs with Schwartz kernels of the form \eqref{eq:kernel-of-pdo1} in any local coordinates. 

\subsubsection{$C^1$ families of $\lambda$-FIOs.}\label{subsec:FIOs}
Consider a $C^1$ family of  exact Lagrange immersions 
\begin{equation}\label{eq:lagrangian-immersions}
\imath_t :\Lambda\to T^\ast M^d,\quad t\in [0,\delta],  
\end{equation}
which means that  the map $[0,\delta]\ni t\mapsto \imath_t\in C^\infty(\Lambda,T^\ast M^d)$ is $C^1$, $\imath_t$ is an immersion  and  the pull-back $\imath_t^\ast(\xi dx)$ of the canonical one-form $\xi dx$ of $T^\ast M^d$ is exact for each  $t\in [0,\delta]$.  Then there exists a $C^1$ mapping $[0,\delta]\ni t\mapsto f_t\in C^\infty(\Lambda)$ such that
\begin{equation}\label{eq:f-t}
\imath_t^\ast(\xi dx) = df_t \quad t\in [0,\delta]. 
\end{equation}
Fix $t\in [0,\delta]$. 
Recall from   \cite{Du} and  \cite{Hor}  that a real valued phase function  $\Phi_t(x,\theta)$ defined in a neighborhood  of a point $(x^0,\theta^0)\in \R^d\times \R^N$ with  $d_\theta\Phi_t(x^0,\theta^0) = 0$ is  nondegenerate at $(x^0,\theta^0)$ if 
\begin{equation}\label{eq:rank-of-phi}
{\rm rank}\, d_{(x,\theta)}d_\theta \Phi_t (x^0,\theta^0)= N.
\end{equation}
Then  there exists a  neighborhood $V\subset \R^d\times \R^N$  of $(x^0,\theta^0)$  such that \eqref{eq:rank-of-phi} holds for any $(x,\theta)\in V$ and 
\[
C_{\Phi_t} := \{(x,\theta)\in V: \,  d_\theta \Phi_t = 0\}
\] 
is a smooth manifold of dimension $d$.  Moreover, the differential of the map
\begin{equation}
\imath_{\Phi_t} : C_{\Phi_t} \ni (x,\theta) \longrightarrow (x,d_x\Phi(x,\theta))
\in \Lambda_{\Phi_t}:=\imath_{\Phi_t} (C_{\Phi_t})
                                             \label{lagrangian1}
\end{equation}
is of rank $d$ and shrinking  $V$ if necessary we obtain that $\Lambda_{\Phi_t}$ is an embedded  Lagrangian submanifold of $T^\ast M^d$.   We say that  the nondegenerate function $\Phi_t(x,\theta)$, $(x,\theta)\in V$,  defines locally the Lagrange immersion $\imath_t:\Lambda\to T^\ast M^d$ if there is an open subset  $W_t \subset  \Lambda $   such that  
\begin{equation}\label{eq:t-immerson} 
  \imath_t : W_t \to \Lambda_{\Phi_t} \quad \mbox{is a diffeomorphism}. 
\end{equation}
We can take $\Phi_t=\Phi_t(x)$ depending only on the coordinates $h$ ($N=0$) if the corresponding Lagrangian manifold is ``horizontal'' which means that the projection to the base is a local diffeomorphism. 
 The collection $(\imath_{t}^{-1}(\Lambda_{\Phi_t}), \imath_{\Phi_t}^{-1}\circ \imath_{t})$ provides the Lagrangian immersion $\imath_t:\Lambda \to T^\ast M^d$ with an atlas of local charts. Given an interval $J\subset [0,\delta]$ we say that a $C^1$  map $J\ni t\to \Phi_t\in C^\infty(V,\R)$ is a $C^1$ family of nondegenerate phase functions generating locally the $C^1$ family of Lagrange immersions \eqref{eq:lagrangian-immersions} in $J$ if for any $t\in J$ the phase function $\Phi_t$ is nondegenerate in $V$ and \eqref{eq:t-immerson} holds. Such $C^1$ families of phase functions can always be constructed locally. 
Consider the function $\Phi_t^\Lambda: = \big({\it k}_t\circ\imath_{\Phi_t}^{-1}\circ \imath_{t}\big)^\ast \big(\Phi_t\big)$ on $\imath_{t}^{-1}(\Lambda_{\Phi_t})$ where ${\it k}_t: C_{\Phi_t} \to  \R^d\times\R^N$ is the inclusion map.
Observe that
$$
d \Phi_t ^\Lambda
=\big(\imath_{\Phi_t}^{-1}\circ \imath_{t}\big)^\ast \Big( {\it k}_t^\ast\Big(\frac{\partial\Phi_t}{\partial x}dx  +
 \frac{\partial\Phi_t}{\partial  \theta}d \theta\Big)\Big)
$$
$$
= \imath_{t}^\ast \Big(\big({\it k}_t\circ \imath_{\Phi_t}^{-1}\big)^\ast\Big(\frac{\partial\Phi_t}{\partial x}dx \Big)\Big) =  
\imath_{t}^\ast\big(\xi dx\big) = df_t  
$$ 
and we choose $\Phi_t$ so that $\Phi_t^\Lambda=f_t$ on $\imath_{t}^{-1}(\Lambda_{\Phi_t})$, where $f_t$ is defined in \eqref{eq:f-t}. 

Given a $C^1$ family of classical amplitudes  $[0,\delta]\ni t\to a_t\in S^0_{\rm cl\, }(V\times {\mathcal D})$ such that $a_t=0$ for $t\notin J$   we consider the $C^1$ family of oscillatory $\frac{1}{2}$-densities 
\begin{equation}\label{oscillatory1}
I_{\Phi_t,a_t}(x,\lambda)\, |dx|^\frac{1}{2}\,  =\,  \Big(\frac{\lambda}{2\pi}\Big)^{m+\frac{d+2N}{4}}\, 
\Big(\int_{{\R}^N}\ e^{i\lambda\Phi_t(x,\theta)} a_t(x,\theta,\lambda)\, d\theta \Big)\, |dx|^\frac{1}{2}\,  
\end{equation}
with the convention that there is no integration when $N=0$.  
Notice that the function $[0,\delta]\ni t\mapsto I_{\Phi_t,a_t}\in C^\infty(M^d)$ is $C^1$  for  each $\lambda\in {\mathcal D}$ fixed. Its oscillation  is detected as $\lambda \to \infty$ by  the corresponding semi-classical wave front set. Integrating by parts one obtains ${\rm WF}_\lambda (I_{\Phi_t,a_t}(\cdot,\lambda)) \subset \Lambda_{\Phi_t}$, where ${\rm WF}_\lambda $ is the frequency set  (or semi-classical $\hbar$-wave-front with $\hbar=1/\lambda$) (cf. \cite{Alex}, \cite{D-S}, \cite{Zw}, \cite{GS}). A (global) $C^1$ family of  oscillatory $\frac{1}{2}$-densities is given by 
\begin{equation}\label{eq:global-FIO}
u_t(x,\lambda) = \sum_{j}  I_{\Phi_{t}^j,a_{t}^j}(x,\lambda)\, |dx|^\frac{1}{2}
\end{equation}
where  $\Phi_{t}^j$ are nondegenerate phase functions in $V_j\subset \R^d\times \R^{N_j}$ such that $\imath_{t}^{-1}(\Lambda_{\Phi_{t}^j})$, $j=1,2,  \ldots $, is a locally finite covering  of $\Lambda$ with open sets for each $t$ fixed.

We denote the class of these oscillatory $\frac{1}{2}$-densities by $I^{m}(M^d, \Lambda_t; \Omega^{\frac{1}{2}}(M^d))$ or simply by $I^{m}(M^d, \Lambda_t)$. In order to simplify the notations we denote the immersion $\imath_t:\Lambda\to T^\ast M^d$ by $\Lambda_t$. 

To any oscillatory integral $u_t(x,\lambda)$ of the form \eqref{oscillatory1} one can associate  a principal symbol
\begin{equation}
e^{i\lambda f_t}\, \sigma_t \quad \mbox{where}\quad \sigma_t= \Big(\frac{\lambda}{2\pi} \Big)^m    \sigma_{1,t}\otimes \sigma_{2,t}
\label{eq:principal-symbol0}
\end{equation}
$t\to \sigma_{1,t}$ is a $C^1$ family of sections of the half-density bundle $\Omega^{\frac{1}{2}}(\Lambda)$ and $\sigma_{2,t}$ is a section of the Keller-Maslov bundle $M(\Lambda_t)$ for each $t$ fixed (cf.  \cite{Du, Mein, GS}).

In any local chart the half-density part $\sigma_{1,t}$ can be written in terms of  the nondegenerate phase functions $\Phi_t$ and the leading part $a_{0,t}$ of the amplitude $t$ in \eqref{oscillatory1} as follows 
\begin{equation}
\big(\imath_t^{-1}\circ\imath_{\Phi_t}\big)^\ast(\sigma_{1,t}) = a_{t,0} \big| d_{C_{\Phi_t}} \big|^{\frac{1}{2}}
\label{eq:principal-symbol-density}
\end{equation} 
(cf.  \cite{Hor}, Sect. 25.3), where $ d_{C_{\Phi_t}}$ is a Leray form on $ C_{\Phi_t}$, i.e.  $ d_{C_{\Phi_t}}={\it k}_t^\ast \big(\widetilde d_{C_{\Phi_t}}\big)$ is the pull-back via the inclusion map ${\it k}_t: C_{\Phi_t} \rightarrow  \R^d\times \R^N$ of a form $\widetilde d_{C_{\Phi_t}}$ such that 
\[
\widetilde d_{C_{\Phi_t}}\wedge d \frac{\partial {\Phi_t}}{\partial \theta_1}\wedge \cdots \wedge d \frac{\partial {\Phi_t}}{\partial \theta_N} = dx_1 \wedge \cdots \wedge dx_d \wedge d\theta_1 \wedge \cdots \wedge d\theta_N \, .
\]
Given for any $t$ a suitable  system of coordinates $\mu=(\mu_1,\ldots,\mu_d)$ on $C_{\Phi_t}$ extended to a neighborhood of $C_{\Phi_t}$ one obtains
\begin{equation}
\left\{
\begin{array}{lcrr}
d_{C_{\Phi_t}}= b_t\,  d\mu_1 \wedge \cdots \wedge d\mu_d\quad {\rm with }\\[0.3cm]
\displaystyle b_t= \frac{dx_1 \wedge \cdots \wedge dx_d \wedge d\theta_1 \wedge \cdots \wedge d\theta_N }
{d\mu_1 \wedge \cdots \wedge d\mu_d\wedge d \frac{\partial {\Phi_t}}{ \partial \theta_1}\wedge \cdots \wedge d \frac{\partial {\Phi_t}}{\partial \theta_N}} = \left|\frac{D(\mu,(\Phi_t)^\prime_\theta)}{D(x,\theta)}\right|^{-1}. 
\end{array}
\right.
\label{eq:principal-symbol-density1}
\end{equation}
 More generally, given a vector bundle $E$ over $M^d$, we denote by $I^{m}(M^d, \Lambda_t; \Omega^{\frac{1}{2}}(M^d)\otimes E)$ the corresponding class of oscillatory $\frac{1}{2}$-densities of order $m$ with values in the space of sections $\Gamma(E)$, and by $S^{m}(\Lambda_t,  \Omega^{\frac{1}{2}}(\Lambda)\otimes M(\Lambda_t)\otimes E_t)$ the corresponding class of symbols, where $E_t$ is the lifting of $E$ to $\Lambda_t$. 

Given two manifolds $M_j$, $j=1,2$, we denote by $\omega_j$ the corresponding canonical symplectic  forms  on $T^\ast(M_j)$ and consider the symplectic manifold 
$T^\ast(M_2)\times T^\ast(M_1)$ equipped with the exact symplectic form $\omega_2-\omega_1$. A $C^1$ family of  (exact) canonical relations ${\mathcal C}_t$, $t\in [0,\delta]$,  ``from $T^\ast(M_1)$ to $T^\ast(M_2)$'' is  given by a $C^1$ family of (exact) Lagrange immersions $\imath_t : {\mathcal C} \to T^\ast(M_2)\times T^\ast(M_1)$. To any $C^1$ family of (exact) canonical relations ${\mathcal C}_t$ one associates    a $C^1$ family of  (exact)  Lagrangian submanifolds  ${\mathcal C}^\prime_t$ of $T^\ast(M_2\times M_1)$ defined  by the exact  Lagrange immersions $\imath^\prime_t : {\mathcal C} \to T^\ast(M_2\times M_1)$ where $\imath'_t = \jmath\circ\imath_t$ and 
\begin{equation} \label{lagrangian-manifold}
\jmath : T^\ast(M_2\times M_1) \to T^\ast(M_2)\times T^\ast(M_1),\quad
\jmath(x_2,x_1,\xi_2,\xi_1)= (x_2,\xi_2,x_1,-\xi_1). 
\end{equation}
We  use the same notations as in \cite{Hor}, Sect. 25,  for the corresponding classes of $\lambda$-FIOs. Given vector bundles $E_j$ on $M_j$ and a $C^1$ family of exact canonical relations  ${\mathcal C}_t$  from $T^\ast(M_1)$ to $T^\ast(M_2)$ we say that 
\[
A_t : C_0^\infty\big(M_1, \Omega^{\frac{1}{2}}(M_1)\otimes E_1\big) \to C^\infty\big(M_2, \Omega^{\frac{1}{2}}(M_2)\otimes E_2\big)
\]
is a $C^1$ family of $\lambda$-FIOs of order $m$ if the family of the corresponding  Schwartz kernels $K_{A_t}$ is a $C^1$ family of oscillatory $\frac{1}{2}$-densities belonging to   $I^{m}\big(M_2\times M_1, {\mathcal C}_t; \Omega^{\frac{1}{2}}(M_2\times M_1)\otimes {\rm Hom}(E_1,E_2)\big)$. 
The composition of $\lambda$-FIOs with exact canonical relations having transversal and more generally a clean composition  can be defined in the same way as in the case of classical FIOs \cite{GS, Mein}. The microlocal calculus is even simpler since the amplitudes are uniformly compactly supported with respect to $\theta$. In particular we have the following analogue of Theorem 25.2.4 \cite{Hor} (see \cite{Du}, \cite{GS})
\begin{Theorem}\label{theo:pseudo-to-integral}
Let $P_t$ be a $C^1$ family of classical $\lambda$-PDOs of order $0$ acting on $\frac{1}{2}$-densities in $M_2$ with principal symbol $ p_t$ and subprincipal symbol $c_t$. Let ${\mathcal C}_t$ be a $C^1$ family of exact canonical relations from  $T^\ast(M_1)$ to $T^\ast(M_2)$ with Schwartz kernels $K_{A_t}\in I^{k}\big(M_2\times M_1, {\mathcal C}_t^\prime; \Omega^{\frac{1}{2}}(M_2\times M_1)\big)$ with principal symbols
$ e^{i\lambda f_t} \sigma_t $. Suppose that $p_t$ vanishes on the projection of ${\mathcal C}_t^\prime$ to $X_2$. Then $P_t A_t$ is a $C^1$ family of $\lambda$-FIOs of order $k-1$ with kernels  $K_{P_t A_t}$ in $I^{k-1}\big(M_2\times M_1, {\mathcal C}_t^\prime; \Omega^{\frac{1}{2}}(M_2\times M_1)\big)$ and principal symbols
\[
e^{i\lambda f_t} \Big( i^{-1}{\mathcal L}_{X_{p_t}}\sigma_t + c_t\sigma_t\Big)
\]
where $X_{p_t}$ is the Hamiltonian vector field of $p_t$ lifted to functions in $T^\ast(M_2\times M_1)$ and ${\mathcal L}_{X_{p_t}}$ is the Lie derivative. 
\end{Theorem}

\subsubsection{Quantization of $C^1$ families of billiard ball maps.}\label{subsubsec:quantization-billiard}

The aim of this section is to construct a family of monodromy operators quantizing   billiard ball maps of a $C^1$ family of billiard tables. The monodromy operators will arise as  boundary values of the  microlocal outgoing parametrizes   
$H_t(\lambda):L^2(\Gamma)\rightarrow  C^\infty( \widetilde X )$, $t\in [0,\delta]$, 
of the Dirichlet problem for the Helmholtz equation. 
We will construct $H_t(\lambda)$ for $t\in [0,\delta]$  as a $C^1$ family of $\lambda$-FIOs satisfying 
asymptotically   the   Helmholtz equation at high frequencies ($|\lambda|\to\infty$), i.e. 
\begin{equation}\label{parametrix1}
\forall\,  N\in \N\, , \quad (\Delta_t-\lambda^2)H_t(\lambda)u = O_N(|\lambda|^{-N})u
\end{equation} 
 in a neighborhood of $X$ in a  smooth extension $(\widetilde X, g_t)$ of the  Riemannian manifold 
 of $(X,g_t)$. Hereafter,
$$
O_N(|\lambda|^{-N})\, :\ L^2(\Gamma)
\to  L^2 (\widetilde X)
$$ 
stands for a $C^1$ family with respect to $t$ of   operators $A_t(\lambda): L^2(\Gamma) \to L^2(\widetilde X)$ depending on $\lambda\in {\mathcal D}$ such that 
\[
\|A_t(\lambda)\|_{L^2} \le C_{N}(1 + |\lambda|)^{-N}  
\]
for each $t$ and $\lambda\in {\mathcal D}$ where 
$C_{N}>0$ is a constant independent of $t$ and of $\lambda$.  Moreover, 
   $u$  are suitable ``initial data'' on $\Gamma$.  Set 
\begin{equation} \label{eq:Lambda} 
\Lambda:= \{(s,\rho)\in \R\times T^\ast\Gamma:\, \rho\in U ,\, -\varepsilon < s \ < T_0( \rho) + 2\varepsilon\}, 
\end{equation}
where $U$ is a compact  subset of the domain of definition  $\widetilde B_0^\ast \Gamma$ of the billiard ball map $B_0$, $0<\varepsilon \ll 1$,  and $T_t:U\to (0,+\infty)$ is the ``return time function'' which assigns  to each $ \rho\in U$ the time of the first impact at the boundary, i.e.  the 
smallest positive  $s=T_t( \rho)$ such that 
$\exp(sX_{ h_t})(\pi^+_{t}( \rho))\in \Sigma_t^-.$ Recall from  Sect. \ref{subsec:billiard-ball} that $ h_t$ is the Hamiltonian corresponding to the Riemannian metric $ g_t$ via the Legendre transform, $X_ { h_t}$ is the corresponding Hamiltonian vector field, and the map $\pi_{t} ^+:B^\ast \Gamma\to \Sigma_t^+$ is defined by \eqref{eq:outgoing-vector}. In particular,  $\exp(s X_ { h_t})(\pi_{t} ^+( \rho))$ lies on the cosphere bundle  
\begin{equation}
\widetilde \Sigma_t:= S_t^\ast \widetilde X= \{(x,\xi)\in T^\ast \widetilde X:\,  h_t(x,\xi)=1\} .
\label{cosphere}
\end{equation}
The FIOs $H_t(\lambda)$, $t\in [0,\delta]$, will be associated to
the $C^1$ family of canonical relations ${\mathcal C}_t$ in  $T^\ast \widetilde X \times T^\ast\Gamma$  given by the $C^1$ family of immersions
\begin{equation}\label{canonical-relation}
\imath_t: \Lambda \to T^\ast \widetilde X \times T^\ast\Gamma,\quad \imath_t(s,\rho)=\big(\exp\big(s X_ { h_t}\big)(\pi_{t} ^+( \rho)),\rho\big).
\end{equation}
Choosing $\delta>0$ sufficiently small we suppose  that the set 
$U$ in \eqref{eq:Lambda} is a connected open subset of $T^\ast\Gamma$ such that 
\begin{itemize}
	\item  $U$ is contained in $\widetilde B_t^\ast \Gamma$ for any $t\in [0,\delta]$;
	\item $T_t(\rho)<T_0( \rho) +\varepsilon$ for any $t\in [0,\delta]$ and $\rho\in \overline U$.
\end{itemize}
 Then $T_t$ is a smooth function on $\overline{U}$,  its image  is a compact  interval and there exist $0<a<b$ such that  $T_t(U)\subset [a,b]$ for any $t$. 
Moreover, Lemma A.1. in \cite{PT4} 
implies that 
\[
\imath_t^\prime = \jmath\circ \imath_t : \Lambda \to T^\ast (\widetilde X \times \Gamma) . 
\]
is a $C^1$ family of exact Lagrangian immersions which will be denoted by ${\mathcal C}_t'$, $t\in [0,\delta]$. 
We choose the corresponding function $f_t$ in \eqref{eq:f-t} to be just the action $f_t(s,\rho)=2s$  on the  bicharacteristic arc associated  with $\imath_t (s,\rho)\in {\mathcal C}_t$, where $(s,\rho)\in \Lambda \subset \R\times T^\ast \Gamma$.

Our aim now is to define  the immersed Lagrangian manifold  ${\mathcal C}_t'$ locally  by a nondegenerate
phase function. Fix $t^0\in [0,\delta]$ and 
take  $  \varrho^0 = (x^0,y^0,\xi^0, -\eta^0)\in {\mathcal C}_{t^0}' $.   Choose a smooth submanifold  $M^0$ of $\widetilde X$ of dimension $n-1$ passing through $x^0$ and  transversal at $x^0$ to the geodesic of $g_{t^0}$ starting from $y^0$ with codirection $(\eta^0)^+$. 
Consider the symplectic map $\chi_t: U^0\to T^\ast M^0$ defined in a neighborhood $U^0$ of $ \rho^0=(y^0,\eta^0)$ in $\widetilde B_t^\ast\Gamma$ by 
$$
\chi_t( \rho) = (\pi\circ\exp(s( \rho)X_{  h_t})\circ \pi_{t}^+)( \rho), \quad \rho\in U^0\subset \widetilde B_t^\ast\Gamma, 
$$
where $s( \rho)>0$ is the arrival time  at $T^\ast X_{|M^0}$ and $\pi(x',x_n,\xi',\xi_n)=(x',\xi')$. If $M^0=\Gamma$, this is just the billiard ball map $B_t$ defined in Sect. \ref{subsec:billiard-ball}. Denote by
${\mathcal C}'_{\chi_t}\subset T^\ast(M\times\Gamma)$ the Lagrangian manifold corresponding to the canonical relation ${\mathcal C}_{\chi_t}:=\{(\chi_t( \rho), \rho):\,  \rho\in U^0\}$.

Let $J\ni t\to x_t=(x_t',x_{t,n})\in C^\infty({\mathcal O}, \R^n)$ be a $C^1$ family of normal coordinates to $M^0$ with respect to the metrics $g_t$, where $J\subset [0,\delta]$ is an interval containing $t^0$ and ${\mathcal O}$ is a sufficiently small neighborhood of $x^0$. For any fixed $t\in J$ we have $M^0\cap {\mathcal O}=\{x_n=0\}\cap {\mathcal O}$  and the normal vector field to $M^0\cap {\mathcal O}$ associated to $g_t$ and determined by $\xi^0(\nu_t(x^0))>0$  becomes $\nu_t=(0,\ldots,0,1)$ in these coordinates.   Then  the Hamiltonian $  h_t$ is of the form  
\begin{equation}\label{eq:namiltonian-in-normal-coordinates}
 h_t(x,\xi) = \xi_n^2 +   r_t(x,\xi')
\end{equation}
in these coordinates, where $J\ni t\to r_t$ is a $C^1$ family of smooth functions in a neighborhood of $(x^{0},\xi^{0\prime})$.   If $x^0\in\Gamma$, we take $M^0$ to be a neighborhood of $x^0$ in $\Gamma$, then $  r_t(x',0,\xi')= h_t^0(x',\xi')$ is the Hamiltonian corresponding to the induced Riemannian metric $g_t^0$ on $\Gamma$. 

Following the proof of H\"ormander \cite{Hor}, Proposition 25.3.3, we can find local coordinates $y\in \R^{n-1}$ in a neighborhood of $y^0$ in $\Gamma$ such that projection ${\mathcal C}'_{\chi_{t^0}}\ni (x',\xi',y,\eta) \to (x',\eta)\in T^\ast \R^{n-1}$ is a local diffeomorphism in a neighborhood of $ (x^{0\prime},\xi^{0\prime},y^{0},\eta^{0})$. Shrinking $J$ if necessary we obtain that the map ${\mathcal C}'_{\chi_{t}}\ni (x',\xi',y,\eta) \to (x',\eta)\in T^\ast \R^{n-1}$ is a local diffeomorphism as well for any $t\in J$. Then there exists a $C^1$ family of smooth functions $\phi_t^0$ defined in a neighborhood $V^0$ of  $(x^{0\prime}, \eta^0)$ in $ \R^{n-1}\times \R^{n-1}$
such that
\begin{equation}
\det \frac{\partial^2\phi_t^0}{\partial x'\partial \eta} (x', \eta) \neq 0 \quad \mbox{for  $(x^{\prime}, \eta)\in V^0$.}
                                             \label{nondegenerate-phase-function0}
\end{equation}
and 
\begin{equation}\label{eq:graph_of_chi}
{\rm graph\, } \big(\chi_t\big)\,  =\,  \{(x',(\phi_t^0)'_{x'}(x',\eta); (\phi_t^0)'_{\eta}(x',\eta),\eta); (x',\eta)\in V^0\} 
\end{equation}
(see \cite{Hor}, Theorem  22.2.18). 
Then solving a suitable Hamilton-Jacobi equation we obtain a $C^1$ family of 
nondegenerate  phase functions 
\begin{equation}
\Phi_t(x,y,  \eta)= \phi_t(x, \eta) - \langle y, \eta\rangle 
                                             \label{phase-function}
\end{equation}
in a neighborhood  of $(x^0,y^0, \eta^0)$  
in $\R^n\times\R^{n-1}\times\R^{n-1}$ generating locally ${\mathcal C}_t'$ in a neighborhood ${\mathcal C}' _{\Phi_t}$ of $  \varrho_0$ 
where $\phi_t(x'0, \eta)=\phi_t^0(x', \eta)$. In particular, we have
\begin{equation}
\det \frac{\partial^2\phi_t}{\partial x'\partial \eta} (x, \eta) \neq 0 \quad \mbox{in a neighborhood $V$ of $(x^0, \eta^0)$.}
                                             \label{nondegenerate-phase-function}
\end{equation}
We summarize this construction as follows. 
\begin{Lemma}\label{lemma:phase-function}  There exists an open interval interval $J\subset [0,\delta]$ containing $t^0$,  local coordinates $y\in\R^{n-1}$ in a neighborhood  of $y^0$ in $\Gamma$ and independent of $t$, a neighborhood $V^0\subset  \R^{n-1}\times\R^{n-1}$ of $(x^{0\prime}, \eta^0)$ and a $C^1$ family of function $\phi_t^0\in C^\infty(V^0)$  satisfying \eqref{nondegenerate-phase-function0}  and such that the following holds 
\begin{enumerate}
	\item  the function $\Phi_t^0(x',y, \eta):=\phi_t^0(x', \eta) - \langle y, \eta\rangle$
is  a local  generating function of the Lagrangian manifold ${\mathcal C}'_{\chi_t}\subset T^\ast(M^0\times\Gamma)$  for every $t\in J$; 
	\item the Lagrangian manifold ${\mathcal C}_t'$  is defined in a neighborhood of $ \varrho_0$ by a  phase function 
\[
\Phi_t(x,y,  \eta)= \phi_t (x, \eta) - \langle y, \eta\rangle, 
\]
where $\phi_t(x, \eta)$ is  a local solution of the Hamilton-Jacobi equation
\begin{equation}
\partial_{x_n} \phi_t(x, \eta)\, =\,  \sqrt{1-  r_t(x,\partial_{x'}\phi_t(x, \eta))}\, , \quad \phi_t(x',0, \eta)=\phi_t^0(x', \eta),
\label{hamilton-jacobi1}
\end{equation}
and $ r_t$ is given by \eqref{eq:namiltonian-in-normal-coordinates}. 
	\item  there exists a neighborhood $V\subset \R^n\times \R^{n-1}$ of $(x^0,\eta^0)$ such that the family of functions $J\ni t\to \phi_t\in C^\infty(V,\R)$ is $C^1$ and \eqref{nondegenerate-phase-function} holds for any $t\in J$. 
\end{enumerate}
\end{Lemma}
In particular,  
\begin{equation}
\Phi_t(x,y,  \eta)  = \phi_t^0(x', \eta)-\langle y, \eta\rangle  + x_n\phi_t^1(x,\eta) ,
\label{hamilton-jacobi}
\end{equation}
where $J\ni t \to \phi_t^1\in C^\infty(V,\R)$ is a $C^1$ map. 
We take  $\phi^0_t(x', \eta)=\langle x', \eta\rangle$ if the image of  $  \varrho_0$ by the involution in \eqref{lagrangian-manifold}   belongs to $\Sigma^+_{t^0}\times U$, which means that $x^0=(y^0,0)$ and $\xi^0 = (\eta^0)^+$.   

In order to obtain  the Maslov part of the principal symbol picked up by the phase functions constructed by the Lemma we need the following. 
Fix $\widetilde t^0\in [0,\delta]$, 
take  $ \widetilde \varrho^0 = (\widetilde x^0,\widetilde y^0,\widetilde \xi^0, -\widetilde \eta^0)\in {\mathcal C}_{\tilde t^0}' $ and denote by   $\widetilde M^0$ the corresponding submanifold   transversal at $\widetilde x^0$ to the geodesic of $g_{\tilde t^0}$ starting from $\widetilde y^0$ with codirection $(\widetilde \eta^0)^+$ and by 
 $\widetilde \chi_t: \widetilde U^0\to T^\ast \widetilde M^0$  the corresponding symplectic map. Let $\widetilde J$ be the corresponding interval about $\widetilde t^0$ and $\widetilde \Phi_t(\widetilde x, \widetilde y, \widetilde \eta)= \phi_t(\widetilde x,  \widetilde \eta) - \langle \widetilde y, \widetilde \eta\rangle$, $t\in \widetilde J$,  the corresponding $C^1$ family of phase functions given by Lemma \ref{lemma:phase-function}. Suppose  that $\tau\in J\cap \widetilde J \neq \emptyset$ and that there exists $\zeta\in \Lambda$ such that 
\[
\imath_\tau^\prime(\zeta) \, \in \,   {\mathcal C}'_{\Phi_t} \cap  {\mathcal C}'_{\widetilde \Phi_t} .
\]
\begin{Lemma}\label{lemma:signature} There exists a neighborhood $I\subset J\cap \widetilde J$ of $\tau$ and a neighborhood $V$ of  $\zeta$ in $\Lambda$ such that the  function $\mu: V\times I \to 2\Z$ defined by
\[
\begin{array}{rcll}
\mu(\varrho,t)\, :=\  {\rm sgn}\,  (\phi_t)''_{ \eta \eta}(x,\eta)- {\rm sgn}\,  (\widetilde \phi''_t)_{\widetilde  \eta\widetilde \eta}(\widetilde x,\widetilde\eta), \\[0.3cm]
\imath_t(\rho)\, =\,  \imath_{\Phi_t} (x, (\phi_t)'_{\eta}(x,\eta), \eta)= \imath_{\widetilde\Phi_t} (\widetilde x, (\widetilde\phi_t)'_{\widetilde\eta}(\widetilde x,\widetilde\eta), \widetilde\eta).
\end{array}
\]
is constant on $V\times I$. 
\end{Lemma}
{\em Proof}. The assurtion follows from an argument in \cite{D-G} using 
  H\"{o}rmander's index $\sigma(M_1,M_2;L_1,L_2)\in\Z$  of four Lagrangian spaces $M_1,M_2,L_1,L_2$ in the Lagrangian Grassmannian $\Lambda(n-1)$, where $L_1$ and $L_2$ are transversal to both $M_1$ and $M_2$ in $T^\ast \R^{n-1}$ (see \cite{Du}, Definition 3.4.2). It is known that $\sigma$ is locally constant and continuous with respect to all the variables  $(M_1,M_2,L_1,L_2)$. Set $\imath_\tau(\zeta)= (u, (\phi_\tau)'_u(u,v),  (\phi_\tau)'_v(u,v),v)\in T^\ast \widetilde X\times T^\ast \Gamma$ where $\phi_t$, $t\in J$, is the phase function in Lemma \ref{lemma:phase-function} corresponding to the coordinates $x_t: {\mathcal O} \to \R^n$. 
Take a section $M$ in $\widetilde X$ passing trough the point $u$ and transversal to the geodesic starting from $\partial_v \phi_\tau(u,v)\in\Gamma$ and havig a codirection $v^+$. We can suppose that  $M= \{ x_n= g(x')\}$ in these coordinates with some smooth function $g$. Let us change the $x$ coordinates in a neighborhood of $M$ by $z'=x'$ and $z_n=x_n- g(x')$ and set  $\psi_t(z,\eta)=\phi_t(z',z_n+g(z'),\eta)$. Consider the (local) symplectic transformation $\chi_t^0:T^\ast \Gamma \to T^\ast M$ defined by 
$$
\chi_t^0( \rho) = \big(\pi^0\circ\exp(s_t^0( \rho)X_{  h_t})\circ \pi_{t}^+\big)( \rho), \quad \rho\in U_0\subset \widetilde B_t^\ast\Gamma, 
$$
for $t$ sufficiently close to $\tau$, where $s_t^0( \rho)>0$ is the arrival time  at $T^\ast X_{|M}$ and $\pi^0(z',z_n,\xi',\xi_n)=(z',\xi')$. Then $(z',\eta)\to \psi_t(z',0,\eta)$ is a generating function of $\chi_t^0$ in the sense of \eqref{eq:graph_of_chi}. 
Given $(x,y,\eta)\in  C_{\Phi_t}$ in a neighborhood of $(u, \partial_u \phi_\tau(v),v)$ with $x\in M$,  
 we obtain as in \cite{D-G} p. 69
\[
\begin{array}{lcrr}
{\rm sgn}\, (\Phi_{t})''_{ \eta \eta}(x',g(x'),y,\eta)=  {\rm sgn}\, (\psi_{t})''_{ \eta \eta}(z',0,\eta)
=  {\rm sgn}\, (V,H_t; (d\chi_t^0)^{-1}(V))
\end{array}
\]
where ${\rm sgn}\, (M_1,M_2; L)$ is defined in \cite{Du}, Definition 3.4.3, $V=\{(\delta z',\delta \xi'):\, \delta z'=0\}$ is the vertical space (the tangent space to the fiber)  and   $H_t$ is the horizontal space $\{(\delta z',\delta \xi'):\, \delta \xi'=0\}$ for the local coordinates $x_t$ in ${\mathcal O}$ used in the construction of $\Phi_{t}$. Repeating this procedure in $\widetilde {\mathcal O}$ for the phase function $\widetilde \Phi_t$ obtained by Lemma \ref{lemma:phase-function} corresponding to the coordinates $\widetilde x_t: {\mathcal O} \to \R^n$, we obtain 
\[
\begin{array}{lcrr}
{\rm sgn}\, (\Phi_{t})''_{ \eta \eta}(x,y,\eta)-{\rm sgn}\, (\widetilde \Phi_{t})''_{\widetilde \eta \widetilde\eta}(\widetilde x,\widetilde y,\widetilde \eta)\\[0.3cm] 
= {\rm sgn}\, (V,H_t; (d\chi_t^0)^{-1}(V))- {\rm sgn}\, (V, \widetilde H_t; (d\chi_t^0)^{-1}(V))
= 2 \sigma(H_t, \widetilde H_t ; (d\chi_t^0)^{-1}(V),V)
\end{array}
\]
where   $\widetilde H_t$ is the horizontal space $\{(\delta z',\delta \xi'):\, \delta \xi'=0\}$ for the corresponding local coordinates $\widetilde x_t: \widetilde {\mathcal O}\to \R^n$ used in the construction of $\widetilde \Phi_{t}$. This shows that $\mu$ is idependent of $t$ and of $\varrho$ in a small neighborhood of $(\tau,\zeta)$. 
\finishproof

Using the phase functions obtained in Lemma \ref{lemma:phase-function} one can define the space of $\lambda$-FIOs  corresponding to the $C^1$ family of canonical relations ${\mathcal C}_t$.  We are  looking for solutions $H_t(\lambda)$ of \eqref{parametrix1} with Schwartz kernels in  $I^{-1/4}(\widetilde X \times\Gamma, {\mathcal C}_t^\prime)$.
To any $C^1$ family of nondegenerate  phase function $\Phi_t(x,y, \eta)$ of the form \eqref{phase-function} generating ${\mathcal C}'_t$ in a neighborhood of a point $  \varrho_0=(x^0,y^0,\xi^0, -\eta^0)\in {\mathcal C}'$ ($\Phi_t$ is given by Lemma \ref{lemma:phase-function}) there is a $C^1$ family of classical amplitudes 
\[
b_t(x, \eta,\lambda)\sim b_{0,t}(x, \eta) + b_{1,t}(x, \eta)\lambda^{-1}+ \cdots
\]  
such that  
the Schwartz kernel of $H_t(\lambda)$ can be written microlocally near  $  \varrho_0$ as a $C^1$ family of oscillatory $\frac{1}{2}$-densities
\begin{equation}\label{oscillatory}
I_{\Phi_t}(x,y,\lambda) = \Big(\frac{\lambda}{2\pi} \Big)^{n-1}\, \Big(
\int_{{\R}^{n-1}}\ e^{i\lambda \Phi_t(x,y, \eta)} b_t(x, \eta,\lambda)  d \eta \Big)|dx|^{1/2}|dy|^{1/2}. 
\end{equation}
(see  \eqref{oscillatory1}).

 Notice that the Hamiltonian $p_t$ in $T^\ast(X\times\Gamma)$ obtained by lifting of the principal symbol $  h_t-1$ of the operator $\lambda^{-2}\Delta_t - 1$  vanishes on ${\mathcal C}_t^\prime$. Thus to compute the principal symbol of $(\Delta_t - \lambda^2)H_t(\lambda)$ we can use Theorem \ref{theo:pseudo-to-integral}. Note also that the corresponding subprincipal symbol is $c_t=0$. 

We are going to define suitable sections $\sigma_{1,t}$ and $\sigma_{2,t}$  of the half density bundle and of the Keller-Maslov bundle of ${\mathcal C}_t^\prime$. 
The lifting of the Hamiltonian vector field $X_{  h_t}$ to $T^\ast X\times T^\ast \Gamma$ is  $Y_t=(X_{  h_t},0)$ and its flow $S_t^\tau$ restricted to ${\mathcal C}_t$ is given by
\begin{equation}
S_t^\tau(x,\xi, y,\eta)=S_t^\tau(\exp(sX_{  h})(\pi^+_\Sigma(y,\eta)),y,\eta)= (\exp((s+\tau)X_{  h_t})(\pi^+_t(y,\eta)),y,\eta)
\label{flow}
\end{equation}
for any $(x,\xi, y,\eta)\in {\mathcal C}_t$. 
The  volume form $\beta_t$ on ${\mathcal C}_t$  given by the pull-back by $\imath_t$  of 
\[
ds\wedge (dy_1 \wedge d\eta_1)\wedge \cdots  \wedge (dy_{n-1} \wedge d\eta_{n-1}) \in \Omega(\Lambda)
\] 
is invariant with respect to the flow $S_t^\tau$ or equivalently, the Lie derivative ${\mathcal L}_{Y_t} \beta_0$ vanishes. 
Then the Lie derivative of the $\frac{1}{2}$-density $\sigma_{0,t}:= |\jmath^\ast(\beta_0)|^{\frac{1}{2}}\in \big|\Omega\big|^{\frac{1}{2}}({\mathcal C}_t^\prime)$ with respect to $X_{p_t}$ is zero. Recall that $\jmath$ is given by \eqref{lagrangian-manifold} and that $X_{p_t}$ is the Hamiltonian vector field of $  h_t-1$   lifted to functions in $T^\ast(\widetilde X\times \Gamma)$. We set 
\begin{equation}
\sigma_{1,t}= \widetilde b_{0,t} \sigma_{0,t}, \quad \widetilde b_{0,t}\in C_0^\infty({\mathcal C}_t'). 
\label{principal-symbol2}
\end{equation}
The relation between $\widetilde b_{0,t}$ and the principal part $b_{0,t}$ of the amplitude $b_{t}$ in \eqref{oscillatory} is obtained  in \cite{PT4}, (A.23). More precisely, let us denote by $\eta_t=\eta_t(x,\xi')$ the local solution of  $\xi'=(\phi_t)'_{x'}(x,\eta)$ obtained by the implicit function theorem and set 
$$
\widetilde b_{0,t}'(x,\xi') = \widetilde b_{0,t}(\pi_1^{-1}(x,\xi',\xi_n)), 
$$
where 
$\pi_1: {\mathcal C}_t' \to \widetilde \Sigma_t \subset T^\ast \widetilde X$ is the  projection  $\pi_1(x,y,\xi,-\eta)=(x,\xi)$.   Then \cite{PT4}, (A.23), yields
\begin{equation}
b_{0,t}(x,\eta_t(x,\xi')) = \frac{\widetilde b_{0,t}'(x,\xi')}{\sqrt{2|\xi_n|}}\left|\det (\phi_t)''_{x' \eta}(x,\eta(x,\xi')) \right|^{\frac{1}{2}} 
\label{principal-symbol3}
\end{equation} 
in a neighborhood of $(x^0,\xi^0)$, where 
$\xi_n =  \sqrt{1-r_t(x,\xi')}. $

The Keller-Maslov bundle $M({\mathcal C}_t')$ of ${\mathcal C}_t'$ admits a natural trivialization by locally constant sections. Recall from  H\"ormander \cite{Hor1}, p. 148,  that a section of the line bundle 
$M({\mathcal C}'_t)$ is given by a family of functions $f_{\Phi_t}:{\mathcal C}'_{\Phi_t}\to \C$, where $\Phi_t$ is a nondegenerate phase function  generating locally ${\mathcal C}'_t$ at ${\mathcal C}'_{\Phi_t}=\imath_{\Phi_t} (C_{\Phi_t})$ (see \eqref{lagrangian1}) such that
$\displaystyle f_{\widetilde {\Phi_t}}= i^{\mu_t} f_{\Phi}$ on ${\mathcal C}'_{\widetilde {\Phi_t}}\cap {\mathcal C}'_{ {\Phi_t}}$. The  function  
$\mu_t=\mu_{\widetilde {\Phi_t} {\Phi_t}}$ is defined  by 
\begin{equation}
\mu_{\widetilde {\Phi_t} {\Phi_t}}( \varrho):=\frac{1}{2}\left(({\rm sgn}\,  (\Phi_t)''_{ \theta \theta}(x,y, \theta)-N)- ({\rm sgn}\,  \widetilde(\Phi_t)''_{\widetilde  \theta\widetilde \theta}(x,y,\theta)-\widetilde N)\right), 
\label{transition-function}
\end{equation}
where 
$\theta\in \R^N$, $\widetilde\theta\in \R^{\widetilde N}$, $(\Phi_t)'_\theta(x,y,\theta)=0$, $\widetilde (\Phi_t)'_{\widetilde\theta}(x,y,\widetilde\theta)=0$ and 
$ \varrho = \imath_{\Phi_t}(x,y,\theta)= \imath_{\widetilde {\Phi_t}}(x,y,\widetilde \theta)\in {\mathcal C}'_{\widetilde {\Phi_t}}\cap {\mathcal C}'_{{\Phi_t}}$. Moreover,  $\mu_{\widetilde {\Phi_t} {\Phi_t}}\in \Z$ and it is constant on each connected component of ${\mathcal C}'_{\widetilde {\Phi_t}}\cap {\mathcal C}'_{{\Phi_t}}$. 
The section will be called ``natural'' if $f_{{\Phi_t}}$ are constant functions taking values in $\{i^k:\, k\in \Z\}$.   
 In our case $M({\mathcal C}_t')$ can be trivialized in a band $|s|<\epsilon$ using the phase functions $\Phi_t$ given by Lemma \ref{lemma:phase-function}, where $\phi_t^0(x', \eta)=\langle x', \eta\rangle$. Then  
 $(\Phi_t)''_{ \eta \eta}(y,0,y, \eta)=0$  in view  \eqref{hamilton-jacobi} and we get  ${\rm sgn}\,  (\Phi_t)''_{ \eta \eta}(y,0,y, \eta)=0$.  
 This yields a natural trivialization of the Keller-Maslov bundle in a band ${\mathcal C}'\cap \{|s|<\epsilon\}$ for some $\epsilon>0$, choosing a locally constant  section which equals $1$ in that band.   In particular, the Lie derivative ${\mathcal L}_{Y_t} \sigma_{2,t}$ vanishes for each $t$. This argument holds whenever ${\mathcal C}'|_{s=0}$ is a conormal bundle of a smooth submanifold (see \cite{Hor1} Sect. 3.3 and  \cite{D-G}, p. 65). 

Using Lemma \ref{lemma:signature} one can obtain 
 a  natural section $\sigma_{2,t}$ of $M({\mathcal C}_t')$ which is independent of $t$ in a small neighborhood of $t^0$ for any $t^0\in [0,\delta]$. 
The section $\sigma_{2,t}$ of $M({\mathcal C}'_t)$ can be described as in \cite{D-G}, \cite{Hor1} and \cite{Mein} as a Maslov index of a suitable path. Let $ \varrho_1=(x^0,y^0,\xi^0, -\eta^0)\in {\mathcal C}_{t^0}'$ and $(x^0,\xi^0)=\exp (T X_{h_{t^0}})(y^0,\eta^0)$. Let $M$ be a submanifold transversal 
to the corresponding geodesic in $x^0$ and
 let $\Phi_t$ be a $C^1$ family of  generating function of ${\mathcal C}_t'$ in a neighborhood of $ \varrho_1$ given by Lemma \ref{lemma:phase-function}.  
Consider the path 
$\widetilde\gamma_t$ on  ${\mathcal C}_t'$ defined by 
$$
\pi_1(\widetilde\gamma_t(s))= \exp(sX_{  h_t})(\pi^+_t(y^0, \eta^0)),\ s\in [0,T].
$$ 
We have 
$$
\widetilde\gamma_{t^0}(0):=  \varrho_0=(y^0,y^0,(\eta^0)^+, -\eta^0)\in {\mathcal C}_{t^0}'|_{s=0}\quad \mbox{and}\quad \widetilde\gamma(T)= \varrho_1 \in {\mathcal C}_{t^0}'.
$$
 Choose   a partition $0=s_0<s_1<\cdots<s_k=T$  and phase functions $\Phi_{t,j}$, $j=1,\ldots,k$, as in Lemma \ref{lemma:phase-function} generating locally 
 ${\mathcal C}_t'$ in a neighborhood of $\widetilde\gamma_t(s_j)$ for $t$ in a small neighborhood of $t^0$ and  such that 
$$
\widetilde\gamma_t([s_{j-1},s_j])\subset {\mathcal C}'_{\Phi_{t,j}},\quad \Phi_{t,k}=\Phi_t \quad \mbox{and}\quad  
\Phi_{t,1}(x',0,y, \eta)=\langle x'-y, \eta\rangle. 
$$
Then  trivializing $M({\mathcal C}_t')$ in a neighborhood of $ \varrho_1$ by the phase function $\Phi_t$ we get
\begin{equation} \label{maslov-index}
\begin{array}{lcrr}
\displaystyle (\sigma_{2t})_{\Phi_t}\,  =\,   i^{\mu(\widetilde\gamma_t)}, \quad \mbox{where}\\[0.3cm] 
\displaystyle \mu(\widetilde\gamma_t)\, :=\, \frac{1}{2}\sum_{j=1}^{k-1}\Big( {\rm sgn}\, (\Phi_{t,j})''_{ \eta \eta}(\widetilde\gamma_t(s_{j}))-{\rm sgn}\, (\Phi_{t,j+1})''_{ \eta \eta}(\widetilde\gamma_t(s_{j}))\Big)\in \Z . 
\end{array}
\end{equation}  
Now Lemma \ref{lemma:signature} implies that $\mu(\widetilde\gamma_t)$ is independent of $t\in I$ where $I$ is a sufficiently small  neighborhood of  $t^0$ in $[0,\delta]$. In other to construct $\sigma_{2,t}$ one  can use finitely  many paths $\widetilde\gamma_t$ since ${\mathcal C}_{t^0}$ is compact, hence $I$ can be chosen to be common for all the paths.   We set
\[
\sigma_t\,  =\,  \sigma_{1,t}\times \sigma_{2,t}\, =\, \widetilde b_{0,t} \, \sigma_{0,t} \times \sigma_{2,t}.
\]
According to Theorem \ref{theo:pseudo-to-integral} the oscillatory integral $(\Delta_t-\lambda^2)K_{H_t}(x,y,\lambda)$ belongs to $I^{3/4}(\tilde X\times \Gamma, {\mathcal C}_t')$ and its principal symbol is just the Lie derivative ${\mathcal L}_{Y_t}\sigma_t$  multiplied by $(\lambda/2\pi)^{3/4}$  since the subprincipal symbol of the Laplace-Beltrami operator is $0$. 
Moreover,  the Lie derivative with respect to $Y$ of the sections $\sigma_{0,t}$ and $\sigma_{2,t}$  vanishes, hence,  the transport  equation  ${\mathcal L}_{Y_t}\sigma_t=0$ becomes 
\begin{equation}
(S^t)^\ast \tilde b_{0,t} = \tilde b_{0,t}. 
\label{transport-equation}
\end{equation}
Multiplying  $\tilde b_0$ with a suitable cut-off function, which  equals $1$ in a neighborhood of 
${\mathcal C}_t'\cap T^\ast(X\times \Gamma)$,  we can suppose that  $\tilde b_0$ has a compact support with respect to $(s,y,\eta)\in \Lambda$.  In this way we obtain a $C^1$ family of $\lambda$-FIOs
 $H_{0,t}(\lambda)$ with Schwartz kernels in   $I^{-1/4}(\widetilde X \times \Gamma, {\mathcal C}_t')$ such that the  Schwartz kernel of 
$(\Delta_t-\lambda^2)H_{0,t}(\lambda)$ belongs to $I^{-1/4}( X\times \Gamma, {\mathcal C}_t')$. Repeating this procedure we get an operator $H_{1,t}(\lambda)$ such that $H_0(\lambda)+H_1(\lambda)$ solves (\ref{parametrix1}) modulo  a $\lambda$-FIO of order $-5/4$ and so on. 
The initial data $\tilde b_0|_{s=0}$  will be determined by Lemma \ref{Lemma:restriction} below. 

Denote by   
$\imath_\Gamma^\ast: C^\infty(\widetilde X) \to  C^\infty(\Gamma)$ the operator of restriction  $\imath_\Gamma^\ast(u)=u_{|\Gamma}$.  We would like to represent $\imath_\Gamma^\ast$ microlocally as a $\lambda$-FIO. To this end, denote by ${\mathcal N}$ the conormal bundle of the graph of the inclusion map $\imath_\Gamma: \Gamma \to \widetilde X$ and by ${\mathcal R}= {\mathcal N}^{-1}$ the corresponding inverse canonical relation. In other words,
\[
{\mathcal R}:= \{(x,\xi;x,\widetilde \xi)\in T^\ast \Gamma\times T^\ast \widetilde X:\, x\in\Gamma, \xi=\widetilde \xi|_{T_x\Gamma}\}. 
\]
The operator $\imath_\Gamma^\ast $ can be considered microlocally  as a  $\lambda$-FIO with Schwartz kernel of the class
$I^{1/4}(\Gamma\times \tilde X,  {\mathcal R}; \Omega^{\frac{1}{2}})$  which means that the composition $\imath_\Gamma^\ast \circ A(\lambda)$ belongs to that class for any classical $\lambda$-PDO $A(\lambda)$ of order $0$.  Moreover, its   principal symbol can be identified with $(\lambda/2\pi)^{1/4}$ modulo the corresponding $\frac{1}{2}$-density (see  \cite{PT4}, Sect. A.1.4). 
In what follows, we shall investigate the composition $\imath_\Gamma^\ast H_t(\lambda)$ of $\lambda$-FIOs. Firstly, notice that the composition ${\cal R}\circ {\mathcal C}_t$ of the corresponding canonical relations is transversal (see \cite{PT4},  Sect. A.1.4). 
Recall that $\pi_1: {\mathcal C}_t' \to T^\ast X$ and $\pi_2: {\mathcal C}_t' \to T^\ast \Gamma$ are given by 
$\pi_1(x,y,\xi,-\eta)= (x,\xi)$ and $\pi_2(x,y,\xi,-\eta)= (y,\eta)$. Denote by $dv( \rho):= dy\wedge d\eta$  the symplectic volume form on $T^\ast\Gamma$.   Recall that $\nu_t(x)\in T_x\widetilde X\big|_{\Gamma}$ is the unit inward normal to $\Gamma$ related to the metric $g_t$ and that $\pi_t^\pm (x,\xi)= (x,\xi_t^\pm)\in \Sigma_t^{\pm} $ for $(x,\xi)\in B_t^\ast\Gamma$. Moreover,   
\[
\langle\xi_t^\pm,\nu_t\rangle(x,\xi):= \langle\xi_t^\pm(x,\xi),\nu_t(x)\rangle =\pm \sqrt{1-r_t(x,\xi')}
\]
in the normal coordinates used in Lemma \ref{lemma:phase-function}. 
Using \eqref{principal-symbol2} and \eqref{principal-symbol3} and the theorem about the composition of $h$-FIOs one obtains
\begin{Lemma}\label{Lemma:restriction}
The composition of canonical relations 
${\cal R}\circ {\mathcal C}_t$ is transversal    and it is a disjoint union $\Delta^0 \sqcup {\mathcal C}_t^0$ of
the diagonal $\Delta^0$ in $U\times U$ (for $s=0$) and  the graph ${\mathcal C}_t^0$ of the billiard ball map
$B_t:U\to B_t(U)$ (for $s=T$). 
 Moreover, 
\begin{equation}
\imath_\Gamma^\ast H_t(\lambda)\,  =\,  P_t(\lambda) + 
G_t(\lambda)
+ O_M(|\lambda|^{-M})
\, ,  
                     \label{boundary-trace1}
\end{equation}
where  $P_t(\lambda)$ is a $C^1$ family of classical $\lambda$-PDOs on $\Gamma$ of order $0$ and $G_t(\lambda)$ is a $C^1$ family of $\lambda$-FIOs with Schwartz kernels in  $I^{0}(\Gamma, \Gamma,  {\mathcal C}_t^{0\prime})$.   The principal symbol of the operator 
$P_t(\lambda)$ can be identified with 
\begin{equation}
\widetilde b_0( \pi_2^{-1}(  \rho))\, (2| \langle \xi_t^+,\nu_t\rangle( \rho) |)^{-1/2}\, |dv( \rho)|^{1/2}, \quad  \rho \in U.  
\label{eq:boundary-trace-symbol}
\end{equation} 
The principal symbol of  $G_t(\lambda)$  can be identified with  
\begin{equation}
\widetilde b_{0,t}( \pi_1^{-1}( \pi_t^-( \rho)))\, |2\langle\xi_t^-,\nu_t\rangle( \rho)|^{-1/2}\, e^{i\lambda A_{\gamma_t(\rho)}}\, |dv( \rho)|^{1/2}\otimes \sigma_t', \quad  \rho \in B_t(U), 
\label{eq:boundary-trace-symbol1}
\end{equation}
where 
$A_{\gamma_t}=
\int_{\gamma_t}\xi dx $ is the action  along the integral curve $\gamma_t$ of the Hamiltonian vector
field $X_{h_t}$  starting from  $\pi^+_t( B_t^{-1}(\rho))$ and  with endpoint  $\pi^{-}_t( \rho)$ and $\sigma_t'$ is a natural section of the Maslov bundle $M({\mathcal C}_t^{0\prime})$. Moreover, for each $t^0\in [0,\delta]$ one can choose $\sigma_t'$ to be independent of $t$ in a neighborhood of $t^0$. 
\end{Lemma}
The Lemma is proved in  \cite{PT4}, Sect. A.1.4. 

Let  $\Psi(\lambda)$ be a classical $\lambda$-PDO of order $0$ with frequency set in $ U$ and principal symbol $\Psi_0( \rho)$, $ \rho\in U$. 
 We take $\Psi(\lambda)$ as initial data of $H_0(\lambda)$ as $s=0$ setting $P_t(\lambda)=\Psi(\lambda)$ in Lemma \ref{Lemma:restriction}. Recall that $\widetilde b_0$ satisfies \eqref{transport-equation}. 
On the other hand 
\[
(x',\xi')=B_t(y,\eta) \quad \mbox{if and only if} \quad \pi_1^{-1}( x',0,\xi^-)=   S^{T_t( y,\eta)}(\pi_2^{-1}( y,\eta))
\]
where $T_t$ is the return time function. Then 
 \eqref{transport-equation} and \eqref{eq:boundary-trace-symbol} imply
\[
\widetilde b_0( \pi_1^{-1}( x',0,\xi^-)) = \widetilde b_0( (\pi_2^{-1}( y,\eta))) = \Psi_0(y,\eta)(2|\langle \eta^+,\nu\rangle(y,\eta)|)^{1/2}.
\]
Then  parameterizing ${\mathcal C}_t^{0\prime}$ by the variables $(y,\eta)\in U$ and using \eqref{eq:boundary-trace-symbol1} we write the principal symbol of   $G_t(\lambda)$  as follows
 \begin{equation}
\sigma(G_t(\lambda))\,  = \,  
\Psi_0(y,\eta) \, \frac{|\langle\eta_t^+,\nu_t\rangle(y,\eta)|^\frac{1}{2}}{|\langle \xi_t^-,\nu_t\rangle(B_t(y,\eta))|^{1/2}} \, e^{i\lambda A_t(y,\eta)}\, 
| d  y\wedge d \eta  |^\frac{1}{2}\otimes \sigma_t' , 
                                                   \label{principal-symbol}
\end{equation}
where   $A_t(y,\eta)=
\int_{\gamma_t}\xi dx $ is the action  along the integral curve $\gamma_t$  of the Hamiltonian vector
field $X_{h_t}$  starting from  $\pi^+_t(y,\eta)$ and  with endpoint  $\pi^{-}_t(B_t(y,\eta))$. 

In the same way, using Lemma \ref{Lemma:restriction} we  determine the initial conditions of $H_{1,t}(\lambda)$ and so on. In this way  we obtain a $C^1$ family of $\lambda$-FIOs 
\begin{equation}\label{eq:parametrix-H} 
H_t(\lambda)=H_{0_t}(\lambda) + H_{1,t}(\lambda) + \cdots 
\end{equation} 
with Schwartz kernels in  $ I^{-1/4}(\tilde X, \Gamma, {\mathcal C}_t^\prime)$ satisfying \eqref{parametrix1} and such that $P_t(\lambda)=\Psi(\lambda)$. From now on, to simplify the notations we drop the corresponding $\frac{1}{2}$-density.
Denote by $E_t(\lambda)$ a $C^1$ family of classical $\lambda$-PDOs of order $0$ on $\Gamma$ with  principal symbols $E_{0,t}\in C_0^\infty(\widetilde B_t^\ast \Gamma)$ such that  
\begin{equation}
E_{0,t}( \rho) = |\langle \xi_t^+,\nu_t\rangle( \rho)|^\frac{1}{2} = |\langle \xi_t^-,\nu_t\rangle( \rho)|^\frac{1}{2}
                    \label{boundary-trace5}
\end{equation}
in a compact neighborhood of $\overline U$ in $\widetilde B_t^\ast \Gamma$.  Then using Egorov's theorem,  \eqref{principal-symbol} and \eqref{boundary-trace5} we obtain $\imath_\Gamma^\ast H_t(\lambda)$. We summarize this construction in the following
\begin{Prop}\label{prop:parametrix} 
The $C^1$ family of $\lambda$-FIOs operators $t\to H_t(\lambda)$ gives for any half density in $f$ in $L^2(\Gamma)$ a family of solution $u_t = H_t(\lambda)f$ of 
\[
(\Delta_t -\lambda^2)u_t = O_N(|\lambda|^{-N})f.
\]
Moreover, 
\begin{equation}
\imath_\Gamma^\ast H_t(\lambda) = \Psi(\lambda) + G_t(\lambda) + O_M(|\lambda|^{-M}) \, ,\quad 
G_t(\lambda) = E_t(\lambda)^{-1}G_t^0(\lambda)E_t(\lambda)
\, ,  
                     \label{boundary-trace3}
\end{equation}
where $E_t(\lambda)$ is a family $\lambda$-PDOs of order 0 which are  of elliptic microlocally in a neighborhood of  $\mbox{WF}_\lambda'(\Psi)$,  the principal symbol of $G_t^0(\lambda)$ can be identified with 
\begin{equation}
\Psi_0( \rho)e^{i\lambda A_t( \rho)}\,  |dv( \rho)|^\frac{1}{2}\otimes \sigma_t' , \quad  \rho \in U,  
                      \label{boundary-trace4}
\end{equation}
and $\sigma_t'$ could be chosen to be independent of $t$ in a small neighborhood of any $t_0$. 
\end{Prop}
In particular,  the frequency 
set $WF'$  of
$G_t(\lambda)$  
is contained in 
$B_t(U)\times U$.

We are going to estimate the $L^2$-norm of $u_t= H_t(\lambda)f$, where $f$ is a $\frac{1}{2}$-density on $\Gamma$. 
Consider the  $L^2$-adjoint operator $H_t(\lambda)^\ast$ of $H_t(\lambda)$ which is well-defined for any $\lambda\in {\mathcal D}$ fixed as an operator
from  $L^2(\widetilde X)$ to $L^2(\Gamma)$. Moreover, it can be considered as a $\lambda$-FIO associated with the canonical relation ${\mathcal C}_t^{-1}$ the Schwartz kernel of which belongs to $ I^{-1/4}(\Gamma, \tilde X,  ({\mathcal C}_t^{-1})^\prime)$. As in  \cite{PT4}, Sect. A.1.4, we obtain 
\begin{Prop}\label{Lemma:continuity}
The family $t\to C_t(\lambda):=H_t(\lambda)^\ast H_t(\lambda): L^2(\Gamma) \to L^2(\Gamma)$ is a $C^1$ family od  classical $\lambda$-PDOs of order $0$. The  principal symbol of $C_t(\lambda)$ can be identified with by 
$$ 
C_{0,t}(y,\eta)\, :=\,  \int_\R |\widetilde b_0 (s,y,\eta)|^2 ds\, ,\quad (y,\eta),\ (x,\eta) \in U.
$$
Moreover, 
$$ 
C_{0,t}(y,\eta)\, \ge\,  T_t(y,\eta) | \langle \xi_t^+,\nu_t\rangle(y,\eta) |\Psi_0(y,\eta)|^2\, ,\quad (y,\eta)\in U, 
$$
where $T_t$ is the return time function. In particular, there exists $C>1$ such that
\[
C^{-1} \|\Psi(\lambda)f\|_{L^2(\Gamma)} \le \|H_t(\lambda)f\|_{L^2(X)} \le C \|\Psi(\lambda)f\|_{L^2(\Gamma)} 
\]
for each $(t,\lambda)\in [0,\delta]\times {\mathcal D}$. 
\end{Prop}

\subsection{Reduction to the boundary}\label{subsec:reduction}
The reduction to the boundary  is a variant of the reflection method for the wave equation. We shall describe it in the case of Dirichlet boundary conditions. In the case of Neumann and more generally of Robin boundary conditions it is done in \cite{PT4}.  

Denote by $(\widetilde X,g_t)$ a $C^\infty$ extension of  $(X,g_t)$ and by $h_t$ the Hamiltonian corresponding to $g_t$ via the Legendre transform. 
Consider a $C^1$ family of Kronecker invariant tori  $[0,\delta]\ni t\to \Lambda_t(\omega) \subset \widetilde {\bf B}_t^\ast \Gamma$ of $B_t$ having frequencies in the set $\Omega_\kappa^0$ of points of positive Lebesgue density in  $\Omega_\kappa=(\Omega-\kappa)\cap D(\kappa,\tau)$, where $\Omega= B(\omega_0,\varepsilon)$,  $0<\kappa<\varepsilon/2\ll 1$ and $0\le \delta\le 1$. For  $0<\delta\ll 1$ such families of Kronecker invariant tori of $B_t$ are provided by  Theorem \ref{Theo:soft-BNF}. 
Denote by
$$
{\mathcal T}_{t}^j\, :\, =\cup\{B_t^j(\Lambda_{t}(\omega))\, :\ \omega\in\Omega_\kappa^0\}  \subset \widetilde B_{t}^\ast \Gamma
$$
the cooresponding union of the invariant tori of $P_t\circ B_t^j$ for $0\le j <m$ and set
\[
{\mathcal T}^j\, :=\, \cup \{{\mathcal T}_{t}^j,\, 0\le j <m.
\]  
Fix $t_0\in [0,\delta]$ and choose   open sets $U_j \subset V_j \subset T^\ast\Gamma$ for  $0\le j \le m$ and a sufficiently small   interval $J\subset [0,\delta]$ around  $t_0$ such that 
\[
{\mathcal T}^j \subset U_j \subset\!\subset V_j \subset\!\subset \widetilde {\bf B}_t \Gamma \quad \mbox{and} \quad B_t(V_j)\subset U_{j+1} 
\]
for each $t\in J$ and $0\le j\le m-1$  and
\[
\overline U_0 \cup \overline U_m \subset V_m  \subset\!\subset \widetilde {\bf B}_t \Gamma
\]
for each $t\in J$. The relation $U \subset\!\subset V$ means here that $\overline U \subset V$ where $\overline U$ is the closure of $U$. 
We suppose that  the $C^1$ family of exact symplectic mappings 
\[
J\ni t \to P_t
\]
admits a $C^1$ family of BNFs in a neighborhood $U$ of $\overline V_0$ in the sense of Definition \ref{Def:BNF} (see also Theorem \ref{Theo:soft-BNF}). In other words, we suppose that there exist $C^1$-smooth with respect to $t\in J$ families of exact symplectic diffeomorphisms 
$\chi_t:\A \to \chi_t(\A) \subset U$   and of real valued functions  $L_t\in C^\infty(D)$ and $R_t^0\in C^\infty(\A)$ where $\A=\T^{n-1} \times D$ and $D=\nabla L_{t_0}^\ast(\Omega)$  such that for  each $t\in J$ the following holds
\begin{enumerate}
\item $\overline V_0 \subset \chi_t(\A)$; 
\item $\Lambda_t(\omega)=\chi_t(\T^{n-1}\times \{I_t(\omega)\})\subset U_0$ for $\omega\in  \Omega_{\kappa}^0$, where $I_t(\omega)$ is given by \eqref{eq:momentum-I}; 
\item The function 
$$
\R^{n-1}\times D\ni (x,I)\mapsto \phi_t(x,I):=\langle x,  I\rangle -L_t(I) - R_t^0(x,I)
$$ 
is a generating function of the exact symplectic map
\[
P_t^0:= \chi_t^{-1}\circ P_t\circ \chi_t: \A\to\A  
\]
in the sense of Definition \ref{def:generating-function}; 
\item $\nabla L_t: D \to \Omega$ is a diffeomorphism and  $L_t=L_{t_0}$  outside $D^1:= \nabla L_{t_0}^\ast(\Omega-\kappa/2)$; 
\item $R_t^0$ is flat at $\T^{n-1}\times E_t^{\kappa}$, where $E_t^{\kappa} =\nabla L_t^\ast( \Omega_{\kappa}^0)$. 
\end{enumerate}
Chose the $\lambda$-PDO $\Psi(\lambda)$ giving the ``initial data'' of the  operators $H_t(\lambda)$ in   \eqref{boundary-trace3} such that 
\begin{equation}\label{eq:initial-data-psi}
\mbox{WF}_\lambda'(\Psi-{\rm Id}))\cap \overline V_j = \emptyset\quad  \forall\, 0 \le j \le m. 
\end{equation}
Recall from Proposition \ref{prop:parametrix} that
\[
\imath_\Gamma^\ast H_t(\lambda) = \Psi(\lambda) + G_t(\lambda) + O_M(|\lambda|^{-M})
\]
where $G_t(\lambda)$ is described in \eqref{boundary-trace3} and \eqref{boundary-trace3}. 
Take now a classical $\lambda$-pseudodifferential operator $\lambda$-PDO $\Psi_0(\lambda)$ such that 
\[
 \mbox{WF}_\lambda'(\Psi_0) \subset V_0 \quad \mbox{and}  \quad \mbox{WF}_\lambda'(\Psi_0 -{\rm Id})\cap \overline U_0=\emptyset. 
\]
Consider the ``outgoing'' solution of 
 the Helmlotz equation 
\begin{equation}\label{eq:parametrix-Helmholtz}
(\Delta_t -\lambda^2)u_t = O_N(|\lambda|^{-N})f
\end{equation}
for $\lambda\in {\mathcal D}$ and  $t\in J$ with ``initial data''   $\Psi_0(\lambda) f$ which is given by $u_t:= H_t^0(\lambda)f$, where 
\[
H_t^0(\lambda)=  H_t(\lambda)\Psi_0(\lambda).
\] 
Recall that
$O_N(|\lambda|^{-N})\, :\ L^2(\Gamma)
\to  L^2 (\widetilde X)$ 
stands here  for a  family  of   operators 
\[ 
A_t(\lambda): L^2(\Gamma) \to L^2(\widetilde X)
\]
such that 
\[
\|A_t(\lambda)\|_{L^2} \le C_{N}(1 + |\lambda|)^{-N}  
\]
for each $t\in J$ and $\lambda\in {\mathcal D}$ where 
$C_{N}>0$ is constant independent of $t$ and of $\lambda$.  Then
\[
\imath_\Gamma^\ast H_t^0(\lambda) = \Psi_0(\lambda) + G_t(\lambda)\Psi_0(\lambda) + O_M(|\lambda|^{-M})
\]
since  $\mbox{WF}_\lambda'((\Psi-{\rm Id})\Psi_0)= \emptyset$ in view of \eqref{eq:initial-data-psi}.

To satisfy the ``boundary conditions'' on $ U^1$ in the case when $m\ge 2$ we use the reflexion method.  Let $\Psi_1(\lambda)$ be a classical $\lambda$-PDO such that 
\[
 \mbox{WF}_\lambda'(\Psi_1) \subset V_1 \quad \mbox{and}  \quad \mbox{WF}_\lambda'(\Psi_1 -{\rm Id})\cap \overline U_1=\emptyset. 
\] 
Set   
\[
H_t^1(\lambda)=  H_t(\lambda)\Psi_1(\lambda)G_t(\lambda)\Psi_0(\lambda)
\] 
and  consider 
\[
u_t(\lambda)= \widetilde{H}_t(\lambda)f:= H_t^0(\lambda)f -  H_t^1(\lambda)f.
\] 
Then $u_t$ satisfies \eqref{eq:parametrix-Helmholtz} and it satisfies microlocally the Dirichlet boundary conditions on $U_1$. Notice that $\mbox{WF}_\lambda(\imath_\Gamma^\ast(u_t)) \subset U_0\cap U_2$. 
Similarly if $m>2$ one can treat  the boundary conditions in  $U_j$ for 
any $0<j<m$ which leads to a solution $u_t= \widetilde H_t(\lambda)f$   satisfying the boundary conditions microlocally in  $U_j$  for each  $0<j<m$. Let $\Psi_j(\lambda)$, $0 \le j \le m-1$ be a classical $\lambda$-PDO such that 
\[
 \mbox{WF}_\lambda'(\Psi_j) \subset V_j \quad \mbox{and}  \quad \mbox{WF}_\lambda'(\Psi_j -{\rm Id})\cap \overline U_j=\emptyset. 
\] 
We set $H_t^0(\lambda)=  H_t(\lambda)\Psi_0(\lambda)$ if $m=1$ and 
\begin{equation}\label{eq:parametrix}
\left\{
\begin{array}{rcll}
\displaystyle \widetilde H_t(\lambda)f &=& \displaystyle  \sum_{j=0}^{m-1}(-1)^j H_t^j(\lambda)f \quad \mbox{where} \\[0.3cm]
\displaystyle  H_t^{j}(\lambda)&=& H_t(\lambda)\Psi_{m-1}(\lambda)G_t(\lambda)\Psi_{m-2}(\lambda) \cdots G_t(\lambda)\Psi_{0} (\lambda) 
\end{array}
\right.
\end{equation}
if $m \ge 2$. Then  $u_t$ satisfies \eqref{eq:parametrix-Helmholtz} and  
\[
\mbox{WF}_\lambda(\imath_\Gamma^\ast(u_t)) \subset U_0\cap B_t(U_{m-1}) \subset\!\subset V_m. 
\] 
More precisely, 
\[
\imath_\Gamma^\ast (u^t)= \Psi_0(\lambda)f - M_t(\lambda)  f + O_N(\lambda^{-N})f 
\]
where $M_t(\lambda):= -G_t(\lambda)\Psi_0(\lambda)$ for $m=1$ and 
\[
\displaystyle  M_t(\lambda):= (-1)^{m-1}G_t(\lambda)\Psi_{m-1}(\lambda)\cdots G_t(\lambda)\Psi_{0}(\lambda)
\]
if $m\ge 2$. Taking into account \eqref{boundary-trace3} we obtain $M_t(\lambda) = E(\lambda)^{-1} M_t^0(\lambda)E(\lambda)$, where 
\begin{equation}\label{eq:the-operator-S}
M_t^0(\lambda):= (-1)^{m-1}Q_t(\lambda)\Psi_{m-1}(\lambda)\cdots Q_t(\lambda)\Psi_{0}(\lambda).
\end{equation}
Moreover, using Proposition \ref{prop:parametrix} and the theorem about the composition of $\lambda$-FIOs (here we use it in the simple case of canonical transformations) and parameterizing ${\rm graph} (P_t)\subset V_m\times V_0$ by its projection on $V_0$ we obtain that for each $t\in J$ the principal symbol of $M_t^0(\lambda)$ is given by
\[
(-1)^{m-1}\exp (i\lambda A_t(x,\xi))|dv( \rho)|^\frac{1}{2}\otimes\sigma'_m
\]  
over  $U_0$, where    
\[
A_t(x,\xi) =
\sum_{j=0}^{m-1}A_t(x_t^j,\xi_t^j), \quad 
(x_t^j,\xi_t^j)= B_t^j (x,\xi),
\]  
is the action along the corresponding broken geodesic and $\sigma'_m$ is a ``natural'' section of the corresponding Keller-Maslov  bundle which can be chosen to be independent of $t\in J$. 

Let $\psi_0(\lambda)$ be a classical $\lambda$-PDO of order zero such that 
\[
\mbox{WF}'(\psi_0(\lambda))\subset U_0\quad  \mbox{and}\quad \mbox{WF}'(\psi_0(\lambda)- {\rm Id})\cap {\cal T}^0 =\emptyset. 
\]
We summarize the above construction  by the following 
\begin{Prop} \label{prop:monodromy1}
Let  $v_{t,\lambda}\in L^2(\Gamma)$ and  $u_{t,\lambda}= \widetilde H_t(\lambda)\psi_0(\lambda)v_{t,\lambda}$ where $(t,\lambda)\in J\times {\mathcal D}$. Then 
\[
\left\{
\begin{array}{rcll}
(\Delta_t-\lambda^2)u_{t,\lambda}\ &=&\ 
O_N(\lambda_q^{-N})u_{t,\lambda} \, , \\
\imath^\ast_\Gamma \,  u_{t,\lambda}  \ &=&\  
 O_N(\lambda_q^{-N})u_{t,\lambda}  
 \end{array}
\right.
\] 
if and only if 
\begin{equation}
(  M_t^0(\lambda) -  \mbox{\rm Id}\, )\psi_0(\lambda)v_ {t,\lambda}\ = \ 
O_N(\lambda_q ^{-N})\, v_ {t,\lambda} .
\label{eq:quasimode-equation}
\end{equation}
\end{Prop}
The  structure of the monodromy operator $M_t^0(\lambda)$ is given by
\begin{Prop} \label{prop:monodromy2}
The canonical relation of 
$  M_t^0(\lambda):=E_t(\lambda) M_t(\lambda)E_t(\lambda)^{-1}$ is given by
the graph ${\rm graph}\, (P_t)\subset V_m\times V_0$ of the symplectic map 
$P_t= B^m_t:V_0 \to V_m$, which is $C^1$ with respect to $t$. The family 
$J\ni t\to M_t^0(\lambda)$ of  classical $\lambda$-FIO 
of order 0 with a
large parameter $\lambda\in {\cal D}$ is  $C^1$ smooth with respect to $t\in J$. 
Parameterizing ${\rm graph}\, P_t$ by its projection on $V_0$ for $t\in J$,  the principal symbol of $M_t^0(\lambda)$ becomes 
$$
\sigma(M_t^0)\, =\, (-1)^{m-1}\exp (i\lambda A_t(x,\xi))\, |dv( \rho)|^\frac{1}{2}\otimes\sigma'_m
$$  over  $U_0$, where    $\sigma'_m$ is a ``natural'' section of the corresponding Keller-Maslov  bundle which does not depend on $t\in J$. 
\end{Prop}

\subsection{Quantum Birkhoff Normal  Form}\label{subsec:QBNF}
Using the $C^1$ family of exact symplectic transformations $\chi_t$ given by Theorem \ref{Theo:soft-BNF}   
 we identify the first cohomology groups 
$H^1(\Lambda_t(\omega),\Z) =  H^1(\T^{n-1},\Z)=\Z^{n-1}$ for $\omega\in \Omega_\kappa^0$ and $t\in J$, 
and we denote by $\vartheta_0\in \Z^{n-1}$ 
the Maslov class of the invariant
tori $\Lambda_t(\omega)$. Notice that $\vartheta_0\in \Z^{n-1}$ does not depend on $t\in J$ and  $\omega\in \Omega_\kappa^0$. 
Consider as in \cite{CV} 
the flat Hermitian line bundle $\LL$ over
$\T^{n-1}$  associated to  the representation $\varrho:\Z^{n-1}\to SU(1)$ of the fundamental one group $\pi_1(\T^{n-1})=\Z^{n-1}$ defined  by
$\varrho(k)=\exp \left(i\frac{\pi}{2}\langle \vartheta_0,k\rangle\right)$, $k\in\Z^{n-1}$ (see \cite{Kob}, Sect. 1.2).  More precisely, $\LL$ is the quotient of $\R^{n-1}\times \C$ by the action of $\Z^{n-1}$ given by $k.(x,z)=(x+2\pi k,\varrho(k)z)$.  
Then sections $s$
of ${\LL}$ can be identified canonically with smooth functions
$\widetilde{s}:\R^{n-1}
\rightarrow \C$ such  that
\begin{equation}
\displaystyle \widetilde{s}(x +2\pi k)\ =\ e^{i\frac{\pi}{2}\langle 
\vartheta_0,k\rangle}
\widetilde{s}(x) \quad \forall \, x\in \R^{n-1},\,  k\in \Z^{n-1}. 
                                              \label{sections}
\end{equation}
An orthonormal basis of $L^2(\T^{n-1},\LL)$ 
is given by
$e_k,\ k\in \Z^{n-1}$,  where 
\[
\widetilde e_k (x)\ =\ 
\exp\left( i \langle k + \vartheta_0/4,x\rangle
\right).
\] 
We quantize the family of exact symplectic  transformations $\chi_t: \A=\T^{n-1}\times D \to T^\ast \Gamma$ as in \cite{CV}, Sect. 5 and  \cite{PT4}, Sect. 3.3.  
Denote by ${\mathcal C}_{\chi_t}$ the graph of $\chi_t$ in $T^\ast \Gamma \times T^\ast \T^{n-1} $ and by ${\mathcal C}'_{\chi_t}= \jmath({\mathcal C}_{\chi_t})$ the corresponding Lagrangian submanifold of $T^\ast (\Gamma \times \T^{n-1})$, where $\jmath$ is defined in \eqref{lagrangian-manifold}. Consider the class of $\lambda$-FIOs $T_t(\lambda):C^\infty(\T^{n-1}, \LL) \rightarrow 
C^\infty(\Gamma,\C)$ of order $0$ associated with the canonical relation ${\mathcal C}_{\chi_t}$. The Schwartz kernel  $K_{T_t(\lambda)}$ of $T_t(\lambda)$ belongs to the class 
$I^{0}(\Gamma\times \T^{n-1}, {\mathcal C}'_{\chi_t}; p_2^\ast( \LL))$, where $p_2:\Gamma\times \T^{n-1}\to \T^{n-1}$ is the projection on the second factor.  
Recall from \cite{CV}, Sect. 5,  that the principal symbol $\sigma(K_{T_t})(\lambda)$ of $K_{T_t(\lambda)}$  can be canonically identified with a smooth function in $T^\ast \T^{n-1}$. Indeed,  $\sigma(K_{_t}T)(\lambda)$ belongs to the symbol class  
$S^0({\mathcal C}_{\chi_t}^\prime, M({\mathcal C}^\prime_{\chi_t})\otimes \pi_{2}^\ast  (\LL'))$, where $\pi_{2}: {\mathcal C}_{\chi_t}^\prime \to \T^{n-1}\times D$ and $\pi_2\circ \jmath : {\mathcal C}_{\chi_t} \to  \T^{n-1}\times D$ is the restriction at ${\mathcal C}_{\chi_t}$ of the projection $T^\ast \Gamma\times T^\ast\T^{n-1} \to T^\ast\T^{n-1} $ on the second factor while $\LL'$ is the dual bundle to $\LL$ (the base manifold of $\LL$ and $\LL'$ here is $\T^{n-1}\times D$ instead of $\T^{n-1}$). On the other hand, $M({\mathcal C}'_{\chi_t})= \pi_{2}^\ast  (\LL)$  and 
using the parametrization of ${\mathcal C}'_{\chi_t}$ given by $\pi_{2}$ we  identify the above  class of symbols  with $S^0(\T^{n-1}\times D, \LL\otimes \LL')$ which can be canonically identified with $C_0^\infty(\T^{n-1}\times D)$ since $\LL\otimes\LL'$ is trivial (cf. \cite{Hirz}, Chapter I,  3.7). This allows us to obtain  a $\lambda$-FIO $T_t(\lambda)$ of order $0$ associated to the canonical relation ${\mathcal C}_{\chi_t}$, which is  microlocally unitary  over $\A^0:=\T^{n-1}\times D^0$, where $D^0$ 
 is a neighborhood of $\cup_{t\in J} E_t^\kappa$ in  $D$ and  
\[
E_t^\kappa = \nabla L_t^\ast(\Omega_\kappa^0) = I_t(\Omega_\kappa^0) 
\]
has been defined in {\it 4}, Theorem \ref{Theo:soft-BNF}. This means that
$$
\mbox{WF}'(T_t(\lambda)^\ast T_t(\lambda) - \mbox{Id})\cap \A^0 =\emptyset. 
$$
Trivializing the $\frac{1}{2}$-density bundle of ${\mathcal C}^\prime_{\chi_t}$ by $\pi_2^\ast|dv|^\frac{1}{2}$, where $dv$ is the symplectic volume form on $T^\ast\T^{n-1}$, we take the principal symbol of $T_t(\lambda)$ to be  equal to one 
in $\T^{n-1}\times D^0$ modulo a Liouville
factor  $\exp(i\lambda \Psi_t(\varphi,I))$, 
where  the function $\Psi_t$ is real valued. Consider the $C^1$ family of $\lambda$-FIOs of order zero 
\[
M_t^1(\lambda):=T_t(\lambda)^{\ast}  M_t^0(\lambda) T_t(\lambda) :   C^\infty(\T^{n-1}, \LL) \rightarrow C^\infty(\T^{n-1}, \LL).
\] 
The corresponding  canonical relation ${\mathcal C}_t$   is just the graph of $P_t^0= {\chi_t}^{-1}\circ P_t \circ {\chi_t}$  i.e. 
 \begin{equation}\label{eq:canonical-relation-C-t}
{\mathcal C}_t:= \{(P_t^0( \rho), \rho):  \rho\in\A \}. 
\end{equation} 
 Denote by ${\mathcal C}_t'$ the corresponding Lagrangian submanifold of $T^\ast(\T^{n-1}\times \T^{n-1})$.  
Using the theorem about the composition of $\lambda$-FIOs in the special case of canonical transformations we obtain   that 
the Schwartz kernel of $M_t^1(\lambda)$ belongs to  
 $I^{0}( \T^{n-1} \times \T^{n-1}, {\mathcal C}_t'; M({\mathcal C}^\prime_{t})\otimes{\rm End}\, (\LL))$. Let us find its principal symbol, parameterizing  
${\mathcal C}'$ by the variables $ \rho=(\varphi,I)\in \A$.  
\begin{Lemma}\label{Lemma:principal-symbol}
The  principal symbol of $M_t^1(\lambda)$ is given by 
\[
\sigma( M_t^1)(\lambda)= (-1)^m\exp (i\lambda f_t)s_{t,0}\otimes\sigma_0\otimes |dv(\rho)|^{1/2}
\]
 where  $s_{t,0}$ is a $C^1$ family of smooth function in $\T^{n-1}\times D$ such that $s_{t,0}(\varphi,I)= 1$ in $\T^{n-1}\times D^0$, $dv(\rho)$ the symplectic volume form on $T^\ast \T^{n-1}$, $\sigma_0$ is a  natural section of the Keller-Maslov bundle $M({\mathcal C}_t')$ independent of $t$ and 
\begin{equation}
\label{liouville-exponent}
f_t(\varphi,I)= A_t(\chi_t(\varphi,I)) + \Psi_t(\varphi,I) -
\Psi_t(P^0(\varphi,I))\, , \quad (\varphi,I)\in \T^{n-1}\times D .
\end{equation}
\end{Lemma}
{\em Proof}. 
 Notice that ${\rm End}\, (\LL) \cong \LL\otimes \LL'$ is trivial as a bundle over $\T^{n-1}\times D$, hence, 
  smooth sections  can be canonically identified with smooth functions  in  $\T^{n-1}\times D$. Then 
 parameterizing   
${\mathcal C}_t'$ by the variables $ \rho=(\varphi,I)\in \A$  and using the   $\lambda$-FIO calculus  and Proposition \ref{prop:monodromy2}  we obtain  the principal symbol of $M^1_t(\lambda)$. To prove \eqref{liouville-exponent} we 
write microlocally the Schwartz kernels of the  corresponding $\lambda$-FIOs in  as oscillatory integrals of the form  \eqref{oscillatory} with suitable phase functions  and then we  evaluate the  phase function of the composition at the stationary points. The claim that $\sigma_0$ is natural and independent of $t$ follows from the fact that  the section $\sigma_m'$ in Proposition \ref{prop:monodromy2} is natural and from  the composition law of FIOs.   \finishproof

Recall that the Lagrangian manifolds ${\mathcal C}_t'$ are generated by the $C^1$ family of functions $\Phi_t(x,y,I) = \phi_t(x,I) - \langle y, I\rangle$, where 
\[
\phi_t(x,I) =\langle x, I\rangle -L_t(I) - R_t^0(x,I)  
\]
satisfies  ${\it 3}$-${\it 5}$ in Sect. \ref{subsec:reduction} (see also Definition \ref{Def:BNF} and Theorem \ref{Theo:soft-BNF}).  
\begin{Prop}\label{operator-T} We have 
$$
T_t(\lambda)^{\ast}  M_t^0(\lambda) T_t(\lambda)
 = e^{i\pi\vartheta/2}
W_t(\lambda)
$$ 
where $\vartheta\in\Z$ is a Maslov's index independent of $t\in J$ and  
\[
J\ni t\to W_t(\lambda): C^\infty(\T^{n-1}, \LL)\to C^\infty(\T^{n-1}, \LL)
\] 
is a $C^1$ family of $\lambda$-FIOs of order zero with  canonical relations given by the graph of $P_t^0$ over $\A$. Moreover, the Schwartz kernel of $W_t(\lambda)$ is of the form 
\begin{equation}
\widetilde{W_t}(x,y,\lambda) |dx|^\frac{1}{2} |dy|^\frac{1}{2}= 
\Big(\frac{\lambda}{2\pi}\Big)^{n-1}\, \Big(\int_{\R^{n-1}} 
\, e^{i\lambda (  \phi_t(x,I) - \langle y,I\rangle)} \, w_t(x,I,\lambda) 
 \, dI \Big)\, |dx|^\frac{1}{2} |dy|^\frac{1}{2}, 
                                  \label{eq:operator}
\end{equation}
where 
$t \to w_t= \sum_{j=0}^\infty w_{t,j} $ is a $C^1$ family of classical amplitudes $2\pi$-periodic with respect to $x$ and 
$w_{t,0}(x,I)= 1$ for  $(x,I)\in \R^{n-1}\times D^0$.  
\end{Prop}
{\em Proof}. \quad The Schwartz kernel of $M_t^1(\lambda)$ can be written in the form \eqref{eq:operator} with a phase function $C+\Phi_t(x,y,I)$, where $C$ is a constant since $\Phi_t$ is a globally defined  generating function of ${\mathcal C}_t'$. We are going to show that $C=0$. Indeed, the exponent on the Liouville factor picked up by these phase functions is 
\[
C + \langle I, \nabla L_t(I)\rangle -L_t(I) + \langle I, \nabla_I R^0_t(\varphi, I)\rangle- R_t^0(\varphi,I) = f_t^0 (\varphi,I),
\]
then  taking $ (\varphi, I)\in E_t^\kappa$ and using Lemma \ref{lemma:beta}, \eqref{liouville-exponent} and the equality $R_t^0\big|_{E^t_\kappa} = 0$ we get $C=0$. 
Trivializing the Maslov bundles $M({\mathcal C}_t')$ by the $
C^1$ family of phase phase functions $\Phi_t$ we get 
$(\sigma_0)_{ \Phi_t}=\exp\left(i \frac{\pi}{2}\vartheta_1\right)$ for some $\vartheta_1\in \Z$ independent of $t$ since  $\sigma_0$ does not depend on $t$. We set $\vartheta=\vartheta + m\pi$. 
Moreover,  by \eqref{eq:principal-symbol-density} and \eqref{eq:principal-symbol-density1} we obtain that $\big| d_{C_{\widetilde \Phi}}\big| = dxdI$. Hence, $w_{t,0}(x, I)=s_{t,0}({\rm pr\,}(x), I)=1$ for each $(x,I)\in \R^{n-1}\times D^0$. \finishproof

In the case of Neumann and Robin boundary conditions we have $\vartheta=\vartheta_1$.

Our aim is to make $w_{t,j}(x,I)$ independent of the angle variable $x$ for $I\in E_t^\kappa$ conjugating ${W_t}(\lambda)$ by a suitable $C^1$ in $t$ family of  $\lambda$-PDOs which are elliptic on $\T^{n-1}\times D^0$. 
\begin{Prop}\label{prop:commutator}
There exists a $C^1$ family of $\lambda$-PDOs $J\ni t\to A_t(\lambda)$ of order $0$ 
 acting on  $C^\infty(\T^{n-1},{\LL})$ and  a $C^1$ family of 
$\lambda$-FIO $J\ni t\to W_t^0(\lambda)$ of the form (\ref{eq:operator}) 
such that 
\begin{equation}
W_t(\lambda)A_t(\lambda)\ =\ A_t(\lambda)W_t^0(\lambda) +
Z_t(\lambda)  \,  ,
\label{eq:equation}
\end{equation}
where 
\begin{enumerate}
\item[(1)] the full symbols of $A_t(\lambda)$ and of $W_t^0(\lambda)$ are
\[
\begin{array}{rcll}
\displaystyle \sigma(A_t)(\varphi,I,\lambda)&:=& \displaystyle a_t(\varphi, I,\lambda) \sim \sum_{j=0}^\infty \lambda^{-j} a_t^j(\varphi,I,\lambda) \ \mbox{and} \\[0.3cm] 
\displaystyle  \sigma(W_t^0)(\varphi,I,\lambda)&:=& \displaystyle p_t(I,\lambda) \sim \sum_{j=0}^\infty \lambda^{-j} p_t^j(I) \, ,
\end{array}
\] 
where $J\ni t\to a_t(\varphi, I,\lambda)$  and $J\ni t\to p_t(I,\lambda)$ are $C^1$ families of classical symbols
and  $a_t^0(\varphi,I)=1$ and $p_t^0(I)=1$  for $I\in D^0$, 
\item[(2)]  
$J\ni t\to Z_t(\lambda)$ is a $C^1$ family of 
$\lambda$-FIOs of order $0$ of the form 
(\ref{eq:operator}) with symbols 
\[
S_t(\varphi,I,\lambda)\sim \sum_{j=0}^\infty \lambda^{-j} S_t^j(\varphi,I)
\] 
such that the functions 
$ S_t^j$, $j\ge 0$,  are flat on $\T^{n-1}\times E_t^\kappa$. 
\end{enumerate}
\end{Prop}
{\em Proof.} The proof of the proposition is similar to that in \cite{C-P} and \cite{PT4}. First, comparing the symbols of order $-j$ of the left and the right hand side of \eqref{eq:equation} we shall derive the corresponding homological equation.  Set
\[
\phi_t^0(x,I) = L_t(I) +R_t(x,I).
\]
We write the  Schwartz kernel of the operator $W_t(\lambda)A_t(\lambda)$ of 
the form (\ref{eq:operator}) with amplitude 
\[
u_t(x,I,\lambda) = 
\left(\frac{\lambda}{2\pi}\right)^{n-1}\, \int_{{\R}^{2n-2}} 
\, e^{i\lambda (\langle x-z,\xi-I\rangle - (\phi_t^0(x,\xi) -\phi_t^0(x,I)))} \, 
w_t(x,\xi,\lambda)
a_t(z,I,\lambda)\, d\xi dz \ ,
\] 
which belongs to $C^\infty ({\T}^{n-1}\times D)$ for each $\lambda$ fixed. 
Changing the variables we write $u_t(x,I,\lambda)$ of the form 
\[
\left(\frac{\lambda}{2\pi}\right)^{n-1}\, \int_{{\R}^{2n-2}} 
\, e^{-i\lambda \langle v,\eta\rangle } \, w_t(x,I + \eta,\lambda)
a_t(v +x + K_t(I,\eta) + H_t(x,I,\eta) ,I,\lambda)\, d\eta dv \, ,
\] 
where $K_t(I,\eta):= \int_0^1 \nabla_I L_t(I +s \eta) d s$ and $ H_t(x,I,\eta) := \int_0^1 \nabla_I R_t(x,I + s \eta) d s$. Note that 
$K_t(I,0)= \nabla L_t(I)$. Moreover,   $\partial_I^\alpha\partial_\eta^\beta H_t(x,I,0)=0$ for each $I\in E_t^\kappa$ and any $\alpha,\beta\in\Z^{n-1}$ since  the function $I \to R_t(x,I )$ is flat at $E_t^\kappa$ for every $x\in\R^{n-1}$ in view of  {\em 5}, Sect. \ref{subsec:reduction}. 
 Using the Taylor formula 
for the amplitude at $v=0$ and integrating by parts we get
\[
u_t(x,I,\lambda) \sim \sum_{j=1}^\infty u_t^j(x,I)\lambda^{-j}
\] 
where
\[
u_t^{j}(x,I)\ :=\ 
\sum_{r+s+|\gamma|=j}\,
\frac{1}{\gamma!}\, \left[ D^\gamma_\eta \left(a_t^r(x,I + \eta)\, 
\partial^\gamma_x\,
a_t^s(x + K_t(I,\eta) +H_t(x,I,\eta) ,I)\right)\right]_{|\eta =0}\, . 
\]
In the same way we write the Schwartz kernel of 
 $A_t(\lambda)W_t^0(\lambda)$ in 
the form (\ref{eq:operator}) with amplitude $q_t(x,I,\lambda)$ 
given by the oscillatory integral 
\[
\left(\frac{\lambda}{2\pi}\right)^{n-1}\,   
p^0(I,\lambda)\, \int_{{\R}^{2n-2}} 
\, e^{i\lambda (\langle x-z,\xi-I\rangle - (\phi_t^0(z,I) - \phi_t^0(x,I)))} \, 
a_t(x,\xi,\lambda)
 d\xi dz \, .
\] 
Changing the variables we obtain $q_t = q_t^0 + q_t^1$, where 
$q_t^0(x,I,\lambda) = a_t(x,I,\lambda) p_t(I,\lambda)$ and
$q_t^1(x,I,\lambda)$ is given by 
\[
\left(\frac{\lambda}{2\pi}\right)^{n-1}\,  p_t(I,\lambda)\,  
\int_{{\R}^{2n-2}} 
\, e^{-i\lambda \langle v,\eta\rangle } 
[a_t(x, \eta  + I + H_t^1(x,v,I), \lambda) - a_t(x, \eta + I, \lambda)]\, 
d\eta dv \ ,
\] 
where $ H_t^1(x,v,I) = \int_0^1 \nabla_x R_t(x+\tau v,I) d \tau$. Moreover, 
all the  derivatives of $ H_t^1(x,v,I)$
vanish for $I\in E_t^\kappa$ since the function $I \to R_t(x,I )$ is flat at $E_t^\kappa$ for every $x\in\R^{n-1}$. In this way we obtain for any $j\ge 1$ that
\begin{equation}
S_t^j(\varphi, I)=a_t^j(\varphi+\nabla L_t(I), I)- a_t^j(\varphi, I) - p_t^j(I) - F^j(\varphi, I,t)\, ,
\label{eq:homological}
\end{equation}
where $F^j$ is a polynomial of $\partial_\varphi^\alpha\partial_I^\beta a_t^l$ and $\partial_I^\beta p_t^l$  for $l<j$ and  $|\alpha|+|\beta|\le 2j$ and of $\partial_I^\beta L_t$ for $|\beta|\le 2j+1$. 

We are looking for functions $a_t^j$ and $p_t^j(I)$ such that $S_t^j(\varphi, I)=0$ on $\T^{n-1}\times E_t^\kappa$.
We shall solve this equation  recursively  with respect to $j$  changing the variables by  $I = I_t(\omega)$, $\omega\in \Omega$, and we consider $\Omega$ as a subset of $\R^{n-1}$. 
Set $f(\varphi, \omega,t) := a_t^j(\varphi, I_t(\omega))$, $c(\omega,t) := p_t^j( I_t(\omega))$ and $F(\varphi, \omega,t) := F^j(\varphi, I_t(\omega),t)$. Then we get the homological equation 
\begin{equation}
f(\varphi+ \omega, \omega,t)- f(\varphi, I,t) = c(\omega,t) + F(\varphi, \omega,t)\, , \quad \omega\in\Omega_\kappa^0\, . 
\label{eq:homological1}
\end{equation}
We are looking for smooth functions  $f$ and $c$ on $\T^{n-1}\times \Omega$ and $\Omega$ respectively, which solve \eqref{eq:homological1} for every $\omega\in \Omega_\kappa^0$. We have the following
\begin{Lemma}\label{Lemma:homological}
Let $J\ni t\to F(\cdot,\cdot,t)\in C^\infty(\T^{n-1}\times \Omega)$ be a $C^1$ family of functions such that $F(\varphi,\omega,t)=0$ for each $\omega\in \Omega_\kappa^0$. Then there exist $C^1$ families functions $J\ni t\to f(\cdot,\cdot,t)\in C^\infty(\T^{n-1}\times \Omega)$ and $J\ni t\to c(\cdot,t)\in C^\infty( \Omega)$ such that the function 
\[
(\varphi,\omega) \to S(\varphi, \omega,t):=  f(\varphi + \omega, \omega,t)- f(\varphi, I,t) - c(\omega,t) - F(\varphi, \omega,t)
\]
is flat at $\T^{n-1}\times \Omega_\kappa^0$ for each  $t$ fixed. 
\end{Lemma}
{\em Proof}. Given $g\in L^1(\T^{n-1})$ we denote by $\hat g_k$, $k\in\Z^{n-1}$, its Fourier coefficients. For any $k\in\Z^{n-1}$ we have 
\[
\hat S_k(\omega,t)=  \hat f_k(\omega,t)\left(e^{ i\langle \omega, k\rangle} -1\right) - c(\omega,t) - \hat F_k(\omega,t). 
\]
We set $c(\omega,t) = - \hat F_0(\omega,t)$, which gives $\hat S_0(\omega,t)=0$. We are going to find $S_k(\omega,t)$ for $k\neq 0$. To this end we 
choose $\phi\in C^\infty_0(\R)$ such that
\[
0\le \phi\le 1, \
\phi(x)=1\ \mbox{for}\ |x|\le \pi/5\  \mbox{and}\ \phi(x)=0\ \mbox{for}\ |x|\ge \pi/4. 
\]
For any $0\neq k\in \Z^{n-1}$ set 
\[
\phi_k(x):= \sum_{j\in\Z}\phi((x - 2\pi j)|k|^{\tau} \kappa^{-1}).
\]
We have $|k|^{\tau} \kappa^{-1}\ge \kappa^{-1} >1$, hence, $\phi((x - 2\pi j)|k|^{\tau} \kappa^{-1})=0$ for $j\neq [x]_\pi$, where $[x]_\pi/2\pi\in \Z$ is the unique integer such that $-\pi\le x-[x]_\pi<\pi$. Then $\phi_k(x)= \phi(\{x\}|k|^{\tau} \kappa^{-1})$, where $\{x\}= x-[x]_\pi$. 
Fix $k\neq 0$ in $\Z^{n-1}$ and consider the smooth  function
\[
\omega \to z_k(\omega) = 1-e^{ i\langle \omega,k\rangle}  
+\ \frac{1}{3} \kappa (1+|k|)^{-\tau}\,  
\phi_k\left( \langle \omega, k\rangle\right)\,. 
\]
\begin{Lemma}\label{lemma:z}
We have 
 \begin{equation}\label{eq:estimate-z}
|z_k(\omega)|\, \ge \frac{1}{3}\kappa  (1+|k|)^{-\tau} \quad \forall\, \omega\in \Omega. 
\end{equation}
Moreover, 
\begin{equation}\label{eq:z}
z_k(\omega)=  1-e^{i\langle \omega,k\rangle}\quad \forall\, \omega\in\Omega_\kappa . 
\end{equation}
\end{Lemma}
{\em Proof}. 
Let $\Omega^1$ be  the set of all $\omega\in\Omega$ such that 
\[
\pi/6\, \le\,   |\{\langle \omega, k\rangle\}|\, |k|^{\tau} \kappa^{-1}\,  \le\,   \pi
\] 
and  $\Omega^2$ the set of all $\omega\in\Omega$ such that  
\[
|\{\langle \omega, k\rangle\}|\, |k|^{\tau} \kappa^{-1}\,  \le\,  \pi/6. 
\]
For every $\omega\in\Omega^1$ we have 
\[
|1 -\exp ( i \langle
k,\omega\rangle )|= 
2|\sin( \frac{1}{2}\{\langle \omega, k\rangle\})|\ge \frac{4}{\pi}|\{\langle \omega, k\rangle\}| \ge  \frac{2}{3}\kappa  (1+|k|)^{-\tau} .
\]
This implies 
\[
|z_k(\omega)|\ge \frac{1}{3}\kappa  (1+|k|)^{-\tau}\quad \forall\, \omega\in\Omega^1.
\]
If $\omega\in\Omega^2$, then $\phi_k\left( \langle \omega, k\rangle\right)= \phi\left( \{\langle \omega, k\rangle\}|k|^{\tau} \kappa^{-1}\right)= 1$, hence,   
\[
{\rm Re\, }(z_k(\omega)) \ge \frac{1}{3}\kappa  (1+|k|)^{-\tau}
\]
which proves \eqref{eq:estimate-z}. 
Moreover, for any $\omega\in\Omega_\kappa$ we have  
$\phi_k\left( \langle \omega, k\rangle\right)= \phi\left(\{\langle \omega, k\rangle\} |k|^{\tau} \kappa^{-1}\right)=0$ in view of  \eqref{eq:sdc}, which implies \eqref{eq:z}. \finishproof

Let us go back to the homological equation \eqref{eq:homological1}. For every $k\neq 0$ we set 
\[
\hat f_k(\omega,t) := - \frac{\hat F_k(\omega,t)}{z_k(\omega)}.
\]
Using Lemma \ref{lemma:z} we obtain that the function 
\[
(\varphi,I)\to f(\varphi,\omega,t):= \sum_{k\in\Z^{n-1}} \hat f_k(\omega,t) e^{i\langle\varphi,k\rangle}
\]
belongs to $C^\infty(\T^{n-1}\times \Omega)$ for any $t$ fixed and the map $J\ni t\to f(\cdot,\cdot,t)\in C^\infty(\T^{n-1}\times \Omega)$ is $C^1$. Hence, the map $J\ni t\to S(\cdot,\cdot,t)\in C^\infty(\T^{n-1}\times \Omega)$ is $C^1$ as well. 
Moreover,  $F_k(\varphi,\omega,t)=0$ for each $\omega\in\Omega_\kappa^0$ and $k\in \Z^{n-1}$ and using \eqref{eq:z} we obtain that   that $S(\varphi,\omega,t)=0 $ for $\omega\in\Omega_\kappa^0$. Now Lemma \ref{Lemma:flat} implies that the function 
$\omega \to S(\varphi,\omega,t)$ is flat at $\Omega_\kappa^0$ for each $\varphi$ and $t$ fixed. Now using 
This completes the proof of  Lemma \ref{Lemma:homological}.
\finishproof

Using Lemma \ref{Lemma:homological} we find $a_t^j$ and  $p_t^j(I)$ such that $\partial_I^\alpha S_t^j(\varphi,I,t)=0$ for every $I\in E_t^\kappa$ and $\alpha\in \Z^{n-1}$. Using Lemma \ref{prop:realization} we find  $C^1$ of realisations $S_t(\varphi, I,\lambda)$ and $p_t(I,\lambda)$ of the formal symbols and  such that 
which completes the proof of Proposition \ref{prop:commutator}. 
\finishproof

We are looking for $C^1$ families of 
solutions  $t\to (\lambda(t), v(t))$  of the equation \eqref{eq:quasimode-equation} of the form 
\[
v(t)=E_t(\lambda)^{-1}T_t(\lambda)A_t(\lambda)e(t)
\] 
for  $t\in J$.  
In view of  Proposition  \ref{prop:commutator}, $e(t)$ should satisfy the equation 
\begin{equation}
\label{eq:equation1}
e^{i\pi\vartheta/2}W_t^0(\lambda)e(t) + e^{i\pi\vartheta/2}S_t(\lambda)e(t)  = e(t) + O_{N}(|\lambda|^{-N})e(t). 
\end{equation} 
Natural candidates for $e(t)$ are the sections $e_k$, $k\in \Z^{n-1}$. 
Since $\lambda\in {\mathcal D}$ may be complex, 
we consider  almost analytic extensions of  order $M \ge 2N+n+2$ 
of the  functions  $\phi_t^0$, $p_t$ and $S_t^j$ in $\zeta =\xi + i\eta$, where $\xi\in D$ and $|\eta| \le C$. The almost analytic extension of $\phi_t^0$ is 
given  by \[
\phi_t^0(x,\xi + i\eta)=L_t(\xi + i\eta)+ R_t(x,\xi + i\eta)
\] 
where
\[
L_t(\xi + i\eta) = \sum_{|\alpha|\le M} \partial_\xi^\alpha L_t(\xi)(i\eta)^\alpha (\alpha !)^{-1}\ \mbox{and}\
R_t(x,\xi + i\eta) = \sum_{|\alpha|\le M} \partial_\xi^\alpha
R_t(x,\xi)(i\eta)^\alpha (\alpha !)^{-1}\,. 
\]
It is easy to see that 
\[
 \partial_\zeta \phi_t^0(x,\xi + i\eta) = O(|\eta|^{M}) .
\]
 Moreover, 
\begin{equation}
\partial_\zeta^\alpha \bar \partial_\zeta^\beta R_t(x,\xi + i\eta) = O_M\left(|\xi - E_t^\kappa|^{M} \right)\, ,\quad |\eta|\le C ,
\label{eq:almost-analytic}
\end{equation}
for $\alpha,\, \beta \in \N^{n-1}$ 
since $R_t$ is flat at $\R^{n-1}\times E_t^\kappa$. 
In the same way we obtain almost analytic extensions $S_t(\varphi, \zeta, \lambda)$  of $S_t^j$ and  $p_t^j(\zeta)$ of $p_t^j$,  $\zeta:=\xi + i\eta$,  such that
\begin{equation}
 \overline \partial_\zeta p_t^j(\xi + i\eta) = O(|\eta|^{M}),\quad 
  \partial_\xi^\alpha\partial_\eta^\beta \overline \partial_\zeta S_t^j(x,\xi + i\eta) = O_{\alpha,\beta}(|\eta|^{M-|\beta|}), 
\label{eq:almost-analytic-p-S}
\end{equation}
for $\alpha,\, \beta\in \N^{n-1}$, $|\beta|\le M$. Moreover,  
\begin{equation}
\partial_\xi^\alpha  \partial_\eta^\beta S_t^j(\varphi,\xi + i\eta) = O\left(|\xi - E_t^\kappa|^{N} \right)\, ,
\label{eq:almost-analytic2}
\end{equation}
for $|\eta|\le C$,  and  
 ${\rm supp\, }_\zeta S_t^j \subset K$, ${\rm supp\, } p_t^j \subset K$ for $j\in\N$, where $K$ is a fixed compact subset of $\R^{n-1}$. We have as well $p_t^0(\xi + i\eta)=1$ whenever $\xi\in D^0$. 
\begin{Prop}  \label{prop:spectral-decomposition}
For each  $t\in J$ we have 
\begin{equation}
\label{eq:newoperator}
\begin{array}{rcll}
W_t^0(\lambda)e_k(\varphi)\, &=&\,  
\exp\Big(-i\lambda \phi_t^0\big(\varphi,(k + \vartheta_0/4)\lambda^{-1}\big)\Big) \\[0.3cm]
&\times&\, \displaystyle \Big(\sum_{j=0}^N \, p_t^j\big((k + \vartheta_0/4)\lambda^{-1}\big)\, \lambda^{-j}\Big)\,   e_k(\varphi)\,  + \, 
O_{N}(|\lambda|^{-N-1}) e_k(\varphi) 
\end{array}
\end{equation}
and 
\begin{equation}
S_t(\lambda)e_k(\varphi)\,  = \, 
O_{N}\left(|\lambda|^{-N-1}+ |E_t^\kappa - (k+\vartheta/4)\lambda^{-1} |^{N+1}\right) e_k(\varphi)
                                  \label{eq:newrest}
\end{equation}
where   $\lambda \in {\cal D}$   and  $k\in \Z^{n-1}$.  
\end{Prop}
{\em Proof}. 
The proof of the proposition is close to that of Proposition  3.11, \cite{PT4} but we give it for the sake of completeness. We have
\[
\begin{array}{lcrr}
\widetilde{W^0_t(\lambda)e_k}(x)\ =\  \widetilde{e_k}(x)\, 
e^{-i\lambda \phi_t^0(x,\xi_k)} \\ [0.3cm]
\displaystyle \times \, 
\lambda^{-j}\, \left(\frac{\lambda}{2\pi}\right)^{n-1}\, \sum_{j=0}^N\, \int_{{\R}^{2n-2}} 
\, e^{i\lambda\langle x-y +  
w_t(x,\xi_k, \eta_k),\eta_k \rangle} \, 
 \, p_t^j(I)\, dI \,  dy \, +\,  \, O_{N}\left(|\lambda|^{-N-1}\right) \widetilde e_k(x)\, ,
\end{array}
\]
where $\lambda\in {\mathcal D}$ and 
\[
w_t(x,\xi,\eta)=\int_0^1\nabla_\xi\phi_t^0(x,\xi+\tau \eta)d\tau,\quad 
\xi_k =  (k + \vartheta_0/4)/\lambda, \quad \eta_k=I-(k + \vartheta_0/4)/\lambda.
\]
If $|k| \ge C_0|\lambda|$ and $C_0\gg 1$ ($C_0$ depends only on the compact set $K\subset \R^{n-1}$ such that ${\rm supp\, } p_t^j \subset K$ for every $j\in \N$) then $|\eta_k|\ge 1$ and we can integrate by parts with respect to $y$ gaining $O_{N}(|\lambda|^{-N-1})$. Suppose now that 
 $|k| \le C_0|\lambda|$.   We have 
\begin{equation}\label{eq:estimate-imaginary-part}
\big|{\rm Im\, }\big((k + \vartheta_0/4)/\lambda\big)\big| \le \frac{C}{|\lambda|} \quad \mbox{for $|k|\le C_0 |\lambda|$ and $\lambda\in {\mathcal D}$. }
\end{equation}
Then deforming the contour of integration  we obtain
\[
\begin{array}{lcrr}
W_t^0(\lambda)e_k(\varphi)\ =\  e_k(\varphi)\, 
e^{-i\lambda \phi_t^0(\varphi,(k + \vartheta_0/4)/\lambda)} \\ [0.3cm]
\displaystyle \times \, 
  \sum_{j=0}^N\,  \lambda^{-j}\, \left(\frac{\lambda}{2\pi}\right)^{n-1}\, \int_{{\R}^{2n-2}} 
\, e^{-i\lambda\langle u,v\rangle } \, 
p_t^j(v + (k + \vartheta_0/4)/\lambda)\, du\,  dv \,  + \, 
O_N(|\lambda|^{-N-1}) e_k(\varphi)\,  ,
\end{array}
\]
which implies  (\ref{eq:newoperator}).

To prove \eqref{eq:newrest} we write $\widetilde{S_t(\lambda)e_k(x)}$ as an oscillatory integral as above, and then for $|k| \le C_0|\lambda|$ we change the contour of integration with respect to $y$ 
by 
$$y \to v = y-x-w_t\big(x,(k + \vartheta_0/4)/\lambda, I - (k + \vartheta_0/4)/\lambda\big)\, 
$$
while for $|k| \ge  C_0|\lambda|$ we integrate by parts to gain $O_N(|\lambda|^{-N-1})$. 
This implies, using \eqref{eq:estimate-imaginary-part}, that
\[
\begin{array}{lcrr}
S_t(\lambda)e_k(\varphi)\ =\  e_k(\varphi)\, 
e^{-i\lambda \phi_t^0(\varphi,(k + \vartheta_0/4)/\lambda)}  \\ [0.3cm]
\displaystyle \times \,   \sum_{j=0}^N\,  \left(\frac{\lambda}{2\pi}\right)^{n-1}\, 
 \int_{{\R}^{2n-2}} 
\, e^{-i\lambda\langle v, I -(k + \vartheta_0/4)/\lambda\rangle }\, 
S_t^j(\varphi, I)\lambda^{-j}\, dI \, dv \, + \, O_N(|\lambda|^{-N-1})e_k(\varphi).
\end{array}
\] 
Since $M>2N+n+2$, taking the Taylor expansion of order $N$  of the function
\[
[0,1]\ni s \to \psi(s):=S_t^j\big(\varphi, (k + \vartheta_0/4)/\lambda) +s(I-(k + \vartheta_0/4)/\lambda)\big)
\]
 at $s=0$  with an integral reminder and using \eqref{eq:almost-analytic-p-S} and \eqref{eq:estimate-imaginary-part} we get
\[
S_t^j(\varphi, I) = \sum_{|\alpha\le N|} \partial^\alpha_{\zeta} S_t^j(\varphi, (k + \vartheta_0/4)/\lambda) (I-(k + \vartheta_0/4)/\lambda)^\alpha /\alpha ! \, 
+\,  T_N (\varphi, I)\, +\,  O(|\lambda|^{-N-n-1})
\]
where the reminder term is 
\[
 T_N (\varphi, I):=(N+1) \sum_{|\alpha|= N+1} \int_0^1(1-s)^N\partial^\alpha_I S_t^j\big(\varphi, I+s(k + \vartheta_0/4)/\lambda\big)\big(I-(k + \vartheta_0/4)/\lambda\big)^\alpha/\alpha !\, ds . 
\]
We have 
\[
\partial^\alpha_\zeta  S_t^j(\varphi, (k + \vartheta_0/4)/\lambda) = O_{N,\alpha,\beta}\left(|E_t^\kappa - (k+\vartheta/4)\lambda^{-1} |^{N+1} \right),\quad \lambda \in {\mathcal D}, 
\]
for every $N\in \N$ and $\alpha\in \N^{n-1}$ in view of \eqref{eq:almost-analytic2}. 

To estmate the reminder
  we integrate $N+1$ times by parts with respect to $v$  in the corresponding oscillatory integral with amplitude 
\[
(N+1) \sum_{|\alpha|= N+1} \int_0^1(1-s)^N\partial^\alpha_I S_t^j\big(\varphi, I+s(k + \vartheta_0/4)/\lambda\big)\big(I-(k + \vartheta_0/4)/\lambda\big)^\alpha/\alpha !\, ds
\]
and we estimate it by $C_N |\lambda|^{-N-1}$. 
This implies   (\ref{eq:newrest}). \finishproof

Proposition \ref{prop:spectral-decomposition} suggests that we should look for pairs $(\lambda, k)\in {\mathcal D}\times \Z^{n-1}$  such that $|\lambda| \gg 1$ and 
\begin{equation}\label{eq:quantization-condition}
|E_t^\kappa - (k+\vartheta/4)\lambda^{-1} | \le \frac{C}{|\lambda|}
\end{equation}
where $C>0$ is a constant. Then \eqref{eq:newoperator} and \eqref{eq:newrest} imply
\[
\begin{array}{rcll}
W_t^0(\lambda)e_k(\varphi)\, &=&\,  
\exp\Big(-i\lambda L_t\big((k + \vartheta_0/4)\lambda^{-1}\big)\Big) \\[0.3cm]
&\times&\, \displaystyle \Big(\sum_{j=0}^N \, p_t^j\big((k + \vartheta_0/4)\lambda^{-1}\big)\, \lambda^{-j}\Big)\,   e_k(\varphi)\,  + \, 
O_{N}(|\lambda|^{-N-1}) e_k(\varphi) 
\end{array}
\]
and 
\[
S_t(\lambda)e_k(\varphi)\,  = \, 
O_{N}\left(|\lambda|^{-N-1}\right) e_k(\varphi)
\]
Thus taking $e=e_k$ in \eqref{eq:equation1} we obtain 
\[
\begin{array}{lcrr}
\displaystyle \exp\Big(-i\lambda L_t\big((k + \vartheta_0/4)\lambda^{-1}\big)+i \pi\vartheta/2\Big) \,  \Big(\sum_{j=0}^N \, p_t^j\big((k + \vartheta_0/4)\lambda^{-1}\big)\, \lambda^{-j}\Big)\,   e_k(\varphi)\\[0.3cm]  
\displaystyle = 
O_{N}(|\lambda|^{-N-1}) e_k(\varphi) 
\end{array}
\]
for every $N\in \N$. Recall that $p_t^0(\xi+i\eta)=1$ if $\xi\in D^0$. Then for $|\lambda|\gg 1$ and $t\in J$ we can write the above equation  as follows
\begin{equation} \label{eq:new-equation}
\begin{array}{rcll}
\displaystyle\lambda \, L_t\Big(\frac{k + \vartheta_0/4}{\lambda}\Big)\,  &=& \, 2\pi k_n\pi  \, +\,  \pi \vartheta/2  \\[0.3cm]  
 &+&\,   \displaystyle\frac{1}{i}\, {\rm Log\, } \Big(1+ \sum_{j=1}^N \, p_t^j\Big(\frac{k + \vartheta_0/4}{\lambda}\Big)\, \lambda^{-j}\Big)\, +\, 
O_{N}(|\lambda|^{-N-1})
\end{array}
\end{equation}
where $k_n\in\Z$ and $ {\rm Log}\, z = \ln |z| + i\,  {\rm arg}\, z,\ -\pi <{\rm arg}\,
z < \pi$.

Hence, to construct quasi-modes we have to find pairs $(\lambda,k)$ satisfying both \eqref{eq:quantization-condition} and \eqref{eq:new-equation}.

\section{$C^1$ families of  quasi-modes and iso-spectral invariants} \label{Sec:Quasi}
Given $t\in J$ and $\omega\in \Omega_\kappa^0$ the formulas \eqref{eq:quantization-condition} and \eqref{eq:new-equation} suggest that the quantization condition of the Lagrangian torus $\Lambda_t(\omega)$ should be of the form
\[
\Big|\lambda I_t(\omega)-(k + \vartheta_0/4)\Big| + \Big|\lambda \, L_t\Big(\frac{k + \vartheta_0/4}{\lambda}\Big) - 2\pi k_n\pi - \pi \vartheta/2 \Big| \le C
\]
for some $C>0$, where $I_t(\omega)\in E_t^\kappa$ is the corresponding action on the torus $\Lambda_t(\omega)$,  $(k,k_n)\in \Z^n$ and $\lambda \in {\mathcal D}$. To obtain iso-spectral invariants from $C^1$-families of quasi-modes we need a stronger quantization condition which will be formulated below. 
\subsection{Quantization condition}\label{Subsec:Quantization}
Fix $t\in J$. The quantization condition corresponding to  a Lagrangian torus $\Lambda_t(\omega)$ with a frequency $\omega\in \Omega_\kappa^0$ will be given by means the following Lemma. 
\begin{Lemma} \label{Lemma:quantization}
Given $t\in J$
there is a set $\Xi_\kappa^t\subset \Omega_\kappa^0$ of full Lebesgue measure in  $\Omega_\kappa^0$ such that the following holds.\\

\noindent
For any $\omega\in \Xi_\kappa^t$ there is an infinite sequence $\widetilde {\mathcal M}(\omega)$  of $(q,\lambda)\in \Z^n \times [1,\infty)$ such that
\begin{equation}
q=(k,k_n)\in \Z^{n-1}\times \Z \, , \quad \lambda=\mu_q^0\ge 1\quad  \mbox{satisfies} \quad  c_0^{-1}|q| \le \mu_q^0 \le c_0 |q|\ \mbox{with}\  c_0>0 ,
                                                     \label{eq:q and mu}
\end{equation} 
and
\begin{equation}
\lim_{|q|\to \infty}\, \Big|\mu_q^0 \, \Big(I_t(\omega), L_t(I_t(\omega))\Big)\  -\ \Big(k+\frac{\vartheta_0}{4}, 2\pi \Big(k_n+\frac{\vartheta}{4}\Big)\Big)\Big|\,  =\,  0.
                                                     \label{eq:quantization}
\end{equation} 
\end{Lemma}
{\em Proof}. Denote by $\Xi_\kappa^t$ the set of all $\omega\in \Omega_\kappa^0$ such that
\begin{equation}
 \quad  2\pi k_nI_t(\omega) \neq  L_t(I_t(\omega)) k \quad  \mbox{for each} \ 0\neq (k,k_n)\in \Z^{n-1}\times \Z .
                                                     \label{eq:non-resonant}
\end{equation} 
We claim that the complement $\Omega_\kappa^0\setminus\Xi_\kappa^t$ of $\Xi_\kappa^t$ in $\Omega_\kappa^0$  is of Lebesgue measure zero. Suppose the contrary. Then there is $ 0\neq (k,k_n)\in \Z^{n-1}\times \Z$ and a set of positive Lebesgue measure $R_t\subset \Omega_\kappa^0$ such that 
\[
 2\pi k_nI_t(\omega) \, = \,  L_t(I_t(\omega))k  \quad  \forall\, \omega\in R_t.
\]
On the other hand, the map $\Omega\ni \omega\to I_t(\omega)= \nabla L_t^\ast(\omega)\in D$ is a local diffeomorphism with inverse $I\to \nabla L_t(I)$ by {\it 4}, Theorem \ref{Theo:soft-BNF}, hence, the set $R_t^0:=\{I_t(\omega):\ \omega\in R_t\}$ is of  positive Lebesgue measure in $\R^{n-1}$. Moreover, 
\begin{equation}\label{eq:k-to-L}
  2\pi k_nI =  L_t(I) k   \quad  \forall\, I\in R_t^0
\end{equation}
and $R_t:=\{\nabla L_t(I):\ I\in R_t^0\}$ by definition. 
Let $I^0\in R_t^0$ be a point of  positive Lebesgue density in $R_t^0$. Set $\omega^0=(\omega^0_1,\ldots,\omega^0_{n-1}):=\nabla L_t(I^0)\in \Omega_\kappa^0$.  Differentiating  \eqref{eq:k-to-L} with respect to $I$ at $I^0$ and using Lemma \ref{Lemma:flat} we get
$2\pi k_n= k_j\omega^0_j$, for $j= 1,\ldots, n-1$, which contradicts  \eqref{eq:sdc}. Hence, the Lebesgue measure of $ \Omega_\kappa^0\setminus \Xi_\kappa^t$ is zero. 
On the other hand, \eqref{eq:non-resonant} implies that for any  $\omega\in\Omega_\kappa^0$ the trajectory
\[
\{\lambda(I_t(\omega), L_t(I_t(\omega)))\ ({\rm mod}\, \Z^n):\,  \lambda\ge 1\, \}\, \subset \, \R^n/\Z^n
\]
is not periodic, hence, it is dense on the torus $\R^n/\Z^n$ which implies that there exists  an infinite sequence $(q_j,\lambda_j)_{j\in \N}$ satisfying \eqref{eq:quantization}. The inequality in \eqref{eq:q and mu} follows from \eqref{eq:quantization} since 
the continuous function  
\[
\omega\mapsto \|(I_t(\omega), L_t(I_t(\omega)))\|=\|(\nabla L_t^\ast(\omega), L_t(\nabla L_t^\ast(\omega)))\|
\]
does not vanish on the compact set $\Omega_\kappa^0$ 
 in view of  \eqref{eq:beta function} and \eqref{eq:beta-alpha}.  \finishproof

\noindent
We point out that the set $\Xi_\kappa^t$ and the sequence $\widetilde {\mathcal M}(\omega)$ may  depend on $t$.\\

\noindent
From now on we fix  $\omega$ in the set $\Xi_\kappa^t$ given by Lemma \ref{Lemma:quantization} and 
denote by  $\mathcal M\subset \Z^n$ the image of the projection of $\widetilde {\mathcal M}(\omega)\subset \Z^n \times  [1,\infty)$ on the first factor. $\mathcal M$ will be the index set of the $C^1$ family of quasi-modes that we are going to construct and  \eqref{eq:quantization} - the quantization condition for $s=t$. To obtain a quantization condition for the tori $\Lambda_s(\omega)$ for $s$  close to $t$ we consider 
for any $q\in {\mathcal M}$ the interval 
\[
J_q:=\left[t,t+2|q|^{-1}  \right] .
\] 
Getting rid of finitely many elements  $q\in {\mathcal M}$ we suppose that  $J_q$ is contained in $J\subset [0,\delta]$ for every $q\in {\mathcal M}$. 
Recall from Theorem \ref{Theo:soft-BNF} that 
the maps  $s\mapsto L_s\in C^\infty(D)$ and 
$s\mapsto I_s\in C^\infty(\Omega;\R^{n-1})$ are $C^1$ on the interval $J$. Then using 
 \eqref{eq:q and mu} and \eqref{eq:quantization} we obtain that there exists a constant $C=C(\omega)>0$ independent of  $q\in {\mathcal M}$ and  $s\in J_q$ such that 
\begin{equation}
 \Big| \mu_q^0 \Big(I_s(\omega), L_s(I_s(\omega))\Big)\  -\ \Big(k+\frac{\vartheta_0}{4}, 2\pi \Big(k_n+\frac{\vartheta}{4}\Big)\Big)\Big|\, \le\,  C \quad \forall\, q\in {\mathcal M},\  s\in J_q. 
                                                     \label{eq:quantization1}
\end{equation} 
The quantization condition \eqref{eq:quantization1} will be used below to construct a $C^1$ quasi-mode with an index set ${\mathcal M}$ for $s\in J_q$, $q\in {\mathcal M}$. 
Condition  \eqref{eq:quantization} is not needed for the the construction of the quasi-mode, but it is essential   for the proof of Lemma \ref{Lemma:main} below. 

\subsection{Construction of $C^1$ families of quasi-modes}\label{sec:construction-quasi-modes}
Fix a positive  integer $M\ge 0$.  For any $q\in {\cal M}$ with $|q|\ge q_0\gg 1$ we are going to construct a family of quasi-modes of order $M$ depending on $s\in J_q$ such that the corresponding family of quasi-eigenvalues $s\mapsto \mu_q(s)^2$ belongs to $C^1(J_q)$. 
\begin{Theorem}\label{Theorem:quasimodes}
For every $q=(k,k_n)\in {\mathcal M}$ and $s\in J_q$ there exists a  quasi-mode $(\mu_q(s)^2,u_{s,q})$ of $\Delta_s$ of order $M$ such that 
\begin{enumerate}
\item[(i)] $ u_{s,q}\in  D(\Delta_s)$ and $\|u_{s,q}\|_{ L^2(X)}=1$;
\item[(ii)] There exists a constant $C_M>0$ such that
\begin{equation}  
\left\{  
\begin{array}{lcr}  
\left\|\Delta\,  u_{s,q}\ - \ \mu_q^2(s)\,  u_{s,q}\right\|\ \le \   
C_M\,  \mu_q^{-M}(s)  \, \quad \mbox{in}\ L^2(X)\, , \\ [0.3cm] 
\displaystyle {\mathcal B}\,  u_{s,q}|_\Gamma \   =\ 0 \,   
\end{array}  
\right.  
                                         \label{eq:thequasimode}
\end{equation} 
for every $q\in {\cal M}$ and  $s\in J_q$;   
\item[(iii)] We have 
\[
\mu_q(s) = \mu_q^0 +c_{q,0}(s) + c_{q,1}(s) \frac{1}{\mu_q^0}  + \cdots + 
c_{q,M}(s) \frac{1}{(\mu_q^0)^M}\, ,\quad  \mbox{where}
\]   
\item[(iv)]
The functions $s\mapsto c_{q,j}(s)$ are real valued and $C^1$ on the interval $J_q$; 
\item[(v)] There exists  a constant $C_M'>0$ such that
$|c_{q,j}(s)| \le C_M'$ for every  $q\in {\mathcal M}$, $0\le j\le M$,  and any $s\in J_q$;
\item[(vi)] There exists $C>0$ such that 
\[
\Big|\mu_q(s) L_s\Big(\frac{k + \vartheta_0/4}{\mu_q(s)}\Big) -  2\pi \Big(k_n +
\frac{\vartheta}{4}\Big)\Big| \, \le\,  \frac{C}{\mu_q(s)}
\]
for every $q\in  {\cal M}$ and $s\in J_q$;
\item[(vii)]  We have
\[
\frac{k + \vartheta_0/4}{\mu_q(t)}\, =\, I_t(\omega) + o\Big(\frac{1}{|q|}\Big)\quad \mbox{as}\ |q|\to \infty. 
\]
\end{enumerate}
\end{Theorem}
{\em Proof}. We are looking for a  perturbation
$\lambda=\mu_q(s)$ of  $\mu^0_q$ satisfying \eqref{eq:new-equation} which means that
\[
\mu_q(s) L_s\Big(\frac{k + \vartheta_0/4}{\mu_q(s)}\Big)\,  -\,  \frac{1}{i}\, {\rm Log}\, \Big(1 +\sum_{j=1}^{M}\, p_s\left(\frac{k + \vartheta_0/4}{\mu_q(s)}\right) \mu_q(s)^{-j}\Big)
\]
\[
\
= \   2\pi\Big( k_n  + \frac{\vartheta}{4}\Big) +
O_{M}\Big(\frac{1}{(\mu_q^0)^{M+1}}\Big)
\]    
uniformly with respect to   $q\in {\mathcal M}$ and   $s\in J_q$.   
Introducing a small parameter $\varepsilon_q =  (\mu_q^0)^{-1}$ we 
are looking for 
\begin{equation}
\label{eq:mu_q}
\left\{
\begin{array}{rcll}
\mu_q(s) \, &=&\, \mu_q^0 + c_{q, 0}(s) + c_{q, 1}(s)\varepsilon_q +
\cdots  c_{q, M}(s)\varepsilon_q^{M}\ , \\[0.3cm] 
\zeta_q (s)\, &=&\, I^0_s(\omega)  +  b_{q, 0}(s)\varepsilon_q+
\cdots  b_{q, M}(s) \varepsilon_q^{M+1}(s) + b_{q, M+1}(s) \varepsilon_q^{M+2}
\end{array}
\right.
\end{equation}
such that
\begin{equation}
\label{eq:system}
\left\{
\begin{array}{rcll}
\mu_q(s)\zeta_q(s) &=& \displaystyle k + \frac{\vartheta_0}{4}  \\[0.3cm] 
\displaystyle \mu_q(s) L_s(\zeta_q(s)) &=&  \displaystyle    2\pi\Big( k_n  + \frac{\vartheta}{4}\Big) + \frac{1}{i}\, {\rm Log}\, \Big(1 + \sum_{j=1}^{M}p^j_s(\zeta_q(s))\mu_q(s)^{-j}\Big) +  
O_{M}(\varepsilon_q^{M+1})\, .
\end{array}
\right.
\end{equation}

We are going to find $\mu_q(s)$. 
Using \eqref{eq:mu_q} we write 
\[
\begin{array}{clr}
\displaystyle \mu_q(s)\zeta_q(s) - k - \vartheta_0/4 &=&  \displaystyle \sum_{j=0}^M \varepsilon_q^j \left[ b_{q, j }(s)+ c_{q, j}(s) I_s(\omega) - W_{q,j}(s)  \right]  \\[0.3cm]
\displaystyle &+& \varepsilon_q^{M+1}\left[(\varepsilon_q \mu_q) b_{q, M+1}(s)- W_{q,M+1}(s)  \right], 
\end{array}
\]
where 
\begin{equation}
\label{eq:W}
\left\{
\begin{array}{lcrr}
\displaystyle W_{q,0}(s) = k + \vartheta_0/4 - \mu_q^0 I_s(\omega)\,  ,\\[0.3cm]
\displaystyle W_{q,j}(s) = -\sum_{r+s=j-1}c_{q, r}(s)b_{q, s }(s)\ \mbox{for}\  1\le j \le M,\ \mbox{and}\\ [0.3cm]
\displaystyle  W_{q,M+1}(s) = -\sum_{l=0}^M\sum_{r= M-l}^Mc_{q, r}(s)b_{q, l }(s).
\end{array}
\right.
\end{equation}
Expanding $ L_s(\zeta_q(s))$ and $p^j_s(\zeta_q(s))$, $1\le j\le M$, in Taylor series at $\zeta=I_s(\omega)$ up to order $M$
we obtain from \eqref{eq:system} 
the following linear  systems
\[
\left\{
\begin{array}{rcll}
b_{q, j}(s) + c_{q, j}(s)I_s(\omega)    &=&\
W_{q,j}(s) \\[0.3cm]
\displaystyle L_s( I_s(\omega))c_{q, j}(s) +   \langle \omega,b_{q, j}(s)\rangle  &=& \displaystyle   
V_{q,j}(s)\, ,
\end{array} 
\right.
\]
for $0 \le j\le M$, and we put
$b_{q, M+1}(s)= (\varepsilon_q \mu_q)^{-1} W_{q,M+1}(s)$, 
where $W_{q,j}(s)$ is given by \eqref{eq:W}, 
and $V_{q,j}(s)$ is a  polynomial of $c_{q, r}(s)$ and  $b_{q,r'}(s)$ with  $0 \le r,r' \le j-1$ and  with $C^1$ with respect to $s$ coefficients. 
By \eqref{eq:beta function} and \eqref{eq:beta-alpha} the corresponding determinant is
$$
D(I_s(\omega)):= L_s(I_s(\omega)) -   \langle I_s(\omega), \omega\rangle =  -\beta_s(\omega)=2 
\int_{\Lambda_s(\omega)} A_s( \rho)d\mu_s\,  > 0
$$
and we obtain a unique
solution $(c_{q, j}(s),b_{q ,j}(s))$, $0\le j \le M-1$.  
More precisely, 
\begin{equation}
\label{eq:quasi-coefficients}
\left\{
\begin{array}{lcrr}
 c_{q, j}(s) = D(I_s(\omega))^{-1}\left[ V_{q,j}(s) - 
 2\pi\langle \omega,W_{q, j}(s)\rangle\right]\\[0.3cm]
b_{q, j}(s) = W_{q,j}(s) -c_{q, j}(s) I_s(\omega) .
\end{array}
\right.
\end{equation} 
We have  
\begin{equation}
\label{eq:first-coefficients}
\left\{
\begin{array}{lcrr}
W_{q,0}(s)=k + \vartheta_0/4 - \mu_q^0 I_s(\omega) = O(1), \\[0.3cm]
V_{q,0}(s) = 2\pi k_n -\pi \vartheta/2- \mu_q^0 L_s(I_s(\omega))= O(1)\, ,\ q\in \cal M\, ,
\end{array}
\right.
\end{equation}
uniformly with respect to $q\in {\cal M}$ and  $s\in J_q$, 
in view  of (\ref{eq:quantization1}). 
Hence, $b_{q,0}(s)$ and $c_{q,0}(s)$, $q\in {\mathcal M}$,  are $C^1$ in $J_q$ and  uniformly bounded. 
 By recurrence we prove that $b_{q, j}(s)$ and $c_{q, j}(s)$, $q\in {\mathcal M}$,  are
$C^1$ in $J_q$ and uniformly bounded with respect to $q\in {\mathcal M}$ and $s\in J_q$. To evaluate $b_{q, M+1}(s)$ observe that $\varepsilon_q \mu_q = 1 + O(\varepsilon_q)$. 

For such $\mu_q(s)$ the quantization condition \eqref{eq:quantization1} gives the estimate
\[
\frac{k + \vartheta_0/4}{\mu_q(s)} = \zeta_q(s)= I_s(\omega) + O\left(\mu_q(s)^{-1}\right)
\]
uniformly with respect to $s\in J_q$.  Then Proposition \ref{prop:spectral-decomposition} for $N=M$ implies   that $(\lambda=\mu_q(s),e_k)$ satisfy \eqref{eq:equation1}  
and we obtain that
\[
\left\{  
\begin{array}{lcr}  
\left\|\Delta\,  u_{s,q}\ - \ \mu_q^2(s)\,  u_{s,q}\right\|\ \le \   
C_M\,  \mu_q^{-M}(s)  \, \quad \mbox{in}\ L^2(X)\, , \\ [0.3cm] 
\displaystyle \left\|{\mathcal B}\,  u_{s,q}\right\|\   \le \ C_M\,  \mu_q^{-M}(s)  \, \quad \mbox{in}\ L^2(\Gamma)\, .   
\end{array}  
\right.  
\]
In order to  prove  the property (i) and to satisfy the boundary conditions in (ii) exactly we follow the proof given  in   \cite{PT4}, Sect. 3.6.3, using Proposition \ref{Lemma:continuity}. The property (vi) follows from the second equation of \eqref{eq:system}. To prove (vii) observe that
\[
\left\{
\begin{array}{lcrr}
W_{q,0}(t)=k + \vartheta_0/4 - \mu_q^0 I_t(\omega) = o(1) \\[0.3cm]
V_{q,0}(t) = 2\pi k_n -\pi \vartheta/2- \mu_q^0 L_t(I_t(\omega))= o(1)
\end{array}
\right.
\]
as $|q|\to \infty$ in view of \eqref{eq:quantization}. Then \eqref{eq:quasi-coefficients} implies that $c_{q,0}(t)= o(1)$ and $b_{q,0}(t)= o(1)$ and we obtain
\[
\frac{k + \vartheta_0/4}{\mu_q(t)}\, =\, \zeta_t(\omega) \, =\, I_t(\omega) + b_{q,0}(t) \varepsilon_q+ O(\varepsilon_q^2)\, =\, I_t(\omega) + o\Big(\frac{1}{|q|}\Big)\quad \mbox{as}\ |q|\to \infty. 
\]
This completes the proof of the Theorem. 
\finishproof

\subsection{From quasi-modes to isospectral invariants}\label{sec:quasi-modes-to-isospectral}
We are going to complete the proof of Theorem \ref{Th:main1}. The items (i) and (ii) have been proven in Sect. \ref{Sec:BBM-BNF}. We are going to prove item (iii) which states that the functions  $\beta_t(\omega)$, $I_t(\omega)$ and  $\alpha_t(I_t(\omega))= L_t(I_t(\omega))$   
are independent of $t\in [0,\delta]$ for any $\omega\in\Xi$ provided that  the billiard tables satisfy the weak isospectral condition $(\mbox{H}_1)-(\mbox{H}_2)$. Recall that the set $\Xi$ is of the form \eqref{eq:Xi}, hence, it suffices to prove the statement for each $\omega$ in $\Omega_\kappa^0 \subset \Xi$. 

Given $\alpha\in\R$ we say that a family of functions $f_q:J_q\to \C$,  $q\in {\mathcal M}$, is $o\left(|q|^\alpha\right)$ as $q\to \infty $ uniformly with respect to $s$ in $J_q$ if 
\[
\displaystyle \lim_{q\to\infty}\, \left(|q|^{-\alpha} \sup_{s\in J_q}|f_q(s)|\right) = 0\, .
\]
We say that ``$f_q=O\left(|q|^\alpha\right)$ uniformly with respect to $s$ in $J_q$'' if there is $C>0$ such that $|q|^{-\alpha} |f_q(s)| \le C$ for any 
$q\in {\mathcal M}$ and $s\in J_q$.
The isospectral condition implies
\begin{Lemma}\label{Lemma:isospectral}
Suppose that $(\mbox{H}_1)-(\mbox{H}_2)$ holds. Fix an integer $M>2d\ge 0$. Then 
\[
\mu_q(s)-\mu_q(t)  = o(1)\quad \mbox{as}\ q\to \infty 
\]
 and
\[
\mu_q(s) = \mu_q(t)\left(1+ o\left(\frac{1}{|q|}\right)\right) \quad \mbox{as}\ q \to \infty
\]
uniformly with respect to $s\in J_q$. 
\end{Lemma}
{\em Proof}. 
It is easy to see that for any $q\in {\mathcal M}$ and $s\in J_q$, the distance from $\mu_q(s)^2$ to the spectrum of $\Delta_s$ can be estimated above by
\[
d_{s,q}:=\left|\, \mbox{Spec}\left(\Delta_s \right)\, -\,  
\mu_q(s)^2\right|\,  \le\,  C_M\, \mu_q(s)^{-M}\,   .
\]
Indeed, if $d_{s,q}\neq 0$  the spectral theorem and \eqref{eq:thequasimode} yield
\[
\frac{1}{d_{s,q}} \ge \|(\Delta_s - \mu_q(s)^2)^{-1}\|  \ge 
\|(\Delta_s - \mu_q(s)^2)u_{s,q}\|^{-1} \ge \frac{\mu_q(s)^M}{C_M}. 
\]
Then Theorem \ref{Theorem:quasimodes} and  \eqref{eq:q and mu} imply that  for any  $q\in {\mathcal M}$, $|q|\ge q_0 \gg 1$, and $s\in J_q$ 
there is $\lambda_{s,q}\in \mbox{Spec}\left(\Delta_{s}\right)$ such that $\lambda_{s,q} \ge \mu_q(s)^2/4 \ge (2c_0)^{-2}|q|^2$ and 
\begin{equation}
\left|\, \lambda_{s,q} \, -\,  
\mu_q(s)^2\right|\, \le\, 
C' \lambda_{s,q}^{-M/2}\,   
\label{eq:quasi-to-spectrum}
\end{equation}
where $C'=2^{M}C_M $. 
 Now  using (H$_2$) we get 
  for any  $q\in {\mathcal M}$ with $|q|\ge q_0\gg 1$ and $s\in J_q$    an integer $k=k(s,q)\ge 1$ such that 
\begin{equation}\label{eq:e:la_in}
\lambda_{s,q}\in [a_k,b_k].
\end{equation}
Fix $\gamma$ so that $M >2\gamma > 2d\ge 0$. Then choosing  $q_0$ sufficiently large we obtain from \eqref{eq:quasi-to-spectrum} and \eqref{eq:e:la_in}   that for any $q\in {\mathcal M}$ with $|q|\ge q_0$ and $s\in J_q$ the quasi-eigenvalue $\mu_q(s)^2$ belongs to  the interval 
\begin{equation}\label{eq:e:the_intervals}
 I_k:=\left[a_k-\frac{c}{2} a_k^{-\gamma},b_k+ \frac{c}{2} a_k^{-\gamma}\right]\, , 
\end{equation}
where $ k=k(s,q)$ and  $c>0$ is the constant of the third assumption of  $(\mbox{H}_1)$. In particular, 
\[
b_{k(q,s)}\ge  \mu_q(t)^2 - \frac{c}{2} a_{k(s,q)}^{-\gamma} \ge C_1|q|^2-C_2,
\]
for some positive constants $C_1$ and $C_2$,  which implies that $\lim k(s,q) = \infty$  as $q\to\infty$ {\em uniformly} with respect to $s\in J_q$. 
On the other hand, using the third assumption of $(\mbox{H}_1)$, the relation $b_k=a_k(1+o(1))$ as $k\to\infty$, which  follows from the first two assumptions in {(\rm H$_1$)}, 
and the inequality $\gamma>d$, we get 
\[
(a_{k+1}-\frac{c}{2}a_{k+1}^{-\gamma})-(b_k+\frac{c}{2}a_k^{-\gamma})=(a_{k+1}-b_k)-\frac{c}{2}a_{k+1}^{-\gamma}-\frac{c}{2}a_k^{-\gamma}\ge
c b_k^{-d}-c a_k^{-\gamma}>0
\]
for any $k\ge k_0$, where $k_0\gg 1$. This shows that the intervals  $I_k$ in \eqref{eq:e:the_intervals} do not intersect each other for  $k\ge k_0$.
Choose $q_0\gg 1$ so that $k(s,q)\ge k_0$ for any $q\in {\mathcal M}$ with $|q|\ge q_0$ and $s\in J_q$ (recall that $k(s,q)\to \infty$ as $|q|\to \infty$ uniformly with respect to $s\in J_q$). 
The function  $\mu_q(s)^2$ is continuous on $J_q$ (even $C^1$), hence,  it can not jump from one interval to another when $|q|\ge q_0$. 
Consequently,    $k(s,q)$ does not depend on $s$ for $|q|\ge q_0$. We have proved that for any $q\in {\mathcal M}$ such that $|q|\ge q_0$ there is $k=k(q)\in\N$ independent
of $s$ such that
\begin{equation}\label{eq:e:mu_in}
\mu_q(s)^2 \in \left[a_k-\frac{c}{2}a_k^{-\gamma},b_k+ \frac{c}{2} a_k^{-\gamma}\right]\quad \forall\, s\in J_q\, .
\end{equation}
Moreover,  $k(q)\to\infty$ as $q\to\infty$  and we obtain  
\[
\mu_q(s)^2 \ge a_{k(q)} - \frac{c}{2}a_k^{-\gamma} \ge \frac{1}{4}a_{k(q)} 
\]
for $|q|\ge q_0 \gg 1$ and  $s\in J_q$. 
Thus for $|q|\ge q_0 \gg 1$ we obtain 
$$
|\mu_q(s) - \mu_q(t)| < \mu_q(t)^{-1}|\mu_q(s)^2 - \mu_q(t)^2|\le 2 a_{k(q)}^{-1/2}\left(\left(b_{k(q)} -a_{k(q)}\right) + ca_{k(q)}^{-\gamma}\right) :=\epsilon_q
$$
where  $C>0$ is independent of $q$ and of $s\in J_q$.  Now $(\mbox{H}_1)$ implies  that
$\epsilon_q \to 0$ as $q\to\infty$.  Hence, $\mu_q(s)-\mu_q(t)  = o(1)$ as $q\to \infty$ uniformly with respect to $s\in J_q$. Moreover, 
\[
\mu_q(s) = \mu_q(t)\left(1 + \frac{o(1)}{\mu_q(t)}\right)= \mu_q(t)\left(1+ o\left(\frac{1}{|q|}\right)\right)\quad \mbox{as}\ q\to \infty
\]
uniformly with respect to $s\in J_q$ since $\mu_q(t) \ge \mu_q^0/2 \ge (2c_0)^{-1}|q|$ for $|q|\ge q_0\gg 1$. 
\finishproof

\noindent

Consider the function $t\to \beta_t(\omega)= \langle \omega, I_t(\omega)\rangle - L_t(I_t(\omega))$.   

\begin{Lemma}\label{Lemma:main}
Suppose that $(\mbox{H}_1)-(\mbox{H}_2)$ holds.  Then $\beta_t(\omega)=\beta_0(\omega)$ for any $t\in [0,\delta]$ and $\omega\in\Omega_\kappa^0$.
\end{Lemma}
{\em Proof}. 
We are interested in the variation 
\[
\dot{\beta}_t(\omega):= \frac{d}{dt}\beta_t(\omega).
\]
Fix $t\in [0,\delta)$ and choose $\omega$ in the set $\Xi_\kappa^t$ given by Lemma \ref{Lemma:quantization}. Consider the quasi-mode of order $N>2d+2\ge 2$  constructed by Theorem \ref{Theorem:quasimodes}. 
Now Lemma \ref{Lemma:isospectral} and Theorem \ref{Theorem:quasimodes}, (vii), imply together that 
\begin{equation}
\zeta_q(s) = \frac{k + \vartheta_0/4}{\mu_q(s)} = \frac{k + \vartheta_0/4}{\mu_q(t)(1 + o(1/|q|))}= \frac{k + \vartheta_0/4}{\mu_q(t)} + o\left(\frac{1}{|q|}\right)= I_t(\omega) + o\left(\frac{1}{|q|}\right)\
\label{eq:variation-of-zeta}
\end{equation}
as $|q|\to \infty$ and uniformly with respect to $s\in J_q$.  On the other hand  Theorem \ref{Theorem:quasimodes}, (vi), yields 
\[
L_s\left(\zeta_q(s)\right) = 2\pi\frac{k_n-\vartheta/4}{\mu_q(s)} 
 +   O(|q|^{-2}) 
\]
uniformly with respect to $s\in J_q$ and  using Lemma \ref{Lemma:isospectral} we obtain as in \eqref{eq:variation-of-zeta} that 
\begin{equation}
L_s\left(\zeta_q(s)\right)= L_t\left(\zeta_q(t)\right)  + o\left(\frac{1}{|q|}\right)\quad \mbox{as}\ q\to \infty
\label{eq:variation-of-L}
\end{equation}
uniformly with respect to $s\in J_q$. Then setting $\eta:=1/|q| \to 0$ we obtain by \eqref{eq:variation-of-zeta} and \eqref{eq:variation-of-L} the equality
\[
\begin{array}{lcrr}
L_{t+\eta}\big(I_t(\omega)\big)=L_{t+\eta}\big(\zeta_q(t+\eta)+ o(\eta)\big)=L_{t+\eta}\big(\zeta_q(t+\eta)\big)+ o(\eta)\\[0.3cm]
= L_{t}\big(\zeta_{q}(t)\big)+ o(\eta)=  
L_t\big(I_t(\omega)\big)  + o(\eta).
\end{array}
\]
We have used also that the map $[0,\delta]\to L_s\in C^\infty (D)$ is $C^1$. Hence, 
\begin{equation}
 \dot{L_t}(I_t(\omega))=\frac{d}{ds}L_s(I_t(\omega))\big|_{s=t} = 0 \quad \forall\,  \omega\in \Xi_\kappa^t\, . 
\label{eq:derivative-of-L}
\end{equation}
On the other hand, $\Xi_\kappa^t$ is dense in $\Omega_\kappa^0$ since any point of $\Omega_\kappa^0$ is of positive Lebesgue density and   $\Omega_\kappa^0\setminus\Xi_\kappa^t$ has measure zero, and by continuity (the function $I\to \dot L_t(I)$ is smooth) we get  \eqref{eq:derivative-of-L} for any  $\omega\in\Omega_\kappa^0$.  The point $t$ has been fixed arbitrary in $[0, \delta)$, hence, \eqref{eq:derivative-of-L} holds true for every $t\in [0, \delta)$ and $\omega\in\Omega_\kappa^0$. 
Now differentiating $\beta_t(\omega)$ with respect to $t$   we obtain
\[
 \dot{\beta}_t(\omega) = \langle \omega, \dot{I}_t(\omega)\rangle - \dot{L}_t(I_t(\omega)) - \langle \nabla L_t(I_t(\omega)), \dot{I}_t(\omega)\rangle = 0 \quad \forall\,  \omega\in \Omega_\kappa^0 
\]
since $\nabla L_t(I_t(\omega))=\omega$. Hence, $\beta_t(\omega)=\beta_0(\omega)$ for every $t\in [0,\delta)$ and $\omega\in \Omega_\kappa^0$.
By continuity we get  the last equality for every $t\in [0,\delta]$ as well. \finishproof

Recall that $\Omega_\kappa^0$ is a set of points of positive Lebesgue density. Differentiating the equality 
\[
\langle \omega, I_t(\omega)\rangle - L_t(I_t(\omega)) =\beta_t(\omega)=\beta_0(\omega)= \langle \omega, I_0(\omega)\rangle - L_0(I_0(\omega))
\]
with respect to $\omega\in \Omega_\kappa^0$ and using Lemma \ref{Lemma:flat} we get $I_t(\omega)=I_0(\omega)$. Then plugging it in the expression of $\beta_t(\omega)$ we obtain $E_{\kappa,t}=E_{\kappa,0}$ as well as  the equality $L_t(I)=L_0(I)$, $I\in E_{\kappa, 0}$.  This completes the proof of Theorem \ref{Th:main1}. \finishproof

\part{KAM theorems and Birkhoff Normal Forms}
\section{KAM theorems}\label{Sec:KAM}
In this Section we prove  KAM theorems and obtain BNF for $C^k$ smooth families of Hamiltonians $t\to H_t$ or exact symplectic maps $t\to P_t$. The main novelty in it can be briefly summarized as follows 
\begin{itemize}
	\item The constant  $\epsilon$ in the smallness condition depends only on the dimension of the configuration space and on the exponent in the Diophantine condition;
	\item $C^k$ smooth families of invariant tori $t\to \Lambda_t(\omega)$ with Diophantine frequencies are obtained;
	\item $C^k$ smooth with respect to the parameter $t$ BNF is obtained around the union of $\Lambda_t(\omega)$;
	\item Uniform estimates in the whole scale of H\"older spaces are obtained.  To this end a new approach to the iterative schema is proposed. The Modified Iterative  Lemma proven in Sect. \ref{Prop:IterativeLemma} provides in a limit smooth functions in the whole domain $\Omega$  (not only on the Cantor set $\Omega_{ \kappa}$) with a good control of the H\"older  norms. In particular, it avoids the Whitney $C^\infty$ extension theorem. 
\end{itemize}
In order to formulate the main results we recall the notion of the Legendre transform.
Let $D\subset \R^d$, $d\ge 1$, be an open set. 
We say that a real valued function $F\in C^\infty(D,\R)$   is {\em non-degenerate} if
\begin{equation}\label{eq:non-degenerate}
\mbox{$\nabla F  : D \longrightarrow  D^\ast:=\nabla F(D) \subset \R^d$ is a diffeomorphism.}
\end{equation}
The Legendre transform $F^\ast$ of $F$ is defined by
\[
F^\ast(\xi)\, :=\,  \mbox{Crit.val.\,}_{x\in D}\{\langle x,\xi\rangle - F(x)\}
\]
which is equivalent to 
\begin{equation}\label{eq:L-transform}
F(x) + F^\ast(\xi) = \langle x, \xi \rangle\, , \ \mbox{where $x\in D$ and $\xi = \nabla F(x)\in D^\ast$.}
\end{equation}
It is easy to see that  $F^\ast\in C^\infty(D^\ast,\R)$ and that $\nabla F^\ast : D^\ast \longrightarrow {D}$ is the inverse  to the map  \eqref{eq:non-degenerate}.
Moreover, $D^{\ast\ast}= D$ and $F^{\ast\ast}= F$.

The real valued function  $F$  defines a non-degenerate (in Kolmogorov sense) completely integrable Hamiltonian in $\T^d\times D$. Hereafter, $\T^d:=\R^d/2\pi\Z^d$.  The corresponding Hamiltonian flow is given by $(s,\theta,r) \to (\theta + s\nabla F(r), r)$. The frequency vector of the restriction of the flow to the invariant torus $\T^{d}\times \{r\}$ is $\omega= \nabla F(r)\in \Omega:=D^\ast$ and the corresponding rotation vector is $\omega/2\pi$. We work here with  frequency vectors instead of rotation vectors because they are more adapted to the Fourier analysis. One can parameterize the invariant tori by their frequency vectors $\omega\in \Omega$ since $F$ is non-degenerate. We are interested below in  families of non-degenerate completely integrable Hamiltonians $F_t$, $t\in [0,\delta]$, with frequency vectors in a fixed open set $\Omega\subset \R^d$. To this end we consider  a family of non-degenerate functions $F^{\ast}_t\in C^\infty(\Omega)$ and define $F_t$ as the Legendre transform of $F^{\ast}_t$ in $D_t:= \nabla F^{ \ast}_t(\Omega)$. The advantage is that the set of frequency vectors is independent of the parameter $t$. In particular,  the Diophantine conditions will be  the same for all $t$. 
The same discussion holds as well for  families of completely integrable exact symplectic maps.

This part is organized as follows. The basic KAM theorem is proved in Sect. \ref{Sec:KAM-parameters}.   

\subsection{KAM theorems for $C^k$ families of Hamiltonians}\label{Subsec:KAM-hamiltonians}
Let $\Omega\subset \R^{n}$ be   an \emph{ open convex  bounded set},  $k\in\{0;1\}$, and $\delta>0$.   Denote by $ \overline{\Omega}$ the closure of $\Omega$ in $\R^n$. Consider  a $C^k$ family $H^{0\ast}$  of real valued   functions 
\[
[0,\delta]\ni t\to H^{0\ast}(\cdot, t)= H_t^{0\ast}(\cdot)\in C^\infty( \overline{\Omega},\R)
\]
satisfying the non-degeneracy condition 
\begin{equation}\label{eq:non-degenerate2}
\mbox{$\nabla H_t^{0\ast} : \Omega \longrightarrow   D_t:=\nabla H_t^{0\ast}(\Omega)$ is a diffeomorphism}.    
\end{equation}
The corresponding  family $H^{0}$ of Legendre transforms 
\[
[0,\delta]\ni t\to H_t^{0} =  H_t^{0\ast\ast}\in C^\infty(D_t,\R)
\]
is $C^k$ as well and 
$ H_t^{0}$ satisfies  \eqref{eq:non-degenerate} on $D_t$ for each $t$.  Consider a $C^k$ family $H$ of perturbations 
\[
[0,\delta] \ni t \to H(\cdot,t) = H_t(\cdot) \in C^\infty(\A_t,\R)
\]
of $H^0$ where $\A_t:=\T^n\times D_t$. 

Let us introduce the arithmetic conditions on the frequency vectors. 
Fix $\kappa >0$ and $\tau>n-1$, and denote by  $\widetilde D(\kappa,\tau)$ the set of all $\omega\in \R^n$ satisfying the 
 $(\kappa,\tau)$-Diophantine  condition 
\begin{equation}  
\forall\, 0\neq k\in \Z^{n}\  :  \quad                    
|\langle\omega,k\rangle |\ \ge \ \frac{\kappa}{|k|^{\tau}}\, . 
                       \label{eq:sdc2}                                        
\end{equation}
Denote by $\Omega_\kappa$ the set of all $(\kappa,\tau)$-Diophantine vectors $\omega\in\Omega$ such that the distance from $\omega$ to the complement $\R^{n}\setminus \Omega$ of $\Omega$ is $\ge \kappa$.  
We will often use the following notation 
\begin{equation}\label{eq:omega-kappa}
\left\{
\begin{array}{llcr}
\Omega +\kappa &:=& \{\omega\in \R^n :\ \mbox{\rm dist} (\omega, \Omega) <\kappa\}
\\[0.3cm]
\Omega -\kappa &:=& \{\omega\in \Omega :\ \mbox{\rm dist} (\omega,\R^{n}\setminus \Omega) >\kappa\}.
\end{array}
\right.
\end{equation}
Then $\Omega_\kappa= \widetilde D(\kappa,\tau)\cap \overline{\Omega -\kappa}$. 

In order to formulate the smallness condition   we need the following notations. Firstly we define  weighted $C^\ell$ H\"older norms as in \cite{Poe}. 
Given   $\ell\ge 0$, $0<\kappa\le 1$, and a domain $D\subset \R^n$, we denote the weighted (with respect to the small parameter $\kappa$)  $C^\ell$-norm of $u\in C^\ell( \T^n\times D,\R^k)$ evaluated at $\T^n\times D$ by 
\begin{equation}
\big\|u\big\|_{\ell,\T^n\times D;\kappa}:= \big\|u\circ \sigma_{\kappa}\big\|_{\,  C^\ell\left(\sigma_{\kappa}^{-1}(\T^n\times D)\right) }
\label{eq:holder-norms}
\end{equation}
where  $\|\cdot\|_{C^\ell}$ is the corresponding H\"older norm (see Sect. \ref{Sec:ApprLemma}) and 
$\sigma_{\kappa}:\T^n\times \R^n\to \T^n\times \R^n$ is  the partial dilation  $\sigma_{\kappa}(\theta,r):= (\theta, \kappa r)$. If $\ell\in\N$ then
\[
\big\|u\big\|_{\ell,\T^n\times D;\kappa}\, = \, \sup_{|\alpha|+|\beta|\le \ell}\, \sup_{(\theta,r)\in \T^n\times D}\, |\partial_\theta^\alpha (\kappa \partial_r)^\beta u(\theta,r)|,
\]
where $|\cdot|$ is the Euclidean norm. 
In the same way we introduce the norm $\|u\|_{\ell,D;\kappa}$ for $u\in C^\ell(D, \R^k)$. 
We set as well
\begin{equation}\label{eq:3-norm}
|\!|\!|u |\!|\!|_{\ell,D;\kappa}\, =\, \sup_{0\le m \le \ell}\|u\|_{\ell-m,D;\kappa}, \quad \mbox{where}\ m\in \N .
\end{equation}
If $D$ is \emph{convex}, then $|\!|\!|u |\!|\!|_{\ell,D;\kappa}\, =\, \|u\|_{\ell,D;\kappa}$. 
Given $\ell\ge 1$ and a family of functions $$u:=\{u_t\in C^\infty(\Omega,\R): \ t\in [0,\delta]\}$$ we set
\begin{equation}
\label{eq:S-sequence}
S_{\ell}(u) := \sup_{0\le t\le \delta }\, \big(1+ \|u_t\|_{C^1(\Omega)}\big)^{\ell-1} \big(1+ \|u_t\|_{C^\ell(\Omega)}\big).
\end{equation}
This expression
 arises   when one evaluates the   $C^\ell$-norms of a composition of functions the form $f_t\circ u_t$  with $u_t=\nabla K_t^{0\ast}$ (see Appendix, Sect. \ref{subsec:interpolation}). If $\Omega$ is {\em convex} then $\|u_t\|_{C^\ell(\Omega)}\le \|u_t\|_{C^\mu(\Omega)}$ for $0\le \ell\le \mu$ and the function $\ell \to S_{\ell}(u)$ becomes increasing in $[1,+\infty)$. 
 
Fix  $\vartheta_0>1$  and set 
\begin{equation}
\ell_0:=2\tau + 2 + 2\vartheta_0\quad \mbox{and} \quad  \ell(m):=  2m(\tau + 1)+ \ell_0 \, , \quad m \ge 0. 
\label{eq:ell(m)}
\end{equation}  
 Given $0<\varrho,\kappa \le 1$ and $m\ge 0$, we denote by ${\mathcal A}_m^0 $ the expression 
\begin{equation}
\label{eq:A-0-sequence}
{\mathcal A}_m^0 \, := \,  \sup_{0\le t\le \delta }\, \left( \varrho^2 |\!|\!| \partial^2 H_t^0 |\!|\!|_{\ell(m), D_t;\kappa} \, +\,  |\!|\!| H_t-H_t^0|\!|\!|_{\ell(m),\A_t; \kappa}\right)
\end{equation}
and set 
\begin{equation}
\label{eq:A-sequence}
{\mathcal A}_m \, =\, S_{\ell(m)+1}(\nabla H^{0\ast})\, {\mathcal A}_m^0 \, . 
\end{equation}
Here $\partial^2 H_t^0(I)$ is the Hessian matrix of $H_t^0$ at $I\in D_t$. The role of the small parameter $\varrho$ is to compensate the norm of the Hessian matrix which  could be very large. The function $m\to {\mathcal A}_m$, $m\ge 0$, is {\em increasing} when $\Omega$ is convex.

Let $\Phi_{t}^s := \exp\left(s X_{H_{t}} \right)$, $s\in\R$, be 
the flow of the Hamiltonian vector field $X_{H_{t}}$ with Hamiltonian $H_{t}$ in $ \A_t= \T^n\times  D_t$.  
Recall that  
for any $\omega\in \Omega$ the map  $R_{\omega}:\T^n\to\T^n$ stands for the translation $R_{\omega}(\varphi) =\varphi + \omega \, \mbox{\rm (mod $2\pi$)}$. 
Fix  $k\in \{0,1\}$. 

\begin{Theorem}\label{Theo:KAM2} There exists $\epsilon = \epsilon(n,\tau,\vartheta_0)>0$ depending only on $n$, $\tau$ and $ \vartheta_0$ such that 
the following holds.\\
\noindent
Let  $ \Omega\subset \R^{n}$ be  an  open convex  bounded set,  $0<\varrho\le\kappa\le 1$ and $\Omega_\kappa\neq\emptyset$. Let $H^{0\ast}$ be a $C^k$ family  of real-valued functions $[0,\delta]\ni t\to H_t^{0\ast}\in C^\infty( \overline\Omega,\R)$   satisfying \eqref{eq:non-degenerate2},  and let  
 $[0,\delta]\ni t\to H_t\in C^\infty( \A_t,\R)$ be a $C^k$ family of Hamiltonians such that
\begin{equation}\label{eq:smallness-condition2}
S_{\ell_0}(\nabla H^{0\ast})\,{\mathcal A}_0^0\,  =\, \sup_{0\le t\le \delta }\, \left( \varrho^2 |\!|\!|\partial^2 H_t^0 |\!|\!|_{\ell_0, D_t;\kappa} \, +\,  |\!|\!| H_t-H_t^0|\!|\!|_{\ell_0,\A_t; \kappa}\right) S_{\ell_0}(\nabla H^{0\ast}) \, 
\le \, \epsilon \varrho \kappa .
\end{equation}
Then there exists  a $C^k$ mapping $[0,\delta]\ni t\to    \Psi_t=(\widetilde  U_t,\widetilde  V_t)\in C^\infty(\T^{n}\times  \Omega; \T^{n}\times  D_t)$ such that
\begin{enumerate}
\item[(i)]
for any $\omega\in \Omega_{\kappa}$,  
$[0,\delta]\ni s\to  \Lambda_t(\omega)=  \Psi_{t,\omega}(\T^{n})$ is a $C^k $  family of Kronecker invariant tori 
of the Hamiltonian vector fields $X_{H_t}$ with a frequency vector $\omega$, where $   \Psi_{t,\omega}:= \Psi_t(\cdot; \omega)$. Moreover, 
for any $\omega\in\Omega_\kappa$ and  $s\in\R$ the following diagram is commutative
\[
\displaystyle{\begin{array}{cccl} 
\displaystyle \T^{n}&\stackrel{R_{s\omega}}{\longrightarrow}&\T^{n}\cr
\downarrow\lefteqn{ \Psi_{t,\omega}}& &\downarrow\lefteqn{ \Psi_{t,\omega}} \cr
\displaystyle \Lambda_t(\omega)&\stackrel{\Phi_{t}^s}{\longrightarrow}&\Lambda_t(\omega)&  
\end{array} }
\]
\item[(ii)] for any $m\in \{0\}\cup [1,+\infty)$ the following estimates hold 
\begin{equation}\label{eq:estimates-tilde-F}
\begin{array}{lcrr}
\displaystyle \big|\partial_\varphi^\alpha (\kappa\partial_\omega)^\beta \big(\widetilde U_t(\varphi;\omega)-\varphi\big)\big| \, \le\,  C_m \frac{{\mathcal A}_m}{\kappa\varrho}\\[0.3cm]
\displaystyle
  \big|\partial_\varphi^\alpha (\kappa\partial_\omega)^\beta   \big(\widetilde V_t(\varphi;\omega)-\nabla H_t^{0\ast} (\omega)\big)\big| 
\, \le\,  C_m \frac{{\mathcal A}_m}{\varrho} \left(1+\frac{{\mathcal A}_1}{\varrho}\right)^{m} 
\end{array}
\end{equation}
for each $t\in [0,\delta]$, $(\varphi,\omega)\in \T^{n}\times \Omega$  and $\alpha,\beta\in \N^{n}$  with  $|\alpha| +|\beta|(\tau+1)\le m(\tau+1) +1$, where the constant $C_m>0$ depends only on $n$, $\tau$, $\vartheta_0$  and $m$. 
\item[(iii)]
${\rm supp\, }\big((\widetilde U_t,\widetilde V_t) -({\rm id}, \nabla H_t^{0\ast} )\big) \subset \T^{n}\times (\Omega-\kappa/2)$.
\end{enumerate}
\end{Theorem}

\begin{Remark}\label{rem:analiticity-1}
	If $P$ is analytic with respect to $t$ in  the disc $B(0,a):= \{t\in \C:\, |t|<a\} $  and \eqref{eq:smallness-condition2} holds for $t\in B(0,a)$, then $\Psi$ and $\phi$ can be chosen to be  analytic with respect to $t$ in $B(0,a)$. Moreover, 
	for any $\alpha,\beta\in \N^n$ of length $|\alpha|+|\beta|(\tau+1)\le m(\tau+1)+ 1$ and $0<a_1<a$,  the estimate \eqref{eq:estimates-tilde-F} holds true for $t\in B(0,a_1)$,  where the supremum with respect to $t$  in the definition of ${\mathcal A}_m$ is taken in $B(0,a)$. 
	\end{Remark}

\noindent
{\em Proof}. 
The idea of the proof is given by P\"oschel in \cite{Poe1}. It can be summarized as follows.  Let us  fix $\omega\in \Omega$ set $r= \nabla H_t^{0 \ast}(\omega) + I$ and apply Taylor's formula up to order two to the function  $I\to H_t^0(\nabla H_t^{0 \ast}(\omega) + I)$ at $I=0$. Then the afine linear with respect to $I$ term is just $N_t(I;\omega)=e_t(\omega)+\langle I,\omega\rangle$ and we  put the quadratic term in the perturbation. Multiplying the perturbation by suitable cut-off functions  we obtain a Hamiltonian $\widetilde H_t(\theta,I;\omega)=N(I;\omega) + P_t(\theta,I;\omega)$, where $t\to P_t\in C^\infty(\T^n\times\R^n\times\Omega)$ is $C^k$ and $P_t$ are compactly supported with respect to $(I;\omega)$.  The smallness condition allows one  to apply Theorem \ref{Theo:A}. The main difficulty in the proof is to obtain the corresponding estimates in $C^\ell$, $\ell\ge 0$.
We devide the proof in several steps.\\

\noindent
{\em Step 1. Construction of the Hamiltonian $\widetilde H_t(\theta,I;\omega)$.}
Given $\omega\in \Omega -\kappa/4$ 
and $I$ in the ball $B^n(0,\varrho)\subset\R^n$ of center $0$ and radius $\varrho$,   we set 
$r=\nabla H^{0\ast}_t(\omega) + I$.  Choosing $\epsilon=\epsilon(n,\tau,\vartheta_0)\le 1/9$ in  \eqref{eq:smallness-condition2} we will show that
\begin{equation} \label{eq:domain-H}
\omega\in \Omega-\kappa/4,\ I\in B^n(0,\varrho) \quad \Longrightarrow \quad  \left\{
\begin{array}{lcrr}
\nabla H^{0\ast}_t(\omega) + I\in D_t \quad \mbox{and}\\[0.3cm]
\displaystyle \nabla  H^{0}_t(\nabla H^{0\ast}_t(\omega) + I)\in \Omega-\frac{1}{8}\kappa  
\end{array}
\right.
\end{equation}
for each $t\in [0,\delta]$.  Fix $t\in [0,\delta]$.  
The smallness condition \eqref{eq:smallness-condition2} implies that $\|\partial^2  H^{0}_t\|_{C^0(D_t)} \le \epsilon \kappa /\varrho$ since $S_{\ell_0}(\nabla H^{0\ast})>1$. Then  for each  $\omega\in \Omega-\kappa/4$  there is a positive number $c\le \varrho$ such that for any  $I\in B^n(0,c)$ the following relation holds
\begin{equation} \label{eq:domain-H-1}
 \left\{
\begin{array}{lcrr}
\nabla H^{0\ast}_t(\omega) + I\in D_t \quad \mbox{and}\\[0.3cm]
\displaystyle |\nabla  H^{0}_t(\nabla H^{0\ast}_t(\omega) + I) -\omega|\,  \le \, \|\partial^2  H^{0}_t\|_{C^0(D_t)}\, |I|\,  \le\,  \epsilon \kappa \frac{c}{\varrho}\,  <\,  \frac{\kappa}{9}  \, .
\end{array}
\right.
\end{equation}
Using the notations \eqref{eq:omega-kappa} one obtains 
\[
\nabla  H^{0}_t(\nabla H^{0\ast}_t(\omega) + I)\in (\Omega-\kappa/4)+\kappa/9=(\Omega-\kappa/8)-\kappa/72
\] 
for any $I\in B^n(0,c)$. If $c<\varrho$, then there exists $c<c'\le\varrho$ such that  \eqref{eq:domain-H-1} still holds for each $I\in B^n(0,c')$. 
This  proves  \eqref{eq:domain-H}. 
Then  Taylor's formula yields
\[
H_t^0(r) = e_t(\omega) + \langle I,\omega\rangle + \int_0^1(1-s)\langle \partial^2H_t^0(\nabla H_t^{0\ast}(\omega) + sI)I,I\rangle\,  ds.
\]
for  $\omega\in \Omega-\kappa/4$ and $I\in B^n(0,\varrho)$,  
where $e_t(\omega) := H_t^0(\nabla H^{0\ast}_t(\omega))$. 

In order to apply Theorem \ref{Theo:A}  we need  suitable cut-off functions. 
\begin{Lemma}\label{lemma:cut-off-omega}
	For any open set $U\subset \R^n$ and $0<\varepsilon\le 1$ there exists a smooth cut-off function  $\psi^U_\varepsilon\in C_0^\infty(\R^{n}, [0,1])$
	such that 
	$\psi^U_\varepsilon=1$ on $U-\varepsilon$, $\psi^U_\varepsilon=0$ on the complement of $U + \varepsilon$, and
	\[
	\|\psi^U_\varepsilon\|_{\ell;\varepsilon} \le C_\ell
	\]
	for any $\ell\ge 0$, where the positive constants $C_\ell=C(\ell,n)$ depend only on $\ell$ and $n$.
\end{Lemma}
The proof of the Lemma is given in the Appendix, Sect. \ref{subsec:composition-inverse}. \\

Denote by $\psi_\kappa\in C_0^\infty(\R^{n}, [0,1])$ the function given by Lemma \ref{lemma:cut-off-omega} with $U=\Omega-\kappa/2$ and $\varepsilon=\kappa/4$. Then  $\psi_\kappa=1$ on $\Omega-3\kappa/4$, $\psi_\kappa=0$ on the complement of $\Omega-\kappa/4$, and
\begin{equation}\label{eq:estimates-psi-kappa}
\| \psi_\kappa\|_{\ell;\kappa} \le C_\ell
\end{equation}
for any $\ell\ge 0$, where the constants $C_\ell=C(\ell,n)$ depend only on $\ell$ and $n$. Let $\widetilde\psi_\varrho$ be the cut-off function given by  Lemma \ref{lemma:cut-off-omega} with $U=B^n(0, 3\varrho/4)$ and  $\varepsilon=\varrho/4$. Then the support of $ \widetilde\psi_\varrho$ is contained in   $B^n(0,\varrho)$, $\widetilde\psi_\varrho=1$ on $B^n(0,\varrho/2)$, and 
\begin{equation}\label{eq:estimates-psi-rho}
\|\widetilde\psi_\varrho\|_{\ell;\varrho} \le C_\ell
\end{equation}
for any $\ell\ge 0$, where the constants $C_\ell=C(\ell,n)$ depend only on $\ell$ and $n$.

It follows from \eqref{eq:domain-H} that the function 
\[
P_t^0(I;\omega):= \int_0^1(1-s)\langle \partial^2H_t^0(\nabla H_t^{0\ast}(\omega) + sI)I,I\rangle \psi_\kappa(\omega)\widetilde \psi_{\varrho}(I)\,  ds,
\]
is well defined and compactly supported in $B^n(0,\varrho)\times \Omega$. 
Setting
\[
P_t^1(\theta,I;\omega):= (H_t-H_t^0)(\theta, \nabla H^{0\ast}_t(\omega) + I)\psi_\kappa(\omega)\widetilde\psi_{\varrho}(I)\quad \mbox{and}  \quad P_t:= P_t^0 + P_t^1, 
\]
we consider the Hamiltonian
\[
\widetilde H_t(\theta,I;\omega)= e_t(\omega) + \langle I,\omega\rangle + P_t(\theta, I;\omega).
\]
We have 
\begin{equation}\label{eq:equality-for-H}
H_t(\theta,\nabla H^{0\ast}_t(\omega) + I)=\widetilde H_t(\theta,I;\omega)
\end{equation}
 for $I\in B^n(0,\varrho/2)$ and $\omega\in \Omega-3\kappa/4$. \\

\noindent
{\em Step 2. H\"older estimates of $P_t$.}
For any $\ell\ge 1$ we are going to evaluate  the weighted norm 
\[
\|P_t\|_{\ell;r,\kappa} := \|P_t\circ \sigma_{r,\kappa}\|_{C^\ell(U)}\, , \quad U:=\sigma_{r,\kappa}^{-1}(\T^n\times B^n(0,\varrho)\times \Omega), 
\]
introduced in \eqref{eq:norm-rho-kappa}. In order to estimate  $\|P_t^0\|_{\ell;r,\kappa}$ we 
set $\Gamma:= B^n(0,\varrho)\times (\Omega-8\kappa/9)$, $\Gamma_t:= B^n(0,\varrho)\times D_t$, 
\[
Q_t^1(I,r)=\int_0^1(1-s)\,  \partial^2H_t^0(r+sI)\, ds 
\]
and $Q_t^0(I,r)= \langle Q_t^1(I,r)I,I\rangle$ for $(I,r)\in \Gamma_t$. 
Then we write 
\[
P_t^0(I,\omega)=\psi_\kappa(\omega)\widetilde\psi_{\varrho}(I)Q_t^0(I, \nabla H_t^{0\ast}(\omega)). 
\] 
The function $Q_t^0\circ ({\rm id}, \nabla H_t^{0\ast}) $ belongs to $C^\infty(\overline \Gamma)$, where $\Gamma$ is convex. Then 
using Remark \ref{rem:interpolation}, \eqref{eq:estimates-psi-kappa} and \eqref{eq:estimates-psi-rho} we obtain
\[
\|P_t^0\|_{\ell;\varrho,\kappa}  \,  \le  \, 
C_{\ell} \|Q_t^0\circ ({\rm id}, \nabla H_t^{0\ast}) \|_{\ell;\varrho,\kappa}.
\]
Proposition \ref{prop:composition-holder1}, {\em 3}, and Remark \ref{rem:anisotrop}
imply for  any  $\ell\ge 1$  the estimate
\[
\|Q_t^0\circ ({\rm id}, \nabla H_t^{0\ast}) \|_{\ell;\varrho,\kappa} \, \le  \, 
C_{\ell} \,  |\!|\!| Q_t^0|\!|\!|_{\ell,\Gamma_t;\varrho,\kappa}  \, 
S_{\ell}(\nabla H^{0\ast}), 
\]
where the norm $ |\!|\!| \cdot |\!|\!|$ is defined in \eqref{eq:3-norm}. 
Using  \eqref{eq:Leibnitz1}  and the inequality $0<\varrho\le\kappa$ we obtain 
\[
\| Q_t^0\|_{\ell,\Gamma_t;\varrho,\kappa} \le  C_{\ell} \, \varrho^2 |\!|\!| Q_t^1|\!|\!|_{\ell, \Gamma_t;\varrho,\kappa}  \le  C_{\ell} \, \varrho^2 |\!|\!| Q_t^1|\!|\!|_{\ell, \Gamma_t;\kappa,\kappa} \le  C_{\ell} \, \varrho^2 |\!|\!| \partial^2 H_t^0 |\!|\!|_{\ell, D_t;\kappa} 
\]
as $|I|\le \varrho$ for $(I,r)\in\Gamma_t$. This inequality implies
\[
|\!|\!| Q_t^0|\!|\!|_{\ell,\Gamma_t;\varrho,\kappa} \le  C_{\ell} \, \varrho^2 |\!|\!| \partial^2 H_t^0 |\!|\!|_{\ell, D_t;\kappa} 
\]
and we obtain
\[
\|P_t^0\|_{\ell;\varrho,\kappa}  \,  \le  \, 
C_{\ell} \, \varrho^2\, |\!|\!| \partial^2H^0 |\!|\!|_{ \ell, D_t; \kappa}  \, 
S_{\ell}(\nabla H^{0\ast}).  
\]
On the other hand $\|P_t^1\|_{\ell;\varrho,\kappa} \le \|P_t^1\|_{\ell;\kappa,\kappa}$  since $0<\varrho\le \kappa\le 1$ and by the same argument we get
\[
\|P_t^1\|_{\ell;\varrho,\kappa} 
\le C_{\ell} |\!|\!|H_t-H_t^0|\!|\!|_{\ell,\A_t; \kappa}
S_{\ell}(\nabla H^{0\ast}).
\]
Finally we obtain
\begin{equation}\label{eq:estimate-P-with-H}
\|P_t\|_{\ell;\varrho,\kappa} 
\le \displaystyle C_{\ell} \left(\varrho^2\, |\!|\!| \partial^2H^0_t |\!|\!|_{ \ell, D_t; \kappa}\, +\, |\!|\!|H_t-H_t^0|\!|\!|_{\ell,\A_t; \kappa}\right)\ 
S_{\ell}(\nabla H^{0\ast}).
\end{equation}

\noindent
{\em Step 3. Application of Theorem \ref{Theo:A} .}
The estimate  \eqref{eq:smallness-condition2} gives 
\begin{equation} \label{eq:0-estimates-P}
\|P_t\|_{\ell_0;\varrho,\kappa} \, \le \, C\,  {\mathcal A}_0 S_{\ell_0}(\nabla H_t^{0\ast}) \, \le\,   C \, \epsilon \varrho \kappa \, .
\end{equation}
This allows us to apply Theorem \ref{Theo:A} to the Hamiltonian $(\theta,I)\to \widetilde H_t(\theta,I;\omega)$. 
Set $ \Psi_t = (\widetilde U_t, \widetilde V_t)$, where $\widetilde U_t= U_t$ and $\widetilde V_t=(\nabla H^{0\ast}_t)\circ \phi_t +V_t$, where $(U_t,V_t,\phi_t)$ are given by Theorem  \ref{Theo:A}. Notice that $\|V_t\|_{C^0} \le c\epsilon \varrho$ in view of the estimates in (ii), Theorem \ref{Theo:A}, where the constant $c$ depends only on  $n$, $\tau$ and $\vartheta_0$, and taking $\epsilon<\min(1,1/c)/2$ we obtain $V_t(\theta;\omega)\in B^n(0,\varrho/2)$ for any $\theta\in \T^n$ and $\omega\in \Omega-\kappa$. In the same way we get $\phi_t(\omega)\in \Omega-3\kappa/4$ for $\omega\in \Omega-\kappa$. In particular, $\widetilde\psi_{\varrho}(V_t(\theta;\omega))=1$ and $\psi_{\kappa}(\phi_t(\omega))=1$ for any  $(\theta,\omega)\in \T^n\times (\Omega-\kappa)$. 

By \eqref{eq:smallness-condition3} and \eqref{eq:main-estimates}, we have 
$$|d_\theta U_t(\theta;\omega)-\mathrm{Id} |\le C_1(n,\tau,\vartheta_0)\epsilon \le 1/2$$ 
for $(\theta,\omega,t)\in {\T}^n\times \Omega\times [0,a]$, choosing $\epsilon $ sufficiently small.  
Now  Remark  \ref{rem:Kronecker} implies  that for any $\omega\in \Omega_\kappa$ and $t\in [0,\delta]$ the Lagrangian manifold 
$ \Lambda_t(\omega):=  \Psi_{t}(\T^n;\omega)$
is a Kronecker invariant torus of $H_t$ of a frequency vector $\omega$ satisfying  (i) in Theorem \ref{Theo:KAM2}. \\

\noindent
{\em Step 4. Estimates of $\widetilde U_t$ and $\widetilde V_t$.}
The estimates (ii), Theorem \ref{Theo:A}, imply (ii) in Theorem \ref{Theo:KAM2} using the estimates of $\displaystyle\| P_t\|_{\ell;\varrho,\kappa}$ given above. To estimate the derivatives of 
\[
\widetilde V_t- \nabla H^{0\ast}_t = V_t+ (\nabla H^{0\ast}_t)\circ \phi_t - \nabla H^{0\ast}_t
\]
we use (ii), Theorem \ref{Theo:A} and the following
\begin{Lemma}\label{lemma:estimates-V}
For any $m\in \{0\}\cup [1,+\infty)$ the following estimate holds
\[
\|(\nabla H^{0\ast}_t)\circ \phi_t - \nabla H^{0\ast}_t\|_{m,\Omega;\kappa} \le C_{m} \frac{{\mathcal A}_m}{\varrho}\left(1+\frac{{\mathcal A}_1}{\varrho} \right)^m.
\]
\end{Lemma}
{\em Proof}. 
The proof of the lemma  is based on higher order  H\"older estimates of a composition of functions given in the Appendix. 
 We write
\[
\nabla H^{0\ast}_t(\phi_t(\omega))- \nabla H^{0\ast}_t(\omega) = u_t(\omega)\cdot v_t(\omega)
\]
where 
\[
u_t(\omega)= \phi_t(\omega)-\omega\in M_{1,n}(\R) \quad \mbox{and} \quad v_t(\omega)=\int_0^1 (\partial^2 H^{0\ast}_t)(\omega + s(\phi_t(\omega)-\omega))\,  ds \in M_{n,n}(\R). 
\]
Theorem \ref{Theo:A} and \eqref{eq:estimate-P-with-H} imply
\[
\|\phi_t-{\rm id}\|_{m,\Omega;\kappa}\,  \le\, C_m\frac{{\mathcal A}_m^0}{\varrho} S_{\ell(m)}(\nabla H^{0\ast}) 
\]
where
\[
S_{\ell(m)}(\nabla H^{0\ast}) = \sup_{0\le t\le \delta}\big(1+ \|\nabla H^{0\ast}_t\|_{1}\big)^{\ell(m)-1}\big(1+ \|\nabla H^{0\ast}_t\|_{\ell(m)}\big) 
\]
is increasing with respect to $m\in\N$ since $\Omega$ is convex. 
Moreover, \eqref{eq:smallness-condition2} yields $\|\phi_t-{\rm id}\|_{0} \le C \varepsilon \kappa$. 
Notice that $v_t$ is well-defined and $C^\infty$ smooth in the convex set $\overline \Omega$. Indeed, the image of $\Omega$ under the map $\omega\mapsto \omega + s(\phi_t(\omega)-\omega)$ is contained in  $\Omega$ when $\varepsilon<1/2C$ since 
${\rm supp\, }(\phi_t -{\rm id})\subset \Omega- \kappa/2$ and $\|\phi_t-{\rm id}\|_0 <\kappa/2$.

Let $m=0$. Then
\[
\|(\nabla H^{0\ast}_t)\circ \phi_t - \nabla H^{0\ast}_t\|_{C^0(\Omega)}\, \le\,  C_0\frac{{\mathcal A}_0^0}{\varrho} S_{\ell_0}(\nabla H^{0\ast}) \|\nabla H^{0\ast}_t\|_{C^1(\Omega)} \, < \, C_0\frac{{\mathcal A}_0^0}{\varrho} S_{\ell_0+1}(\nabla H^{0\ast}) \, = \, 
C_{0} \frac{{\mathcal A}_0}{\varrho}
\]
since $\Omega$ is convex. 

Let $m\ge 1$. Using  Remark \ref{rem:interpolation} we get 
\[
\|\nabla H^{0\ast}_t\circ \phi_t- \nabla H^{0\ast}_t\|_{m,\Omega;\kappa}\,  \le\, C_m \big( \|\phi_t-{\rm id}\|_{m,\R^n;\kappa}\|v_t\|_{C^0(\Omega)} 
+\|\phi_t-{\rm id}\|_{C^0(\R^n)}\|v_t\|_{m,\Omega;\kappa}\big).
\]
We obtain as above
\[
\begin{array}{rcll}
\displaystyle \|\phi_t-{\rm id}\|_{m,\R^n;\kappa}\|v_t\|_{C^0(\Omega)} \, \le\, C_m\frac{{\mathcal A}_m^0}{\varrho} S_{\ell(m)}(\nabla H^{0\ast}) \|\nabla H^{0\ast}_t\|_{C^1(\Omega)}\\[0.3cm] 
\displaystyle  < \, C_m\frac{{\mathcal A}_m^0}{\varrho} S_{\ell(m)+1}(\nabla H^{0\ast})  = C_m\frac{{\mathcal A}_m}{\varrho}
  \end{array}
\]
since $\Omega$ is convex. 
On the other hand, Proposition \ref{prop:composition-holder-interpolation}, {\em 3}, applied to $f=({\rm id} +s(\phi_t-{\rm id}))\circ \sigma_\kappa$ and $g=\partial^2 H^{0\ast}_t\in C^\infty(\overline \Omega, M_{n,n}(\R))$ yields
\[
\begin{array}{rcll}
\|v_t\|_{m,\Omega;\kappa} &\le& C_m \big(1+  \|\phi_t-{\rm id}\|_{1;\kappa}^{m-1}\big)\\[0.3cm]
 &\times& \displaystyle  \, \big(\|\nabla H^{0\ast}\|_{m+1}(1+ \|\phi_t-{\rm id}\|_{1})+\|\nabla H^{0\ast}\|_{2}\|\phi_t-{\rm id}\|_{m;\kappa} \big)
 \end{array}
\]
where the corresponding  $C^\ell$ norms of $\nabla H^{0\ast}$  are evaluated on $\Omega$. Then
\[
\begin{array}{lcrr}
\displaystyle\|\phi_t-{\rm id}\|_{C^0(\R^n)}\|v_t\|_{m,\Omega;\kappa}\, \le\, C_m \frac{{\mathcal A}_0^0}{\varrho}\Big(1+\frac{{\mathcal A}_1}{\varrho}\Big)^{m-1}S_{\ell_0}(\nabla H^{0\ast}) \\[0.3cm]
 \displaystyle \times\, \left(\|\nabla H^{0\ast}\|_{m+1}\Big(1+\frac{{\mathcal A}_1}{\varrho}\Big) +  \|\nabla H^{0\ast}\|_{2} \frac{{\mathcal A}_m^0}{\varrho} S_{\ell(m)}(\nabla H^{0\ast})  \right)
\end{array}
\]
The interpolation inequalities \eqref{eq:interpolation} and  Remark \ref{rem:interpolation}  yield  for any $ r, s \ge 1$ and $f\in C^{r+s-1}(\overline\Omega)$ the estimate 
$\|f\|_r\|f\|_s \le C_{r,s} \|f\|_1\|f\|_{r+s-1}$, where $C_{r,s}>0$ depends only on $r,s$ and  $n$. Applying this inequality to $f= \nabla K_t^\ast$ and using the convexity of $\Omega$ we obtain
\[
S_{\ell_0}(\nabla H^{0\ast})\|\nabla H_t^{0\ast}\|_{m+1} < C_m S_{\ell_0+m}(\nabla H^{0\ast})< C_mS_{\ell(m)+1}(\nabla H^{0\ast})
\]
and 
\[
S_{\ell(m)}(\nabla H^{0\ast})\|\nabla H_t^{0\ast}\|_{2} \le C_m S_{\ell(m)+1}(\nabla H^{0\ast}).
\]
On the other hand ${\mathcal A}_0^0S_{\ell_0}(\nabla H^{0\ast}) \le \varepsilon \kappa\varrho$ by \eqref{eq:smallness-condition2} and ${\mathcal A}_0^0 \le {\mathcal A}_m^0$ for $m\ge 0$ by definition 
and  we obtain 
\[
\|\phi_t-{\rm id}\|_{C^0(\R^n)}\|v_t\|_{m,\Omega;\kappa}\, \le\, C_m \frac{{\mathcal A}_m}{\varrho}\Big(1+\frac{{\mathcal A}_1}{\varrho}\Big)^{m}
\]
This completes the proof of the Lemma. 
 \finishproof
 
\noindent 
Statement  (iii) follows from the definition of $\widetilde U$ and $\widetilde V$ using Theorem  \ref{Theo:A}.  \finishproof

\subsection{KAM theorem with parameters for  symplectic maps}\label{Subsec:KAM-maps-parameters}
We are going to prove an analogue of Theorem \ref{Theo:A} for symplectic maps. More precisely, given a $C^k$ family of ``small'' exact symplectic  perturbations $(\theta,r) \to P_t(\theta,r;\omega)$ of the translation  $(\theta,r)\to R_\omega(\theta,r):= (\theta+\omega,r)$ with a Diophantine frequency $\omega$, we are going to find a $C^k$ family of Kronecker invariant tori of $P_t(\theta,r;\phi_t(\omega))$, where $t\to \phi_t$ is a $C^k$ family of diffeomorphisms close to the identity map. Moreover, we will estimate the $C^m$, $m\in\N$,  norm of $\phi_t -{\rm id}$, and of the displacement of the Kronecker tori with respect to the corresponding inperturbed tori.

Let $\Omega \subset [0,2\pi)^{n-1}$ be an open convex set.  We identify $\Omega$ with an open convex subset of $\T^{n-1}$. Fix $k\in \{0, 1\}$ and $\varrho_0>0$ and  consider a $C^k$ family of exact symplectic maps 
\[
[0,\delta]\ni t\to P_t(\cdot;\omega)\in C^\infty(\T^{n-1}\times B^{n-1}(0,\varrho_0), \T^{n-1}\times \R^{n-1})
\]
 depending smoothly on a parameter $\omega\in \Omega$. 
 We suppose that  $P_t(\cdot;\omega)$ is defined by a generating functions $\widetilde G_t(\cdot;\omega)$ of the form
\begin{equation} 
\label{eq:generating-function-of P}
\widetilde G_t(\theta,r;\omega):= \langle \theta,r\rangle -\langle \omega,r\rangle -G_t(\theta,r;\omega)
\end{equation}
i.e. 
\[
P_t(\theta-\omega-\nabla_{r}G_t(\theta,r;\omega),r;\omega)=(\theta, r-\nabla_{\theta}G_t(\theta,r;\omega),r)
\]
for any $(\theta,r)\in \T^{n-1}\times B^{n-1}(0,\varrho_0)$ and $\omega\in \Omega$ and  that the $C^2$ norm of $G_t(\cdot;\omega)$ is  sufficiently small. 
 If $G_t=0$ then $P_t(\cdot;\omega)$ becomes a translation with $\omega$, namely, 
$R_\omega(\theta,r)=(\theta+\omega,r)$.  In general we consider $G_t$ as a small perturbation  depending smoothly on the frequency $\omega$ as well. 
Thus the perturbation is a  real valued function $(\theta,I;\omega,t) \mapsto G(\theta,I;\omega,t)$ defined in $\A\times \Omega\times [0,\delta]$, where $\A:= \T^n\times B(0,\varrho_0)$. Hereafter, we assume that
\begin{equation}
G\in C^k\left([0,\delta];C^\infty(\A\times  \Omega)\right)
\label{eq:families-of-perturbations-maps}
\end{equation}
with $k=0$ or $k=1$, 
i.e. the map 
  $t\to G_t:=G(\cdot,t)\in C^\infty(\T^{n-1}\times B(0,\rho_0)\times  \Omega)$ is  $C^k$-smooth on the interval $[0,\delta]$.  

Given  $\ell>0$ and $0<\varrho,\kappa\le 1$ we denote by $\|G_t\|_{\ell;\varrho,\kappa}$ the  weighted H\"older norm 
\begin{equation}\label{eq:norm-rho-kappa-maps}
\|G_t\|_{\ell;\varrho,\kappa} := \|G_t\circ \sigma_{\varrho,\kappa}\|_{C^\ell(\sigma_{\varrho,\kappa}^{-1}(\A\times  \Omega))}
\end{equation}
where  $\sigma_{\varrho,\kappa}$ is  the partial dilation  $\sigma_{\varrho,\kappa}(\varphi,I;\omega):= (\varphi, \varrho I; \kappa \omega)$ and the  H\"older norms are defined in Section \ref{Sec:ApprLemma}. Note that the function $\ell \to \|G_t\|_{\ell;\varrho,\kappa}$ is {\em increasing} in the interval $[0,+\infty)$ since the set $B(0,\rho_0)\times  \Omega$ is {\em convex}. 

Fix $\tau > n-1$ and $\kappa\in (0,1)$, and 
define   $\Omega_\kappa$ as  the set of all $(\kappa,\tau)$-Diophantine vectors  satisfying \eqref{eq:sdc} and such that ${\rm dist}\, (\omega,\T^{n-1}\setminus \Omega)\ge \kappa$, i.e.
\begin{equation}
\label{eq:sdc3}
\Omega_\kappa =  D(\kappa,\tau)\cap \overline{\Omega -\kappa}.   
\end{equation}
Recall that $\vartheta_0$, $\ell_0$ and $\ell(m)$ are defined in \eqref{eq:ell(m)}. The following statement is a counterpart of Theorem \ref{Theo:A} in the case of exact symplectic mappings. 
\begin{Theorem}\label{Theo:A-for-maps}
 There is a positive constant $\epsilon= \epsilon(n,\tau,\vartheta_0)$  depending only  on $n$, $\tau$ and on $\vartheta_0$ such that for any $\delta>0$, $0<\kappa<1$ and $0<\varrho \le \rho_0$ the following holds.\\
 
\noindent 
Let $P_t$ be a  $C^k$ family of exact symplectic maps with generating functions  $\widetilde G_t$ satisfying \eqref{eq:generating-function-of P} and  \eqref{eq:families-of-perturbations-maps}   and such that 
\begin{equation}
   \sup_{t\in [0,\delta]} \, \| G_t\|_{\ell_0+1;\varrho,\kappa} \  \le \ \epsilon \varrho \kappa  \, .
\label{eq:smallness-condition-maps}
\end{equation}
Then there exists a $C^k$ family of maps 
\[
[0,\delta]\ni t \mapsto \phi_t\in C^\infty(\Omega;\Omega), \quad \ [0,\delta]\ni t \mapsto \Psi_t=(U_t,V_t)\in C^\infty(\T^{n-1} \times \Omega;\T^{n-1}\times B^{n-1}(0,\varrho)) \ 
\] 
such that   
\begin{enumerate}
\item [(i)]  ${\rm supp}\, (\phi_t-{\rm id})\subset \Omega-\kappa/2$, ${\rm supp}\, \big((U_t,V_t)-({\rm id},0)\big) \subset  \T^{n-1}\times (\Omega-\kappa/2)$; 
\item[(ii)]  
For each $\omega \in \Omega_\kappa$ and $t\in [0,\delta]$ the map 
$\Psi_{t,\omega} := \Psi_t(\cdot,\omega):{\T}^{n-1}  \rightarrow
 {\T}^{n-1}\times B(0,\varrho)$ 
is a smooth embedding, 
$\Lambda_t(\omega)  := \Psi_{t,\omega}({\T}^{n-1})$ is an embedded Lagrangian 
torus invariant  with respect to the exact symplectic map given by 
$ P_{t,\phi_t(\omega)}(\theta,I) := P(\theta, I; \phi_t(\omega),t)$, and 
\[
P_{t,\phi_t(\omega)} \circ \Psi_{t,\omega} = 
\Psi_{t,\omega} \circ R_{\omega} \quad \mbox{on}\  {\T}^n ;  
\]
\item[(iii)]
For any $m\in \N$ there is $C_m>0$ depending only on $m$, $n$, $\tau$ and $\vartheta_0$ such that
for any $\alpha,\beta\in \N^n$ with  $|\alpha|+|\beta|(\tau+1)\le m(\tau+1) +1$  the following estimate holds
\begin{equation}\label{eq:estimates-F}
\begin{array}{lrc}
\left|\partial_\theta^{\alpha}(\kappa \partial_\omega)^{\beta}(U_t(\theta;\omega) - \theta)\right|\,  + \,
\varrho ^{-1}
\left|\partial_\theta^{\alpha}(\kappa \partial_\omega)^{\beta} V_t(\theta;\omega)\right|  
+ \kappa^{-1} \left|(\kappa \partial_\omega)^{\beta}(\phi_t(\omega) - \omega)\right|
\\    [0.5cm]  
\displaystyle \leq\    C_{m}\, 
 (\varrho \kappa   )^{-1}\, \,\| G_t\|_{\ell(m)+1;\varrho,\kappa} 
\end{array}
\end{equation}
 uniformly in  $(\theta,\omega,t)\in {\T}^{n-1}\times \Omega_\kappa\times [0,\delta]$. 
\end{enumerate}
If $P_t$ is analytic with respect to $t$ in a disc $B(0,a)$ then so are $U$, $V$ and $\phi$.   
\end{Theorem}	
{\em Proof.} To simplify the notations we fix $k=1$. The first step in the proof will be to modify $P_t$ by multiplying $G_t$ by suitable cut-off functions without changing the corresponding estimates. This allows us to work with compactly supported functions. 
To this end we use the cut-off functions $\widetilde\psi_{\varrho}$ and $\psi_{\kappa}$ constructed in the previous sub-section by means of  Lemma \ref{lemma:cut-off-omega}. We set $G^1_t(\theta,I;\omega):= G_t (\theta,I;\omega) \widetilde\psi_{\varrho}(I) \psi_{\kappa}(\omega)$ and denote by $P_t^1$ the exact symplectic map with generating function 
\[
\widetilde G_t^1(\theta,r;\omega):= \langle \theta,r\rangle -\langle \omega,r\rangle -G_t^1(\theta,r;\omega). 
\]
The function $\ell \to \|G_t\|_{\ell;\varrho,\kappa}$ is {\em increasing} in the interval $[0,+\infty)$ since the set $B(0,\rho_0)\times  \Omega$ is {\em convex}. Then using \eqref{eq:Leibnitz1} one obtains
\[
\|G_t^1\|_{\ell;\varrho,\kappa} \le C_\ell \|G_t\|_{\ell;\varrho,\kappa} ,
\]
where $C_\ell$ depends only on $n$ and $\ell$. In particular, 
it follows from \eqref{eq:smallness-condition-maps} that
\[
\|\sigma_\varrho^{-1}{\rm sgrad}\,  G_t^1\|_{1;\varrho,\kappa} \le C_1\|\sigma_\varrho^{-1}{\rm sgrad}\,  G_t\|_{\ell_0-1;\varrho,\kappa} \le C_1\epsilon \kappa
\]
which allows one to apply Lemma \ref{rem:with-parameters}. Hereafter
\[{\rm sgrad\, }G_t^1(\theta,r;\omega):= (\nabla_r G_t^1(\theta,r;\omega), -\nabla_\theta  G_t^1(\theta,r;\omega))\]
is the simplectic gradient of $G_t^1$. 
We have 
$$
P_t^1(\theta,r;\omega)=(\theta+\omega,r) \quad \mbox{for} \quad (r,\omega)\notin  B^{n-1}(0,7\varrho/8)\times (\Omega-\kappa/4).
$$
Moreover, 
$$
P_t^1(\theta,r;\omega)= P_t(\theta,r;\omega) \quad \mbox{for} \quad (\theta,r;\omega)\in \T^{n-1}\times B^{n-1}(0,\varrho/2)\times (\Omega-\kappa)
$$ 
since $G_t^1=G_t$ on that set.  
This allows us to replace $G_t$ by $G_t^1$ and $P_t$ by $P_t^1$ in the theorem. 
From now on we suppose that 
\[
{\rm supp}\, G_t\subset \A':=\T^{n-1}\times B^{n-1}(0,\varrho)\times (\Omega-\kappa/4). 
\]
Set 
\[
{\mathcal N}_{\omega,\omega_n}(r,r_n):=  \langle \omega, r\rangle + \omega_nr_n, \quad (\omega,\omega_n)\in \Omega\times I, \quad I=(\pi, 3\pi),
\] 
Using an argument of Douady \cite{Dou} (see also Theorem 1.1 \cite{P2} and Theorem 3.1 \cite{P3}), 
we are going  to find a $C^1$ family of  Hamiltonians 
\begin{equation}\label{eq:inclusion-map-hamiltonian}
(\theta,\theta_n,r,r_n)\to H_t(\theta,\theta_n,r,r_n; \omega,\omega_n) =  \langle \omega, r\rangle + \omega_nr_n + h_t(\theta,\theta_n,r;\omega)
\end{equation}
 in $\T^n\times\R^n$ depending smoothly on parameters $(\omega,\omega_n)\in \Omega\times I$ with the following properties - $h_t$ is ``small'' and   the corresponding   Poincar\'e map is given by $P_t$ at any energy surface. Then we will apply Theorem \ref{Theo:A} to the family $H_t$.   
We set $y=(y',y_n)\in \R^n$, $(\theta,\theta_n)=\mathrm{pr}(y)\in \T^n$,  $\eta=(\eta',\eta_n)=(r,r_n)$ and $\widetilde \omega:= (\omega,\omega_n)\in \widetilde \Omega:= \Omega\times I$. We shall denote by $\Sigma_c$ the energy surface $ \Sigma_c:=\{H_t=c\} \subset \T^n\times\R^n$ for $c\in\R$ and by $\imath_c:\T^{n-1}\times\R^{n-1}\to\ \Sigma_c\cap \{\theta_n=0\}$ the corresponding inclusion map, i.e. 
$$\imath_c(\theta,r)=(\theta,0,r,(c-\langle \omega,r\rangle- h_t(\theta,0,r;\omega))/\omega_n).$$

\begin{Prop}\label{Prop:H}
There is a $C^1$ family of  Hamiltonians 
$$
H_t(\theta,\eta;\widetilde \omega):=  {\mathcal N}_{\widetilde \omega}(\eta) + h_t(\theta,\eta';\omega)\, ,\quad h_t\in C^\infty (\T^n\times\R^{n-1}\times \Omega)\, ,\  \widetilde \omega\in \Omega\times (\pi, 3\pi),
$$  
such that 
\begin{itemize}
	\item[(i)] $P_t = \imath^{-1}_{c}\circ \Phi_t^{\frac{2\pi}{\omega_n}}\circ\imath_{c}$ for any $c\in\R$, where 
	$\imath_c$ is the corresponding inclusion map 
	and $\Phi_t^{x}$,  $x\in\R$, is the Hamiltonian flow of $H_t$, 
	\item[(ii)] $\mathrm{supp\, } (h_t) \subset A'':= \T^n\times B^{n-1}(0,\varrho)\times \Omega$ and 
	the following estimate holds 
	\[
	\|h_t\|_{\ell,\A'';\varrho,\kappa} \le C_{\ell}  \|G_t\|_{\ell+1,\A' ;\varrho,\kappa} 
\]	 
	for any  $\ell\ge 2$ and $0<\varrho,\kappa\le 1$, where $\T^n\times\R^{n-1}\times \Omega$,  $C_\ell>0$ depends only on $\ell$ and $n$. 
	\item[(iii)] $h_t(\theta,\eta';\omega)= 0$ for $|\theta_n|\le \pi/2$.  
\end{itemize}
\end{Prop}
{\em Proof.}
Choose $\eta\in C^\infty(\R)$ such that $0\le \eta \le 1$, $\eta(s)=0$ for $|s|\le 1/4$ and $\eta(s)=1$ for $|s|\ge 1/2$.
 Consider the family of exact symplectic maps $s\to P_t^s$ in $\A=\T^{n-1}\times\R^{n-1}$ for $s\in \R$ having  generating functions of the form
\[
\widetilde G_t^s(\theta,r;\omega):= \langle \theta,r\rangle - s\langle \omega,r\rangle   - \eta(s)G_t(\theta,r;\omega),\ (\theta,r;\omega)\in \R^{n-1}\times\R^{n-1}\times \Omega.
\] 
We have $P_t^s=Q^s\circ W_t^s$, where $Q^s(\theta,r;\omega)=(\theta +s\omega,r)$ and the generating function of $W_t^s$ is 
$ \langle \theta,r\rangle   - \eta(s)G_t(\theta,r;\omega)$. 
Set $G(\theta,r;\omega,t)=G_t(\theta,r;\omega)$ and denote the symplectic gradient of $G$ with respect to $(\theta,r)$ by
\[
{\rm sgrad\, }G(\theta,r;\omega,t):= (\nabla_r G(\theta,r;\omega,t), -\nabla_\theta  G(\theta,r;\omega,t)).
\]
Notice that
\begin{equation}\label{eq:P-s}
\begin{array}{rcll}
P_t^s(\theta,r;\omega) &=& (\theta +s\omega,r)\quad \forall\, s\in [-1/4,-1/4]\\[0.3cm]
P_t^s(\theta,r;\omega) &=& P_t(\theta,r;\omega) \quad \forall\, s\in (-\infty,-1/2]\cup [1/2,+\infty) .
\end{array}
\end{equation}
Denote by 
$$\xi_t^s \, := \, \frac{d P_t^s}{ds}  \circ (P_t^s)^{-1}, \quad s\in \R$$  
the corresponding  vector field and set $v(\theta,s,r;\omega) =(\omega,0)$. Then
\[
\xi_t^s(\theta,r;\omega) = v(\theta,s,r;\omega) = (\omega, 0)\quad \forall\, s\in (-\infty,-1/2]\cup [-1/4,1/4]\cup [1/2,+\infty) .
\]
We set $\xi_t(\theta,s,r;\omega)=\xi_t^s(\theta,r)$ where $(\theta, s)\in\T^{n-1}\times [0,1]$. 
\begin{Lemma}\label{lemma:zeta}
We have $\xi_t^s(\theta,r;\omega) = (\omega, 0)$ for $s\in [0,1/4]\cup [1/2,1]$ and  
$$\mathrm{supp\, }(\xi_t-v) \subset \A':= \T^{n-1}\times  (0 , 1)  \times B^{n-1}(0,\varrho)\times \Omega.$$ 
Moreover, the following estimates hold
\[
 \|\sigma_\varrho^{-1}(\xi_t-v)\|_{\ell, \A';\rho,\kappa} \le  C_{\ell} \, 
\|\sigma_\varrho^{-1}{\rm sgrad\, }G_t(\cdot;\omega)\|_{\ell,\A';\rho, \kappa} 
\]
where $\ell\ge 1$ and $C_{\ell}>0$ depends only on $\ell$ and $n$.
\end{Lemma}
{\em Proof.}    We have $\xi_t^s =v + dW_t^s/ds \circ (W_t^s)^{-1}$. By  Lemma \ref{lemma:generating-functions} the support of $W_t^s -\mathrm{id}$ is contained in $\A'$. 
The estimates of $dW_t^s/ds$ and $(W_t^s)^{-1}$ follow from Proposition \ref{prop:composition-holder1}, Lemma \ref{lemma:generating-functions} and Lemma \ref{rem:with-parameters} taking into account \eqref{eq:smallness-condition-maps}. To obtain the corresponding estimates for the composition $dW_t^s/ds\circ (W_t^s)^{-1}$ we use Proposition \ref{prop:composition-holder1}.   \finishproof

\noindent
Notice  that the one-form $\imath(\xi_t^s)d\theta\wedge dr$ is exact, where $d\theta\wedge dr=\sum_{j=1}^{n-1}d\theta_j\wedge dr_j$ is the standard symplectic two-form in $T^\ast \T^{n-1}$ and $\imath(v)$ is the inner product with the vector field $v$. This follows from the identity 
$$\imath(\xi_t^s)d\theta\wedge dr=(P_t^s)_\ast(\frac{d}{ds}\, (P_t^s)^\ast ( r d\theta)-d(\imath(\xi_t^s)r d\theta)) $$
 since $(P_t^s)^\ast ( r d\theta) -  r d\theta$ is exact. 
Let $h_t^s$ be a primitive of 
\begin{equation}\label{eq:primitive-h}
\left\{
\begin{array}{lcrr}
\imath(\xi_t^s-v)d\theta\wedge dr= - dh_t^s\\[0.3cm]
h_t^s(0,0;\omega)=	0.
\end{array}
\right.
\end{equation}
The first equality in \eqref{eq:primitive-h} means that
\[
\begin{pmatrix}
\partial_{\varphi}h_t^s\cr
\partial_{r}h_t^s
\end{pmatrix}
=
\begin{pmatrix}
0 & I\cr
-I & 0
\end{pmatrix}
(\xi_t^s -v) .
\]
The second one combined with  Lemma \ref{lemma:zeta}   implies that $\mathrm{supp\, }(h_t) \subset \A'$ and
\[
\begin{array}{lcrr}
\displaystyle \|h_t\|_{\ell,\A';\varrho,\kappa}  \le C \big(\|\partial_\theta h_t\|_{\ell,\A';\varrho,\kappa} + \|(\varrho \partial_I )h_t\|_{\ell,\A';\varrho,\kappa}  \\[0.3cm]
\displaystyle \le\,  C_\ell^\prime\big(\|\partial_\theta h_t\|_{\ell,\A';\varrho,\kappa} + \|(\varrho \partial_I )h_t\|_{\ell-1,\A';\varrho,\kappa}\big)\\[0.3cm]
\displaystyle \le  C_{\ell} \|G_t\|_{\ell+1,\A';\varrho,\kappa} .
\end{array}
\]

We denote as above
$y=(y',y_n)\in \R^n$, where $y_n = 2\pi s$, $(\theta,\theta_n)= \mbox{pr}(y) \in \T^n$, $\eta=(\eta',\eta_n)=(r,r_n)$,   and set 
\[
h_t(y,\eta';\omega):=h_t^{y_n/2\pi}(y',\eta';\omega),  \ H_t(y,\eta;\omega,\omega_n):= \langle \omega, \eta'\rangle+ \omega_n \eta_n + h_t(y,\eta';\omega). 
\]
Let 
${\mathcal N}(\eta;\omega,\omega_n)= \langle \omega,\eta'\rangle + \omega_n\eta_n$ be the corresponding normal form. 
The  H\"{o}lder norms of $H_t-{\mathcal N}=h_t$ have been estimated above which proves (ii).

Notice that $
h_t^s(\theta,r)= 0$ for  $s\in [0,1/4]\cup [1/2,1]$ 
since $\xi_t^s(\theta,r;\omega)=(\omega,0)$ there. 
Then 
\[
h_t(\theta,\eta_n,r,\omega)= 0 \quad \forall\, \theta_n \in  [0,\pi/2]\cup [\pi,2 \pi]
\]
which gives (iii). To obtain (i), consider the Hamiltonian vector field $X_{H_t}$ given by 
\[
\imath(X_{H_t})dy\wedge d\eta= - dH_t(y,\eta).
\]
By \eqref{eq:primitive-h} we get 
\[
(\partial_r H_t(\theta,y_n,r,r_n), - \partial_\theta H_t(\theta,y_n,r,r_n)) = (\partial_r h_t^{y_n/2\pi}(\theta,r)+\omega, - \partial_\theta h_t^{y_n/2\pi}(\theta,r)) = \xi_t^{y_n/2\pi}(\theta,r).
\]
Moreover, $\partial_{r_n} H_t(\theta,y_n,r,r_n)=\omega_n$. Setting $P_t^s(\theta,r)= (p_t^s(\theta,r),  q_t^s(\theta,r))$ we obtain that the Hamiltonian flow $\Phi_t^s$ of $H_t$ is given by 
\[
\Phi_t^x(\theta,y_n,r,r_n)= (p_t^{x\omega_n/2\pi}(\theta,r),y_n+x\omega_n, q_t^{x\omega_n/2\pi}(\theta,r), q_n^{x\omega_n/2\pi}(\theta,y_n,r,r_n))
\]
for $x\in \R$, where $t\to q_n$ is a $C^1$ family of smooth functions. 
This yields
\[
\Phi_t^{\frac{2\pi}{\omega_n}}(\theta,0,r,r_n)= (p_t^{1}(\theta,r),2\pi, q_t^{1}(\theta,r), q_n^{1}(\theta,0,r,r_n)). 
\]
This implies (i) and  completes the proof of the proposition. \finishproof

In particular, we have
\[
\|h_t\|_{\ell_0,\A'';\varrho,\kappa} \le C \epsilon \kappa \varrho
\]
for $t\in [0,a]$, where $C>0$ depends only on $\ell_0$ and $n$. Then 
Proposition  \ref{Prop:H} allows us to apply Theorem \ref{Theo:A} to the perturbation $H_t$ of ${\mathcal N}$ in 
$\T^n\times  B^n(0,\varrho)\times \tilde \Omega$, where $\tilde \Omega:=\Omega\times (\pi,3\pi)$.   Denote by 
\[
\tilde \Psi_t=(\tilde U_t,\tilde V_t,\tilde \phi_t): \T^n\times  B^n(0,\varrho)\times \tilde \Omega\to \T^n\times  B^n(0,\varrho)\times \tilde \Omega
\] 
the map given by Theorem \ref{Theo:A}. 
To obtain the map $\Psi_t$ from $\tilde \Psi_t$ we use an argument of  Douady \cite{Dou}.  Fix $\omega\in \Omega_\kappa$ and $t\in [0,\delta]$. 
Then $\tilde \omega:=(\omega,2\pi)\in \tilde\Omega_\kappa:=\widetilde D(\kappa,\tau)\cap (\tilde \Omega-\kappa)$, where $\widetilde D(\kappa,\tau)$ is given by \eqref{eq:sdc2}. By  Theorem \ref{Theo:A} there exists  a $C^1$ family of Kronecker invariant tori $\tilde\Lambda_t(\tilde\omega)$  of the flow $\Phi^s_{t}(\cdot; \tilde \phi_t(\tilde\omega))$ of the Hamiltonian $H_t(\cdot; \tilde \phi_t(\tilde\omega))$ with a frequency  vector $\tilde\omega$ lying on a certain energy surfaces  $\Sigma_{t}:=\{  H_t = c_t\}$. Moreover, 
\begin{equation}
\left(\tilde U_t(y + x \tilde\omega; \tilde \omega),\tilde V_t(y+x \tilde\omega; \tilde\omega)\right) \, =\, \Phi^x_{t}\left(\tilde U_t(y,\tilde\omega), \tilde V_t(y, \tilde\omega);\tilde \phi_t(\tilde\omega)\right)\,  
\label{eq:g}
\end{equation}
for any $x\in \R$ and $y\in \T^n$. On the other hand,  $\Phi^x_{t}$ is of the form
\[
\Phi^x_{t}(y',y_n,\eta;\tilde \phi_t(\tilde\omega))= (z_t'(x,y,\eta;\tilde \omega),y_n + x\omega_n, \zeta_t(x,y,\eta;\tilde \omega))
\]
since $\partial H_t/\partial \eta_n=  \omega_n$, hence,   the last coordinate 
$\tilde U_{t,n}$ of $\tilde U_t$ satisfies the equality
$$\tilde U_{t,n}(y + x \tilde \omega;\tilde \omega)= \tilde U_{t,n}(y;\tilde \omega) +x\omega_n $$ 
for any $y\in \T^n$ and $x\in \R$. 
Then $\nabla_y\tilde U_{t,n}(y + x \tilde \omega;\tilde \omega)= \nabla_y\tilde U_{t,n}(y;\tilde \omega)$, and since the flow $x\mapsto x \tilde \omega$ is ergodic on $\T^n$ (recall that $\tilde \omega$ is a Diophantine frequency) we get 
$$\nabla_y\tilde U_{t,n}(y;\tilde \omega)=\nabla_y\tilde U_{t,n}(0;\tilde \omega):=a_t(\tilde \omega)$$
 for any $y\in\T^n$. The function $\T^n \ni y\mapsto U_{t,n}(y ;\tilde \omega)\in \T$ is  determined   up to a translation and we fix it by $U_{t,n}(0;\tilde \omega)=0$. Then    $\tilde U_{t,n}(y;\tilde \omega) = \langle a_t(\tilde \omega),y \rangle$ and $a_t(\tilde \omega)\in\Z^n$. On the other hand,   $|\nabla_y\tilde U_t(y;\tilde \omega)- \mathrm{Id}| \le C(n,\tau,\vartheta_0)\epsilon$ in view of \eqref{eq:main-estimates}  and \eqref{eq:smallness-condition3}. 
 Taking $\epsilon=\epsilon(n,\tau,\vartheta_0)$ small enough we obtain  that $a_t=(0,\ldots,0,1)$  and  $\tilde U_{t,n}(y; \tilde \omega) = y_n$. Now we set 
 \[
 \Psi_{t}(\theta,r) =(U_t(\theta,\omega),V_t(\theta,\omega) ) :=  \mathrm{p} (\tilde U_t(\theta,0,\tilde\omega), \tilde V_t(\theta,0,\tilde\omega)),\quad \phi_t(\omega)=\tilde \phi_t(\tilde\omega),\quad  \tilde\omega=(\omega, 2\pi), 
 \]
 where $\mathrm{p}$ is the projection given by $\mathrm{p}(\theta,\theta_n,r, r_n)= (\theta,r)$.

Using Proposition \ref{Prop:H}  one  obtains that for each $ \omega\in \Omega_\kappa$, the torus 
 \[
 \Lambda_t(\omega):= \mathrm{p} (\tilde\Lambda_t(\tilde\omega))  
 \]   
 is a Kronecker invariant torus of $P_{t,\phi_t(\omega)}$ with a frequency vector $\omega$ and  we obtain (ii). Moreover, \eqref{eq:estimates-tilde-F} implies \eqref{eq:estimates-F}. 
 To prove the analyticity with respect to $t$ we use Cauchy theorem at any step of the construction.
 This completes the proof of Theorem \ref{Theo:KAM}. \finishproof
\\

\subsection{KAM theorems for $C^k$ families of symplectic maps}\label{Subsec:KAM-maps}
The aim of this section is to obtain $C^k$ families of Kronecker invariant tori of $C^k$ families of exact  maps close to the family  $(\theta,r) \to (\theta +\nabla K_t(r),r)$, $t\in [0,\delta]$. 

Let $\Omega \subset [0,2\pi]^{n-1}$ be an open convex set which we  identify with an open convex subset of $\T^{n-1}$.
Fix $k\in\{0,1\}$ and consider  a $C^k$-family  of real-valued functions 
\[
[0,\delta]\ni t\to K_t^\ast\in C^\infty( \overline \Omega,\R)
\]  
 satisfying the non-degeneracy condition 
\begin{equation}\label{eq:non-degenerate1}
\nabla K_t^\ast : \Omega \longrightarrow  D_t:=\nabla K^\ast_t(\Omega) \ \mbox{is a diffeomorphism.}
\end{equation}
where $\overline \Omega$ is the closure of $\Omega$. 
For any $t\in [0,\delta]$ the Legendre transform $K_t:  D_t\to  \Omega$ of $K_t^\ast$ is in $C^\infty(  D_t,\R)$ and it satisfies the non-degenerate condition \eqref{eq:non-degenerate}. Moreover, the corresponding family of functions    $[0,\delta]\ni t\to K_t\in C^\infty( D_t,\R)$ is $C^k$ smooth. 
We set $\A_t:=\T^{n-1}\times D_t$ and   denote by
\begin{equation}\label{eq:Q}
Q_t:\A_t\to \A_t, \quad Q_t(\theta,r)= (\theta +\nabla K_t(r),r)
\end{equation}
the corresponding family of exact symplectic maps on $\A_t$.   The frequency vector of $Q_t$ on the invariant torus $\T^{n-1}\times \{r\}$ is $\omega=\nabla K_t(r)\in \Omega$.

We consider  a $C^k$-family   of exact symplectic maps 
\begin{equation}\label{eq:P-t-family-1}
[0,\delta]\ni t\to P_t \in C^\infty( \A_t, \T^{n-1}\times \R^{n-1})
\end{equation}
close $Q_t$.   We suppose  that  $P_t$ is defined by  a generating function $\widetilde G_t$ of the form
\begin{equation} 
\label{eq:generating-function-of-P-1}
\widetilde G_t(\theta,r):= \langle \theta,r\rangle - K_t(r) -G_t(\theta,r)
\end{equation}
which  means that 
\[
P_t(\theta-\nabla K_t(r)-\nabla_{r}G_t(\theta,r),r)=(\theta, r-\nabla_{\theta}G_t(\theta,r))
\]
for any $(\theta,r)\in \A_t$. We suppose as well that
the map
\[
[0,\delta] \ni t\to G_t \in C^\infty (\T^{n-1}\times D_t, \R)
\]
is $C^k$ smooth with $k=0$ or $k=1$. 

We assume as well that the $C^2$ norm of $G_t$ is  sufficiently small. Then  the inverse function theorem implies that the map $\theta \to\theta-\nabla_{r}G_t(\theta,r)$ is a diffeomorphism of  $\T^{n-1}$ for any fixed $r\in D_t$ and $P_t$ is well defined. 

Denote by $R_t$ the exact symplectic map with generating function $(\theta,r)\to \langle \theta,r\rangle  -G_t(\theta,r)$, i.e. $R_t(\theta-\nabla_{r}G_t(\theta,r),r)=(\theta, r-\nabla_{\theta}G_t(\theta,r))$.  One can show that 
 $P_t=R_t\circ Q_t$ on $\A_t$ (see Lemma \ref{Lemma:decomposition}).

Given $\varrho,\kappa\in (0,1)$ and $m\ge 0$ we set
\begin{equation}
\label{eq:B-0-sequence}
{\mathcal B}_m^0 \, := \,   \sup_{0\le t\le \delta }\, \left( \varrho^2 |\!|\!| \partial^2 K_t |\!|\!|_{\ell(m)+1,  D_t;\kappa} \, +\,  |\!|\!|  G_t|\!|\!|_{\ell(m)+1,\A_t; \kappa}\right), 
\end{equation}
where  $\|G_t\|_{\ell, \A_t;\kappa} := \|G_t\circ \sigma_{\kappa}\|_{C^\ell(\sigma_{\kappa}^{-1}(\A_t))}$, 
 $\sigma_\kappa$ is defined by $\sigma_\kappa(x,\xi)=(x,\kappa \xi)$, and
\[
|\!|\!| u |\!|\!|_{\ell,D;\kappa}\, =\, \sup_{0\le m\le \ell}\|u\|_{\ell-m,D;\kappa}.
\]
 Recall that $\vartheta_0$, $\ell_0$ and $\ell(m)$ are defined in \eqref{eq:ell(m)}. The sequence ${\mathcal B}_m^0$, $m\ge 0$, is increasing by definition. 
 
  To formulate the smallness condition in the KAM theorem below we need as well  the notation 
\begin{equation}\label{eq:S-K-sequence}
S_{\ell}(\nabla K^{\ast}) := \sup_{0\le t\le \delta }\, \big(1+ \|\nabla K^{\ast}_t\|_{C^1(\Omega)}\big)^{\ell-1} \big(1+ \|\nabla K^{\ast}_t\|_{C^\ell(\Omega)}\big)
\end{equation}
introduced  in \eqref{eq:S-sequence}. This in an increasing sequence with respect to $\ell\in [1,+\infty)$ since $\Omega$ is convex. 
For any $m\ge 0$ we set
\begin{equation}
\label{eq:B-sequence}
{\mathcal B}_m \, := \,   {\mathcal B}_m^0 S_{\ell(m)+2}(\nabla K^{\ast}). 
\end{equation}

\begin{Theorem}\label{Theo:KAM} There exists $\epsilon = \epsilon(n,\tau,\vartheta_0)>0$ depending only on $n$, $\tau$ and $\vartheta_0$ such that 
the following holds.\\

\noindent
Let  $\Omega\subset \T^{n-1}$ be an open convex  set,  $[0,\delta]\ni t\to K_t^\ast\in C^\infty( \overline \Omega,\R)$  a $C^k$ family  satisfying \eqref{eq:non-degenerate1} and let $\varrho>0$ and $\kappa>0$ be such that $0<\varrho\le \kappa\le 1$ and $\Omega_\kappa\neq\emptyset$. 
Consider a  $C^k$ family  of  exact symplectic maps  
$[0,\delta]\ni P_t\in C^\infty( \A_t,  \T^{n-1}\times\R^{n-1})$  with generating functions $\widetilde G_t$ of the form \eqref{eq:generating-function-of-P-1} such that
\begin{equation}\label{eq:smallness-condition-maps-2}
{\mathcal B}_0^0S_{\ell_0+1}(\nabla K^{\ast})\, =\,   \sup_{0\le t\le \delta }\, \left( \varrho^2 |\!|\!| \partial^2 K_t |\!|\!|_{\ell_0+1,  D_t;\kappa} \, +\,  
|\!|\!| G_t|\!|\!|_{\ell_0+1,\A_t; \kappa}\right)\, 
S_{\ell_0+1}(\nabla K^{\ast})\ 
\le \ \epsilon \varrho \kappa .
\end{equation}
Then there is a $C^k$ family $[0,\delta]\ni t\to    f_t=( u_t, v_t)\in C^\infty(\T^{n-1}\times \Omega; \T^{n-1}\times  D_t)$ such that 
\begin{enumerate}
\item[(i)]
for any $\omega\in \Omega_{\kappa}$,  
$[0,\delta]\ni s\to \Lambda_t(\omega)=  f_{t,\omega}(\T^{n-1})$ is a $C^k $  family of Kronecker invariant tori 
of of $P_t$ with a frequency vector $\omega$, where $  f_{t,\omega}:=  f_t(\cdot; \omega)$, and 
the following diagram is commutative
\[
\displaystyle{\begin{array}{cccl} 
\displaystyle \T^{n-1}&\stackrel{R_{\omega}}{\longrightarrow}&\T^{n-1}\cr
\downarrow\lefteqn{ f_{t,\omega}}& &\downarrow\lefteqn{  f_{t,\omega}} \cr
\displaystyle \Lambda_t(\omega)&\stackrel{P_t}{\longrightarrow}&\Lambda_t(\omega)&  
\end{array} }
\]
\item[(ii)] for any $m\in \{0 \}\cup [1,+\infty)$ the following estimates hold 
\begin{equation}\label{eq:estimates-F-1}
\begin{array}{lcrr}
\displaystyle \big|\partial_\varphi^\alpha (\kappa\partial_\omega)^\beta \big(u_t(\varphi;\omega)-\varphi)\big)\big| \, \le\,  C_m \frac{{\mathcal B}_m}{\kappa\varrho}\\[0.3cm]
\displaystyle
  \big|\partial_\varphi^\alpha (\kappa\partial_\omega)^\beta   \big(v_t(\varphi;\omega)-\nabla K_t^{\ast} (\omega)\big)\big| 
\, \le\,  C_m \frac{{\mathcal B}_m}{\varrho} \left(1+\frac{{\mathcal B}_1}{\varrho}\right)^{m} 
\end{array}
\end{equation}
for any $(\varphi,\omega)\in \T^{n-1}\times \Omega_\kappa$, $t\in [0,\delta]$,   and multi-indices $\alpha,\beta\in \N^{n-1}$ with  $|\alpha|+|\beta|(\tau+1)\le m(\tau+1)+1$, where  the constant $C_m>0$ depends only on $n$, $\tau$,  $\vartheta_0$ and $m$,  
\item[(iii)]
${\rm supp\, }\big((u_t,v_t) -({\rm id}, \nabla K_t^\ast \big) \subset \T^{n-1}\times (\Omega-\kappa/2)$.
\end{enumerate}
If $P$ is analytic with respect to $t$ in a disc $B(0,a)$ then so are $u$ and $v$.   
\end{Theorem}
We note that the $C^1$ families of invariant tori $\Lambda_t(\omega)$, $t\in [0,\delta]$, given by the theorem are uniquely defined. \\

\noindent
{\em Proof}. The proof is similar to that of Theorem \ref{Theo:KAM2} and we give only the main steps in it. We are going to apply Theorem \ref{Theo:A-for-maps}. To this end we will firstly construct the function $\widetilde G_t(\theta,r;\omega)$ in \eqref{eq:generating-function-of P}. \\

\noindent
{\em Step 1. Construction of the generating function $\widetilde G_t(\theta,r;\omega)$.}
Given $\omega\in \Omega -\kappa/4$ 
and $I$ the ball $B^{n-1}(0,\varrho)$,   we set 
$r=\nabla K_t^{\ast}(\omega) + I$. 
Choosing $\epsilon=\epsilon(n,\tau,\vartheta_0) < 1/9$ in  \eqref{eq:smallness-condition-maps-2} we obtain as in the proof of of Theorem \ref{Theo:KAM2} the following relation
\begin{equation} \label{eq:domain-G}
\forall\,  \omega\in \Omega-\kappa/4,\ \forall\, I\in B^{n-1}(0,\varrho), \quad  \left\{
\begin{array}{lcrr}  
\nabla K^{\ast}_t(\omega) + I\in D_t \quad \mbox{and}\\[0.3cm] 
\displaystyle \nabla  K_t(\nabla K^{\ast}_t(\omega) + I)\in \Omega-\frac{1}{8}\kappa  . 
\end{array}
\right.
\end{equation}
By Taylor's formula we obtain
\begin{equation}\label{eq:generating-function-G}
\widetilde G_t(\theta,r) = \langle \theta, \nabla K_t^\ast(\omega)\rangle + \langle \theta, I\rangle - K_t(\nabla K_t^\ast(\omega))  - \langle \omega, I\rangle - G_t^0(I;\omega) - G_t^1(\theta,I;\omega),
\end{equation}
where
\[
G_t^0(I;\omega) = \int_0^1(1-s)\langle \partial^2K_t(\nabla K_t^{\ast}(\omega) + sI)I,I\rangle\,  ds 
\]
and 
\[
G_t^1(\theta,I;\omega):= G_t(\theta, \nabla K^{\ast}_t(\omega) + I). 
\]
It follows from \eqref{eq:domain-G} that the functions  $G_t^0$ and $G_t^1$ are well defined  for $I\in B^{n-1}(0,\varrho)$ and $\omega\in \Omega-\kappa/4$.
 Denote by $P_{t,\omega}:\A_t\to \T^{n-1}\times \R^{n-1}$ the exact symplectic map defined by means of the  generating function 
\[
\begin{array}{rcll}
(\theta,I)\to \widetilde G_{t,\omega}(\theta,I)&:=& \langle \theta, I\rangle - \langle \omega, I\rangle - G_t^0(I;\omega) - G_t^1(\theta,I;\omega)\\[0.3cm]
 &=& \widetilde G_t(\theta,\nabla K_t^\ast(\omega)+I) -\langle \theta, \nabla K_t^\ast(\omega)\rangle - K_t(\nabla K_t^\ast(\omega))
\end{array}
\]
and set $\psi_\omega(\theta,I)= (\theta, \nabla K_t^\ast(\omega)+ I)$. 
\begin{Lemma}\label{lemma:equality-for-P}  For any $\omega\in \Omega-\kappa/4$ the map $P_{t,\omega}:\A_t\to \T^{n-1}\times \R^{n-1}$ is well defined and 
\[
 P_{t,\omega} = \psi_\omega^{-1}\circ P_t \circ \psi_\omega
\]
on $\T^{n-1}\times B^{n-1}(0,\varrho)$ provided that the constant  $\epsilon = \epsilon(n,\tau,\vartheta_0)$ in \eqref{eq:smallness-condition-maps-2} is sufficiently small. 
\end{Lemma}
{\em Proof}. The smallness condition \eqref{eq:smallness-condition-maps-2} implies that   
$\|\nabla_\theta \nabla_I \widetilde G_{t,\omega}(\theta,I)-{\rm Id}_n\| \le C\epsilon$  for $(\theta,I;\omega)\in \T^{n-1}\times B^{n-1}(0,\varrho)\times (\Omega-\kappa/4)$, where $C=C(n)$ depends only on the dimension $n$ and ${\rm Id}_n\in M_n(\R^{n-1})$ is the identity matrix. Choosing   $\epsilon = \epsilon(n,\tau,\vartheta_0)$ sufficiently small we obtain that the map $\theta\to \varphi=\nabla_I \widetilde G_{t,\omega}(\theta,I)$ is a diffeomorphism on $\T^{n-1}$ for any fixed $I\in B^{n-1}(0,\varrho)$ and $\omega\in \Omega-\kappa/4$, hence, $ P_{t,\omega}$ is well-defined. Notice that
$\nabla_I \widetilde G_{t,\omega}(\theta,I)=\nabla_r \widetilde G_t(\theta,\nabla K_t^\ast(\omega)+I)$ and $\nabla_\theta \widetilde G_{t,\omega}(\theta,I)=\nabla_\theta \widetilde G_t(\theta,\nabla K_t^\ast(\omega)+I)- \nabla K_t^\ast(\omega)$. Then
\[
 (\psi_\omega \circ P_{t,\omega})\big(\nabla_r \widetilde G_t(\theta,\nabla K_t^\ast(\omega)+I),I\big)=\big(\theta, \nabla_{\theta}\widetilde G_t(\theta,\nabla K_t^\ast(\omega)+I)\big).
\]
On the other hand,
\[
\begin{array}{rcll}
( P_t\circ\psi_\omega) \big(\nabla_r \widetilde G_t(\theta,\nabla K_t^\ast(\omega)+I),I\big)&=&P_t\big(\nabla_r \widetilde G_t(\theta,\nabla K_t^\ast(\omega)+I),\nabla K_t^\ast(\omega)+I\big)\\[0.3cm]
&=&\big(\theta, \nabla_{\theta}\widetilde G_t(\theta,\nabla K_t^\ast(\omega)+I)\big)
\end{array}
\]
and we obtain that $ \psi_\omega \circ P_{t,\omega} =  P_t \circ \psi_\omega$ since the map $\theta\to \nabla_r \widetilde G_t(\theta,r)$ is a diffeomorphism. \finishproof
\\

\noindent
{\em Step 2. H\"older estimates of $\widetilde G_t$.}
We are going to apply Theorem \ref{Theo:A-for-maps} to the family of exact symplectic maps $ P_t(\cdot;\omega):=\widetilde P_{t,\omega}(\cdot)$. 
To this end  we evaluate  the weighted norms of  $G_t^0$ and $G_t^1$. We have  $G_t^0=Q_t^0\circ ({\rm id}, \nabla K_t^{\ast})$, where
\[
Q_t^0(I;r)=\int_0^1(1-s)\langle \partial^2K_t(r+sI)I,I\rangle\, ds 
\]
is well defined and smooth in $\overline\Gamma$, $\Gamma:= B^{n-1}(0,\varrho)\times (\Omega-\kappa/4)$. 
As in the proof of Theorem \ref{Theo:KAM2} we obtain that 
\[
 |\!|\!|  Q_t^0 |\!|\!| _{\ell,\Gamma;\varrho,\kappa} \le  C_{\ell} \, \varrho^2 |\!|\!| \partial^2 K_t |\!|\!|_{\ell, D_t;\kappa} . 
\]
Since $\Omega$ is convex and $K^{\ast}_t\in C^\infty(\overline \Omega)$ using  Proposition \ref{prop:composition-holder1}, {\em 3.}, and Remark \ref{rem:anisotrop} as in the proof of Theorem \ref{Theo:KAM2}we obtain  for  any  $\ell\ge 1$ the estimate
\[
\|G_t^0\|_{\ell;\varrho,\kappa} \le 
C_{\ell} \, \varrho^2\, |\!|\!| \partial^2K_t |\!|\!|_{ \ell, D_t; \kappa}
S_{\ell}(\nabla K^{\ast}).  
\]
Moreover,  $\|G_t^1\|_{\ell;\varrho,\kappa} \le \|G_t^1\|_{\ell;\kappa,\kappa} \le C_\ell |\!|\!|G_t|\!|\!|_{\ell;\kappa}S_{\ell}(\nabla K^{\ast})$  since $0<\varrho\le \kappa\le 1$ and we obtain
\begin{equation}\label{eq:estimate-G-0-1}
\begin{array}{rcll}
\displaystyle\|G_t^0\|_{\ell;\varrho,\kappa} + \|G_t^1\|_{\ell;\varrho,\kappa}
\le \displaystyle C_{\ell} \left(\varrho^2\, |\!|\!| \partial^2 K_t |\!|\!|_{ \ell, D_t; \kappa}\, +\, |\!|\!| G_t |\!|\!|_{\ell,\A_t; \kappa}\right)\ 
S_{\ell}(\nabla K^{\ast}).
\end{array}
\end{equation}
\\

\noindent
{\em Step 3. Applying  Theorem \ref{Theo:A-for-maps}}. 
Now \eqref{eq:smallness-condition-maps-2} gives 
\[
\|G^0_t\|_{\ell_0+1;\varrho,\kappa} + \|G^1_t\|_{\ell_0+1;\varrho,\kappa} \le {\mathcal B}_0 S_{\ell_0+1}(\nabla K^{\ast})\le \epsilon \varrho \kappa .
\]  
This allows us to apply Theorem \ref{Theo:A-for-maps} to  the family of exact symplectic maps $ P_t(\cdot;\omega)=\widetilde P_{t,\omega}(\cdot)$. 
Set $ \Psi_t = (u_t, v_t)$, where $u_t= U_t$, $v_t=(\nabla K^{\ast}_t)\circ \phi_t +V_t$ and  $(U_t,V_t,\phi_t)$ are given by Theorem  \ref{Theo:A-for-maps}. Notice that $\|V_t\|_{C^0} \le c\epsilon \varrho$ in view of the estimates in (ii), Theorem \ref{Theo:A-for-maps}, where the constant $c$ depends only on  $n$, $\tau$ and $\vartheta_0$, and taking $\epsilon<\min(1,1/c)$ we obtain $V_t(\theta;\omega)\in B^{n-1}(0,\varrho)$ for any $\theta\in \T^{n-1}$ and $\omega\in \Omega-\kappa$. In the same way we get $\phi_t(\omega)\in \Omega-\kappa/4$ for $\omega\in \Omega-\kappa$. 

Now  Lemma \ref{lemma:equality-for-P} implies that for any $\omega\in \Omega_\kappa$ and $t\in [0,\delta]$ the Lagrangian manifold 
$ \Lambda_t(\omega):=  \Psi_{t}(\T^{n-1};\omega)$
is a Kronecker invariant torus of $P_t$ of a frequency vector $\omega$ satisfying  (i) in Theorem \ref{Theo:KAM}.  \\

\noindent
{\em Step 4. Estimates of $u_t$ and $v_t$.}

The estimates (ii), Theorem \ref{Theo:A-for-maps}, imply (ii) in Theorem \ref{Theo:KAM}. To estimate the derivatives of 
\[
v_t- \nabla K^{\ast}_t = V_t+ (\nabla K^{\ast}_t)\circ \phi_t - \nabla K^{\ast}_t
\]
we use (ii), Theorem \ref{Theo:A-for-maps} 
and the following Lemma which is an analogue of Lemma \ref{lemma:estimates-V}.
\begin{Lemma}\label{lemma:estimates-v}
For any $m\in \N$ the following estimate holds
\[
\|(\nabla K^{\ast}_t)\circ \phi_t - \nabla K^{\ast}_t\|_{m,\Omega;\kappa} \le C_{m} \frac{{\mathcal B}_m}{\varrho}\left(1+\frac{{\mathcal B}_1}{\varrho} \right)^m.
\]
\end{Lemma}
To prove (iii) we use suitable cut-off functions in $\omega$ given by Lemma \ref{lemma:cut-off-omega}. 
\finishproof

\subsection{Birkhoff Normal Forms for $C^k$-families of symplectic maps}\label{Subsec:BNF-for-maps}
Let $\Omega \subset [0,2\pi]^{n-1}$ be an open convex set which we  identify with an open convex subset of $\T^{n-1}$. Fix $\tau>n-1$ and denote by $\kappa_0(\Omega)$ the supremum of all $0<\kappa\le 1$ such that the set $\Omega_\kappa= D(\kappa,\tau)\cap \overline{\Omega-\kappa}$ is of  positive Lebesgue measure. Given $0<\kappa<\kappa_0(\Omega)$ we denote by $\Omega_\kappa^0$ the set of points of $\Omega_\kappa$ of {\em positive Lebesgue density}. Recall that $\omega\in \Omega_\kappa^0$ if  the Lebesgue measure of $\Omega_\kappa\cap U$ is positive for any neighborhood $U$ of $\omega$. Then $\Omega_\kappa\setminus\Omega_\kappa^0$ is a set of measure zero.  Recall that $\ell_0$ and $\ell(m)$ defined  in \eqref{eq:ell(m)}, i. e. $\ell_0:=2\tau + 2 + 2\vartheta_0$ and $\ell(m):=  2m(\tau + 1)+ \ell_0$ for $m \ge 0$. 

We are going to use as well the notations ${\mathcal B}_{m}$  and $S_{\ell}(\nabla K_t^{\ast})$  introduced  in  \eqref{eq:B-sequence}  and  \eqref{eq:S-K-sequence}. 

To construct a BNF we have to deal with the second differential $d^2K_t$ of the Legendre transform $K_t$ of $K_t^\ast$. 
We denote by $\partial^2 K_t(I)$ its Hessian matrix. Its norm  could be very large, as in the case of the  billiard ball map close to the boundary,  and to measure it we introduce a parameter $\lambda\ge 1$. More precisely, we suppose below that 
\begin{equation}\label{eq:non-degeneracy-BNF}
  \sup_{t\in [0,\delta]} \|\partial^2 K_t\|_{2,D_t; \kappa}\,  \le \, \lambda \, < \, \infty\, , 
\end{equation}
where $\lambda\ge 1$. 
 
\begin{Theorem}\label{Theo:BNF}  There exists $\epsilon = \epsilon(n,\tau,\vartheta_0)>0$ depending only on $n$, $\tau$, and $\vartheta_0$ such that 
the following holds. 
\\

\noindent
Let $\Omega\subset \T^{n-1}$ be an  open convex set and  $0<\varrho<\kappa< \kappa_0(\Omega)$.
Let    $[0,\delta]\ni t\to K_t^\ast\in C^\infty( \overline {\Omega},\R)$ be a $C^k$ family  satisfying the non-degeneracy condition \eqref{eq:non-degenerate1} and suppose that its Legendre transform $K_t$ satisfies \eqref{eq:non-degeneracy-BNF}.
Let  
$[0,\delta]\ni P_t\in C^\infty( \A_t,  \A)$,   be a  $C^k$ family  of  exact symplectic maps defined by generating functions $\widetilde G_t(\theta,r)=\langle \theta,r\rangle -K_t(r)-G_t(\theta,r)$
such  that  
\begin{equation}\label{eq:smallness-condition-1}
{\mathcal B}_2\ \le \ \epsilon \varrho \kappa\lambda^{-4} . 
\end{equation} 
Then 
\begin{enumerate}
\item[(i)] there exist  $C^k$-smooth with respect to $t\in [0,\delta]$ families of exact symplectic maps 
$\chi_t:\A_t \to \A_t$ and of real valued functions  $L_t\in C^\infty(D_t)$ and $R_t^0\in C^\infty(\A_t)$  such that 
\begin{enumerate}
\item $(\varphi,I)\mapsto \langle\varphi, I\rangle -L_t(I) - R_t^0(\varphi,I)$ is a generating function of $P_t^0:= \chi_t^{-1}\circ P_t\circ \chi_t$;
\item $\nabla L_t: D_t \to \Omega$ is a diffeomorphism, $L_t=K_t$  outside $D_t^1:= \nabla K_t^\ast(\Omega-\kappa/2)$, and  $\nabla L_t^\ast(\omega)=I_t(\omega)$ is given by \eqref{eq:momentum-I} for each $\omega\in \Omega_{ \kappa}^0$;
\item $R_t^0$ is flat at $\T^{n-1}\times \nabla L_t( \Omega_{ \kappa}^0$; 
\end{enumerate}
\item[(ii)] For any $t\in [0,\delta]$ and $m\in \N$ the following estimates hold
\begin{equation}\label{eq:estimates-chi}
\begin{array}{lcrr}
\|\sigma_\kappa^{-1}(\chi_t - {\rm id})\|_{m,\A_t;\kappa} + \|\sigma_\kappa^{-1}(\chi^{-1}_t - {\rm id})\|_{m,\A_t;\kappa}\\[0.3cm]
\displaystyle
\le \frac{C_m}{\varrho\kappa} \,  {\mathcal B}_{m+1} \lambda^{2m}  \big(\lambda +\|\partial^2 K_t\|_{m+1,D_t;\kappa}\big)
\end{array}
\end{equation}
and
\begin{equation}\label{eq:estimates}
\begin{array}{llcrr}
 \|\nabla L_t - \nabla K_t\|_{m,D_t;\kappa}   + \|\sigma_\kappa \nabla R_t^0\|_{m,D_t;\kappa}  \\[0.3cm]
 \le\,  \displaystyle \frac{C_m}{\varrho } \,  {\mathcal B}_{m+1}  \lambda^{2m}   \big(\lambda+\|\partial^2 K_t\|_{m+1,D_t; \kappa}\big)\,   , 
\end{array}
\end{equation}
where the constant $C_m>0$ depends only on $n$, $\tau$, $\vartheta_0$ and $m$.
\end{enumerate}
If $\widetilde G$ is analytic with respect to $t$, then so are $\chi$, $L$ and $R^0$.  
\end{Theorem}
Before proving the Theorem we observe that
\begin{Remark}\label{rem:inv-tori} (Birkhoff Normal Form) Let $k=1$. Then 
for any $\omega\in \Omega_\kappa^0$ the map 
\[
[0,\delta]\ni t\to \Lambda_t(\omega):=\chi_t(\T^{n-1}\times \{I_t(\omega)\})
\] 
provides  a $C^1$ family of invariant tori of $P_s$ with a frequency vector $\omega$ and taking into account Lemma \ref{Lemma:flat} we obtain
\[
P_t^0(\varphi,I)= (\varphi + \nabla L_t(I),I) + O_N(|I-\nabla L_t^\ast(\omega)|^N)
\] 
for any $N\in\N$ . Moreover,  the last formula can be differentiated $N-1$ times with respect to $(\varphi,I)$. Hence, Theorem \ref{Theo:BNF} gives a simultaneous Birkhoff Normal Form  of $P_t$ on the invariant tori $\Lambda_t(\omega)$, where $t\in [0,\delta]$ and the frequency vectors $\omega$ are in $\Omega_\kappa^0$. Recall as well that the complement of  $\Omega_\kappa^0$ in  $\Omega_\kappa$ is of Lebesgue measure zero. 
\end{Remark}
\noindent
{\em Proof of Theorem \ref{Theo:BNF}.} Without loss of generality we consider only the case when $k=1$. We devide the proof in several staps.\\ 

\noindent
{\em Step 1. Writing $\Lambda_t(\omega)$ as graphs.}\quad  
 For any $\omega\in \Omega_\kappa$ and $t\in [0,\delta]$ the Lagrangian manifold 
\[
 \Lambda_t(\omega):= \left\{\left(u_t(\theta;\omega), v_t(\theta;\omega)\right):\, \theta\in\T^n \right\}
\]
given by Theorem \ref{Theo:KAM}
is a Kronecker invariant torus of $P_t$ of a frequency vector $\omega$ satisfying (i) of Theorem \ref{Theo:KAM}. 
Firstly we will solve the equation $\varphi=u_t(\theta,\omega)$ with respect to $\theta$ and get the respective estimates of the solution. 
To this end we consider the map $w_t: \A\to\A$ defined by $w_t(\theta,\omega)=(u_t(\theta,\omega),\omega)$, where $\A=\T^{n-1}\times \R^{n-1}$. Recall from Theorem \ref{Theo:KAM}, $(iii)$, that ${\rm supp}\, (w_t-{\rm id})\subset \T^{n-1}\times (\Omega-\kappa/2)$.
It follows from \eqref{eq:estimates-F-1} and \eqref{eq:smallness-condition-1}  that 
\[
\|\sigma_\kappa^{-1}(w_t-{\rm id})\|_{1;\kappa} = \|w_t-{\rm id}\|_{1;\kappa} \le C_1\epsilon <(2n-2)^{-1}
\]
choosing $\epsilon = \epsilon(n,\tau, \vartheta_0)$ small enough (recall that $C_1$ depends only on $n$, $\tau$ and $\vartheta_0$). Then 
applying Proposition \ref{prop:inverse-holder1} we obtain a solution $\theta_t(\varphi,\omega)=\varphi+\psi_t(\varphi,\omega)$ of the equation  $\varphi=u_t(\theta,\omega)$, where $(\varphi,\omega)\in\T^{n-1}\times \Omega$. Then ${\rm supp}\, \psi_t\subset \T^{n-1}\times (\Omega-\kappa/2)$ and 
\begin{equation}\label{eq:derivatives-theta}
  \|\psi_t\|_{m;\kappa} \, \le \,  \frac{C_m}{\kappa\varrho }\,  {\mathcal B}_m 
\end{equation}
for any $m\in\N$, where the constant $C_m>0$ depends only on $n$, $\tau$, $\vartheta_0$ and $m$. Setting $\widetilde F_t(\varphi,\omega)= v_t(\theta_t(\varphi,\omega),\omega)$ and 
\begin{equation}
\label{eq:F}
F_t(\varphi,\omega):=  -\nabla K_t^{\ast}(\omega) + \widetilde F_t(\varphi,\omega)) 
\end{equation}
for $(\varphi,\omega)\in \T^{n-1}\times \Omega$ we write
\begin{equation}\label{eq:kam-tori-F}
\Lambda_t(\omega)=\{(\varphi, \widetilde F_t(\varphi,\omega)):\, \varphi\in \T^{n-1}\}=\{(\varphi, \nabla K_t^{\ast}(\omega) +  F_t(\varphi,\omega)):\, \varphi\in \T^{n-1}\},\ \omega\in \Omega_\kappa .
\end{equation}
Notice that that ${\rm supp}\, F_t\subset \T^{n-1}\times (\Omega-\kappa/2)$. 
We are going to prove that
\begin{equation}\label{eq:derivatives-F}
 \|F_t\|_{m;\kappa}   \, \le \,  \frac{C_m}{\varrho}  \, {\mathcal B}_m ,
\end{equation}
for any $m\in\N$, where the constant $C_m>0$ depends only on $n$, $\tau$, $\vartheta_0$ and $m$. 
To this end we write
\[
F_t(\varphi,\omega) = \big(v_t(\varphi,\omega)- \nabla K_t^{\ast}(\omega)\big) + \big(v_t(\theta_t(\varphi,\omega),\omega)-v_t(\varphi,\omega)\big). 
\]
The estimate of $v_t(\varphi,\omega)- \nabla K_t^{\ast}(\omega)$ follows directly from \eqref{eq:estimates-F-1} using the inequality ${\mathcal B}_1\le {\mathcal B}_2 \le \varepsilon\kappa\varrho\le \varrho$. To obtain the estimate of the second term of \eqref{eq:F} we write
\[
v_t(\theta_t(\varphi,\omega),\omega)-v_t(\varphi,\omega) \, =\, \displaystyle \int_0^1 (d_\theta v_t)(\varphi+s\psi_t(\varphi,\omega))\, \psi_t(\varphi,\omega)\,  ds .
\]
Then one  uses \eqref{eq:interpolation-leibnitz}, Proposition \ref{prop:composition-holder1}, {\em 2},  and  \eqref{eq:smallness-condition-1} as well.

Denote by $p:{\R}^{n-1} \to  {\T}^{n-1}$ the natural
projection. 
\begin{Lemma}\label{Lemma:averaging}
There is a $C^1$ family of real-valued functions 
$h_t\in  C^\infty( {\R}^{n-1}\times \Omega)$ and 
$I_t\in   C^\infty(\Omega)$ in $t\in [0,\delta]$ such that 
$h_t^0(x,\omega):= h_t(x,\omega)-\langle x,I_t(\omega)\rangle$ is $2\pi$-periodic with respect to $x$ and 
\begin{itemize}
\item[(i)] $\forall (x,\omega)\in{\R}^{n-1}\times \Omega_\kappa$, 
 $\nabla_x h_t(x,\omega)= \widetilde F_t(p(x),\omega)$,
\item[(ii)] $\nabla_x h_t^0(x,\omega) = 0$ and $I_t(\omega)=\nabla K_t^\ast(\omega)$ for $\omega\notin \Omega-\kappa/2$, 
\item[(iii)]   
$\displaystyle
 \|h_t^0\|_{m;\kappa}  \, +\, \|I_t - \nabla K_t^\ast\|_{m;\kappa} 
\, \le \,  \frac{C_m}{\varrho }   \, {\mathcal B}_m$
\end{itemize} 
for $m\in\N$ , where $C_m$ is a positive constant depending only on $n$, $\tau$, $\vartheta_0$ and $m$. 
\end{Lemma}
{\em Proof.}
To obtain $h_t$ we consider the function
$$ 
\widetilde h_t (x,\omega)\  =\  \int_{\gamma_x}^{}\, \sigma\ 
=\ \int_{0}^{1}\, \langle \widetilde F_t(p(sx),\omega),\, x\rangle\, ds\ =\ \langle \nabla K_t^\ast(\omega),x\rangle + \int_{0}^{1}\, \langle  F_t(p(sx),\omega),\, x\rangle\, ds\quad 
$$ 
for  $(x,\omega)\in {\R}^{n-1}\times \Omega$, 
where $\gamma_x = \{(sx, \widetilde  F_t(p(sx),\omega)): 0\leq s\leq 1\}$ and $\sigma =\xi dx$ is the canonical one-form on $T^\ast{\R}^{n-1}$. In view of 
\eqref{eq:derivatives-F},  the function 
$Q_t^0(x,\omega):=\widetilde h_t(x, \omega) - \langle \nabla K_t^\ast(\omega),x \rangle$
satisfies the 
estimate
\begin{equation}\label{eq:estimate-Q-0}
 \| Q_t^0\|_{m;\kappa} 
\, \le \,  \frac{C_m}{\varrho }   \, {\mathcal B}_m
\end{equation}
for  $m\in \N$.  We set 
\[
 I_{tj} (\omega) = \widetilde h_t( 2\pi e_j,\omega)/2\pi, \quad \omega\in \Omega, 
\]
where 
$\{e_1,\ldots, e_{n-1}\}$ is the canonical  basis in ${\R}^{n-1}$. Then \eqref{eq:estimate-Q-0} implies
\begin{equation}\label{eq:estimate-I-K}
 \| I_t-\nabla K_t^\ast\|_{m;\kappa}
\, \le \,  \frac{C_m}{\varrho }   \, {\mathcal B}_m
\end{equation}
for $m\in \N$ and $t\in [0,\delta	]$.

As $\Lambda_t(\omega)$ 
is a Lagrangian torus for $\omega\in \Omega_\kappa$, we get 
\begin{equation}\label{eq:integral-for-h}
\forall\, y\in {\R}^{n-1} ,  \quad \widetilde h_t (x+y,\omega) -\widetilde h_t(x,\omega)\ =\ 
\int_{l_t(x,y)}\, \sigma
\end{equation}
where $l_t(x,y) = \{(x+sy,\widetilde F_t(p(x+sy),\omega)):\, 0\leq s\leq 1\}$ and $\sigma$ is the pull-back to $\Lambda_t(\omega)$ of the fundamental one-form $Idx$.  The 
integral in \eqref{eq:integral-for-h} is equal to 
$$ 
\int_{l_t(x,y)}\, \sigma=\int_{0}^{1}\, \langle \widetilde  F_t(p(x+sy),\omega), y\rangle ds\ =\ 
\langle \widetilde F_t(p(x),\omega),y\rangle + O(y^2), 
$$ 
and we obtain  $\nabla_x \widetilde h(x,\omega)= \widetilde F_t(p(x),\omega)$  for any  $\omega\in \Omega_\kappa$. 
In particular, the function 
$\nabla_x \widetilde h_t(x,\omega)$ is $2\pi$-periodic with
respect to $x$ and we get
\[
\forall\, \omega\in\Omega_\kappa,\ \forall\, \alpha \in {\Z}^{n-1}, \quad \widetilde h_t (x+2\pi \alpha,\omega) - \widetilde h_t (x,\omega) = \widetilde h_t (2\pi \alpha,\omega) - \widetilde h_t (0,\omega)
=\langle 2\pi  \alpha, I_t(\omega)\rangle.
\] 
Consider the function 
\[
\widetilde h_t^0(x,\omega) = 
\widetilde h_t(x,\omega) - \langle x, I_t(\omega) .
\]
It  is $2\pi$-periodic with respect to
$x$ for $\omega\in \Omega_\kappa$ and $\widetilde h_t^0$ satisfies the  estimates \eqref{eq:estimate-Q-0} in 
$[0,\delta]\times \R^{n-1}\times \Omega$. 
We are going to average  $\widetilde h_t^0$ on ${\T}^{n-1}$ using the following 
\begin{Lemma}\label{Lemma:averaging function} There exists 
$f\in C^\infty({\R}^{n-1})$ with ${\rm supp}\, f \subset [\pi,
7\pi]^{n-1}$ such that 
\[
\sum_{k\in {\Z}^{n-1}} f(x -  2\pi k) = 1
\]
for each $x\in {\R}^{n-1}$. 
\end{Lemma}
Consider the function 
\[
h_t^0(x,\omega) =
\sum_{k\in {\Z}^{n-1}} (f \widetilde h_t^0)(x -  2\pi k,\omega).
\] 
It is $2\pi$-periodic with respect to $x$ by
construction  and 
$h_t^0(x,\omega) = \widetilde h_t^0(x,\omega)$ for $(x,\omega)\in 
{\R}^{n-1}\times \Omega_\kappa$.
Moreover,  $h_t^0$ satisfies \eqref{eq:estimate-Q-0} in 
$[0,\delta]\times\R^{n-1}\times \Omega$.  
 We set
\[
h_t(x,\omega)=h_t^0(x,\omega) + \langle x, I_t(\omega)\rangle .
\] 
Recall that  
\[
\mbox{dist}\, (\Omega_\kappa, {\R}^{n-1}\setminus
\Omega) \ge  \kappa. 
\]
Then multiplying $h_t^0$ and $I_t- \nabla K_t^\ast$ by
a suitable 
cut-off function  which is equal to  one on 
$\Omega-3\kappa/4$ and zero outside $\Omega-\kappa/2$, we can
assume that $h_t(x,\omega) =  \langle x, \nabla K_t^\ast(\omega)\rangle$ and $I_t(\omega) = \nabla K_t^\ast(\omega)$ outside $\Omega-\kappa/2$  (see Lemma \ref{lemma:cut-off-omega}). This proves $(ii)$.  The statement $(iii)$ follows from \eqref{eq:estimate-Q-0}, the definition of $h_t^0$ and \eqref{eq:estimate-I-K}.\finishproof
\\

\noindent
{\em Step 2. Inverting $I_t$.}
\begin{Lemma}\label{Lemma:action} Choosing $\epsilon=\epsilon(n,\tau,\vartheta_0)>0$ small enough one has  the following for each $t\in [0,\delta]$. 
\begin{enumerate}
\item The map $I_t:\Omega \rightarrow D_t$ is a diffeomorphism  and its inverse $\omega_t: D_t\to \Omega$ satisfies the estimates 
$\|\omega_t -\nabla K_t\|_{1,D_t;\kappa} \le C_0\epsilon\kappa$ and 
\begin{equation}\label{eq:omega}
 \| \omega_t -\nabla K_t\|_{m, D_t;\kappa}\, 
 \le \,  \frac{C_m}{\varrho}  \,{\mathcal B}_{m} \lambda^{m}\big(\lambda+\|\partial^2 K_t\|_{m,D_t;\kappa}\big)
\end{equation}
for any $m\in\N_\ast$. Moreover, $\omega_t=\nabla K_t$ outside the set $\overline D_t^1:= \nabla K_t^\ast(\Omega-\kappa/2)$. 
\item For any $x\in \R^{n-1}$ the map $\Omega\ni\omega\rightarrow \nabla_x h_t(x,\omega)\in D_t$ is a diffeomorphism.
\end{enumerate}
\end{Lemma}
{\em Proof}. We are going to show that the map $\Omega\ni \omega\to I_t(\omega)$ a diffeomorphism. To this end we write
\[
I_t = ({\rm id} +\varphi_t)\circ \nabla K_t^\ast,\quad \varphi_t:= (I_t -\nabla K_t^\ast)\circ \nabla K_t.
\]
Moreover, Lemma \ref{Lemma:averaging}, \eqref{eq:non-degeneracy-BNF} and  \eqref{eq:smallness-condition-1} imply
\[
\|\varphi_t\|_{1} \le  \|I_t -\nabla K_t^\ast\|_{0} + \|I_t -\nabla K_t^\ast\|_{1}(1+\|\nabla K_t\|_{1}) \le C \epsilon \lambda^{-2}
\]
where $C=C(n,\tau,\vartheta_0)>0$. Recall that $I_t=\nabla K_t^\ast$ outside $\Omega-\kappa/2$, hence ${\rm supp\,}\varphi_t\subset D_t$. Then choosing $0<\epsilon<1/(2C)$ and  applying Proposition \ref{prop:inverse-holder1} we obtain that ${\rm id} +\varphi_t:\R^{n-1}\to \R^{n-1}$ is invertible and that
$({\rm id} +\varphi_t)^{-1} ={\rm id} +\psi_t$, where ${\rm supp\,}\psi_t\subset D_t$. Hence, $I_t(\Omega)= D_t$ and 
 $I_t:\Omega\to D_t$ is a diffeomorphism with inverse 
 \[
 \omega_t=\nabla K_t\circ ({\rm id} +\psi_t):D_t\to \Omega.
 \]
\begin{Lemma}\label{Lemma:estimates-phi} 
For any $m\in\N_\ast$ there exists $C_m>0$ depending only on $m$, $n$, $\tau$ and $\vartheta_0$ such that
\[
\|\varphi_t\|_{m,\kappa} \le \frac{C_m}{\varrho}  \,{\mathcal B}_{m}  \lambda^{m-1}   \big(\lambda+\|\partial^2 K_t\|_{m-1,D_t;\kappa}\big). 
\]
\end{Lemma}
{\em Proof}. The support of $\varphi_t$ is contained in the closure of  $D_t^1:= \nabla K_t^\ast (\Omega-\kappa/2)$. Set $2r_t:= {\rm dist}\, (D_t^1, \R^{n-1}\setminus D_t)$ and fix $I^0$ in $\overline D_t^1$. Applying Remark \ref{rem:interpolation} to the restriction of $K_t$ to the closed ball 
$\overline B(I^0, r_t)$ as well as   Proposition \ref{prop:composition-holder1},    \eqref{eq:non-degeneracy-BNF} and  \eqref{eq:smallness-condition-1} one obtains 
\[
\begin{array}{cclr}
\displaystyle |(\kappa\partial_I)^\beta\varphi_t(I)|\le 
 C_m^\prime\left(1+ \|\partial^2 K_t\|_{C^0(D_t)}^{m-1}\right) \\[0.3cm]
\displaystyle \times \left(\|I_t -\nabla K_t^\ast\|_{m,\kappa}\| \partial^2 K_t\|_{{C^0(D_t)}}+\|I_t -\nabla K_t^\ast\|_{C^1}\|\partial^2 K_t\|_{m-1,D_t;\kappa}\right)\\[0.3cm]
\displaystyle \le \frac{C_m}{\varrho}  \,{\mathcal B}_{m}   \lambda^{m-1}   \big(\lambda+\| \partial^2 K_t\|_{m-1,D_t;\kappa}\big),
\end{array} 
\]
for $(I,t)\in B(I^0, r_t)\times [0,\delta]$ and $\beta\in \N^{n-1}$, $ |\beta|=m\in\N_\ast$, 
where $C_m^\prime>0$ depends only on $m$ and $n$ and $C_m>0$ depends only on $m$, $n$, $\tau$ and $\vartheta_0$. 
On the other hand, $\varphi_t=0$ outside $D_t^1$ which completes the proof of the Lemma. This argument will be used many times in the sequel.
\finishproof

\noindent
Proposition \ref{prop:inverse-holder1} applied to ${\rm id} + \kappa^{-1}\varphi_t \circ\sigma_\kappa$  implies that ${\rm supp}\, \psi_t\subset \overline D_t^1$, $\|\psi_t\|_{C^0}\le C_0\epsilon/\lambda$  and 
\[
\|\psi_t\|_{m,\kappa} \le \frac{C_m}{\varrho}  \,{\mathcal B}_{m}   \lambda^{m-1}  \big(\lambda+\|\partial^2 K_t\|_{m-1,D_t;\kappa}\big)
\]
for $m\in\N^\ast$. In particular, $\|d \psi_t\|_{C^0} \le C_1\epsilon$. 
Consider
\[
\omega_t(I) -\nabla K_t(I) = \int_0^1   d(\nabla K_t)(I + s\psi_t(I)) \,  \psi_t(I)\,  ds, \quad I\in D_t. 
\]
The support of $\omega_t -\nabla K_t$ is contained in $\overline D_t^1$. Moreover, $\|\omega_t -\nabla K_t\|_{C^0(D_t)} \le C_0\epsilon\kappa$, and using Remark \ref{rem:interpolation} 
and Proposition \ref{prop:composition-holder1} 
one obtains as in the proof of Lemma \ref{Lemma:estimates-phi} the estimate
\[
\begin{array}{cclr}
\displaystyle
\|(\kappa \partial)^\alpha(\omega_t -\nabla K_t)\|_{C^0} \le C_m \|\partial^2 K_t\|_{C^0(D_t)}\|\psi_t\|_{m,\kappa} \\[0.3cm]
\displaystyle +\,  C_m\|\psi_t\|_{C^0}\big(1+ \|d \psi_t\|_{C^0}^{m-1}\big)\big(\|\partial^2 K_t\|_{1,D_t;\kappa}\|d\psi_t\|_{m-1,\kappa} +
\|\partial^2 K_t\|_{m,D_t;\kappa}\|d\psi_t\|_{C^0(D_t)}\big)
\end{array}
\]
 for any $m\in \N_\ast$ and $\alpha\in\N^{n-1}$ with $|\alpha|=m$ . 
Then 
using \eqref{eq:non-degeneracy-BNF}, \eqref{eq:smallness-condition-1} and the previous estimates we obtain \eqref{eq:omega}. The estimate  
$\|d\omega_t- d\nabla K_t\|_{C^0(D_t)} \le C \epsilon\lambda$ follows from \eqref{eq:omega} with $m=1$ and \eqref{eq:smallness-condition-1}. 
We are going to prove that for each $x\in \R^{n-1}$ the map $\Omega\ni\omega\rightarrow \nabla_x h_t(x,\omega)$ is a diffeomorphism. 
To this end we fix $x$ and we write the map $\omega\to \nabla_x h_t(x,\omega)$ as follows 
\[
\nabla_x h_t= ({\rm id} +\varphi_t^1)\circ \nabla K_t^\ast, \quad
\varphi_t^1:= (\nabla_x h_t^0 + I_t -\nabla K_t^\ast)\circ \nabla K_t .
\]
Then ${\rm supp\,}\varphi_t^1\subset \overline D_t^1$  and 
\[
\|\varphi_t^1\|_{C^1} \le (\|\nabla_x h_t^0\|_{C^1}+\|I_t -\nabla K_t^\ast\|_{C^1})(1+\|\partial^2 K_t\|_{C^0}) \le C_1 \epsilon/\lambda
\]
and we complete the proof of {\em 2} as above. 
\finishproof
\\

\noindent
{\em Step 3. Construction of $\chi_t$.}\quad 
The second statement of the Lemma implies that there is 
a $C^\infty$-foliation of ${\T}^{n-1}\times D_t$ 
by  Lagrangian tori 
\[
\Lambda_t(\omega)\ =\ \{(p(x),\nabla_x h_t(x,\omega)) :\, 
x\in{\R}^{n-1}\},\ \omega\in \Omega,
\]
which is a smooth extension of the family of the Kronecker invariant  tori \eqref{eq:kam-tori-F} of $P_t$. 
Notice that  $I_t(\omega)$ is the action along the basis of cycles $[\gamma_{t,j}(\omega)], \ldots, [\gamma_{t,n-1}(\omega)]$ of   $H_1(\Lambda_t(\omega),\R)$, where $\gamma_{t,j}(\omega) = 
\{(p(s2\pi e_j),\, \nabla_x h_t(s2\pi e_j,\omega))\, :\, 
0\leq s\leq 1\}$. Indeed, it follows from the definition of $h_t$ that $I_t(\omega)=\nabla_x h_t(x,\omega)-\nabla_x h_t^0(x,\omega)$, where $h_t^0$ is $2\pi$-periodic in $x$ and we obtain
\begin{equation}\label{eq:momentum-I-2}
I_{t}(\omega)\ =\ \left(\int_{\gamma_{t,1}(\omega)}\, \sigma ,\ldots, \int_{\gamma_{t,n-1}(\omega)}\, \sigma \right)
\end{equation}
for  $\omega\in \Omega$.   
Now we set $\Phi_t(x,I) = h_t(x,\omega_t(I))$. Then  
\[
\Phi^0_t(x,I):=\langle x, I\rangle -\Phi_t(x,I)  = - h^0_t(x,\omega_t(I))
\]
is $2\pi$-periodic with
respect to $x$, and  it has a compact support in $\T^{n-1}\times D_t$. Moreover, it follows from Lemma \ref{Lemma:averaging},  Lemma \ref{Lemma:action}, Remark \ref{rem:interpolation} and Proposition \ref{prop:composition-holder-interpolation} that
\[
\begin{array}{lcrr}
\|\Phi^0_t\|_{m,\A_t;\kappa}\, 
\le \, C_m \left(1+ \|d \omega_t\|_{C^0(D_t)}^{m-1} \right) \Big(\|h_t^0\|_{m,\kappa}\|d \omega_t\|_{C^0(D_t)}\\[0.3cm]
\displaystyle
  + \|h_t^0\|_{1,\kappa}\|\partial^2 K_t\|_{m-1,D_t;\kappa}+ \|h_t^0\|_{1,\kappa}\|(d \omega_t- d\nabla K_t)\|_{m-1,D_t;\kappa}\Big) \\[0.3cm]
\displaystyle
\le \frac{C_m}{\varrho}  \,{\mathcal B}_{m}  \lambda^{2m-2}  \big(\lambda +\|\partial^2 K_t\|_{m,D_t;\kappa}\big)
\end{array}
\]
for  $t\in [0,\delta]$ and  $m\in\N_\ast$, where $C_m$ depends only on $m$, $n$, $\tau$ and $\vartheta_0$. This implies 
\begin{equation}\label{eq:estimates-Phi-0}
\|\sigma_\kappa^{-1}{\rm sgrad}\Phi^0_t\|_{m,\A_t;\kappa}\, 
\le \,  \frac{C_m}{\kappa\varrho}  \,{\mathcal B}_{m+1}  \lambda^{2m}  \big(\lambda +\|\partial^2 K_t\|_{m+1,D_t;\kappa}\big)
\end{equation}
for  $t\in [0,\delta]$ and  $m\in\N$, where $C_m$ depends only on $m$, $n$, $\tau$ and $\vartheta_0$.
In particular, one obtains by means of \eqref{eq:non-degeneracy-BNF} and \eqref{eq:smallness-condition-1} that
\begin{equation}\label{eq:estimates-Phi-0-2}
\|\sigma_\kappa^{-1}{\rm sgrad}\, \Phi^0_t\|_{1,\A_t;\kappa}  \le c\epsilon/\lambda  
\end{equation}
for $t\in [0,\delta]$, where $c=c(n,\tau,\vartheta_0)>0$. Using  Lemma \ref{lemma:generating-functions} we obtain
\begin{Lemma}\label{Lemma:chi} Choosing $\epsilon=\epsilon(n,\tau,\vartheta_0)>0$ small enough we have the following
\begin{enumerate}
\item $\Phi_t$ is a generating function of a symplectic transformation $\chi_t:\T^{n-1}\times D_t \to \T^{n-1}\times D_t$ and the map $[0,\delta]\ni t\to \chi_t\in C^\infty(\A_t,\A_t)$ is $C^1$; 
\item $\chi_t(\Lambda_t(\omega))=\T^{n-1}\times \{I_t(\omega)\}$ for any  $\omega\in \Omega$ and $t\in [0,\delta]$;
\item $\chi_t-{\rm id}$ and $\chi_t^{-1}-{\rm id}$ are compactly supported in $\T^{n-1}\times \overline D_t^1$, where $D_t^1=\nabla K_t^\ast(\Omega-\kappa/2)$ and they satisfy the estimates \eqref{eq:estimates-chi}. Moreover,
\[
\|\sigma_\kappa^{-1}(\chi -{\rm id})\|_{1,\A_t;\kappa} + \|\sigma_\kappa^{-1}(\chi -{\rm id})\|_{1,\A_t;\kappa} \, \le\,   c\epsilon/\lambda  . 
\] 
\end{enumerate}
\end{Lemma}
{\em Proof}. Using  Lemma \ref{lemma:generating-functions}  one  
obtains a symplectic transformation $\chi_t: {\T}^{n-1}\times D_t\to {\T}^{n-1}\times \R^{n-1}$ defined  by 
\[
\chi_t(\nabla_I\Phi_t(\theta,I), I)= (\theta, \nabla_\theta\Phi_t(\theta,I)), \quad (\theta,I)\in  {\T}^{n-1}\times \R^{n-1}. 
\]
Notice that the map
\[
D_t \ni I \to \nabla_\theta\Phi_t(\theta,I)= I + \nabla_\theta h_t^0(\theta,\omega_t(I))\in D_t
\] 
is a diffeomorphism since the map 
$\Omega \ni \omega \to I_t(\omega) + \nabla_\theta h_t^0(\theta,\omega)= \nabla_\theta h_t(\theta,\omega)\in D_t$ is a diffeomorphism in view of Lemma \ref{Lemma:action}, {\em 2}, hence,  $\chi_t(\A_t)=\A_t$. 
For any $\omega\in
\Omega$ and any $\theta$ we
have 
\[
(\theta, \nabla h_t(\theta,\omega)) = (\theta, \nabla_\theta\Phi_t(\theta,I_t(\omega)))
= \chi_t(\nabla_I\Phi_t(\theta,I_t(\omega)), I_t(\omega)),
\] 
hence, $\Lambda_t(\omega) =
\chi_t({\T}^{n-1}\times \{I_t(\omega)\})$. Moreover,  $\chi_t(\varphi, I_t(\omega))= (\varphi, I_t(\omega))= (\varphi,\nabla K_t^\ast(\omega))$ if  $\mbox{dist}\, (\omega, {\R}^{n-1}\setminus \Omega) \le \kappa/2$, hence, the support of both $\chi_t-{\rm id}$ and  $\chi_t^{-1}-{\rm id}$ is contained in $\T^{n-1}\times \overline D_t^1$. The estimate    \eqref{eq:estimates-chi} follows from \eqref{eq:estimates-Phi-0} and Lemma \ref{lemma:generating-functions}. 
\finishproof
\\

\noindent
{\em Step 4. Estimates.}\quad 
Consider the exact symplectic map $P_t^0 = \chi_t^{-1} P_t\chi_t$. Using Lemma \ref{Lemma:decomposition} we write $P_t$  as a  composition $P_t=W_t Q_t$, where $W_t$ is the exact symplectic map defined by the  generating function $(x,r)\to \langle x, r\rangle - G_t(x,r)$ and $Q_t(\theta,r)=(\theta +\nabla K(r),r)$. 
Then  
\[
P_t^0 = W_t^0 Q_t,
\] 
where  
\[
W_t^0 = \chi_t^{-1} W_t + \chi_t^{-1} W_t Q_t (\chi_t^{-1} - {\rm id}) Q_t^{-1}.
\]
\begin{Lemma}\label{Lemma:W-1} The exact symplectic map $W_t^0$, $t\in [0,\delta]$,  admits a generating function of the form 
\[
(x,I) \to \langle x, I\rangle - G_t^0(x,I)
\]
 such that the map $[0,\delta]\ni t \to G_t^0\in C^\infty(\A)$ is $C^1$, ${\rm supp}\, (dG_t^0) \subset \T^{n-1}\times \overline D_t^1$ and 
\[
\|\sigma_\kappa^{-1}{\rm sgrad}\, G^0_t\|_{m,\A_t;\kappa} \, \le\,  
\frac{C_m}{\varrho\kappa } \,{\mathcal B}_{m+1}  \lambda^{2m} \big(\lambda+\|\partial^2 K_t\|_{m+1,D_t;\kappa}\big)
\]
for any $m\in\N$, where  $C_m$ depends only on $m$, $n$, $\tau$ and $\vartheta_0$.
\end{Lemma}
{\em Proof}. 
We have 
\begin{equation}\label{eq:identity-for-W-1}
\begin{array}{lcrr}
W_t^0 - {\rm id}\,  =\,  \big(W_t - {\rm id} + (\chi_t^{-1}- {\rm id}) W_t\big)\, +\,  Q_t(\chi_t- {\rm id})Q_t^{-1} \\[0.3cm]
+ \big( W_t - {\rm id} + (\chi_t^{-1}- {\rm id}) W_t\big)\, \circ\, \big(Q_t(\chi_t- {\rm id})Q_t^{-1}\big).
\end{array}
\end{equation}
We estimate the $C^m$ norms of it term by term. 
Notice that the support of each term is contained in $\T^{n-1}\times \overline D_t^1$. 

Lemma \ref{lemma:generating-functions} and \eqref{eq:smallness-condition-1} imply $\|\sigma_\kappa^{-1}(W_t - {\rm id})\|_{1,\A_t;\kappa}< C\epsilon/\lambda$ and 
\[
\|\sigma_\kappa^{-1}(W_t - {\rm id})\|_{m,\A_t;\kappa}\, \le\,  \frac{C_m}{\varrho\kappa}  \,{\mathcal B}_{m}  . 
\]
The last estimate, Lemma \ref{Lemma:chi}, {\em 3}, and Lemma \ref{prop:composition-holder1}, {\em 2}, imply
\[
\|\sigma_\kappa^{-1}(\chi_t^{-1}- {\rm id}) W_t\|_{m,\A_t;\kappa}\, \le\,  
\frac{C_m}{\varrho\kappa }  \,{\mathcal B}_{m+1} \lambda^{2m} \big(\lambda+\|\partial^2 K_t\|_{m+1,D_t;\kappa}\big).   
\]
and 
\[
\|\sigma_\kappa^{-1}(\chi^{-1} -{\rm id})W_t\|_{1,\A_t;\kappa}  \, \le\,   c\epsilon/\lambda  . 
\] 
Using the argument in the proof of Lemma \ref{Lemma:estimates-phi} first to $ (\chi_t- {\rm id})Q_t^{-1}$ and then to  $Q_t(\chi_t- {\rm id})Q_t^{-1}$ we obtain the estimate
\[
\|\sigma_\kappa^{-1}Q_t(\chi_t- {\rm id})Q_t^{-1}\|_{m,\A_t;\kappa}\, \le\,  
\frac{C_m}{\varrho\kappa } \,{\mathcal B}_{m+1} \lambda^{2m}\big(\lambda+\|\partial^2 K_t\|_{m+1,D_t;\kappa}\big) .  
\]
This yields the  estimate of the first line of \eqref{eq:identity-for-W-1}. To obtain the estimates for the composition in the second line one uses the preceding estimates and the argument of Lemma \ref{Lemma:estimates-phi}. It remains to apply Lemma \ref{lemma:generating-functions} in order to complete  the proof of the Lemma.  \finishproof
\\

\noindent
{\em Step 4. Proof of $(i)$.}\quad  
Now Lemma \ref{Lemma:decomposition} implies that the function 
$$\widetilde G^0_t(x,I)= \langle x, I\rangle -  K_t(I) - G_t^0(\varphi,I)$$ 
is a generating function of $P_t^0$. The function  $\widetilde G^0_t$ is uniquely defined modulo a constant depending only on $t$ which is chosen  appropriately in order to obtain a $C^1$-smooth  the map $t\to \widetilde G^0_t$. Set $L_t(I):= K_t(I)- G^0_t(0, I)$ and $R_t(\theta, I)= G^0_t(\theta, I)- G^0_t(0, I)$. We have 
\[
P_t^0({\T}^{n-1}\times \{I_t(\omega)\}) = (\chi_t^{-1}\circ P_t) (\Lambda_t(\omega))= \chi_t^{-1}(\Lambda_t(\omega))= {\T}^{n-1}\times \{I_t(\omega)\}
\]
for any $\omega\in \Omega_\kappa$, which implies 
\[
I_t(\omega) - \nabla_x R_t^0(x,I_t(\omega))= \nabla_x \widetilde G_t^0(x,I_t(\omega))=I_t(\omega)
\]
for any such $x\in\R^{n-1}$. On the other hand, $R_t^0(0,I)=0$, hence, $R_t^0(\theta,I)=0$ on ${\T}^{n-1}\times E_{t,\kappa}$, where  $E_{t,\kappa}:= I_t(\Omega_\kappa)$. 
Now Lemma \ref{Lemma:flat} implies that $\partial_I^\beta R_t^0(\theta,I)=0$ for any  $(\theta, I)\in {\T}^{n-1}\times E_{t,\kappa}^0$, where  $E_{t,\kappa}^0:= I_t(\Omega_\kappa^0)$ is the set of points of   positive Lebesgue density in $E_{t,\kappa}$. 
Then
$$P_t^0(\theta + \nabla L_t(I), I)= (\theta,I) \quad \mbox{on} \quad {\T}^{n-1}\times E_{t,\kappa}^0, $$
and we obtain that $\nabla_I L_t(I_t(\omega))=\omega$ for each $\omega\in  \Omega_{\kappa}^0$.  Hence $\nabla L_t^\ast(\omega)= I_t(\omega)$ for each $\omega\in  \Omega_{\kappa}^0$, where  $I_t(\omega)$  is given by   \eqref{eq:momentum-I} for such $\omega$,  according to \eqref{eq:momentum-I-2}.  
The estimates of the derivatives of $L_t$ and $R_t^0$ follow from that of Lemma \ref{Lemma:W-1}. 
\finishproof

\section{KAM theorem with parameters}\label{Sec:KAM-parameters}
The  theorems formulated above follow from a KAM theorem with parameters. A complete and very comprehensive proof of it has been given by P\"oschel \cite{Poe1} and Kuksin  \cite{Kuk} in the analytic case. It can be  extended to the case of smooth Hamiltonians using suitable approximation lemma. In the case of Gevrey Hamiltonians this has been done in \cite{P4}. The advantage of this approach is that frequencies are separated from action variables which makes it easier to obtain smoothness with respect to them. Moreover, it allows one to prove   H\"older estimates of the transformations putting  the Hamiltonian to a normal form. Here, the normal form of the Hamiltonian is  $N(I;\omega):= \langle\omega,I\rangle$.
The perturbation is a  real valued function $(\theta,I;\omega,t) \mapsto P(\theta,I;\omega,t)$ defined in $\A^n\times \Omega\times [0,a]$, where $\A^n:= \T^n\times B(0,\rho_0)$, $B(0,\rho_0)\subset \R^n$ is the ball centered at $I=0$ with radius $\rho_0\in (0,1]$ and $\Omega$ is a bounded domain in $\R^n$. Hereafter, we assume that
\begin{equation}
P\in C^k\left([0,a];C_0^\infty(\A^n\times\Omega)\right), \quad k\in \{0;1\},
\label{eq:families-of-perturbations}
\end{equation}
i.e. the map
  $t\to P_t:=P(\cdot,t)\in C_0^\infty(\T^n\times B(0,\rho_0)\times \Omega)$ is  $C^k$-smooth on the interval $[0,a]$.  This means that  the support of the function   $(I,\omega) \to P(\theta,I;\omega,t)$ is contained in a fixed  compact subset of  $B(0,\rho_0)\times \Omega$ independent of  $(\theta,t)\in \T^n\times [0,a]$ and 
that the maps  
\[
t\to P_t:=\partial_t^q P(\cdot,t)\in C^j(\T^n\times B(0,\rho_0)\times \Omega), \quad 0\le q\le k, 
\]
are continuous in $t\in [0,a]$ for  $j\in \N$. 
Given  $\ell>0$ and $0<r,\kappa\le 1$, $r\le \rho_0$,  we denote by $\|P_t\|_{\ell;r,\kappa}$ the  weighted H\"older norm 
\begin{equation}\label{eq:norm-rho-kappa}
\|P_t\|_{\ell;r,\kappa} := \|P_t\circ \sigma_{r,\kappa}\|_{C^\ell(\sigma_{r,\kappa}^{-1}(\A^n\times\Omega))}
\end{equation}
where  $\sigma_{r,\kappa}$ is  the partial dilation  $\sigma_{r,\kappa}(\varphi,I;\omega):= (\varphi, rI; \kappa \omega)$. The H\"older norms are defined in Section \ref{Sec:ApprLemma} (see also \cite{Poe}).

Fix  $1\le \vartheta_1<\vartheta_0 < \tau +1$ 
and set 
\begin{equation}
\ell_0:=2\tau +  2+ \vartheta_0  \quad \mbox{and} \quad  \ell(m):=  2m(\tau + 1)+ \ell_0 \, , \quad m \in \N. 
\label{eq:ell-def}
\end{equation}

Denote the Hamiltonian vector field associated to  the Hamiltonian $(\theta,I)\to N(I;\omega)=\langle\omega,I\rangle$ by ${\mathcal L}_\omega := X_N(\cdot,\omega)=\langle \omega, \partial/\partial \theta\rangle$.  We consider $C^k$-families, $ k\in \{0;1\}$,  of Hamiltonians $t \to  H_{t,\omega}$, $t\in [0,\delta]$, where
\[
H_{t,\omega}(\theta,I):= H(\theta,I; \omega,t) = N(I;\omega) + P(\theta,I; \omega,t)
\]
and $P$ satisfies \eqref{eq:families-of-perturbations}. Recall that for given $0<\kappa\le 1$ and $\tau>n-1$, the set $\Omega_{ \kappa}=\widetilde D (\kappa,\tau)\cap \overline{\Omega-\kappa}$ consists of all $(\kappa,\tau)$-Diophantine frequencies  $\omega$ in $\Omega$ ($\omega$ satisfies \eqref{eq:sdc2}) such that the distance from $\omega$ to the complement of $\Omega$ in $\R^n$ is greater or equal to $\kappa$. 
Set
 \begin{equation}\label{eq:main-estimates1}
 \begin{array}{rcll}
\displaystyle \left\langle P \right\rangle_{\ell(m);r,\kappa} ^{(0)}&=&  \displaystyle  \sup_{t\in [0,a]} \, \|P_t\|_{\ell(m);r,\kappa}  \\ [0.5cm]
\displaystyle   \left\langle P \right\rangle_{\ell(m);r,\kappa} ^{(1)}&=&  \displaystyle   \frac{\left\langle P \right\rangle_{\ell(m);r,\kappa} ^{(0)}}{\left\langle P \right\rangle_{\ell(0);r,\kappa} ^{(0)}}\, \sum_{0\le p\le 1}\sup_{t\in [0,a]} \, \|\partial_t^p  P_t\|_{\ell(m);r,\kappa}. 
\end{array}
\end{equation}
The following result is an analogue of Theorem A in \cite{Poe1}. 
\begin{Theorem}\label{Theo:A}
 There exists a  positive constant $\epsilon= \epsilon(n,\tau,\vartheta_0, \vartheta_1)>0$  depending only  on $n$, $\tau$, $\vartheta_0$ and $\vartheta_1$  such that, for any $a>0$, $0<\kappa<1$, $0<r < \rho_0$ and  any real valued Hamiltonian  $H=N+P$, where $N(I;\omega)=\langle\omega,I\rangle $ and  $P$ satisfies   \eqref{eq:families-of-perturbations}  and the smallness hypothesis
\begin{equation}
   \sup_{t\in [0,a]} \, \| P_t\|_{\ell_0;r,\kappa} \  \le \ \epsilon \kappa r \, ,
\label{eq:smallness-condition3}
\end{equation}
 the following holds.

\vspace{0.3cm}
\noindent
There exist $C^k$ families of maps 
\[
[0,a]\ni t \mapsto \phi_t\in C^\infty(\Omega;\Omega) \quad \mbox{and} \quad [0,a]\ni t \mapsto \Psi_t=(U_t,V_t)\in C^\infty(\T^n \times\Omega;\T^n\times B(0,r)) \ 
\] 
such that   ${\rm supp\,}(\phi_t-{\rm id})\subset \Omega-\kappa/2$, ${\rm supp\,}\big((U_t,V_t)-({\rm id},0)\big)\subset \T^n\times (\Omega-\kappa/2)$,  and 
\begin{enumerate}
\item[(i)]
For each $\omega \in \Omega_\kappa$ and $t\in [0,a]$ the map 
$\Psi_{t,\omega} := \Psi_t(\cdot,\omega):{\T}^n  \rightarrow
 {\T}^n\times B(0,r)$ 
is a smooth embedding, 
$\Lambda_t(\omega)  := \Psi_{t,\omega}({\T}^n)$ is an embedded Lagrangian 
torus invariant  with respect to the Hamiltonian flow of 
$ H_{t,\phi_t(\omega)}(\theta,I) := H(\theta, I; \phi_t(\omega),t)$, and 
\[
X_{H_{t,\phi_t(\omega)}} \circ \Psi_{t,\omega} = 
D \Psi_{t,\omega} \cdot {\mathcal L}_\omega \quad \mbox{on}\  {\T}^n ,  
\]
\item[(ii)]
For any $m\ge 0$ there is $C_m>0$ depending only on  $n$, $\tau$, $\vartheta_0$,  $\vartheta_1$, and $m$, such that
for any $\alpha,\beta\in \N^n$ of length $|\alpha|+|\beta|(\tau+1)\le m(\tau+1)+ \vartheta_1$ and $0\le q\le k$ the following estimate holds
\begin{equation}\label{eq:main-estimates}
\begin{array}{lrc}
\left|\partial_\theta^{\alpha}(\kappa \partial_\omega)^{\beta}\partial_t^q
(U_t(\theta;\omega) - \theta)\right|\,  + \,
r ^{-1}
\left|\partial_\theta^{\alpha}(\kappa \partial_\omega)^{\beta} \partial_t^q V_t(\theta;\omega)\right|  
\\    [0.5cm]  
\displaystyle  +\,  \kappa^{-1} \left|(\kappa \partial_\omega)^{\beta} \partial_t^q
(\phi_t(\omega) - \omega)\right| \, \leq\,    C_{m}\, 
 (\kappa  r )^{-1}\,  \left\langle P \right\rangle_{\ell(m);r,\kappa} ^{(q)}
\end{array}
\end{equation}
 uniformly in  $(\theta,\omega,t)\in {\T}^n\times \Omega\times [0,a]$. 
\end{enumerate}
\end{Theorem}

\begin{Remark}\label{rem:analiticity-t}
If $P$ is analytic with respect to $t$ in  the disc $B(0,a):= \{t\in \C:\, |t|<a\} $  and \eqref{eq:smallness-condition3} holds for $t\in B(0,a)$, then $\Psi$ and $\phi$ can be chosen to be  analytic with respect to $t$ in $B(0,a)$. Moreover, 
for any $\alpha,\beta\in \N^n$ of length $|\alpha|+|\beta|(\tau+1)\le m(\tau+1)+ \vartheta_1$ and $0<\delta<a$,  the following estimate holds
\[
\begin{array}{lrc}
\left|\partial_\theta^{\alpha}(\kappa \partial_\omega)^{\beta}
(U_t(\theta;\omega) - \theta)\right|\,  + \,
r ^{-1}
\left|\partial_\theta^{\alpha}(\kappa \partial_\omega)^{\beta} V_t(\theta;\omega)\right|  
+ \kappa^{-1} \left|(\kappa \partial_\omega)^{\beta} (\phi_t(\omega) - \omega)\right|
\\    [0.5cm]  
\displaystyle \leq\    C_{m,\delta}\, 
(\kappa  r )^{-1}\,  \sup_{t\in B(0,a)} \,\|\partial_t^p P_t\|_{\ell(m);r,\kappa} 
\end{array}
\]
uniformly in  $(\theta,\omega,t)\in {\T}^n\times \Omega\times B(0, a-\delta)$, with  $C_{m,\delta}>0$ depending only on $n$, $\tau$, $\vartheta_0$,$\vartheta_1$, $\delta$, $m$. 
\end{Remark} 
Before starting the proof of Theorem \ref{Theo:A} and Remark \ref{rem:analiticity-t}
we are going to list several comments. 
For each $t\in[0,\delta]$ and $\omega\in \Omega$ 
denote by 
\[
\Phi_{t,\omega}^s := \exp\left(s X_{H_{t,\omega}} \right) ,\ s\in\R,
\]
the flow of the Hamiltonian vector field $X_{H_{t,\omega}}$ of the Hamiltonian $H_{t,\omega}$ and set
\[
g_\omega^s(\theta) = \theta + p(s \omega) ,\ \theta \in\T^n,\ s\in \R, \ \omega\in \Omega,
\] 
where $p:\R^n\to \T^n$ is the canonical projection. By \eqref{eq:smallness-condition3} and \eqref{eq:main-estimates}, we have 
$$|d_\theta U_t(\theta;\omega)-\mathrm{Id} |\le C_1(n,\tau,\vartheta_0)\epsilon \le 1/2$$ 
for $(\theta,\omega,t)\in {\T}^n\times \Omega\times [0,a]$, choosing $\epsilon $ sufficiently small and we obtain 
\begin{Remark}\label{rem:Kronecker} 
The assertion (i) of Theorem \ref{Theo:A} means that for each $\omega\in\Omega_\kappa$ the family $[0,\delta]\ni t \rightarrow \Lambda_t(\omega)$ is a $C^k$ family of Kronecker invariant tori with respect to the flow $\Phi_{t,\widetilde\omega}^s$, where $\widetilde\omega=\phi_t(\omega)$. More precisely, for each $t\in[0,\delta]$, $\omega\in\Omega_\kappa$, and $s\in\R$, the following diagram is commutative 
\[
\displaystyle{\begin{array}{cccl} 
\displaystyle \T^{n}&\stackrel{g_{\omega}^s}{\longrightarrow}&\T^{n}\cr
\downarrow\lefteqn{\Psi_{t,\omega}}& &\downarrow\lefteqn{\Psi_{t,\omega}} \cr
\displaystyle \Lambda_t(\omega)&\stackrel{\Phi_{t,\widetilde\omega}^s}{\longrightarrow}&\Lambda_t(\omega)&  
\end{array} } 
\]
\end{Remark}
\begin{Remark}\label{rem:kam} - \\ 
1. The Theorem could  be obtained  for any $k\in \N$ (then $C_m$ depends on $k$ as well). We suppose here that $ k\in \{0;1\}$ to simplify the proof. 
\\
2. We point out that the parameter  $\varepsilon>0$ does not depend on the parameters $\kappa$ and $r$, the domain $\Omega$, the annulus $\A^n= \T^n\times B(0,\rho_0)$, nor on the interval $[0,a]$.\\
3.  (ii) still holds if $P\in C^k([0,1];C_0^{\ell(M)}(\A^n\times\Omega))$ with $M\ge 0$ (see Theorem \ref{Theo:Holder}). \\
\end{Remark}
\begin{Remark}\label{rem:kappa}
	Without loss of generality one can assume that $\kappa=r=1$. Indeed, consider the $C^k$-family of Hamiltonians 
	$$
	\widetilde H_t= (\kappa r)^{-1}(N+P_t)\circ\sigma_{\kappa,r} = N+ \widetilde P_t,
	$$
	where 
	$$
	\widetilde P_t(\theta,I;\omega) = (\kappa r)^{-1}P_t(\theta, rI;\kappa\omega), \quad (\theta,I;\omega)\in \T^n\times \R^n\times (\kappa^{-1}\Omega).
	$$
	If $P_t$ satisfy \eqref{eq:smallness-condition3}, then so do $\widetilde P_t$ with $\kappa=r=1$. Let $\widetilde \phi_t$ and $\widetilde \Psi_t=(\widetilde U_t, \widetilde V_t)$ be the  family of maps obtained  by Theorem \ref{Theo:A} for the family of Hamiltonians $\widetilde H_t$  with $\kappa=\rho=1$. 
	Then taking 
	$$
	\phi_t:=\kappa \widetilde \phi_t\circ\sigma_{\kappa}^{-1}, \quad  \Psi_t:=(\widetilde U_t, r \widetilde V_t)\circ\sigma_{\kappa,r}^{-1}
	$$ 
	we obtain items (i)-(iii) in Theorem  \ref{Theo:A}  for $H_t$ and for   $0<\kappa\le 1$, $0<r \le 1$. 
\end{Remark} 	
In order to avoid the repeating use of the parameters $\kappa$ and and $\rho$, we suppose from now on that 
\begin{equation}\label{eq:kappa-rho}
\kappa=\rho=1.
\end{equation}

\noindent
{\em Idea oh the Proof.}
The proof of Theorem \ref{Theo:A} and Remark \ref{rem:analiticity-t} is organized as follows. In Sect. \ref{Sec:kamstep} we prove the KAM Lemma and choose the parameters for the next iteration. The KAM Lemma is close to that of P\"oschel in \cite{Poe1} but one needs additional arguments to estimate the derivatives with respect $t$. To this end we give a complete prove of it skipping some details. 
 In Sect. \ref{Sec:Iteration} we iterate the KAM Step infinitely many times. The choice of the parameters leads to an exponentially converging scheme. Additional efforts are needed to  get convergence for the derivatives with respect to $t$ and to obtain the corresponding estimates in the Iterative Lemma.  The iteration procedure is convergent in a Whitney sense only on the Cantor set $\Omega_{ \kappa}$ and one can not hope to get the global (in $\Omega$) estimates \eqref{eq:main-estimates} using Whitney's extension theorem for $C^\infty$ jets. For this reason we propose a new method in Sect. \ref{Sec:modified-iteration-lemma}. Using suitable almost analytic extensions in Gevrey classes, we prove a Modified Iterative Lemma which provides a convergent scheme  over the whole domain $\Omega$ and yields the desired estimates. The almost analytic extensions is obtained in Sect. \ref{Sec:Gevrey}.

\section{Proof of Theorem \ref{Theo:A}}\label{Sec:Proof-Theorem A}
\subsection{The KAM Step} \label{Sec:kamstep}
\subsubsection{\it The KAM Lemma}\label{Sec:kamlemma} 
Given two domains $D_j\subset \C^{n_j}$, $j=1,2$,  we denote by ${\mathcal A}(D_1,D_2)$ the space of  analytic maps $f:D_1\to D_2$ equipped by the inductive topology generated  by  sup-norms on compact sets of $D_1$,  and by $C^k([0,a], {\mathcal A}(D_1,D_2))$, $k\in\N$,  the corresponding space of the $C^k$ functions.   If $D_2=\C$ we write ${\mathcal A}(D_1):={\mathcal A}(D_1,\C)$.
Recall that an analytic function $f\in {\mathcal A}(D_1)$ is said to be real analytic if $D_1\cap \R^{n_1}\neq \emptyset$ and $f(D_1\cap \R^{n_1})\subset \R$. 
Introduce the complex domains 
$$ 
D_{s,r} = \{\theta\in {\C}^n/ 2\pi {\Z}^n:\ |{\rm Im}\,  \theta|<s\} 
          \times \{I\in {\C}^n:\ |I|<r\}, 
$$ 
$$ 
O_{h} = \{\omega\in {\C}^n:\ |\omega - \Omega_1 |<h\}. 
$$ 
Hereafter, $\displaystyle |v|=|(v_1,\ldots,v_n)|=\sup_j |v_j|$ is the sup-norm of $v\in\C^n$. 
The sup-norm of functions in ${\cal V}:=D_{s,r}\times O_h$ will be denoted by 
$|\cdot |_{s,r,h}$ and the corresponding space of analytic functions in ${\cal V}$ by ${\mathcal A}({\mathcal V})$. 
We state below a variant of the KAM Lemma  following   P\"oschel \cite{Poe1}. It involves a small parameter $\varepsilon >0$ and several parameters $\sigma, s, r, \eta, K$ such that 
\begin{equation}
0<s,r<1,\ 0<\eta<1/8,\ 0< 5\sigma < s < 1,\ K\geq 1, 
\label{eq:constants}
\end{equation}
 as well as a positive $c_0=c_0(n,\tau)\le 1$ depending only on $n$ and $\tau$. We suppose that the following inequalities are satisfied 
\begin{itemize}
\item[{\rm (a)}]
$\varepsilon \le c_0    \eta r \sigma^{\tau + 1},$ 
\item[{\rm (b)}]
$\varepsilon \le c_0 h r\,  ,$ 
 \item[{\rm (c)}]  $\displaystyle  
 h \leq \frac{1} {2K^{\tau + 1}}$.  
 \end{itemize}
Moreover, we will require below the  inequality 
 \begin{itemize}
  \item[{\rm (d)}] $2h \le   \sigma^{\tau +1}$
\end{itemize}
which follows from (c) provided that $K\sigma \ge 1$. 
Fix $k\in \{0;1\}$. 
\begin{Prop}[KAM Step Lemma] \label{Prop:kam-step}
	
There is a positive  $c_0=c_0(n,\tau) < 1$ depending only on $n$ and $\tau$ such that,  for any   $\sigma, s,h, r,\eta,K$, $a>0$ and $\varepsilon>0$ satisfying \eqref{eq:constants} and (a)-(c) and for every real valued Hamiltonian $H=N+P$, where 
\[
N(I;\omega,t) = e(\omega,t) + \langle \omega, I\rangle \quad \mbox{and} \quad P\in C^k([0,a], {\mathcal A}(D_{s,r}\times O_h))
\]
satisfies  the estimate
\begin{equation} 
\label{eq:smallness}
 \sup_{0\le p\le k} \,  \sup_{t\in [0,a]} \, |\partial_t^pP_t|_{s,r,h}\ \leq\  \varepsilon, 
\end{equation}
the following holds. 
\begin{itemize}
\item[{\rm (1)}]                       
There exists a $C^k$ family of real analytic transformation ${\cal F}= (\Phi,\phi)$, where
\[
\Phi \in C^k([0,a], {\mathcal A}( D_{s-5\sigma,\eta r}\times O_{h/4} , D_{s,r}))\quad \mbox{and}\quad \phi\in C^k([0,a], {\mathcal A}( O_{h/4} , O_h))
\]
such that $H\circ {\cal F} = N_+ + P_+$ with 
\[
N_+(I;\omega,t)=e_+(\omega,t) + \langle \omega, I\rangle \quad \mbox{and} \quad P_+\in C^k ([0,a],{\cal A}(D_{s-5\sigma,\eta r}\times O_{h/4}))
\]
satisfying the estimate
\begin{equation} 
|\partial_t^p P_+(\cdot,t)|_{s-5\sigma,\eta r,h/4} 
\le C_0 \left(\frac{\varepsilon ^2}{  r\sigma^{\tau + 1}} 
+ (\eta^2 +\sigma^{-n} e^{-K\sigma})\varepsilon\right)  
                                   \label{eq:new-error} 
\end{equation} 
for any $t\in [0,a]$ and $0\le p\le k$, where $C_0=C_0(n,\tau)>0$  depends only on $n$ and $\tau$;
\item[{\rm (2)}] $\Phi(\theta,I;\omega,t) = (U(\theta;\omega,t),
V(\theta,I;\omega,t))$, where $V$ is affine linear with respect to $I$
and  the transformation $(\theta,I)\to \Phi(\theta,I ;\omega,t)$ 
is canonical  for each $(\omega,t)$ fixed. Moreover, for any $0\le p \le k$, $\alpha,\beta\in\N^n$, and $|\gamma|\le 1$ the maps 
$\Phi$ and $\phi$ satisfy the estimates 
\[
\begin{array}{lcrr}
\displaystyle |W\partial_t^p(\sigma\partial_\theta)^\alpha(r\partial_I)^\beta(\Phi(\theta,I;\omega,t) - (\theta,I) )| \le C_{\alpha,\beta}
 \frac{\varepsilon }{  r\sigma^{\tau + 1}}\, , \\[0.3cm]
\displaystyle   |(h\partial_\omega)^\gamma \partial_t^p(\phi_t - id )| \le C  \frac{ \varepsilon }{r} \, ,
\end{array}
\]
uniformly on $D_{s-5\sigma,\eta r}\times O_{h}\times [0,a]$ 
and  $O_{ h/4}\times [0,a]$, 
respectively, where 
$W = {\rm diag \, }\left(\sigma^{-1}{\rm Id \, }, r^{-1}{\rm Id \, }\right)$, $C_{\alpha,\beta}>0$ depends only on $n$, $\tau$, $\alpha$, $\beta$,  and  $C>0$ depends only on  $n$ and  $\tau$. 
\end{itemize}
\end{Prop} 
\begin{Remark}\label{Rem:estimates}      Set $\overline W =  
{\rm diag \, }\left(\sigma^{-1}{\rm Id \, }, r^{-1}{\rm Id \, },   h^{-1}{\rm Id \, }\right)$ and suppose that (d) holds, i.e. $2h \le  
\sigma^{\tau +1}$. Then 
\[
|\overline W \partial_t^p(D{\cal F}(\cdot,t)  - {\rm Id} )
\overline W^{-1}| \le C_0  \frac{\varepsilon }{r h}  
\] 
on $D_{s-5\sigma,\eta r}\times O_{h/4}\times [0,a]$, where $D{\cal F}(\cdot,t)$ stands for the Jacobian of ${\cal F}(\cdot,t)$. Moreover, 
   (2) and the  Cauchy estimate of the derivatives of ${\cal F}$  with respect to
$\omega$ yield for $0\le p\le 1$  and any $\alpha,\beta,\gamma\in\N^n$ the estimate
\[
\begin{array}{lcrr}
\displaystyle |\overline W\partial_t^p(\sigma\partial_\theta)^\alpha(r\partial_I)^\beta(h\partial_\omega)^\gamma({\cal F}(\theta,I;\omega,t) - (\theta,I;\omega ))|\ \le C_{\alpha,\beta,\gamma} \frac{\varepsilon }{r h}  
\end{array}
\]
on $D_{s-5\sigma,\eta r}\times O_{h/6}\times [0,a]$, where $C_{\alpha,\beta,\gamma}>0$  depends only on $n$, $\tau$, $\alpha$, $\beta$ and $\gamma$. 
\end{Remark}
\begin{Remark}\label{rem:analytisity-t-KAMstep}
If $P$ is analytic with respect to $t$ in the disc $B(0,a)\subset \C$  and \eqref{eq:smallness} holds in $B(0,a)$  for $k=0$, then $\Psi$ is analytic with respect to $t$ in $B(0,a)$  and items (1) and (2) hold for $t\in B(0,a)$ with $p=0$. 
\end{Remark}
\begin{Remark}\label{Rem:Cauchy} 
Hereafter we use the  Cauchy estimates for analytic functions  in $\C^n$ (see for example Theorem 2.2.7, \cite{Hor2} and Appendix A in \cite{Poe1} ). More precisely, let $D$  be a domain in $\C^n$ and $D_r:=\{z\in \C^n:\, |z-D|<r\}$ the corresponding polydisc. Then for any analytic function $f$ in $D_r$ with a bounded sup-norm $|f|_r:=\sup_{z\in D_r}|f(z)|$ and any $0\le \rho<r$ and $\alpha\in\N^n$ one has
\begin{equation}
|\partial^\alpha f|_{\rho} \ \le\ \alpha! (r-\rho)^{-|\alpha|-1} |f|_{r} . 
\label{eq:Cauchy1}
\end{equation}
\end{Remark}
We recall  as well the standard estimates of the Fourier coefficients 
\[
f_k = \int_{\T^n}f(\theta)e^{-2\pi i \langle k,\theta\rangle} d\theta\ ,\quad k \in \Z^n, 
\]
of an analytic function $f$ in a strip $\T^n + s:=\{\theta\in \C^n/ \Z^n:\ |{\rm Im\, } \theta| < s\}$, $s>0$,  with a bounded sup-norm $|f|_s$, namely,  
\begin{equation}
|f_k| \ \le\ e^{-|k|s} |f|_{s} \, , 
\label{eq:Fourier-coeff}
\end{equation}
where $|k| = \sum_{j=1}^n |k_j|$. \\

\noindent
{\em Proof of Proposition \ref{Prop:kam-step}}. 
For Hamiltonians independent of $t$ the proposition is formulated and proved in
 \cite{Poe1}. It follows easily from \cite{Poe1} in the case $k=0$.   The proof of the corresponding estimates of $\partial_tP_t$ requires additional efforts.
 For this reason   we 
give a complete  proof in the case $k=1$.\\

\noindent
{\em Step 1. Truncation. }  Consider the linear part   of $P$ with respect to $I$
\[
Q(\theta,I;\omega,t):= P(\theta,0;\omega,t) + \langle \nabla_I P(\theta,0;\omega,t), I\rangle .
\]  
 Given a positive integer $K$ we denote by 
\[
R(\theta,I;\omega,t) := \sum_{|k|\le K} R_k(I;\omega,t) e^{i\langle k,\theta\rangle}
\]
the trigonometric polynomial of degree $K$ in the Fourier series expansion of $Q$ with respect to $\theta$. By \eqref{eq:smallness} and   the Cauchy inequalities \eqref{eq:Cauchy1} one  obtains the following estimates 
\[
|\partial_t^p Q|_{s,r} < C_0 \varepsilon \, ,\quad |\partial_t^p(P-Q)|_{s,2\eta r} < C_0 \eta^2 \varepsilon
\] 
for $0\le p\le 1$ uniformly with respect to $(\omega,t)\in O_h\times [0,a]$ (recall that $0<\eta<1/8$). 
Moreover, estimating the Fourier coefficients of $\partial_t^p Q$ by \eqref{eq:Fourier-coeff} one obtains
\begin{equation}
|\partial_t^p(Q-R)|_{s-\sigma,r} < C_0 \sigma^{-n} e^{-K\sigma} \varepsilon \quad \mbox{and} \quad |\partial_t^p R |_{s-\sigma,r} < C_0  \varepsilon . 
\label{eq:estimates-R}
\end{equation} 
The Cauchy estimates imply
\begin{equation}
\begin{array}{lcrr}
|\partial_t^p (\sigma\partial_\theta)^\alpha(r\partial_I)^\beta R |_{s-2\sigma,r/2} \le C_{\alpha,\beta}  \varepsilon \\[0.3cm]
|\partial_t^p (\sigma\partial_\theta)^\alpha(r\partial_I)^\beta (P-R) |_{s-2\sigma,2\eta r} \le C_{\alpha,\beta} (\eta^2 + \sigma^{-n} e^{-K\sigma})\varepsilon
\end{array}
\label{eq:nabla-R}
\end{equation}
for $0\le p\le 1$ uniformly in $(\omega,t)\in O_h \times [0,a]$. Hereafter  $C_0\ge 1$ stands for a constant depending only on $n$ and $\tau$ and we denote by   $C_{\alpha,\beta}$  a positive constant depending   only on $n$, $\tau$, $\alpha$ and $\beta$.  \\

\noindent
{\em Step 2. Homological equation. } The idea is to  put $\partial_t^p(P-R)$  in the error term and to  to kill $\partial_t^p R$ by means of a canonical transformation $\Phi$ which is  the time-one-map of
a Hamiltonian vector field $X_F= (\nabla_I F, - \nabla_\theta F)$. More precisely,  consider the Hamiltonian flow 
\[
(x,\theta,I) \to \exp(x X_F)(\theta,I)= (u(x,\theta,I), v(x,\theta,I)) 
\] 
and set 
\[
\Phi= (U,V):= \exp(X_F), \quad  \mbox{where}\  U(\cdot)=u(1,\cdot)\  \mbox{and}\ V(\cdot)=v(1,\cdot). 
\]
The corresponding Hamiltonian system is
\begin{equation}
\left\{
\begin{array}{lcrr}
\displaystyle \frac{d u}{dx} = \nabla_v F(u,v;\omega,t)\\[0.3cm]
\displaystyle \frac{d v}{dx} = - \nabla_u F(u,v;\omega,t)\\[0.3cm]
u(0)=\theta\, ,\ v(0)= I\, . 
\end{array}
\right.
\label{eq:hamiltonian-system}
\end{equation}
The Lie method is based on the identity
\[
\frac{d}{d x}(f\circ \exp(x X_F) )= \{f,F\}\circ \exp(x X_F)\, ,
\]
where $\{f, F \}= \langle \nabla_I f, \nabla_\theta F \rangle - \langle \nabla_\theta f, \nabla_I F \rangle$ is the  Poisson bracket. Using Taylor's formula with respect to $x$ at $x=0$ and  the above identity  one gets
\begin{equation}
\begin{array}{lcrr}
\displaystyle  (N+R)\circ \Phi = N\circ \exp (x X_F)|_{x=1} + R\circ \exp (x X_F)|_{x=1} \\[0.3cm] 
\displaystyle  = N + \{N,F\} +R + \int_0^1\{(1-x)\{N,F\} + R,F\}\circ \exp(x X_F)\,  dx .
\end{array}
\label{eq:taylot-to-homological}
\end{equation}

We are looking for a trigonometric polynomial $F$ of degree $K$  and for a function $\hat N$ depending only on $(I,\omega,t)$ solving the homological equation
\begin{equation}
\{N, F \}  + R =  \hat N . 
\label{eq:homological2}
\end{equation}
Recall that $ N(I;\omega,t) = e(\omega,t) + \langle \omega, I\rangle$. Then \eqref{eq:homological2} becomes
${\cal L}_\omega F = \hat N - R$, 
where ${\cal L}_\omega =\langle \omega, \partial/\partial \theta\rangle$. 
Take
\[
\hat N(I;\omega,t) := R_0(I;\omega,t) = \int_{\T^n}R(\theta,I;\omega,t) d\theta\, ,
\]
which is affine linear in $I$. Then the zero order term of the trigonometric polynomial $\hat N-R$ is zero which is a necessary condition for solving the  above equation.  On the other hand, 
the Diophantine condition \eqref{eq:sdc2} with $\kappa=1$  and (c) imply
\begin{equation}                
|\langle\omega,k\rangle |\ \ge \ \frac{1}{2|k|^{\tau}}\quad \mbox{for all} \  \omega\in O_h \ \mbox{and}\ 0 \neq |k|\le K ,
                       \label{eq:sdc3}                   
\end{equation}
where $|k| = \sum_{j=1}^n |k_j|$.   Denote by  $\mathcal{H}_K$ the space 
 of trigonometric polynomials in $\theta\in \T^n$ of degree $\le K$ with zero order terms equal to $0$. This space is generated by the functions $\exp(i\langle k,\theta\rangle)$, where $k\in\Z^n$ and $0<|k|=|k_1|+ \cdots + |k_n|\le K$. It follows from \eqref{eq:sdc3} that 
  the map  ${\cal L}_\omega: \mathcal{H}_K \to \mathcal{H}_K$ is an automorphism. 
 Denote by ${\cal L}_\omega^{-1}: \mathcal{H}_K \to \mathcal{H}_K$ the inverse map and set
 \begin{equation} \label{eq:F} 
 F:= {\cal L}_\omega^{-1}(R-\hat N). 
 \end{equation}
 The Fourier coefficients of 
$F$ are $F_0=0$, $F_k= (i\langle \omega,k\rangle)^{-1}R_k$ for $0<|k|\le K$ and $F_k=0$ for $|k|>K$. Hence, $F$
is well defined, it solves \eqref{eq:homological2} and  is affine linear in $I$. Moreover, it is uniquely defined by
\begin{equation}
\int_{\T^n}\, F(\theta,I;\omega,t)\, d\theta \ = F_0(I;\omega,t)= 0 . 
\label{eq:uniqueness-of-U}
\end{equation}
Now \eqref{eq:taylot-to-homological} reads
\begin{equation}
(N+R)\circ \Phi 
= N + \hat N  + \int_0^1\{(1-x)\hat N + xR,F\}\circ \exp(x X_F) dx . 
\label{eq:taylot-to-homological1}
\end{equation}
Moreover, \eqref{eq:Fourier-coeff},  \eqref{eq:estimates-R} and \eqref{eq:sdc3}  imply 
\begin{equation}
|\partial_t^p F |_{s-2\sigma,r} \le C_0 \frac{|\partial_t^p R |_{s-\sigma,r}}{  \sigma^\tau} < C_0  \frac{\varepsilon}{  \sigma^\tau}. 
\label{eq:estimates-F0}
\end{equation}
Using the Cauchy estimates one gets 
\begin{equation}
|\partial_t^p (\sigma\partial_\theta)^\alpha(r\partial_I)^\beta F | \le C_{\alpha,\beta}  \frac{\varepsilon}{   \sigma^{\tau}} 
\label{eq:nabla-F}
\end{equation}
uniformly in $D_{s-3\sigma,r/2}\times O_h \times [0,a]$ and for any $0\le p\le 1$ and $\alpha,\beta \in \N^n$, $|\beta|\le 1$.  By   \eqref{eq:Fourier-coeff},  \eqref{eq:estimates-R} and Cauchy one has as well
\begin{equation}
|\partial_t^p \hat N|_{r} = |\partial_t^p R_0|_{r} \le |\partial_t^p R|_{s-\sigma,r} \le  C_0  \varepsilon \quad \mbox{and}\quad |\partial_t^p (r\partial_I) \hat N|_{r/2} \le  C_0  \varepsilon 
\label{eq:estimates-derivatives-N}
\end{equation}
for $0\le p\le 1$ uniformly with respect to $(\omega,t)\in O_h\times [0,a]$.  The derivatives of   $F$ and $\hat N$ with respect $I$ of order bigger than one are all zeros since the functions are affine linear in $I$. \\

\noindent
{\em Step 3. Canonical transformation. } The solution $(u,v)$  of the Cauchy problem \eqref{eq:hamiltonian-system} are real analytic in $(x,\theta,I,\omega)$ and $C^1$ in $t$. 
Consider the canonical transformation $\Phi=(U,V)$, where $U(\cdot)=u(1,\cdot)$ and $V(\cdot)=v(1,\cdot)$ are defined in Step 2.  Since $F$ is affine linear in $v$ one observes that $u$ is independent of $I$ and $v$ is affine linear in $I$. In particular, $U$ is independent of $I$ and $V$ is affine linear as a function of $I$. 
 Moreover, \eqref{eq:nabla-F}  and condition (a) imply for $p\in \{0;1\}$  the inequality
\begin{equation}
|\partial_t^p \nabla_\theta F |\le \eta r \le r/8\, \quad \mbox{and}\quad |\partial_t^p \nabla_I F | \le \sigma  
\label{eq:estimates-nabla-F}
\end{equation}
in $D_{s-3\sigma,r/2}\times O_h \times [0,a]$  choosing the constant $c_0=c_0(n,\tau)<1$ in (a) sufficiently small. Then 
\begin{equation}
\exp(x X_H): D_{s-4\sigma,r/4} \to D_{s-3\sigma,r/2}
\label{eq:range-exponential-map}
\end{equation}
for every $(\omega,t)\in O_h\times [0,a]$ and $0\le x\le 1$. In particular, 
$\Phi(\cdot;\omega,t)=(U(\cdot;\omega,t),V(\cdot;\omega,t))$ is a well defined real analytic  map
\[
\Phi(\cdot;\omega,t): D_{s-4\sigma,r/4} \to D_{s-3\sigma,r/2}
\]
for every $(\omega,t)\in O_h\times [0,a]$ and we get $\Phi\in C^1([0,a], {\mathcal A}(D_{s-4\sigma,r/4}\times O_h, D_{s-3\sigma,r/2}))$. 
We are going to show that
\begin{equation}
\left\{
\begin{array}{lcrr}
\displaystyle \frac{1}{\sigma}|\partial_t^p (u(x,\theta,I;\omega,t)-\theta)| \le C_0  \frac{\varepsilon}{   r\sigma^{\tau+1}}\\[0.3cm]
\displaystyle \frac{1}{r}|\partial_t^p (v(x,\theta,I;\omega,t)-I)|  \le C_{0}  \frac{\varepsilon}{  r \sigma^{\tau+1}}
\end{array}
\right.
\label{eq:estimates-u-v}
\end{equation}
in $[0,1]\times D_{s-4\sigma,r/4}\times O_h\times [0,a]$. For $p=0$ it follows directly from \eqref{eq:hamiltonian-system} and \eqref{eq:nabla-F}. Let $p=1$. Set 
\[
\overline{u}(\theta,I;\omega,t):= \frac{1}{\sigma} \sup_{0\le x\le 1}|\partial_t u(x,\theta,I;\omega,t)| \quad \mbox{and} \quad \overline{v}(\theta,I;\omega,t):= \frac{1}{r} \sup_{0\le x\le 1}|\partial_t v(x,\theta,I;\omega,t)| .
\] 
Differentiating \eqref{eq:hamiltonian-system} with respect to $t$ and using \eqref{eq:nabla-F} one gets  
\[
\begin{array}{lcrr}
\displaystyle \overline{v}\ \le \ \frac{1}{\sigma r} \sup_{0\le x\le 1}|\partial_t  (\sigma \nabla_\theta) F | + \frac{1}{\sigma r}\sup_{0\le x\le 1}|  (\sigma\partial_\theta)  (\sigma\nabla_\theta)F |\,  \overline{u} + \frac{1}{\sigma r} \sup_{0\le x\le 1}|  (r \partial_I) (\sigma\nabla_\theta F |\,  \overline{v} \\[0.5cm]
\displaystyle \le \ C \left( \frac{\varepsilon}{  r \sigma^{\tau+1}} + \frac{\varepsilon}{  r \sigma^{\tau+1}}\,  \overline{u} + \frac{\varepsilon}{  r \sigma^{\tau+1}}\,  \overline{v}\right)
\end{array}
\]
where $C=C(n,\tau)>0$ and by (a)  one obtains
\[
\displaystyle \overline{v}\ \le \ C \frac{\varepsilon}{  r \sigma^{\tau+1}} + Cc_0(\overline{u} + \overline{v}) 
\]
in $D_{s-4\sigma,r/4}\times O_h\times [0,a]$. 
Choosing $c_0\le (4C)^{-1}$ this gives
\[
\displaystyle \overline{v}\ \le \ C \frac{\varepsilon}{  r \sigma^{\tau+1}} + \frac{1}{4}(\overline{u} + \overline{v}). 
\]
The same estimate holds for $\overline{u}$ and we get
 \eqref{eq:estimates-u-v}. 
By Cauchy this implies 
\begin{equation}
\left\{
\begin{array}{lcrr}
\displaystyle \frac{1}{\sigma}|\partial_t^p (\sigma\partial_\theta)^\alpha(r\partial_I)^\beta (u(x,\theta,I;\omega,t)-\theta)| \le C_{\alpha,\beta}  \frac{\varepsilon}{   r\sigma^{\tau+1}}\\[0.3cm]
\displaystyle \frac{1}{r}|\partial_t^p (\sigma\partial_\theta)^\alpha(r\partial_I)^\beta (v(x,\theta,I;\omega,t)-I)| \le   C_{\alpha,\beta}  \frac{\varepsilon}{  r \sigma^{\tau+1}}
\end{array}
\right.
\label{eq:derivatives-u-v}
\end{equation}
 in $[0,1]\times D_{s-5\sigma,r/8}\times O_h\times [0,a]$.  Since $\eta<1/8$ this proves the estimates of $\Phi$ in statement (2) of the KAM step. By \eqref{eq:estimates-nabla-F} we get
\[
\begin{array}{lcrr}
\displaystyle
|U(\theta,I;\omega,t)-I| \le \sup_{0\le x\le 1} |\nabla_I F(\cdot;\omega,t |_{s-3\sigma, r/2}  \le \sigma , \\[0.3cm]
\displaystyle |V(\theta,I;\omega,t)-I| \le \sup_{0\le x\le 1} |\nabla_\theta F(\cdot;\omega,t |_{s-3\sigma, r/2}  \le \eta r
 \end{array}
\]
on $D_{s-5\sigma,\eta r}\times O_h\times [0,a]$. 
This implies   that $\Phi(\cdot;\omega,t)$ maps $D_{s-5\sigma,\eta r}$ to $D_{s-4\sigma,2\eta r}$,  and that 
\begin{equation}
\Phi\in C^1([0,a], {\mathcal A}(D_{s-5\sigma,\eta r}\times O_h, D_{s-4\sigma,2\eta r})). 
\label{eq:class-of-Phi}
\end{equation}

\noindent
{\em Step 4. New error term. } The identity \eqref{eq:taylot-to-homological1} yields
\begin{equation}
H\circ \Phi=  (N+R)\circ \Phi + (P-R)\circ \Phi = N_+ + P_+
\label{eq:transformed-hamiltonian}
\end{equation}
where the Hamiltonian $N_+=N +\hat N$ is independent of $\theta$ and affine linear in $I$ and 
\begin{equation}
P_+ =  \int_0^1\{(1-x)\hat N + xR,F\}\circ \exp(x X_F) dx + (P-R)\circ \Phi
\label{eq:new-error-term}
\end{equation}
is the new error term.

We are going to prove \eqref{eq:new-error}. In the case $p=0$ it follows from the corresponding estimates in \cite{Poe1}. 
Take $p=1$ and consider firstly 
\[
\partial_t( (P-R)\circ \Phi) = (\partial_t(P-R))\circ \Phi + (D(P-R)\circ \Phi) \, \partial_t \Phi,
\]
where $D$ stands for the differential with respect to $(\theta,I)$. 
By \eqref{eq:class-of-Phi} and \eqref{eq:nabla-R} we have 
\[
|\partial_t^p (P-R)\circ \Phi|_{s-5\sigma, \eta r} \le |\partial_t^p (P-R)|_{s-4\sigma, 2\eta r} 
\le C_0 (\eta^2 + \sigma^{-n} e^{-K\sigma})\varepsilon .
\]
Moreover, \eqref{eq:nabla-R} implies 
\[
|D(P-R)W^{-1}|_{s-4\sigma,2\eta r} \le C_0 (\eta^2 + \sigma^{-n} e^{-K\sigma}), 
\]
while \eqref{eq:derivatives-u-v} gives 
$|W\partial_t \Phi| \le C_0 \varepsilon/  r \sigma^{\tau+1}$ on $D_{s-5\sigma,r/8}$, and we get
\[
\begin{array}{lcrr}
|(D(P-R)\circ \Phi) \, \partial_t \Phi|_{s-5\sigma,\eta r} \le |D(P-R)W^{-1}|_{s-4\sigma,r/2}|W\partial_t \Phi|_{s-5\sigma,r/8}\\[0.3cm] 
\displaystyle \le  C_0 (\eta^2 \varepsilon + \sigma^{-n} e^{-K\sigma}\varepsilon)\frac{\varepsilon}{  r \sigma^{\tau+1}} 
\end{array}
\]
uniformly with respect to $(\omega,t)\in O_h\times [0,a]$.
To evaluate  the derivative with respect to $t$ of the first term  in \eqref{eq:new-error-term}, we consider 
\[
G:= \partial_t(\{R,F\}\circ \exp(x X_F))=G_1+G_2+G_3,
\]
where
\[
 G_1:= \{\partial_tR,F\}\circ \exp(x X_F)\, ,\quad G_2:= \{R,\partial_t  F\}\circ \exp(x X_F) 
\]
and
\[
G_3:= \big(D\{R,F\}\circ \exp(x X_F)\big)\,   . \partial_t \exp(x X_F). 
\]
Using \eqref{eq:range-exponential-map} one obtains 
\[
|\{\partial_t^p R, \partial_t^q F\}\circ \exp(x X_F)|_{s-5\sigma, \eta r} 
\le |\{\partial_t^p R, \partial_t^q F\}|_{s-3\sigma, r/2}
\]
for $0\le p,q\le 1$. Now  \eqref{eq:nabla-R} and \eqref{eq:nabla-F}  imply
\[
\begin{array}{lcrr}
\displaystyle |\{\partial_t^p R, \partial_t^q F\}|_{s-3\sigma, r/2}  \le 
 |\partial_t^p \nabla_I R||\partial_t^q \nabla_\theta  F| + |\partial_t^p \nabla_I F||\partial_t^q \nabla_\theta  R|\\[0.3cm] 
 \displaystyle \le C_0\left(\frac{\varepsilon}{r}\cdot \frac{\varepsilon}{  \sigma^{\tau+1}} + \frac{\varepsilon}{  r \sigma^{\tau}}\cdot \frac{\varepsilon}{ \sigma} \right)
 =2 C_0 \frac{\varepsilon^2}{  r \sigma^{\tau+1}}
\end{array}
\]
uniformly with respect to $(\omega,t)\in O_h\times [0,a]$, 
which gives the desired estimate for $G_1$ and $G_2$. 
By the same argument one obtains 
\[
 |(\sigma\partial_\theta)^\alpha(r\partial_I)^\beta\{ R, F\}|_{s-3\sigma, r/2} 
 \le C_{\alpha,\beta}\left(\frac{\varepsilon}{r}\cdot \frac{\varepsilon}{  \sigma^{\tau+1}} + \frac{\varepsilon}{  r \sigma^{\tau}}\cdot \frac{\varepsilon}{ \sigma} \right)
 =2 C_{\alpha,\beta} \frac{\varepsilon^2}{  r \sigma^{\tau+1}}.
\label{eq:estimates-poisson-bracket}
\]
Using \eqref{eq:derivatives-u-v} and the preceding estimate one gets
\[
\begin{array}{lcrr}
\displaystyle 
|G_3|_{s-5\sigma, \eta r} \le |(\sigma D_\theta)\{R,F\}|_{s-3\sigma, r/2}\,  |\sigma^{-1}\partial_t u|_{s-5\sigma,\eta r}\\[0.3cm]
 \displaystyle+ |(r D_I)\{R,F\}|_{s-3\sigma, r/2} \, |r^{-1} \partial_t v|_{s-5\sigma,\eta r} 
 \le C_0 \frac{\varepsilon^2}{  r \sigma^{\tau+1}}
\end{array}
\]
where $D_\theta$ and $D_I$ are the partial differentials with respect to $\theta$ and $I$ respectively. 
The function $\partial_t(\{\hat N,F\}\circ \exp(x X_F))$ can be evaluated in the same way using \eqref{eq:estimates-derivatives-N}. 
This proves \eqref{eq:new-error}. \\

\noindent
{\em Step 5. Transforming the frequencies. } Consider 
\[
N(I;\omega,t)+ \hat N(I;\omega,t) = e_+(\omega,t) + \langle \omega + (\nabla_I \hat N)(\omega,t),I\rangle = e_+(\omega,t) + \langle \omega + (\nabla_I R_0)(\omega,t),I\rangle. 
\]
Following P\"oschel \cite{Poe1}, Sec. 4, we 
obtain a real analytic inverse  $\phi_t: O_{h/4} \to O_{h/2}$ of the map 
\[
\omega\to \omega_+:=\omega + (\nabla_I R_0)(\omega,t)
\] 
i.e. 
\begin{equation}\label{eq:inverse-map-phi}
\phi_t(\omega) + (\nabla_I R_0)(\phi_t(\omega),t)=\omega  \, , \quad \omega\in O_{h/4},\ t\in [0,a]. 
\end{equation}
Moreover, the following estimate is true
\[
 |\phi_t - id | +h |D\phi_t - {\rm Id}|
\le C \frac{ \varepsilon }{r}\, ,
\]
on $O_{h/4}$. We set $\phi(\cdot,t)=\phi_t$ and $N_+=(N+\hat N)\circ \phi$.   

We are going to estimate $\partial_t\phi$ on $O_{h/4}$. Using \eqref{eq:estimates-derivatives-N} and Cauchy we obtain the estimate
\begin{equation}\label{eq:estimate-N}
|(h\partial_\omega)^\gamma\partial_t^p \nabla_I R_0(\omega,t)| \le C_\gamma  \frac{ \varepsilon }{r} \quad \mbox{in} \  (\omega,t)\in O_{h/4}\times [0,a]
\end{equation}
for each $0\le p\le 1$ and $\gamma\in \N^n$ (recall that $R_0$ is affine linear in $I$). In particular, using (b) we obtain
\begin{equation}\label{eq:estimate-N1}
| D \nabla_I R_0(\omega,t)| \le C(n,\tau) \frac{ \varepsilon }{hr} \le C(n,\tau) c_0(n,\tau) <\frac{1}{2}  \quad \mbox{in} \  (\omega,t)\in O_{h/4}\times [0,a]
\end{equation}
for $c_0$ small enough. 
Differentiating \eqref{eq:inverse-map-phi} we get
\begin{equation}\label{eq:inverse-map-dphi}
\partial_t \phi_t(\omega) + (D \nabla_I R_0)(\phi_t(\omega),t)\, .\partial_t \phi_t(\omega)+ (\nabla_I \partial_t R_0)(\phi_t(\omega),t)=0
\end{equation}
Using \eqref{eq:estimate-N} and \eqref{eq:estimate-N1}  we get the  estimate 
\[
|\partial_t \phi_t(\omega)|\le C(n,\tau)\frac{\varepsilon}{r} \le C_0(n,\tau)h  \quad \mbox{in} \  (\omega,t)\in O_{h/4}\times [0,a]
\]
Differentiating \eqref{eq:inverse-map-dphi} with respect to $\omega$ we obtain
\[
(h\partial_{\omega_j})\partial_t \phi_t(\omega) + (D \nabla_I R_0)(\phi_t(\omega),t)\, .(h\partial_{\omega_j})\partial_t \phi_t(\omega) = Q_t(\omega) 
\]
where 
\[
\begin{array}{rcl}
\displaystyle 
Q_t =& -&(hD)^2 (\nabla_I \partial_t^p R_0)(\phi_t,t)[\partial_{\omega_j}\phi_t, h^{-1}\partial_t \phi_t] - (h\partial_{\omega_j})(hD )\nabla_I R_0(\phi_t,t)\, .h^{-1}\partial_t \phi_t \\
       &-& (h\partial_{\omega_j})(\nabla_I \partial_t R_0)(\phi_t,t) - hD(\nabla_I \partial_t R_0)(\phi_t,t)\, . \partial_{\omega_j}\phi_t
\end{array}
\]
Hereafter, $D^2f[\cdot,\cdot ]$ stands for the quadratic form representing the second differential of $f$. Using  \eqref{eq:estimate-N},  \eqref{eq:estimate-N1}  and (b),  we obtain 
\[
|(h\partial_{\omega_j}) \partial_t^p(\phi_t - id )| \le C  \frac{ \varepsilon }{r} \quad \mbox{in} \ O_{h/4}.
\]
 This completes the proof of the KAM Step Lemma.  \finishproof

\noindent
The analyticity with respect to $t$ in Remark \ref{rem:analytisity-t-KAMstep} follows from the theorem of Cauchy. \finishproof

Does the transformation ${\mathcal F}$  obtained by the KAM Step Lemma depend on the choice of the parameters $K$, $\sigma$, $h$, $r$, $\eta$ and how? Following the construction of  ${\mathcal F}$ we obtain the following 

\begin{Remark}[Uniqueness by construction in the KAM Lemma]\label{rem:uniqueness-KAM}
	The transformation ${\mathcal F}$, the new normal form $N_+$ and the error term $P_+$  depend on the choice of $K$ via the truncation in Step 1. 
	If  $K$ is fixed,  then  they do not depend on the choice of the other parameters $\sigma$, $h$, $r$ and $\eta$ in the following sense. 
	Let $\sigma'$, $h'$, $r'$ and $\eta'$ be another admissible choice of the parameters and ${\mathcal F}'$,  $N'_+$ and $P'_+$,  be the corresponding transformation,  normal form and  error term. 
	Then ${\mathcal F}'={\mathcal F}$,  $N'_+= N_+$  and $P'_+= P_+$ on the intersection of their domains of definition. 
\end{Remark}

\subsubsection{\it Preparing next iteration.}
\label{Sec:next-iteration}
 
We are going to prepare the next iteration. 
Choose 
a ``weighted error'' $E$ satisfying
\begin{equation}
\label{eq:E-and-eta}
0<E\le \eta^2<1/64  
\end{equation}
fix $0< \hat\varepsilon \le 1$ and set 
\begin{equation}
\varepsilon = \hat\varepsilon
      r \sigma^{\tau + 1} E. 
\label{eq:def-of-epsilon}
\end{equation} 
where  $0<\sigma<1/5$. 
Define $K$ and $h$ by
\begin{equation}
\label{eq:K}
K=\sigma^{-1} \ln^2(\sigma), \  2h = \frac{1}{K^{\tau +1}} =   \big(\sigma/\ln^{2}(\sigma)\big)^{\tau +1} . 
\end{equation}
\begin{Lemma}\label{lemma:conditions-a-d}
There exists  $E^0=E^0(n,\tau)>0$ depending only on $n$ and $\tau$ such that the hypothesis (a)-(d) of Proposition \ref{Prop:kam-step} are satisfied for   $0<\sigma<1/5$ and $0<E\le E^0(n,\tau)$ provided that
\begin{equation}
\displaystyle   2  \ln^{2\tau+2}(\sigma)E   \le c_0. 
\label{eq:condition-b}
\end{equation}
\end{Lemma}
{\em Proof}.
Firstly, a) follows from the definition of $\varepsilon$ choosing $E\le E^0\le c_0(n,\tau)^2$ and c) follows from the definition  $h$, while (d) follows from the inequality $K\sigma = \ln^2(\sigma)>\ln^2(5)>1$. The hypothesis  (b)  follows  from the inequality 
\begin{equation}\label{eq:b-again}
\frac{\varepsilon}{hr} \le \frac{   r \sigma^{\tau + 1} E}{hr} =2 \ln^{2\tau+2}(\sigma)E   \le c_0
\end{equation}
in view of \eqref{eq:def-of-epsilon} and \eqref{eq:K}, which yields 
 (b) in Proposition \ref{Prop:kam-step}. 
\finishproof

We are going to fix $\eta$ and determine the parameters   $s_+$, $\sigma_+$, $r_+$, $\eta_+$, $K_+$, $h_+$, $\varepsilon_+$, and the weighted error $E_+$ for the next iteration. 
Suppose that 
\begin{equation} \label{eq:sigma-eta}
\sigma^{-n} \exp(-K\sigma) = \sigma^{-n} \exp\left(-\ln^2(\sigma)\right) \le  \eta^2 . 
\end{equation}
Then using \eqref{eq:E-and-eta} and \eqref{eq:def-of-epsilon}, one obtains from  \eqref{eq:new-error}    the following inequality
\begin{equation}
\label{eq:new-error1}
|\partial_t^p P_+|_{s-5\sigma,\eta r,h/4}  \le  C_0   \hat\varepsilon   r \sigma^{\tau + 1} E \left( E
+ \eta^2 +\sigma^{-n} e^{-K\sigma} \right) < 3 C_0 \hat\varepsilon \eta^2     r \sigma^{\tau + 1} E :=\frac{1}{2}\varepsilon_+,
\end{equation}
where $0\le p \le k$ and $C_0=C_0(n,\tau) > 1$  depends only on $n$ and $\tau$. 
Set 
\begin{equation}
\label{eq:r-and-s}
r_+ = \eta r,\ s_+ = s-5\sigma,\ \sigma_+ = \delta\sigma,\ s = \frac{5}{1-\delta}\sigma,
\end{equation}
where  $0<\delta< 1/6$ will be fixed below and put $\varepsilon_+:=   \hat\varepsilon    r_+\sigma_+^{\tau + 1}  E_+$. 
Plugging the expression of $\varepsilon_+$ in  \eqref{eq:new-error1} and using \eqref{eq:r-and-s} we get 
\[
E_+ =  \left(6 C_0(n,\tau)    \delta^{-\tau -1}\right)\eta E
\]
which leads to an exponentially converging iteration scheme if  $\displaystyle 6C_0(n,\tau)\delta^{-\tau -1}\eta <1$. Now we fix
\begin{equation}
\label{eq:delta+eta}
 0<\vartheta < \min(\vartheta_0/4, 1),  \quad \delta := (6 C_0(n,\tau))^{-\frac{1}{\vartheta}},   \quad  \eta := \delta^{\tau +1 +\vartheta+\nu}, 
\end{equation}
where  $\nu$  is a positive number which will be determined   in Sect. \ref{Sec:Iteration} and $\vartheta_0>1$ is fixed in \eqref{eq:ell-def}. 
In particular,  \eqref{eq:delta+eta} implies  that 
\[
0<\delta  < 1/(6C_0)< 1/6 \quad \mbox{and} \quad  \eta < \delta^{2} < 1/6,
\] 
since $C_0>1$ and $\tau>n-1\ge 1$.  Moreover, 
\begin{equation}
\label{eq:E+}
E_+ =  \delta^{\nu} E. 
\end{equation}
Set $\eta_+=\eta \delta^{\nu_+-\nu}$ with certain $\nu_+\ge \nu$ which will be determined by the next iteration  and put $h_+=(1/2)K_+^{-\tau -1}$, where  $K_+=\sigma_+^{-1}\ln^2(\sigma_+)$.
Notice that 
\[
s_+ = s-5\sigma =\delta s, \quad s_+-5\sigma_+ = \delta (s-5\sigma) >0,
\]
and one obtains that   $\sigma_+, s_+,  r_+,\eta_+$ and $K_+$ satisfy  \eqref{eq:constants}. 
Moreover, 
\begin{equation}
\label{eq:h+}
\frac{h_+}{h} < \left(\frac{\sigma_+}{\sigma}\right)^{(\tau +1)} < \delta^{\tau +1} < \frac{1}{6}, 
\end{equation}
and \eqref{eq:new-error1} implies
\begin{equation}
|\partial_t^p P_+|_{s_+,r_+,h_+}\ \leq\ \frac{1}{2} \varepsilon_+ = \frac{1}{2}   \hat\varepsilon 
   r_+\sigma_+^{\tau + 1}  E_+  
                           \label{eq:new-error2}
\end{equation}
for $0\le p \le k$. 
We have prepared  the next iteration.

\subsection{Iteration}
\label{Sec:Iteration}
\subsubsection{\it Choice of the small parameters.}
\label{Sec:choice-parameters}
As in \cite{Poe1} we are going to iterate the KAM step 
infinitely many times choosing appropriately the parameters $0<s,r,\sigma, 
h,\eta<1$ and so on. Our goal is to get a convergent scheme. We are going
to define suitable strictly decreasing sequences of positive numbers 
$\{s_j\}_{j=0}^\infty$, 
$\{r_j\}_{j=0}^\infty$ and  $\{h_j\}_{j=0}^\infty$, tending to $0$.   
Set 
\begin{equation}
\label{eq:s}
s_j = s_0 \delta^j\, ,\ \sigma_j = \sigma_0 \delta^j\, ,\ 
s_0 = 5 \sigma_0(1-\delta)^{-1}\in (0,1) ,
\end{equation}
where $\delta=\delta(n,\tau,\vartheta) < 1/6$ is given by \eqref{eq:delta+eta}.

Given $m\ge 0$, we define an increasing sequence $ \nu(m):=(\nu_j(m))_{j\in\N}$ as follows. We set 
% and we define the sequences $ \nu(m):=(\nu_j(m))_{j\in\N}$ and $ (\ell_j)_{j\le J(m)}$  by
\begin{equation}
\nu_j(m) = \left\{
\begin{array}{lcrr}
\vartheta_0-\vartheta \quad \mbox{for}\quad j < J(m)  \\[0.3cm]
m(\tau +1) + \vartheta_0-\vartheta\quad  \mbox{for}\ j \ge J(m),
\end{array}
\right.
\label{eq:secuence-nu-m}
\end{equation} 
where $0<\vartheta <\min (\vartheta_0/4, 1)$, and 
\begin{equation}
\label{eq:J(m)}
J(m)\ge m(\tau+1)\vartheta^{-1}
\end{equation}
is an integer which will be determined in Sect. \ref{sec:nu+epsilon}.  If $m=0$,  we have  $\nu_j(0)= \vartheta_0-\vartheta$ for any $j\in\N$ and  we set $J(0)=0$.

%\begin{equation}
%\label{eq:nu}
%\left\{
%\begin{array}{lcrr}
%\vartheta = \vartheta_0/6=\nu_0\le \nu_1\le \cdots\,  ,\\[0.3cm]
%\displaystyle 2(\nu_j-\nu_0) \, \le \, \Sigma_{q=0}^{j-1}  \, \nu_q,  \quad j\ge 1, \\[0.3cm]
%\nu_j   \le C(n,\tau)(j +1),  \quad j\ge 1, 
%\end{array}
%\right.
%\end{equation}
%where $C(n,\tau)\ge 1$ may depend only on $n$ and $\tau$. 
%The first inequality will be used to obtain \eqref{eq:E-and-eta} for $j\ge 1$. 
Taking into account  \eqref{eq:r-and-s},  \eqref{eq:delta+eta} and  \eqref{eq:E+}, we define the sequences  $\{r_j(m)\}_{j\in\N}$, $\{\eta_j(m)\}_{j\in\N}$  and   $\{E_j(m)\}_{j\in\N}$ as follows.
Fix 
$$
r_0=s_0<1,\ \eta_0 =\delta^{\tau +1 +\vartheta + \nu_0} = \delta^{\tau +1 +\vartheta_0}   , 
$$
and set for $j\ge 1$ 
\begin{equation}
\label{eq:r}
\left\{
\begin{array}{lcrr}
\displaystyle \eta_j=\eta_j(m):= \delta^{\nu_j-\nu_{j-1} }\eta_{j-1} =\delta^{\nu_j-\nu_{0} }\eta_{0} =\delta^{\nu_j+\tau+1+\vartheta } , \\[0.3cm]
\displaystyle r_j = r_j (m):=  \eta_{j-1} r_{j-1}=\delta^{p_j}r_0,  \quad p_j= j(\tau+1+\vartheta )+ ( \nu_0+ \cdots +\nu_{j-1} ),
\end{array}
\right. 
\end{equation}   
and 
\begin{equation}
\label{eq:E}
E_j =E_j(m):= \delta^{\nu_{j-1}} E_{j-1}= \delta^{\nu_0 + \cdots +\nu_{j-1}}E_0. 
\end{equation}
%To indicate the dependence of $E_j$ of the sequence $\nu=(\nu_j)_{j\in\N}$, we write $E_j=E_j(\nu)$. 
Take the positive number $E_0=E_0(n,\tau,\vartheta_0)$ sufficiently small so that 
$$ 
E_0< \eta_0^2= \delta^{2\tau +2+2\vartheta_0} . 
$$
The inequality  \eqref{eq:J(m)} implies that 
\begin{equation}\label{eq:nu-j-inequality}
2\nu_j-2\nu_0\le \nu_0 + \cdots + \nu_{j-1}, \quad j\ge 1.
\end{equation}
Indeed, if $m=0$ then $\nu_j=\nu_0= \vartheta_0-\vartheta>0$ for each $j$. Let $m\ge 1$. 
For $j<J(m)$  have $\nu_j=\nu_0= \vartheta_0-\vartheta>0$ and for $j\ge J(m)$ we get 
$$
2\nu_j-2\nu_0 =2 m(\tau+1)\le J(m)\vartheta \le 2 j\vartheta < j(\vartheta_0 - \vartheta)\le \nu_0 + \cdots + \nu_{j-1}
$$ 
which yields  the inequality for any $m\in \N$ and $ j\in \N$. 
Now, \eqref{eq:r}, \eqref{eq:E} and \eqref{eq:nu-j-inequality} yield
\begin{equation}\label{eq:E-eta-j}
0<E_j<\eta_j^2\le \eta_0^2<1/64, \quad j\in\N.
\end{equation}
Taking into account   \eqref{eq:def-of-epsilon},  we put
\begin{equation}
\label{eq:epsilon}
\varepsilon_j:=  \hat\varepsilon    r_j\sigma_j^{\tau + 1}  E_j = 
  \hat\varepsilon    r_0\sigma_0^{\tau + 1}  E_0 \delta^{q_j}, 
\end{equation}
where $q_0=0$ and $q_j=q_j(m)$ is given for $j\ge 1$ by
\begin{equation}\label{eq:q-j}
\begin{array}{rcl}
q_j &:=&p_j+  j(\tau +1) + (\nu_0+\cdots+\nu_{j-1})\\
&=& j(2\tau +2+\vartheta) + 2(\nu_0+\cdots+\nu_{j-1}) .
\end{array}
\end{equation}
The parameter  $0<\hat\varepsilon \le 1 $   will be  chosen later. 
Finally, taking into account \eqref{eq:K} we set
\begin{equation}
\label{eq:h-and-K}
K_j= \sigma_j^{-1}\ln^2(\sigma_j)\quad \mbox{and}\quad  2 h_j =   K_j^{-\tau-1} =  \big(\sigma_j/\ln^2(\sigma_j)\big)^{\tau+1},\  j\in\N.
\end{equation}
We have
\begin{equation}
\label{eq:recurrence-relations}
\displaystyle s_{j+1} = s_j - 5\sigma_j,\ 
\sigma_j = 5^{-1}(1-\delta)s_j \quad \mbox{and}\quad h_{j+1}/h_j< \delta^{\tau+1}< 1/6
\end{equation}
in view of \eqref{eq:h+} and \eqref{eq:s}. Moreover, $ s_{j+1}-5s_{j+1}=\delta (s_{j}-5s_{j})$ and  \eqref{eq:constants} holds for each $j\in\N$. 
We are going to show that  \eqref{eq:sigma-eta} and hypothesis (a) - (d) of Proposition \ref{Prop:kam-step} are satisfied for any $j\in \N$. 
%More generally, we have the following 
\begin{Lemma}\label{Lemma:a b and c}  
	There exist   constants  
	\[
	0<\widetilde \sigma_0=\widetilde \sigma_0(n,\tau,\vartheta_0,\vartheta)<(1-\delta)/5, \quad 0< \widetilde E_0=\widetilde E_0(n,\tau,\vartheta_0,\vartheta)<1/64, 
	\]
	depending only on $n$, $\tau$, $\vartheta_0$ and $\vartheta$, such that \eqref{eq:constants}, \eqref{eq:sigma-eta}, \eqref {eq:E-eta-j},  and the hypothesis (a)-(d) are satisfied for any $j\in \N$,  provided that      
	\[
	0<  \sigma_0 \le \widetilde \sigma_0, \quad 0<E_0\le \widetilde E_0, \quad 2\ln^{2\tau +2}(\sigma_0)E_0\le c_0 .  
	\]
\end{Lemma}
\noindent
{\em Proof.}  We have already obtained \eqref{eq:constants} and \eqref {eq:E-eta-j} for $j\in\N$. Choosing $E_0 \le \widetilde E_0(n,\tau,\vartheta_0)\le c_0^2$ we get (a) for any $j\in\N$, while \eqref{eq:h-and-K} implies (c). Moreover, (d) holds since $K_j \sigma_j  =\ln^2(\sigma_j)>1$. 
On the other hand, (b)  holds if $E_j$ and $\sigma_j$ verify  
\eqref{eq:condition-b}. By \eqref{eq:secuence-nu-m} and since $\vartheta_0>4\vartheta$, we obtain
\[
2\ln^{2\tau +2}(\sigma_j)E_j = 2\ln^{2\tau +2}(\sigma_0\delta^j)\delta^{\nu_0+\cdots +\nu_{j-1}}E_0 < 2 \ln^{2\tau +2}(\sigma_0\delta^j)\delta^{2j\vartheta}E_0= f(\delta^j),
\]
where the function 
\[
x\mapsto f(x)=  2\ln^{2\tau +2}(\sigma_0x) x^{2\vartheta}E_0
\]
is increasing in the interval $(0,1]$, provided that $0<\sigma_0\le \widetilde \sigma_0^\prime :=\exp\left(-(\tau+1)\vartheta^{-1}\right)$. Then we have 
\[
2\ln^{2\tau +2}(\sigma_j)E_j < f(\delta^j) \le f(1)=2 \ln^{2\tau +2}(\sigma_0)E_0 \le c_0
\]
for  $0<\sigma_0 \le \widetilde \sigma_0^\prime$.

We are going to prove  \eqref{eq:sigma-eta}. For $j=0$ this means that
\[
 \sigma_0^{-n} \exp\left(-\ln^2(\sigma_0)\right)\delta^{-2\tau -2 -2\vartheta_0} \le 1.
\]
The function $x\mapsto  x^{-n} \exp\left(-\ln^2 x\right)$ is increasing in the interval $(0,e^{-\sqrt{n}}]$ and $\delta$ depends only on $n$, $\tau$ and $\vartheta$, hence, 
there exists a positive constant $ \widetilde \sigma_0'' =\widetilde \sigma_0''(n,\tau,\vartheta)\le e^{-\sqrt{n}}$ such that the inequality is satisfied for any $0<\sigma_0 \le \widetilde \sigma_0''$.

%Consider the function
%\[ f(\sigma):= \sigma^{-n} \exp(-K\sigma) = \sigma^{-n} \exp\left(-\ln^2(\sigma)\right). \]
%It is strictly increasing in an  interval $\big(0, \exp(-n/2)\big]$. 
 %We have 
 %\[  f'(\sigma) = \sigma^{-1} f(\sigma)\big(\vartheta (\tau+1)^{-1}\sigma^{-\vartheta (\tau+1)^{-1}}-n\big),\] 
 %We take  $\widetilde \sigma_0$ in that interval so  that 
% $f(\widetilde\sigma_0)\le \eta_0^2= \delta^{ \tau+2+2\vartheta_0-2\vartheta } $
 %which implies \eqref{eq:sigma-eta} for $j=0$. 
 Suppose now that $j\ge 1$. 
 Notice that 
 $$\nu_j(m)\le  j\vartheta +\vartheta_0 -\vartheta$$ 
  in view of \eqref{eq:secuence-nu-m} and \eqref{eq:J(m)}, hence, 
  \[
  \eta_j \ge \delta^{ j\vartheta +\tau+1+\vartheta_0 }= \delta^{ j\vartheta}\eta_0.
  \]
  This implies 
 \[  
  \sigma_j^{-n} \exp\left(-\ln^2(\sigma_j)\right) \eta_j^{-2} \le  \sigma_0^{-n}  \delta^{-jn} \exp\left(-\ln^2(\sigma_0\delta^j)\right)\delta^{-2j\vartheta}\eta_0^2 := g(\delta_j^{-n}).
 \]
 The function 
 $$x\mapsto g(x):=x^{-n-2\vartheta} \exp\left(-\ln^2(\sigma_0 x)\right)  \sigma_0^{-n}  \eta_0^2$$
 is increasing in the interval $0<x\le 1$ for $0<\sigma_0\le  \widetilde \sigma_0''':=e^{-\sqrt{n+2\vartheta}}$ and we get
  \[  
 \sigma_j^{-n} \exp\left(-\ln^2(\sigma_j)\right) \eta_j^{-2} \le g(\delta_j^{-n}) \le g(1)=  \sigma_0^{-n} \exp\left(-\ln^2(\sigma_0)\right) \eta_0^{-2} \le 1.
 \]
 for  $0<\sigma_0\le \widetilde \sigma_0'''(n,\tau,\vartheta)$. 
This yields \eqref{eq:sigma-eta} for any $j\in\N$.

\finishproof 

\noindent
We fix   $0<\sigma_0< (1-\delta)/5$ 
ones forever by 
\begin{equation}\label{eq:sigma-0}
0<\sigma_0 = \sigma_0(n,\tau,\vartheta_0, \vartheta):=  \min\big(\frac{1}{37}, \widetilde \sigma_0', \widetilde \sigma_0'', \widetilde \sigma_0'''\big),
\end{equation}
and then  choose $\widetilde E_0$ in Lemma \ref{Lemma:a b and c}  such that $2\widetilde E_0\le c_0\ln^{-2\tau -2}(\sigma_0)$. Then (b) holds for any $0<E_0\le \widetilde E_0$ and $j\in\N$. 
The choice of $\sigma_0 $ is motivated by the previous Lemma and by \eqref{eq:u-0}. 

Using  the proof of (b) in Lemma \ref{Lemma:a b and c}  we obtain the inequality
\begin{equation}
\frac{\varepsilon_j}{r_jh_j} \le 2 \ln^{2\tau +2}(\sigma_j)E_j      \le 3 \sigma_0^{-3\vartheta} \ln^{2\tau +2}(\sigma_j)\sigma_j^{3\vartheta}E_0 \le 
 C(n,\tau,\vartheta_0, \vartheta) \sigma_j^{2\vartheta}E_0 , 
\label{eq:epsilon-by-rh}
\end{equation}
since $\vartheta<\vartheta_0/4$. 
\begin{Remark}\label{rem:dependence-of-m}
The sequences of $\eta_j=\eta_j(m)$, $r_j=r_j(m)$ and weighted errors $E_j=E_j(m)$ depend on the choice of $m\in\N$, but $\sigma_j$, $h_j$ and $K_j$ do not depend on $m$. 
\end{Remark}

\subsubsection{\it Analytic smoothing of $P_t$.}
\label{Sec:analytic-smoothing}
The Hamiltonian $P_t$ is not analytic and one can not apply directly the KAM step to it. We are going to approximate it by real analytic functions. To this end 
we recall some facts about the analytic smoothing technique in Section \ref{Sec:ApprLemma}. 
We are going to apply the Approximation Lemmas \ref{Lemma:approximation-lemma}  to the real valued Hamiltonian  $P\in C^k\left([0,a];C_0^L(\A^n\times\Omega)\right)$, where $a>0$ and  $0<L\le \infty$. 
Set  
\begin{equation}
\label{eq:u-v-w-sequences}
u_j = u_0 \delta^j \, ,\ j\in\N,
\end{equation}
where 
\begin{equation}
\label{eq:u-0}
0<u_0= 6 s_0= 30 (1-\delta)^{-1} \sigma_0 \le 36 \sigma_0< 1,
\end{equation}
the small parameter $0<\delta=\delta(n,\tau,\vartheta)<1/6$ is given by \eqref{eq:delta+eta} and $\sigma_0$ is fixed in \eqref{eq:sigma-0}. 
%Note that
%\[
%\sup\{s_j/u_j,  r_j/v_j,h_j/w_j\} \le \frac{1}{6}   ???????
%\]
%for $j\in\N$. 
Let us denote by  ${\cal U}_j$ the complex strips  in  
${\C}^n/2\pi {\Z}^n \times  {\C}^n \times  {\C}^n$ consisting of all
$(\theta, I; \omega)$ such that 
\begin{equation} 
|{\rm Im\, } \theta| \, ,\  
|{\rm Im\, } I| \, ,\ 
|{\rm Im\, } \omega| \ < \  u_j\,  ,
\label{eq:u-v-w}
\end{equation}
and by 
$A({\cal U}_j)$ the set of all real-analytic  bounded functions in ${\cal
	U}_j$ equipped with the sup-norm $|\cdot|_{u_j}$. 
Define
\begin{equation}\label{eq:approximation-of-P}
P_t^j:= S_{u_j}P_t\, ,\quad j\in\N,
\end{equation}
by means of the Approximation Lemma \ref{Lemma:approximation-lemma}. 
This is  a $C^k$ family with respect to $t\in[0,a]$ of real analytic  in ${\C}^n/2\pi {\Z}^n \times  {\C}^n \times  {\C}^n$  functions.    In view of \eqref{eq:approximation2-1},  for each finite $\ell\le L$ and $0\le \ell'\le \ell$, the following inequality is true
\begin{equation}\label{eq:P-j}
\|P_t^j  -  P  \|_{\ell'}\,  \le \, C(n,\ell)\,  u_j^{\ell-\ell'}\|P\|_\ell
\end{equation}
in the corresponding H\"{o}lder norms on $\T^n\times\R^n\times\Omega$. 
On the other hand, the inequality \eqref{eq:approximation2}  with $\rho=u_j$ and  $\tilde\rho=u_{j-1}=\delta^{-1} u_j$, yields  the estimate 
\[
\left|\partial_t^p(P_t^{j} -  P_t^{j-1} )  \right|_{u_j} \le C_0 u_{j-1}^{\ell}\|\partial_t^p P_t\|_{\ell} =C u_j^{\ell}\|\partial_t^p P_t\|_{\ell}
\]
for each finite $\ell$, $0\le \ell\le L$ and $0\le p\le k$, where $C=C(\ell,n,\tau, \vartheta_0)=C_0(n,\ell)\delta^{-\ell}$ is a positive constant depending only on $\ell,n,\tau, \vartheta_0$. Moreover, 
\[
\left|\partial_t^p P_t^{0}  \right|_{u_0} \le C_0 \|\partial_t^p P_t\|_{0} \le \widetilde C u_0^{\ell}\|\partial_t^p P_t\|_{\ell} 
\]
where $\widetilde C = C_0(n)u_0^{-\ell}$. The positive constant $C_\ell := \max(c, \widetilde c)=C_\ell(n,\tau,  \vartheta_0)$ depends only on $\ell, n,\tau$ and $\vartheta_0$. 
Hence, 
\begin{equation}
\label{eq:approximation-estimates1}
\left|\partial_t^p P_t^{0}    \right|_{u_0} \le \widetilde\varepsilon_{\ell,0,p}\quad \mbox{and} \quad \left|\partial_t^p P_t^{j} -  \partial_t^p P_t^{j-1}   \right|_{u_j} \le \widetilde\varepsilon_{\ell,j,p} \quad \mbox{for} \quad j\ge 1
\end{equation}
for any finite $\ell$, $0\le \ell\le L\le \infty$, where 
\begin{equation}
\label{eq:approximation-estimates2}
\widetilde\varepsilon_{\ell,j,k}\, :=\,  C_\ell u_j^{\ell}\, \sum_{p=0}^k \, \sup_{0\le t \le a} \|\partial_t^p P_t\|_{\ell} 
\end{equation}
and $C_\ell=C(\ell,n,\tau,\vartheta_0)>0$ depends only on $\ell$, $n$, $\tau$ and $\vartheta_0$. 

We would like to deal with $P^{j}$ at the $j$-th iteration putting $P^{j}-P^{j-1}$ in the error term. 
To this end   we need for $0\le p\le k$ the following inequalities
\begin{equation}
\label{eq:approximation6-0}
\left|\partial_t^p P_t^{0}    \right|_{u_0} \le \frac{\varepsilon_{1}}{2} \quad \mbox{and} \quad \left|\partial_t^p P_t^{j} -  \partial_t^p P_t^{j-1}   \right|_{u_j} \le \frac{\varepsilon_{j+1}}{4} \quad \mbox{for} \quad j\ge 1. 
\end{equation}
These inequalities will be obtained in Sect. \ref{sec:nu+epsilon}, choosing appropriatelly the sequence $\nu$ and the small constants $\epsilon$ and $\hat{\varepsilon}$. 

Using the notations introduced in the beginning of Sect. \ref{Sec:kamlemma} we set 
\begin{equation}\label{eq:complex-V}
D_j := D_{s_j,r_j}\, ,\ O_j:=O_{h_j} \, ,\ {\cal V}_j := D_j\times O_j .
\end{equation}
Moreover, given an integer $1\le q \le 3$ we set
\begin{equation}\label{eq:D-j-q}
D_{j}^q:= D_{\frac{q}{4} s_j, \frac{q}{4}r_j}, \   O_{j}^q :=O_{\frac{q}{4}h_j}  \ \mbox{and} \ {\cal V}_j^q := D_j^q\times O_j ^q. 
\end{equation}
We have 
\[
D_{j+1}\times O_{j+1} \subset D_j^1\times O_j ^1 
\]
since \( {\rm sup\, }\{s_{j+1}/s_{j},\, r_{j+1}/r_{j}, \, h_{j+1}/ h_{j}\}\le \delta <1/6. \) 

\subsubsection{\it Iterative Lemma.}
\label{Sec:iteration-lemma}
We are ready to make the iterations. Consider the real analytic in 
${\cal U}_j$ Hamiltonian  
\[
H_t^j(\varphi,I;\omega)= H^j(\varphi,I;\omega,t) := N_0(I;\omega) + P_t^j(\varphi,I;\omega),
\] 
where 
$ N_0(I;\omega) := \langle \omega,I\rangle$ and ${\cal U}_j$ is defined by \eqref{eq:u-v-w}.   
Let us  denote by ${\cal U}_j^0$ the subset of 
${\C}^n/2\pi{\Z}^n \times  {\C}^n \times  {\C}^n$ consisting of all
$(\theta, I; \omega)$ such that 
\[
|{\rm Im\, } \theta| \, ,\  
|{\rm Im\, } I| \, ,\ 
|{\rm Im\, } \omega| \  <\   \frac{1}{2} u_j\, .
\]
We have $2s_j < u_j$, which yields 
$ D_j\times O_j\subset {\cal U}_j^0$. 
Using the notations introduced in \eqref{eq:D-j-q} we  obtain
\begin{equation}\label{eq:inclusion}
D_{j+1}\times O_{j+1} \subset D_{j}^2\times O_{j}^2 \subset  D_j\times O_j \subset {\cal U}_j^0\subset {\cal U}_j   
\end{equation}
since ${\rm sup\, }\{s_{j+1}/s_{j},\, r_{j+1}/r_{j}, \, h_{j+1}/ h_{j}\}\le \delta <1/6$. 
For any $j\in {\N}$,
let us  denote by ${\cal D}_j$ the class of  real-analytic 
diffeomorphisms 
\[
{\cal F}_{j}: D_{j+1}\times  O_{j+1} \to D_{j}^2\times  O_{j}^2
\]
of the form  
\begin{equation}
\label{eq:transformation-F}
{\cal F}_j(\theta,I;\omega)\ =\  (\Phi_j(\theta,I;\omega),\phi_j(\omega))\, ,\ 
\Phi_j(\theta,I;\omega)\ =\ (U_j(\theta;\omega),V_j(\theta,I;\omega))\, ,                            
\end{equation}
where $V_j(\theta,I;\omega)$ is affine  linear  with respect to $I$,  
and  
$(\theta, I) \to \Phi_j(\theta,I;\, \omega)$ 
is a canonical transformation   for  any fixed $\omega$. To simplify the notations we
denote the sup-norm of functions $f:D_{j}\times  O_{j}\to \C$ by $|f |_j=|f|_{s_j,r_j,h_j}$. 
%Obviously,  $D_{j}\times  O_{j} \subset {\cal U}_j$. 
Fix $k\in \{0;1\}$.

\begin{Prop}[Iterative Lemma] \label{Prop:IterativeLemma} 
Let  $P^j\in C^k([0,a], {\mathcal A}({\mathcal U}_j))$, $j\in \N$,  be a $C^k$ family of  real analytic Hamiltonians in ${\mathcal U}_j$ satisfying \eqref{eq:approximation6-0} and $H^j=N_0+P^j$.  
Then for each $j\in \N$ there is a normal form 
$N_j(I;\omega,t) 
= e_j(\omega,t) + \langle \omega,I\rangle$ and  a $C^k$ family of  real analytic 
transformations 
\begin{equation}
{\cal F}^j\in  C^k([0,a], {\mathcal A}( D_{j}\times O_{j}, (D_0\times O_0) \cap {\cal
U}_{j}^0)), \quad {\cal F}^j_t = {\cal F}^j(\cdot,t),
\label{eq:range-F-0}
\end{equation}
such that
\begin{enumerate}
	\item 
${\cal F}^0 = {\it
id}$ and  ${\cal F}^{j+1}_t = {\cal F}_{t,0}\circ\cdots\circ {\cal F}_{t,j}$, for $ j\ge 0$, where 
\begin{equation}
{\cal F}_j\in C^k([0,a], {\mathcal A}( D_{j+1}\times O_{j+1}, D_{j}^2\times O_{j}^2))\quad \mbox{and}\quad {\cal F}_{t,j}(\cdot,t):={\cal F}_{j}(\cdot,t)\in {\cal D}_j ;  
\label{eq:range-F}
\end{equation}
	\item 
$H^j\circ {\cal F}^{j+1} = N_{j+1} + R_{j+1}$ and  
$\displaystyle | \partial_t^p R_{j+1} |_{j+1}  \le \varepsilon_{j+1}/2$  for $0\le p\le k$;
	\item The following estimates hold
\begin{equation} 
\displaystyle  |\overline{W}_j \partial_t^p({\cal F}_{t,j} - {\rm id} )|_{j+1} \ +\ 
|\overline{W}_j \partial_t^p(D {\cal F}_{t,j} - {\rm Id} )
\overline{W}_j^{\, -1}|_{j+1} \ <   
 \frac{C_0 \varepsilon_j}{r_jh_j} , 
                         \label{eq:estimates-F1} 
\end{equation}
\begin{equation} 
\displaystyle  |\partial_t^p({\cal F}^{j+1}_t - {\cal F}^j_t)|_{j+1} \, 
< \frac{C_0 \varepsilon_j}{r_jh_j} ,
                          \label{eq:estimates-F2} 
\end{equation} 
for $0\le p\le k$ and uniformly with respect to $t\in [0,a]$,  where $C_0 =C_0(n,\tau,\vartheta_0,\vartheta)>0$,   $D {\cal F}^{j}_t$
stands for the Jacobian of ${\cal F}^{j}_t$ with respect to
$(\theta,I;\omega)$, and \\
$\overline{W}_j = {\rm diag}\, \left({\sigma_j}^{-1}{\rm Id}, 
{r_j}^{-1}{\rm Id},{h_j}^{-1}{\rm Id}\right)$.   
\end{enumerate}
\end{Prop} 
{\em Proof}. For $k=0$ the proof is similar to that of the Iterative Lemma in 
\cite{Poe1} and it is done in \cite{P4} in the case of Gevrey Hamiltonians independent of $t$. Additional efforts are required for the proof of the estimates \eqref{eq:estimates-F1} and \eqref{eq:estimates-F2} in the case when $p=k=1$.

Consider firstly the Hamiltonian $H^0= N_0+P^0$. It satisfies the hypothesis of Proposition \ref{Prop:kam-step} in $D_0\times O_0$ for $t\in [0,a]$ in view of \eqref{eq:approximation6-0} and Lemma \ref{Lemma:a b and c}. 
Hence, applying the KAM Step Lemma to the Hamiltonian $H^0$  we find 
${\cal F}^{1} = {\cal F}_{0}$ such that 
$H^0 \circ {\cal F}^1 =  N_1 + R_{1}$, where $R_{1}(\cdot,t)$ is real analytic
in $D_1\times O_1$ and $| \partial_t^p R_{1}(\cdot,t)|_{1} \leq \varepsilon_1/2$. Moreover, \eqref{eq:estimates-F1} holds for $j=0$.

Given $j\ge 1$ we suppose that the Proposition holds for all indexes $0\le l\le j-1$. We are going to prove it for $l=j$. We are looking for a transformation 
${\cal F}^{j+1} = {\cal F}^{j}\circ {\cal
	F}_{j}$, where ${\cal F}_j$  belongs to ${\cal D}_j$. 
By the inductive
assumption we have 
\[
H^{j-1} \circ {\cal F}^{j} = N_j + R_{j},
\] 
where 
$N_j(I;\omega,t) = e_j(\omega,t) + \langle \omega, I\rangle$, $R_{j}(\cdot,t)$ is real analytic   in $D_j\times
O_j$, and  $ |  \partial_t^p R_{j} (\cdot,t)|_{j}  \leq  \varepsilon_{j}/2$. 
Then we write  
$$
H^j \circ {\cal F}^{j+1} = (N_0 + P^{j-1})\circ {\cal F}^{j+1} + 
(P^{j} - P^{j-1})\circ {\cal F}^{j+1} 
$$
$$
= \left(H^{j-1} \circ {\cal F}^{j}\right)\circ {\cal F}_{j} 
+ (P^{j} - P^{j-1})\circ {\cal F}^{j+1} 
$$
$$
= \left(N_j + R_{j}+ (P^{j} - P^{j-1})\circ {\cal
F}^{j}\right) \circ {\cal F}_{j} .     
$$
Consider the Hamiltonian $\widetilde H_j= N_j+ R_j + (P^{j} - P^{j-1})\circ {\cal
	F}^{j}$   in $D_j\times O_j$ for $t\in [0,a]$ and set $R_j^1= (P^{j} - P^{j-1})\circ {\cal
	F}^{j}$. Using \eqref{eq:approximation6-0}  we get
\[
|(P_t^{j} - P_t^{j-1})\circ {\cal F}^{j}(\cdot,t)|_{j} 
\leq |P_t^{j} - P_t^{j-1}|_{{\cal U}_{j}^0} 
\leq \frac{\varepsilon_{j}}{4} .
\]
On the other hand, by the inductive assumptions  \eqref{eq:range-F-0} we obtain 
\[
|\partial_t((P_t^{j} - P_t^{j-1})\circ {\cal F}^{j}(\cdot,t))|_{j} \le 
 |\partial_t P_t^{j} - \partial_t  P_t^{j-1}|_{{\cal U}_{j}^0} + |(D(P_{j} - P_t^{j-1})\circ {\cal F}^{j} (\cdot,t))\, .\partial_t {\cal F}^{j}(\cdot,t)|_{j} .
\]
The firs term of the right hand side is estimated by $\frac{\varepsilon_{j}}{4} $ in veiw of \eqref{eq:approximation6-0}. Using \eqref{eq:range-F-0} we estimate the second one by
\[
|D(P_t^{j} - P_t^{j-1})W_j^{-1}|_{{\cal U}^0_{j}}\, |W_j\partial_t {\cal F}^{j}(\cdot,t)|_{j} 
\]
(here we consider $\overline W_{j}$ as a linear operator acting on  $\C^{3n}$). 
Now Cauchy estimates (see Remark \ref{Rem:Cauchy} ) and \eqref{eq:approximation6-0} yield 
\begin{equation}\label{eq:Cauchy-P}
\begin{array}{rcll}
\displaystyle \left|D(P_t^{j} - P_t^{j-1})W_j^{-1}\right|_{{\cal U}^0_{j}}  = \left|\left(s_j\nabla_{\theta}, r_j\nabla_{I},h_j\nabla_{\omega}\right)(P_t^{j} - P_t^{j-1})\right|_{{\cal U}^0_{j}} \\ [0.4cm]
\le \displaystyle \sup\{s_j,r_j,h_j \} \, \frac{2}{u_j} \, |P_t^{j} - P_t^{j-1}|_{{\cal U}_{j}} 
\le \displaystyle 2\times \frac{1}{6} \times\frac{\varepsilon_{j}}{4} <  \frac{\varepsilon_{j}}{4} .  
\end{array}
\end{equation}
Moreover, \eqref{eq:estimates-F1} and \eqref{eq:epsilon-by-rh} imply that $|W_j\partial_t {\cal F}^{j}(\cdot,t)|_{j} \le 1$ for $E_0=E_0(n,\tau,\vartheta_0,\vartheta)>0$ small enough.  Finally, we obtain 
$$
|\partial_t^p R_j + \partial_t^p R_j^0|_{j} < \varepsilon_{j}\, , \quad p\in \{0;1\}. 
$$
We apply the KAM Step  Lemma - Proposition \ref{Prop:kam-step} - to the $C^k $ family of
Hamiltonians $\widetilde H_j$. 
 Using Remark \ref{Rem:Cauchy}, \eqref{eq:new-error2} and  \eqref{eq:epsilon-by-rh} as well,  
we find a $C^k$ family of real-analytic maps 
\[
{\cal F}_{j}(\cdot,t): 
D_{j+1} \times O_ {j+1} \to  D_{j}^2 \times O_{j}^2
\]
which belong to the class ${\cal D}_{j}$, satisfy \eqref{eq:estimates-F1} and such that 
$(N_j + R_{j})\circ {\cal F}_{j} = N_{j+1} + R_{j+1}$, where 
$$ 
|\partial_t^p R_{j+1}|_{j+1} \leq 
 \frac{1}{2} 
 \hat\varepsilon_j    r_{j+1}\sigma_{j+1}^{\tau +1}
E_{j+1}   \le \frac{1}{2} 
 \hat\varepsilon_{j+1}    r_{j+1}\sigma_{j+1}^{\tau +1}
E_{j+1} = \frac{\varepsilon_{j+1}}{2} .  
$$ 
%Moreover,  ${\cal F}_{j}$ satisfies
%(\ref{eq:estimates-F1}) in view of Remark \ref{Rem:estimates} and Lemma \ref{Lemma:a b and c}. ????? 
We are going to show that 
\begin{equation}
{\cal F}^{j+1} : D_{j+1} \times O_{j+1}\  \longrightarrow \ {\cal
U}_{j}^0 .
                                  \label{eq:map-F}
\end{equation}
%This implies
%$$
%|\partial_t^p(P_{j} - P_{j-1})\circ {\cal F}^{j+1}|_{j+1} 
%\leq |\partial_t^p(P_{j} - P_{j-1})|_{{\cal U}_{j}} 
%\leq \frac{\varepsilon_{j+1}}{2} ,   
%$$
%and  we obtain  
%$H^j\circ {\cal F}^{j+1} = N_{j+1} + R_{j+1}$, where
%$|\partial_t^p R_{j+1}|_{j+1} \le \varepsilon_{j+1}$ for $0\le p\le k \le 1$. 
To prove (\ref{eq:map-F}) we estimate the norm of the linear operator  $W_q \overline W_{q+1}^{\, -1}$. We have
\[
|\overline W_q \overline W_{q+1}^{\, -1}| = 
\sup \left\{s_{q+1}/s_q,r_{q+1}/r_q ,h_{q+1}/h_q
\right\} =  s_{q+1}/s_q = \delta, 
\]
since $r_{q+1}/r_q \le \delta ^{\tau +1} < \delta$ and $h_{q+1}/h_q \le 
\delta ^{\tau +1} < \delta$ for any $q\in {\N}$ by \eqref{eq:r} and \eqref{eq:recurrence-relations}. 
Recall that  $\delta$ and $E_0$ depend only on $n$, $\tau$, $\vartheta_0$ and $\vartheta$. 
Then using \eqref{eq:epsilon-by-rh} and  the inductive assumption (\ref{eq:estimates-F1}), 
we  estimate the Jacobian of ${\cal
F}^{j+1}$ in $D_{j+1}\times O_{j+1}$ as follows (see also \cite{Poe1})
\[
\begin{array}{rcl}
\left|\overline W_0 D{\cal F}^{j+1}\overline W_{j}^{\, -1}\right|_{j+1}
\,  &=& \, 
\left|\overline W_0 D ({\cal F}_{0}\circ \cdots \circ 
{\cal F}_{j})\overline W_{j}^{\, -1}\right|_{j+1}       \\[0.3cm]
&\le& \displaystyle \prod_{q=0}^{j-1}
\left(\left|\overline W_q D {\cal F}_{q}\overline W_q ^{\, -1}\right|_{q+1}\, 
\left|\overline W_q \overline W_{q+1}^{\, -1}\right|\right)
\left|\overline W_j D {\cal F}_{j}\overline W_{j}^{\, -1}\right|_{j+1}\, 
 \\[0.3cm]
&\le& \displaystyle\delta^{j}  \prod_{k=0}^\infty\, 
\left(1 +  \frac{C \varepsilon_k}{r_kh_k} \right) \le \displaystyle \delta^{j} 
\exp\left(\sum_{k=0}^\infty 
 \frac{C \varepsilon_k}{r_kh_k}\right)\\[0.3cm] 
 &<& \displaystyle  \delta^{j}\exp\left(C(1-\delta^{2\vartheta}) ^{-1}E_0\right) , 
\end{array}
\]
where $C=C(n,\tau,\vartheta_0,\vartheta)$ stands for different positive constants depending only on $n$, $\tau$, $\vartheta_0$ and $\vartheta$. 
Choosing the parameter $E_0 =E_0(n,\tau,\vartheta_0,\vartheta)>0$ sufficiently small 
we obtain  
\begin{equation}
 \left|\overline W_0 D{\cal F}^{j+1}\overline W_{j}^{\, -1}\right|_{j+1}  <  \delta^{j}, \quad j\in \N.
 \label{eq:estimate-differential-F}
\end{equation}
Set 
\[
z=(\theta,I,\omega) =x +iy\in D_{j+1}\times  O_{j+1},
\]
where  $x$
and $y$ are respectively the real and the imaginary part of $z$. Then $|\overline W_{j+1} \,  y| \le 1$, where $|\cdot|$ stands for the sup-norm. We have 
\[
\begin{array}{lcrr}
\displaystyle {\cal F}^{j+1}(x+iy) \, =\,  {\cal F}^{j+1}(x) +  i
\overline W_0^{\, -1} T_{j+1}(x,y)\overline W_{j}\,  y \, ,\\ [0.3cm]
\displaystyle  T_{j+1}(x,y) = \displaystyle {  \int_0^1 \,  \overline W_0 
D{\cal F}^{j+1}(x + ity) \overline W_{j}^{\, -1} \, dt \,  } 
\end{array} 
\]
(we consider $\overline W_{j}$ as a linear operator acting in $(\R^{3n}, |\cdot|)$). 
Moreover,   
$|T_{j+1}(x,y)| <  \delta^{j}$ and since 
$|\overline W_{j} \,  y| \le \delta|\overline W_{j+1} \,  y| \le \delta\le  1/6$,  we
get 
\[
|T_{j+1}(x,y)\overline W_{j}\,  y| <  \frac{1}{2}  \, \delta^{j}\,  
,\quad x+iy\in D_{j+1}\times  O_{j+1} \, . 
\]
Denote by $Z_{j+1}(x,y)$ the imaginary  part of ${\cal F}^{j+1}(x+iy) $.  
Since ${\cal F}^{j+1}(x)$ is
real valued, $Z_{j+1}(x,y)$ is equal to  the real  part of $\overline W_0^{\, -1} T_{j+1}(x,y)\overline W_{j}\,  y $. Then we get
\[
u_j^{-1}  | Z_{j+1}(x,y)| \le \delta^{-j}  u_0^{-1}  | \overline W_0^{\, -1} ||T_{j+1}(x,y)\overline W_{j}\,  y|  <  \frac{1}{2}  \,  
,\quad x+iy\in D_{j+1}\times  O_{j+1} \, ,
\]
and we obtain \eqref{eq:map-F}.

It remains to prove \eqref{eq:estimates-F2}. In the case when $p=0$ it follows from the arguments in \cite{Poe1}. Suppose now that $p=k=1$.  
Denote by $D{\cal F}^{j}(z)$  the differential of ${\cal F}^{j}$ with respect to $z=(\theta,I,\omega)$ acting on vectors $\eta\in \C^{3n}$ by 
$\eta\to D{\cal F}^{j-1}(z)\, .\, \eta$. 
%  and by  $(\xi,\eta)\to D^2 {\cal F}^{j}(z)[\xi,\eta]$  the symmetric bilinear form corresponding to the second differential of $ {\cal F}^{j}$ at $z$. 
Consider
\[
\begin{array}{lcrr}
\displaystyle \partial_t({\cal F}^{j+1} - {\cal F}^{j}) = 
\partial_t( {\cal F}^{j}\circ{\cal F}_{j} - {\cal F}^{j})  \, \\ [0.3cm]
\displaystyle = 
(D{\cal F}^{j}\circ{\cal F}_{j} )\, . \partial_t {\cal F}_{j}  + (\partial_t{\cal F}^{j})\circ{\cal F}_{j} -\partial_t {\cal F}^{j} = \Sigma_1 + \Sigma_2
\end{array} 
\]
where 
\[
\begin{array}{lcrr}
\displaystyle \Sigma_1:=   (D{\cal F}^{j}\circ{\cal F}_{j} )\, \cdot \partial_t {\cal F}_{j}  ,\\ [0.3cm]
\displaystyle \Sigma_2:=  \int_0^1 \, (D \partial_t {\cal F}^{j} )(x{\cal F}_{j} +(1-x) {\,  id\, }) \, . \, ({\cal F}_{j} - {\it id\, })\ dx . 
\end{array} 
\]
We are going to estimate $\Sigma_l$, $1\le l\le 2$. By  \eqref{eq:range-F}, \eqref{eq:estimates-F1} and \eqref{eq:estimate-differential-F} we get for any $j\ge 1$ 
\begin{equation}\label{eq:estimate-Sigma1}
 \left| \Sigma_1 \right|_{j+1} < h_0^{-1}\left|\overline W_0 D{\cal F}^{j}\overline W_{j}^{\, -1}\right|_{j} \left|\overline W_{j}\partial_t\, {\cal F}_{j} \right|_{j+1}  < C_0   \frac{\varepsilon_j}{r_jh_j} .
\end{equation}

Consider $  \Sigma_2$ now. Set ${\cal F}^{0}={\cal F}_{-1}= id$ and put ${\cal F}^{q,j}= {\cal F}_{q}\circ \cdots \circ  {\cal F}_{j-1}$ for $q\le j-1$ and ${\cal F}^{j,j}= id$. 
For $j\ge 1$  we have 
\[
\partial_t {\cal F}^{j} = 
\partial_t \left( {\cal F}_{0}\circ \cdots \circ  {\cal F}_{j-1}\right) = \sum_{q=0}^{j-1} (D{\cal F}^{q}\circ{\cal F}^{q,j})\, .\, 
((\partial_t  {\cal F}_{q})\circ{\cal F}^{q+1,j}) \, .
\]
Using \eqref{eq:range-F}, \eqref{eq:estimates-F1},  \eqref{eq:estimate-differential-F} and \eqref{eq:epsilon-by-rh}  we get as above 
\[
\begin{array}{lcrr}
 \displaystyle  \left| \partial_t {\cal F}^{j} \right|_{j} < h_0^{-1}  \sum_{q=0}^{j-1} \left|\overline W_0  D{\cal F}^{q}W_{q}^{\, -1}\right|_{q} \, 
\left|\overline W_q \partial_t  {\cal F}_{q}\right|_{q+1} \\ [0.3cm]
\displaystyle  \le C \sum_{q=0}^{j} \frac{\varepsilon_q}{r_qh_q } < C(1-\delta^{2\vartheta}) ^{-1}E_0
\end{array} 
\]
where $C$ stands for different constants depending only on $n$, $\tau$, $\vartheta_0$ and $\vartheta$. By Cauchy this implies 
\begin{equation}\label{eq:estimate-D-t-F}
 \left| D\partial_t {\cal F}^{j} \overline W_{j} \right| \le C
 \end{equation}
uniformly on $ D_j^2 \times O_{j}^2\times [0,a]$, 
and we get 
\begin{equation}\label{eq:estimate-Sigma2}
 \left|  \Sigma_2 \right|_{j+1}\,  \le\,  \sup_{ D_j ^2\times O_{j}^2\times [0,a]}\, \left|\overline W_0 D \partial_t{\cal F}^{j}\overline W_{j}^{\, -1}\right|  \left|\overline W_{j}( {\cal F}_{j} - id) \right|_{j+1}  \le C   \frac{\varepsilon_j}{r_jh_j} 
\end{equation}
where $C=C(n,\tau,\vartheta_0,\vartheta)$ stands for different positive constants depending only on $n$, $\tau$, $\vartheta_0$  and $\vartheta$.

This proves \eqref{eq:estimates-F2} for $p=k=1$. In the case when $p=0$ we use the same arguments. 
This completes the proof of Proposition \ref{Prop:IterativeLemma}. 
\finishproof

\begin{Remark}[Uniqueness by construction in the Iterative Lemma]\label{rem:uniqueness-Iteration}
	The transformations ${\mathcal F}_{t,j}$, the normal forms $N_{t,j}:= N_{j}(\cdot,t)$ and the error terms $P_{t,j}:=P_{j}(\cdot,t)$   do not depend on the choice of $m\ge 0$ in \eqref{eq:secuence-nu-m} in the following sense. 
	Let $m'\ge 0$  and let ${\mathcal F}_{t,j}'$,  $N'_{t,j}$ and $P'_{t,j}$,  be the corresponding transformations,  normal forms and  error terms. 
	Then ${\mathcal F}'_{t,j}={\mathcal F}_{t,j}$,  $N'_{t,j}= N_{t,j}$  and $P'_{t,j}= P_{t,j}$ on the intersection of their domains of definition. 
\end{Remark}
Remark \ref{rem:uniqueness-Iteration} follows from Remark \ref{rem:uniqueness-KAM} and Remark \ref{rem:dependence-of-m} by induction whit respect to $j\in \N$. 

The Iterative Lemma provides a convergent schema  giving in a limit a $C^\infty$ function on  $\T^n\times \Omega_1$ in a Whitney sense. To avoid inconveniences arising  from the Whitney extension theorem, we propose a modified Iterative Lemma in the next section. 

%-----------------------------------------------------------------------------------

\subsubsection{\it Modified Iterative Lemma.}
\label{Sec:modified-iteration-lemma}

We are going to modify ${\mathcal F}^j_{t}$ multiplying ${\mathcal F}^j_{t}- id$ by a suitable almost analytic cut-off function in $\omega\in \C^n$. 

\vspace{.3cm}
\noindent
{\bf 1.}    {\em  Construction of almost analytic cut-off functions}. 

\vspace{.3cm}
\noindent
We say that a function $f:\C^n\to \C^n$,  given by $x+iy\mapsto f(x+iy):=f(x,y)$ for $x,y\in\R^n$, is $\R$-smooth, or $C^\infty$ in a real sens, if the function 
$\R^n\times \R^n \ni (x,y)\mapsto f(x,y)$ 
is $C^\infty$-smooth.  
As usually we denote by  $\bar{\partial}_l$,  $1\le l \le n$, the operators
\[
 {\bar{\partial}}_l =\frac{\partial}{\partial \bar{z}_l} = \frac{1}{2}\left(\frac{\partial}{\partial x_l} +i \frac{\partial}{\partial y_l}  \right), 
\]
and we set  $\bar\partial= (\bar\partial_1, \ldots, \bar\partial_n)$. A $\R$-smooth function $f:\C^n\to \C^n$ is called almost-analytic if the vector-function 
$(x,y)\mapsto \bar{\partial}f(x+iy)$ is flat at $\R\times \{ 0\}$, in the sense that 
\[
\partial_y^\beta\bar{\partial} f(x,y)|_{y=0} = 0 \quad \mbox{for any}\ \beta\in \N^n. 
\]
Such a function is ``very small'' for $y$ small. If $f$ is an almost analytic  Gevrey functions, then it is even exponentially small.  
Given $\rho>1$ and $L\ge 1$, we say that $f$ belongs to the Gevrey class $\mathcal{G}^\rho_L(\C^n)$ if it is $\R$-smooth and 
\[
\|f\|_{L} := \sup_{\alpha,\beta\in{\N}^n}\, 
\sup_{(x,y)\in{\R}^n \times \R^n  } \, 
\left(|\partial_x^\alpha \partial_y^\beta  f(x,y)|\,  L^{-|\alpha|-|\beta|}
\alpha !^{-\rho} \beta !^{-\rho}\right)\ <\ \infty \, ,
\]
where 
$|\alpha| = \alpha_1 + \cdots + \alpha_n$ and $\alpha ! =
\alpha_1! \cdots \alpha_n!$ for $\alpha=(\alpha_1,\ldots,\alpha_n)
\in {\N}^n$.  We say that $f$ is Gevrey-${\mathcal G}^\rho$ function.  If the function $f\in \mathcal{G}^\rho_L(\C^n)$ is almost-analytic, then there exist positive constants $C=C(n,\rho)$ and $c=c(n,\rho)$ depending only on $n$ and $\rho$, such that
\[
|\partial_x^\alpha\partial_{y}^{\beta} \bar{\partial_l}f(x+iy)|\ \leq\ 
C \|f\|_L\,  L^{|\alpha|+|\beta|} 
\alpha !\,^\rho  \beta !\, ^{\rho}\,  
\exp\left(- c
(L|y|)^{\, -\frac{1}{\rho-1}}\right)  
\]
for any $\alpha, \beta\in \N^n$ and $1\le l\le n$.

\begin{Prop}\label{prop: almost-analytic} 
For any $n\ge 2$ and $\rho>1$ there exist positive constants $C= C(n,\rho)$, $L=L(n,\rho)$ and $c=c(n,\rho)$, and a family of almost-analytic  functions $\chi_j\in \mathcal{G}^{\rho}_{L/h_{j+1}}(\C^n)$, $j\in\N$, with the following properties
\begin{enumerate}
	\item [(i)] $\mbox{\rm supp\, }(\chi_j) \subset O_{j+1}^3$ and $\chi_j=1$ on ${O}_{j+1}^2$;
	\item[(ii)]  $\|\chi_j\|_{L/h_{j+1}}\le C$ for $j\in \N$; 
	\item[(iii)]  the following estimate holds
	\[
	|\partial_x^\alpha\partial_{y}^{\beta} \bar{\partial}\chi_j(x+iy)|\ \leq\ 
	C  (L/h_{j+1})^{|\alpha|+|\beta|+1} 
	(\alpha !\,  \beta !)\, ^{\rho}\,  \exp\left(-  c
	(|y|/h_{j+1})^{\, -\frac{1}{\rho-1}}\right)  
	\]
	on $C^n$ for any $j\in\N$ and $\alpha, \beta\in \N^n$. 
\end{enumerate}

\end{Prop}

The proposition will be proved in Section \ref{Sec:Gevrey}.

\vspace{.3cm}
\noindent
{\bf 2.}    {\em  Modified transformations}. 

\vspace{.3cm}
\noindent
From now on we take $\rho=2$ in Proposition \ref{prop: almost-analytic}. 
We define the modified transformations 
\[
 {\mathcal H}_{t,j} : D_{j+1}\times \C^n \to {\C}^n/ 2\pi {\Z}^n\times \C^n\times \C^n
 \]
  by 
  \begin{equation}
 {\mathcal H}_{t,j}(z) := z +\chi_j(\omega)\left({\mathcal F}_{t,j} (z)-z\right), \quad z=(\theta,I;\omega)\in D_{j+1}\times \C^n .
  \label{eq:modified-transformations}
  \end{equation}
  Setting
  ${\mathcal H}_{t,j} =(\widetilde \Phi_{t,j}, \widetilde \phi_{ t,j})$, this means that
\[
\left\{
\begin{array}{lcrr}
 \widetilde \Phi_{t,j}(\theta,I;\omega) &=& (\theta,I) +\chi_j(\omega)\left(\Phi_{t,j} (\theta,I;\omega)-(\theta,I)\right),  \\
  \widetilde \phi_{j,t}(\omega) &=& \omega +\chi_j(\omega)\left(\phi_{t,j} (\omega)-\omega\right), 
  \end{array}
  \right.
\]
for  $(\theta,I;\omega)\in D_{j+1}\times \C^n$. 
\begin{Lemma}\label{lemma:rnage-of-H} The following relations hold  for any $j\in \N$ provided  $E_0=E_0(n,\tau,\vartheta_0,\vartheta)$ is sufficiently small 
	\begin{enumerate}
		\item[(1)] ${\mathcal H}_{t,j} : D_{j+1}\times \C^n \to D_{j}^2\times \C^n$, 
		\item[(2)] $\widetilde \phi_{t,j}(O_{j+1})\, \subset\, O_{j}^2$ and ${\mathcal H}_{t,j} : D_{j+1}\times O_{j+1} \to D_{j}^2\times O_{j}^2$. 
	\end{enumerate}

\end{Lemma}

\noindent
{\em Proof}.
\quad {\em (1)}  \
Recall from \eqref{eq:estimates-F1}  that
\[
|{W}_j (\Phi_{t,j} - {\rm id} )|_{D_{j+1}}\ <   
\frac{C_0 \varepsilon_j}{r_jh_j} <C_0 C(n,\tau,\vartheta_0,\vartheta) E_0
\]
by \eqref{eq:estimates-F1} and \eqref{eq:epsilon-by-rh}, where ${W}_j = {\rm diag}\, \left({\sigma_j}^{-1}{\rm Id}, 
{r_j}^{-1}{\rm Id} \right)$.   Moreover,
\[
  |\chi_j|_{\C^n}<\|\chi_j\|_{1/h_{j+1}}\le C(n)
\]  
   in view of Proposition \ref{prop: almost-analytic}, (ii). This yields
\[
|{W}_j (\widetilde \Phi_{t,j} - {\rm id} )|_{D_{j+1}}\ <C_1(n,\tau,\vartheta_0,\vartheta) E_0 \le \frac{1}{8}
\]
choosing $E_0=E_0(n,\tau,\vartheta_0,\vartheta)$ sufficiently small, and we obtain {\em (1).}
%\begin{equation}\label{eq:image-of-H}
%{\mathcal H}_{t,j} : D_{j+1}\times \C^n \to D_{j}^2\times \C^n. 
 %\end{equation}
 
 {\em (2)}  \  Let $\omega\in    O_{j+1}$. Then there exists $\omega'\in \Omega_1$ such that    $|\omega-\omega'|   \le h_{j+1}$ and        we get  as above 
 \[
 \begin{array}{rcll}
 \displaystyle |\widetilde \phi_{t,j}(\omega)-\omega'| &\le&  \displaystyle |\omega-\omega'| + \| \chi_j\|_{1/h_{j+1}} | \phi_{t,j}(\omega)-\omega'| \le h_{j+1} + C(n)C_0\frac{\varepsilon_j}{r_j}  \\ [0.3cm]
 &<&  \displaystyle  h_{j+1} + C(n,\tau,\vartheta_0,\vartheta)E_0 h_j < h_{j+1} + \frac{1}{3}h_j < \frac{1}{2}h_j,
 \end{array}
 \]
 for $E_0=E_0(n,\tau, \vartheta_0,\vartheta)>0$ sufficiently small, hence, $\widetilde \phi_{t,j}(\omega)\in  O_{j}^2$. \finishproof

Let  us define ${\mathcal H}_{j}=  (\widetilde \Phi_{j}, \widetilde \phi_{ j})$ and ${\mathcal H}^{j}$ by $ {\mathcal H}_{j}(\cdot,t) = {\mathcal H}_{t,j}$ and $ {\mathcal H}^{j}(\cdot,t)={\mathcal H}^{j}_t$ for $t\in [0,a] $, where 
\[
   {\mathcal H}^0_t  =   {\it 	id}, \  {\mathcal H}^{j+1}_t :=  {\mathcal H}_{t,0} \circ \cdots \circ {\mathcal H}_{t,j} : D_{j+1}\times \C^n \to D_{0}\times \C^n.
\]
%Set
%\[\widetilde O_j  := \{ \omega\in \C^n :\, n|{\rm Im\,}(\omega_\alpha)| < h_j \ln^{-2}(\sigma_j),\, \alpha=1,\ldots,n \} . \]
We set $\bar{\partial}_l  = \frac{\partial }{\partial \bar{\omega}_l}$ for $1\le l\le n$ and $\bar{\partial}=(\bar{\partial}_1, \ldots, \bar{\partial}_n)$ . 
We are going to use as well  the convention $\frac{1}{+0}=+\infty$ and $\exp(-\infty)=0$. 

\begin{Prop}[Modified Iterative Lemma] \label{Prop:ModifiedIterativeLemma} 
Under the assumptions of Proposition 	\ref{Prop:IterativeLemma}, the transformations $ {\mathcal H}^{j}$ are well defined on $D_j\times \C^n\times [0,a]$ and 
have the following properties
	\begin{enumerate}
		\item [(i)]  $ {\mathcal H}^{j} \in C^k \left([0,a], C^\infty( D_{j}\times \C^n, D_0\times \C^n) \right)$ and $ {\mathcal H}_t^{j}=  {\mathcal F}_t^{j}$ on $D_{j}\times O_j^2$ for $t\in [0,a]$. Moreover, 
		\[
		 {\rm supp\, }({\mathcal H}_t^{j+1}-{\mathcal H}_t^{j})\subset D_{j+1}\times  O_{j+1}^3	\quad  \mbox{and}\quad {\mathcal H}_t^{j+1}-{\mathcal H}_t^{j} = {\mathcal F}_t^{j} \circ {\mathcal H}_{t,j} - {\mathcal F}_t^{j} ; 
		\]
		\item [(ii)] \quad $ \displaystyle  |\partial_t^p({\cal H}^{j+1}_t - {\cal H}^j_t)(z)| \, 
		< \frac{C_0 \varepsilon_j}{r_jh_j} $ for $z=(\theta,I;\omega)\in D_{j+1}\times \C^n$,  $t\in [0,a]$, and $0\le p\le k$ ,  where $C_0 =C_0(n,\tau,\vartheta_0,\vartheta)>0$; 
		\item [(iii)]  ${\mathcal H}^{j+1}-{\mathcal H}^{j}$ is analytic with respect to $(\theta,I)\in D_j$ and  almost analytic and Gevrey-${\mathcal G} ^2$ with respect to $\omega$. Moreover,  for any $0\le p\le k$   the following estimate holds
		\begin{equation} 
		\begin{array}{lcrr}
		\displaystyle  \left|\bar{\partial} \partial_t^p \left({\mathcal H}^{j+1}-{\mathcal H}^{j}\right)(\theta,I;\omega,t)\right| \, \le \, 
		C  h_{j+1}^{-1}\,  \exp\left(- c\, \frac{h_{j+1}}{|{\rm Im\,}(\omega)|}\right)   \, \frac{\varepsilon_j}{r_jh_j},\\ [0.3cm]
		\mbox{for}   \ (\theta,I;\omega)\in D_{j+1}^3\times \C^n, \ t\in [0,a], 
		\end{array}
		\label{eq:estimates-H1} 
		\end{equation}
		where $C=C(n,\tau,\vartheta_0, \vartheta)$ and  $c=c(n)$ are positive constants;
		\item [(iv)]   the following estimate is true
		\begin{equation} 
		\displaystyle  \left|\partial_t^p \partial_\theta^\beta\partial_\omega^\gamma \left({\mathcal H}^{j+1}-{\mathcal H}^{j}\right)(\theta,I;\omega,t)\right| \, \le \, 
		C_{\beta,\gamma} \frac{\varepsilon_j}{r_jh_j}\sigma_{j+1}^{-|\beta|}h_{j+1}^{-|\gamma|}  \ln^{2|\gamma|+2}(\sigma_{j+1})
		\label{eq:estimates-H} 
		\end{equation}
		for $ \ (\theta,I;\omega)\in D_{j+1}^2\times \R^n$, $t\in [0,a]$,  and $0\le p\le k$,  $\beta,\gamma\in\N^n$, where \\
		$C_{\beta,\gamma}=C_{\beta,\gamma}(n,\tau,\vartheta_0, \vartheta)>0$.
	\end{enumerate}
\end{Prop} 
{\em Proof}. \quad {\em (i)} \quad 
%Item {\em 2} in Lemma \ref{lemma:rnage-of-H} implies for any $0\le l\le j-1$ the relation 
%\[   \widetilde \phi_{t,l}\circ \cdots \circ \widetilde \phi_{t,j}: O_{j+1}\to O_{l}^2.   \]
%On the other hand $\chi_l=1$ on $O_{l+1}^2$ and we obtain
%\begin{equation}\label{eq:chi-phi}
%\chi_{j-1}\circ\widetilde \phi_{t,j}(\omega)= \chi_{j-2}\circ\widetilde \phi_{t,j-1}\circ\widetilde \phi_{t,j}(\omega)= \cdots =\chi_0\circ\widetilde \phi_{t,1}\circ\cdots\circ\widetilde \phi_{t,j}(\omega) =1 \quad \forall\, %\omega \in O_{l+1}.
%\end{equation}
Recall that $\chi_{j-1}= 1$ on $O_j^2$, hence, ${\cal H}_{t, j-1}= {\cal F}_{t, j-1}$ on $D_{j}\times O_j^2$. Moreover, $ {\cal F}_{t, j-1}:D_{j}\times O_j \to D_{j-1}^2\times O_{j-1}^2$ by Proposition \ref{Prop:IterativeLemma} and $\chi_{j-2}= 1$ on $O_{j-1}^2$, which implies 
\begin{center}
\( ({\mathcal H}_t^{j-1}\circ {\cal H}_{t, j-1})(z)=( {\mathcal H}_t^{j-1}\circ {\cal F}_{t, j-1})(z)= ({\mathcal F}_t^{j-1}\circ {\cal F}_{t, j-1})(z)\) 
\end{center}
for any $z\in D_{j}\times O_j^2$. 
Repeating this argument we obtain  the equality
\[
{\mathcal H}_t^{j}(z)=( {\mathcal H}_t^{j-1}\circ {\cal F}_{t, j-1})(z)= \cdots = ( {\mathcal F}_{t,0}\circ \cdots\circ {\cal F}_{t, j-1})(z) = {\mathcal F}_t^{j}(z).
\]

If $\omega\notin O_{j+1}^3$, then $\chi_j(\omega)=0$, ${\cal H}_{t, j}(\theta,I;\omega)=(\theta,I;\omega)$ and $ {\mathcal H}_t^{j+1}(\theta,I;\omega)= {\mathcal H}_t^{j}(\theta,I;\omega)$, hence, 
${\rm supp\, }({\mathcal H}_t^{j+1}-{\mathcal H}_t^{j})\subset D_{j+1}\times  O_{j+1}^3$. 

Let $z=(\theta,I;\omega)\in D_{j+1}\times O_{j+1}^3 \subset  D_{j+1}\times O_{j+1}$. 
Lemma \ref{lemma:rnage-of-H} implies that ${\cal H}_{t, j}(z)\in D_{j}^2\times O_{j}^2$, hence, $\chi_{j-1}(\widetilde \phi_j(\omega))=1$ and we obtain 
 \[({\cal H}_{t, j-1}\circ{\cal H}_{t, j})(z)= ({\cal F}_{t, j-1}\circ {\cal H}_{t, j})(z) \in D_{j-1}^2\times O_{j-1}^2. \] 
On the other hand $\chi_{j-2}=1$ on $O_{j-1}^2$ and repeating this argument we get
\[
\begin{array}{rcll}
 ({\mathcal H}_t^{j} \circ {\cal H}_{t, j})(z)&=& ( {\cal H}_{t, 0} \circ\cdots \circ {\cal H}_{t, j-1} \circ {\cal H}_{t, j})(z)= ( {\cal H}_{t, 0} \circ\cdots \circ {\cal F}_{t, j-1} \circ {\cal H}_{t, j})(z)\\
 &= &
 \cdots =( {\cal F}_{t, 0} \circ\cdots \circ {\cal F}_{t, j-1} \circ {\cal H}_{t, j})(z) = ({\mathcal F}_t^{j} \circ {\cal H}_{t, j})(z) . 
 \end{array}
\]
Moreover, $ {\mathcal H}_t^{j}(z)= {\mathcal F}_t^{j}(z)$ since $z\in  D_{j+1}\times O_{j+1}\subset   D_{j}^2\times O_{j}^2$, and we obtain the equality
\[
{\mathcal H}_t^{j+1}-{\mathcal H}_t^{j} = {\mathcal F}_t^{j} \circ {\mathcal H}_{t,j} - {\mathcal F}_t^{j}
\]
on $D_{j+1}\times O_{j+1}^3 $. On the other hand, both sides of it vanish at any $z\in D_{j+1}\times (\C^n\setminus O_{j+1}^3 )$. Hence, the above equality is true on $D_{j+1}\times \C^n $.
This completes the proof of (i).

%by induction with respect to $j$ that $ {\mathcal H}_t^{j+1}-{\mathcal H}_t^{j} = {\mathcal F}_t^{j} \circ {\mathcal H}_{t,j} $. 

\vspace{0.3cm}
{\em (ii)}  \quad Let $p=0$. Using (i) we obtain
\[
{\mathcal H}^{j+1}-{\mathcal H}^{j} = {\mathcal F}^{j} \circ {\mathcal H}_{j} - {\mathcal F}^{j} =  \chi_j \int_0^1 \, (D  {\cal F}^{j} )(x{\cal H}_{j} +(1-x) {\,  id\, }) \, . \, ({\cal F}_{j} - {\it id\, })\, dx .
\]
Moreover, for any $z\in  {\rm supp\, }({\mathcal H}_t^{j+1}-{\mathcal H}_t^{j})\subset D_{j+1}\times  O_{j+1}\subset  D_j^2\times O_j^2  $ we have 
\begin{equation}\label{eq:inclusion-H}
x{\cal H}_{j}(z) +(1-x)z\in D_j^2\times O_j^2 \quad   \mbox{for}\quad  0\le x \le 1,
\end{equation}
and \eqref{eq:estimates-F1} and \eqref{eq:estimates-F2} imply
\[
|{\mathcal H}^{j+1}-{\mathcal H}^{j}|_{j+1} \le  C(n ) h_0^{-1}\left|\overline W_0 D{\cal F}^{j}\overline W_{j}^{\, -1}\right|_{j} \left|\overline W_{j}( {\cal F}_{j}- {\,  id\, } )\right|_{j+1}  < C_0   \frac{\varepsilon_j}{r_jh_j} .
\]

Let $p=1$. By the chain rule we get
\[
\begin{array}{lcrr}
\partial_t ( {\mathcal H}^{j+1}-{\mathcal H}^{j}) =\partial_t ( {\mathcal F}^{j} \circ {\mathcal H}_{j} - {\mathcal F}^{j}) \\ [0.3cm]
\displaystyle = 
(D{\cal F}^{j}\circ{\cal H}_{j} )\, . \partial_t {\cal F}_{j}  + (\partial_t{\cal F}^{j})\circ{\cal H}_{j} -\partial_t {\cal F}^{j} = \Sigma_1 + \Sigma_2
\end{array} 
\]
where 
\[
\begin{array}{lcrr}
\displaystyle \Sigma_1:=   (D{\cal F}^{j}\circ{\cal H}_{j} )\, \cdot \partial_t {\cal F}_{j}  ,\\ [0.3cm]
\displaystyle \Sigma_2:=  \int_0^1 \, (D \partial_t {\cal F}^{j} )(x{\cal H}_{j} +(1-x) {\,  id\, }) \, . \, ({\cal F}_{j} - {\it id\, })\ dx . 
\end{array} 
\]
Using \eqref{eq:inclusion-H},  we estimate $\Sigma_1$ and $\Sigma_2$ as in \eqref{eq:estimate-Sigma1} and  \eqref{eq:estimate-Sigma2}. 

\vspace{0.3cm}
{\em (iii)} \quad  
Let $p=0$. 
 %If $|\alpha| =1$, then \eqref{eq:estimates-H}  follows from  (ii) and \eqref{eq:estimates-F2}, taking into account  the inequality  \eqref{eq:estimate-h-sigma}.
%We are going to prove \eqref{eq:estimates-H}  for any $\alpha\in \N^n$ with $|\alpha|\ge 1$. To this end we will first generalize  (ii) as follows. 
Recall that 
 ${\mathcal F}^{j}$ is analytic on  $ D_{j}\times  O_{j}$, $ {\cal F}_{j}$ is analytic with respect to $z=(\theta,I;\omega)\in D_{j+1}\times  O_{j+1}$, and that ${\mathcal H}_{j}(z,t)= z + \chi_j(\omega)({\cal F}_{j}(z,t)-z)$. 

Differentiating the identity in (i) we obtain for any 
$1\le l\le n$ the following one
\begin{equation}\label{eq:bar-partial-H}
\begin{array}{crll}
\displaystyle  \bar{\partial}_l  ({\mathcal H}^{j+1}-{\mathcal H}^{j})(z,t)=\frac{\partial}{\partial\bar{\omega}_{l}  } ( {\mathcal F}^{j} \circ {\mathcal H}_{j} - {\mathcal F}^{j})(z,t)\\ [0.3cm]
 \displaystyle = \,   \bar\partial_l\chi_j(\omega )D{\mathcal F}_t^{j}({\mathcal H}_{j}(z,t))\, .\left({\mathcal F}_{j}(z ,t)-z\right)
\end{array}
\end{equation}
for   each  $z\in  {\rm supp\, }({\mathcal H}_t^{j+1}-{\mathcal H}_t^{j})$ and $t\in [0,a]$. According to Proposition \ref{prop: almost-analytic} we have
 \begin{equation}\label{eq:estimate-R}
 \sup_{\omega\in\C^n} |\bar\partial_l\chi_j(\omega )| \le C h_{j+1}^{-1} \,  \exp\left(- c\, \frac{h_{j+1}}{|{\rm Im\,}(\omega)|}\right)   .
 \end{equation}
 Morreover,  $z\in  D_{j+1}\times  O_{j+1}^3	$ and 
${\mathcal H}_{j}(z,t)\in D_{j}^2\times  O_{j}^2$ 
in view of item (i), and arguing as in the proof of \eqref{eq:estimate-Sigma1} we obtain (iii) for $p=0$.

\vspace{0.3cm}
Let $p=k=1$. 
Differentiating the identity \eqref{eq:bar-partial-H} with respect to $t$  we obtain  
\[
\begin{array}{crll}
  \partial_t \bar\partial_l ({\mathcal H}^{j+1}-{\mathcal H}^{j})(z,t) \\ [0.3cm] 
 = \displaystyle  \bar\partial_l\chi_j(\omega ) \partial_tD{\mathcal F}^{j}_t({\mathcal H}_{j}(z,t))\, . \left({\mathcal F}_{j}(z ,t)-z\right)\\ [0.3cm]
 + \displaystyle  \bar\partial_l\chi_j(\omega ) D{\mathcal F}^{j}_t({\mathcal H}_{j}(z,t))\, . \partial_t{\mathcal F}_{j}(z ,t)\\ [0.3cm]
 +  \displaystyle \chi_j(\omega) \bar\partial_l\chi_j(\omega ) D^{2}{\mathcal F}^{j}_t({\mathcal H}_{j}(z,t))\left[\partial_t{\mathcal F}_{j}(z ,t),   {\mathcal F}_{j}(z,t) -z)\right]
\end{array}
\]
for    $z\in  {\rm supp\, }({\mathcal H}_t^{j+1}-{\mathcal H}_t^{j})$. 

Consider  for any $z=(\theta,I;\omega)\in D_j\times O_j$ the symmetric bilinear form  
\[
(\xi,\eta)\to D^2 {\cal F}_t^{j}(z,t)[\xi,\eta], \quad \xi,\eta\in \C^{3n}, 
\] 
representing  the second differential  of $ {\cal F}_t^{j}$ at $z$.  We have
\begin{Lemma}\label{lemma:estimate-bilinear}
	There  exists $C=C(n,\tau,\vartheta_0,\vartheta) >0$ such that 
	\begin{equation}\label{eq:estimate-D2-F}
	\left|  D^2  \partial_t^p {\cal F}_t^{j}(z)[\overline{ W}_{j}^{\, -1}\xi, \overline{ W}_{j}^{\, -1}\eta]\right| \,   \le \,  C |\xi|_{\C^{3n}}|\eta|_{\C^{3n}}  , \quad \xi,\eta\in \C^{3n}, 
	\end{equation}
	for any  $z\in  D_{j}^2\times O_{j}^2$,  $t\in [0,a]$,  $0\le p\le k$, and $j\in \N$. 
\end{Lemma}
{\em Proof}. \quad The expression $D^2  \partial_t^p {\cal F}^{j}(z,t)[\overline{ W}_{j}^{\, -1}\xi,\overline{ W}_{j}^{\, -1}\eta]$ is a sum of monomials of the form
\[
(\sigma_j \partial_\theta)^\alpha (r_j \partial_I)^\beta (h_j\partial_\omega) ^\gamma  \partial_t^p {\cal F}^{j}(\theta,I;\omega,t)\,  \xi_l\eta_m, 
\]
where $|\alpha|+ |\beta| + |\gamma|=2$ and $1\le l,m\le n$. 
The estimate follows from \eqref{eq:estimates-F2} and the Cauchy inequalities . \finishproof

For each $z\in  {\rm supp\, }({\mathcal H}_t^{j+1}-{\mathcal H}_t^{j})\subset  D_{j+1}\times  O_{j+1}$ we have ${\mathcal H}_t^{j}(z)\in  D_{j}^2\times  O_{j}^2$  in view of Lemma \ref{lemma:rnage-of-H}. 
Then the  estimate \eqref{eq:estimates-H}  follows from \eqref{eq:estimate-R}, 
Lemma \ref{lemma:estimate-bilinear},    and  \eqref{eq:estimates-F1} as in the case $p=0$.
This proves {\em (iii)}. 

%Recall that $\omega\in O_{j+1}^3$ on the support of $R_q$. 

%Differentiating \eqref{eq:simplified-chain} we obtain
%\[
%\begin{array}{rcll}
%\displaystyle \frac{\partial}{\partial \bar{\omega}_\alpha}\partial_t\left({\mathcal H}^{j+1}-{\mathcal H}^{j}\right)  &=& 
% \displaystyle\frac{\partial \chi_j}{\partial \bar{\omega}_\alpha}\, \cdot \left(D{\mathcal F}^{j}\circ  {\mathcal H}_{j}\right)\, \cdot \partial_t({\mathcal F}_{j}-{\rm id}) \\ [0.3cm]
%& +& \displaystyle \frac{\partial \chi_j}{\partial \bar{\omega}_\alpha}\, \cdot \left(D\left(\partial_t{\mathcal F}^{j}\right)\circ  {\mathcal H}_{j}\right)\, \cdot ({\mathcal F}_{j}-{\rm id}) \\ [0.3cm]
% &+& \displaystyle \frac{\partial \chi_j}{\partial \bar{\omega}_\alpha} \chi_j\, \cdot \left(D^2{\mathcal F}^{j}\circ  {\mathcal H}_{j}\right)\left[\partial_t({\mathcal F}_{j}-{\rm id}), ({\mathcal F}_{j}-{\rm id})\right].
%\end{array}
%\]

\vspace{0.3cm}
{\em (iv)} \quad We are going to use Cauchy formula for almost analytic functions. Let $f:\C\to\C$ be a $\R$-smooth almost analytic function. Denote by $D_r(x)$ the open disc $\{z\in\C:\,  |z-x|<r\}$ of radius $r>0$ and by $\partial D_r(x)$ its boundary oriented counter clockwise. For any $x\in\R$ and $\beta\in \N$, the following Cauchy integral formula is valid
\begin{equation}\label{eq:Cauchy-almost-analytic}
\displaystyle \frac{d^\beta f}{d x^\beta}(x)\, =\, \frac{\beta !}{2\pi i}\, \left\{   \int_{\partial D_r(x)}^{} (z-x)^{-\beta-1}f(z)  \,  dz \, +\,   \iint_{ D_r(x)}^{} (z-x)^{-\beta-1}\bar{\partial}f(z)\,  dz\wedge d\bar{z}  \right\}
\end{equation}
(see e.g. \cite{Hed}, Proposition 1.1). Notice that the second integral is well defined since the function $\bar{\partial} f$ is flat at $\R$.

Set 
\[
R_j = h_{j+1}\ln^{-2}(\sigma_{j+1})
\]
and 
\[
I_{\beta,\gamma_1} (\theta,I;\omega,t):=\partial_t^p \partial_\theta^\beta\partial_{\omega_1}^{\gamma_1}\ \left({\mathcal H}^{j+1}-{\mathcal H}^{j}\right)(\theta,I;\omega,t)
\]
for  $(\theta,I;\omega)\in D_{j+1}^2\times \R^n, \ t\in [0,a]$. Recall that ${\mathcal H}^{j}$ is analytic with respect to $(\theta,I)\in  D_{j+1}$. Applying first the Cauchy inequality with respect to $\theta$,  we get
\[
|I_{\beta,\gamma_1} (\theta,I;\omega,t)| \le  (2\sigma_{j+1})^{-|\beta|-1} |I_{0,\gamma_1} (\theta,I;\omega,t)| .
\]
for $(\theta,I;\omega)\in D_{j+1}^2\times \R^n, \ t\in [0,a]$. We are going to estimate $ |I_{0,\gamma} (\theta,I;\omega,t)|$.

%Then $ B_{R_j}^q(\omega) \subset \widetilde O_{j+1}$. 
Applying Cauchy formula \eqref{eq:Cauchy-almost-analytic} to the variable $x=\omega_1\in \R$, keeping $\omega'=(\omega_2,\ldots,\omega_n)$ fixed in $\R^{n-1}$,  we obtain
\begin{equation}\label{eq:integrals}
\begin{array}{rcll}
I_{0,\gamma_1} (\theta,I;\omega,t) \, &=&\, \displaystyle \frac{\gamma_1 !}{2\pi i}\,   \int_{D_{R_j}(\omega_1)}  \,  \frac{\partial_t^p \bar{\partial}_1({\mathcal H}^{j+1}-{\mathcal H}^{j})(\theta,I;z,\omega',t)}
{(z-\omega_1)^{\gamma_1+1}}\, dz \, \\ [0.5cm]
&+&  \displaystyle\frac{\gamma_1!}{2\pi i}\, \int_{D_{R_j}(\omega_1)}\, \frac{\partial_t^p \bar{\partial}_1({\mathcal H}^{j+1}-{\mathcal H}^{j})(\theta,I;z,\omega',t)}
{(z-\omega_1)^{\gamma_1+1}}\, dz\wedge d\bar{z}
\end{array}
\end{equation}
for  $(\theta,I;\omega)\in D_{j+1}^2\times \R^n$ and $t\in [0,a]$.

Using (ii) we estimate the first integral by 
\[
C_{\beta,\gamma_1}  \frac{\varepsilon_j}{r_jh_j} h_{j+1}^{-\gamma_1-1} \ln^{2\gamma_1+2}(\sigma_{j+1})
\]
for $(\theta,I;\omega)\in D_{j+1}^2\times \R^n$ and $t\in [0,a]$. 

In order to estimate the second integral we are going to use the  following 
 estimate 
\[
\begin{array}{crll}
\displaystyle  \frac{1}{ |z-\omega_{1}|^{\gamma_{1}+1}}\,\exp\left(- c\, \frac{h_{j+1}}{|{\rm Im\,}(z)|}\right)  h_{j+1}^{-1}\\ [0.3cm]
\displaystyle \le 
\frac{1}{  |{\rm Im\,}(z)|^{\gamma_{1}+1}}\,\exp\left(- \frac{c}{2}\, \frac{h_{j+1}}{ |{\rm Im\,}(z)| }\right) 2^{\tau+1} \sigma_{j+1}^{-\tau-1}\ln^{2\tau+2}(1/\sigma_{j+1}) e^{- \frac{c}{2}\ln^{2}(\sigma_{j+1})}\\ [0.3cm]
\displaystyle \le 
C_{\gamma_1}   \frac{1}{h_j^{\gamma_{1}+1} }
\end{array}
\]
for $z\neq \omega_1$. 
Using    (iii) and the estimate above, we estimate  the second integral by
\[
  C_{\gamma}  \frac{\varepsilon_j}{r_jh_j} h_{j+1}^{-\gamma_1-1}
\]
for $(\theta,I;\omega)\in D_{j+1}^2\times \R^n$ and $t\in [0,a]$. 	

This implies
\[
|I_{0,\gamma_1} (\theta,I;\omega,t)| \le 
C_{\gamma_1}  \frac{\varepsilon_j}{r_jh_j} h_{j+1}^{-\gamma_1-1}\ln^{2\gamma_1+2}(1/\sigma_{j+1})
\]
for $(\theta,I;\omega)\in D_{j+1}^2\times \R^n$ and $t\in [0,a]$. 	

Finally we obtain 
\[
|I_{\beta,\gamma_1} (\theta,I;\omega,t)| \le 
C_{\beta, \gamma_1}  \frac{\varepsilon_j}{r_jh_j} h_{j+1}^{-\gamma_1-1} \ln^{2\gamma_1+2}(1/\sigma_{j+1})\sigma_{j+1}^{-|\beta|}
\]
for $(\theta,I;\omega)\in D_{j+1}^2\times \R^n$ and $t\in [0,a]$. 	This proves (iv) in the case when $\gamma =(\gamma_1, 0,\ldots,0)$. By a permutation of the indexes, we obtain it as well for 
$\gamma =(0, \ldots,0,\gamma_l, 0,\ldots,0)$. It remains to prove the estimate for the mixed derivatives with respect to $\omega$. To this end we shall use the following 
\begin{Lemma}\label{lemma:directional-derivatives}
For any $\gamma\in \N^n$ of length $N=|\gamma|\ge 2$ there exist  $(N+1)^{n-1}$ vectors $\vec{v}_m$ and constants $c_m$ such that
\[
\partial_\omega^\gamma = \sum c_m {\mathcal L}_{\vec{v}_m}^N
\]
where ${\mathcal L}_{\vec{v}_m}$ stands for the directional derivative ${\mathcal L}_{\vec{v}_m}f(\omega)= \frac{d}{ds}|_{s=0}f(\omega +s \vec{v}_m)$. 
\end{Lemma}	
{\em Proof}. \quad We proceed by induction with respect to $n\ge 2$. Let $n=2$ and $\gamma=(\gamma_1,\gamma_2)$ with $\gamma_2\neq 0$. Denote by $\vec{e}_1$ and $\vec{e}_2$ the canonical basis of $\R^2$ and set $\vec{v}= \vec{e}_1 +\lambda \vec{e}_2$, where $\lambda>0$. We  have
\[
{\mathcal L}_{\vec{v}}^N =\sum_{l=0}^N \lambda^{l} L_{l}, \quad  L_l:= \frac{N!}{l! (N-l)!} \partial_1^{N-l}\partial_2^{l} .
\]
Choosing  $\lambda_m=m/(N+1)$ and $\vec{v}_m= \vec{e}_1 +\lambda_m\vec{e}_2$ for $1\le m \le N+1$, we obtain the linear system
\[
\sum_{l=0}^N \lambda_m^{l} L_{l} = {\mathcal L}_{\vec{v}_m}^N, \quad m=0,\ldots,N.
\]
This system has a unique solution with respect to $L_l$, $0\le l \le N$,  since the corresponding determinant is just the Vandermonde determinant. Then we use induction with respect to $n$. \finishproof

Applying the preceding argument for each derivative  ${\mathcal L}_{\vec{v}_m}^N$, we 
complete  the proof of (iv). \finishproof

\begin{Remark}[Uniqueness  in the Modified Iterative Lemma]\label{rem:uniqueness-Modified}
	The transformations ${\mathcal H}^{j}$  do not depend on the choice of  $m\ge 0$ in \eqref{eq:secuence-nu-m} in the sense of Remark \ref{rem:uniqueness-Iteration}. 
\end{Remark}
%-----------------------------------------------------------------------------------

\subsubsection{\it Choice of  the  sequence $\nu $ and the small parameters $\epsilon$ and $\hat{\varepsilon}$.}\label{sec:nu+epsilon}

Given  $m\in\N$ we consider the sequence $(\nu_j(m))_{j\in \N}$ introduced in \eqref{eq:secuence-nu-m} and set 
\begin{equation}
\label{eq:secuence-ell-m}
\ell_j(m)= \left\{
\begin{array}{lcrr}
\ell_0 = 2\tau +2 + 2\vartheta_0 \ &\mbox{for}&\  j<J(m), \\ [0.3cm]
\ell_m= 2m(\tau +1) +\ell_0 \ &\mbox{for}&\  j\ge J(m), 
\end{array}
\right.
\end{equation}
where $J(0)=0$ and $J(m)$, $m>0$,  will be  a suitable integer satisfying \eqref{eq:J(m)}.

Consider  the family of functions $P_t^j$, $j\in\N$,  defined by \eqref{eq:approximation-of-P}. In order to apply the Iterative Lemma  and the Modified Iterative Lemma to that family, we have to show  that it satisfies \eqref{eq:approximation6-0}, To this end we will choose appropriately  the small constants $\epsilon$ and $\hat{\varepsilon}$ as well as  the integer $J(m)$ for each $m>0$.  
By  \eqref{eq:approximation-estimates1} and \eqref{eq:approximation-estimates2}, it suffices to prove for each $j\in \N$ that 
\begin{equation}
\label{eq:approximation6-1}
C_{\ell} u_j^{\ell}\, \sum_{p=0}^k \, \sup_{0\le t \le a} \|\partial_t^p P_t\|_{\ell} \, \le \, \frac{\varepsilon_{j+1}}{4} \, 
=\,  \frac{1}{4}\hat\varepsilon    r_{j+1}(m)\sigma_{j+1}^{\tau + 1}  E_{j+1} (m)  \quad \mbox{with}\    \ell=\ell_j(m). 
\end{equation}
Here $r_j(m)$ and $E_j(m)$  are given by \eqref{eq:r} and \eqref{eq:E}, respectively, $\hat\varepsilon\in (0,1]$.
%, which will be determined below. To satisfy \eqref{eq:approximation6-1}, we shall choose $\hat\varepsilon\in (0,1]$ appropriately. 
In view of  \eqref{eq:u-v-w-sequences} and \eqref{eq:r}-\eqref{eq:q-j}, the relation  \eqref{eq:approximation6-1} becomes
\begin{equation}
\label{eq:approximation6-2}
\sum_{p=0}^k \, \sup_{0\le t \le a}\|\partial_t^p P_t\|_{\ell}\,  \le \,   \hat\varepsilon    r_0\sigma_0^{\tau + 1}   \frac{\left(\delta/6 s_0\right)^{\ell}}{4C_{\ell}} \delta^{M_j(m)}   \quad \mbox{with}\    \ell=\ell_j(m), 
\end{equation}
where 
\begin{equation}
\label{eq:M(nu)}
\begin{array}{rcll}
M_j(m)\, & :=& \, q_{j+1}- (j+1)\ell_j(m)  \\
&=&\, ( j+1)(2\tau +2+\vartheta-\ell_j(m)) + 2(\nu_0(m)+\cdots+\nu_{j}(m)).
\end{array}
\end{equation}
Since  $r_0=s_0 r>\sigma_0 r$, 
the inequality  \eqref{eq:approximation6-2} will follow from the following one 
\begin{equation}
\label{eq:approximation6-3}
\sup_{0\le t \le a}\|\partial_t^p P_t\|_{\ell_j(m)}\,  \le \,    \hat\varepsilon \, \epsilon_j(m)    \delta^{M_j(m)} , \quad 0\le p\le k, 
\end{equation}
where
\begin{equation}\label{eq:epsilon-ell}
\epsilon_j(m):= \sigma_0^{\tau + 2}  E_0 \frac{\left(\delta/ 6 s_0\right)^{\ell_j(m)}}{8C_{\ell_j(m)}} .
\end{equation}
We have
\[
\epsilon_j(m) 
= \left\{
\begin{array}{lcrr}
\epsilon_0 =  \sigma_0^{\tau + 2}  E_0 \frac{\left(\delta/ 6 s_0\right)^{\ell_0}}{8C_{\ell_0}} \quad  &\mbox{for}& \ j<J(m); \\ [0.3cm]
\epsilon_m=  \sigma_0^{\tau + 2}  E_0 \frac{\left(\delta/ 6 s_0\right)^{\ell_m}}{8C_{\ell_m}}   \quad  &\mbox{for}&\ j\ge J(m). 
\end{array}
\right.
\]
%depends only on $n,\tau, \vartheta, \vartheta_0$ and $m$. 
%If   $m=0$, then we  have 
%$M(j,\ell_0,(\nu(0)))=-(\vartheta_0-2\vartheta)(j+1)$ for any $j\in\N$,  and the estimates \eqref{eq:approximation6-3}  become
%\begin{equation}
%\label{eq:smallness0}
%\sup_{0\le p \le k} \,  \sup_{0\le t \le a} \| \partial_t^p P_t\|_{\ell_0} \le   \hat\varepsilon \, \epsilon(0)  \delta^{-(\vartheta_0-2\vartheta)(j+1)}\,   , \quad  \forall j\in\N.
%\end{equation}
For $m=0$, taking into account \eqref{eq:secuence-nu-m} and \eqref{eq:secuence-ell-m}, we obtain
\[
M_j(0) = - (j+1)\vartheta.
\]
For $m>0$ we obtain in the  same way
\begin{equation}\label{eq:M(m)}
M_j(m) 
= \left\{
\begin{array}{lcrr}
-(j+1)\vartheta \,  ,  &j&< J(m),\\ [0.5cm]
-j(2m(\tau+1)+\vartheta) -\vartheta  \, ,  &j&=J(m),\\ [0.5cm]
M_{J(m)}(m) - (j-J(m))\vartheta \, ,   &j&\ge J(m).
\end{array}
\right.
\end{equation}
%The relation  \eqref{eq:approximation6-3} becomes 
%\begin{equation}
%\label{eq:smallness0}
% \left\{
%\begin{array}{lcrr}
%\displaystyle \sup_{0\le p \le k} \, \sup_{0\le t \le a} \|\partial_t^p P_t\|_{\ell_0} \le  \hat\varepsilon \, \epsilon_0 \delta^{-(j+1)(\vartheta_0-2\vartheta)}  , \  &j&<J(m) ;\\ [0.3cm]
%  \displaystyle   \sup_{0\le p \le k} \, \sup_{0\le t \le a} \|\partial_t^p P_t\|_{\ell_m} \le   \hat\varepsilon \, \epsilon_m  \delta^{-(j+1)(m(\tau+1)+\vartheta_0-2\vartheta)}\delta^{2m(\tau+1)}   , \ & j&= J(m); \\ [0.3cm]
 % \displaystyle   \sup_{0\le p \le k} \, \sup_{0\le t \le a} \|\partial_t^p P_t\|_{\ell_m} \le   \hat\varepsilon \, \epsilon_m  \delta^{-(j+1)(m(\tau+1)+\vartheta_0-2\vartheta)}\delta^{2m(\tau+1)}   , \ & j&\ge J(m) .
%\end{array}
%\right.
%\end{equation}

Our aim now is to satisfy  \eqref{eq:approximation6-3} for each $m\in\N$, $j\in \N$, and $0\le p\le k\le 1$, choosing  appropriately  $\hat\varepsilon\in (0,1]$ and $J(m)$.  

\vspace{0.3cm}
\noindent
{\em\large (1) \quad The case when $p=0$.}
\vspace{0.2cm}

\noindent
Suppose firstly that  $m=0$. Then 
%To illustrate the construction, we firstly  take  $m=0$,  set 
%This case will be used to prove item  3 of Remark \ref{rem:kam}. $\ell =\ell_0=\ell(0)= \tau + 1+ \vartheta_0$,  and  consider the sequence  $\nu=(\vartheta)$, where  $\nu_j=\vartheta=\vartheta_0/6$, $ j\in\N$. 
%We have 
 \eqref{eq:approximation6-3}    becomes
\begin{equation}
\label{eq:smallness1}
\sup_{0\le t \le a} \|P_t\|_{\ell_0} \le   \hat\varepsilon \, \epsilon_0      \delta^{-\vartheta j}\,   , \quad  \forall j\in\N.
\end{equation}
%To satisfy \eqref{eq:smallness1} for $p=0$   we fix the small  constant $\epsilon$ in \eqref{eq:smallness-condition3} and $\hat\varepsilon_0$ by 
Let us set 
\begin{equation}
\label{eq:small-constant}
\left\{
\begin{array}{rcll}
\epsilon \ &:=&\  \displaystyle \epsilon_0   = \sigma_0^{\tau + 2}  E_0 \frac{\left(4 s_0\right)^{-\ell_0}}{4C_{\ell_0}}   , \\ [0.3cm]
\hat\varepsilon\ &:=&\  \displaystyle  \varepsilon^{-1}\sup_{0\le t \le a} \|P_t\|_{\ell_0}. 
\end{array}
\right.
\end{equation}
Then  \eqref{eq:smallness1} holds for any $j\in\N$.  Moreover,  \eqref{eq:smallness-condition3} just means that  
$$
0<\hat\varepsilon_0\le 1. 
$$ 
Notice that $\epsilon $ depends only on $n$, $\tau$, $\vartheta_0$ and $\vartheta$ since $c_{\ell_0}$, $\sigma_0$, $s_0$ and $E_0$ depend only on $n$, $\tau$, $\vartheta_0$ and $\vartheta$ by  Lemma \ref{Lemma:a b and c},  \eqref{eq:recurrence-relations} and \eqref{eq:sigma-0}. Hence, one can apply the Iterative Lemma. 
%In this way one can obtain a family of continuous with respect to $t$ invariant tori of finite H\" older regularity. 
%In order to get $C^\infty$ families (with respect to $\omega$)  of regular tori and to prove  item (ii) of Theorem \ref{Theo:A},    we set 
%\begin{equation}\label{eq:secuence-ell-m}
%\ell_j=\ell_0 \ \mbox{for}\  j<J(m), \quad  \ell_m=\ell(m)=m(\tau +1) +\ell_0, 
%\end{equation}
%where $J(m)$ is a suitable integer satisfying \eqref{eq:J(m)}.  

%set  $J(q)=0$ for $q<m$ and we define the sequences $ \nu(m):=(\nu_j(m))_{j\in\N}$ and $ (\ell_j)_{j\le J(m)}$  by
%and  $(\ell_j)_{j\le J(m)}$ 
%Given an integer $J(m)\ge 1$ which will be fixed below, we set 
%\begin{equation}
%\ell_j = \left\{
%\begin{array}{lcrr}
%\ell_0 \quad \mbox{for}\quad j < J(m), \\[0.3cm]
%2m(\tau + 1) + \ell_0\quad  \mbox{for}\ j \ge J(m)
%\end{array}
%\right.
%\label{eq:l-j}
%\end{equation} 
%and  
%Notice that  \eqref{eq:approximation6-3} is satisfied for  $0\le  j< J(m)$ in view of \eqref{eq:smallness1} and \eqref{eq:small-constant}.  

Suppose now that $m>0$. If  $j<J(m)$, then $\ell_j(m)=\ell_0$ and $M_j(m)=-(j+1)\vartheta$, and \eqref{eq:approximation6-3} for $p=0$ reduces to \eqref{eq:smallness1} with  $\epsilon$ 
and   $\hat\varepsilon$ given by \eqref{eq:small-constant}. 

 On the other hand, for  $j = J(m)$ and any $k\in\{0;1\}$,  the inequality   \eqref{eq:approximation6-3}  becomes 
\begin{equation}
\label{eq:smallness-A}
A_j^k(m):= C_0(m) \delta^{j(2m(\tau+1) +\vartheta)}  \sum_{p=0}^k \ \sup_{0\le t \le a}\|\partial_t^p P_t\|_{\ell(m)} \le    \hat\varepsilon ,
\end{equation}
where $C_0(m) = \delta^{\vartheta} \epsilon_m^{-1}$. 

The sequence $(A_j^0(m))_{j\in \N}$ is decreasing and it tends to zero. 
Let $J(m)$ be the smallest integer 
\[
j \ge m(\tau+1)\vartheta^{-1}  
\]
such that $A_{j}^0(m)\le \hat\varepsilon$.  Then \eqref{eq:approximation6-3}  holds for $j=J(m)$ and $p=0$. 
Moreover,  $J(m)$ satisfies \eqref{eq:J(m)} by definition.  For $j\ge J(m)$ we have 
\[
\sup_{0\le t \le a}\| P_t\|_{\ell_j(m)}\,  =  \sup_{0\le t \le a}\| P_t\|_{\ell_m}\,  \le \,    \hat\varepsilon \, \epsilon_m    \delta^{M_{J(m)}(m)} \le   \hat\varepsilon \, \epsilon_m    \delta^{M_{j}(m)} 
\]
in view of \eqref{eq:M(m)}, hence,  \eqref{eq:approximation6-3}  is satisfied for each $j\in\N$ when $p=0$. 
\begin{Lemma}\label{lem:A-m}
There exist 
$$
\widetilde C_m=\widetilde  C_m(n,\tau,\vartheta,\vartheta_0)> 0,
$$
depending only on $n,\tau,\vartheta,\vartheta_0$ and $m$, such that
\begin{equation}\label{eq:A-m}
\widetilde  C_m \hat\varepsilon   \le A_{J(m)}^0(m)  \le \hat\varepsilon . 
\end{equation} 
\end{Lemma}
{\em Proof}. \quad 
If $J(m)-1\ge m(\tau+1)\vartheta^{-1}$, then $A_{J(m)-1}^0(m) >\hat\varepsilon$, and we get 
$$
\hat\varepsilon \ge A_{J(m)}^0(m)= \delta^{2m(\tau+1) +\vartheta}A_{J(m)-1}^0(m) >\hat\varepsilon  \delta^{2m(\tau+1) +\vartheta}.
$$
If $J(m)<m(\tau+1)\vartheta^{-1} + 1$, then 
$$
\begin{array}{rcll}
\displaystyle \hat\varepsilon  &\ge&   \displaystyle A_{J(m)}^0(m)   \ge 
C_0(m)   \delta^{J(m)(2m(\tau+1) +\vartheta)}  \sup_{0\le t \le a}\|P_t\|_{\ell_0}\\
 &\ge&      \displaystyle C_0(m  )\delta^{b(m)}  \epsilon  \,   \hat\varepsilon  , 
 \end{array}
$$
where $b(m)= (m(\tau+1)\vartheta^{-1} + 1)(2m(\tau+1) +\vartheta)$, 
and we obtain \eqref{eq:A-m} since $ C_0(m)$, $ \epsilon$  and $\delta$ depend only on  $n,\tau,\vartheta, \vartheta_0$ and $m$. 
\finishproof

\vspace{0.2cm}
\noindent
{\em\large (2)\quad The case when $p=k=1$.}

\vspace{0.2cm}
\noindent
Choosing  $\hat\varepsilon$ and  $J(m)$ as in the case (1),  we obtain that $P_t$ satisfies   \eqref{eq:approximation6-3}  for $p=0$. 
To satisfy  \eqref{eq:approximation6-3}  for  $p=1$,  we need an additional argument. 
 We rescale $t$ by setting  $\tilde t = t T(m)\in [0,\tilde a(m)]$, where $\tilde a(m) = a T(m)$ and 
\begin{equation}\label{eq:T}
\displaystyle T(m)\, :=\,   \frac{1}{\epsilon\hat \varepsilon} \,   \sup_{p\in\{0;1\}}  \,   \sup_{0\le t \le a} \|\partial_t^p P_t\|_{\ell_m} = 
 \frac{\displaystyle  \sup_{p\in\{0;1\}}  \,   \sup_{0\le t \le a} \|\partial_t^p P_t\|_{\ell_m} }{\displaystyle  \sup_{0\le t \le a} \|P_t\|_{\ell_0} }    \ge  1.  
\end{equation}
Then $\widetilde P_{\tilde t}$ defined by $\widetilde P(\cdot,\tilde t):= P(\cdot,\tilde t/T(m)$ satisfies  \eqref{eq:approximation6-3} for $p=0$. Moreover, 
%the derivative of $\widetilde P_{\tilde t} :=  P_{\tilde t/T}$ with respect to $\tilde t$ satisfies \eqref{eq:smallness0} by construction, hence,   $\widetilde P_{\tilde t}$ verifies \eqref{eq:smallness0} for $k=1$. 
%Indeed,
\[
\|\partial_{\tilde t} \widetilde P_{\tilde t}\|_{\ell_m}  \le   \frac{1}{T(m)} \sup_{0\le t \le a} \|\partial_{t} P_{ t}\|_{\ell_m} \le \sup_{0\le t \le a}   \|P_{t}\|_{\ell_0} 
\]
%and
%\[\|\partial_{\tilde t} \widetilde P_{\tilde t}\|_{\ell_m} = \|P_{t}\|_{\ell_0} \le  \|P_{t}\|_{\ell_m} \]
and we obtain  \eqref{eq:approximation6-3} for $p=1$. 
Replacing $P_t$ with $\widetilde P_{\tilde t}$, we can  
apply Proposition \ref{Prop:kam-step} to $\widetilde P_{\tilde t}^j$ for $\tilde t\in [0,\widetilde a]$ at each iteration 
(recall that the constants in Proposition \ref{Prop:kam-step} do not depend on $a$). 
%In this way we will obtain a convergent scheme for $\widetilde P_{\tilde t}$. 
%Then scaling back with respect to $t$ we will obtain the corresponding  for $P_t$.  ??????
%To get (ii) for $P_t$ and $q=1$ we use the inequality $T\ge 1$.  
%In this way, 
%replacing  $P_t$ by  $\widetilde P_{\widetilde t}$  and  $a$ by $\widetilde a(m)$,  we may assume that \eqref{eq:smallness1} holds for $0 \le p\le k\le 1$ with $\epsilon$ and $\hat\varepsilon$ given by \eqref{eq:smallness0}. \\

We summarize the above construction by the following
\begin{Lemma}\label{Lemma:approximation}
	Fix the positive constants $\epsilon$ and $\hat\varepsilon\le 1$ by \eqref{eq:small-constant}. Then 
	\begin{itemize}
	\item[(i)] If 	$k=0$  then 
	for each  $m\in\N$  the sequence  $(\widetilde\varepsilon_{\ell_j(m),j,0})_{j\in\N}$ defined in \eqref{eq:approximation-estimates2}    satisfies \eqref{eq:approximation6-1}  and 
	the Iteration Lemma as well as the Modified Iteration Lemma hold for any $m\in \N$; 
	\item[(ii)]  If $k=1$, then  $\tilde P_{\tilde t}$, $\tilde t\in [0, \tilde a(m)]$, satisfies  \eqref{eq:approximation6-1}  and the Iteration Lemma as well as the  Modified Iteration Lemma  hold for any $m\in \N$. 
\end{itemize}
\end{Lemma}
How do the maps ${\mathcal F}_t^j$,$j\in\N$,   constructed by the Iteration Lemma and ${\mathcal H}_t^j$ given by the Modified Iteration Lemma, depend on $m\in\N$? The answer of this question is given in Remark \ref{rem:uniqueness-Iteration} and Remark \ref{rem:uniqueness-Modified}  and we summarize it by the following 
\begin{Lemma}[Uniqueness by construction] \label{Lemma:dependence-of-m}
We have 	the following:
\begin{itemize}
	\item[(i)] The transformations ${\mathcal F}_t^j$ and ${\mathcal H}_t^j$, $j\in\N$, do not depend on $m$ in the following sense. If  $m'\in\N$  and ${\mathcal F}_t^{\prime j}$ and ${\mathcal H}_t^{\prime j}$,  are the corresponding transformations, 
	then ${\mathcal F}_t^j={\mathcal F}_t^{\prime j}$ and ${\mathcal H}_t^j= {\mathcal H}_t^{\prime j}$ on the intersection of their domains of definition. 
	\item[(ii)] Let $k=1$ and $\widetilde {\mathcal F}_{\tilde t}^j$ and $\widetilde {\mathcal H}_{\tilde t}^j$, $j\in\N$, be the transformations corresponding to $\widetilde P_{\tilde t}$, where $\tilde t = t T(m)\in  [0,T(m)a]$. Let  ${\mathcal F}_t^j$ and ${\mathcal H}_t^j$ be the transformations corresponding to $P_t$,  $t\in [0,a]$. Then $\widetilde {\mathcal F}_{\tilde t}^j={\mathcal F}_t^j$ and $\widetilde {\mathcal H}_{\tilde t}^j={\mathcal F}_t^j$. 
\end{itemize}
\end{Lemma}	
Item (ii) means that the map $[0,a] \to {\mathcal H}_t^j$ is $C^1$, if $k=1$, and that
\[
 \partial_t {\mathcal H}_t^j = T(m)\partial_{\tilde t} \widetilde {\mathcal H}_{\tilde t} ^j |_{\tilde t=t T(m)}.
\]
In order to prove (ii), Theorem \ref{Theo:A}, we need  the following
\begin{Lemma}\label{lem:estimate-eps-b}
	For each   $m\ge 0$  there exists 	$C_m=C_m(n,\tau,\vartheta, \vartheta_0)>0$ depending only on $m,n,\tau,\vartheta,\vartheta_0$ such that 
	\[
	\hat\varepsilon E_{j}(m)\,  \le\,   C_m \sigma_{j+1}^{m(\tau+1)+ \vartheta_0-\vartheta} \,   \sup_{0\le t \le a}\|P_t\|_{\ell(m)} .
	\]
\end{Lemma}
{\em Proof.}\quad Let  $m=0$. We have 
\[
 E_{j}(0)= \delta^{\nu_0(0)+\cdots + \nu_{j-1}(0)}E_0 = \delta^{j(\vartheta_0-\vartheta)}E_0 < (\sigma_0\delta)^{-\vartheta_0}\sigma_{j+1}^{\vartheta_0-\vartheta}E_0
\]
and the estimate holds in view of the choice of  $\hat\varepsilon$ in \eqref{eq:small-constant}. 

Suppose now that $m >0$. 
Using Lemma \ref{lem:A-m} and \eqref{eq:smallness-A}  we obtain 
\[ 
\begin{array}{rcll}
	\displaystyle \hat\varepsilon E_{j}(m) &\le& \displaystyle \widetilde  C_m ^{-1} A_{J(m)}^0(m)E_{j}(m) \\ [0.3cm]
%		&=&  \displaystyle \widetilde  C_m ^{-1}  C_0(m) \delta^{(J(m)+2)(m(\tau+1) +\vartheta_0-2\vartheta)} E_{j}(m)\sup_{0\le t \le a}\| P_t\|_{\ell(m)} \\  [0.3cm]
		 &=&  \displaystyle C_m ^{-1} C_0(m)  \delta^{  F_m(j) } \sup_{0\le t \le a}\| P_t\|_{\ell(m)} ,
\end{array}	
 \]
 where
 \[ 
 F_m(j) := J(m)(2m(\tau+1) +\vartheta) +\nu_0(m) + \cdots +\nu_{j-1}(m).
  \]
  Let  $ j\le J(m)$. Using \eqref{eq:secuence-nu-m}, we get 
  $$
  \begin{array}{rcll}
  F_m(j) &=&   J(m)(2m(\tau+1) +\vartheta) + j(\vartheta_0-\vartheta )\\ 
  &\ge &  j(2m(\tau+1) +\vartheta_0) . 
  \end{array}	
  $$ 
 If $ j\ge J(m)+1$, we obtain by  \eqref{eq:secuence-nu-m} the inequality
  $$
  \begin{array}{rcll}
  \displaystyle 
  F_m(j) &=& J(m)(2m(\tau+1) +\vartheta) 
+J(m) (\vartheta_0-\vartheta)\\
&+& (j-J(m))(m(\tau+1)+\vartheta_0-\vartheta)\\
 &\ge& j(m(\tau+1) +\vartheta_0-\vartheta).
  \end{array}	
  $$
  Choosing $C_m :=  \displaystyle \widetilde  C_m ^{-1}  C_0(m)(\sigma_0\delta)^{-m(\tau+1) -\vartheta_0}$ we complete the proof of the Lemma.
  \finishproof

%-----------------------------------------------------------------------------------------------------------------------

\noindent

We are ready to prove Theorem  \ref{Theo:A}. Fix the parameter $m$ and set 
\begin{equation}\label{eq:B-p-m}
 \left\langle P \right\rangle_{\ell(m)}^{(p)}  :=\left\langle P \right\rangle_{\ell(m);1,1} ^{(p)}, \quad 0\le p\le k, 
\end{equation}
 using  the notations in \eqref{eq:main-estimates1}. If $p=1$ we scale back with respect to $t$ by $T(m)$. 
 Combining \eqref{eq:estimates-H}  and Lemma \ref{lem:estimate-eps-b} and using Lemma \ref{Lemma:dependence-of-m},  (ii),  in the case when $p=1$,  we obtain the estimate
\[
\begin{array}{lcrr}
\displaystyle  \left|\partial_t^p \partial_\theta^\beta\partial_\omega^\gamma \left({\mathcal H}^{j+1}-{\mathcal H}^{j}\right)(\theta,I;\omega,t)\right|
\\ [0.3cm]
 \, \le \, 
\displaystyle C_{\beta,\gamma} \frac{\varepsilon_j}{r_jh_j}\sigma_{j+1}^{-|\beta|}h_{j+1}^{-|\gamma|}  \ln^{2|\gamma|}(1/\sigma_{j+1}) T(m)^p\\ [0.3cm]
\displaystyle  \le 	 C _{m,\beta,\gamma} \,  \sigma_{j+1}^{m(\tau+1)+ \vartheta_0-\vartheta}\sigma_{j+1}^{-|\beta|}h_{j+1}^{-|\gamma|}  \ln^{2|\gamma|}(1/\sigma_{j+1}) \,   \left\langle P \right\rangle_{\ell(m)} ^{(p)}\\ [0.3cm]
\displaystyle  \le 	 C _{m,\beta,\gamma} \,  \sigma_{j+1}^{(m-|\gamma|)(\tau+1)-|\beta|+ \vartheta_0-\vartheta} \ln^{4|\gamma|}(1/\sigma_{j+1}) \,   \left\langle P \right\rangle_{\ell(m)}^{(p)} \\ [0.3cm]
\displaystyle  \le 	 C _{m,\beta,\gamma} \,  \sigma_{j+1}^{(m-|\gamma|)(\tau+1)-|\beta|+ \vartheta_0-2\vartheta}  \,    \left\langle P \right\rangle_{\ell(m)}^{(p)}
	\end{array}
\]
on $D_{j+1}\times \R^n \times [0,a]$, where $ C _{m,\beta,\gamma}$ stands for possibly different positive constants.    Let us fix $\vartheta=(\vartheta_0-\vartheta_1)/4$. Then 
\begin{equation}\label{eq:estimate-H1}
\begin{array}{rcll}
 \left|\partial_t^p \partial_\theta^\beta\partial_\omega^\gamma \left({\mathcal H}^{j+1}-{\mathcal H}^{j}\right)(\theta,I;\omega,t)\right| 
& \le &	 C _{m,\beta,\gamma} \,  \sigma_{j+1}^{(m-|\gamma|)(\tau+1)-|\beta|+ \vartheta_1 +2\vartheta}  \,    \left\langle P \right\rangle_{\ell(m)}^{(p)} \\ 
&\le &  C _{m,\beta,\gamma} \,  \sigma_{j+1}^{ 2\vartheta}  \,    \left\langle P \right\rangle_{\ell(m)}^{(p)}
  	\end{array}
\end{equation}
on $D_{j+1}\times \R^n \times [0,a]$ provided that $|\beta|+|\gamma|(\tau+1) \le m(\tau+1) + \vartheta_1$.  
Set 
\[
%{\cal F}^{j}(\theta,I;\omega,t)= (\Phi^j(\theta,I;\omega,t),\phi^j(\omega,t)) \quad \mbox{and} \quad 
{\cal H}^{j}(\theta,I;\omega,t)= ( \Phi^j(\theta,I;\omega,t),\phi^j(\omega,t))
\] 
where  $\Phi^j(\theta,I;\omega,t)= (U^j(\theta;\omega,t),V^j(\theta,I;\omega,t))$ and $V^j$ is affine linear in $I$ by construction.  Set
\begin{equation}
\label{eq:limits}
\left\{
\begin{array}{lcrr}
\displaystyle \Psi(\theta;\omega,t)=\Psi_t(\theta;\omega)=(U(\theta;\omega,t),V(\theta;\omega,t))=\lim_{j\to\infty} \Phi^{j}(\theta,0;\omega, t),\\[0.3cm]
\displaystyle \phi(\omega,t)=\phi_t(\omega)=\lim_{j\to\infty} \phi^{j}(\omega, t),\ (\theta,\omega,t)\in \T^n \times \R^n\times [0,a]. 
\end{array}
\right.
\end{equation}
Lemma \ref{Lemma:dependence-of-m} implies that the transformations $\Psi$ and $\phi$ do not depend on the choice of $m$. Then 
it follows from \eqref{eq:estimate-H1} that the function $[0,a]\ni t\mapsto (\partial_\theta^\beta\partial_\omega^\gamma\Psi_t, \partial_\omega^\gamma\phi_t) \in C^{\vartheta_1}$ is  $C^k$ for  for any $m\ge 0$ and $\alpha,\beta\in \N^n$ such that $|\alpha| + |\beta| (\tau+1) \le m(\tau +1)+\vartheta_1$. Moreover, the estimates in Theorem \ref{Theo:A}, (ii), hold (here $\kappa=\rho =1$). 
 
We are going to prove (i). To this end we use the identity 
\begin{center}
$ {\cal H}^j = {\cal F}^j$ on $D_{j}^2\times O_{j}^2\times [0,a]$ 
\end{center}
given in  Proposition \ref{Prop:ModifiedIterativeLemma}, (i).
As in Sect. 5.d, \cite{Poe1}, we obtain that
\[
\left| X_{H^j} \circ {\cal F}^j - D \Phi^j \cdot X_N\right| \le 
\frac{c\varepsilon_j}{r_jh_j}
\]  
on ${\T}^n \times \{0\}\times \Omega_1$ for all $j\ge 0$,
where $ X_{H^j}$ and $X_{N}={\mathcal L}_\omega$ stand for the Hamiltonian vector fields of 
$H^j(\theta,I;\omega,t)$ and $N(\theta,I;\omega)=\langle \omega,I \rangle $,
respectively. On the other hand, $\nabla_{(\theta,I)} H^j$ converges uniformly to 
$\nabla_{(\theta,I)} H$ as $j\to \infty$ in view of the estimate \eqref{eq:P-j}, with $\ell=\ell_0$ and $\ell'=\vartheta_1$ 
%\eqref{eq:approximation-estimates1} and \eqref{eq:approximation-estimates2}, 
hence, 
$$
X_{H(\cdot;\phi(\omega,t),t)} \circ \Psi(\cdot;\omega,t) = D \Psi(\cdot;\omega,t)  \cdot {\mathcal L}_\omega
$$
on ${\T}^n \times \{0\}\times \Omega_1$. Moreover,  \eqref{eq:estimate-H1}  implies that
$$
U(\cdot,\omega,\cdot)= \lim U_{j}(\cdot,\omega,\cdot)\in C^k([0,a], C^{1+\vartheta}(\T^n)) , \quad \mbox{\rm for each}\ \omega \in \Omega, 
$$ 
and 
\[
\|U(\cdot;\omega,t) - {\rm id} \|_{1+\vartheta}  <
 C(n,\tau, \vartheta_0,\vartheta_1)\left\langle P \right\rangle_{\ell(m)}^{(0)}   \le C(n,\tau, \vartheta_0,\vartheta_1)\epsilon <  1/2 ,
\]
choosing $\epsilon $ small enough in a function of $n$, $\tau$,$\vartheta_0$ and $\vartheta_1$, 
hence, $U(\cdot,\omega,t)$ is  an embedding. 
Then 
\[
t\to \{\Psi(\theta;\omega,t):\, \theta\in {\T}^n\}
\] 
is a $C^k$ family of  embedded  invariant
tori  
of the Hamiltonians $(\theta,I) \to H(\theta,I;\phi(\omega,t),t)$
with frequency $\omega \in \Omega_1$. They are Lagrangian by
construction (see also \cite{Her}, Sect. I.3.2). 
%In particular, we prove Remark \ref{rem:kam}, 3,(i).

Using Remark \ref{rem:analytisity-t-KAMstep} and Cauchy one obtains 
\begin{Remark}\label{rem:analyticity-t-iteration}
	If $P^j$ are analytic with respect to $t$ in $B(0,a)$ and satisfy \eqref{eq:approximation6-0} for $t\in B(0,a)$, then ${\cal F}_j$ are analytic with respect to $t$ in $B(0,a)$ and the estimates \eqref{eq:estimates-F1} and \eqref{eq:estimates-F2} hold for $p=0$ and $t\in B(0,a)$. Moreover, $\Psi$ and $\phi$ are analytic in $t$ in $B(0,a)$. 
\end{Remark}

\subsection{KAM theorem with parameters in H\"older classes.}   
\label{Sect.ii}

%Fix $m\in \N$ and $0<\vartheta_1< \vartheta_0 (\tau+1)^{-1}$. Recall that $\ell_0=\tau+2+\vartheta_0$, and $\ell(m)= m(\tau+1) + \ell_0$, $m\in \N$, where $0<\vartheta_0<\tau+1$. We suppose that $\ell(m)\notin \N$ for each $m\in \N$.   ????  
%Fix $k\in\{0;1\}$. 

Better  H\"older estimates  of the transformations $\Psi_t$ and $\phi_t$ then those in (ii) Theorem \ref{Theo:A}  can be obtained by means of the anisotropic H\"older spaces $C^{\rho(m)}(\T^n\times \Omega)$ introduced by P\"oschel \cite{Poe}, where 
$$
\rho(m)=((m(\tau+1)+\vartheta_1,m+ \vartheta_2) , \ m\ge 0,\ 1 < \vartheta_1<\vartheta_0, \ \vartheta_2=  (\vartheta_0-\vartheta_1)/(4\tau+4).
$$
Denote the corresponding weighted H\"older norms by $\|\cdot\|_{\rho(m);\kappa}$.

\begin{Theorem}\label{Theo:Holder}
	There exists a  positive constant $\epsilon= \epsilon(n,\tau,\vartheta_0, \vartheta_1)>0$  depending only  on $n$, $\tau$, $\vartheta_0$ and  $\vartheta_1$ such that, for any 
	$a>0$, $0<\kappa<1$, $0<r < \rho_0$ and $M\ge 0$, and  any real valued Hamiltonian  $H=N+P$, where the perturbation  $P\in C^k\big([0,a];C_0^{\ell(M)}(\A\times\Omega)\big)$ satisfies   \ the smallness condition
	\begin{equation}
	\sup_{t\in [0,a]} \, \| P_t\|_{\ell_0;r,\kappa} \  \le \ \epsilon \kappa r \, ,
	\label{eq:smallness-condition-4}
	\end{equation}
and $N(I;\omega)=\langle\omega,I\rangle $ is the normal form, 	the following holds.
	
	There exist families of maps 
	\begin{equation}\label{eq:maps-Holder}
	[0,a]\ni t \mapsto \phi_t\in C^{M+\vartheta_2}(\Omega;\Omega)\, ,\ [0,a]\ni t \mapsto \Psi_t=(U_t,V_t)\in C^{\rho(M)}(\T^n \times\Omega;\T^n\times B(0,r)) \ 
	\end{equation}
	such that   ${\rm supp\,}(\phi_t-{\rm id})\subset \Omega-\kappa/2$, ${\rm supp\,}\big((U_t,V_t)-({\rm id}_{\T^n},0)\big)\subset \T^n\times (\Omega-\kappa/2)$ and   item (i) of Theorem \ref{Theo:A} holds true.  
	
Moreover, for any  $0\le m\le M$ there is $C_m>0$ depending only on  $n$, $\tau$, $\vartheta_0$, $\vartheta_1$,  and $m$,  such that
		\begin{equation}\label{eq:estimatesKAM-Holder}
		\begin{array}{lrc}
		\|\partial_t^q(U_t - {\rm id}_{\T^n})\|_{\rho(m);\kappa}\,  + \,
		r ^{-1}
		\|\partial_t^q V_t\|_{\rho(m);\kappa}
	+ \kappa^{-1} \|
		\partial_t^q(\phi_t - {\rm id})\|_{m+\vartheta_2;\kappa}  	\\    [0.3cm]  
		\displaystyle  	\leq\    C_{m}\, 
		(\kappa  r )^{-1}\, \sup_{0\le p\le q} \,  \sup_{t\in [0,a]} \,\|\partial_t^p P_t\|_{\ell(m);r,\kappa} 
		\end{array}
		\end{equation}
	for each $m\in [0,M]$ and  $t\in [0,a]$. These estimates hold for each $m\in [0,+\infty)$ if $P\in C^k\big([0,a];C_0^{\infty}(\A\times\Omega)\big)$.
	
	 If $P$ is analytic with respect to $t$ in an open  disc  $B(0,a)\subset \C$ of radius $a$ and \eqref{eq:smallness-condition-4} holds for any $t\in B(0,a)$ , then $\phi$ and $\Psi$ are analytic in $t\in B(0,a)$, and the inequalities \eqref{eq:estimatesKAM-Holder}  hold uniformly  in $t\in B(0,a')$, $0<a'<a$,  where  $q=k=0$, the interval $[0,a]$ is replaced by the disc  $B(0,a)$ in the right hand side of \eqref{eq:estimatesKAM-Holder},   and the constant $C_m$ depends on $a'$ as well. 
\end{Theorem}

The estimates follow from \eqref{eq:estimate-H1}, the properties of the norms $\|\cdot\|_{\rho;\kappa}$ for anisotropic H\"older spaces obtained in \cite{Poe} and the Inverse Approximation Lemma obtained by P\"oschel in  \cite{Poe}. 

\begin{Remark}
Can $\ell(m)=  2m(\tau + 1)+ \ell_0$ be  replaced by $ m(\tau+1)+ \ell_0$ ? The loss of  $m(\tau+1)$ derivatives in the estimates (ii) is due to the fact that  we take only the affine linear approximation $Q$ of $P$ with respect to $I$ in the KAM Step Lemma below. Using the approximation proposed by R\"{u}ssmann  in Theorem 7.2 \cite{Rus} as  Bounemora \cite{Boun}, one could prove (ii) with $\ell(m)$ replaced by $m(\tau + 1)+ \ell_0$. This needs  additional efforts and will be done elsewhere. 
\end{Remark}

\appendix
\section{Appendix.}
\setcounter{equation}{0}
\renewcommand{\theequation}{A.\arabic{equation}}

\subsection{Approximation Lemma}\label{Sec:ApprLemma}
The Hamiltonian $P_t$ is not analytic and one can not apply directly the KAM step to it. We are going to approximate it by real analytic functions. To this end 
we recall some facts about the analytic smoothing technique invented by Moser \cite{M1}, \cite{M2},  and developed in different situations by Zehnder \cite{Ze}, P\"oschel \cite{Poe}, Salamon \cite{Sa} and Salamon and Zehnder \cite{Sa-Ze}. The Approximation Lemma and the Inverse Approximation Lemma characterize H\"older classes of differentiable functions in terms of quantitative estimates of approximating sequences of analytic functions. 

 Let $m\in \N$, $0 < \mu \le 1$, and let $U\subset \R^n$ be an open  set.  The H\"older space $C^{m,\mu}(U)$ consists of all $f\in C^m(U)$ such that 
\begin{equation}
\|f\|_{C^{m,\mu}(U)}\,  :=\,  \sup\, (\|f\|_{C^m(U)},  H_{m,\mu}(f)) < \infty,
\label{eq:holder-norm}
\end{equation}
where 
\begin{equation}
\|f\|_{C^m(U)}\,  := \, \sup_{|\alpha|\le m}\, \sup_{x\in U}\, |\partial^\alpha f(x)| 
\label{eq:C-k-norm}
\end{equation}
is the $C^m$ norm of $f$ and 
\begin{equation}
H_{m,\mu}(f)\,  :=\,  \sup\, \frac{|\partial^\alpha f(x)-\partial^\alpha f(y)|}{|x-y|^\mu}
\label{eq:holder-norm1}
\end{equation}
where the supremum is taken over all  $x,y\in U$ such that
$x\neq y$ and all $\alpha=(\alpha_1,\ldots,\alpha_n)\in \N^n$ of length $|\alpha|=\alpha_1+\cdots+\alpha_n=m$. 

Given a  non negative number $\ell = m +\mu\notin \N$, where $m=[\ell]\in \N$ is the entire part of $\ell$ and $0\le \mu=\{\ell\}<1$ the residual one, we set $C^{\ell}(U)=C^{m,\mu}(U)$.
To simplify the notations we set  
\begin{center}
$\|f\|_\ell=\|f\|_{\ell, U}=\|f\|_{C^\ell(U)}$. 
\end{center}

Denote by 
\[
(T^{m}_x f)(y):= \sum_{|\alpha|\le m} \partial^\alpha f(x)y^\alpha/\alpha !
\]
the Taylor polynomial of $f$ up to order $m$. Given $0<\rho \le \infty$ and an open set $U\subset \R^n$ we denote by $U_\rho$ the strip of all $x+iy\in \C^n$ such that $x,y\in \R^n$, $x\in U$ and $|y| < \rho$. Recall that ${\mathcal A}(U_\rho)$ is the  space of analytic  functions on $U_\rho$. We denote by $|\cdot|_\rho$ the sup-norm on $U_\rho$. The function  $f$ is said to be real analytic in $U_\rho$ if $f$ is analytic on $U_\rho$ and real valued on $U$. In this Section we take $U=\R^n$. The space of entire functions ${\mathcal A}(\C^n)$ is endowed by the inductive topology generated by the sup-norms on compact sets of $\C^n$. 
\begin{Lemma}\label{Lemma:approximation-lemma} {\bf(Approximation Lemma)} (\cite{Sa},\cite{Ze}). 
There exists  an entire  function $K\in {\mathcal A}(\C^n)$ generating a family of convolution operators 
\begin{equation}
S_\rho f(x):= \rho^{-n} \int_{\R^n}\, K\left(\rho^{-1}(x-y)\right) f(y) dy\, ,\quad 0<\rho\le 1,
\label{eq:convolution-operators}
\end{equation}
from $C^0(\R^n)$ to ${\mathcal A}(\C^n)$ with the following properties.
\begin{enumerate}
\item
For any $\ell=m+\mu \ge  0$, where $m\in \N$ and $0\le \mu \le 1$,  there is a constant $C=C(n,\ell)>0$ such that,  for every $f\in C^{m,\mu}(\R^n)$, any $\alpha\in\N^n$ of length $|\alpha|\le \ell$ and any $x=u+iv\in\C^n$, $u,v\in\R^n $ with $|v|<\rho$, we have
\begin{equation}
\left|\partial^\alpha S_\rho f (u+iv) -  (T_u^{[\ell]-|\alpha|}\partial^\alpha f)(iv)          \right| \le C \rho^{\ell-|\alpha|}\|f\|_{C^{m,\mu}},
\label{eq:approximation1}
\end{equation}
and in particular for any $0<\rho<\tilde\rho\le 1$, and $f\in C^\ell(\R^n)$,
\begin{equation}
\left|\partial^\alpha S_\rho f  -  \partial^\alpha S_{\tilde\rho} f   \right|_\rho \le C \tilde \rho\, ^{\ell-|\alpha|}\|f\|_\ell . 
\label{eq:approximation2}
\end{equation}
\item The restriction of $S_\rho f $ to $\R^n$ satisfies
\begin{equation}
\left\|S_\rho f  -  f   \right\|_{s} \le C \rho^{\ell-s}\|f\|_{C^{m,\mu}},\quad 0\le s < \ell. 
\label{eq:approximation2-1}
\end{equation}
\item
$K(\R^n )\subset \R$ and in particular the  function $ S_\rho f$ is real analytic whenever $f$ is real valued. Moreover, if  $f$ is periodic in some variables then so is $S_\rho f$ in the same variables.
\item If $[0,1]\ni t\mapsto f_t \in C^\ell(\R^n)$   is a $C^k$ family, then for any $\rho>0$ fixed, the family 
\[
[0,1]\ni t\mapsto S_\rho f_t\in {\mathcal A}(\C^n)
\] 
is  $C^k$ as well and $\left(\frac{d}{dt}\right)^pS_\rho f_t= S_\rho \left(\frac{d}{dt}\right)^pf_t$ for $0\le p \le k$. Moreover, if $t\to f_t(x)$ is analytic in a complex neighborhood $V$ of $t=0$ for each  $x\in \R^n$, then so is $S_\rho f_t(x)$,   and \eqref{eq:approximation1}-\eqref{eq:approximation2-1} are satisfied for $t\in V$. 
\end{enumerate}
\end{Lemma}
A complete proof of the claims 1.-3.  is given for example in \cite{Sa}, Lemma 3, and in \cite{Ze}. In the case of anysotrop H\"older spaces the lemma has been obtained by P\"oschel in \cite{Poe}.  The claim 4. follows easily from the properties of $K$. In order to obtain item {\em 1}, one uses  Taylor's formula with integral remainder, which yields  the estimate
\[
| \partial^\alpha f (u+v) -  (T_u^{[\ell]-|\alpha|}\partial^\alpha f)(v)  \le c    \|f\|_{C^{m,\mu}} |v|^{\ell-|\alpha|} 
\]
for $|\alpha|\le m$ (see (3.4) in \cite{Sa}). \finishproof

Using item {\em 3.} one obtains as in \cite{Sa}, Lemma 5, the interpolation and product estimates.
More precisely,  let $r,s,\ell$ be positive numbers such that  $0\le r<s<\ell$, and $\ell = m+\mu$, $0\le \mu \le 1$. If $\mu<1$, then  $C^{m,\mu}= C^\ell$.  Set $\nu:=(\ell-s)/(\ell-r)$. 
Then there is a constant $c=c_{r,\ell}>0$ depending only on $\ell$ and $r$ such that for any  compactly supported $f\in C_0^{m,\mu}(\R^n)$  the following estimate holds 
\begin{equation}
\|f\|_{C^s} \, \le \, c_{r,\ell} \, \|f\|_{C^r}^\nu\, \|f\|_{C^{m,\mu}}^{1-\nu}.  
\label{eq:interpolation}
\end{equation}  
Moreover, given $f,g\in C_0^{m,\mu}(\T^n\times \R^n)$ one can estimate the $C^\ell$-norm of the product 
\begin{equation}
\|fg\|_{C^{m,\mu}} \, \le \,  \|f\|_{C^0}\|g\|_{C^{m,\mu}} + \|f\|_{C^{m,\mu}}\|g\|_{C^0}.  
\label{eq:interpolation-leibnitz}
\end{equation}

\begin{Remark}\label{rem:interpolation}
	Let $D\subset \R^m$ be an open bounded convex set. Then  
	\[
	\|f\|_{C^s} \, \le \, c_{r,\ell,m} \, \|f\|_{C^r}^\nu\, \|f\|_{C^{\ell}}^{1-\nu}
	\]
 for each $f\in C^\infty(\overline D)$ or $f\in C^\infty(\T^n\times \overline D)$.  Moreover, there exists $C_\ell$ depending only on $\ell$ and on the dimensions $n$ and $m$ such that if $f,g\in C^\infty(\overline D)$ or $f,g\in C^\infty(\T^n\times \overline D)$ then  
 \[
 \|fg\|_{\ell} \, \le \, C_{\ell} (\|f\|_{0}\|g\|_{\ell} + \|f\|_{\ell}\|g\|_{0}).
 \]
\end{Remark}
{\em Proof}. Firstly  we apply  Whitney's extension theorem  (see e.g. \cite{St}, Chapter VI, Theorem 4)  to $f\in C^{k,\mu}(\overline{D})$, where $\overline{D}$ is compact, $k\in\N$, $0<\mu\le 1$.
We obtain an extension 	$\widetilde f\in C^{k,\mu}(\R^m)$  of $f$ such that	  
\[
\|\widetilde f\|_{C^{k,\mu}(\R^m)} \le C_\ell \| f\|_{C^{k,\mu}(D)} , 
\]
where $C_\ell = C_\ell (m) >0$ depends only on $\ell=k+\mu$ and $m$.  Moreover, 
\[
\| f\|_{C^{k,1}(D)} \le   \| f\|_{C^{k+1}(D)}, 
\]
since $f$ is $C^\infty$-smooth on the compact $\overline{D}$ and $D$ is convex. 
Then we  apply the interpolation inequalities \eqref{eq:interpolation} to the extension $\widetilde f\in C^{k,\mu}(\R^m)$  of $f$. 
In the same way we prove the product estimate.
\finishproof

Consider a subdivision $x=(x^{(1)}, \ldots, x^{(p)}) \in \R^{n_1}\times \cdots \times \R^{n_p}$, $n= n_1+\cdots +n_p$, where $1\le p\le n$. 
Given  $a=(a_1,\ldots,a_p)$, where $a_j$ are positive numbers for $1\le  j\le p$, we 
denote by $\sigma_a:\R^n\to\R^n$ the  dilation $\sigma_a(x) = (a_1 x^{(1)}, \ldots, a_p x^{(p)})$. More generally, for any $d\in \N$ we denote by 
$\sigma_a:\T^d\times \R^n\to\T^d\times \R^n$ the partial dilation $\sigma_a(\theta,x) = (\theta, a_1 x^{(1)}, \ldots, a_p x^{(p)})$ (by convention $\T^0=\{0\}$) and 
define the ``$a$-weighted'' H\"older norm of $f\in C^\ell (\T^d\times \R^n)$ by
\begin{equation}
\big\|f\big\|_{\ell;a}:= \big\|f\circ \sigma_{a}\big\|_{\,  C^\ell } .
\label{eq:holder-norms-A}
\end{equation}

\subsection{Almost analytic Gevrey extensions}\label{Sec:Gevrey}
{\em Proof of Proposition \ref{prop: almost-analytic} }. 
Fix $\rho>1$. Let  $\varphi$ be a real valued 
compactly supported  Gevrey function belonging to the class $\mathcal{G}^\rho_\lambda(\R^n)$ for some $\lambda>0$, which means that
\begin{equation}\label{eq:Gevrey-lambda}
\|\varphi\|_{\lambda} := \sup_{\alpha\in{\N}^n}\, 
\sup_{x\in{\R}^n  } \, 
\left(|\partial_x^\alpha  \varphi(x)|\,  \lambda^{-|\alpha|}
\alpha !^{-\rho} \right)\ <\ \infty \, ,
\end{equation}
where 
$|\alpha| = \alpha_1 + \cdots + \alpha_n$ and $\alpha ! =
\alpha_1! \cdots \alpha_n!$ for $\alpha=(\alpha_1,\ldots,\alpha_n)
\in {\N}^n$.  We suppose as well that the support of $\varphi$ is contained in the unit ball $B_1^n(0)=\{x\in\R^n: |x|<1\}$ in $\R^n$ and that 
\[
\int_{\R^n}^{}\, \varphi(x)\, dx \, =\, 1.
\]
Set $U_j^0:=O_j^2\cap \R^n$ and $U_j^q:=\{x\in \R^n:\ {\rm dist\, }(x,U_j^0)<q/16 \}$, $q\in\{0;1;2;3;4\}$, in particular, $U_j^4= O_j^3\cap \R^n$. 
Denote by  $\1_j$ the characteristic function of the set $U_j^2$ in $\R^n$ and consider for any $j\in \N$ the function $f_j$ defined by the convolution
\[
f_j(x)\ =\  (16/h_j)^{n}\, \int_{\R^n}^{}\, \1_j(x-u) \varphi\left( 16u/h_j)\right)\, du \ =\  (16/h_j)^{n}\, \int_{\R^n}^{}\, \1_j(u) \varphi\left(16(x-u)/h_j\right)\, du . 
\]
These functions have the following properties
\begin{equation}\label{eq:properties-f}
\begin{array}{lcr}
	  (1)\quad f_j\in \mathcal{G}^\rho_{\lambda_j}(\R^n) \ \mbox{with} \ \lambda_j= 16 \lambda/h_j \  \mbox{and } \  \|f_j\|_{\lambda_j} \le{\rm vol\, }(B_1^n(0))  \|\varphi\|_{\lambda} ;\\
	  (2)\quad  {\rm supp\, }  f_j \subset U_j^3\  \mbox{and } \    f_j=1 \  \mbox{on } \  U_j^1,
\end{array}
\end{equation}
where the positive constant $\lambda$ is given in \eqref{eq:Gevrey-lambda}. 
We are going to obtain a Gevrey-$\mathcal{G}^\rho$ almost analytic extension of $f_j$ in $\C^n$ which is equal to one on $O_j^2$ and has a support in $O_j^3$.

To  this end we introduce a family of the compact sets in $\R^n\times \R^n$ given by 
\[K_j:= (\R^n\times \{0\} )\cup \overline{O_j^2}\cup (\overline{O_j}\setminus O_j^3), \]
where the set $O_j\subset \C^n$ is identified with the corresponding open set in $\R^n\times \R^n$ via the map $\C^n\ni x+iy \mapsto (x,y)\in \R^n\times \R^n$ 
and $\overline{O_j}$ stands for the closure of $O_j$. Let us extend $f_j$ to a continuous function  with support in $K_j$ by 
\[
\tilde f_j(x,y):=
\left\{
\begin{array}{rcll}
f_j(x)\quad &\mbox{if}& \ x\in \R^n, \ y=0; \\
1  \quad &\mbox{if}& \ (x,y)\in  \overline{O_j^2}; \\
0  \quad &\mbox{if}& \ (x,y)\in  (\R^n\times \R^n )\setminus O_j^3.
\end{array}
\right.
\]
It is easy to see that a formal almost analytic extension of $\tilde f_j$ is given by the power series 
\begin{equation}\label{eq:power-series}
 \sum_{\gamma\in\N^n}^{}\, (iy)^\gamma \partial_x^\gamma f_j(x)/\beta! ,
 \end{equation}
which means that the operators $\bar{\partial}_k$, $k=1, \ldots,n$, annihilate it. 
The  corresponding Taylor series centered at $(x,y)\in K_j$ is 
\[ \sum_{(\beta,\beta')\in\N^n\times \N^n}^{}\, (x'-x)^\beta (iy'-iy)^{\beta'}  \partial^{\beta +\beta'} f_j(x)/(\beta!\beta' !). \]

The family of jets $F_j = \left( f_j^ {(\beta,\beta')}\right)_{\beta,\beta'\in\N^{n}}$ corresponding to the power series given above is 
  defined  for any $ (\beta,\beta')\in \N^n\times\N^n$ and $(x,y)\in K_j$ by
	\begin{equation}\label{eq:jet-F}
	f_j^{(\beta,\beta')}(x,y)\, =\,
	 i^{|\beta'|}\partial^{\beta+\beta'}  f_j(x).
	\end{equation}	
	\begin{Remark}\label{rem:y neq 0}
		If $(x,y)\in K_j$ and  $y\neq 0$, then either $x\in U_j^0$ or $x\notin U_j^4$. On the  other hand, $f_j=1$ on $U_j^1$ and  $f_j=0$ on $\R^n\setminus U_j^3$. Then $f_j^{\alpha}(x,y)=0$ for $\alpha =(\beta,\beta')\neq 0$ and $f_j^{0}(x,y)=1$ if $(x,y)\in U_j^0=O_j^2$, $f_j^{0}(x,y)=0$ if $(x,y)\notin U_j^4=O_j^3$. In particular, 
		\begin{center}
			$f_j^{\alpha}(x,y)=f_j^{\alpha}(x,0)$  for each $\alpha$.
		\end{center}
		\end{Remark}
	
We are going to extend the jet $F_j$ to a Gevrey-$\mathcal{G}^\rho$ function using a Whitney extension theorem in Gevrey classes. 

Let us first  recall the notion of Gevery smoothness of Whitney jets. Let $K$ be a compact set in $\R^d$, $d\ge 1$, and $F=  \left( f^\beta\right)_{\beta\in\N^{d}}$ a jet of continuous functions $f^\beta \in C(K)$.   For each $N\in \N$  we denote by $T^N_u F$ the  formal Taylor polynomial of order $N$ centered at $u\in K$, i.e. 
\[
T^N_u F(z)\,  :=\, \sum_{|\beta|\le N}^{}\, f^\beta(u)(z-u)^\beta/\beta ! \, , \quad z\in \R^d.
\]
Given $\alpha\in \N^d$ we denote by  $F^{(\alpha)}$ the jet $ \left( f^{\alpha+\beta}\right)_{\beta\in\N^{d}}$. Then for $|\alpha|\le N$,  the partial derivative $ \partial^\alpha_{}$ of the Taylor polynomial is given by
 \[
 \partial^\alpha_{z}T^N_u F(z)\,  =\, T^{N-|\alpha|}_u F^{(\alpha)}(z)\, =\,  \sum_{|\beta|\le N-|\alpha|}^{}\, f^{\alpha+\beta}(u)(z-u)^\beta/\beta ! \, .
 \]
 For each $\alpha\in \N^d$, the corresponding Taylor remainder is  defined by
  \[
   \begin{array}{rcll}
\displaystyle   R^N_u F^{(\alpha)}(z)\,  &=&\,  \displaystyle   f^\alpha(z)-\partial^\alpha_{z}T^N_u(F)(z)\\
  &=&\,   f^\alpha(z)-  T^{N-|\alpha|}_u F^{(\alpha)}(z)\\ 
 &=&\,  \displaystyle  f^\alpha(z)- \sum_{|\beta|\le N-|\alpha|}^{}\, f^{\alpha+\beta}(u)(z-u)^\beta/\beta! \, .
  \end{array}
 \]
 Recall from Stein \cite{St} p. 177  the following identity 
 \begin{equation}\label{eq:identity-Stein}
 \partial^\alpha_{z}T^N_v F(z) - \partial^\alpha_{z}T^N_u F(z) \, =\, \sum_{|\beta|\le N}^{}(z-v)^\beta R_uF^{(\alpha+\beta)}(v)/\beta !
\end{equation}
  for any $N\in \N$, $|\alpha|\le N$, $u,v\in K$ and $z\in \R^d$ . 

 Let $L>0$. 
 The jet $F=(f^\beta)_{\beta\in
 	{\N}^d}$ is said to belong to the Whitney space $W{\mathcal G}_L^\rho(K)$ of Gevrey jets  if there exists $A>0$ such that
 \begin{equation}
 \begin{array}{lcrr}
\displaystyle  (1) \quad  |f^\beta(u)| \le A L^{|\beta|} (\beta !)^\rho  \quad  \mbox{for}\  \beta\in
 \N^d\, ,\ u\in K; \\ [0.3cm]
 \displaystyle (2) \quad |R^N_{u} F^{(\gamma)}(z)| \le  A L^{N+1}  ((N+1)!) ^\rho
 |z-u|^{N-|\gamma|+1}/ (N-|\gamma|+1)!  \\
 \quad \quad \ \mbox{for}\  |\gamma|\le N\, ,\ u, z\in  K .
 \end{array}
 \label{eq:gevrey-K}
 \end{equation}
 The corresponding norm of $F$ is defined by $\|F\|_L := \inf A$. The space $W\mathcal{G}_L^\rho(K)$ equipped with this space is a Banach space. We recall the Whitney extension theorem of Bruna \cite{Bru} as it  has been  presented in \cite{P4}, Theorem 3.8.
 \begin{Theorem} \label{thm:Bruna}
 	There exist positive constants $A_0=A_0(d,\rho)$ and 
 	$ C_0=C_0(d,\rho)$ such that the following holds.
 	 
 	 For any compact subset  $K$ of $\R^d$ and
 	 jet $F=(f^\beta)_{\beta\in {\N}^d}\in W\mathcal{G}_L^\rho(K)$, satisfying (\ref{eq:gevrey-K}) on
 	$K$ with some $ L >0$,  there exists 
 	 $ f\in {\cal G}_{C_0L}^{\rho}({\R}^d)$ such that
 	\begin{enumerate}
 		\item[(i)]    $\partial^{\beta} f = f^\beta$ on $K$ for any $\beta$;
 		\item[(ii)]   $\|f\|_{C_0L}  \le A_0\|F\|_L$. 
 		\end{enumerate}
 \end{Theorem}
 We are going to prove that $F_j$ are Whitney jets. 
 Using \eqref{eq:properties-f}  we obtain 
\begin{Lemma}\label{lemma:jet-F}
For each $j\in\N$ the jet $F_j$ belongs to the Whitney space $W{\mathcal G}_{L_j}^\rho(K_j)$ with $L_j= 385 n^2 2^{n\rho} \lambda /h_j$ and $\|F_j\|_{L_j} \le {\rm vol\, }(B_1^n(0))  \|\varphi\|_{\lambda} $. 
\end{Lemma}
{\em Proof.}\quad 
Set $A= {\rm vol\, }(B_1^n(0))  \|\varphi\|_{\lambda} $. 
Using   item (1) of \eqref{eq:properties-f} we obtain for each $j\in\N$ the estimate
\[
|f_j^{(\beta,\beta')}(x,y)| \, <  A \,  \lambda_j^{|\beta|+|\beta'| } ((\beta +\beta' )!)^\rho \,   \le  \, A(2^{n\rho}  \lambda_j)^{|\beta|+|\beta'| } (\beta !)^\rho  (\beta' !)^\rho \
\]
for $ (\beta, \beta')\in
\N^n\times \N^n\, ,\ (x,y)\in K_j$. 

We are going to prove (ii).
For any $u=(x,y)\in K_j$ the formal Taylor polynomial of $F_j$ of order $N$ which is centered at $u$ and  evaluated at $z=(x',y')$,  is given by
\[
T^N_{u}F_j(z)\,  =\, \sum_{|\beta|+|\beta'|\le N}^{}\, i^{|\beta'|} \partial^{\beta+\beta'}f_j(x)(x'-x)^\beta (y-y')^{\beta'}(\beta!\beta' ! )\, .
\]
We consider separately the following two cases.

\vspace{0.3cm}
\noindent
{\em 1.}\quad Let $z=(x',0)$. Suppose at first that $u=(x,0)$. 
Then, setting $M= N+1-|\alpha|-|\alpha'|$, we obtain by Taylor's formula
\[
\begin{array}{lcrr}
R^N_{u} F_j^{(\alpha,\alpha')}(z) = f_j^{(\alpha,\alpha')}(z) -\partial^{(\alpha,\alpha')}_zT^N_{u}F_j(z)\\ [0.3cm]
\displaystyle  =\,  i^{|\alpha'|}\partial^{\alpha+\alpha'}  f_j(x')-\sum_{|\beta|\le N-|\alpha|-|\alpha'|}^{}\,   i^{|\alpha'|} \partial^{\alpha +\alpha'+\beta}f_j (x)(x'-x)^\beta /(\beta! ) \\ [0.5cm]
\displaystyle  =\,  i^{|\alpha'|}  M \sum_{|\beta|=M}^{}\, \frac{(x'-x)^\beta}{\beta !} \int_0^1 \, (1-t)^{M-1} \,\partial^{\alpha +\alpha'+\beta}f_j(x+t(x'-x)) \, dt .
\end{array}
\]
Now  item (1) of \eqref{eq:properties-f} yields
\[
|R^N_{u} F_j^{(\alpha,\alpha')}(z) \le   A L_j^{N+1} ((N+1)!)^\rho |x'-x|^M \sum_{|\beta|=M}^{}\frac{M}{\beta !}\, .
\]
On the other hand
\[
\sum_{|\beta|=M}^{}\frac{M}{\beta !} = \frac{M}{M!}\sum_{|\beta|=M}^{} \frac{( \beta_1+\cdots+\beta_n)!}{\beta_1 !\cdots \beta_n !} =  \frac{M}{M!} n^M<  \frac{(2n)^M}{M!} \, .
\]
Setting $\gamma=(\alpha,\alpha')\in \N^n\times\N^n $ and $\widetilde \lambda_j=  2n \lambda_j =32n\lambda /h_j$ we obtain
\begin{equation}\label{eq:estimate-remainder}
\begin{array}{lcrr}
|R^N_{u} F_j^{(\gamma)}(z) < A  \widetilde \lambda_j^{N+1} ((N+1)!)^\rho  |z-u|^{N-|\gamma|+1}/ (N-|\gamma|+1)! \\ [0.3cm]
 \mbox{for} \ z=(x',0), u=(x,0) \in K_j. 
\end{array}
\end{equation}
Let $u=(x,y)\in K_j$ and $y\neq 0$.  By Remark \ref{rem:y neq 0}  we have $f^\gamma(x,y)=f^\gamma(x,0)$ for each $\gamma\in \N^n\times \N^n$, hence, $R^N_{u} F^{(\gamma)}(z)=R^N_{(x,0)} F^{(\gamma)}(z)$. Then using \eqref{eq:estimate-remainder} we obtain the same estimate since $|z-(x,0)| < |z-u|$. 
Hence \eqref{eq:estimate-remainder} is true for any $z=(x',0)$ and $u=(x,y)$ in $K_j$.

\vspace{0.3cm}
\noindent
{\em 2.}\quad Let $z=(x',y')\in K_j$, $y'\neq 0$ and $u\in K_j$. Set $v=(x',0)$. 
Remark \ref{rem:y neq 0} implies 
\[
f_j^\gamma(z) =  \partial^\gamma_{z}T^N_v F_j(z)
\]
and by means of  \eqref{eq:identity-Stein} we obtain 
\[
R^N_{u} F_j^{(\gamma)}(z) =  \partial^\gamma_{z}T^N_v F_j(z) -\partial^\gamma_{z}T^N_u F_j(z) \, =\, \sum_{|\beta|\le N}^{}(z-v)^\beta R_u^NF_j^{(\gamma+\beta)}(v)/\beta ! .
\]
Now applying \eqref{eq:estimate-remainder} to $ R_u^NF_j^{(\gamma+\beta)}(v)$ we get
\[
|R^N_{u} F_j^{(\gamma)}(z)\ < A  \widetilde{\lambda}_j^{N+1} ((N+1)!)^\rho  \sum_{|\beta|\le N}^{} \frac{|y'|^{|\beta|}|v-u|^{N-|\beta|-|\gamma|+1}}{(N-|\beta|-|\gamma|+1)!\beta! }.
\]
Remark \ref{rem:y neq 0} implies  that  $|x'-x|\ge h_j/16$ on the support of the function 
$u\mapsto R^N_{u} F^{(\gamma)}(z)$ and we get
\[
|y'|\le \frac{2}{3}h_j \le \frac{32}{3}|x'-x|\le  \frac{32}{3}|z-u| \quad \mbox{and} \quad |v-u|\le |x'-x|+|y| \le  \frac{35}{3}|z-u|. 
\]
Setting $L_j = 385 n^2 >
35\widetilde \lambda_j/3$ we obtain as above 
\[
|R^N_{u} F^{(\gamma)}(z) < A   L_j^{N+1} ((N+1)!)^\rho  |z-u|^{N-|\gamma|+1}/ (N-|\gamma|+1)! .
\]
This completes the proof of the lemma. \finishproof

Lemma \ref{lemma:jet-F} enables us to apply Theorem \ref{thm:Bruna} and we denote by $\chi_j$  the corresponding extension of the jet $F_{j+1}\in W{\mathcal G}_{L_{j+1}}^\rho(K_{j+1})$. By construction,  the  Taylor series of $\chi_{j-1}$ at $(x,0)$ coincides with the power series \eqref{eq:power-series}, which implies that  $\chi_j$ is almost analytic. The function $\chi_j$ satisfies item (i) of Proposition 
\ref{prop: almost-analytic} since it coincides with $\tilde f_{j+1}$ on $K_{j+1}$. It satisfies (ii), taking $L= 385 n^2 2^{n\rho}  \lambda C_0$, where $\lambda$ is introduced in \eqref{eq:Gevrey-lambda} and $C_0$ is the constant in Theorem \ref{thm:Bruna}. Fixing for any $n$ the function $\varphi$ and the constants $\lambda$ and  $C_0$, we may suppose that $L=L(n,\rho)$ depends only on $n$ and $\rho$. 

It remains to prove (iii). Expanding $\partial_x^\alpha \partial_y^\beta \chi_j(x,y)$ in
Taylor series at $y=0$,  we obtain
for any $\alpha,\beta \in {\N}^n$ and $m\in {\N}$
$$
|\partial_x^\alpha\partial_{y}^{\beta} \chi_j(x,y)|\ \leq\ 
 A\, ( L/h_{j+1})^{|\alpha|+ |\beta|} 
\left( ( L/h_{j+1})\right)^{m}\, 
\alpha !\,^\rho  \beta !\, ^{\rho}\, m !\, ^{\rho-1}|y|^m\, , 
\quad  (x,y)\in \R^n\times \R^n .
$$ 
Using Stirling's formula 
we minimize the right-hand side  with respect to 
$m\in {\N}$. An optimal choice for $m$ is given by 
$$m\sim  ( L|y|/h_{j+1})^{\, -\frac{1}{\rho-1}},$$ 
which leads to  
\[
|\partial_x^\alpha\partial_y^{\beta} \chi_j(x,y)|\ \leq\ 
C_0 A\,   ( L/h_{j+1})^{|\alpha|+ |\beta|} 
\alpha !\,^\rho  \beta !\, ^{\rho}\,  
\exp\left(- \frac{1}{2}( L|y|/h_{j+1})^{\, -\frac{1}{\rho-1}}\right)  
\]
for any $\alpha,\beta \in {\N}^n$ and $(x,y)\in \R^n\times \R^n$, $0<|y|\le 1, $ with $C_0\ge 1$.   \finishproof

\subsection{Borel's Theorem}\label{Sec:Whitney}

\noindent
{\em Proof of Proposition \ref{prop:realization}}.\quad We follow the standard proof of the Borel's theorem (see \cite{Zw}, Theorem 4.15 and \cite{Ma}, Proposition 2.3.2). To simplify the notations set $z=(x,\xi)\in T^\ast\R^d$.  Consider the increasing sequence $(\eta_j)_{j\in\N}$ where
\[
\eta_j\, :=\, 1\, +\,  \sup_{0\le k\le 1}\sup_{|\alpha|\le j} \sup_{m\le j} \sup_{(t,z)\in I\times T^\ast\R^d} \big|\partial_t^k\partial_z^\alpha a_{t,m}(z) \big|.
\] 
Choose $\chi\in C^\infty(\R)$ such that  $\chi =1$ on $(-\infty,1/2]$, $0< \chi < 1$ on $(1/2,1)$, and  $\chi =0$ on the interval $[1,+\infty)$ and set 
\[
a_t(z,\lambda)= \sum_{j\ge 0} \chi\Big(\frac{\eta_j}{|\lambda|}\Big) a_{t,j}(z) \lambda^{-j}
\]
where $\lambda \in {\mathcal D}$. Notice the the sum is finite for any $\lambda \in {\mathcal D}$ fixed. 
Moreover, $\displaystyle \frac{\eta_j}{|\lambda|} \chi\Big(\frac{\eta_j}{|\lambda|}\Big) \le 1$, hence, for any $k\in\{0,1\}$ and $\alpha\in \N^{2d}$ with $ |\alpha| \le j$ we have  
\begin{equation}\label{eq:inequality-for-a}
\chi\Big(\frac{\eta_j}{|\lambda|}\Big) \big|\partial_t^k\partial_z^\alpha a_{t,j}(z) \big| <  \chi\Big(\frac{\eta_j}{|\lambda|}\Big) \eta_j  \le  |\lambda|.
\end{equation}
On the other hand $|\lambda| \le 2\eta_j$ whenever $\displaystyle \chi\Big(\frac{\eta_j}{|\lambda|}\Big) < 1$, hence, 
\begin{equation}\label{eq:inequality-for-a1}
\Big(1-\chi\Big(\frac{\eta_j}{|\lambda|}\Big)\Big)\big|\partial_t^k\partial_z^\alpha a_{t,j}(z) \big|  \le  C_{p,j,\alpha}|\lambda|^{-p}
\end{equation}
for any $p\in\N$. For any $N\in\N$  we  have
\[
\begin{array}{lcrr}
\displaystyle \Big|\partial_t^k\partial_x^\alpha \Big(a_t(z,\lambda) - \sum_{j=0}^{N-1} a_{t,j}(z)\lambda^{-j}\Big)\Big| \\[0.5cm]
\displaystyle \le
\sum_{j=0}^{N-1}  \Big(1-\chi\Big(\frac{\eta_j}{|\lambda|}\Big)\Big)|\partial_t^k\partial_x^\alpha a_{t,j}(z)||\lambda|^{-j} +\sum_{j\ge N} \chi\Big(\frac{\eta_j}{|\lambda|}\Big)|\partial_t^k\partial_x^\alpha a_{t,j}(z)||\lambda|^{-j}.
\end{array}
\]
Using \eqref{eq:inequality-for-a1} we estimate the first sum by $C_{N,\alpha} |\lambda|^{-N}$. If $N\ge |\alpha|$ taking into account \eqref{eq:inequality-for-a} we estimate the second sum by
$\eta_N |\lambda|^{-N} + \sum_{j\ge N+1}|\lambda|^{-j+1} \le C_{N,\alpha} |\lambda|^{-N}$ for $|\lambda|\ge 2$. Finally, if $N<|\alpha|$ we estimate the first $|\alpha|-N$ terms of the second sum  by $C_{N,\alpha} |\lambda|^{-N}$ and for the other terms we apply the preceding argument. This proves \eqref{eq:amplitude}.   The proof shows as well that for any $k\in \{0,1\}$ and $\alpha\in\N^{2d}$ the function $(t,z)\to \partial_t^k\partial_x^\alpha a_{t}(z,\lambda)$ is a sum of a normally convergent series of functions, hence, it is a continuous  for any $\lambda$ fixed. Thus the map $J\ni t\to a_t$ is a $C^1$ family of symbols of order zero.
\finishproof

\subsection{ Higher order H\"{o}lder estimates of a composition and of the inverse function.}\label{subsec:composition}
\subsubsection{ Estimates of the composition and  the inverse function.}\label{subsec:composition-inverse}
The aim of this section is to obtain estimates of the composition and the inverse function in certain H\"{o}lder norms with constants depending only on the dimension of the spaces and the H\"{o}lder exponent. 
We start by introducing  certain semi-norms as follows. Hereafter we denote the Euclidean norm of $x\in\R^n$ by $|x|$.
Recall that  for any $\ell\in \R$, $[\ell]\in \R$ stands for the entire part of $\ell$ and  $\{\ell\}:= \ell - [\ell]\in [0,1)$ for the residual one. Given $f\in C^\ell(U)$, where $U\subset \R^n$ is an open set,  we introduce a $\ell$-semi-norm of $f$ by 
\begin{equation}\label{eq:semi-norm-l}
 |  f | _{\ell, U} := \left\{ 
\begin{array}{lcrr}
\sup_{|\alpha|=\ell,  x\in U}|\partial^\alpha f(x)| &,& \ \mbox{if}\ \ell\in\N; \\
H_{[\ell], \{\ell\}}(f) &,& \ \mbox{if}\ \ell\notin\N ,
\end{array}
\right.
\end{equation}
where the semi-norm $H_{k, \mu}(f)$ is defined by \eqref{eq:holder-norm1}. To simplify the notations we often write $|f|_{\ell}:= |f|_{\ell,U}$. This notation should not be confused with the sup-norm in Sect. \ref{Sec:ApprLemma}.
The $\mu$-semi-norm of the product of two functions  $f,g\in C^\mu(U)$ with $0<\mu<1$ can be estimated by
\begin{equation}\label{eq:Leibnitz}
|  f g | _{\mu} \le  |  f | _{\mu} |   g | _{0} + |   f | _0  |   g | _\mu. 
\end{equation}
More generally, for any $m\in\N$, $0\le \mu<1$, and $f,g\in C^{m+\mu}(U)$, the Leibniz formula implies 
\begin{equation}\label{eq:Leibnitz1}
|fg|_{m+\mu}\, \le\, C_m\sum_{k=0}^m \Big(|f|_{k+\mu} |g|_{m-k} + |f|_k|g|_{m-k+\mu}\Big),   
\end{equation}
where $C_m$ depends only on $m$ and on the dimension $n$. 

Consider now the composition of H\"older functions.
\begin{Lemma}\label{lemma:composition-mu}
	Let $U_j\subset \R^{n_j}$, $j=1,2$, be open sets. Suppose that $U_1$ is convex. Let $0<\mu<1$ and $f_1\in C^1(U_1)$,  $f_2\in C^\mu(U_1)$. Then
	$|f_2\circ f_1|_{\ell,U_1} \le |f_2|_{\ell,U_2} |f_1|_{1,U_1} ^\ell$.
\end{Lemma}

To estimate higher H\"older norms of the composition of two functions we apply the   Faa di Bruno formula. To this end we introduce the following notations. 
Given $\ell\ge 1$, $f\in C^\ell(U)$  and an integer $1\le m \le \ell$ we set
\begin{equation}
\label{eq:P-l-m}
{\mathcal P}^{\ell,m}_U (f) \ :=\ \displaystyle \sum_{\Delta(\ell,m)}\, \sum_{j=1}^{m} |f|_{k_1} \cdots  |f|_{k_{j-1}} |f|_{k_j +\{\ell\}} |f|_{k_{j+1}}  \cdots |f|_{k_m}
\end{equation}
where the $|f|_{\alpha}:= |f|_{\alpha,U}$ and the index set $\Delta(\ell,m)$  consists of all  
$(k_1,\ldots , k_m)\in \N^m$ 
such that 
\[
k_1 + \cdots + k_m=[\ell]\, , \quad \min_{1\le j\le m}\, k_j \ge 1 .
\] 
For any $\ell\ge 1$ we set
\[{\mathcal P}^{\ell}_U (f):= \sum_{m=1}^{[\ell]}{\mathcal P}^{\ell,m}_U (f).
\]

\begin{Prop}\label{prop:composition}
Let $f_j\in C^\ell(U_j,U_{j+1})$, $j=1,2$, where $\ell \ge 1$ and $U_j\subset \R^{n_j}$ are open sets such that $f_1(U_1)\subset U_2$. Then the following holds:
\begin{enumerate}
	\item If  $\ell\in \N_\ast$ then 
\[
|f_2\circ f_1|_{\ell,U_1}\,  \le\,  C_\ell \sum_{m=1}^\ell |f_2|_{m,U_2}\,  {\mathcal P}^{\ell,m}_{U_1}(f_1)
\]
where the constant $C_\ell$ depends only on $\ell$ and on the dimensions  $n_j$, $j=1,2,3$. 
\item If  $\ell\in [1,\infty)$  and $U_1$ is convex, then 
\[
\begin{array}{lcrr}
|f_2\circ f_1|_{\ell,U_1}\,  \le\,  C_{\ell} \left(1+|f_1|_{1,U_1}^{\{\ell\}}\right)\\[0.3cm] 
\displaystyle \times \sum_{m=1}^{[\ell]}\left(|f_2|_{m+{\{\ell\}},U_2}\,  {\mathcal P}^{[\ell],m}_{U_1}(f_1)+|f_2|_{m,U_2}\,  {\mathcal P}^{\ell,m}_{U_1}(f_1)\right) 
\end{array}
\]
where the constant $C_{\ell}$ depends only on $[\ell]$ and on the dimensions  $n_j$, $j=1,2,3$. 
\end{enumerate}
\end{Prop} 
{\em Proof.} Statement {\em 1.}  follows directly from the Faa di Bruno formula
\[
\partial^\alpha (f_2\circ f_1) = \sum \frac{(\partial^\beta f_2)\circ f_1}{\beta !}\,  (\partial^{\alpha_1} f_1)  \cdots  (\partial^{\alpha_m} f_1) \, \frac{\alpha_1 !\cdots \alpha_m !}{\alpha !}
\]
where the summation is over all the indices 
\[
m\in \N_\ast, \quad \beta\in \N^{n_2}, \quad (\alpha_1, \ldots, \alpha_m)\in \underbrace{\N^{n_1}\times \cdots \times \N^{n_1}}_{m}= \N^{n_1 m}
\]
such that 
\[
1 \le |\beta|=m\le |\alpha|,\quad \alpha_1 +\ldots + \alpha_m = \alpha, \quad \min_{1\le j \le m}|\alpha_j| \ge 1. 
\]
Here $|\beta|=\beta_1 +\cdots +\beta_{n_2}$ stands for the length of  $\beta\in \N^{n_2}$. \\

We are going to prove {\em 2.} Suppose now that  $U_1$ is convex and that $\ell=m+\mu$,  where $m\in \N_\ast$ and $0<\mu<1$.  Firstly we 
apply  \eqref{eq:Leibnitz} to the product in the Faa di Bruno formula.  
Then we  estimate  $|(\partial^\beta f_2)\circ f_1|_{\mu,U_1}$ by means of 
 Lemma \ref{lemma:composition-mu}, which yields
\[
|(\partial^\beta f_2)\circ f_1|_{\mu,U_1} \le |\partial^\beta f_2|_{\mu,U_2} |f_1|_{1,U_1}^{\mu}. 
\] 
This implies {\em 2}. 
\finishproof

Similar inequalities can be proven for  compensated domains $U_1$ (\cite{L-O}, Theorem 5.4) but then the constants depend on $U_1$.

\begin{Prop}\label{prop:inverse}
Let $\ell > 0$ and let $f\in C^{\ell+1}(U,V)$ be a difeomorphism with inverse $g=f^{-1}$, where $U$ and $V=f(U)$ are open subsets of $\R^n$, $n\ge 1$. 
Then the following holds:
\begin{enumerate}
	\item If  $\ell\in \N_\ast$ then 
\[
\|dg\|_{\ell,V}\,  \le \,  C_\ell\left(1+\|(df)^{-1}\|_{C^0(U)}^{3\ell}\right) \,   {\mathcal P}_U^{\ell}(df) 
\]
where  $C_\ell>0$ depends  only on $\ell$ and on the dimension $n$;
	\item  Let $0<\ell<1$ and let $V=f(U)$ be convex. Then 
\[
\|dg\|_{C^\ell(V)} \,  \le \, \|(df)^{-1}\|_{C^0(U)}^\ell  \|(df)^{-1}\|_{C^\ell(U)}\,  \le \,  \|(df)^{-1}\|_{C^0(U)}^{\ell+2} \,  \|df\|_{C^\ell(U)}  ;
\]
	\item  Let $\ell =m+\mu$, where $m\in\N_\ast$ and $0<\mu<1$.  Let $V=f(U)$ be convex. Then 
\[
\begin{array}{rcl}
\|dg\|_{C^\ell(V)}\,  
                              &\le& \,  C_\ell\left(1+\|(df)^{-1}\|_{C^0(U)}^{3\ell+2}\right) \,  \left( {\mathcal P}_U^{\ell}(df) + \|df\|_{C^{\mu}(U)}{\mathcal P}_U^{m}(df) \right) ,
\end{array}
\]
where  $C_\ell>0$ depends  only on $\ell$ and on the dimension $n$. 
\end{enumerate}
\end{Prop} 
{\em Proof}. To prove the statement one  needs a sort of  Lagrange inversion formula for  higher derivatives of  the Jacobian matrix $Dg\in C^\infty(U, M_n(\R))$. Hereafter   $M_n(\R)= M_{n,n}(\R)$ is the space  of real $n\times n$ matrices equipped with the corresponding sup-norm and we denote by $A\cdot B$  the product of  two matrices $A$ and $B$.  
Denote by $L_p: M_{n}(\R) \to M_{1,n}(\R)$ the linear operator which assigns to each matrix $A\in M_{n}(\R)$ its $p^{th}$-line and by $C_q: M_{n}(\R) \to M_{n,1}(\R)$ the linear operator which assigns to each matrix $A\in M_{n}(\R)$ its $q^{th}$-column. We identify $L_p$ and $C_q$ with the corresponding matrices in $M_{1,n}(\R)$ and  $M_{n,1}(\R)$, respectively. 
Since $(Dg)(f(x))= Df(x)^{-1}$, where $Df(x)^{-1}$ is the inverse of the matrix $Df(x)$, we get for any $1\le p,q\le n$ the equality
\begin{equation}
\left(\frac{\partial}{\partial y_q} g_p\right)(y)=  L_p\cdot Df(x)^{-1} \cdot C_q  \in \R \quad \mbox{as} \quad y=f(x).
\label{eq:derivative-Dg}
\end{equation}
Then differentiating the identity $Dg(y)=(Df)^{-1}(g(y))$ with respect to $y_q$ we obtain
\[
\left(\frac{\partial}{\partial y_q}  Dg\right)(y) = \sum_{p=1}^n \frac{\partial}{\partial x_p}(Df(x)^{-1}) \cdot (L_p \cdot Df(x)^{-1} \cdot C_q)I_n \quad \mbox{as}\  x=g(y). 
\]
Denote by ${\mathcal A}$ the set of automorphisms  of $M_n(\R)$ generated under composition by the identity map  and the automorphisms  of the form 
\[
A\mapsto (L_p\cdot A \cdot C_q) I_n, \quad  A\in M_n(\R).
\]  
Using  \eqref{eq:derivative-Dg} we obtain by induction  with respect to  $m\in \N$ the following relation
\begin{equation}
\partial_y^\alpha Dg(y) = \sum_{\gamma\in\Delta_\alpha} c_\gamma \prod_{1\le j \le m+1} {\mathcal L}_j ( \partial_x^{\gamma_j}(Df(x)^{-1})) \quad \mbox{as}\   x=g(y) \, , 
\label{eq:derivative-Dg-1}
\end{equation}
for any $\alpha\in \N^n$ with $|\alpha| =m$. The index set $\Delta_\alpha$ consists of all  $\gamma=(\gamma_1,\ldots,\gamma_{m+1})\in (\N^n)^{m+1}$ such that 
\[
|\gamma_1| + \cdots + |\gamma_{m+1}|=m=|\alpha|,
\] 
 $c_\gamma\in \R$ are universal constants, and ${\mathcal L}_j\in {\mathcal A}$. 
 Consider the derivatives of the inverse matrix $Df(x)^{-1}$ of $Df(x)$. 
One can easily show that
\begin{equation}
\frac{\partial}{\partial x_p}  (Df(x)^{-1}) = - Df(x)^{-1} \cdot \left(\frac{\partial}{\partial x_p}  Df(x)\right) \cdot  Df(x)^{-1} .
\label{eq:derivative-Df}
\end{equation}
 This equality  implies by induction that for any $0\neq \alpha \in\N^n$, 
\begin{equation}
\partial_x^\alpha  (Df(x)^{-1}) = \sum_{\beta\in\Delta_\alpha^1} c_\beta \prod_{1\le j \le |\alpha|} Df(x)^{-1} \cdot \partial_x^{\beta_j} Df(x) \cdot Df(x)^{-1} , 
\label{eq:derivative-Df-1}
\end{equation}
where $c_\beta\in \R$ are universal constants and the index set $\Delta_\alpha^1$ consists of all  \[
\beta=(\beta_1,\ldots,\beta_{|\alpha|})\in \N^n\times\cdots\times \N^n= (\N^n)^{|\alpha|}
\]
such that 
\[
|\beta_1| + \cdots + |\beta_{|\alpha|}|=|\alpha|.
\]
Now statement {\em 1.} follows easily from \eqref{eq:derivative-Dg-1} and \eqref{eq:derivative-Df-1}.

We are going to prove {\em 2} and   {\em 3} . 
Since $V$ is convex, we have in view of Lemma \ref{lemma:composition-mu}
\[
\|h\circ g\|_{C^\mu(V)} \le \|h\|_{C^\mu(U)}\|Dg\|_{C^0(V)}^\mu = \|h\|_{C^\mu(U)}\|(Df)^{-1}\|_{C^0(U)}^\mu 
\]
for any $h\in C^\mu(U)$ and $0<\mu<1$. Moreover,
\[
\|(Df)^{-1}\|_{C^\mu(U)}\le \|Df\|_{C^\mu(U)}\|(Df)^{-1}\|_{C^0(U)}^2.
\]
We have $Dg=(Df)^{-1}\circ g$, hence, taking $h=(Df)^{-1}$ we obtain {\em 2}. 
Statement    {\em 3} follows from   \eqref{eq:derivative-Dg-1} and \eqref{eq:derivative-Df-1} as in  the proof of Proposition \ref{prop:composition}. 
\finishproof

\noindent
{\em Proof of Lemma \ref{lemma:cut-off-omega}. }
One can take the convolution  $\psi^U_\varepsilon= \1_U \ast\chi_{\varepsilon}$,  where $\1_U$ is the characteristic function of $U$, $\chi_\varepsilon(x) = \varepsilon^{-n}\chi(x/\varepsilon)$ and $\chi\in C_0^\infty(\R^n)$ is a test-function such that 

\begin{center}
$\chi(x)>0$ if $|x|<1$, $\chi(x)=0$ if $|x|\ge1$ and $\displaystyle \int_{\R^n}\chi(x)dx=1$. 
\end{center}
More preciselly, we first define a smooth function $\widetilde \chi$ by  
$\widetilde \chi(x)=\exp(-(1-|x|^2)^{-1})$ for $|x|<1$ and by $\widetilde \chi(x)=0$ for $|x|\ge1$ and then we set $\chi(x)=\widetilde \chi(x)/\|\widetilde \chi\|_{L^1}$. 
For any $\ell\ge 0$ we have $\|\chi_\varepsilon\|_{\ell;\varepsilon}=\varepsilon^{-n}\|\chi\|_{\ell}$. Set $\phi:= \psi^B_1$ and $\phi_\varepsilon(x) = \varepsilon^{-n}\phi(x/\varepsilon)$, where $B=B^n(0,2)$. Then $\phi(x)=1$ for $|x|\le 1$, hence, $\phi(x/\varepsilon)=1$ for $x$ in  the support of $\chi_\varepsilon$ and  we easily obtain the inequality 
\[
|\partial_x^\alpha\chi_\varepsilon(x-y)-\partial_z^\alpha\chi_\varepsilon(z-y)|\,  \le\,  \big(\phi((x-y)/\varepsilon)+\phi((z-y)/\varepsilon)\big) |\partial_x^\alpha\chi_\varepsilon(x-y)-\partial_z^\alpha\chi_\varepsilon(z-y))|
\]
for any $x,y,z\in\R^n$ and $0<\varepsilon\le 1$. 
This implies
\[
\|\psi^U_\varepsilon\|_{\ell;\varepsilon} \le 2\varepsilon^n\|\1_U\ast\phi_{\varepsilon}\|_{L^1}\|\chi_\varepsilon\|_{\ell;\varepsilon} \le 2\|\phi_{\varepsilon}\|_{L^1}\|\chi\|_{\ell} = 2\|\phi\|_{L^1}\|\chi\|_{\ell}=C_\ell
\]
which proves the Lemma. 
\finishproof

\subsubsection{ Higher order H\"{o}lder estimates  and Interpolation inequalities.}\label{subsec:interpolation}
The above estimates can be simplified considerably  if the domain of definition of the functions  is the whole space or an open convex bounded set. We set $\A:= \T^n\times D$, where $D$ is an open set in $\R^d$. We shall use the convention $\A:= \T^n$ if $d=0$ and $\A:= D$ if $n=0$.

Using the interpolation inequalities as in \cite{L-O}, Proposition 5.5,  one obtains
\begin{Prop}\label{prop:composition-holder-interpolation} Let  $f\in C^\infty(\A_1,\A_2)$ and $g\in C^\infty(\A_2,\R)$, where $\A_1=\T^{n_1}\times \R^{d_1}$,  $\A_2=\T^{n_2}\times D_2$  and 
	$D_2\subset \R^{d_2}$ is an open set.   Then the following holds
\begin{enumerate}
	\item
For any  $\ell\ge 1$, 
\[
\begin{array}{lcrr}
\big|g\circ f \big|_{\ell,\A_1}   \,  \le\,  C_{\ell} (1+\|d f \|_{0,\A_1}^{\ell})  \\[0.3cm] 
\displaystyle \times 
\sum_{m=1}^{[\ell]}\left(|g|_{m,\A_2}\|d f\|_{\ell-m,\A_1} + |g|_{m+{\{\ell\}},\A_2} \|d f\|_{[\ell]-m,\A_1}\right)
\end{array}
\]
where  $C_{\ell}>0$ depends only on $\ell\ge 1$ and on the dimensions $n_j$ and $d_j$. 
\item Let $\A_j=\T^{n_j}\times \R^{d_j}$, $j=1,2$.  Then
\begin{equation}\label{eq:composition-interpolation}
\begin{array}{lcrr}
\big|g\circ f \big|_{\ell,\A_1}   \,  \le\,  C_{\ell}  (1+\|d f \|_{0,\A_1}^{\ell})  \\[0.3cm] 
\displaystyle \times 
\left(\|g\|_{\ell,\A_2}\| d f\|_{0,\A_1} + \|g\|_{1,\A_2}\|  d f\|_{\ell-1,\A_1} \right)
\end{array}
\end{equation}
where  $C_{\ell}>0$ depends only on $\ell\ge 1$ and on the dimensions $n_j$ and $d_j$.
\item Let  $D_j$, $j=1,2$, be  open convex subsets of $\R^{d_j}$ and   $\overline \A_j=\T^{n_j}\times \overline D_j$. Then \eqref{eq:composition-interpolation} holds for any 
$f\in C^\infty(\overline \A_1,\overline \A_2)$ and $g\in C^\infty(\overline \A_2,\R)$ with a constant $C_\ell$ depending only on $\ell\ge 1$ and the dimensions $n_j$ and $d_j$.  
\end{enumerate}
\end{Prop}
{\em Proof.} To prove the first statement we make use of  Proposition \ref{prop:composition} and of the interpolation inequalities. Consider a typical term of \eqref{eq:P-l-m} given by
\[
A:=  |f|_{k_1} \cdots  |f|_{k_{j-1}} |f|_{k_j +\{\ell\}} |f|_{k_{j+1}}  \cdots |f|_{k_m}
\]
where $k_1+ \cdots +k_m= [\ell]$ and $k_p \ge 1$,  for any $p$. Set $r=0$, $s=s_p=k_p+ \delta_p-1$ and $t=\ell-m$, where $\delta_p=0$ if $p\neq j$ and $\delta_j= \{\ell \}$.    By means of the interpolation inequalities we get
\[
|f|_{k_p+ \delta_p} \, \le \, \|d f\|_{s} \, \le \, c_{\ell}\,  \|d f\|_{0}^{\frac{t-s}{t}}\, \|d f\|_{t}^{\frac{s}{t}}
\]
which  implies 
\[
A \, \le \, c_{\ell}\,  \|d f\|_{0}^{m-1}\, \|d f\|_{\ell-m} , 
\]
where $c_{\ell} >0$ depends only on $\ell$. Hence,
\begin{equation}\label{eq:estimate-P}
 {\mathcal P}^{\ell,m}_{U}(f)\, \le \, C_{\ell}\,  \sum_{m=1}^{[\ell]} \|d f\|_{0}^{m-1}\, \|d f\|_{\ell-m}.
 \end{equation}
 Notice also that $|f|_{1,\A_1}^{\{\ell\}}\le 1 + |df|_{0}$. 
Using Proposition \ref{prop:composition} one obtains {\em 1}.

To prove {\em 2}  one uses the  interpolation inequalities with respect to both functions $f$ and $g$. Namely, 
given $0<s<t$ and  $u,v\in C^t$ one obtains
\begin{equation}\label{eq:u-v-Holder}
\|u\|_{s}\|v\|_{t-s}\, \le\,  c_t \|u\|_{0}^{\frac{t-s}{t}}\|u\|_{t}^{\frac{s}{t}}\|v\|_{0}^{\frac{s}{t}}\|v\|_{t}^{\frac{t-s}{t}}
\, < \,  c_t \big(\|u\|_{t}\|v\|_{0}+\|v\|_{t}\|u\|_{0}\big)
\end{equation}
by means of \eqref{eq:interpolation} and Young's inequality 
\[
xy\, \le\,  \frac{1}{p}x^p + \frac{1}{q}y^q \, <\,  x^p + y^q ,
\] 
where $x$ and $y$ are non-negative  and $p$ and $q$  are  positive numbers such that  $\frac{1}{p} + \frac{1}{q} =1$. Putting $s=m$, $t=\ell$, $u=g$, $v=dg$, and then  $s=m-1+\{\ell\}$, $t=\ell-1$, $u=dg$, $v=dg$, and using {\em 1}, we obtain  {\em 2}. 
The inequality in {\em 2} has been proven for more general domains in (\cite{L-O}, Proposition 5.5)  but the constants there depend on the domains.
Statement {\em 3} follows from {\em 2} and Remark \ref{rem:interpolation}. 
\finishproof

In the same way we obtain
\begin{Prop}\label{prop:inverse1}
Let $f\in C^\infty(\overline\A_1,\overline\A_2)$ be a difeomorphism with inverse $g=f^{-1}\in C^\infty(\overline\A_2,\overline\A_1)$, where $\overline\A_j=\T^{n}\times \overline D_j$ and 
	$D_j\subseteq \R^{d}$, $j=1,2$,  are open   sets. Then the following holds.
 \begin{enumerate}
   \item Let $D_1$ be convex. Then for any  positive integer $\ell\in\N_\ast$,  
\begin{equation} \label{eq:inverse-interpolation}
\|dg\|_{\ell,\A_2}  \le C_\ell\left(1+\|(df)^{-1}\|_{0,\A_1}^{3\ell +1}\right)\,   \sum_{m=1}^{[\ell]} \|d f\|_{0,\A_1}^{m-1}\, \|d f\|_{\ell-m,\A_1}
\end{equation}
where  $C_\ell$  depends  only on $\ell$ and on the dimensions $n$ and $d$. 
   \item Let  both $D_1$ and $D_2$ be convex. Then \eqref{eq:inverse-interpolation} holds for any $\ell\ge 1$.  
\end{enumerate}
\end{Prop}
{\em Proof.}
 The inequality \eqref{eq:inverse-interpolation} follows for $\ell\in\N_\ast$ from Proposition \ref{prop:inverse} and Remark \ref{rem:interpolation} using \eqref{eq:estimate-P}  as in the proof of Proposition \ref{prop:composition-holder-interpolation}. To prove it for any  $\ell\ge 1$ we use \eqref{eq:u-v-Holder} as well. 
 \finishproof

Given $n,d\in \N$ we set $\A= \T^n\times \R^d$. 
Using the interpolation inequalities and Proposition \ref{prop:inverse} we obtain
\begin{Prop}\label{prop:inverse-holder1}
	Let $u={\rm id} + \phi$, where $\phi\in C^\infty(\A,\A)$, $\A= \T^n\times \R^d$,  and
	\[
	(n+d)\| \phi\|_{1}   \le \varepsilon_0 <1. 
	\]
	Then $u$ is a diffeomorphism homotope to the identity with inverse $u^{-1}= {\rm id} + \psi$, where $\psi\in C^\infty(\A,\A)$ and   for any  $\ell\ge 0$ we have 
	\[
	\|\psi\|_{\ell} \,  \le\,  C_{\ell}\,  \| \phi\|_{\ell}
	\]
	where  $C_{\ell} = C_{\ell} (\varepsilon_0,n,d)>0$ depends only on $\ell$, $\varepsilon_0$,  and on the dimensions $n$ and $d$. If  ${\rm supp\, }\phi \subset \T^n\times K$ then ${\rm supp\, }\psi \subset \T^n\times K$ as well. Moreover, if the map $[0,\delta]\ni t\to \phi_t\in C^\infty(\A,\A)$ is $C^k$ then so is the map $t\to \psi_t$.
	The same holds  if $\A= \T^n$ ($d=0)$ or $\A= \R^d$ ($n=0$). 
\end{Prop}
{\em Proof.} Notice that $\| d \phi\|_{0} < (n+d)\| \phi\|_{1} \le  \varepsilon_0 <1$. 
The inverse function theorem implies that $u$ is a diffeomorphism with inverse $u^{-1}= {\rm id} + \psi$, where $\psi\in C^\infty(\A,\A)$. Moreover, 
\begin{equation}\label{eq:identity-psi-phi}
\psi=-\phi\circ ({\rm id} + \psi)^{-1},
\end{equation}
which implies $\|\psi\|_{0}  = \| \phi\|_{0}$. If  $0\le \ell <1$, then 
\[
\|\psi\|_\ell
\le \| \phi\|_{0} + \|\phi\|_{C^\ell}\|({\rm Id} + d\phi)^{-1}\|_{0}^\ell 
\le C(1-\varepsilon_0)^{-1} \|\phi\|_{\ell}. 
\]
The estimate of $\psi$ in the $C^\ell$ norms with $\ell \ge 1$ follows from \eqref{eq:identity-psi-phi} using Proposition \ref{prop:composition-holder-interpolation}  and Proposition \ref{prop:inverse1}. 
To prove the assertion about the support notice that $u^{-1}={\rm id}$ on $\T^n\times (\R^d\setminus K)$. 
\finishproof

\subsubsection{ Weighted H\"older norms and interpolation inequalities.}\label{subsec:weighted-interpolation}
Given $0<\kappa\le 1$, $\ell\ge 0$,  and $f\in C^\ell(\T^n\times D)$, where $D\subset \R^d$ is an open set  we define the corresponding weighted $C^\ell$ norm by
\[
\|f\|_{\ell,\T^n\times D;\kappa}\, := \, \|f\circ \sigma_\kappa\|_{\ell,\sigma_\kappa^{-1}(\T^n\times D)}, 
\]
where $\sigma_\kappa(\theta,r)=(\theta,\kappa r)$. 
We set as well 
$$|f|_{\ell,\T^n\times D;\kappa}\, := \, |f\circ \sigma_\kappa|_{\ell,\sigma_\kappa^{-1}(\T^n\times D)}.$$
In particulat, if $\ell\in\N$, then 
$$\|f\|_{m,\T^n\times D;\kappa}=\sup_{0\le m\le \ell}|f|_{m,\T^n\times D;\kappa},$$
where
\[
|f|_{m,\T^n\times D;\kappa}\, := \, \sup_{|\alpha|+ |\beta|= m}\, \|\partial_\theta^\alpha (\kappa\partial_r)^\beta f\|_{C^0(\T^n\times D)}.
\]
Applying \eqref{eq:interpolation} to $f=u\circ \sigma_\kappa$ with $0<\kappa\le 1$ one gets the interpolating inequalities for $\|u\|_{t;\kappa}:=\|u\|_{t,\A;\kappa}$, where $\A=\T^{n}\times \R^{d}$.

We list below several estimates which follow directly from Proposition \ref{prop:composition-holder-interpolation}. Set
\[
|\!|\!| u |\!|\!|_{\ell,D;\kappa}\, =\, \sup_{0\le m\le \ell}\|u\|_{\ell-m,D;\kappa},   
\]
where $m$ are integers. 
\begin{Prop}\label{prop:composition-holder1} Fix $0<\kappa\le 1$. 
\begin{enumerate}
	\item  Let  $u\in C^\infty(\A_1,\A_2)$ and $v\in C^\infty(\A_2,\R)$, where $\A_1=\T^{n_1}\times \R^{d_1}$, $\A_2=\T^{n_2}\times D_2$  and 
	$D_2\subset \R^{d_2}$ is an open set.   Then for any  $\ell\ge 1$, 
\[
\begin{array}{lcrr}
\big|v\circ u \big|_{\ell,\A_1;\kappa}   
\le   C_{[\ell]} \left(1+\|d(\sigma_\kappa^{-1}\circ u\circ \sigma_\kappa)\|_{C^0}^{\ell-1}\right)
\\[0.3cm]
\displaystyle \times \sum_{m=1}^{[\ell]}\left(\|v\|_{m,\A_2;\kappa}\|d(\sigma_\kappa^{-1}\circ u\circ \sigma_\kappa)\|_{C^{\ell-m}(\A_1)} + \|v\|_{m+{\{\ell\}},\A_2;\kappa} \|d(\sigma_\kappa^{-1}\circ u\circ \sigma_\kappa)\|_{C^{[\ell]-m}(\A_1)}\right)
\end{array}
\]
where  $C_{\ell}>0$ depends only on $\ell$ and on the dimensions $n_j$ and $d_j$.
\item Let $\A_1 =\T^{n_2}\times D_1$, where $D_1$ is an open convex subset of $\R^{d_1}$, $\A_2 =\T^{n_2}\times \R^{d_2}$, and $u\in C^\infty(\overline\A_1,\A_2)$. Then
\[
\begin{array}{lcrr}
|v\circ u|_{\ell,\A_1;\kappa}\,  \le \,  C_\ell \left(1+  \|d (\sigma_\kappa^{-1} \circ u\circ \sigma_\kappa)\|_{C^0}^{\ell-1}\right)
\\[0.3cm]
\displaystyle \times \big(\|v\|_{\ell;\kappa}\|d(\sigma_\kappa^{-1}\circ u\circ \sigma_\kappa)\|_{C^0} + \|v\|_{1;\kappa}\| d(\sigma_\kappa^{-1}\circ u\circ \sigma_\kappa)\|_{C^{\ell-1}} \big)
\end{array}
\]
for any $\ell\ge 1$, where  $C_{\ell}>0$ depends only on $\ell$ and on the dimensions $n_j$ and $d_j$.
In particular, if  $\A_2 =\T^{n_2}$ ($d_2=0$) then 
\[
\|v\circ u\|_{\ell;\kappa}\,  \le\,  \|v\circ u\|_{C^0}+C_\ell \left(1+  \|u\|_{1;\kappa}^{\ell-1}\right)\big(\|v\|_{\ell;\kappa}\|u\|_{1;\kappa} + \|v\|_{1;\kappa}\| u\|_{\ell;\kappa} \big).
\]
\item Let $D_j\subset \R^{d_j}$, $j=1,2$, be open sets in $\R^{d_j}$ and let $D_1$ be  convex. Then for any $u\in C^\infty (\overline{D}_1, \R^{d_2})$ with $u(D_1)\subset D_2$ and $v\in C^\infty(\A_2)$ with $\A_2=\T^{n_2}\times D_2$ we have
\[
\begin{array}{lcrr}
|v\circ u |_{\ell,D_1;\kappa} \,  \le \,    
C_{\ell} \left(1+\|d u\|_{C^0(D_1)}^{\ell-1}\right)
\\[0.3cm]
\displaystyle \times \, 
\sum_{m=1}^{[\ell]}\left(|v|_{m,\A_2;\kappa}\,  (1+\|d u\|_{\ell-m,D_1; \kappa})  + |v|_{m+{\{\ell\}},\A_2;\kappa}\,  (1+\|d u\|_{[\ell]-m,D_1; \kappa})\right) 
\\[0.3cm]
\displaystyle
\le C_{\ell} \left(1+\|d u\|_{C^0(D_1)}^{\ell-1}\right) \,  \big(1+\|d u\|_{C^\ell(D_1)} \big)\, |\!|\!| v |\!|\!| _{\ell,\A_2;\kappa}\, 
\end{array}
\]
for any $\ell\ge 1$, where the constant $C_{[\ell]}>0$ depends only on $[\ell]$ and on the dimensions $d_j$ and $n_2$. 
\end{enumerate}
\end{Prop}
{\em Proof}. To prove {\em 1} we put $f=\sigma_\kappa^{-1}\circ u \circ \sigma_\kappa$ and $g=v\circ \sigma_\kappa$ in Proposition \ref{prop:composition-holder-interpolation}. The statement {\em 2} follows from Proposition \ref{prop:composition-holder-interpolation}, {\em 3}. 
To prove {\em 3}  we use Proposition \ref{prop:composition}, {\em 3},  and apply  Remark \ref{rem:interpolation} to $u\in C^\infty(\overline D_1)$  using the interpolation inequalities for the extension of $u$ as in the proof of Proposition \ref{prop:composition-holder-interpolation}. To prove the second inequality in {\em 3} notice   that $\|d u\|_{\mu,D_1; \kappa} \le \|d u\|_{\ell,D_1; \kappa}$ for $\mu\le\ell$ since $D_1$ is convex. 
\finishproof
\begin{Remark}\label{rem:anisotrop}
The estimates in  Proposition \ref{prop:composition-holder1} hold when $\kappa=(\kappa',\kappa'')$ with $0<\kappa',\kappa''\le 1$,  $D_j =  D_j'\times D_j'$, $\sigma_\kappa(\theta,x',x'')= (\theta,\kappa'x',\kappa''x'')$ and $|u|_{\ell,\A_1;\kappa}:= |u\circ \sigma_\kappa|_{C^\ell(\sigma_\kappa^{-1}(\A_1))}$. 
\end{Remark}

\subsubsection{ Symplectic transformations and generating functions.}\label{subsec:symplectic}

Consider a $C^1$ family of  exact symplectic maps $W_t: \A\to \A$, $t\in [0,\delta]$, where $\A:= \T^d\times \R^d$. Suppose that $W_t-{\rm id}$ is compactly supported for any $t$. We are looking for a $C^1$ family of generating functions 
\begin{equation}\label{eq:generating-function1}
\widetilde G_t(\varphi,r)=\langle \varphi, r\rangle - G_t(\mathrm{pr}(\varphi),r), \quad (\varphi,r)\in \R^{d}\times \R^{d},
\end{equation}
of $W_t$ such that the function 
 $G_t\in C^\infty(\T^d\times \R^d)$ is  compactly supported with respect to $r$ and 
\begin{equation}\label{eq:generating-function2}
W_t\left(\nabla_r \widetilde G_t(\varphi,r),r\right)=\left(\varphi, \nabla_\varphi \widetilde G_t(\varphi,r)\right),\quad (\varphi,r)\in \A
\end{equation}
(see Definition \ref{def:generating-function}). Slightly abusing the notations, we will identify  below a $2\pi$-periodic function  with the corresponding functions on $\T^d$. 
Given a smooth function $G$ in $\A:=\T^d\times D$,   we denote  by 
\[
{\rm sgrad\, }G(\theta,r):= (\nabla_r G(\theta,r), -\nabla_\theta  G(\theta,r))
\]
its symplectic gradient. Notice that $\|\sigma_\kappa^{-1}\,  {\rm sgrad\, }G_t\|_{\ell,\A;\kappa} \le \kappa^{-1} \| G_t\|_{\ell+1,\A;\kappa}$. 
\begin{Lemma}\label{lemma:generating-functions} \emph{1.} Let $[0,\delta]\ni t\to W_t: \A\to \A$ be a $C^1$ family of exact symplectic mappings. Suppose that  
\begin{equation}
\label{eq:small-estimate-W}
2d	\|\sigma_\kappa^{-1}(W_t-{\rm id})\|_{1,\A;\kappa}\le  \varepsilon_0<1
\end{equation}
for $t\in [0,\delta]$.
	Then there exists a $C^1$ family of generating functions $\widetilde G_t$ of $W_t$ given by \eqref{eq:generating-function1} such that for any $\ell\ge 0$ the following estimate hold true
\begin{equation}\label{eq:estimate-G}
\|\sigma_\kappa^{-1}\,  {\rm sgrad\, }G_t\|_{\ell,\A;\kappa} \le C_\ell  \|\sigma_\kappa^{-1} (W_t-{\rm id})\|_{\ell,\A;\kappa}
\end{equation}
for $t\in [0,\delta]$, where  $C_{\ell} = C_{\ell} (\varepsilon_0,d)>0$  depends only on $\ell$, $\varepsilon_0$ and $d$. Moreover, the relation ${\rm supp\, }(W_t -{\rm id}) \subset \T^d \times K$ implies  ${\rm supp\, } ({\rm sgrad\, } G_t) \subset \T^d \times K$ as well. \\
\emph{2.}  Conversely, let $G_t$ be a $C^1$ family of functions  such that  
\begin{equation}
\label{eq:small-estimate-G}
 	2d \|\sigma_\kappa^{-1}\,  {\rm sgrad\, }G_t\|_{1,\A;\kappa} \le \varepsilon_0<1 
\end{equation}
for $t\in [0,\delta]$.
Then $\widetilde G_t$ given by \eqref{eq:generating-function1} is  a $C^1$ family of generating functions of symplectic maps $W_t:\A\to\A$ and for any  $\ell\ge 0$ we have
\begin{equation}\label{eq:estimate-W}
\displaystyle\begin{array}{lcrr}
\|\sigma_\kappa^{-1}(W_t-{\rm id})\|_{\ell,\A;\kappa} + \| \sigma_\kappa^{-1}(W_t^{-1}-{\rm id})\|_{\ell,\A;\kappa}\\[0.3cm] 
\displaystyle \le C_\ell  \|\sigma_\kappa^{-1} {\rm sgrad\, }G_t \|_{\ell,\A;\kappa}
\end{array}
\end{equation}
where $t\in [0,\delta]$ and  $C_{\ell} = C_{\ell} (\varepsilon_0,d)>0$  depends only on $\ell$, $\varepsilon_0$ and $d$.  Moreover, if ${\rm supp\, } ({\rm sgrad\, }G_t) \subset \T^d \times K$ then ${\rm supp\, }(W_t -{\rm id}) \subset \T^d \times K$ as well. 
\end{Lemma}
{\em Proof.} \emph{1.}   Set $W_t=(U_t,V_t): \A \to \A$. It follows from \eqref{eq:small-estimate-W} that the map $\theta\to U_t(\theta,r)-\theta$ can be identified with a  $2\pi$-periodic vector function on $\R^d$. Consider  the map 
$$f_t={\rm id} + g_t:\A\to\A, \quad  \mbox{where}\quad   
g_t(\theta,r)= (U_t(\theta,r)-\theta,0). $$
By \eqref{eq:small-estimate-W} one obtains 
\[
2d	\, \|g_t\|_{1,\A;\kappa} <2 d  \|\sigma_\kappa^{-1}(W_t-{\rm id})\|_{1,\A;\kappa}\le  \varepsilon_0 <1
\]
	for any $t\in [0,\delta]$. The inverse  function theorem  (Proposition \ref{prop:inverse-holder1}) implies that $f_t:\A \to \A$
is a diffeomorphism homotope to the identity. In particular, the equation  $\varphi=U_t(\theta,r)$ has a unique smooth solution  $\theta=\varphi+\phi_t(\varphi,r)$, where $\phi_t$ can be identified with a  $2\pi$-periodic with respect to $\varphi\in \R^d$ function and the map 
$[0,\delta]\ni t \to \phi_t\in C^\infty(\A,\A)$ is $C^1$. Then 
$$f_t^{-1}={\rm id} + h_t , \quad  \mbox{where}\quad  h_t= (\phi_t,0).$$
Proposition \ref{prop:inverse-holder1}    applied to 
$$\sigma_{\kappa}^{-1} \circ f_t \circ \sigma_{\kappa}= {\rm id}  + g_t \circ \sigma_{\kappa} \quad \mbox{and}\quad 
\sigma_{\kappa}^{-1}\circ  f_t^{-1}\circ  \sigma_{\kappa}  = {\rm id}  + h_t \circ \sigma_{\kappa}$$
yields
\begin{equation}\label{eq:estimate-phi-W}
\|\phi_t\|_{\ell;\kappa}\le  C_\ell  \|g_t\|_{\ell;\kappa} \le C_\ell \|\sigma_\kappa^{-1}(W_t-{\rm id})\|_{\ell,\A;\kappa} 
\end{equation}
where  $C_\ell>0$ depends only on $\ell$, $\varepsilon_0$  and $d$. On the other hand, the map $W_t$ is exact symplectic and close to the identity and there exists a $C^1$-family of generating functions  $\widetilde G_t$ such that $G_t$ is compactly supported and 
\[
\nabla_r G_t(\varphi,r) = -\phi_t(\varphi,r)\, ,\quad \nabla_\varphi G_t(\varphi,r) = V_t(\varphi+\phi_t(\varphi,r),r)-r.
\]
We are going to prove \eqref{eq:estimate-G}. The estimate of $\nabla_r G_t$ follows from \eqref{eq:estimate-phi-W}. To prove the estimate of $\kappa^{-1}\nabla_\varphi G_t$ we write
\[
\nabla_\varphi G_t(\varphi,r) = (V_t(\varphi,r)-r) + \int_0^1 d_\theta V_t(\varphi+s\phi_t(\varphi,r),r)\,  \phi_t(\varphi,r) \, ds\, ,
\]
where  $d_\theta$ is the partial differential with respect to the first variables $\theta$. 
Notice that $\|\phi_t\|_{1;\kappa} \le C_1 \varepsilon_0$ in view  of \eqref{eq:small-estimate-W} and \eqref{eq:estimate-phi-W}. 
Then using  \eqref{eq:interpolation-leibnitz}, Proposition \ref{prop:composition-holder1}, {\em 2},  and \eqref{eq:small-estimate-W} we complete the proof of \eqref{eq:estimate-G}.
Suppose now that ${\rm supp\, }(W_t -{\rm id}) \subset \T^d \times K$. Then ${\rm supp\, }(f_t -{\rm id}) \subset \T^d \times K$ which implies that ${\rm supp\, }(f_t^{-1} -{\rm id}) \subset \T^d \times K$. Hence, $\phi_t(\varphi,r)=0$ for $\varphi\in\T^n$ and $r\notin K$ and 
 ${\rm supp\, } ({\rm sgrad\, }  G_t) \subset \T^d \times K$. 

\emph{2.} In the same way we prove  the second part of the Lemma. Suppose that \eqref{eq:small-estimate-G} holds. 
Using the inverse  function theorem  given by Proposition \ref{prop:inverse-holder1} one  solves as above  the equation 
\[
\theta=\varphi -  \nabla_rG_t(\varphi,r)
\]
 with respect to $\varphi\in \T^{d}$. The corresponding solution has the form $\varphi= \theta  + \psi_t(\theta,r)$ and 
\[
W_t(\theta,r) = (\theta  + \psi_t(\theta,r), r + \nabla_\theta G_t(\theta  + \psi_t(\theta,r ),r)). 
\]
Moreover, Proposition \ref{prop:inverse-holder1}   yields as above the estimate
\begin{equation}\label{eq:implicit-psi}
\|\psi_t\|_{\ell,\A;\kappa}
 \le   C_\ell \|\sigma_\kappa^{-1}{\rm sgrad\, }G_t \|_{\ell,\A;\kappa}
\end{equation}
for any $t\in [0,\delta]$, where $\ell\ge 0$ and $C_{\ell}>0$ depends only on $\ell$, $\varepsilon_0$   and $n$.  Then using  \eqref{eq:generating-function2} and Proposition \ref{prop:composition-holder1}  we estimate of  $W_t-{\rm id}$. 
We get the same estimates for $(W_t)^{-1}-{\rm id}$, where
 \[
(W_t)^{-1}(\varphi,r-\nabla_{\varphi}G_t( \varphi,r))=( \varphi-\nabla_{r}G_t( \varphi,r), r).
\]
To this end we first solve the equation $r-\nabla_{\varphi}G_t( \varphi,r)=I$ with respect to $r$ and then we proceed as above. 
 \finishproof

The estimates \eqref{eq:estimate-G} and \eqref{eq:estimate-W} are still valid if we add additional parameters $s\in \T^p$ and $\omega\in \R^q$. 
Set $\A := \A_1\times \A_2$, where $\A_1= \T^{d}\times \R^{d}$ and $\A_2=\T^{p}\times \R^{q}$. 
Given $\mu=(\varrho,\kappa)$ with $0<\varrho,\kappa\le 1$ and $f\in C^\ell(\A_1\times \A_2)$ we set $\|f\|_{\ell, \A;\mu} = \|f\circ\sigma_{\mu}\|_{C^\ell(\A)}$, where $\sigma_{\mu}(\theta,r;s,\omega)=(\theta,\varrho r; s, \kappa \omega)$. We consider the symplectic gradient of the function $(\theta,r) \to G_t(\theta,r; s,\omega)$ for $(s,\omega)$ fixed. Following the proof of  Lemma \ref{lemma:generating-functions} we obtain 
\begin{Lemma}\label{rem:with-parameters}
Suppose that the map $[0,\delta]\ni t\to G_t\in C^\infty(\A,\R)$ is $C^1$ and
\[
 (2d+p+q)	\, \|\sigma_{\varrho}^{-1}({\rm sgrad\, }G_t(\cdot;s,\omega)-{\rm id}_{\A_1}(\cdot))\|_{1,\A_1;\varrho}\, \le\,  \varepsilon_0\, <\, 1,
\]
for $t\in [0,\delta]$  and $(s,\omega)\in \A_2$, where   ${\rm id}_{\A_1}$ is the identity map on $\A_1$. Then for any $(s,\omega)\in\A_2$ fixed, the function $(\theta,r)\to \widetilde G_t(\theta,r;s,\omega) =\langle \theta,r\rangle -G_t(\theta,r; s,\omega)$ is a generating function of an exact symplectic map $W_t(\cdot;s,\omega)$ in $\A_1$, the map $[0,\delta]\ni t\to W_t\in C^\infty(\A,\A_1)$ is $C^1$ and  
\begin{equation}\label{eq:estimate-W-s}
\| \sigma_{\varrho}^{-1}(W_t-{\rm id}_{\A_1})\|_{\ell, \A;\mu} + \| \sigma_{\varrho}^{-1}(W_t^{-1}-{\rm id}_{\A_1})\|_{\ell, \A;\mu}
 \le C_\ell \|\sigma_{\varrho}^{-1} {\rm sgrad\, }G_t \|_{\ell, \A;\mu} 
\end{equation}
where   $W_t^{-1}(\cdot;s,\omega)$ is the inverse of $W_t(\cdot;s,\omega)$ in $\A_1$ with $(s,\omega)\in \A_2$ fixed and  $C_\ell>0$ depends only on $\ell$, $\varepsilon_0$ and on the dimensions $d$,$p$ and $q$.  
\end{Lemma}

Let $D\subset \R^{d}$ be an open set and $\A=\T^{d}\times D$. 
Consider a   function
$\widetilde G\in C^\infty( \A,  \A)$   of the form 
\[
\widetilde G(\theta,r)=\langle \theta,r\rangle -K(r)-G(\theta,r), 
\]
Recall that the map $Q$ defined by  $Q(\theta,r)= (\theta +\nabla K(r),r)$ is a symplectic map with generating function $(\theta,r)\to \langle \theta,r\rangle -K(r)$. 
\begin{Lemma}\label{Lemma:decomposition}
Let  ${\rm supp}\, G \subset \T^{n-1}\times F$, where  $F\subset D$ is a compact and let 
 $G$ satisfy \eqref{eq:small-estimate-G}. Then the function $(\theta,r)\to \langle \theta,r\rangle -G(\theta,r) $ is a generating functions of a  symplectic transformation $W: \A \to \A$,  $\widetilde G$ is a  generating functions of a  symplectic transformations $P:\A \to \A$, the support of $W- {\rm id}$ is contained in $\T^{n-1}\times F$ and $P=W\circ Q$. 
\end{Lemma}
{\em Proof}. The assertions about $W$ follow from Lemma \ref{lemma:generating-functions}. As in the proof of Lemma \ref{lemma:generating-functions}  one obtains  that the map $\theta \to \theta - \nabla_r G(\theta,r)$ is a diffeomorphism of $\T^{d}$ homothope to the identity mapping. Then comparing the identities 
\[
P(\theta - \nabla_r K(r)-\nabla_r G(\theta,r),r) = (\theta, r-\nabla_\theta G(\theta,r))
\]
and 
\[
(W\circ Q) (\theta - \nabla_r K(r)-\nabla_r G(\theta,r),r) =W(\theta -\nabla_r G(\theta,r),r) = (\theta, r-\nabla_\theta G(\theta,r))
\]
we obtain the relation $P=W\circ Q$. \finishproof

\section{Appendix.}

\subsection{Invariant characterization of Liouville billiards}\label{sec:Appendix B}

\noindent Here we prove the following invariant characterization of Liouville billiard tables defined 
in \cite[Sec. 2]{PT1}.

\begin{Th}\label{th:invariant_definition}
Let $(X,g)$ be a smooth oriented compact and connected Riemannian manifold of dimension two with connected boundary 
$\Ga\equiv\partial X$. Assume that
\begin{itemize}
\item[$(a)$] There exists a smooth quadratic in velocities integral of the geodesic flow $I : TX\to \R$ that is 
invariant with respect to the reflection at the boundary $TM|_{\Ga}\to TM|_{\Ga}$, $\xi\mapsto \xi - 2g(\nu,\xi)$, 
where $\nu$ is the outward unit normal to $\Ga$. In addition, we assume that the metric $g$ does {\em not} allow global 
Killing  symmetries;
\item[$(b)$] There is no point $x_0\in\Ga$ and a constant $c\in\R$ such that $g_{x_0}(\xi,\xi)= c I_{x_0}(\xi,\xi)$
for any  $\xi\in T_{x_0}X$.
\end{itemize}
Then $(X,g)$ is isometric to a Liouville billiard table.\footnote{In particular, $X$ is diffeomorphic to the
unit disk ${\mathbb D}^2$ in $\R^2$.}
Conversely, any Liouville billiard table satisfies the properties stated above.
\end{Th}
\begin{Remark}
The assumption that $g$ does not allow global Killing symmetries is needed for excluding the case when $(X,g)$
is a surface of revolution. Condition (b) can be replace by a similar condition but can {\em not} be avoided.
One can easily see this by considering the billiard table on the surface
of the ellipsoid $\big\{\frac{x^2}{a^2}+\frac{y^2}{b^2}+\frac{z^2}{c^2}=1\big\}$, $0<a<b<c$, defined by the 
condition $y\ge 0$. This billiard table is completely integrable but it is {\em not} a Liouville billiard table.
Its boundary is the geodesic that corresponds to the intersection of the coordinate plane
$O_{xz}$ with the ellipsoid. In particular, this curve is {\em not} locally geodesically convex and it contains 
the four umbilics of the ellipsoid. One can easily see that the billiard table defined this way satisfies all 
conditions of Theorem \ref{th:invariant_definition} except $(b)$.
Condition $(b)$ is also needed to ensure that the integral is non-trivial, i.e. $I\not\equiv c g$ where
$c$ is a real constant.
\end{Remark}
\noindent As a consequence of Theorem \ref{th:invariant_definition} we see that there exists a double
covering map with two branched points,
\[
\tau : C\to X,
\]
where $C$ denotes the cylinder ${(\R/\Z)}\times [-N, N]$, $N>0$, coordinatized by the variables $x$ and $y$
respectively, so that the metric $\tau^*(g)$ and the integral $\tau^*(I)$ have the following form on $C$,
\begin{eqnarray}\label{eq:the_metric}
dg^2&=&\big(f(x)-q(y)\big)(dx^2+dy^2)\\
dI^2&=&\alpha\,dF^2+\beta\,dg^2\nonumber
\end{eqnarray}
where $\alpha\ne 0$ and $\beta$ are real constants and
\begin{equation}\label{eq:the _integral}
dF^2:=\big(f(x)-q(y)\big)\big(q(y)\,dx^2+f(x)\,dy^2\big)\,.
\end{equation}
In other words, the integral $dI^2$ belongs to the pencil of $dg^2$ and $dF^2$.
Here $f\in C^\infty(\R)$ is 1-periodic, $q\in C^\infty([-N,N])$, and
\begin{itemize}  
\item[(i)] $f$ is  even, $f>0$ if $x\notin\frac{1}{2}{\Z}$, and        
$f(0)=f(1/2)=0$;  
\item[(ii)] $q$ is even, $q<0$ if $y\ne 0$,  $q(0)=0$ and $q^{''}(0)<0$;  
\item[(iii)] $f^{(2k)}(l/2)=(-1)^kq^{(2k)}(0)$,  $l=0,1$,    
for every natural $k\in{\N}$.  
\end{itemize}  
In particular, if  $f\sim \sum_{k=1}^{\infty}\ f_kx^{2k}$ is the Taylor expansion of $f$  
at $0$,  then,   by  (iii),  the Taylor expansion of $q$ at $0$  
is $q\sim \sum_{k=1}^{\infty}\ (-1)^k f_kx^{2k}$.  

\begin{Remark}
The branched points of the covering correspond to the points $(0,0)$ and $(1/2,0)$ of the cylinder $C$.
The metric \eqref{eq:the_metric} and the integral \eqref{eq:the _integral} on $C$ vanish at these points.
\end{Remark}

\noindent{\em Proof of Theorem \ref{th:invariant_definition}.}
Consider a tubular neighborhood $V\equiv V(\Ga)\subseteq X$ of the boundary $\Ga$ in $X$ that is 
diffeomorphic to the strip $(\R/\Z)\times[-\ep,0]$, $\ep>0$, and assume that the boundary $\Ga$ 
corresponds to the circle $(\R/\Z)\times\{0\}$. By gluing two 2-dimensional closed disks along the 
boundaries of this strip and then by extending the Riemannian metric $g$ to a smooth Riemannian metric 
$g$ on the corresponding 2-sphere, we obtain an isometrical embedding of our tubular neighborhood
$V$ of the boundary $\Ga$ into a Riemannian manifold diffeomorphic to the unit 2-sphere ${\mathbb S}^2$ 
in $\R^3$. Using the metric $g$ on ${\mathbb S}^2$ and passing to isothermal charts we obtain a 
complex atlas on ${\mathbb S}^2$, that transforms  ${\mathbb S}^2$ into a Riemann surface. 
Then, by the Riemann mapping theorem, this Riemann surface is biholomorphically equivalent to the 
standard Riemann sphere that we identify with the complex projective plane $\C P^1$. 
Taking a point $N$ on $\C P^1$ that does not lie in the image of the strip $V$
and then applying stereographic projection $\C P^1\setminus\{N\}\to\C$ we obtain an embedding
of the strip $V$ into the complex plane. By construction, the push-forward of the metric $g$ 
is conformally equivalent to the Euclidean metric on $\C$. Let $\{(x,y)\}$ 
and $z=x+i y$ be the coordinates in $\C$. For simplicity, we will identify the metric $g$, the integral
$I$, and  the neighborhood $V$ and $\Ga$ with their corresponding push-forward images. Then we have,
\begin{equation*}
dg^2=\frac{1}{2}\la(x,y)(dx^2+dy^2)
\end{equation*}
and $V$ is a closed domain in $\C$ diffeomorphic to the annulus $(\R/\Z)\times[-\ep,0]$. 
By construction, $g$ is extended to a smooth Riemannian metric 
$d{\tilde g}^2=\frac{\la(x,y)}{2}(dx^2+dy^2)$ on the whole of $\C$.
Let $\{(p_1,p_2,x,y)\}$ be the standard coordinates on $T^*\R^2$ where $\R^2$ is identified with 
$\C$. Applying the Legendre transform, then passing to complex notations and introducing 
the complex impulses $p:=\frac{1}{2}(p_1-i p_2)$, $\bar p:=\frac{1}{2}(p_1+i p_2)$, and 
$\partial_p:=\partial_{p_1}+i\partial_{p_2}$, $\partial_{\bar p}:=\partial_{p_1}-i\partial_{p_2}$,
we obtain complex coordinates $\{(p,z)\}$ on $T^*\R^2$ so that
\begin{equation}\label{eq:the_metric'}
H=\frac{2 p{\bar p}}{\la(z,{\bar z})}
\end{equation}
and
\begin{equation}\label{eq:the_integral'}
I= A(z,{\bar z})p^2+B(z,{\bar z})p{\bar p}+\overline{A(z,{\bar z})}\,{\bar p}^2
\end{equation}
with $\overline{B(z,{\bar z})}=B(z,{\bar z})$.\footnote{${\bar I}=I$ as $I$ is real-valued.} 
In view of condition $(b)$ of the theorem, the coefficient $A(z,{\bar z})$ does {\em not} vanish on 
$\Ga$, i.e., 
\begin{equation}\label{eq:condition(b)}
A(z,{\bar z})\ne 0\,\,\,\forall z\in\Ga\,.
\end{equation}
In the coordinates $\{(z,p)\}$ on $T^*\R^2$ the canonical symplectic structure $\om$ takes the 
form $\om=dp\wedge dz+d{\bar p}\wedge d{\bar z}$.
Hence,
\begin{equation}\label{eq:bracket}
\{H,I\}=(H_p I_z-H_z I_p)-(H_{\bar p} I_{\bar z}-H_{\bar z} I_{\bar p})=
2 \mathop{\tt Re}(H_p I_z-H_z I_p)\,.
\end{equation}
As $I$ is a first integral of the geodesic flow of $g$, we have
\begin{equation}\label{eq:poisson_bracket_vanish}
\{H,I\}|_V\equiv 0\,.
\end{equation}
Using \eqref{eq:bracket} one sees that  equation \eqref{eq:poisson_bracket_vanish} is equivalent to 
the following system of equations
\begin{equation}
\left\{
\begin{array}{l}\label{eq:the_system}
A_{\bar z}=0,\\
\la A_z+2 \la_z A+B\la_{\bar z}+\la B_{\bar z}=0\,.
\end{array}
\right.
\end{equation}
In particular, we see that the coefficient $A(z,{\bar z})$ in front of $p^3$ in the formula for the integral 
\eqref{eq:the_integral'} is holomorphic in $z\in V\subseteq\C$, $A=A(z)$.
Take $z_0$ in the interior of $V$ and consider a biholomorphic change of the variable $w=w(z)$ in an open
neighborhood of $z_0$ in the interior of $V$. Then the expression for the integral \eqref{eq:the_integral'} implies that,
\begin{equation}\label{eq:the_form}
{\tilde A}(w)=A(z)\Big(\frac{dw}{dz}\Big)^2,
\end{equation}
where ${\tilde A}(w)$ is the coefficient in front of $({\tilde p})^3$ in the expression for the integral
$I$ in the chart corresponding to $w$. Here ${\tilde p}$ is the complex impulse in the chart corresponding
to $w$.\footnote{Note that $p={\tilde p}\frac{dw}{dz}$.}
\begin{Remark}
In fact, \eqref{eq:the_form} implies that the bivector field,
\[
\Omega_{\{z\}}:=A(z)\,\partial_z\otimes\partial_z,
\] 
when written in an isothermal atlas, will correspond to a globally defined holomorphic section of the boundle 
$T^{2,0}_{\C} X\subseteq T_{\C}X$ . As the integral $I$ is non-trivial (condition $(b)$),
the holomorphic bivector $\Omega$ vanishes only at finitely many points in the interior of $X$. 
If $X$ were a closed surface, then by Hopf theorem, $\deg (\Omega)=2\chi(X)$, where $\chi(X)=2-2g$ is the Euler 
characteristic of $X$ and $\deg(\Omega)$ is the number of zeros (counted with multiplicities) of $\Omega$.
This would imply that $g=0, 1$, and therefore $X$ would be diffeomorphic to the 2-sphere or the 2-torus
(\cite{Kol,Koz}).
\end{Remark}
\noindent Next, we want to simplify \eqref{eq:the_system} by passing to a new complex variable $w=w(z)$,
with $w(z)$ holomorphic in the interior of $V$, so that ${\tilde A}(w)\equiv 1$. In view of \eqref{eq:the_form}, 
this amounts to  solving the differential equation $1=A(z)(\frac{dw}{dz})^2$. In the case when $A(z_0)\ne 0$ and 
$z_0$ lies in the  interior of $V$, the later equation can we solved explicitly in a sufficiently small open disk 
centered at $z_0$,
\[
w(z)=w(z_0)+\int_{z_0}^z \frac{d\la}{\sqrt{A(\la)}}\,,
\]
where the path of integration connecting $z_0$ with $z$ is $C^1$-smooth and lies in the small disk
centered at $z_0$. The square root $\sqrt{A(\la)}$ is holomorphic in the considered disk and is defined up to the
choice of the sign. As it was mentioned above, condition $(b)$ of the theorem implies \eqref{eq:condition(b)}. Hence,
by shrinking the strip $V$ if necessary, we can ensure that $\Ga\subseteq V$ and $A(z)\ne 0$ for any $z\in V$. 
Now, take $z_0\in\Ga$ and consider the map,
\begin{equation}\label{eq:Phi}
\Phi : z\mapsto w(z):=\int_{z_0}^z \frac{d\la}{\sqrt{A(\la)}},\quad V\to\C\,,
\end{equation}
where the path of integration connecting $z_0$ with $z$ is $C^1$-smooth and is contained in $V$.
Clearly, the map $\Phi$ above is well-defined on $V$ and holomorphic in the
interior of $V$. Moreover, it follows from \eqref{eq:Phi} that the directional derivatives of $\Phi$ of all orders exist and 
are continuous up to the boundary of $V$. This allows us to extend $\Phi$ to a smooth map defined in some
open set $\widetilde V\supseteq V$. Next, let us consider the image $\Phi(\Ga)$ of the boundary $\Ga$.
Take $z_1\in\Ga$. By the inverse function theorem, there exist an open neighborhood $U(z_1)$ of $z_1$ 
in $\C$ and an open neighborhood ${W}(w_1)$ of $w_1:=\Phi(z_1)$ in $\C$ so that 
$\Phi|_{U(z_1)} : U(z_1)\to W(z_1)$ is a diffeomorphism. Let $w=u+i v$. Then, as ${\tilde A}(w)=1$ for all
$w\in\Phi(V)$ we conclude from \eqref{eq:the_metric'} and \eqref{eq:the_integral'} that ${\tilde g}:=\Phi_*(g)$ and 
${\tilde I}:=\Phi_*(I)$ are diagonal in the coordinates $\{(u,v)\}$ on $W(w_1)\cap\Phi(V)$ and
non-proportional at all points of $\Phi(V)$. In other words, the coordinate vector fields $\partial_v$ and $\partial_u$
on $W(w_1)\cap\Phi(V)$ coincide with the principle directions of the quadratic forms ${\tilde g}$ and ${\tilde I}$.
As by assumption the integral ${\tilde I}$ is invariant with respect to the reflections at the boundary $\Phi(\Ga)$ we
conclude that $\Phi(\Ga)\cap W(w_1)$ is a coordinate line. As $z_1\in\Ga$ was chosen arbitrarily, we see that
$\Phi(\Ga)$ is a straight line. By shrinking the strip $V$ so that $\Ga\subseteq V$ onece more if necessary and 
by rotating the target copy of $\C$ we get that for some $\de>0$,
\[
\Phi(V)=\{w=u+ i v\,|\,-\de \le v\le 0\},
\]
and $\Phi : V\to\Phi(V)$ is a smooth covering map. The boundary $\Phi(\Ga)$ coincides with the real line.
This proves that there exist a tubular neighborhood $V$ of $\Ga$ in $X$ and $\de, l>0$ such that $V$ is 
diffeomorphic  to the cylinder 
\[
Z:=\{z=x+ i y\,|\,x\in \R/l \Z , -\de \le y\le 0\}, 
\]
with $\Ga$ corresponding to $(\R/l \Z )\times\{0\}$,
\begin{equation}\label{eq:H and I}
H=2 p{\bar p}/\la,\quad I= p^2+B p{\bar p}+{\bar p}^2,\quad {\bar B}=B,
\end{equation}
and
\begin{equation}\label{eq:B-equation}
B\la_{\bar z}+\la B_{\bar z}=-2\la_z\,.
\end{equation}
Note that \eqref{eq:B-equation} is equivalent to the vanishing of the Poisson bracket $\{H,I\}$.
Equation \eqref{eq:B-equation} is also equivalent to $(\la B)_{\bar z}=-2\la_z$. By separating its real and 
the imaginary parts we get
\[
\left\{
\begin{array}{l}
(\la B)_x=-2 \la_x,\\
(\la B)_y=2 \la_y\,.
\end{array}
\right. 
\]
This implies that 
\[
\la(B+2)=\big(\la(B+2)\big)|_{(y,0)}=q(y)/4,\,\,\,\la(B-2)=\big(\la(B-2)\big)|_{(0,x)}=f(x)/4
\]
where $f\in C^\infty(\R)$ is 1-periodic and $q\in C^\infty([-\delta,0])$. By subtracting these two equations
we see that,
\[
\la(x,y)=q(y)-f(x)>0.
\]
This, together with \eqref{eq:H and I} and $p=(p_1-i p_2)/2$ implies that
\begin{equation}\label{eq:H}
H=\frac{1}{2}\frac{p_1^2+p_2^2}{q(y)-f(x)}
\end{equation}
and
\begin{eqnarray}
2 I&=&(B+2) p_1^2+(B-2) p_2^2\nonumber\\
&=&\frac{q(y) p_1^2+f(x) p_2^2}{q(y)-f(x)}\label{eq:I}\,.
\end{eqnarray}
\begin{Remark}\label{rem:involutivity}
Our arguments also show that for any choice of functions $f\in C^\infty(\R)$ 1-periodic 
and $q\in C^\infty([a,b])$, $a, b\in\R$, so that $q(y)-f(x)>0$ the functions $H$ and $I$ defined by \eqref{eq:H} and 
\eqref{eq:I} are in involution with respect to  the canonical symplectic structure $\om=dp_1\wedge dx+dp_2\wedge dy$
on the cotangent bundle to the cylinder $(\R/l \Z)\times [a,b]$.
\end{Remark}
\noindent This Remark allows us to extend the metric $g$ and its first integral $I$ to a larger cylinder,
\[
{\widetilde Z}:=\{z=x+ i y\,|\,x\in \R/l \Z , -\de \le y\le \de\}, 
\]
that contains the boundary $\Ga=(\R/l \Z )\times\{0\}$ in its interior. In order to do this we extend the function
$q$ to a function $q\in C^\infty[-\de,\de]$ so that $q(y)-f(x)>0$ on $\widetilde Z$ and
\begin{equation}\label{eq:q_flat}
\forall k\ge 1,\,\,\,q^{(k)}(\de)=0\,.
\end{equation}
Then we use \eqref{eq:H} and  \eqref{eq:I} to extend  the metric $g$ and $I$ to smooth quadratic forms on 
${\widetilde Z}$. By Remark \ref{rem:involutivity}, $I$ continues to be a quadratic integral of the Riemannian metric 
$g$ on ${\widetilde Z}$.
In this way we extend the Riemannian manifold $(X,g)$ to a smooth Riemannian manifold $({\widetilde X},{\tilde g})$
with  connected boundary $\widetilde\Ga$, so that 
$X\subseteq{\widetilde X}$, $\Ga$ is in the interior of ${\widetilde X}$,
${\tilde g}|_X=g$, ${\tilde I}|_X=I$, and $\tilde I$ is a quadratic first integral of ${\tilde g}$.
In addition, a collar neighborhood of $\widetilde\Ga$ in ${\widetilde X}$ can be coordinatized  by the
cylinder ${\widetilde Z}$ so that the Legendre transforms of the metric and the integral are given by 
\eqref{eq:H} and \eqref{eq:I}.

Our final step is to take two copies of $({\widetilde X},{\tilde g})$ and glue them along their boundaries
by a diffeomorphism that,  in the coordinates $\{(x,y)\}$, corresponds to the identity,
\[
(x,\de)\mapsto(x,\de),\quad(\R/l \Z)\times\{\de\}\to (\R/l \Z)\times\{\de\}\,.
\]
In this way we obtain a closed Riemannian manifold $({\hat X}, {\hat g})$. In view of the flatness condition 
\eqref{eq:q_flat} on $q$, the metric $\hat g$ and the corresponding quadratic form $\hat I$ are smooth, 
the Riemannian manifold  $(X,g)$ is isometrically embedded into $({\hat X}, {\hat g})$, and ${\hat I}|_X=I$.
Moreover, by construction, ${\hat I}$ is a quadratic integral of the geodesic flow of ${\hat g}$.
Finally, Theorem \ref{th:invariant_definition} follows from the classification theorem for
Liouville surfaces (see e.g. \cite{Kol}).
\noindent $\Box$\\

\subsection{Kolmogorov Nondegeneracy  of  the bouncing ball map for Liouville billiards}\label{sec:Appendix B2}

In this Appendix we show that the Poincar\'e map of the Liouville billiard tables on 
the surfaces of constant curvature is non-degenerate at the elliptic fixed point.

\medskip

Let $(X,g)$ be a Liouville billiard table of classical type.
Then there exists a double covering with two branched points
\begin{equation}\label{eq:tau}
\tau : C\to X
\end{equation}
where $C$ denotes the cylinder $\big(\R/\Z\big)\times [-N,N]$, $N>0$, coordinatized by the variables
$x$ and $y$ respectively, so that the pull-back of the Riemannian metric on $X$ and 
the corresponding quadratic in velocities first integral take the form
\begin{equation}\label{eq:metric}
dg^2=\big(f(x)-q(y)\big)(dx^2+dy^2)
\end{equation}
\begin{equation}\label{eq:integral}
dF^2=\big(f(x)-q(y)\big)\big(q(y)dx^2+f(x)dy^2\big)
\end{equation}
where $f\in C^\infty(\R)$ is 1-periodic, $q\in C^\infty\big([-N,N]\big)$, and
the hypotheses (i)$\div$(v) in the definition of Liouville billiard tables of classical type hold.
In addition, we will assume that $f$ has a {\em Morse singularity} at $x=1/4$ which amounts 
to $f''(1/4)<0$. Note that the line (taken twice) on the cylinder $C$ corresponding to $x=1/4$ is an elliptic 
closed broken geodesic of $(X,g)$ with two vertices. Let 
\begin{equation}\label{eq:f_Taylor}
f(x)=\alpha_0+\alpha_1(x-x_0)^2+\alpha_2(x-x_0)^4+O\big((x-x_0)^6\big)
\end{equation}
where $\alpha_0>0$ and $\alpha_1<0$ be the Taylor's expansion of $f$ at $x_0=1/4$.
Let $\{(I,\theta)\}$ be action-angle variables in an open neighborhood in ${\bf B}^*\Gamma$ of the elliptic fixed point of 
the billiard ball map of $(X,g)$ normalized so that $I=0$ at the elliptic point. 
We have the following

\begin{Theorem}\label{th:hessian_elliptic_point}
	Denote by $K$ the Hamiltonian that generates the billiard ball map in the action-angle coordinates $\{(I,\theta)\}$.
	Then
	\[
	\frac{dH}{dI}(0)=-\frac{\sqrt{-\alpha_1}}{\pi}\int_{-N}^N\frac{dy}{\sqrt{\alpha_0-q(y)}}
	\]
	and
	\[
	\frac{d^2 K}{dI^2}(0)=\frac{\alpha_1}{4\pi^2}\left(2\int_{-N}^N\frac{dy}{\big(\alpha_0-q(y)\big)^{3/2}}-
	\frac{3\alpha_2}{\alpha_1^2}\int_{-N}^N\frac{dy}{\sqrt{\alpha_0-q(y)}}\right).
	\]
\end{Theorem}

Integrable billiard tables on surfaces of constant curvature are examples of Liouville billiard tables of classical type -- 
see \cite[\S 3]{PT1}. In the case of elliptic billiard tables we have that
\begin{equation}\label{eq:ellipse}
f(x)=4\epsilon^2\pi^2\sin^22\pi x\quad\text{and}\quad q(y)=-4\epsilon^2\pi^2\sinh^22\pi x
\end{equation}
where $\epsilon>0$ is the distance between the center of the ellipse and one of the focuses (see \cite[\S 3.1]{PT1}).
As a consequence we obtain

\begin{Coro}\label{coro:ellipse}
	For any $\epsilon>0$ and for any $N>0$ we have that $-1<\frac{dK}{dI}(0)<0$ and $\frac{d^2 K}{dI^2}(0)<0$. In particular,
	the Poincar\'e map of the elliptic billiard ball map is {\em non-degenerate} (twisted) at the elliptic fixed point. Moreover, 
	it is {\em $4$-elementary} except for five different values of the parameter $N>0$.
\end{Coro}

\begin{Remark}
	Similar results can be proved for the Liouville billiard tables on the surfaces of constant curvature.
\end{Remark}

\begin{proof}[Proof of Theorem \ref{th:hessian_elliptic_point}]
	Let $\{(x,y,p_1,p_2)\}$ be the standard coordinates on the cotangent bundle $T^*C$.
	By the Legendre transform 
	\[
	H=\frac{p_1^2+p_2^2}{f(x)-q(y)},\quad F=\frac{q(y)p_1^2+f(x)p_2^2}{f(x)-q(y)}
	\]
	where 
	\begin{equation}\label{eq:impulses}
	p_1=\big(f(x)-q(y)\big){\dot x},\quad p_2=\big(f(x)-q(y)\big){\dot y}
	\end{equation}
	and $\dot x$ and $\dot y$ denote the components of the velocity vectors in $TC$.
	For 
	\[
	0<h\le \max f=\alpha_0
	\] 
	consider the invariant with respect to the geodesic flow on $T^*C$ surface 
	\[
	Q_h:=\big\{H=1, F=h\big\}\subseteq T^*C.
	\]
	Since the variables separate one easily sees that $Q_h$ is characterized by the set of equations
	\begin{equation}\label{eq:alternative}
	p_1^2=f(x)-h,\quad p_2^2=h-q(y).
	\end{equation}
	One concludes from \eqref{eq:alternative} and the hypothesis (i)$\div$(v) on the functions $f$ and
	$q$ that for $0<h<\alpha_0$ the surface $Q_h$ consists of two copies of $(\R/\Z)\times[-N,N]$.
	The billiard reflection map at the boundary of $C$ preserves the boundary of $Q_h$ and can be used
	to ``glue'' the two components of $Q_h$ into a single Liouville torus ${\widetilde Q}_h$ of the broken geodesic flow on $C$.
	By integrating the Liouville form $\kappa=p_1dx+p_2dy$ along the two cycles on ${\widetilde Q}_h$ that are ``parallel'' to
	the coordinate lines on $C$, we obtain from \eqref{eq:alternative} that
	\begin{equation}\label{eq:K}
	K(h)=2\int_{-N}^N\sqrt{h-q(y)}\,dy
	\end{equation}
	and
	\begin{equation}\label{eq:I}
	I(h)=2\int_{x_h'}^{x_h''}\sqrt{f(x)-h}\,dx
	\end{equation}
	where  $0<x_h'\le 1/4\le x_h''<1/2$ are the two zeros of the equation $f(x)=h$.
	Note that $x_h'=x_h''=1/4$ if and only if $h=\alpha_0$.
	Since $f$ has a Morse singularity at $x_0=1/2$ there exists an orientation preserving
	change of variables $p : U(0)\to V(0)$ from an open neighborhood of zero $U(0)$ 
	onto an open neighborhood of zero $V(0)$ such that
	\begin{equation}\label{eq:u}
	x-x_0=p(u),
	\end{equation}
	and
	\begin{equation}\label{eq:morse}
	f(x)-h=(\alpha_0-h)-u^2.
	\end{equation}
	It follows directly from \eqref{eq:f_Taylor}, \eqref{eq:u}, and \eqref{eq:morse} that
	\begin{equation}\label{eq:p_Taylor}
	p(u)=\frac{1}{\sqrt{-\alpha_1}}\,y+\frac{\alpha_2}{2\alpha_1^2\sqrt{-\alpha_1}}\,y^3+O(y^5).
	\end{equation}
	In view of \eqref{eq:I}, \eqref{eq:u}, \eqref{eq:morse}, and \eqref{eq:p_Taylor} we obtain
	\begin{eqnarray*}
		I(h)&=&2\int_{x_h'}^{x_h''}\sqrt{f(x)-h}\,dx\\
		&=&2\int_{-\sqrt{\alpha_0-h}}^{\sqrt{\alpha_0-h}}p'(u)\sqrt{(\alpha_0-h)-y^2}\,dy\\
		&=&2(\alpha_0-h)\int_{-1}^1p'\big(u\sqrt{\alpha_0-h}\big)\sqrt{1-u^2}\,du\\
		&=&-\frac{\pi}{\sqrt{-\alpha_1}}(h-\alpha_0)+\frac{3\alpha_2\pi}{8\alpha_1^2\sqrt{-\alpha_1}}(h-\alpha_0)^2
		+O\big((h-\alpha_0)^3\big).
	\end{eqnarray*}
	Hence,
	\begin{equation}\label{eq:I-derivatives}
	I(\alpha_0)=0,\quad\frac{dI}{dh}(\alpha_0)=-\pi/\sqrt{-\alpha_1},\quad
	\frac{d^2I}{dh^2}(\alpha_0)=3\pi\alpha_2/4\alpha_1^2\sqrt{-\alpha_1}.
	\end{equation}
	It follows from \eqref{eq:K} and the fact that $\alpha_0>0$ that $K$ is a $C^\infty$-smooth function of $h$ in an open
	neighborhood of $h=\alpha_0$. By combining this with \eqref{eq:I-derivatives} we obtain
	\begin{equation}\label{eq:rotation_limit}
	\frac{dK}{dI}(0)=\frac{dK}{dh}(\alpha_0)\Big\slash\frac{dI}{dh}(\alpha_0)=
	-\frac{\sqrt{-\alpha_1}}{\pi}\int_{-N}^N\frac{dy}{\sqrt{\alpha_0-q(y)}}
	\end{equation}
	and
	\begin{eqnarray}
	\frac{d^2K}{dI^2}(0)\!\!&=&\!\!\left( \frac{d^2K}{dh^2}(\alpha_0)-\frac{dK}{dI}(0)\frac{d^2I}{dh^2}(\alpha_0)\right)
	\Big\slash\Big(\frac{dI}{dh}(\alpha_0)\Big)^2\nonumber\\
	&=&\!\!\frac{\alpha_1}{4\pi^2}\left(2\int_{-N}^N\frac{dy}{\big(\alpha_0-q(y)\big)^{3/2}}-
	\frac{3\alpha_2}{\alpha_1^2}\int_{-N}^N\frac{dy}{\sqrt{\alpha_0-q(y)}}\right).\label{eq:rotation'_limit}
	\end{eqnarray}
\end{proof}

\begin{proof}[Proof of Corollary \ref{coro:ellipse}]
	The Corollary follows directly from Theorem \ref{th:hessian_elliptic_point} and \eqref{eq:ellipse}. In fact,
	it follows from \eqref{eq:ellipse} that 
	\[
	f(x)=4\epsilon^2\pi^2\Big(1-4\pi^2(x-1/4)^2+\frac{16\pi^4}{3}(x-1/4)^4+O\big((x-1/4)^6\big)\Big).
	\]
	Hence,
	\begin{equation}\label{eq:coefficients_ellipse}
	\alpha_0=c^2,\quad\alpha_1=-4\pi^2 c^2,\quad\alpha_2=\frac{16\pi^4}{3}\,c^2,
	\end{equation}
	where we set for simplicity $c:=2\pi\epsilon$.
	Then, in view of Theorem \ref{th:hessian_elliptic_point} we obtain that
	\begin{equation*}
	\frac{dK}{dI}(0)=-2\int_{-N}^N\frac{dy}{\sqrt{1+\sinh^22\pi y}}
	\end{equation*}
	and
	\begin{equation*}
	\frac{d^2 K}{dI^2}(0)=-\frac{1}{2\epsilon\pi}\left(2\int_{-N}^N\frac{dy}{\big(1+\sinh^22\pi y\big)^{3/2}}-
	\int_{-N}^N\frac{dy}{\sqrt{1+\sinh^22\pi y}}\right).
	\end{equation*}
	By passing to the variable $v=\sinh2\pi y$ in the integrals above one obtains that
	\begin{equation}\label{eq:rotation_ellipse'}
	\frac{dK}{dI}(0)=-\frac{2}{\pi}\arctan\big(\sinh 2\pi N\big)<0
	\end{equation}
	and 
	\begin{equation}\label{eq:rotation'_ellipse'}
	\frac{d^2 K}{dI^2}(0)=-\frac{1}{2\epsilon\pi^2}\,\frac{\sinh 2\pi N}{\cosh^2 2\pi N}<0
	\end{equation}
	which completes the proof of the first two statements of the Corollary. 
	It is clear that the spectrum of the Poincar\'e map of the elliptic billiard ball map at the elliptic fixed point
	is equal to $\{e^{\pm i\varphi}\}$ where
	\[
	\varphi\equiv\pm 2\pi\frac{dK}{dI}(0)\!\mod\pi.
	\]
	The last statement of the Corollary then follows from the definition of the $4$-elementary Poincar\'e maps and
	formula \eqref{eq:rotation_ellipse'}.
\end{proof}

% -----------------------------------------------------------------------------------------------------------------

\noindent
\textbf{Acknowledgments.}
{\em Part of this work has been done at the Institute of Mathematics, Bulgarian Academy of Sciences, and we would like to thank the colleagues there for the stimulating discussions. }

\vspace{0.5cm} 
\noindent 
G. P.: 
Universit\'e de Nantes,  \\
Laboratoire de math\'ematiques Jean Leray,
CNRS: UMR 6629,\\
2, rue de la Houssini\`ere, \\
BP 92208,  44072 Nantes 
Cedex 03, France \\

\vspace{0.5cm} 
\noindent 
P.T.: Department of Mathematics,\\
 Northeastern University,\\ 
360 Huntington Avenue, Boston, MA 02115

\end{document}